\documentclass[12pt]{amsart}
\usepackage{amsfonts,amssymb,amscd,amsmath,epsfig}
\usepackage{latexsym}
\usepackage{mathrsfs}
\usepackage{tocvsec2}

\usepackage[]{youngtab}

\usepackage[
hyperindex=true,pagebackref=true,bookmarks=true,
colorlinks=true,linkcolor=blue,citecolor=red]{hyperref}
\def\~{{\rm --}}

\begin{document}

\renewcommand{\tilde}{\widetilde}
\renewcommand{\hat}{\widehat}

\newcommand{\BR}{{\mathbb R}}
\newcommand{\BQ}{{\mathbb Q}}
\newcommand{\BC}{{\mathbb C}}
\newcommand{\BP}{{\mathbb P}}
\newcommand{\BZ}{{\mathbb Z}}
\newcommand{\BN}{{\mathbb N}}
\newcommand{\BS}{{\mathbb S}}

\newcommand{\cH}{{\mathcal H}}
\newcommand{\cA}{{\mathcal A}}
\newcommand{\cB}{{\mathcal B}}
\newcommand{\ccF}{{\mathfrak F}}
\newcommand{\cD}{{\mathcal D}}
\newcommand{\cL}{{\mathcal L}}
\newcommand{\cF}{{\mathcal F}}
\newcommand{\cP}{{\mathcal P}}
\newcommand{\cX}{{\mathcal X}}
\newcommand{\cY}{{\mathcal Y}}
\newcommand{\cS}{{\mathcal S}}
\newcommand{\cSol}{\hbox{$\mathcal Sol$}}
\newcommand{\cT}{\hbox{$\mathcal T$}}

\newcommand{\Z}{{\mathbb Z}}
\newcommand{\Q}{{\mathbb Q}}
\newcommand{\N}{{\mathbb N}}
\newcommand{\C}{{\mathbb C}}
\newcommand{\R}{{\mathbb R}}
\newcommand{\X}{{\mathbb X}}
\newcommand{\Y}{{\mathbb Y}}

\newcommand{\CH}{{\mathcal H}}
\newcommand{\CA}{{\mathcal A}}

\def\HH{\mbox{${\mathcal H}$\kern-5.2pt${\mathcal H}$}}

\newcommand{\binomial}[2]{\genfrac{(}{)}{0pt}{}{ #1 }{ #2 }}
\newcommand{\qbinomial}[2]{\genfrac{[}{]}{0pt}{}{ #1 }{ #2 }_q }
\newcommand{\qbinom}[3]{\genfrac{[}{]}{0pt}{}{ #1 }{ #2 }_{ #3 } }


\def\der{\partial}
\def\tensor{\otimes}
\def\gam{\gamma} \def\Gam{\Gamma}
\def\del{\delta} \def\Del{\Delta}
\def\kap{\kappa}
\def\lam{\lambda} \def\Lam{\Lambda}
\def\Comp{{\mathbb C}}
\def\sM{{\mathcal M}}

\newtheorem{theorem}{Theorem}[section]
\newtheorem{maintheorem}[theorem]{Main Theorem}
\newtheorem{proposition}[theorem]{Proposition}
\newtheorem{definition}[theorem]{Definition}
\newtheorem{lemma}[theorem]{Lemma}
\newtheorem{corollary}[theorem]{Corollary}
\newtheorem{notation}[theorem]{Notation}
\newtheorem{remark}[theorem]{Remark}
\newtheorem{example}[theorem]{Example}

\newtheorem{theorem }{Theorem}[section]
\newtheorem{maintheorem }[theorem]{Main Theorem}
\newtheorem{proposition }[theorem]{Proposition}
\newtheorem{definition }[theorem]{Definition}
\newtheorem{lemma }[theorem]{Lemma}
\newtheorem{corollary }[theorem]{Corollary}
\newtheorem{notation }[theorem]{Notation}
\newtheorem{remark }[theorem]{Remark}
\newtheorem{example }[theorem]{Example}

\newtheorem{ maintheorem }[theorem]{Main Theorem}
\newtheorem{ theorem}{Theorem}[section]
\newtheorem{ proposition}[theorem]{Proposition}
\newtheorem{ definition}[theorem]{Definition}
\newtheorem{ lemma}[theorem]{Lemma}
\newtheorem{ corollary}[theorem]{Corollary}
\newtheorem{ notation}[theorem]{Notation}
\newtheorem{ remark}[theorem]{Remark}
\newtheorem{ example}[theorem]{Example}

\newtheorem{thm}{Theorem}[section]
\newtheorem{prop}[thm]{Proposition}
\newtheorem{lem}[thm]{Lemma}
\newtheorem{cor}[thm]{Corollary}
\newtheorem{conj}[thm]{Conjecture}
\newtheorem{con}[thm]{Conjecture}
\newtheorem{dfn}[thm]{Definition}
\newtheorem{df}[thm]{Definition}
 \newcommand{\rem}{{\bf Comment.\ }}
 \newcommand{\rmk}{{\bf Comment.\ }}
 \newcommand{\exmp}{{\bf Example.\ }}
 \newcommand{\ex}{{\bf Example.\ }}
 \newcommand{\prob}{{\bf Problem.\ }}

\newtheorem{note}{Note} 
\renewcommand{\thenote}{}
\newtheorem*{acka}{Acknowledgments}
\newtheorem{ack}{Acknowledgments}
\renewcommand{\theack}{}
\renewcommand{\appendixname}{\bf Appendix}
\renewcommand{\proof}{{\em Proof.\ }}

\hyphenation{
ap-pen-dix as-ymp-tot-ic at-trib-uted at-trib-ut-able
Bry-li-n-sky com-mu-ta-tion de-ge-ne-rate
de-riv-a-tive dis-trib-ute equi-vari-ant ex-tra-or-di-nary  
geo-met-ric griev-ance griev-ous grad-ed ho-lo-no-my ho-mo-thetic
in-fin-ite-ly in-fin-i-tes-i-mal Ha-rish Cha-n-dra mul-ti-plic-able 
non-euclid-ean non-iso-mor-phic non-smooth par-a-digm 
par-a-bol-ic pa-rab-o-loid pa-ram-e-trize phe-nom-e-non 
post-script pseu-do-dif-fer-en-tial pseu-do-fi-nite 
qua-drat-ics quad-ra-ture Han-kel rec-tan-gle semi-def-i-nite 
set-up wide-spread Euler-ian Feb-ru-ary Gauss-ian Grothen-dieck 
Hamil-ton-ian Her-mi-t-ian her-mi-t-ian Jan-u-ary 
Japan-ese Ka-shi-wa-ra Kor-te-weg Le-gendre No-vem-ber Rie-mann-ian 
Sep-tem-ber Za-mo-lo-d-chi-kov Kni-zh-nik quan-tum Op-dam
Mac-do-nald Ca-lo-ge-ro Su-ther-land Mo-ser 
Ol-sha-net-sky  Pe-re-lo-mov in-de-pen-dent ope-ra-tors 
cy-clo-to-mic ra-tio-nal de-gen-er-a-tion 
in-ter-est-ing de-for-ma-tions de-for-ma-tion pro-ce-dure 
fol-lows ope-ra-tors  pre-serve suf-fices ap-proach 
for-mu-las con-sider its com-ple-tion cor-re-spond-ing 
au-to-mor-phism be-cause pro-por-tional fi-nal-ly let-ting 
equi-v-a-lence ge-n-er-al-ized Mac-do-nald iden-ti-ties 
cor-re-s-pond sub-dia-grams par-ti-tion na-t-u-ral-ly 
or-dered stan-dard de-for-ma-tion ar-gu-ment com-bined 
sphe-r-i-cal rep-re-sen-ta-tions tri-go-no-me-t-ric
ge-n-er-al-ly speak-ing pri-m-it-ive ir-re-du-cible 
sum-ma-tion  rep-re-sen-ta-tives pro-por-ti-o-na-li-ty
ultra-sphe-ri-cal Ro-gers}

\def\ffor{\quad\hbox{ for }\quad}
\def\wwhen{\quad\hbox{ when }\quad}
\def\wwhere{\quad\hbox{ where }\quad}
\def\aand{\quad\hbox{ and }\quad}
\def\for{\  \hbox{ for } \ }
\def\iif{ \ \hbox{ if } \ }
\def\when{ \ \hbox{ when } \ }
\def\where{\  \hbox{ where } \ }
\def\and{\  \hbox{ and } \ }
\def\and{\  \hbox{ and } \ }
\def\oor{\  \hbox{ or } \ }
\def\proof{{\em Proof. \  }}

\def\equal{\stackrel{\,\mathbf{def}}{= \kern-3pt =}}

\def\la{\lambda}
\def\La{\Lambda}
\def\om{\omega}
\def\Om{\Omega}
\def\Th{\Theta}
\def\th{\theta}
\def\al{\alpha}
\def\be{\beta}
\def\ga{\gamma}
\def\ep{\epsilon}
\def\up{\upsilon}
\def\Up{\Upsilon}
\def\de{\delta}
\def\De{\Delta}
\def\ka{\kappa}
\def\kapp{\hbox{\bf \ae}}
\def\si{\sigma}
\def\Si{\Sigma}
\def\Ga{\Gamma}
\def\ze{\zeta}
\def\io{\iota}
\def\bio{b^\iota}
\def\aio{a^\iota}
\def\twio{\tilde{w}^\iota}
\def\hwio{\hat{w}^\iota}
\def\gio{\g^\iota}
\def\Bio{B^\iota}

\def\del{\delta}
\def\pa{\partial}
\def\vp{\varphi}
\def\ve{\varepsilon}
\def\inf{\infty}

\def\vph{\varphi}
\def\vps{\varpsi}
\def\vPh{\varPhi}
\def\vep{\varepsilon}
\def\vpi{{\varpi}}
\def\vth{{\vartheta}}
\def\vsi{{\varsigma}}
\def\vrh{{\varrho}}

\def\bph{\bar{\phi}}
\def\bsi{\bar{\si}}
\def\bvp{\bar{\varphi}}

\newcommand{\bS}{{\mathbf S}}
\newcommand{\bH}{{\mathbf H}}
\newcommand{\bF}{{\mathbf F}}
\newcommand{\bE}{{\mathbf E}}

\def\tal{\tilde{\alpha}}
\def\tbe{\tilde{\beta}}
\def\tde{\tilde{\delta}}
\def\tpi{\tilde{\pi}}
\def\txi{\tilde{\xi}}
\def\tPi{\tilde{\Pi}}
\def\tPhi{\tilde{\Phi}}
\def\tV{\tilde{V}}
\def\tJ{\tilde{J}}
\def\tla{\tilde{\lambda}}
\def\tga{\tilde{\gamma}}
\def\tGa{\tilde{\Gamma}}
\def\tvs{\tilde{{\varsigma}}}
\def\tu{\tilde{u}}
\def\tU{\tilde{U}}
\def\tw{\widetilde w}
\def\tW{\widetilde W}
\def\tB{\tilde B}
\def\tv{\tilde v}
\def\tV{\tilde V}
\def\tz{\tilde z}
\def\tb{\tilde b}
\def\ta{\tilde a}
\def\tih{\tilde h}
\def\trh{\tilde {\rho}}
\def\tx{\tilde x}
\def\tf{\tilde f}
\def\tg{\tilde g}
\def\tG{\tilde G}
\def\tk{\tilde k}
\def\tl{\tilde l}
\def\tL{\tilde L}
\def\tD{\tilde D}
\def\tR{\tilde R}
\def\tP{\tilde P}
\def\tH{\tilde H}
\def\tp{\tilde p}

\def\hH{\hat{H}}
\def\hh{\hat{h}}
\def\hR{\hat{R}}
\def\hY{\hat{Y}}
\def\hX{\hat{X}}
\def\hP{\hat{P}}
\def\hT{\hat{T}}
\def\hV{\hat{V}}
\def\hG{\hat{G}}
\def\hF{\hat{F}}
\def\hw{\widehat{w}}
\def\hW{\widehat{W}}
\def\hu{\hat{u}}
\def\hs{\hat{s}}
\def\hv{\hat{v}}
\def\hb{\hat{b}}
\def\hB{\widehat{B}}
\def\hze{\hat{\zeta}}
\def\hsi{\hat{\sigma}}
\def\hrh{\hat{\rho}}
\def\hth{\hat{\theta}}
\def\hy{\hat{y}}
\def\hx{\hat{x}}
\def\hz{\hat{z}}
\def\hg{\hat{g}}
\def\he{\hat{e}}
\def\hE{\widehat{E}}

\def\B{\mathbf{B}}
\def\I{\mathbf{I}}
\def\P{\mathbf{P}}
\def\G{\mathbf{G}}
\def\S{\mathbf{S}}
\def\F{\mathbf{F}}
\def\one{\mathbf{1}}
\def\Sn{\mathbf{S}_n}
\def\0{\mathbf{0}}
\def\H{\mathbf{H}}
\def\V{\mathbf{V}}

\def\f{\mathcal{F}}
\def\çF{\mathcal{F}}
\def\o{\mathcal{O}}
\def\t{\mathcal{T}}
\def\r{\mathcal{R}}
\def\l{\mathcal{L}}
\def\m{\mathcal{M}}
\def\k{\mathcal{K}}
\def\n{\mathcal{N}}
\def\d{\mathcal{D}}
\def\p{\mathcal{P}}
\def\cP{\mathcal{P}}
\def\a{\mathcal{A}}
\def\h{\mathcal{H}}
\def\c{\mathcal{C}}
\def\y{\mathcal{Y}}
\def\e{\mathcal{E}}
\def\v{\mathcal{V}}
\def\z{\mathcal{Z}}
\def\x{\mathcal{X}}
\def\s{\mathcal{S}}
\def\g{\mathcal{G}}
\def\u{\mathcal{U}}
\def\w{\mathcal{W}}
\def\i{\mathcal{I}}
\def\j{\mathcal{J}}
\def\b{\mathcal{B}}

\def\lan{\langle}
\def\llb{(\!(}
\def\ran{\rangle}
\def\rrb{)\!)}
 \def\dim{{\hbox{\rm dim}}_{\mathbb C}\,}
\def\lng{\hbox{\rm{\tiny lng}}}
\def\sht{\hbox{\rm{\tiny sht}}}
\def\sph{\hbox{\rm{\tiny sph}}}
\def\inv{\hbox{\rm{\tiny inv}}}

\def\br#1{\langle #1 \rangle}

\def\rank{\hbox{rank}}
\def\gl{\mathfrak{gl}_N}

\newcommand{\Aut}{\operatorname{Aut}}
\newcommand{\Hom}{\operatorname{Hom}}
\newcommand{\End}{\operatorname{End}}
\newcommand{\Ind}{\operatorname{Ind}}
\newcommand{\ad}{\operatorname{ad}}
\newcommand{\pr}{\operatorname{pr}}
\newcommand{\aweyl}{\tilde{\mathbb S}_n}
\newcommand{\hec}{{\mathcal H}^t_n}
\newcommand{\Func}{{\mathcal F}({\mathbb C}^n,{\mathcal H}^t_n)}
\newcommand{\tr}{\operatorname{tr}}
\newcommand{\Out}{\operatorname{Out}}
\newcommand{\Rad}{\operatorname{Rad}}
\newcommand{\Spec}{\operatorname{Spec}}
\newcommand{\id}{\operatorname{id}}
\newcommand{\Int}{\operatorname{Int}}
\newcommand{\ct} {\operatorname{ct}}

\newcommand{\rat}{{\mathbb Q}}
\newcommand{\real}{{\mathbb R}}
\newcommand{\cplx}{{\mathbb C}}
\newcommand{\zint}{{\mathbb Z}}

\newcommand{\sq}{\phantom{1}\hfill$\qed$}
\newcommand{\Rea}{\Re}
\newcommand{\Ima}{\Im}

\newcommand{\st}{\bowtie}
\newcommand{\modd}{\mbox{\,mod\,}}
\newcommand{\lr}{\langle}
\newcommand{\rr}{\rangle}
\newcommand{\eps}{\varepsilon}
\newcommand{\phk}{\phi^{(k)}}
\newcommand{\psk}{\psi^{(k)}}
\newcommand{\Res}{\mbox{Res}\;}
\newcommand{\sgn}{\mbox{sgn}}
\newcommand{\mn} {\left\{ \begin{array}{c}m\\
n\end{array}\right\}}

\def\sX{\mathscr{X}}
\def\sH{\mathscr{H}}
\def\sY{\mathscr{Y}}
\def\TT{\mathfrak{T}}
\def\JJ{\mathfrak{J}}
\def\HH{\mathfrak{H}}
\def\FF{\mathfrak{F}}
\def\GG{\mathfrak{G}}
\def\CC{\mathfrak{C}}
\def\LL{\mathfrak{L}}

\def\BB{\mathfrak{B}}
\def\AA{\mathfrak{A}}
\def\ZZ{\mathfrak{Z}}
\def\HH{\hbox{${\mathcal H}$\kern-5.2pt${\mathcal H}$}}
\def\HHH{\hbox{${\mathbb H}$\kern-4.2pt${\mathbb H}$}}
\def\tHH{\widetilde{\HH\ }}

\font\smm=msbm10 at 12pt 
\def\symbol#1{\hbox{\smm #1}}
\def\lsmash{{\symbol n}}
\def\rsmash{{\symbol o}}
\def\#{\sharp}

\font\tenbf=cmbx10
\font\tenrm=cmr10
\font\tenit=cmti10
\font\ninebf=cmbx9
\font\ninerm=cmr9
\font\nineit=cmti9
\font\eightbf=cmbx8
\font\eightrm=cmr8
\font\eightit=cmti8
\font\sevenrm=cmr7
\font\sevenbf=cmbx7


\title [DAHA approach to iterated torus links]
{DAHA approach to iterated torus links}
\author[Ivan Cherednik]{Ivan Cherednik $^\dag$}
\author[Ivan Danilenko] {Ivan Danilenko}

\begin{abstract}
We extend the construction of DAHA-Jones polynomials for
any reduced root systems and
DAHA-superpolynomials in type $A$  from  
iterated torus knots (our previous paper) to links, including 
arbitrary algebraic links. Such a passage essentially corresponds 
to the usage of the products of Macdonald polynomials
and is directly connected to the splice diagrams (Neumann et al.).
The specialization $t=q$ of our superpolynomials 
results in the HOMFLY-PT polynomials. 
The relation of our construction to the stable Khovanov-Rozansky 
polynomials and the so-called ORS-polynomials of the corresponding
plane curve singularities is expected for algebraic links
in the uncolored case. These 2 connections are less
certain, since the Khovanov-Rozansky theory for links is
not sufficiently developed and the ORS polynomials are quite
involved. However we provide some confirmations. 
For Hopf links, our 
construction produces the DAHA-vertex, similar to the
refined topological vertex, which is an important part of our 
work.
\end{abstract}

\thanks{$^\dag$ \today.
\ \ \ Partially supported by NSF grant
DMS--1363138}

\address[I. Cherednik]{Department of Mathematics, UNC
Chapel Hill, North Carolina 27599, USA\\
chered@email.unc.edu}

\address[I. Danilenko]{Department of Mathematics, Columbia
University, 2990 Broadway, New York, NY 10027, USA\\
danilenko@math.columbia.edu\\
ITEP, 25 Bolshaya Cheremushkinskaya, Moscow 117218, Russia}

 \def\sht{\raisebox{0.4ex}{\hbox{\rm{\tiny sht}}}}
 \def\bysame{{\bf --- }}
 \def\~{{\bf --}}
 \def\rr{{\mathsf r}}
 \def\ss{{\mathsf s}}
 \def\mm{{\mathsf m}}
 \def\pp{{\mathsf p}}
 \def\ll{{\mathsf l}}
 \def\aa{{\mathsf a}}
 \def\bb{{\mathsf b}}
 \def\NS{\hbox{\tiny\sf ns}}
 \def\ssum{\hbox{\small$\sum$}}
\newcommand{\comment}[1]{}
\renewcommand{\tilde}{\widetilde}
\renewcommand{\hat}{\widehat}
\renewcommand{\V}{\mathbb{V}}
\renewcommand{\S}{\mathbb{S}}
\renewcommand{\F}{\mathbb{F}}
\newcommand{\dagx}{\hbox{\tiny\mathversion{bold}$\dag$}}
\newcommand{\ddagx}{\hbox{\tiny\mathversion{bold}$\ddag$}}
\newtheorem{conjecture}[theorem]{Conjecture}
\newcommand*\toeq{
\raisebox{-0.15 em}{\,\ensuremath{
\xrightarrow{\raisebox{-0.3 em}{\ensuremath{\sim}}}}\,}
}
\newcommand{\unknot}{\hbox{\tiny\!\raisebox{0.2 em}
{$\bigcirc$}}\!}
\newcommand{\mmu}{\hbox{\mathversion{bold}$\mu$}}
\newcommand{\lla}{\hbox{\mathversion{bold}$\lambda$}}
\newcommand{\dde}{\hbox{\mathversion{bold}$\delta$}}

\newcommand\rightthreearrow{\hbox{\tiny
        $\mathrel{\vcenter{\mathsurround0pt
         \ialign{##\crcr
         \noalign{\nointerlineskip}$\rightarrow$\crcr
         \noalign{\nointerlineskip}$\rightarrow$\crcr
         \noalign{\nointerlineskip}$\rightarrow$\crcr
                }}}$ }}
\newcommand\rightfourarrow{\hbox{\tiny
        $\mathrel{\vcenter{\mathsurround0pt
         \ialign{##\crcr
         \noalign{\nointerlineskip}$\rightarrow$\crcr
         \noalign{\nointerlineskip}$\rightarrow$\crcr
         \noalign{\nointerlineskip}$\rightarrow$\crcr
         \noalign{\nointerlineskip}$\rightarrow$\crcr
                }}}$ }}

\newcommand\rightdotsarrow{\hbox{\small
        $\mathrel{\vcenter{\mathsurround0pt
         \ialign{##\crcr
         \noalign{\nointerlineskip}$\,\rightarrow$\crcr
         \noalign{\nointerlineskip}$\cdots$\crcr
         \noalign{\nointerlineskip}$\,\rightarrow$\crcr
                }}}$ }}
\newcommand\rightdotsarrowtiny{\hbox{\tiny
        $\mathrel{\vcenter{\mathsurround0pt
         \ialign{##\crcr
         \noalign{\nointerlineskip}$\,\rightarrow$\crcr
         \noalign{\nointerlineskip}$\cdots$\crcr
         \noalign{\nointerlineskip}$\,\rightarrow$\crcr
                }}}$ }}

\newcommand{\twoone}
{\hbox{\rm
$\circ\!$\raisebox{-2.6pt}{$\rightarrow$}
\raisebox{-2.6pt}{$\!\!\circ\!\!\rightarrow$}
\kern-33pt\raisebox{+2.6pt}{$\rightarrow$}
\kern+14pt}}
\newcommand{\twoonetiny}
{\hbox{\tiny
$\circ\!$\raisebox{-2.pt}{$\rightarrow$}
\raisebox{-2.pt}{$\!\!\circ\!\!\rightarrow$}
\kern-23pt\raisebox{+2.pt}{$\rightarrow$}
\kern+10pt}}
 
\newcommand{\twotwo}
{\hbox{\rm $\circ\!\rightrightarrows\!
\raisebox{2.5pt}{\hbox{\small $\circ$}}
\kern-5.5pt\raisebox{-2.5pt}{\hbox{\small $\circ$}}
\!\rightrightarrows\,$}}
\newcommand{\twotwotiny}
{\hbox{\tiny $\circ\!\rightrightarrows\!
\raisebox{2.pt}{\hbox{\tiny $\circ$}}
\kern-4.2pt\raisebox{-2.pt}{\hbox{\tiny $\circ$}}
\!\rightrightarrows$}}
\newcommand{\tax}{\hbox{\sf[r,s]}}
\newcommand{\lxi}{\raisebox{0.5pt}{${}^\xi$}\!}
\vskip -0.0cm

\maketitle
\smallskip

{\tiny
\centerline{{\bf MSC} (2010): 14H50, 17B22, 17B45, 20C08,
20F36, 33D52, 30F10, 55N10, 57M25}
}
\smallskip

\vskip -0.5cm
\renewcommand{\baselinestretch}{1.2}
{\vbadness=10000\textmd
{\small \tableofcontents}
}
\renewcommand{\baselinestretch}{1.0}
\vfill\eject

\renewcommand{\natural}{\wr}

\setcounter{section}{-1}
\setcounter{equation}{0}
\section{\sc Introduction}
\vbadness=3000
\hbadness=3000
We extend the construction from \cite{ChD} of the 
{\em DAHA-Jones polynomials\,} (any reduced root systems)
and the {\em DAHA 
superpolynomials\,} in type $A$ from knots
to {\em iterated torus links\,}, including all algebraic links.
There is solid evidence that
under $\,t\!=\!q\,$ our polynomials become correspondingly
the Quantum Group (WRT) invariants of such links. 
This was checked in particular cases in \cite{CJ,ChD} 
using  \cite{Ste} and proper variants of the 
{\em Rosso-Jones formula\,}; see \cite{RJ, Mo, ChE}.
For the {\em HOMFLY-PT  polynomials\,}, this coincidence
was recently verified via the approach from \cite{MoS}
(they proved it for iterated torus {\em knots}).
This paper is devoted to the DAHA-Jones theory (for generic 
$q,t$).

Our main technical tool is the switch from  
the Macdonald polynomials in the DAHA-Jones theory from \cite{CJ}
and further works to their products, 
considered in terms of $X$ and in terms of $Y$.
The projective action of $GL_2(\Z)$ in DAHA is the key here.
The {\em Hopf\, links\,} play
an important role in our work; they directly lead to
the {\em DAHA-vertex\,}, similar to the refined
{\em topological vertex\,}. See  \cite{GIKV,AKMV,AFS}, but 
no exact relation to these papers is expected 
because we do a different version of the Macdonald theory
(fixed root systems).
Though the topological {\em $2$\~vertex} from \cite{AKMV}
is connected to our one. 

\smallskip
Due to the novelty of the direction we present in this work,
explicit examples are its very important part.
We note that the products
of Macdonald polynomials in the context of superpolynomials
of links appeared in physics works. See for instance \cite{GIKV}
(the uncolored Hopf $2$\~link) and especially
Section 4.3 from \cite{DMMSS}; our superpolynomials 
match those suggested there (for simple multiple torus knots). 

\smallskip

\subsection{\bf Overview}
Let try to put our work into perspective.
 
\subsubsection{\sf Topological connections}
The relation of our construction to the 
{\em HOMFLY-PT homology\,} \cite{KhR1,KhR2,Kh,Ras}
is not too certain, since the theory of the latter 
for links and with colors is not sufficiently developed. 
Nevertheless, we conjecture the coincidence of our superpolynomials
(type $A$) with the {\em stable 
Khovanov-Rozansky polynomials\,} for algebraic uncolored
links in the unreduced setting. 

We also establish some relation to the {\em reduced\,}
Khovanov polynomials (for $sl_2$) for small iterated links,
not only algebraic. The reduction to $sl_2$ is actually
a special case of  the approach from \cite{DGR} and 
\cite{CJ,ChD}, which is based on the theory of
differentials from \cite{Kh,Ras}. We can
recover certain (reduced) Khovanov-Rozansky
polynomials for small {\em links\,} directly  
from the corresponding DAHA-Jones polynomials.

There is a fundamental connection of our construction 
with the {\em splice diagrams\,} from \cite{EN};
we mainly use them in the form of pairs
of {\em incidence graphs\,}, which are collections of 
labeled trees. The trees with
strict positivity of $\rr,\ss$ (and under certain further
conditions in the case of {\em twisted unions\,})
describe the links of germs of {\em any\,} plane curve 
singularities. As for non-algebraic links, the simplest 
(and important) example is the so-called Hopf $\kappa$\~link, 
which corresponds to the tree with one vertex labeled 
by $[1,-1]$ and $\kappa$ arrowheads (with arbitrary colors).
The topological symmetries
of links and the corresponding splice diagrams appeared equivalent 
to the invariance of the polynomial representation of DAHA under 
the automorphisms $\tau_-, \eta, \si^2, \iota\,$
and the $\vph$\~invariance of the {\em DAHA-coinvariant\,}.
See Section \ref{sect:Aut} for the definitions.
\smallskip

For the Hopf $3$\~links ($3$\~fibers of the Hopf fibration),
the corresponding DAHA-Jones polynomials and superpolynomials
constitute the {\em topological DAHA-vertex\,}. They 
satisfy important associativity identities, 
generalizing those from paper
\cite{ChF} on the nil-DAHA approach to the Rogers-Ramanujan sums
(upon the limit $t\to 0$).  

This associativity we establish is an interesting
$q,t$\~generalization (with an additional parameter
$a$ for superpolynomials)
of the relations for the so-called $3j$\~symbols.
The Macdonald polynomials times the powers of {\em Gaussians\,}
are the key here. This power is the {\em level\,} of 
the theory; it is $0$ for the $q,t$\~deformation
of the $3j$\~symbols.
The construction possesses the permutation 
invariance, which is important in {\em TQFT\,}, and
the associativity (\ref{ass-s-three}). To be exact, 
the latter holds for the Hopf 
$3$\~links with the linking numbers $\,1,1,-1\,$.
The standard Hopf $3$\~links with the linking
numbers $\,-1,-1,-1\,$ result in the  
$\mathbf{S}_3$\~invariance; 
and then the associativity must be properly adjusted.  
\smallskip

Not all {\em solvable\,} links in $\S^3$
can be obtained by our construction, but  
all algebraic links are reached \cite{EN}. 
The basic algebraic example 
is $T(\kappa\,\rr,\kappa\,\ss)$, described by the tree 
with one vertex labeled by 
$[\rr,\ss]$ and $\kappa$ arrowheads colored by dominant 
weights. Our labels are  {\em Newton pairs\,} 
$[\rr,\ss]$, which however allowed to be negative. 
They naturally emerge in our construction vs. the ``topological" 
$\aa,\rr$\~pairs.
\smallskip

\subsubsection{\sf Algebraic links}
For algebraic links, we conjecture in the {\em uncolored case\,}
(i.e. for the fundamental 
representation) that {\em unreduced\,} DAHA-superpolynomials,
defined as the series
$\hat{\h}^{min}/(1-t)^{\kappa-1}$ for reduced ones 
$\hat{\h}^{min}$
are proportional to $\overline{\mathscr{P}}_{\hbox{\tiny alg}}$ 
from \cite{ORS} defined in
terms of the nested {\em Hilbert schemes\,} of 
(the germs of) plane curve
singularities $\c$. This is not simple to verify;
the {\em weight filtration\,,} which is the key ingredient of
the {\em ORS polynomials}, is generally
quite involved.
\smallskip

Generalizing Conjecture 2 in \cite{ObS} 
(extended to arbitrary colors and proved in \cite{Ma}), 
the {\em ORS conjecture\,} from \cite{ORS}
claims the connection of their polynomial to the Poincar\`e
polynomial 
$\overline{\mathscr{P}}_{\hbox{\tiny KhR}}
(\hbox{Link of\, } \c)$
of the {\em unreduced\,} triply-graded HOMFLY-PT homology of
Khovanov and Rozansky. This was stated there for any algebraic
uncolored links. 
\smallskip

Thus we have three independent mathematical constructions
of link superpolynomials, presumably coinciding 
(up to normalization) with each other in the
case of uncolored algebraic links.
The DAHA construction is certainly the
simplest theoretically and practically and it works
for any root systems and with arbitrary colors. 
Finding $\overline{\mathscr{P}}_{\hbox{\tiny KhR}}$
is an involved task, even for simple knots (and quite
a challenge for links and with non-trivial colors).
The polynomials $\overline{\mathscr{P}}_{\hbox{\tiny alg}}$ 
are sophisticated too, though the super-duality 
and some other important symmetries were checked for them 
\cite{ORS}, which match those in the DAHA theory. There is 
also the $4${\small th} important approach based on the 
rational DAHA \cite{Gor,GORS}, but it is restricted to torus 
knots.

We mainly omit in this work the {\em Jacobian factors\,}
of plane curve singularities.
The conjectural relation of the reduced DAHA superpolynomials
under $a=0,q=1$ to the
Betti numbers of Jacobian factors of {\em unibranch\,} $\c$
played an important role in \cite{ChD}. The passage from 
algebraic knots to links does not add 
too much new here. The DAHA-superpolynomials 
at $q=1$ 
(the case of the trivial center charge) are products of those over 
the branches; the corresponding product formula for the Betti 
numbers is geometrically obvious.
\smallskip

\subsubsection{\sf Our construction}
The data are the pairs of incidence graphs $\,\{\l,\,'\!\l\}\, $,
a unions of subtrees, which are labeled, colored and with 
arrowheads; see Section \ref{sec:Splice} for the splice
interpretation.
The labels are Newton pairs, lifted to $\ga\in GL_2(\Z)$ 
and interpreted as 
DAHA-automorphisms $\hat{\ga}$.
We check that the final formulas do not depend on
the flexibility here and, moreover, depend
only on the corresponding link, which is quite a confirmation
of the connections with topology. 

The tree structure and labels determine which 
products of the corresponding 
polynomials must be considered and which $\ga$ must be
applied; the resulting DAHA operators ($\ga$ act only in
$\HH$)
are always projected onto its polynomial representation
before the next step
\cite{ChD}. In the presence of  $\,'\!\l\,$,
we replace $X$ by $Y^{\mp 1}$ in the 
{\em pre-polynomial\,} for $\,'\!\l, \,'\!\l^\vee$ and apply 
this operator to the pre-polynomial 
for $\l$. The {\em DAHA-coinvariant\,} is applied to the
last pre-polynomial.
 
In the simplest case of the $\kappa$\~fold
torus knot $T(\kappa\rr,\kappa\ss)$, there is no $\,'\!\l$ and
the matrix $\ga$ is an
arbitrary lift of $(\rr,\ss)^{tr}$ considered as its first
column.
It is then lifted to the corresponding DAHA automorphism
$\hat{\ga}$ and is applied to the product of
Macdonald polynomials for the dominant weights assigned to
the $\kappa$ arrowheads. 
\smallskip

The usage of $Y^{\mp 1}$ for $\,'\!\l,\,'\!\l^\vee$ is connected 
with adding the {\em meridian\,} in topology. The sign corresponds
to the orientation.
Our positivity claim from Part $(ii)$ of the
Connection Conjecture \ref{CONCONJ}
and formula (\ref{deg-a-j})
for deg$_a$ from Theorem \ref{STABILIZ} are expected to hold
only for $Y^{+1}$ there (i.e. with $\vee$). The 
classical Hopf links are with 
the pairwise linking numbers $-1$, i.e. for $Y^{-1}$; the 
positivity of the corresponding
series $\hat{\h}^{min}/((1\!-\!t)(1\!-\!q))^M$ does not hold 
for them.
\smallskip 

\subsubsection{\sf DAHA theory}
The algebraic properties of the DAHA-Jones polynomials
and DAHA-superpolynomials from
\cite{CJ,GoN,CJJ,ChD} can be fully extended to 
arbitrary iterated torus {\em links\,}, 
which includes the polynomiality, duality, special values
at $q=1$, the color exchange and further symmetries.
The main new feature (vs. knots) is that the construction now 
depends on the choice of the $q,t$\~integral form of the Macdonald
polynomials.
The normalization is the division by the $LC\!M$ of all
evaluations of the Macdonald polynomials involved.
Upon this division, the notation is  $J\!D^{min}$
and $\hat{\h}^{min}$ (for the superpolynomials). 
The integral form does not matter if only one color is
involved, say in the case of knots.

Algebraically, $\hat{\h}^{min}$ is very reasonable because 
all colors are on equal footing in its definition and this
is a natural setting for our polynomiality theorem.
However $\hat{\h}^{min}$  seems generally too ``small" to be
interpreted  topologically. Namely,  
the {\em reduced\,}
HOMFLY-PT polynomial ($t=q$) requires picking one of the branches
and then the division by the evaluation of the corresponding
Macdonald polynomial, not by the total $LC\!M$. 
The {\em unreduced\,} HOMFLY-PT polynomials are  
without any divisions at all. They are topological, but
have poles for links with respect to $q$; so must be 
generally considered as $q$\~series.   Let us comment a little 
on the HOMFLY-PT polynomials in the case $t=q$.


The Macdonald $P$\~polynomials becomes Schur functions;
they do not depend on $q$ as $t\!=\!q$. So no integral form
is needed in the DAHA approach to 
{\em QG/WRT}-invariants for $t\!=\!q$.  The unreduced HOMFLY-PT 
polynomials are then uniquely recovered from the corresponding
$QG$\~invariants for (all) $A_n$ upon the substitution
$a=q^{n+1}$ for any $n$ (or any infinite sequence of them).
This can be considered as a definition of
HOMFLY-PT polynomials; another approach is
based on the theory of the {\em Skein}. The corresponding 
$QG$\~invariants are $q$\~polynomials (without denominators).
However their
$a$\~stabilizations are generally {\em rational} in terms 
of $q$. Thus we need to switch to the (infinite) $q$\~expansion 
to make the coefficients of HOMFLY-PT polynomials  meaningful
geometrically and topologically simply because of the
$a$\~stabilization. Such  {\em unreduced
setting} is not necessary for knots.
\smallskip

Generally we use the $J$\~polynomials in this paper
instead of $P$\~polynomials, which
provides the absence of the denominators
in the DAHA superpolynomials.
We note that the corresponding $\hat{\h}^{min}$ are not always 
irreducible as polynomials of $a,q,t^{\pm1}$; they are 
really {\em minimal\,} possible (irreducible) 
only for sufficiently small dominant weights.
 

The $J$\~polynomials have no exact counterparts for arbitrary
root systems. However the $QG$\~invariants are 
well defined for any root systems. Therefore the refined theory
of such invariants can be expected to involve
certain theory of $J$\~polynomials. Unfortunately 
the $q,t$\~integral forms of the Macdonald $P$\~polynomials are 
not generally settled for arbitrary root systems. Following 
\cite{C103}, we prove that the {\em spherical\,} $P$\~polynomials
become $q,t$\~integral upon multiplication by certain products 
of $q,t$\~binomials {\em of multiplicity one\,}
(with a minor reservation for the root systems $D_{2m}$). 
This is for nonsymmetric and symmetric Macdonald polynomials,
which  provides an important technical tool in the 
theory of the decomposition of the polynomial DAHA
representation and is of independent algebraic
interest. For instance, this can be used
to simplify and generalize the theory from 
\cite{C103} among other applications.
Some relations to \cite{HHL},\cite{RY},\cite{OS} can be
expected.
\smallskip

The corresponding $q,t$\~integral polynomials seem as close as 
possible to the Macdonald $J$\~polynomials. The
latter are actually not minimal as $q,t$\~integral ones; the 
focus of their theory is the stabilization and positivity 
(not only the $q,t$\~integrality). Say for $t=q$, they are 
greater than the Macdonald $P$\~polynomials, which are Schur 
functions.
The main deviation of our modified $J$\~polynomials, denoted 
by  $\tilde{P}$ in this work, from the classical 
$J$\~theory in type $A$ is that our construction is naturally 
invariant with respect to the automorphism $\iota=-w_{0}$ for the 
element of maximal length
$w_0\in W$. This is not the case for the $J$\~polynomials;
their stabilization property is incompatible with $\iota$.
\medskip

\subsection{\bf Brief history}
Due to the novelty of the DAHA-Jones theory, a brief account
of the history of this (very recent) direction can be useful 
to the readers. We will not discuss other instances of using 
Macdonald polynomials and DAHA in algebra, geometry, topology
and physics; there are quite a few now (including String Theory, 
rational DAHA, Hilbert schemes and so on).
The focus of this section will be on using $q,t$\~DAHA in 
Knot Theory and related theory of plain curve singularities, 
especially on the applications to the {\em refined\,} WRT 
invariants, HOMFLY-PT polynomials and HOMFLY-PT homology 
(stable Khovanov-Rozansky polynomials) and on using the
{\em splice diagrams}. 
\smallskip

\subsubsection{\sf Torus knots via DAHA} 
Using the Macdonald polynomials instead of Schur functions in 
the so-called {\em knot opertators\,} was 
suggested in \cite{AS}, which directly influenced \cite{CJ}. 
These operators naturally appear in the approach to the 
Jones and WRT invariants of torus knots via Conformal Field
Theory. See \cite{Ste, CJ, GoN, CJJ}.
The main difficulty of the construction of \cite{AS}
was the usage of roots of unity, similar to that 
in Verlinde algebra (when $t=q$). Even the simplest ingredient 
needed there, the {\em refined
Verlinde $S$ operator}, requires
knowing {\em all\,} Macdonald polynomials at roots of unity
that occur in the corresponding perfect modules. They
are quite non-trivial beyond $A_1$; one generally 
needs {\em any} $A_n$
to find the superpolynomial. A related problem
is that $t$ must be generally an integral power
of $q$ to ensure the existence of  the $PSL(2,\Z)$\~
action in the corresponding DAHA-module.  

There are 
some formulas for the (coefficients of) 
Macdonald polynomials for $A_n (n>1)$, but they are
very involved theoretically and practically; see e.g. \cite{HHL}. 
In spite of such obvious difficulties, 
 the authors demonstrated that their calculations for
the simplest torus knots match known or conjectured 
formulas (upon their restriction to roots of unity)
for the superpolynomials from 
topology or physics papers. The superpolynomials are
relatively simple  for the simplest torus knots, say 
$1+qt+aq$ for the uncolored trefoil; they are
generally very non-trivial apart from
the family $T(2\mm+1,2)$ even  in the 
absence of colors. 

\smallskip

Let us mention that the {\em refined\,} 
$S$ operator taken alone, already quite involved, results only 
in the unknot (with the trivial superpolynomial).  The whole 
projective unitary action 
of $PSL_2(\Z)$ in {\em perfect} modules, also
called refined (generalized) Verlinde algebras, is 
needed in this approach, 
which is due to Cherednik and A. Kirillov Jr.
This makes using perfect DAHA modules here
difficult,
even if the latest software
(like SAGE) is employed for the Macdonald polynomials. 

The lift of the formulas at roots of unity $q,t$ 
to generic $q,t$ is always quite a challenge (including
the Verlinde algebras, i.e. in the unrefined theory), 
unless the existence of the superpolynomials and the bounds 
for the degrees of $q,t,a$
are  {\em a priori\,} known. The rank stabilization (associated
with $a$) has its own challenges too, even for generic $q,t$. 
Also, we do any colors and non-torus iterated knots/links.
\smallskip

These problems were resolved (actually bypassed) 
in \cite{CJ}, which eliminated the usage of roots
of unity and was written for any torus knots, arbitrary 
reduced root systems and dominant weights (colors).
It was a significant development and the beginning of the 
DAHA-Jones theory. The construction of \cite{AS} results from
that from \cite{CJ} simply because the projective
action of $PSL_2(\Z)$ in DAHA reduces to that in its 
perfect modules (a standard theorem from \cite{C101});
no other justification of this connection with \cite{AS}
is necessary.
\smallskip

The coincidence with the HOMFLY-PT
polynomials (as $t\!=\!q$) for torus knots was proven in 
\cite{CJ} and the connections were conjectured  
with the (colored) WRT-invariants for any root systems  
and reduced uncolored stable Khovanov-Rozansky polynomials
(the Poincar\`e polynomials of the HOMFLY-PT homology). 
The latter are quite involved; not many explicit formulas are
known for them. Quite a few uncolored and colored DAHA 
superpolynomials were calculated in this paper, as well as the 
refined Quantum Group (WRT) invariants 
for the classical and exceptional root systems. 
 
As one of the applications, the refined $QG$ invariants of types 
$E_{6,7,8}$ for minuscule and quasi-minuscule weights and
simplest torus knots were 
calculated in \cite{CJ}, conjecturally coinciding with
the corresponding {\em QG} invariants at $t=q$
(confirmed for $E_6$ by Ross Elliot
via the Rosso-Jones formula). 
\smallskip

Let us mention that the (projective) action of $PSL_2(\Z)$ in DAHA 
is a generalization of the action of $PSL_2(\C)$ in the 
Heisenberg algebra; DAHA is its certain deformation. 
The key obstacle in classical theory of Fourier-Hankel 
transform is that there is no action of $PSL_2(\C)$ in the Fock 
representation. This remains the same in DAHA theory,
unless upon the reduction to perfect modules (generalized
Verlinde algebras). Fortunately, such a reduction appeared
unnecessary for DAHA-Jones theory and DAHA superpolynomials.
\smallskip

Concerning the current status of conjectures from \cite{CJ},
practically all ``intrinsic" ones about the existence
and the structure of the DAHA-Jones polynomials 
and DAHA-superpolynomials (the $A$-case) were
proved in \cite{CJJ,GoN}. The only
conjecture in type $A$ from this paper
that remains open by now is the positivity 
of the DAHA superpolynomials for rectangle Young diagrams.
There are recent geometric developments here (but no
proof so far). 

The existence of
the {\em DAHA-hyperpolynomials} in types $B,C,D$ and their
symmetries remain open,
though the approach from \cite{CJ} (based on \cite{SV})
can be extended to this case and there is a sketch of the 
proof of {\em hyper-duality\,} in \cite{CJJ}. In type $D$, 
the {\em Kauffman polynomials\,} (instead of the HOMFLY-PT ones)
were conjectured in \cite{CJ} to coincide with the DAHA-Jones
polynomials  as $t\!=\!q$; this is justified by now.
\smallskip
 
We note that some of the symmetries of superpolynomials were 
suggested by physicists. The most interesting one is the 
{\em super-duality} conjectured in \cite{GS}.
This conjecture became rigorous and 
for any Young diagram in the
DAHA setting \cite{CJ}. 
In String Theory, the super-duality is related
to the approach to the refined Chern-Simons-Witten theory via the
$M_5$\~ theory; the action of $\C^*\times \C^*$ there
is naturally associated with parameters $q,t$.
See also e.g. \cite{DGR,AS,DMMSS}
concerning various physics aspects and formulas.

The DAHA super-duality for torus knots was proven in 
\cite{GoN} (see \cite{CJJ} for an approach
via roots of unity) together with a justification of the
stabilization of the DAHA-Jones polynomials in type $A$,
which was announced in \cite{CJ} as a theorem. The switch 
from the Macdonald polynomials of type $A_n$ to the so-called
$J$\~polynomials (the stable ones) in \cite{GoN} is 
important. The $J$\~polynomials can be avoided for 
knots, but are absolutely necessary 
for {\em links} (our present work). 

\subsubsection{\sf From torus knots to iterated links}
The main demerit of the DAHA-Jones
theory after \cite{AS,CJ,GoN,CJJ} was its restrictions 
to torus knots.
Arbitrary algebraic knots  and links (not only
torus ones) are very much needed because of many reasons. 
\smallskip

First of all, this generalization is necessary to employ the 
technique 
of physically inspired theory of the {\em resolved conifold\,}, 
and its monoidal-type transformations, used in \cite{Ma} to prove 
the {\em OS conjecture\,} from \cite{ObS}. Second, there are 
significant links to the Fundamental Lemma and related 
algebraic geometry of rational planar curves, 
with various implications within and beyond the Langlands Program. 

Third, the topological reasons for the switch from
torus knots to any {\em iterated torus links} are quite obvious.
This class of links is closed with respect to the {\em cabling\,},
one of the major operations in Knot Theory. We recall that
all algebraic links are cables of unknot. However 
iterated torus knots and links are generally non-algebraic. It is 
known which ones are algebraic (see \cite{EN}), but this is far 
from trivial, especially for links. The present paper is mostly 
written for any iterated torus links (not only algebraic),
though the positivity conjectures for links generally require the
algebraic ones (and the unreduced setting).
\smallskip

The topological perspective is very important here since the 
superpolynomials  (of any origin, geometric, algebraic or
physical) presumably coincide with 
the stable  Khovanov-Rozansky polynomials. The theory of such 
polynomials and HOMFLY-PT homology is mainly restricted 
now to {\em uncolored knots\,}
(unless for the Khovanov $sl_2$ homology). There are 
recent developments, including \cite{WW} and the approach
based on the Howe duality,
that allow in principle to incorporate colors, but technical
difficulties remain very serious. The Khovanov-Rozansky
theory and {\em categorification} are of course for {\em all}
knots, not only for {\em cables}, the level DAHA have
reached so far. 
\smallskip

The passage to arbitrary torus iterated knots from 
 torus knots was an important step of the DAHA approach. Let us 
mention that we expected our theory in \cite{ChD} to be connected 
with paper \cite{Sam}, but this did not materialize. We failed 
to understand the approach used in \cite{Sam} for $q\neq t$. 
In contrast to \cite{ChD}, this paper is for $A_1$ only, but
nevertheless his polynomials $J_n$ have significant 
$q,t$\~denominators, which is not what can be expected. 
The polynomiality of DAHA-Jones polynomials is the key in 
\cite{CJ} and further works, including this one.
The examples of his $J_2$\~polynomials 
for $C\!ab(\pm 5,2)T(3,2)$ from Section 5.2 (ibid.) are very
different from our ones for the same cables as $q\neq t$ 
(and we do not understand how they were obtained and cannot 
reproduce them ourselves). See Section 4 of \cite{ChD}.
\smallskip

\subsubsection{\sf Using the Skein}
Let us discuss paper \cite{MoS}. As $t\!=\!q$, the authors
establish the connection of the 
{\em Skein} of the torus with the Elliptic Hall algebra 
\cite{SV}.  This implies (through quite a few
technical steps) the coincidence of the DAHA superpolynomials 
of arbitrary torus iterated {\em knots} at $t\!=\!q$  with the 
corresponding HOMFLY-PT  polynomials. This was established for 
any iterated torus knots for $A_1$ in \cite{ChD} and for any 
torus knots in \cite{CJ} in the case of arbitrary $A_n$ (and 
therefore for the HOMFLY-PT polynomials). It is 
now checked for the Kauffman polynomials too (this is
the DAHA $D$\~case;
see \cite{CJ}).
\smallskip 

We note that the exact framing 
factor (a power of $q,a$) is missing in \cite{MoS} in contrast
to \cite{ChD} (for $A_1$). Also, a nonstandard framing in 
\cite{MoS}, which makes their construction
directly in terms of Newton's pairs (vs. the ``topological" 
$\aa,\rr$-parameters), potentially creates
problems with the isotopy equivalence (they refer to our 
algebraic proof of the topological invariance). 
\smallskip

Using the {\em Seifert framing\,}, which extends the framing used 
in \cite{MoS} (Definition 7.1) from knots to links, the 
splice diagrams \cite{EN}, and the interpretation of the 
DAHA-multiplication via link operations, 
we can extend Section 7.1 of \cite{MoS} from iterated
torus {\em knots} considered there to arbitrary (colored) torus 
iterated {\em links\,}. For instance, the multiplication by 
$J_\la(X)$ for a Macdonald polynomial is essentially adding an 
unknot colored by $\la$; also, applying $J_\la(Y)$ to a
{\em pre-polynomial\,} is interpreted as adding a {\em meridian}. 
These are standard facts in the theory of skein (Morton and others) 
and they are actually the key for the passage from knots to 
links (from \cite{ChD} to this paper). 
Thus, Conjecture \ref{CONCONJ}, $(i)$ below (for $t\!=\!q$) 
becomes a theorem if our present paper is combined with 
\cite{MoS}; we will provide the details elsewhere.
\smallskip

This coincidence is the $A$\~case of the conjectured 
connections of the DAHA-Jones polynomials at $t=q$ with the 
general {\em WRT-invariants\,} for 
arbitrary root systems and for any dominant weights. 
They do not seem very difficult to justify (including torus
iterated {\em links\,}) using the Rosso-Jones cabling
formula \cite{RJ} and the
DAHA shift operator, as it was demonstrated in \cite{CJ,ChD,ChE}.
We have a sketch of such a justification for $A_n$ and Ross Elliot
can essentially do this for general root systems (unpublished).
Another approach here is via CFT and Verlinde algebras, but
it seems more involved. Using the skein provides a 
justification ``without calculations" (thought there were
quite a few steps in \cite{MoS} to verify),  but this is 
restricted only to the case of $A_n$.
\smallskip

{\sf Summary.} Let us briefly summarize the main points discussed 
above. The passage from knots to arbitrary torus 
iterated links in our present paper is an important step 
in the DAHA-Jones theory. This is absolutely necessary from
the topological-geometric perspective, for instance for reaching
the invariants of $3$-folds (via framed links). Also,
adding colors in Knot Theory is closely related to the passage 
from knots to links. The multi-brunch plane curve singularities 
are very important in its own right (including the Fundamental
Lemma). One more direction to be mentioned,
is the physical and mathematical theory of topological vertex 
based on Hopf links. Algebraically, the DAHA theory of
torus iterated links we present 
here seems reasonably complete; we use it to approach the 
$q,t$\~Skein at the end of the paper.

\subsubsection{\sf Splice diagrams}
There is a long history of using combinatorial/graphic
presentations of knots, links and related $3$-folds.
Algebraic links are of great interest here, since
they are in one-to-one correspondence  with the plane 
curve singularities. {\em Splice diagrams\,} suggested by 
Neumann and systematically developed in \cite{EN} appeared 
a great tool for us. They were aimed at Seifert manifolds 
(and their plumbing). This generality is beyond our paper;  
we assume that the links are in $\mathbb{S}^3$, not in an 
arbitrary homology 3-sphere. 

The main result of
\cite{EN} is that the splice diagrams are  isotopy
invariants of the corresponding links and 
describe all of them. See Theorems 9.2, 9.4 in \cite{EN},
Section \ref{sec:Splice} below and 
the Appendix \ref{TopologySection}
on splice diagrams.
We note that our 
method does not produce all splice diagrams of solvable type;
but {\em all\,} algebraic links can be reached. 

The key operation from \cite{EN} is
{\em splicing\,}; it provides large 
families of links. Other operations can be mainly considered 
as its special cases. The {\em cabling} and
{\em twisted unions} play the major role in the DAHA
approach to torus iterated links. The other operations
(in our context) are erasing components, orientation reversion,
and disjoint sums. The isotopy invariance of our construction 
is stated in Theorem \ref{MAINTHM},$(ii)$ for 
DAHA-Jones polynomials (any root systems) and 
in  Theorem \ref{STABILIZ}, $(ii)$ for the superpolynomials.

The justification of the topological symmetries from 
Theorem \ref{MAINTHM},$(ii)$ is
essentially parallel to Theorem 1.2 from
\cite{CJJ}. Namely, it is proven there that the DAHA-Jones
polynomials $J\!D$ (and superpolynomials) are the same 
for the torus knots $T(\rr,\ss)$ and $T(\ss,\rr)$.
This coincidence, the triviality of these polynomials
for $T(\rr,1)$ and the effect of {\em mirroring\,} 
(changing $q,t$ to $q^{-1},t^{-1}$ when
$\rr\mapsto -\rr$) are the only topological
symmetries in the small universe of torus knots.
These facts were deduced in \cite{CJJ}
from the properties of the DAHA involutions.
We note that the symmetry $T(\rr,\ss)\leftrightarrow
T(\ss,\rr)$ can be far from obvious in other
approaches to superpolynomials (say, 
when rational DAHA is used for obtaining such polynomials).

When switching to torus iterated knots, one also needs to check
that the torus knot  $T(\rr,m\rr+\ss)$ results 
in the same DAHA-Jones and superpolynomials
as the iterated knots corresponding to the Newton pairs 
$[1,m],[\rr,\ss]$ (the latter is employed the first). 
This is the key new symmetry 
here. Topologically, $T(\rr,m\rr+\ss)$ is obviously isotopic 
to  $C\!ab(m\rr+\ss,\rr)T(1,m)$ since
$T(1,m)$ is unknot. The corresponding
relation for the $J\!D$\~polynomials readily 
follows from the commutativity $\tau_-^m$ (associated with
$[1,m]$)
with the projection $\Downarrow$ onto the polynomial
representations, which is equivalent
to the fact that $\tau_-$  acts in this
representation \cite{C101}.  Other symmetries are due to applying 
$\eta, \iota, \si^2$ inside the DAHA coinvariant;
they are compatible with the projection $\Downarrow$ as well. 
We also constantly use the $\vph$\~invariance of the
coinvariant.

Generalizing the theory from \cite{CJJ,ChD} 
to arbitrary torus iterated {\em links}, the major
(additional) fact needed for the isotopy invariance of
our construction is
Theorem \ref{EVALTAUM}, which states that
\begin{align*}
\bigl\{\tau_-^{-1}(f),\tau_-^{-1}(g)\bigr\}_{ev}=
\bigl\{\tau_-^{-1}(fg)\bigr\}_{ev}=
\bigl\{\tau_-^{-1}\bigl(f(X)g(X)\bigr)\bigr\}_{ev}
\end{align*}
for arbitrary $f=f(X),g=g(X)$ and the evaluation pairing
$\{f,g\}_{ev}=\{f(Y^{-1})\bigl(g(X)\bigr)\}_{ev}$,
where $\{\cdot\}_{ev}$ is the DAHA coinvariant.

All other symmetries are either straightforward or follow
from this theorem and the theorems mentioned above
for torus iterated knots (based on the action of 
$\vph,\eta,\tau_-$ and $\si^2$). For instance, applying 
$\vph$ provides switching the components in the pairs of trees 
$\{\l,'\!\l\}$; see  (\ref{lsymlprime}) for the exact
relation.

The pairs of trees give a very natural way to encode the
DAHA invariants, which corresponds to the consideration of
splice diagrams with {\em marked edges\,} in the terminology 
of \cite{EN}. The marked edge shows the place where the DAHA 
coinvariant will be applied to the last 
{\em pre-polynomial\,}.
I.e. this gives the last step of our calculation.
In the case of a single tree, the last edge is 
marked (connecting the
last vertex with one of the leaves from it in the language of
splice diagrams); the tree for an iterated {\em knot} is
simply a path with an arrow at its end and {\em leaves\,} at 
the vertices \cite{ChD}. 

The isotopy invariance of DAHA-Jones
polynomials (and superpolynomials) includes
the proof of their independence of the 
choice of the marked edge; this is the key new
feature when doing links. Assuming that
the graph is connected, it suffices to
check that adjacent marked edges give the same. And this 
follows directly from the symmetries discussed above, 
including Theorem \ref{EVALTAUM} and
the equivalence of $T(\rr,\ss)$ and  $T(\ss,\rr)$
(the transposition of a leaf and an edge from the same
vertex in the terminology of \cite{EN}).
\smallskip

A typical and instructional example of using 
Theorem \ref{EVALTAUM} is the coincidence of
the superpolynomials for the link from (\ref{T6-4-y}) and
that for the cable 
$\bigl(C\!ab(2,3)C\!ab(-1,0)\bigr)T(1,-1)$ from (\ref{T2-3--1-2}). 
We provide there a direct DAHA  deduction of the corresponding
symmetry  and its interpretation
in terms of splice diagrams. 
There are many examples of this kind in the paper. For instance,
the link $\bigl(C\!ab(8,3)C\!ab(0,1)\bigr)T(2,1)$
from (\ref{T8-3-2-1}) and (\ref{T8-3-2-1X}) 
is considered in full detail.
\smallskip

\subsubsection{\sf Conclusion}
As we tried to demonstrate, the extension of 
\cite{CJ,GoN,CJJ,ChD} from {\em knots\,} to torus iterated 
{\em links\,} in this paper (any reduced systems and weights)
is an important step with various possible
applications in Knot Theory and theory of plane curve 
singularities. Using links seems inevitable here for
any self-consistent theory. Even the simplest of them
are quite valuable; say, the Hopf links govern the 
new theory of DAHA-vertex. 

One of the applications of links is the approach to the
toric $q,t$\~Skein we suggest at the end of the paper,
based on {\em generalized 
twisted unions}; see Section \ref{sec:skein}.
The DAHA {\em knot operators\,} 
are certain (symmetric) elements in $\HH$ corresponding to 
links $\l$. Their {\em matrix elements\,} 
are essentially the 
superpolynomials of the links that are obtained from  given
$\l$ by adding
two probe links ``on the left and on the right". 
This direction will be continued in our further works.

Due to the novelty
of the presented theory, we  provide many examples, 
including detailed analysis of algebraic and
topological symmetries. 
The numerical formulas are important. We selected
the most instructional ones;
each and every serves some purpose. The theory of 
algebraic/iterated links is quite ramified and we need to 
demonstrate its range (and related DAHA features) theoretically 
and practically. The DAHA approach provides a unique way of 
obtaining the superpolynomials for arbitrary colored iterated 
links; we think that explicit formulas can help researches from
neighboring fields, including various divisions
of Macdonald theory and DAHA theory. 
\smallskip

The intrinsic algebraic theory of DAHA-Jones polynomials
and the theory of DAHA superpolynomials in type $A$ 
is sufficiently well developed by now (we think, better than 
other approaches to superpolynomials). This includes the key
$q,t^{\pm1},a$\~polynomiality of DAHA-Jones polynomials,
the super-duality and quite a few other (proven) properties. 
We expect that the passage to links in this paper is an important 
step toward the theory of refined invariants of Seifert 
$3$\~manifolds and related toric-type surfaces,  hopefully 
including significant applications in Number Theory.  
\medskip

\setcounter{section}{0}
\setcounter{equation}{0}
\section{\sc Double Hecke algebras}
\subsection{\bf Affine root systems}
Let $R=\{\al\}   \subset \R^n$ be a root system of type
$A_n,\ldots,\!G_2$
with respect to a euclidean form $(z,z')$ on $\R^n
\ni z,z'$,
$W$ the Weyl group 
generated by the reflections $s_\al$,
$R_{+}$ the set of positive  roots
corresponding to fixed simple 
roots $\al_1,...,\al_n;$ $R_-=-R_+$. 
The form is normalized
by the condition  $(\al,\al)=2$ for 
{\em short\,} roots. 
The weight lattices are
$P=\oplus^n_{i=1}\Z \om_i$, 
where $\{\om_i\}$ are fundamental weights:
$ (\om_i,\al_j^\vee)=\de_{ij}$ for the
coroots $\al^\vee=2\al/(\al,\al).$
The root lattice is $Q=\oplus_{i=1}^n \Z\al_i$.
Replacing $\Z$ by $\Z_{\pm}=\{m\in\Z, \pm m\ge 0\}$, we obtain
$P_\pm,Q_\pm$. See  e.g., \cite{Bo} or \cite{C101}. 

Setting 
$\nu_\al\equal (\al,\al)/2$,
the vectors $\ \tal=[\al,\nu_\al j] \in
\R^n\times \R \subset \R^{n+1}$
for $\al \in R, j \in \Z $ form the
{\em twisted affine root system\,}
$\tR \supset R$, where $z\in \R^n$ are identified with $ [z,0]$.
We add $\al_0 \equal [-\vth,1]$ to the simple
roots for the {\em maximal short root\,} $\vth\in R_+$.
The corresponding set
$\tR_+$ of positive roots is 
$R_+\cup \{[\al,\nu_\al j],\ \al\in R, \ j > 0\}$.

The set of the indices of the images of $\al_0$ by all
automorphisms of the affine Dynkin diagram will be denoted by 
$O$ ($O=\{0\} \for E_8,F_4,G_2$). 
Let $O'\equal\{r\in O, r\neq 0\}$.
The elements $\om_r$ for $r\in O'$ are  
{\em minuscule weights\,}, defined by the inequalities 
$(\om_r,\al^\vee)\le 1$ for all $\al \in R_+$. We set $\om_0=0$
for the sake of uniformity.
\smallskip

{\sf Affine Weyl groups}
Given $\tal=[\al,\nu_\al j]\in \tR,  \ b \in P$, let
\begin{align}
&s_{\tal}(\tz)\ =\  \tz-(z,\al^\vee)\tal,\
\ b'(\tz)\ =\ [z,\ze-(z,b)]
\label{ondon}
\end{align}
for $\tz=[z,\ze] \in \R^{n+1}$.
The
{\em affine Weyl group\,} $\tW=\lan s_{\tal}, \tal\in \tR_+\ran)$ 
is the semidirect product $W\lsmash Q$ of
its subgroups $W=$ $\lan s_\al,
\al \in R_+\ran$ and $Q$, where $\al$ is identified with
\begin{align*}
& s_{\al}s_{[\al,\,\nu_{\al}]}=\
s_{[-\al,\,\nu_\al]}s_{\al}\for
\al\in R.
\end{align*}

The {\em extended Weyl group\,} $ \hW$ is $W\lsmash P$, where
the corresponding action is 
\begin{align}
&(wb)([z,\ze])\ =\ [w(z),\ze-(z,b)] \for w\in W, b\in P.
\label{ondthr}
\end{align}
It is isomorphic to $\tW\lsmash \Pi$ for $\Pi\equal P/Q$. 
The latter group consists of $\pi_0=$id\, and the images $\pi_r$
of minuscule $\om_r$ in $P/Q$. 
Note that 
$\pi_r^{-1}$ is $\pi_{r^\iota}$
and $u_r^{-1}$ is $u_{r^\iota}$,  where $\iota$ is 
the standard involution 
(sometimes trivial) of the {\em nonaffine\,} Dynkin diagram,
induced by $\al_i\mapsto -w_0(\al_i)$.
Generally
$\iota(b)=-w_0(b)=b^\iota$, where $w_0$ is the longest element 
in $W$.

The group $\Pi$
is naturally identified with the subgroup of $\hW$ of the
elements of the length zero; the {\em length\, } is defined as 
follows:
\begin{align*}
&l(\hw)=|\La(\hw)| \for \La(\hw)\equal\tR_+\cap \hw^{-1}(-\tR_+).
\end{align*}
One has $\om_r=\pi_r u_r$ for $r\in O'$, where $u_r$ is the 
element $u\in W$ of {\em minimal\,} length such that 
$u(\om_r)\in P_-$,
equivalently, $w=w_0u$ is of {\em maximal\,} length such that  
$w(\om_r)\in P_+$. The elements $u_r$ are very explicit.
Let $w^r_0$ be the longest element
in the subgroup $W_0^{r}\subset W$ of the elements
preserving $\om_r$; this subgroup is generated by simple
reflections. One has:
\begin{align}\label{ururstar}
u_r = w_0w^r_0 \hbox{\,\, and\,\, } (u_r)^{-1}=w^r_0 w_0=
u_{r^\iota} \for r\in O.
\end{align}
\smallskip

Setting $\hw = \pi_r\tw \in \hW$ for $\pi_r\in \Pi,\, \tw\in \tW,$
\,$l(\hw)$ coincides with the length of any reduced decomposition
of $\tw$ in terms of the simple reflections
$s_i,\, 0\le i\le n.$ Thus $\Pi$ is a subgroup of $\hW$ of the
elements of length $0$. 
\smallskip

\subsection{\bf Definition of DAHA}
We follow \cite{CJJ,CJ,C101}.
Let $\mm$ be the least natural number
such that  $(P,P)=(1/\mm)\Z.$  Thus
$\mm=|\Pi|$ unless 
$\mm=2 \for D_{2k}$ and $\ \mm=1 \for B_{2k},C_{k}.$ 

The {\em double affine Hecke algebra, DAHA\,}, depends
on the parameters
$q, t_\nu\, (\nu\in \{\nu_\al\})\,$ and will be defined
over the ring
$\Z_{q,t}\equal\Z[q^{\pm 1/\mm},t_\nu^{\pm 1/2}]$
formed by
polynomials in terms of $q^{\pm 1/\mm}$ and
$\{t_\nu^{\pm1/2}\}.$ Note that the coefficients of the 
Macdonald polynomials will belong to 
$
\Q(q,t_\nu).
$

For $\tal=[\al,\nu_\al j] \in \tR,\ 0\le i\le n$, we set
\begin{align}\label{talqal}
&   t_{\tal} =t_{\al}=t_{\nu_\al}=q_\al^{k_\nu} ,\ \, 
q_{\tal}=q^{\nu_\al}, \ \, t_i = t_{\al_i},  q_i=q_{\al_i},
\end{align}

Also, using here (and below) {\em\small\, sht,\ lng\,} instead 
of $\nu$, we set
\begin{align*}
\rho_k\equal \frac{1}{2}\!\sum_{\al>0} k_\al \al=
k_{\sht}\rho_{\sht}\!+\!k_{\lng}\rho_{\lng},\ \,
\rho_\nu=\frac{1}{2}\!\sum_{\nu_\al=\nu} \al=
\!\!\sum_{\nu_i=\nu,\, i>0}  \om_i.
\end{align*}

For pairwise commutative $X_1,\ldots,X_n,$
\begin{align}
& X_{\tb}\ \equal\ \prod_{i=1}^nX_i^{l_i} q^{ j}
\iif \tb=[b,j],\ \hw(X_{\tb})\ =\ X_{\hw(\tb)},
\label{Xdex}\\
&\hbox{where\ } b=\sum_{i=1}^n l_i \om_i\in P,\ j \in
\frac{1}{ m}\Z,\ \hw\in \hW.
\notag \end{align}
For instance, $X_0\equal X_{\al_0}=qX_\vth^{-1}$.
\medskip

Recall that 
$\om_r=\pi_r u_r$ for $r\in O'$ (see above) and that
$\pi_r^{-1}=\pi_{\iota(i)}$, where
$\,\iota\,$ here is the standard involution  of the nonaffine 
Dynkin diagram.
Finally, we set $m_{ij}=2,3,4,6$
when the number of links between $\al_i$ and $\al_j$ in the affine 
Dynkin diagram is $0,1,2,3$.

\begin{definition}
The double affine Hecke algebra $\HH\ $
is generated over $\Z_{q,t}$ by
the elements $\{ T_i,\ 0\le i\le n\}$,
pairwise commutative $\{X_b, \ b\in P\}$ satisfying
(\ref{Xdex})
and the group $\Pi,$ where the following relations are imposed:

(o)\ \  $ (T_i-t_i^{1/2})(T_i+t_i^{-1/2})\ =\
0,\ 0\ \le\ i\ \le\ n$;

(i)\ \ \ $ T_iT_jT_i...\ =\ T_jT_iT_j...,\ m_{ij}$
factors on each side;

(ii)\ \   $ \pi_rT_i\pi_r^{-1}\ =\ T_j \iif
\pi_r(\al_i)=\al_j$;

(iii)\  $T_iX_b \ =\ X_b X_{\al_i}^{-1} T_i^{-1} \iif
(b,\al^\vee_i)=1,\
0 \le i\le  n$;

(iv)\ $T_iX_b\ =\ X_b T_i\ $ if $\ (b,\al^\vee_i)=0
\for 0 \le i\le  n$;

(v)\ \ $\pi_rX_b \pi_r^{-1}\ =\ X_{\pi_r(b)}\ =\
X_{ u^{-1}_r(b)}
 q^{(\om_{\iota(r)},b)},\  r\in O'$.
\label{double}
\end{definition}

Given $\tw \in \tW, r\in O,\ $ the product
\begin{align}
&T_{\pi_r\tw}\equal \pi_r T_{i_l}\cdots T_{i_1},\where
\tw=s_{i_l}\cdots s_{i_1} \for l=l(\tw),
\label{Twx}
\end{align}
does not depend on the choice of the reduced decomposition
Moreover,
\begin{align}
&T_{\hv}T_{\hw}\ =\ T_{\hv\hw}\  \hbox{ whenever\,}\
 l(\hv\hw)=l(\hv)+l(\hw) \for
\hv,\hw \in \hW. \label{TTx}
\end{align}
In particular, we arrive at the pairwise
commutative elements 
\begin{align}
& Y_{b}\equal
\prod_{i=1}^nY_i^{l_i} \iif
b=\sum_{i=1}^n l_i\om_i\in P,\ 
Y_i\equal T_{\om_i},b\in P.
\label{Ybx}
\end{align}
When acting in the polynomial representation
$\v$ (see below), they are called {\em difference
Dunkl operators}.

\subsection{\bf The automorphisms}\label{sect:Aut}
The following maps can be (uniquely) extended to
automorphisms of $\HH\,$, where 
$q^{1/(2\mm)}$ must be added to $\Z_{q,t}$
(see \cite{C101}, (3.2.10)--(3.2.15))\,:
\begin{align}\label{tauplus}
& \tau_+:\  X_b \mapsto X_b, \ T_i\mapsto T_i\, (i>0),\
\ Y_r \mapsto X_rY_r q^{-\frac{(\om_r,\om_r)}{2}}\,,
\\
& \tau_+:\ T_0\mapsto  q^{-1}\,X_\vth T_0^{-1},\
\pi_r \mapsto q^{-\frac{(\om_r,\om_r)}{2}}X_r\pi_r\
(r\in O'),\notag\\
& \label{taumin}
\tau_-:\ Y_b \mapsto \,Y_b, \ T_i\mapsto T_i\, (i\ge 0),\
\ X_r \mapsto Y_r X_r q^\frac{(\om_r,\om_r)}{ 2},\\
&\tau_-(X_{\vth})= 
q T_0 X_\vth^{-1} T_{s_{\vth}}^{-1};\ \
\si\equal \tau_+\tau_-^{-1}\tau_+\, =\,
\tau_-^{-1}\tau_+\tau_-^{-1},\notag\\
&\si(X_b)=Y_b^{-1},\   \si(Y_b)=
T_{w_0}^{-1}X_{b^{\,\iota}}^{-1}T_{w_0},\ \si(T_i)=T_i (i>0).
\label{taux}
\end{align}
These automorphisms fix $\ t_\nu,\ q$
and their fractional powers, as well as the
following {\em anti-involution\,}:
\begin{align}
&\vph:\ 
X_b\mapsto Y_b^{-1},\, Y_b\mapsto X_b^{-1},\
T_i\mapsto T_i\ (1\le i\le n).\label{starphi}
\end{align} 

The following anti-involution results directly from
the group nature of the DAHA relations:
\begin{align}\label{star-conj}
H^\star= H^{-1} \for H\in \{T_{\hw},X_b, Y_b, q, t_\nu\}.
\end{align}
To be exact, it is naturally extended to the fractional
powers of $q,t$:
$$
\star:\ t^{\frac{1}{2\mm}}_\nu \mapsto t_\nu^{-\frac{1}{2\mm}},\
q^{\frac{1}{2\mm}}\mapsto  q^{-\frac{1}{2\mm}}.
$$
This anti-involution serves the inner product in the theory 
of the DAHA polynomial representation.

We will also need the involution:
\begin{align}
\eta:\ &T_i\!\mapsto\! T_i^{-1},\, X_b\!\mapsto\! X_b^{-1},\,
\pi_r\!\mapsto\! \pi_r ,\ 
\,t^{\frac{1}{2\mm}}_\nu\!\mapsto\! t_\nu^{-\frac{1}{2\mm}},\, 
q^{\frac{1}{2\mm}}\!\mapsto\!  q^{-\frac{1}{2\mm}}.
\label{etatxpi}
\end{align}
where $\,0\le i\le n,\, r\in O',\, b\in P$. 
Its actions on $Y_b$ is not that
uniform:
\begin{align}
\eta:\ Y_r\mapsto \pi_{\iota(r)}^{-1} T_{u_{\iota(r)}}^{-1},\ \,
Y_\vth\mapsto  T_0^{-1} T_{s_{\vth}}^{-1},
\label{etatxpiy}
\end{align}
where $\iota$ is the involution of the {\em nonaffine\,}
Dynkin diagram; see (\ref{ururstar}).
\smallskip

The involution $\eta\,$ extends the
{\em Kazhdan\~Lusztig involution\,} in the affine Hecke theory;
see \cite{C101}, (3.2.19--22). Note that 
\begin{align}\label{eta-tau}
\vph\tau_\pm\vph=\tau_\mp=\si\tau_\pm^{-1}\si^{-1},\ 
\eta\tau_{\pm}\eta=\tau_{\pm}^{-1},\ 
\vph\si\vph=\si^{-1}=\eta\si\eta.
\end{align}

Let us list the matrices corresponding to the automorphisms and
anti-automorphisms above upon the natural projection 
onto $GL_2(\Z)$, corresponding to  
$\,t^{\frac{1}{2\mm}}_\nu=1=q^{\frac{1}{2\mm}}$. 
The matrix {\tiny 
$\begin{pmatrix} \al & \be \\ \ga & \de\\ \end{pmatrix}$}
will then represent the map $X_b\mapsto X_b^\al Y_b^\ga,
Y_b\mapsto X_b^\be Y_b^\de$ for $b\in P$. One has:
\smallskip

$\!\!\!\tau_+\rightsquigarrow$ 
{\tiny 
$\begin{pmatrix}1 & 1 \\0 & 1 \\ \end{pmatrix}$},\ 
$\tau_-\rightsquigarrow$ 
{\tiny 
$\begin{pmatrix}1 & 0 \\1 & 1 \\ \end{pmatrix}$},\
$\si\rightsquigarrow$ 
{\tiny 
$\begin{pmatrix}0 & 1 \\-1 & 0 \\ \end{pmatrix}$},\
$\vph\rightsquigarrow$ 
{\tiny 
$\begin{pmatrix}0 & -1 \\-1 & 0 \\ \end{pmatrix}$},\
$\eta\rightsquigarrow$ 
{\tiny 
$\begin{pmatrix}-1 & 0 \\0 & 1 \\ \end{pmatrix}$}.\
\smallskip

{\sf Enhanced projective $GL_2(\Z)$.}
The {\em projective\,} $GL_2(\Z)$ is 
the group generated by $\tau_{\pm}, \eta$ subject to the
relations $\tau_+\tau_-^{-1}\tau_+=
\tau_-^{-1}\tau_+\tau_-^{-1},\ \eta^2=1$ and 
$\eta\tau_{\pm}\eta=\tau_{\pm}^{-1}.$ The notation will
be $GL_{\,2}^{\wedge}(\Z)$.
The span 
of $\tau_{\pm}$ is the projective $PSL_{2}(\Z)\,$
(due to Steinberg), which is isomorphic to the braid 
group $B_3$.

Let us enrich these groups by the
following automorphisms of $\HH$.
For the pair of arbitrary characters $u,v$ of $\Pi=P/Q,$  
\begin{align}\label{zetauv|}
\ze_{u,v}(X_a T_w Y_b)=u(a)v(b)X_a T_w Y_b \for
a,b\in P,\, w\in W,
\end{align} 
where the order of $X,T,Y$ does not matter here and $T_w$
can be replaced by $T_{\tw}$ for any $\tw\in \tW$. 
The map $X_a T_w Y_b\mapsto u(a)v(b)$ can be 
readily extended to a character of $\HH$ (its one-dimensional
representation).

These automorphisms satisfy the following relations:
\begin{align}\label{zetauvrel}
&\vph\, \ze_{u,v}\,\vph\,=\,\ze_{v^{-1},u^{-1}},\ 
\eta\, \ze_{u,v}\,\eta\,=\,\ze_{u^{-1},v},\\
&\tau_+^{-1} \ze_{u,v}\tau_+=\ze_{u,uv},\, 
\tau_-^{-1} \ze_{u,v}\tau_-=\ze_{\,uv,v},\,
\si^{-1}\ze_{u,v}\si=\ze_{v^{-1},u},\notag\\
\label{zetauvrel1}
&\tau_+^u (\tau_-^v)^{-1} \tau_+^u=\ze_{u,v}\,\si\,\ze_{u,v}^{-1}
\equal \si^{u,v}=
(\tau_-^v)^{-1} \tau_+^u (\tau_-^v)^{-1} \hbox{\,\,for}\\
&\tau_+^u\equal \ze_{u,v}\,\tau_+\,\ze_{u,v}^{-1}=\tau_+\,
\ze_{1,u},\ \,
\tau_-^v\equal \ze_{u,v}\,\tau_-\,\ze_{u,v}^{-1}=\tau_-\,\ze_{v,1}.
\notag
\end{align}

The action of $\tau_{\pm}$ on $\ze_{u,v}$ by conjugation
is dual to the natural
action of $SL_2(\Z)$ in $\Pi^2$. Recall that $u,v$ are 
arbitrary characters of $\Pi$.
Formulas (\ref{zetauvrel1}) readily
follow from (\ref{zetauvrel}). We call the group generated
by $GL_{\,2}^{\wedge}(\Z)$ and all $\ze_{u,v}$ the 
{\em enhanced projective $GL_2(\Z)$} and use the
notation $GL_{\,2}^{\wedge}(\Z)^\ze$. 
\smallskip

{\sf The coinvariant.}
The projective $GL_2(\Z)$ and the coinvariant
are the main ingredients of our approach. Let us define the
latter. 
Any $H\in \HH$ can be uniquely represented in the form 
\,$\sum_{a,w,b} c_{a,w,b}\, X_a T_{w} Y_b$\, for $w\in W$,
$a,b\in P$ (the PBW theorem, see \cite{C101}). Then the
{\em coinvariant\,} is a functional $\HH\to \C$ 
uniquely defined via the following
substitution in such sums:
\begin{align}\label{evfunct}
\{\,\}_{ev}:\ X_a \ \mapsto\  q^{-(\rho_k,a)},\ 
Y_b \ \mapsto\  q^{(\rho_k,b)},\ 
T_i \ \mapsto\  t_i^{1/2}. 
\end{align}
The key symmetries of the coinvariant are:
\begin{align}\label{evsym}
&\{\,\vph(H)\,\}_{ev}\!=\!\{\,H\,\}_{ev},\, 
\{\,\eta(H)\,\}_{ev}\!=\!\{\,H\,\}_{ev}^\star,\, 
\{\,\iota(H)\,\}_{ev}\!=\!\{\,H\,\}_{ev}.
\end{align}
We use here that $\iota$ naturally
acts in $\HH$:
\begin{align}\label{iotaXY}
\iota(X_b)=X_{\iota(b)},\ 
\iota(Y_b)=Y_{\iota(b)},\ T_i^\iota=T_{\iota(i)},\, 1\le i\le n. 
\end{align}

One has $\{H T_w Yb\}_{ev}=\{H\}_{ev} \chi (T_w Y_b)$,
where $\chi$ is the standard character (one-dimensional 
representation) of the algebra
$\h_Y$, by definition generated by $T_w, Y_b$ for 
$w\in W, b\in P$, which sends 
$Y_b\mapsto q^{(\rho_k,b)}, T_i\mapsto t_i$. 
Therefore $\{\ldots\}_{ev}$ acts via the projection 
$H\mapsto H\!\Downarrow \equal H(1)$ of $\HH\,$
onto the {\em polynomial representation \,}$\v$, which is
the $\HH$\~module induced from $\chi$; 
see \cite{C101,CJ,CJJ} and the next section.
\medskip

\subsection{\bf Macdonald polynomials}
We will begin with the explicit construction
of the polynomial representation, denoted by $\v$
in this work. 

{\sf Polynomial representation.}
It was already defined above as an induced representation.
In detail, it
is isomorphic to $\Z_{q,t}[X_b]$
as a vector space and the action of $T_i(0\le i\le n)$ there 
is given by 
the {\em Demazure-Lusztig operators\,}:
\begin{align}
&T_i\  = \  t_i^{1/2} s_i\ +\
(t_i^{1/2}-t_i^{-1/2})(X_{\al_i}-1)^{-1}(s_i-1),
\ 0\le i\le n.
\label{Demazx}
\end{align}
The elements $X_b$ become the multiplication operators 
and  $\pi_r (r\in O')$ act via the general formula
$\hw(X_b)=X_{\hw(b)}$ for $\hw\in \hW$. Note that $\tau_-$ 
$\eta$ and $\,\iota\,$ naturally act in the polynomial 
representation. 
See formula (1.37) from \cite{CJJ} and (\ref{taumineb})
below for $\tau_-$. We will use the notation $\dot{\tau}_-$
for this action; it is explicitly defined as follows:
\begin{align}\label{taumin-poly}
\dot{\tau}_-(g)\equal\Bigl(\tau_-\bigl(g(X)\bigr)\Bigr)(1) 
\for g\in \v.
\end{align}

As for $\eta$:
\begin{align}\label{eta-conj}
\eta(f)\!=\!f^\star, \hbox{\, where\, } X_b^\star\!=\!X_{-b}, 
(q^\upsilon)^\star\!=\!q^{-\upsilon}, (t^v)^\star\!=\!t^{-v} 
\hbox{\, for\, } \upsilon\in \Q.
\end{align}

Also, the automorphisms $\ze_{u,1}$ from (\ref{zetauv|})
act in $\v$. They can be represented by certain translations
in $x\in \C^n$ for $X=q^x$. Since they preserve $Y_b$, we obtain 
that simple $Y$\~eigenvectors in $\v$ are also
$\ze_{u,1}$\~eigenvectors for any character $u:\Pi\to \C^*$. 
Thus they are also invariant
(up to proportionality) under the action of
$\tau_-^u=\tau_-\ze_{u,1}$ from (\ref{zetauvrel1}).

\smallskip

{\sf Symmetric Macdonald polynomials.}
The standard notation for them is  $P_b(X)$ for $b\in P_+$
(they are due to Kadell
for the classical root systems). The definition is
as follows. Let $c_+$ be such that $c_+\in W(c)\cap P_+$
(it is unique);
recall that $Q_+=\oplus_{i=1}^n \Z_+\al_i$.
 For $b\in P_+$,
\begin{align*}
&P_b\! -\!\!\!\!\sum_{a\in W(b)}\!\!\! X_{a} 
\in\, \oplus_{b\neq c_+\in b-Q_+}\Q(q,t_\nu) X_c 
\hbox{\, and\, }
\lan \,P_b X_{c^{\iota}}\mmu(X;q,t)\,\ran\!=\!0
\\
&\hbox{for such $c$,\, where\, } 
\mmu(X;q,t)\!\equal\!\prod_{\al \in R_+}
\prod_{j=0}^\infty \frac{(1\!-\!X_\al q_\al^{j})
(1\!-\!X_\al^{-1}q_\al^{j+1})
}{
(1\!-\!X_\al t_\al q_\al^{j})
(1\!-\!X_\al^{-1}t_\al^{}q_\al^{j+1})}\,.
\end{align*}
Here and further $\lan f\ran$ is the {\em constant term\,}
of a Laurent series or polynomial $f(X)$;
$\mmu$ is considered
a Laurent series of $X_b$ with 
the coefficients expanded in terms of
positive powers of $q$. The coefficients of
$P_b$ belong to the field $\Q(q,t_\nu)$.
One has (see (3.3.23) from \cite{C101}):
\begin{align}\label{macdsym}
&P_b(X^{-1})\,=\,P_{b^{\iota}}(X)\,=\,
P_b^\star(X),\ \, P_{b}(q^{-\rho_k})=
P_{b}(q^{\rho_k})\\
\label{macdeval}
=\,&(P_{b}(q^{-\rho_k}))^\star=
q^{-(\rho_k,b)}
\prod_{\al>0}\prod_{j=0}^{(\al^{\!\vee},b)-1}
\Bigl(
\frac{
1- q_\al^{j}t_\al X_\al(q^{\rho_k})
 }{
1- q_\al^{j}X_\al(q^{\rho_k})
}
\Bigr).
\end{align}
Recall that $\iota(b)=b^{\,\iota}=-w_0(b)$ for $b\in P$.

DAHA provides an important alternative (operator) 
approach to the $P$\~polynomials; namely, they satisfy
the (defining) relations
\begin{align}\label{macdopers}
L_f(P_b)=f(q^{-\rho_k-b})P_b,\  L_f\equal f(X_a\mapsto Y_a)
\end{align}
for any symmetric ($W$\~invariant) polynomial 
$f\in \C[X_a,a\in P]^W$. Here $b\in P_+$ and the coefficient
of $X_b$ in $P_b$ is assumed $1$.
\smallskip

{\sf Nonsymmetric Macdonald polynomials.}
The $P$\~polynomials are $t$-symmetrizations of the
{\em nonsymmetric Macdonald polynomials\,} $E_b\in \v$.
For any $b\in P$, we define them as follows:
\begin{align}\label{YaEb}
Y_a&(E_b)\,=\,q^{-(a,b+w_b^{-1}(\rho_k))}E_b, \hbox{ where }
w_b\in W \hbox{\, is \,} \\
\hbox{a unique} &\hbox{ element of maximal length 
such that\, } w_b(b)\in P_+.\notag
\end{align}
The element $w_b$ for $b\in P_+$ here is the element of maximal 
length (an involution) in the 
centralizer of $b\,$ in $W$. The normalization of $E$ in
(\ref{YaEb}) is by 
the condition that the coefficient of $X_b$ in $E_b$ is $1$. 
For $b\in P_+$, one has that
$E_b-X_b\in \oplus_{b_+\neq c_+\!\in b-Q_+}\,\Q(q,t_\nu) X_c$.
See \cite{Mac} (for $k_{\sht}=k_{\lng}\in \Z_+$) and (6.14) 
from \cite{C103} or (3.3.14) from \cite{C101};
the differential version is due to Opdam with a participation
of Heckman. The Macdonald
conjectures for them were extended from the symmetric
case and justified in \cite{Ch4}; as a matter of fact, these
conjectures become significantly simpler in the nonsymmetric 
setting.

We note that all monomials in $E_b$
or  $P_b$ are in the form $X_{b+a}$ for $a\in Q$,
i.e. have coinciding images in $\Pi=P/Q$. This readily
follows from the commutativity of $Y_c$ with the automorphisms
$\ze_{u,1}$ from (\ref{zetauv|}).
\smallskip

\subsection{\bf Evaluation formula}
One of the key formulas in this work is the following evaluation:
\begin{align}\label{macde-eval}
&E_{b}(q^{-\rho_k})=
q^{-(\rho_k,b)}
\prod_{\al>0}\prod_{j=1}^{(\al^{\!\vee},b)-1}
\Bigl(
\frac{
1- q_\al^{j}t_\al X_\al(q^{\rho_k})
 }{
1- q_\al^{j}X_\al(q^{\rho_k})
}
\Bigr) \for b\in P_+.
\end{align}

For any $b\in P$ and $w_b$ from (\ref{YaEb}),
\begin{align}
E_{b}(q^{-\rho_k}) \ =\ q^{-(\rho_k,b_+)}&
\prod_{\{\al,j\}\in \La^+_b}
\Bigl(
\frac{
1- q_\al^{j}t_\al X_\al(q^{\rho_k})
 }{
1- q_\al^{j}X_\al(q^{\rho_k})
}
\Bigr),
\label{ebebs}\\
\La^+_b\ \equal\
 \{ \{\al>0,j>0\}\ \mid\
&( b_+, \al^\vee )>j> 0
\hbox{ for }  w_b^{-1}(\al)\in R_+,\notag \\
&( b_+, \al^\vee )\ge j> 0
\hbox{ for }  w_b^{-1}(\al)\in R_-\}.\notag
\end{align}

Formula (\ref{ebebs}) is the Macdonald
{\em evaluation conjecture\,}
in the nonsymmetric variant from \cite{Ch4};
see also formulas (6.33), (7.15) from 
\cite{C103} and (3.3.45) from \cite{C101}.
The set $\La_b^+$ appears here due to the
following (ibid.):
\begin{align}\label{lapibdef}
&\{\tilde{R}_+\ni [-\al,j\nu_\al]\,\mid\, \{\al,j\}\in \La_b^+\}
=\La(\pi_b), \where\\
\pi_b&\equal b w_b^{-1} w_0,\ 
\La(\hw)=\tilde{R}_+\cap \hw^{-1}(-\tilde{R}_+) \for 
\hw\in \hW.\notag
\end{align}
Then the DAHA intertwining operators
are used to justify (\ref{ebebs}).
\smallskip

{\sf Spherical normalization.}
We call $E_b^\circ\equal E_b/E_b(q^{-\rho_k})$ for $b\in P$
nonsymmetric {\em spherical polynomials\,}. Accordingly,
$P_b^\circ\equal P_b/P_b(q^{-\rho_k})$ for $b\in P_+$, and
\begin{align}
&P_{b}(q^{-\rho_k})=\Pi_R^b\,E_{b}(q^{-\rho_k}),\  \,
\Pi_R^b\equal \!\!\!\prod_{\al>0,(\al,b)>0}
\frac{
1\!-\! t_\al X_\al(q^{\rho_k})
 }{
1\!-\! X_\al(q^{\rho_k})
},\label{Pbrho}\\
&P_{b}(q^{-\rho_k})=
q^{-(\rho_k,b)}
\prod_{\al>0}^{(\al,b)>0}\,\prod_{j=0}^{(\al^{\!\vee}\!,b)-1}
\Bigl(
\frac{
1- q_\al^{j}t_\al X_\al(q^{\rho_k})
 }{
1- q_\al^{j}X_\al(q^{\rho_k})
}
\Bigr).\label{Pbrho'}
\end{align}

The product $\Pi_R^b$ 
becomes over all $\al>0$
if {\em all\,} $\om_i$ appear in $b\in P_+$; it is then the 
{\em Poincar\'e polynomial\,} $\Pi_R$. See formulas 
(7.15)-(7.19) from \cite{C103} concerning these polynomials 
with different root lengths. Recall that $t_\al=q_\al^{k_\al}$ 
here (and there); see (\ref{talqal}) above.
\smallskip

{\sf The symmetrization.}
The symmetrization
relation between $E_b^\circ$ and $P_{b_+}^\circ$ is as follows.
For any $b\in P$,
\begin{align}\label{nontosym}
P_{b_+}^\circ\!=\!(\Pi_R)^{-1}\mathscr{P}_+
(E^\circ_b)\, \hbox{\, for\, }\, 
\mathscr{P}_+\!\equal\!\sum_{w\in W}
t_{\sht}^{l_{\sht}(w)/2}
t_{\lng}^{l_{\lng}(w)/2}T_w,
\end{align}
where $l_{\sht}, l_{\lng}$ count correspondingly the number 
of $s_i$ for short and long $\al_i$ in any reduced decomposition 
$w=s_{i_l}\cdots s_{i_1}$. We check that the right-hand side here
is proportional to $P_{b_+}$ and then evaluate at 
$\,q^{-\rho_k}$ by
applying the general formula
$\{T_i (f)\}_{ev}=t_i^{1/2}\{f\}_{ev}\,$ for $i>0$ and
any $f(X)$.

Let us provide the following particular case of (\ref{nontosym})
for the standard Macdonald polynomials: 
\begin{align}\label{nontosymdom}
&\prod_{\al>0,(\al,b)=0}
\frac{
1- t_\al X_\al(q^{\rho_k})
 }{
1- X_\al(q^{\rho_k})
}
\,P_{b}=\mathscr{P}_+(E_{b}) \for b\in P_+,
\end{align}  
For generic dominant $b$ (when all $\om_i$
occur in its decomposition), the coefficient of 
proportionality here is $1$.
\smallskip

Let us also give the
formulas for the natural action of $\tau_-$ in $\v$
in the basis of $E$\~polynomials. See e.g., formula (1.37) 
from \cite{CJJ}. Recall that we use the notation 
$\dot{\tau}_-$ for this action:
\begin{align}\label{taumineb}
\dot{\tau}_-(E_b)=q^{-\frac{(b_+\,,\,b_+)}{2}-
(b_+\,,\,\rho_k)}\,E_b
\for b\in P.
\end{align}
This readily gives that 
$\dot{\tau}_-(P_b)=q^{-(b,b)/2-(b,\rho_k)}P_b$ for $b\in P_+$;
for instance, use (\ref{nontosymdom}).
We assume that $E_b,P_b$ are  well defined. 
\medskip

\setcounter{equation}{0}
\section {\sc Integral forms}
\subsection{\bf A preliminary version}
The technique of intertwining operators applied to the 
$E$\~polynomials
results in the following {\,\em existence criterion\,}.
Given $b\in B$ and assuming that $q$ is not a root of unity,
the following polynomials are {\em $q,t$\~integral\,}:
\begin{align}
&\n_b\, E^\circ_b \for \n_b=\prod_{\{\al,j\}\in \La^+_b}\,
\bigl(
1- q_\al^{j}t_\al X_\al(q^{\rho_k})\bigr),
\label{ebebnze}\\
&\d_b\,E_b \for \d_b=\prod_{\{\al,j\}\in \La^+_b}\,
\bigl(
1- q_\al^{j}X_\al(q^{\rho_k})\bigr).\label{ebebdze}
\end{align}
See \cite{C1}, Corollary 5.3 or (6.13) from \cite{C103}.
In the case of $A_n$, this observation is due to F.~Knop.
By {\em integrality\,} (we will mainly omit $q,t$ here),
we mean that their coefficients belong to 
$\Z_{q,t}=\Z[q^{\pm 1},t^{\pm 1}]$. Recall that
$t_\al=q_\al^{k_\al}=q^{\nu_\al k_\al}$.

These two polynomials obviously coincide because
$\n_b/\d_b=E_{b}(q^{-\rho_k})$; see (\ref{macde-eval}).
So we actually have only one statement here.
Switching to $P_b$ for $b\in P_+$, one obtains
the integrality of $\n_b\,P_b^\circ$, $\d_b\,P_b$. 
Note that all binomials in $\n_b,\, \d_{b}\,$ contain 
$q$ and $t$ (both).

Using \cite{C103}, the existence criterion for $E_b^\circ$
can be made much sharper. Namely, the multiplicities 
of the binomials 
$(1-q^\bullet t^\bullet)$ in $\n_b$ (generally large)
can be reduced to $1$, with minor adjustments for $D_{2m}$.
Also, some of these binomials can be further reduced
to their proper factors. We will justify this mainly following 
\cite{C103}. Such
an improvement was not done there; it is helpful
in decomposition theory of $\v$. 
\smallskip

\subsection{\bf The integrality theorem}
Assuming that $q,t_{\lng},t_{\sht}$ are generic, the evaluation 
formula for $E_{b}(q^{-\rho_k})$
stabilizes to an infinite product as $b\in P_+$ 
of {\em integral polynomials\,} as 
$b$ becomes sufficiently large. Namely, it approaches the 
product of $\Pi_{\tR}$ from Part $(iii)$, 
Theorem 11.8 of \cite{C103}, which is extended
to all  $k_{\sht}+\Z_+, k_{\lng}+\Z_+$ and
then divided by the Poincar\`e polynomial $\Pi_R$. 

This product has a canonical presentation 
$E_{\infty}(q^{-\rho_k})=\n_\infty/\d_\infty$,
where the factors of $\n^\infty$ and $\d^\infty$ are products of
binomials in the form 
\begin{align}\label{binaffexp}
&(1-q^{\kapp}) \,\for\, \kapp=j+j_{\sht}k_{\sht}+
j_{\lng}k_{\lng},\,\where\\ 
&j\in \N,\,  j_{\sht}/\nu_{\sht}\in \Z_+
\ni j_{\lng}/\nu_{\lng}\hbox{\, and\, }
j_{\sht}+j_{\lng}>0.\notag
\end{align}

The exponents $\kapp=j+j_{\sht}k_{\sht}+j_{\lng}k_{\lng}$
for the binomials of $\n_\infty$ are called 
{\em affine exponents\,}. There is a deviation from
the usage of this term in \cite{C103}. Namely,
we omit the binomials with $j=0$  and 
(more importantly) 
consider the {\em stable case\,} i.e. do not bound 
$j$. Thus the exact name must be  
{\em affine stable $q$-dependent exponents\,},
but we will simply call them {\em affine}.
We call an affine exponent $\kapp\,$  {\em rational\,} if it 
is not divisible by any exponents $\kapp'$ that occur in 
$\d_\infty$.
\smallskip

In the case $AD\!E$, he denominator $\d_\infty$ is very
simple: $\prod_{j=1}^\infty
(1-tq^j)$. For these and arbitrary root systems, each non-rational  
$\kapp$ canonically corresponds to its divisor  
$\kapp'$ in $\d_\infty$ ($\kapp'$  can have multiplicities there).
We will provide below explicit lists. This
correspondence will be denoted: 
$\kapp\rightsquigarrow \kapp'$.
Importantly, the multiplicity of each $\kapp$ in the list of
affine exponents is one (including $BC\!F\!E$)
unless for the series 
$\kapp=2mk+j$ in the case of $D_{2m}$, which is of
multiplicity $2$. 
 
Note that the binomials
in (\ref{binaffexp}) are in terms of integral (non-negative)
powers of $q,t_{\sht},t_{\lng}$, so we 
use $\kapp$ (in terms of $k_\nu$)  for the sake of simplicity
of notations. For $b\in P$ and $\al\in R_+$, we set
$\de_\al(b)= 1 \hbox{ if } w_b^{-1}(\al)\in R_-$ and $\,0$ 
otherwise; if $b\in P_-$, then $\de_\al(b)=1$ for any 
$\al\in R_+$.

\begin{theorem}\label{REGTHM}
We fix $b\in P$. Let us
consider only the affine exponents $\kapp$ that occur in
$\n_b$, which means the existence of $\al\in R_+$ such that
$\kapp=k_\al+(\al,\rho_k)+j$ and 
$( b_+, \al^\vee )+\de_\al(b) >j>0$.
If $b\in P_+$, the inequality $( b_+, \al^\vee )>0$ and the
relation $b=b_+$ imply that $\de_\al(b)=0$ here.
Define $\tilde{\n}_b$ as the product 
over all such $\kapp$ of the following factors:

(a) the binomials \, $(1-q^{\kapp})$\, for the rational exponents
$\kapp$,
 
(b) and the ratios $\,(1-q^{\kapp})/(1-q^{\kapp'})$ for 
non-rational $\kapp$, 

\noindent
where we use the above correspondence 
$\kapp\rightsquigarrow \kapp'$. There is a reservation
for $\kapp=2mk+j$ in the case 
of $D_{2m}$; we take the 
corresponding factor from $(a)$ or $(b)$ twice if there exist at 
least two $\al>0$ such that $(\al,\rho)=2m-1$ subject to 
$( b_+, \al)+\de_\al(b) >j>0$. 

Then we set 
$\tilde{\n}_b\equal\hbox{gcd}(\n_b,\n_\infty/\d_\infty)$,
and claim that 
the coefficients of the polynomials $\tilde{E}_b\equal
\tilde{\n}_b E_b^\circ$ are $q,t$\~integral for $b\in P$,
as well as those of  $\tilde{P}_b\equal
\tilde{\n}_b P_b^\circ$ for $b\in P_+$. Here we can assume
that $j\le N$ in $\n_\infty/\d_\infty$
for sufficiently large $N$.\sq
\end{theorem}

\subsection{\bf Using the duality}
The justification will be based on the DAHA-duality.
We can assume that the parameters $q,t_{\sht},
t_{\lng}$ are generic. Indeed, the denominators of the
coefficients of $\tilde{E}_b$ can be only in the form
(\ref{binaffexp}) for certain $\kapp$, so it suffices
to establish the absence/reduction of some of them for 
generic parameters only.

Assuming that $q$ is not a root of unity, the key fact is
that for any given values of $t_{\sht},t_{\lng}$ in $\C^*$,
the polynomials $E_c$ become regular (i.e. their coefficients
become regular) for sufficiently large $c\in P$. The latter
means that $(c_+,\al_i)>C$ for proper $C\gg 0$
and all $1\le i\le n$.
 
See formulas (9.16) and (9.17) from \cite{C103} (right before
Section 10 ``The structure of $\v_c$"). It is demonstrated there
that the {\em generalized\,} $Y$\~eigenspace in
$\v$ containing $E_c$ becomes 
$1$-dimensional for $c$ satisfying certain inequalities
for its coefficients; this readily results in the
existence of the corresponding $E_c$ (i.e. in its 
$q,t$\~integrality).

Actually Theorem 10.3 there (in Subsection 10.3
``The semisimple submodule") establishes that $\v$ has 
a canonical nonzero semisimple submodule. However it
was proven in this work under certain restriction and this
(much stronger) fact is not necessary for our proof.
Vice versa, our result here is useful for the decomposition
of the polynomial representation, which we will try to
address in other publications.
\smallskip

Then we apply the {\em duality theorem\,} from \cite{Ch4};
see also (6.30) from \cite{C103}). It states that
for any $b,c\in P:$
\begin{align}\label{ebdual}
E_b^\circ(q^{c_{\#}})=E_c^\circ(q^{b_{\#}}),\  
b_\# \equal b+w_b^{-1}(\rho_k),\ Y_b(E_c)=q^{-(c_\#,b)}E_c&,\\
P_b^\circ(q^{c+\rho_k})= P_c^\circ(q^{b+\rho_k}),\, \,
P_b(Y^{-1})(P_c)=P_b(q^{c+\rho_k})P_c \for b,c\in P_+&.\notag
\end{align}
Note that 
$w_b=w_0$ for $\,b\in P_-$ (including
$b=0$) and therefore $b_\#=b-\rho_k$ for such $b$.
Also, $b_\#=b+\rho_k$ for generic $b\in P_+$ 
(such that $(b,\al_i)>0$ for $i=1,\ldots, n$).
  
The duality is the main motivation
of the {\em spherical normalization\,}. 
{\em All\,} Macdonald's conjectures can be 
deduced from the duality. See Proposition 6.6 in \cite{C103} 
and the end of Section 7 there concerning the deduction of the
Evaluation Conjecture, Norm Conjecture and
Constant Term Conjecture (in this order) from the duality. 
One can add to this list
the Pieri rules and the difference Mehta-Macdonald
Conjecture from \cite{C4}, which also result from the duality. 

\subsection{\bf Concluding the proof}
Fixing $q,t,b\,$ and assuming that $c$ is sufficiently large 
(and that $q$ is not a root of unity), we conclude
that $(\n_c/\d_c) E_b^\circ(q^{c_\#})$ is regular for such $c$. 
This regularity formally results in the 
regularity of the coefficients of this polynomial at such $q,t$; 
let us demonstrate this.
 
Not any set of values can be generally used for such an
implication, even if it is arbitrarily large (infinite) and
sufficient to ``catch" all coefficients of our polynomial. 
Say, the non-integral polynomial $F(X)=(1-X)/(1-q)$ has 
integral values at any $X=q^n (n\in \Z)$; such 
values are of course sufficient to recover (the coefficients of)
{\em any\,} polynomial
$F(X)$. The Lagrange interpolation formula readily 
clarifies which sets are good for this.

Let us consider one variable $X$ and assume that a 
Laurent polynomial $F(X)$ can be uniquely 
recovered from its values at 
certain pairwise distinct points $X=q^{a_i}$.  Then 
the denominators of the
coefficients of $F(X)$ can be only divisors of 
$q^{a_i}-q^{a_j}=0$, where we treat $a_i$ as parameters. 
Thus the regularity of $F(q^{a_i})$ implies that for
the coefficients of $F$ if there are no common factors
of {\em all\,} differences $q^{a_i}-q^{a_j}$.
\smallskip
 
For instance, let us establish the regularity
of $(\n_c/\d_c) E_b^\circ$ for the root system $A_1$
and sufficiently large $c$. 
We set $X=X_{\om}$ for $\om=\omega_1$ and take 
$a_j=(j\om)_\#=\sgn(j) (k+|j|)/2\,$ for $i\in \Z$ with
sufficiently large $|j|$. Note that $\sgn(0)=-1$ here,
but we need only large $|j|$. Thus 
$q^{a_i}-q^{a_j}$ for $i>j\,$ is  
\begin{align*}
&\hbox{either\, }
t^{\hbox{\tiny\rm sgn}(i)/2}(q^{i/2}-q^{j/2}) \for i j>0\,\\ 
&\hbox{or\, } t^{1/2}q^{i/2}-t^{-1/2}q^{-j/2} \for i>0>j.
\end{align*}

Let the set $\j=\{j\}$ be with sufficiently large $|j|$
and generic enough for the recovery of the 
coefficients of $E_b^\circ$ from its values
at $q^{c_\#}$ for $c=j\om$. Let us pick $j$ here
with the same $\sgn(j)$; then the corresponding differences 
in the case $ij>0$ are relatively prime to 
any binomials from (\ref{binaffexp}) due to 
$j_{\sht}=1$ there. The series $i>0>j$
can be used here too provided $i,j$ are 
sufficiently large (depending on the set of
$\kapp$). 
\smallskip

Given $b$, such a choice of $\j$ proves that the 
coefficients of $(\n_c/\d_c) E_b^\circ$ are regular.
Also, since $|j|$ are sufficiently large, we can replace 
$\n_c/\d_c$ by $\n_\infty/\d_\infty$,  
which is its stabilization as $|j|\to \infty$. Since 
it is $q,t$\~integral (a product of $q,t$\~polynomials),
we conclude that $(\n_\infty/\d_\infty) E_b^\circ$
is integral too.  Finally, 
we can restrict ourselves only to the binomials in
$\n_\infty$ that appear in $\n_b$; we use that all
factors from $(a,b)$ are irreducible for
generic $t$ in type $A_1$. 
\smallskip

The case of any root system is actually no different,
but we must avoid using the Lagrange interpolation.
Let us take $c$ such that the polynomials $E_c$ exist for any 
zeros of $\n_b$. We can assume that $c\in P_-$
and therefore $c_\#=c-\rho_k$. Let us substitute
$X\mapsto Xq^{-\rho_k}$, $X_a\mapsto q^{-(\rho_k,a)}X_a$;
this will not change the (non)regularity of the
coefficients. Upon this rescaling,  values of
$F(X)=(\n_c/\d_c) E_b^\circ (X q^{-\rho_k})$
at $q^c$ with sufficiently large $c$ are
all regular for such $c$. We pick sufficiently many
of them ensuring that the coefficients of $F(X)$ can be
recovered from its values in this set (assume its
minimality).
\smallskip
 
Generally the coefficients of $F$ are recovered by 
applying the inverse of the matrix transforming the 
coefficients of $F$ to its values. The entries 
of this matrix are powers of $q$. Therefore the zeros 
of its (Vandermonde-type) determinant cannot be zeros 
of any binomials in (\ref{binaffexp})
because $\kapp$ there always contain $k_{\sht}$ or $k_{\lng}$. 

Then we proceed as for $A_1$
and obtain that the multiplication of $E_b^\circ$
by $\n_\infty/\d_\infty$ is actually sufficient to make this 
polynomial $q,t$\~integral.  We also know that $\n_b E_b^\circ$
is integral. Thus we can keep in $\n_\infty/\d_\infty$ 
only the factors that come from the binomials
from $\n_b$, which reduces $\n_\infty/\d_\infty$ to   
$\tilde{\n}_b$.
\smallskip

Here we can consider a sufficiently
large truncation $(\n_\infty/\d_\infty)^\dag$
of $\n_\infty/\d_\infty$ and treat it and $\n_b$ as
polynomials of one variable $q$ with the coefficients
that are rationals
in terms of $t_{\sht},t_{\lng}$. Given any fixed
$t$\~part (i.e. for any
fixed $j_{\sht},j_{\lng}$ in 
$j_{\sht}k_{\sht}+j_{\lng}k_{\lng}),$ the corresponding
factors from $(a,b)$ in the theorem are relatively 
prime to each other. Indeed, such binomials/ratios
must have then different powers $q^j$ and we can use that
$gcd$\, of the Kummer-type polynomials
$(1-uq^l)$ and $(1-u q^m)$ for $l\neq m$
is $1$ if $u$ is not a power of $q$. 

Let us take here a maximal pair $\{j_{\sht},j_{\lng}\}$
(in the {\em poset\,} of pairs). 
For any divisor of the corresponding factor from
$(a,b)$ that divides \,gcd$(\n_b,(\n_\infty/\d_\infty)^\dag)$  
and divides $\tilde{\n}_b$, the whole factor from $(a,b)$
must divide both of them. Therefore one can multiply
$E_b^\circ$  by the 
product of all such factors for a given maximal 
$\{j_{\sht},j_{\lng}\}$ and then proceed by (double) 
induction with respect to 
the (remaining) $\{j_{\sht},j_{\lng}\}$. This concludes the 
proof of the theorem.
\sq

\subsection{\bf Affine exponents}
We will provide here the lists of affine (stable 
and  $q$\~dependent) exponents. We follow \cite{C103}, 
where the formulas must be naturally adjusted to the stable 
setting. 
We use the {\em Coxeter exponents\,}  $m_i$ (see e.g., \cite{Bo}), 
determined from the relation 
\begin{align}\label{tdegr}
\prod_{\al\in R_+} \frac{1- t^{1+(\al,\,\rho^\vee)}}
{1-t^{(\al,\,\rho^\vee)}}=
\prod_{i=1}^n \frac{1-t^{\,m_i+1}}
{1-t},\ \ 2\rho^\vee=\sum_{\al\in R_+}\,\al^\vee.
\end{align}
This product is the {\em classical Poincar\`e polynomial\,}, 
which is for $t_{\sht}=t=t_{\lng}$ in terms of
our $\Pi_R$; see (\ref{Pbrho}).
One has that $m_1+\ldots +m_n=|R_+|$ 
and $(m_i+1)\cdots(m_n+1)=|W|$; $(m_i+1)$ are called the
{\em degrees} of the root system $R$. 

{\sf The $AD\!E$\~case.}
We set $t_{\sht}=t=t_{\lng},\, t=q^k$ 
in the simply-laced case. Then
\begin{align}\label{ntilde-an}
&\frac{\n_\infty}{\d_\infty}=
\prod_{\al\in R_+}\prod_{ j=1}^{\infty}
\Bigl(\frac{
1- q^{j} t^{1+(\al,\rho)}}
{1- q^{j}t^{(\al,\rho)}}\Bigr)=
\prod_{i=1}^n \prod_{j=1}^{\infty} 
\Bigl(\frac{1-q^{j} t^{\,m_i+1}}
{1-q^j t}\Bigr).
\end{align}

Therefore the list of affine exponents is
\begin{align}\label{affexade}
AD\!E:\ \ &\{ \kapp=(m_i+1)k+j+1\, \mid\, i=1,\ldots, n,\ \,
j\in\Z_+\}, \\
\hbox{unless\, } &\{(m_i+1)(k+j+1)\rightsquigarrow\, \kapp'=k+j+1,
\, j\in \Z_+\}.\notag
\end{align}

Recall that non-rational exponents $\kapp$ (the second line) 
result in the corresponding quotients in $\tilde{\n}_\infty$
with respect to the map $``\rightsquigarrow"$ from
$\kapp$ to their (canonical) divisors $\kapp'$ from the set 
of exponents of $\d_\infty$. 

{\sf The case of $B_n$.}
The are two series of rational exponents in this case:
\begin{align}\label{redbdegrees}
&\{\kapp=2mk_{\lng}+2j+2\,\mid\, 2\le m\le n, j+1\not\in m\Z_+\},\\
&\{\kapp=2mk_{\lng}+2k_{\sht}+2j+1\,\mid\, 0\le m<n,\ \, j\ge 0\}.
\notag
\end{align}
Here and below $j$ will be always from $\Z_+$. 
The list of the non-rational $\kapp$ (with their
canonical divisors $\kapp'$ from $\d_\infty$) is as follows:

{\small
\begin{align}\label{bexpnonrat}
&\{4mk_{\lng} +2k_{\sht}+2j+2\, \rightsquigarrow\, 
2mk_{\lng} +k_{\sht}+j+1 \,\mid\, 1 \le m < n\},\\
&\{m(2k_{\lng} +2j+2)\, \rightsquigarrow\, 
\kapp'=2k_{\lng} +j+1 \,\mid\, 2 \le m \le n\}.  
\end{align}
}

{\sf The case of $C_n$.} There are
two rational series as well: 
\begin{align}\label{redbdegreesc}
&\{\kapp=mk_{\sht}+j+1\,\ \mid\,\ 2\le m\le n \for 
j+1\,\not\in\, m\Z_+\},\\
&\{\kapp=\de_m(2mk_{\sht}+2k_{\lng}+2j+1)\,\mid\, 
0\le m<n,\, j\ge 0\},\notag
\end{align}
where $\de_m=2\,$ if $\,m<n/2\,$ and $1$ otherwise.
\smallskip
 
Note that up to proportionality, the rational $C_n$\~exponents
coincide with those of type $B_n$ upon the transposition 
$k_{\sht}\leftrightarrow k_{\lng}$. This is related to
our usage of ``rational" here; in the limit to {\em rational
DAHA\,}, the rational exponents (and only them)
can be seen and the proportionality does not matter. So 
the affine exponents for $B$ and $C$ coincide in this 
limit. (We mention that we somewhat fixed (8.3) 
from \cite{C103}; the main formula (8.6) 
for {\em all\,} affine exponents in the $C$\~case is
correct here).

The non-rational exponents and their divisors $\kapp'$
(from $\d_\infty$) are as follows in the case of $C_n$:

{\small
\begin{align}\label{cexpnonrat}
&\{4k_{\lng} +2mk_{\sht}+2\ep_m(j+1)\,\rightsquigarrow\,
2k_{\lng} +mk_{\sht}+ \ep_m(j+1)\},\\
&\,\, \hbox{\normalsize where\, } 1 \le m < n,\ \  
\ep_m=1 \hbox{\normalsize\, for odd\,\,}\, m
\hbox{\normalsize\,\, and \,{\small $2$}\, otherwise},\notag\\
&\{\kapp=m(k_{\sht}+j+1)\rightsquigarrow\,\kapp'=k_{\sht}+j+1\}
\hbox{\normalsize\, for\, $2\le m\le n$}.
\notag
\end{align}
}

{\sf The case of $G_2$}. There are
$3$ series of  rational exponents:
\begin{align}\label{redgdegrees}
\{2 k_{\sht}+2j+1\},\,
\{6 k_{\lng}+6j+3\},\, 
\{3 k_{\lng} + 3 k_{\sht}+3j+1,2\},
\end{align}
where $j\in \Z_+$ as in all formulas here.
By {\em series\,}, we mean here the sequences of $\kapp$
with coinciding $k$\~components.

The non-rational series are:
{\small
\begin{align}\label{g2affined}
\{2 (k_{\sht}+j+1)\rightsquigarrow k_{\sht}+j+1\},\ 
\{6 (k_{\lng}+j+1)\rightsquigarrow 3 (k_{\lng}+j+1)\}&,\notag\\
\{9k_{\lng}+ 3k_{\sht}+3j+3 \,\rightsquigarrow\,
3 k_{\lng}+ k_{\sht}+j+1\}\, \hbox{\normalsize 
\,\,for\,\,} j\ge 0&.
\end{align}
}
\smallskip

{\sf The case of $F_4$.}
In this case, the 
rational affine exponents are:
\begin{align}\label{redfdegrees}
&\{2k_{\sht}+2j+1\},\, \{3k_{\sht}+3j+1,2\},\notag\\
&\{4k_{\lng}+4j+2\},\, \{6k_{\lng}+2j+2,4\},\notag\\ 
&\{2k_{\lng} + 4k_{\sht}+2j+1\}, \,
\{4k_{\lng} + 4k_{\sht}+2j+1,2,3\},\notag\\
&\{6k_{\lng} + 6k_{\sht}+2j+1\},\,
\{8k_{\lng} + 4k_{\sht}+4j+2\}.
\end{align}

The non-rational $F_4$\~series are:
{\small
\begin{align}\label{f4affined}
\{2(k_{\sht}+j+1)\, &\rightsquigarrow \ \ \ k_{\sht}+j+1\},\,\ 
\{3(k_{\sht}+j+1)\ \, \rightsquigarrow\ \, 
k_{\sht}+j+1\},\notag\\
\{4(k_{\lng}+j+1)&\,\rightsquigarrow 2(k_{\lng}+j+1)\},\,
\{6(k_{\lng}+j+1)\,\rightsquigarrow 2(k_{\lng}+j+1)\},\notag\\
\{4k_{\lng}& + 4k_{\sht}+4j+4\rightsquigarrow
2k_{\lng} + 2k_{\sht}+2j+2\},\notag\\
\{8k_{\lng}& + 2k_{\sht} +2j+2\ \,\rightsquigarrow\ \,  
4k_{\lng} + k_{\sht} +j+1\},
\notag\\ 
\{8k_{\lng}& +4k_{\sht}+4j+4\ \,\rightsquigarrow\ \, 
2k_{\lng} +k_{\sht}+j+1\},
\notag\\  
\{12k_{\lng}\!& + 6k_{\sht}+2j+2\, \rightsquigarrow\ \,
6k_{\lng} + 3k_{\sht}+j+1\}.
\end{align}
}

Note that the exponents $\kapp'=2k_{\lng}+k_{\sht}+j+1$ 
from $\d_\infty$ can potentially divide 
$8k_{\lng} +4k_{\sht}+4j+4$ and 
$12k_{\lng} + 6k_{\sht}+2j+2$, but $6k_{\lng} + 3k_{\sht}+j+1$
can be used only to divide the latter. Thus the correspondence
$\rightsquigarrow$ is still unique in this case. 
The series  $\kapp'=k_{\sht}+j+1$
and $\kapp'=2(k_{\lng}+j+1)$ are of the same kind.
Such examples exist for other root systems, 
but all instances come from the
Poincar\`e polynomial $\Pi_R$ (upon translations
of $k_{\sht},k_{\lng}$ by $Z_+$). 

\smallskip

\subsection{\bf 
\texorpdfstring{{\mathversion{bold}$J$}}{J}-polynomials in 
type \texorpdfstring{{\mathversion{bold}$A$}}{A}}\label{sec:A-n}
In the case of $A_n$, let us calculate $\,\tilde{\n}_b\,$ for
dominant $\,b=\sum_{i=1}^n b_i \om_i$, i.e. with
$\,b_i\ge 0$ for all $\,i$. The usage of Young diagrams
is standard (and convenient) here; they are defined as
follows:
$$
\la=\la(b)=\{\la_1=b_1+\ldots+b_n,
\la_2=b_2+\ldots+b_n,\ldots, \la_n=b_n\}.
$$

We set $R_+\ni \al=\vep_{lm}=\vep_l-\ep_m$ for
$1\le l<m\le n+1$ in the standard basis
$\{\vep_l\}\in \R^{n+1}$. Then 
\begin{align}\label{omviavep}
&\om_i=\vep_1+\ldots+\vep_i-\frac{i}{n+1}(\vep_1+\ldots+\vep_{n+1})
\for i=1,\ldots,n,\\ 
&b=\sum_{i=1}^n \la_i \vep_i-
|\la|\,\bigl(\vep_1+\ldots+\vep_{n+1}\bigr)/(n+1)
\for |\la|=\sum_{i=1}^n \la_i,\notag\\
&\rho=\om_1\!+\ldots+\!\om_n=\frac{1}{2}\bigl(
(n\!-\!1)\vep_1+(n\!-\!3)\vep_2+\ldots+(1\!-\!n)
\vep_n\bigr).\notag
\end{align}
One has: $(b,\vep_{lm})=b_l+\ldots+b_{m-1}$,
$(b,\rho)=|\la|/2-\la_1/2$. Also, 
$b^2\equal(b,b)=\sum_{i=1}^n \la_i^2-|\la|^2/(n+1)$.


Let us calculate the  set
$\La_b^+$ from (\ref{ebebs}). Since $b\in P_+$,
the root  $w^{-1}(\al)$ there can be only
from  $R_+$. Thus:
\begin{align*}
&\La^+_b\ \equal\
\{ \{R_+\ni\al=\vep_{lm},j>0\}\ \mid\
b_l+\ldots+b_{m-1}>j> 0\}.
\end{align*}

{\sf Our construction.}
In Theorem \ref{REGTHM}, we need to know
all possible $\kapp=(\al,k\rho)+j=
k(m-l)+j\,$ for such $\{\al,j\}$. For $l+p-1\le n$,
they are  as follows:
\begin{align}\label{kapp-an}
&\{kp+j\,\mid\, 1\le p\le n,\, 0\le j <j_p\equal
\max_l (b_l+\ldots+
b_{l+p-1})\};\\ \label{kapp-an1}
&\hbox{if\, \,} b_n\ge b_{n-1}\ge\ldots\ge b_1,
\hbox{\, then\, }
j_p=\la_{n-p+1} \for 1\le p\le n,\\
&\hbox{if\, \,} b_1\ge b_1\ge\ldots\ge b_{n},
\hbox{\, then\, }
j_p\,=\,\la_1-\la_{p+1} \for \la_{n+1}=0\notag.
\end{align}

The sequences  $\{j_p\}$ are obviously the same for  
$b$ and $b^\iota$. The coincidence of the
factors $\tilde{n}_b$ for $b$ and $b^\iota$ is a general
fact for any root systems, 
reflecting that  $P_b$ and $P_{b^{\iota}}$ 
(as well as $E_b$ and $E_{b^{\iota}}$ for $b\in P$) have the same 
$q,t$\~integrality properties. Note
that the Young diagram $\la(b^{\,\iota})$ for $b^\iota$
is the complement of $\la=\la(b)$ in the $(n+1)\times 
\la_1$\~rectangle (naturally containing $\la$)
rotated by $180^\circ$.  
  
The inequalities in (\ref{kapp-an1})
are sufficient but of course not necessary
for the relations $j_p=b_n+\ldots+b_{n-p+1}$. For instance,
take any non-increasing sequences $B^{(s)}=
\{b^{(s)}_m\ge\ldots\ge b^{(s)}_1\}$ of the same size $m$ 
and such that $B^{(s)}\ge B^{(s+1)}$ element-wise.
Then $j_p=\la_{n-p+1}$ for all $p$ in the  
sequence $\{B^{(1)},B^{(2)},\ldots\}$ of the
size $n=m+m+\ldots$. 
\smallskip

Finally, using (\ref{ntilde-an},\ref{affexade})
and the formula $m_p=p$ for the Coxeter exponents
of type $A_n$, we obtain that
\begin{align}\label{regep-an}
&\tilde{\n}_b=
\prod_{p=1}^n \prod_{j=1}^{j_p-1} 
\Bigl(\frac{1-q^{j} t^{p}}
{1-q^j t}\Bigr).
\end{align}
Our theorem states that $\tilde{E}_b=\tilde{\n}_b E_b^\circ$ and
$\tilde{P}_b=\tilde{\n}_b P_b^\circ$ are $q,t$\~integral
for any $b\in P_+$ in this case.
\smallskip

{\sf Classical $J$\~polynomials.}
Let us provide their definition;
we emphasize that only the $A_n$\~version (for $sl_{n+1}$) of the 
$J$\~polynomials will be given and used (not that for $gl_{n+1}$):
\begin{align}\label{P-arms-legs}
J_\la\equal h_\la P_b \for \la=\la(b),\, 
h_\la=\prod_{\Box\in\la}
(1-q^{arm(\Box)}t^{leg(\Box)+1}).
\end{align}
This polynomial is $q,t$\~integral. 
\smallskip

Here $arm(\Box)$ is the {\em arm number\,}, which is
the number of boxes in the same row as $\Box$
strictly after it; $leg(\Box)$ is the {\em leg number\,},
which is the number of boxes in the column of $\Box$
strictly below it in the \Yboxdim7pt  
$\yng(2,1)\,$\~type presentation. Namely,   
$\la_1\!\ge \la_2,\ldots,\la_{n-1}\!\ge \la_n\,$ are the numbers of
boxes in the corresponding rows and the $i${\tiny th} row
is above the  $(i+1)${\tiny th}.
\smallskip

Equivalently:
\begin{align}\label{j-polynom}
&J_\la= t^{-(\rho,b)}
\prod_{p=1}^n \prod_{j=0}^{\la_{p^*\,}-1} 
\Bigl( 1\!-\!q^{j} t^{\,p+1}\Bigr) P_b^\circ, \ \,
p^*= n\!-\!p\!+\!1,\, b\in P_+.\ 
\end{align}
See, for instance, Theorem 2.1 from \cite{GoN}.
This makes $J_\la$ very similar to our 
$\tilde{P}_b$\~polynomials, but our multipliers are 
generally greater than those in $J_\la$; see below.
Note that the {\em arms and\, legs\,} do not 
appear here (in the 
approach via $P_b^\circ$).
Let us identify (\ref{j-polynom}) and (\ref{P-arms-legs})
for the sake of completeness and because we are doing the
$sl_{n+1}$\~case.
\smallskip

The extra binomial factors necessary to
obtain the classical $J_\la$ from $P_b^\circ=P_b/P_b(t^{-\rho})$
are as follows. For any box in the $(n-p+1)${\footnotesize th} 
row of $\la$, the $t$-power of the corresponding
binomial factor from $P_b^\circ$ to the classical
$J_\la$ must be $(n+1)-(n-p+1)+1=p+1$. 
The range of the powers of $q$ in this row is 
$\,\{0,1,\ldots,\la_{n-p+1}-1\}$. So this is what
formula (\ref{j-polynom}) does.
\smallskip

The product in (\ref{j-polynom}) before $P_b^\circ$ is
$\{J_\la\}_{ev}=J_\la(t^{-\rho}).$ It is important
that if one introduces one more variable $a=-t^{n+1}$,
then this evaluation depends only on $a$ and
the Young diagram $\la$. Namely, 
\begin{align}\label{stabeval}
J_\la(t^{-\rho})=t^{-(\rho,b)}
&\prod_{p=1}^n \prod_{j=0}^{\la_{p\,}-1} 
\Bigl( 1-q^{j} t^{n-p+2}\Bigr)\notag\\
=(a^2)^{-\la_1/4}\,t^{|\la|/2}\,
&\prod_{p=1}^n \prod_{j=0}^{\la_{p\,}-1} 
\Bigl( 1\,+\,q^{j}\, a\, t^{-p+1}\Bigr).
\end{align}
This is one of the main points of the approach
from \cite{GoN}.

\smallskip
{\sf Comparing $J_\la$ and $\tilde{P}_b$.}
First of all, Macdonald's {\em arms and\, legs\,} 
$a$({\small$\square$}),
$l$({\small$\square$}) disappear in the approach from 
(\ref{j-polynom}) due to our 
division of $P_b$ by its evaluation at $t^{-k\rho}$. Technically,
we use here 
the so-called {\em co-arms and co-legs\,}, but they are simply
the coordinates of the boxes in $\la$, which can be
directly spelled in terms of $b$ and the root system.

Our usage of $P_b^\circ$ as the starting point for the
definition of the $J$\~polynomials is of course not
just to make their definition less combinatorial. It
is very natural because of the construction of the 
DAHA-polynomials for iterated torus knots. Also,
the spherical normalization $P_b^\circ$ is the key for
the difference Mehta-Macdonald formulas from \cite{C5},
which will be considered below. 
\smallskip
 
Recall that in our construction of $\tilde{P}_b$,
we omit $j=0$ in (\ref{j-polynom}), i.e.
drop the binomials containing only $t$, and also divide
the binomials $\bigl(1- (q^j t)^{m_i+1}\bigr)\,$ by
$(1-q^j t)\,$ in our construction of the integral form. 

Ignoring  the binomials without $q$ in the formula for $J_\la$ 
obviously does not influence the $q,t$\~integrality; this 
makes the integral form smaller (but may impact the symmetries,
including the super-duality of the DAHA superpolynomials).
The second reduction is more interesting; it generally diminishes
$J$. For instance, our construction correctly
gives the minimal denominator of $P_n$ in the case of $A_1$, 
which is 
$(1-tq^{n-1})(1-tq^{n-2})\cdots (1-tq^{n-[n/2]})$. The 
$J$\~regularization requires the multiplication of $P_n$ by
the whole product $\prod_{j=1}^n (1-tq^{n-j})$.

With these reservations, the exact match between $\tilde{P}$
and $J$ is upon 
the relations $j_p=\la_{n-p+1}$  or $j_p=\la_1-\la_{p+1}$
(for all $p$). Then our 
$\tilde{P}_b=\tilde{\n}_b P_b^\circ$ essentially become
$J_\la$ for $\la=\la_b$; see (\ref{kapp-an1}).
This conditions are restrictive;
generally  the binomials from $\tilde{\n}_b$  have
powers of $q$ larger than those  
needed in $J_\la$ (with the same $t$\~powers). 
\smallskip

Such an increase of  our $\tilde{P}_b$ vs. $J_\la$
is mainly because our construction is uniform 
and does not use specific features of the DAHA intertwining 
operators in type $A_n$ (quite special). For instance,
the root system $A_m$ with $m$ equal to the 
number of rows of $\la=\la(b)$ is sufficient to find the 
denominator of $P_\la$.  
See Theorem 3.5.1 and Corollary 3.5.2 from \cite{HHL}, where 
the integral forms of $E$\~polynomials were effectively 
calculated using the intertwiners, and
Corollary 5.1.2 there concerning the $J$\~polynomials
(which are not always optimal as integral forms).  
\smallskip

We note that the factors $\tilde{\n}_b$ coincide for
$b$ and $b^\iota$ in our approach, since the integrality 
properties of the $P$\~polynomials are unchanged
under $\iota\,$. This 
cannot be true in the theory of $J$\~polynomials,
which is designed for undetermined ranks (and aimed at
the positivity). By the way, this readily results in a 
sharper integrality for a fixed $A_n$ and
sufficiently large diagrams $\la\,$
if the inequality $j<\la_{p^*}$ in 
(\ref{j-polynom}) for $p\ge 2$ is replaced by
$\min\{\la_{p^*},\la_1-\la_{p+1}\}$; 
recall that $\la_{p^*}=b_n+\ldots+b_{n+1-p}$ and 
$\la_1-\la_{p+1}=b_1+\ldots+b_p=\la_{p^*}(b^{\,\iota})$.
\smallskip

Because of the symmetry $b\mapsto b^{\,\iota}$ our integral
forms $\tilde{P}_b$ are generally larger than $J_\la$ for 
$\la=\la(b)$, but not always.
Let us provide examples. First of all,
$\tilde{P}_b=\tilde{\n}_b P_b^\circ$ for $b=m\om_n$ 
is actually smaller than $J_\la$ for  $\la=\la(m\om_n)$.
Indeed,  we do not need the binomials without $q$ and
replace those for non-rational affine exponents
by their certain factors; the $A_1$\~case was 
discussed above. However, $\tilde{\n}_{m\om_1}$
coincides with $\tilde{\n}_{m\om_n}$ in our approach, which
makes the former generally greater than 
$J_{\la^*}/P_{b^{\iota}}^\circ$\, 
for $b^\iota=m\om_1$ and $\la^*=\la(b^{\,\iota})$. This ratio 
is a product of $m$ binomials vs. approximately 
$n(m-1)$ factors in our 
$\tilde{\n}_{m\om_1}$, a significant increase.

Concerning the non-$A$ root systems,
the relations of our integral forms to formulas from  
\cite{RY,OS} can be expected. However our
$\tilde{\n}_b$ cannot be obtained (at least directly)   
from the formulas there; we use the spherical 
normalization. Both approaches
are based on the construction of the $E$\~polynomials 
via the DAHA intertwiners from
\cite{C1} (due to Knop and Sahi for $A_n$) and connections
are likely. We note that the $J$\~type integral forms
of $P_b$ can be defined for {\em classical\,} root systems,
which is similar to type $A$, but this will not be discussed here.  
\smallskip

To recapitulate,
there are two outcomes of our considerations:
 
a) $J_\la$ can be naturally given 
entirely via the root system $A_n$, if one 
introduce them in terms of the spherical normalization 
$P_b^\circ$ of $P_b$, which allows avoiding the usage of Young 
diagrams for this, including Macdonald's {\em arms and legs\,}, 
necessary if $P_b$ is used as the starting point;

b) certain counterparts of $J$\~polynomials can be defined 
for any root systems (reasonably close to them, but without 
the $n$\~stabilization), essentially matching the
$J_\la$ for $\la=\la(b), b\in P_+$ (in type $A_n$)  such that 
$\max_{i}
\{b_i+\ldots+b_{i-p}\}$ is reached at  
$i=n\,$ for any $\,p\ge 0$, where $i>p$.
\medskip

\setcounter{equation}{0}
\section{\sc Topological vertex}
\subsection{\bf Theta-functions}
The DAHA-vertex is closely related to
the difference Mehta-Macdonald formulas,
which require the {\em theta-functions\,}. They are
defined for the root system $R$ and will depend
in this work on the choice
of the character $u: \Pi=P/Q\to \C^*$. We denote this
group of characters by $\Pi'$ and the trivial character 
by $1'$. Let
\begin{align}
&\theta_u(X)\ \equal\ 
\sum_{b\in P} u(b) q^{(b,b)/2}X_b=
\ze_{u,1}(\theta) \for \theta\equal \theta_{1'}.
\label{gauser}
\end{align}

The characters $u$ play here the role
of the classical {\em theta- characteristics\,}; they were
important in \cite{ChF} (necessary for the
level-rank duality), though we used a somewhat
different setting there. Namely, they were
introduced using the partial summations $\th_{\varpi}$,
where the images of $b$ are taken from subsets  
$\varpi\subset\Pi$, 
instead of our using the characters $u$ here. 
The difference is that we obtain
another basis in the same space (of theta-functions); 
using $\theta_u$ is more convenient in the present work.

The immediate DAHA-motivation of this definition is
the following lemma, where the parameters $q,t_\nu$ are assumed
generic. 

\begin{lemma}\label{THTAU+}
The formal conjugation by $\theta_u^{-1}$ in a proper completion 
of $\HH$ or that 
for End$(\v)$, where $\v$ is the polynomial representation,
induces the following DAHA\~automorphism: 
$\tau_+^u=\tau_+\,\ze_{1,u}=\ze_{1,u}\tau_+$.
More exactly,
\begin{align}\label{thetatau+}
&(\theta_u)^{-1} H \theta_u= \tau_+^u(H) \for H\in \HH, u\in \Pi'. 
\end{align}
\end{lemma}

{\it Proof.} This is standard for $u=1'$ (see \cite{C101}),
let us outline the {\em enhancement\,} due to $u$.
We note that for practical calculations,
it is convenient to replace here $\theta$ by $q^{-x^2/2}$ for 
$X=q^x$; all formulas remain the same.

The theta-functions above are $W$\~invariant,
so it suffices to check that in a proper completion of $\v$,
$$ 
T_0(\theta_u)=\tau_+^u(T_0)(\theta_u) \hbox{\, and\, }
Y_r(\theta_u)=\tau_+^u(Y_r)(\theta_u) \for r\in O.
$$ 
Recall that  $\tau_+^u(T_0)=\tau_+(T_0)$, so
only the latter relation must be checked: 
\begin{align}
&Y_r(\theta_u)=\!\sum_{b\in P} u(b) q^{-(b,\om_r)}\,X_b q^{b^2/2}
\!=q^{-\om_r^2/2} u(\om_r)X_{\om_r}\theta_u\!=
\tau_+\bigl(\ze_{1,u}(Y_r)\bigr).\notag
\end{align}
\vskip -1.cm \sq

Combining this lemma with the action $(\dot{\tau}^u_-)^{-1}$ in
$\v$, we obtain that the semigroup generated by 
$(\tau_{\pm}^u)^{-1}$ for $u\in \Pi'$ acts in the space 
linearly generated by the products
$X_b\,\theta_{u_1}\cdots \theta_{u_l}$
for any {\em levels\,} $l$. We add here $q^{1/(2\mm)}$ to
the ring of constants $\Z_{q,t}$, as in Section \ref{sect:Aut}.
The products
$\theta_{u_1}\cdots \theta_{u_l}$ are considered 
as Laurent series
with the coefficients in $\Z[t_\nu^{1/2}][[q^{1/(2\mm)}]]$. 
Thus we switch from $\v$ to a certain subspace 
$\tilde{\v}_l$ of Laurent series. Considering 
$\oplus_l \tilde{\v}_l$, one can define there 
the (projective) action of $GL_2(\Z)$.

This can be used to obtain interesting DAHA\~invariants of 
iterated torus links,
but generally of more involved nature. For instance, we
obtain a connection of the {\em DAHA-vertices\,} for any
levels $l\,$ (see below) with the DAHA-Jones 
polynomials. 
\medskip

\subsection{\bf Mehta-Macdonald identities}
We follow Theorem 3.4.5 from \cite{C101}, enhancing
it by characters $u\in \Pi'$. It is worth mentioning
that the passage from the $E$\~polynomials to the
$P$\~polynomials is not always fully
clarified in \cite{C101}; we will
do this in this work.

We denote  the {\em constant term\,}
of Laurent series $f(X)$
by $\langle  f \rangle.$
Let $\mmu^\flat\equal \mmu/\langle \mmu \rangle$,
where
\begin{align}
&\langle\mmu\rangle\ =\ \prod_{\al \in R_+}
\prod_{i=1}^{\infty} \frac{ (1- q^{(\rho_k,\al)+i\nu_\al})^2
}{
(1-t_\al q^{(\rho_k,\al)+i\nu_\al})
(1-t_\al^{-1}q^{(\rho_k,\al)+i\nu_\al})
}.
\label{consterm}
\end{align}
Recall that $q_\al=q^{\nu_\al},\ q^{(z,\al)}=q_\al^{(z,\al^\vee)},
\ t_\al=q_\al^{k_\al}.$
This formula is equivalent to the celebrated
Macdonald constant term conjecture (see \cite{C101}).

The pairing $\,\lan f,g\ran\equal\lan f^\star\, g\,
\mmu^\flat  \ran\,$ for
$f,g\in \v\,$ plays a major role in the theory;
it is served by the anti-involution $\star$ of $\HH$
from (\ref{star-conj}).
Thus the operators $T_{\hw},X_a,Y_b$ are all unitary with
respect to it.

We will use Theorem 5.1 from \cite{C5} (the first formula) 
and Theorem 3.4.5 from \cite{C101}. There are
other Mehta-Macdonald-type formulas in these works, including 
the ones in terms of the
{\em Jackson integral\,}; we will not use them in this
work. Recall that 
\begin{align*}
&E_b^\circ(q^{c_{\#}})=E_c^\circ(q^{b_{\#}}),\,  
b_\# \equal b+w_b^{-1}(\rho_k),\, w_b(b_\#)=b_+ +\rho_k,\,
\end{align*}
and 
$P_{b_+}^\circ(q^{c_+ +\rho_k})=P_{c_+}^\circ(q^{b_+ +\rho_k})$, 
where the latter follows from the former
considered for $b_-=w_0(b_+)$; we use that $(b_-)_\#=
b_- -\rho_k=w_0(b_++\rho_k).$

See (\ref{ebdual}). We will also need (\ref{nontosymdom}):
\begin{align}\label{nontosymdomx}
&\Pi_R P_{b}^\circ=\mathscr{P}_+(E_{b}^\circ)\for
\mathscr{P}_+\!\equal\!\sum_{w\in W}
t_{\sht}^{l_{\sht}(w)/2}
t_{\lng}^{l_{\lng}(w)/2}T_w.
\end{align}  
Recall that $\Pi_R$ is the 
Poincar\`e polynomial (with $t_{\sht},t_{\lng}$).

\begin{theorem}\label{EPEP}
For $b,c\in P$ in the first following formula and
$\,b,c\in P_+$ in the second, 
\begin{align}
u(b+c)\langle E^\circ_b E^\circ_c \theta_u\mmu^\flat  
\rangle & \,=\,
q^{(b_\#,b_\#)/2+(c_\#,c_\#)/2 -(\rho_k,\rho_k)}
E^\circ_c(q^{b_\#})\langle \theta\mmu^\flat  \rangle,
\label{epep}
\\
u(b+c)\langle P^\circ_b P^\circ_c \theta_u\mmu^\flat  
\rangle & \,=\,
q^{b^2/2+c^2/2 +(b+c,\rho_k)}
P^\circ_c(q^{b+\rho_k})\langle \theta\mmu^\flat  \rangle.
\label{pepep}
\end{align}
Here the coefficients of
$\mmu^\flat$ are naturally expanded in terms of
positive powers of $q$ and the proportionality factor is 
a $q$\~generalization of the Mehta -Macdonald integral:
\begin{align}
&\langle \theta\mmu^\flat  \rangle\ =\
\prod_{\al\in R_+}\prod_{ j=1}^{\infty}\Bigl(\frac{
1-t_\al^{-1} q_\al^{(\rho_k,\al^\vee)+j}}{
1-      q_\al^{(\rho_k,\al^\vee)+j} }\Bigr).
\label{mehtamu}
\end{align}
\end{theorem}
{\it Proof.}
The first formula is obtained from that
in \cite{C5} using that 
\begin{align}\label{zetheE}
\theta_u=\ze_{u,1}(\theta),\, \ze_{u,1}(E_b)=u(b) E_b,
\hbox{\, and\, } \lan \ze_{u,1}(f(X))\ran=\lan f(X) \ran
\end{align}
for $u\in \Pi'$ and any Laurent series $f(X)$.
It results in the following important lemma.

\begin{lemma}\label{GAUS-EVAl}
For an arbitrary Laurent series $f_c(X)$ such that
$\ze_{u,1}(f_c)=u(c) f_c$ and a $W$\~invariant one  $g_c(X)$
with the same symmetry,  
\begin{align}\label{epepf}
u(b+c)\langle E^\circ_b\, \dot{\tau}_-(f_c)\, 
\theta_u\mmu^\flat  
\rangle\ & =\
q^{(b_\#,b_\#)/2-(\rho_k,\rho_k)/2}
f_c(q^{b_\#})\langle \theta\mmu^\flat  \rangle,\\
u(b +c)\langle P^\circ_{b}\, \dot{\tau}_-(g_c)\, 
\theta_u\mmu^\flat  
\rangle\ & =\
q^{(b,b)/2+(\rho_k,b)}
g_c(q^{b+\rho_k})\langle \theta\mmu^\flat  \rangle,
\label{epepfsym}
\end{align}
provided the convergence. Here we use the action of
$\dot{\tau}_-$ in the polynomial representation $\v$;
see (\ref{taumin-poly}). We assume 
that $b\in P$ in the
first formula and $b\in P_+$ in the second; note that $u(b)=
u(w(b))$ for any $w\in W, u\in \Pi'$.
\end{lemma}
{\it Proof.}
The polynomials $E_c,P_c$ (or those with
$\circ$) satisfy the condition 
from the lemma on $f_c,g_c$. Namely,
\begin{align}\label{zeonEP}
\ze_{u,1}(E_c)=u(c) E_c,\, \ 
\ze_{u,1}(P_c)=u(c) P_c.
\end{align} 

Formula (\ref{taumineb}) gives that (with or without $\circ$)\,:
\begin{align}\label{tauminebx}
\dot{\tau}_-(E_c)=q^{-\frac{(c_+\,,\,c_+)}{2}-
(c_+\,,\,\rho_k)}\,E_c
= q^{-(c_\#^2/2+\rho_k^2/2)}\,E_c
\for c\in P.
\end{align}
Then we substitute an arrive at (\ref{epepf}).
Formula (\ref{epepfsym}) is 
obtained by the symmetrization inside $\lan\ldots\ran$
in the left-hand side of (\ref{epepf}) for symmetric $g_c$ and
$b=b_-\in P_-$. Namely, let us use the relations
$$
\lan T_w(f)\,,\, g \ran=
\lan f \,,\,((T_w)^\star(g))\ran=
\lan f\,,\, \bigl(T_{w^{-1}})^{-1}(g)\bigr)\ran
$$
for any $f(X),g(X)$. Then (\ref{nontosymdomx}) above
considered for  $b=b_-=w_0(b_+)\in P_-$ gives that
\begin{align*}
&u(b+c)\,\Pi_R\,\langle
P^\circ_{b_+}\, \dot{\tau}_-(g_c)\,\theta_u\,
\mmu^\flat  \rangle\, =\,
u(b+c)\,\langle \mathscr{P}_+
\bigl(E^\circ_{b}\, \dot{\tau}_-(g_c)\,\theta_u\bigr)\,\mmu^\flat  
\rangle\\ 
=&u(b+c)\langle (\mathscr{P}_+^\star(1))^\star\,
E^\circ_{b}\, \dot{\tau}_-(g_c)\,\theta_u\,\mmu^\flat  
\rangle\!=\!u(b+c)\Pi_R\,\langle 
E^\circ_{b}\, \dot{\tau}_-(g_c)\,\theta_u\,\mmu^\flat  
\rangle \\
=&\,\Pi_R\,q^{(b_\#,b_\#)/2-(\rho_k,\rho_k)/2}\,
g_c(q^{b_\#})\,\langle \theta \mmu^\flat  \rangle
=\Pi_R\, q^{b_+^2/2+(b_+,\rho_k)}\,
g_c(q^{b_+ +\rho_k})\,\langle \theta\mmu^\flat  \rangle.
\end{align*}
Dividing by $\Pi_R$, we arrive at (\ref{epepfsym}). Formula
(\ref{pepep}) from the theorem is its particular case,
which concludes the justification of the theorem. \sq
\smallskip

\subsection{\bf Norm-formulas} 
Let us provide the norm formula
for the Macdonald $E$\~polynomials (see e.g., \cite{Ch4}),
the Main Theorem):

\begin{align}\label{normepols}
&\lan E_b,E_c\ran\!=\!
\de_{bc}\!\!\prod_{\{\al,j\}\in \La^+_b}
\Bigl(
\frac{
(1-q_\al^j t_\al^{-1} X_\al(q^{\rho_k}))
(1-q_\al^j t_\al X_\al(q^{\rho_k}))
}{
(1-q_\al^j X_\al(q^{\rho_k})) (1-q_\al^j X_\al(q^{\rho_k}))
}
\Bigr).
\end{align}
The corresponding formula for the symmetric polynomials
$P_b (b=b_+)$,
the Macdonald norm conjecture, reads as:
\begin{align}\label{normppols}
&\lan P_{b},P_{c}\ran\!=\!
\de_{bc}\!\!\!\!&\prod_{\al>0}\!\!\prod_{j=0}^{(\al^\vee\!,b)-1}\!
\frac{
(1-q_\al^{j+1}t_\al^{-1} X_\al(q^{\rho_k}))
(1-q_\al^j t_\al  X_\al(q^{\rho_k}))
}{
(1-q_\al^j X_\al(q^{\rho_k})) (1-q_\al^{j+1} X_\al(q^{\rho_k}))
}.
\end{align}

Note that $\lan P_{b},P_{b}\ran=\lan P_{b^{\iota}},
P_{b^{\iota}}\ran=
\lan P_{b},P_{b}\ran^*$, as well as for 
$\lan P^\circ_b,P^\circ_b \ran$ used below. 

We will need formula (3.4.1) from \cite{C101}
and its symmetric variant:
\begin{align}
 \lan E^\circ_b,E^\circ_c\ran\,=
\de_{bc}\frac{\mmu^{-1}(q^{b_{\#}})}
{\mmu^{-1}(q^{-\rho_k})}\,=\,
\de_{bc} \!\prod_{\{\al,j\}\in \La^+_b}
\Bigl(
\frac{
t_\al -q_\al^j X_\al(q^{\rho_k})}{
1-q_\al^jt_\al X_\al(q^{\rho_k})}
\Bigr)&, \label{normhats}\\
 \lan P^\circ_b,P^\circ_c\ran= 
\de_{bc}\prod_{\al>0}\Bigl(\frac{
1 - X_\al(q^{\rho_k})}{
1-\,X_\al(q^{\rho_k+b})}\Bigr)
\!\!\prod_{j=0}^{(\al^\vee\!,b)-1}
\Bigl(
\frac{
t_\al -q_\al^{j+1} X_\al(q^{\rho_k})}{
1-q_\al^j t_\al X_\al(q^{\rho_k})}
\Bigr)&\notag\\
=\de_{bc}\frac{\dde^{-1}(q^{b+\rho_k})}
{\dde^{-1}(q^{\rho_k})} \for
\dde(X)\!\equal\!\mmu(X) \prod_{\al>0}
\frac{1-X_\al^{-1}}{1\!-\!t_{\al}X_\al^{-1}}&.
\label{normphsym}
\end{align}

Respectively, one takes $b,c\in P$ in the first formula
and $b,c\in P_+$ in the second.
The latter formula follows from the former
for $b=b_+\in P_+$. To see this, one can
proceed as follows:
\begin{align*}
\lan E_b^\circ, P_c^\circ\ran\,&=\,
\lan\, E_b^\circ\, (P_c^\circ)^\star\,,\,
\frac{\mathscr{P}_+^\star}{\Pi^\star_R}(1)\,\ran\,
=\,\lan\, \frac{\mathscr{P}_+}{\Pi_R}
\Bigl(E_b^\circ\,(P_c^\circ)^\star\Bigr)\,,\,1\,\ran\\
&=\,\lan \,\frac{\mathscr{P}_+}{\Pi_R}
\bigl(E_b^\circ\bigr)\,,\,P_c^\circ\,\ran\,=\,
\lan P_b^\circ, P_c^\circ\ran. 
\end{align*}
Then we switch
from $\mmu$ to $\dde$ or expand $P_c$ in terms 
of the $E$\~polynomials and use their orthogonality
relations; see formula (3.3.15) in \cite{C101}, which
is direct from the technique of DAHA intertwiners.

\comment{
\begin{align*}
&\lan P_b^\circ, P_c^\circ\ran=\frac{1}{\Pi_R \Pi_R^\star}
\lan \mathscr{P}_+(E_b^\circ),\mathscr{P}_+(E_c^\circ)\ran\\
=
&\,\frac{1}{\Pi_R \Pi_R^\star}
\lan E_b^\circ, (\mathscr{P}_+^\star \mathscr{P}_+)E_c^\circ\ran=
\lan E_b^\circ, P_c^\circ\ran, 
\end{align*}
using $\mathscr{P}_+^\star=(\Pi_R^\star/\Pi_R)
\mathscr{P}_+\!\!=\!
t_{\sht}^{-l_{\sht}(w_0)}t_{\lng}^{-l_{\lng}(w_0)}\mathscr{P}_+$
and $\mathscr{P}_+^2\!=\!\Pi_R \mathscr{P}_+$. 
}

Formula (\ref{normphsym}) can be of course deduced 
from  (\ref{normppols}) and (\ref{Pbrho}):
\begin{align}\label{normphsymy}
 &\lan P^\circ_b,P^\circ_c\ran\!=\! 
\de_{bc}\prod_{\al>0}\!\!
\!\!\prod_{j=0}^{(\al^\vee\!,b)-1}
\!\!\Bigl(\frac{
1\!-\!q_\al^{j+1}t_\al^{-1}X_\al(q^{\rho_k})}{
1\!-\!q_\al^{j+1} X_\al(q^{\rho_k})}\Bigr)\!
\Bigl(
\frac{
1\!-\!q_\al^{j} X_\al(q^{\rho_k})}{
t_\al^{-1}\!-\!q_\al^j X_\al(q^{\rho_k})}
\Bigr), 
\end{align}
though the simplest justification of formula 
(\ref{normppols}) for $\lan P_b^\circ,P_b^\circ\ran$ is via
$\lan E^\circ_b,E^\circ_b\ran$, which is conceptually
connected with $\mmu^{-1}(q^{b_\#})$. 

\subsection{\bf The case of 
\texorpdfstring{{\mathversion{bold}$A_n$}}{A}}
Continuing for any root system,
we cancel coinciding factors in (\ref{Pbrho'}), but do not 
perform any further divisions: 
$$
P_b(q^{-\rho_k})=q^{-(b,\rho_k)}\prod 
\frac{(1-q_\al^j t_\al^i)}{(1-q_\al^{j\,{}^\flat} 
t_\al^{i^\flat})},
$$
Then $\lan P^\circ_b,P^\circ_c\ran = 
\de_{bc}q^{2(\rho_k,b)}/(f_b\, f'_b)$ for 
\begin{align}\label{normtoeval}
f_b=
\prod \frac{(1-q_\al^j t_\al^i)}
{(1-q_\al^{j\,{}^\flat} t_\al^{i^\flat})},
\ f'_b=
\prod \frac{(1-q_\al^{j+1} t_\al^{i-1})}
{(1-q_\al^{j\,{}^\flat+1} t_\al^{i^\flat-1})}.\end{align}

Let us now apply this formula to (\ref{stabeval})
in the case of $A_n$ using
\begin{align}\label{P-arms-legsx}
P_\la=J_\la/h_\la \for 
h_\la=\prod_{\Box\in\la}
(1-q^{arm(\Box)}t^{leg(\Box)+1}).
\end{align}
Recall that $arm(\Box)$ is
the number of boxes strictly after $\Box$ (in the same row)
and $leg(\Box)$ is the number of boxes 
strictly below it. See (\ref{P-arms-legsx}).

\begin{proposition}\label{MACD-NORMS-A}
For $\la=\la(b),\nu=\la(c),\, b,c\in P_+$\, in type
$A_n$, we obtain that 
\begin{align}\label{norm-arms-legs}
\lan P^\circ_\la,&\,P^\circ_\nu\ran = \de_{\la,\nu}\, 
\frac{(a^2)^{\la_1/2}}{t^{|\la|}}\,
\frac{h_\la h_\la'}{g_\la g_\la'} 
\,\for\, h'_\la=\prod_{\Box\in\la}
(1-q^{arm(\Box)+1}t^{leg(\Box)}),\notag\\
&g_\la=\prod_{p=1}^n \prod_{j=0}^{\la_{p\,}-1} 
\Bigl( 1\!+\!q^{j}\, a\, t^{-p+1}\Bigr),\ 
g_\la'=\prod_{p=1}^n \prod_{j=0}^{\la_{p\,}-1} 
\Bigl( 1\!+\!q^{j+1}\, a\, t^{-p}\Bigr).
\end{align}
\end{proposition}
{\it Proof.}
We combine (\ref{normtoeval}) and (\ref{stabeval}):
\begin{align}\label{stabevalx}
J_\la(t^{-\rho})\!=\!(a^2)^{-\frac{\la_1}{4}}\,
t^{\frac{|\la|}{2}}\,
&\prod_{p=1}^n \prod_{j=0}^{\la_{p\,}-1} 
\Bigl( 1+q^{j}\, a\, t^{-p+1}\Bigr) \for a=-t^{n+1}.
\end{align}
\vskip -1cm \sq

Note that the formulas for our norms of $P_b$ are
similar to (but do not coincide with) those for 
Macdonald's {\em stable\,} pairing; see e.g., 
\cite{GoN}, formula (6). Under his pairing, the denominator
will be $g_\la^2$:
\begin{align*}
\lan P_\la^\circ, P_\nu^\circ\ran_{\hbox{\tiny stab}}\,=\,
\de_{\la,\nu}\,\frac{g_\la'}{g_\la}\,\lan P^\circ_\la,\,
P^\circ_\nu\ran,\where g_\la'=g_\la (a\mapsto aq/t).
\end{align*}

{\sf $A_n$\~stabilization.}
Let us summarize the stabilization properties
of $P_b\,(b\in P_+)\,$ and the formulas above; this will later 
result in the $a$\~stabilization of the DAHA-vertex and 
DAHA-superpolynomials. 

First of all, $P_b$ can be naturally lifted to the root systems 
$gl_{n+1}$. Namely, to $\hat{P}_\la$ that is a polynomial in terms
of positive powers of the variables $X_{\vep_i}$ (with the
natural action of $W=\mathbf{S}_{n+1}$ by permutations of $i$)  
obtained from $P_b$ by the substitutions  
$$
P_b\ni X_c\mapsto X_{\vep_1}^{c_1+\ldots+c_n}
X_{\vep_2}^{c_2+\ldots+c_n}\cdots   
X_{\vep_n}^{c_n}\, \bigl(X_{\vep_1}\cdots 
X_{\vep_{n+1}}\bigr)^{c_{0}},
$$   
where $c=\sum_{i=1}^n c_i\om_i$ and
$c_0\in \Z_+$ is adjusted to ensure that the degrees of all 
monomials in $\hat{P}_\la$ are $|\la|$.  Note that $c_0=0$ for 
the leading monomial $X_b$ in $P_b$. In the opposite direction,
we obtain $P_\la$ by  
\begin{align}\label{hatPinP}
P_b=\bigl(X_{\vep_1}\cdots 
X_{\vep_{n+1}}\bigr)^{-\frac{|\la|}{n+1}}
\,\hat{P}_\la \for \la=\la(b),
\end{align}
where we express the right-hand side
in terms of $X_i=X_{\om_i}$ for $\,i=1,\ldots, n\,$ using
(\ref{omviavep}).
 
Then the passage from $gl_{n+m+1}$ to $gl_{n+1}$ 
(accordingly, from $A_{n+m}$ to $A_n$) is simply
by letting $X_{\ep_i}=0$ for $i>n+1$. The same 
stabilization holds for the $E$\~polynomials (one can
use the DAHA-intertwiners for them), but we will
not use this here. 

Similarly,  $\theta_u$ for $gl_{n+1}$ are as follows:
$$
\hat{\theta}_u=\sum_{\{m_i\}\in \Z^{n+1}} q^{\sum m_i^2/2}\, 
u(m_1+\ldots+ m_{n+1})\,X_{\vep_1}^{m_1}
\cdots X_{\vep_{n+1}}^{m_{n+1}},
$$
where $u$ is any homomorphism $u:\Z\to \C^*$. This definition
is stable with respect to $n$ in the same sense as for $P_\la$.
Then all formulas above can be transformed to $gl_{n+1}$,
The definition of $\mmu$ remains unchanged and switching to
$\hat{\om}_i=\sum_{j=1}^i\vep_j$ from $\om_i$ is natural.
The latter is actually necessary mainly to avoid 
the denominator $(n+1)$ in $b^2=
\sum_{i=1}^n \la_i^2-|\la|^2/(n+1)$ for $\la=\la(b)$;
see Section \ref{sec:A-n}. Then $b^2$ will become 
$\sum_{i=1}^n \la_i^2$,
i.e. will depend only on $\la$. 

We do not 
actually need DAHA theory for the root system $gl_{n+1}$
in this work. It is convenient to address the
$a$\~stabilization, but we can state and check almost all
stabilization claims within $sl_{n+1}$, which makes it
better compatible with our consideration of arbitrary
(reduced irreducible) root systems $R$. 
The passage to $gl_{n+1}$  may change the formulas only by 
some $q,t$\~monomial factors, which are generally ignored 
in the constructions below, as well as in \cite{CJ,CJJ,ChD}.  

\begin{proposition}\label{PROP:stab-values}
Given two Young diagrams $\la$ and $\mu$,\, 
the values $P_\la(q^{\mu+k\rho})$ are $a$\~stable, which means that 
there is a universal expression in terms 
of $q,t,a$ such that its value at $a=-t^{n+1}$ coincides with
 $P_b(q^{c+k\rho})$ for $\la=\la(b),\mu=\la(c)$,
$b,c\in P_+$ for $A_n$ with $n$ no smaller than the
number of rows in $\la$ and in $\mu$. Up to powers of $a^{1/2}$
and $t^{1/2}$,
they are rational function in  $q,t,a$. Also, 
$P_\la^\circ(q^{\mu+k\rho})$,
$\lan P_\la,P_\la\ran$
and $\lan P^\circ_\la,P^\circ_\la\ran$ are $a$\~stable
(in the same sense). 
\end{proposition}
{\it Proof.} The claim about the $a$\~stability of  
$P_b(q^{c+k\rho})$ 
is direct from the stability of $\hat{P}_\la$ for $\la=\la(b)$.
The correction factor from (\ref{hatPinP}) does not contribute
to these values. One needs only to know the 
$a$\~stability of $M_b(q^{c+k\rho})$ for the standard
symmetric monomials $M_b=\sum_{c\in W(b)} X_c$. 
Other claims follow from this and the 
explicit formulas provided above.

\vskip -0.5cm
\sq

\subsection{\bf High--level 3j--symbols}
We are back to an arbitrary root system $R$.
We will provide in this work only the symmetric
version of DAHA-vertex. The case of $E$\~polynomials is quite
parallel, but we do not see at the moment its 
applications (beyond the DAHA theory). For $l\in Z_+$,
called the {\em level\,} of theta-function, and for an
unordered set 
$\mathbf{u}=(u_i,1\le i\le l)\subset \Pi'=\hbox{Hom}(\Pi,\C^*)$,
we put:

\begin{align}\label{P-level}
&P_b^{\mathbf{u},\circ}\equal P_b^\circ \,\theta_{\mathbf{u}} \for 
\theta_{\mathbf{u}}=\prod_{i=1}^l\theta_{u_i},\,\ 
P^\circ_aP^\circ_b\,\theta_{\mathbf{u}}=\!\!\sum_{c\in P_+}
C^{c,\mathbf{u}}_{a,b} P^\circ_c,
\end{align}
where $a,b,c\in P_+$.
Here $\theta_{\mathbf{u}}$ are considered as Laurent series
in terms of $X_b$ with the coefficients that are formal series in
terms of positive powers of $q$. Analytically, the convergence
of their coefficients is granted for $|q|<1$. Here and below
we use that $\{P_c^{\mathbf{u},\circ}, c\in P_+\}$ is 
a basis in the space  $\v^W \theta_{\mathbf{u}}$ of $W$\~invariant 
vectors in $\v \theta_{\mathbf{u}}$, naturally considered over
the ring $\C[[\,q^{1/(2\mm)}\,]]\,[t_\nu^{1/2}]$.
\smallskip

Similarly, $C^{c,\mathbf{u}}_{\mathbf{b}}$
will be defined from the relation
\begin{align}\label{P-multi}
&P^{\mathbf{u},\circ}_{\mathbf{b}}\equal
P^\circ_{\mathbf{b}}\,\theta_{\mathbf{u}}=
P^\circ_{b_1}\ldots P^\circ_{b_k}
\,\theta_{\mathbf{u}}=\!\!\sum_{\mathbf{c}}
C^{c,\mathbf{u}}_{\mathbf{b}} 
P^\circ_{c}, 
\end{align}
where
$\mathbf{b}=(b_i, 1\le i\le k)\subset P_+\ni c.$

We put $\mathbf{u}=\emptyset$ for $l=0$ and, similarly,
use $\emptyset$ if the corresponding $\theta_u$ is
missing in $\theta_{\mathbf{u}}$. The coefficients 
$C_{a,b}^{c,\emptyset}$ are a generalization of 
the classical $3j$\~symbols (as $t_\nu=q_\nu$). 

The key example for us will be $l=1$; then
we write $P_b^{u,\circ}$ and $C^{c,u}_{a,b}$
for $u\in \Pi'$. Arbitrary levels can be reduced to $l=0,1$, 
as we will see. Actually $l=1$ is formally 
sufficient due to the relation
\begin{align}\label{lev-zero-C}
C_{a,b}^{c,u}=\sum_{d\in P_+} C_{a,b}^{d,\emptyset} C_{d,0}^{c,u}
\for a,b,c\in P_+,\, u=1'\in \Pi'.
\end{align}

To be more exact, $C_{a,b}^{d,\emptyset}$ can be expressed 
in terms of $C_{a,b}^{c,1'}$ if the inverse of the matrix
$\bigl(C_{d,0}^{c,u},\, c,d\in P_+\bigr)$ is known;
recall that $1'$ is the trivial character of $\Pi=P/Q$. 
The latter matrix essentially controls the DAHA Fourier transform,
so its inverse can be presented in a similar form; 
we will not discuss this here. Relation (\ref{lev-zero-C})
and similar relations for any levels readily follow 
from the fact that $\{P_c^{\mathbf{u},\circ}, c\in P_+\}$ is 
a basis in $\v^W \theta_{\mathbf{u}}$.


\begin{proposition}\label{LEV-ONE-C}
For any $b,c\in P_+$, using (\ref{pepep}) and 
(\ref{normphsym}):
\begin{align}\label{lev-1-two}
&C_{b,0}^{c,u}\!=\!\langle P^\circ_b 
\frac{P^\circ_{c^{\iota}}}{\lan P^\circ_c,P^\circ_c\ran} 
\theta_u\mmu^\flat  
\rangle\!=
\frac{q^{b^2/2+c^2/2 +(b\!+\!c,\rho_k)}}{u(b\!-\!c)
\lan P^\circ_c,P^\circ_c\ran}
P^\circ_b(q^{-c-\rho_k})\langle \theta\mmu^\flat  \rangle,
\end{align}
\begin{align}\label{normphsymx}
&\hbox{where\,\, } 
\lan P^\circ_c,P^\circ_c\ran=
\lan P^\circ_{c^{\iota}},P^\circ_{c^{\iota}}\ran=
\dde^{-1}(q^{c+\rho_k})/
\dde^{-1}(q^{\rho_k})\\
&=\prod_{\al>0}\Bigl(\frac{
1 - X_\al(q^{\rho_k})}{
1-\,X_\al(q^{\rho_k+c})}\Bigr)
\!\!\prod_{j=0}^{(\al^\vee\!,c)-1}
\Bigl(
\frac{
t_\al -q_\al^{j+1} X_\al(q^{\rho_k})}{
1-q_\al^j t_\al X_\al(q^{\rho_k})}
\Bigr).\notag
\end{align}
Similarly, for 
$\mathbf{b}=(b_i)\subset P_+\ni c$ and $P_{\mathbf{b}}^\circ\equal
\prod_i P_{b_i}^\circ$, 
\begin{align}\label{lev-1-abc}
&C_{\mathbf{b}}^{c,u}\!=\!\langle P^\circ_{\mathbf{b}} 
P^\circ_{c^{\iota}}
\theta_u\mmu^\flat\rangle/\lan P^\circ_c,P^\circ_c\ran\\
&=\frac{\langle P^\circ_{\mathbf{b}}
P^\circ_{c^{\iota}}
\theta\mmu^\flat\rangle}{u(\Si_i b_i-c)\,
\lan P^\circ_c,P^\circ_c\ran}
=\frac{\dot{\tau}_-^{-1}(P^\circ_{\mathbf{b}} P^\circ_{c^{\iota}})
(q^{-\rho_k})
\langle\theta\mmu^\flat\rangle}
{u(\sum_i b_i-c)\,\lan P^\circ_c,P^\circ_c\ran}.\notag
\end{align}
Therefore $C_{\mathbf{b}}^{c,u}$
belong to $\C(q^{1/(2\mm)},t_\nu^{1/2})
\langle \theta\mmu^\flat  \rangle$; the
denominators are the products of binomials 
from (\ref{binaffexp}), i.e. in terms of
$\bigl(1-q^\bullet t_{\sht}^{\bullet}t_{\lng}^{\bullet}\bigr)$
with non-negative powers (strictly positive for $q$). 
\end{proposition}
\vskip -0.7cm
\sq

The following {\em associativity theorem\,} is direct
from the definition of $C_{\mathbf{b}}^{c,\mathbf{u}}$;
we use that $P_c^\circ$ are linearly independent.
 
\begin{theorem} \label{ASSOC-C}
(i) Let us represent the sets $\mathbf{b},\mathbf{u}$ as unions
$$\tilde{\mathbf{b}}=
\bigl((\mathbf{b}^{1}),\ldots,(\mathbf{b}^m)\bigr),\
\tilde{\mathbf{u}}=
\bigl((\mathbf{u}^{1}),\ldots,(\mathbf{u}^m)\bigr),
$$
where the decompositions
$\tilde{\mathbf{b}}$ and $\tilde{\mathbf{u}}$
are obtained from $\mathbf{b}$ and $\mathbf{u}$ by adding
a proper number of $\,0${\small s} for the first and 
$\,\emptyset${\small s} for 
the second when necessary (at any places) to make them of the 
same length $\,m$; the size of $\mathbf{b}^{i}$ can be
different from that of $\mathbf{u}^{i}$.
Then
\begin{align}\label{ass-c-gen}
&C_{\mathbf{b}}^{c,\mathbf{u}}=
\Xi_{\tilde{\mathbf{b}}}^{c,\tilde{\mathbf{u}}}\equal
\\
&\sum_{c_1,c_2,\ldots,c_{m-1}\in P_+}
C^{c_1,\mathbf{u}^1}_{\mathbf{b}^1} 
C^{c_2,\mathbf{u}^2}_{(c_1,\mathbf{b}^2)} 
C^{c_3,\mathbf{u}^3}_{(c_2,\mathbf{b}^3)}
\cdots 
C^{c,\mathbf{u}^m}_{(c_{m-1},\mathbf{b}^m)}.\notag
\end{align}
For instance, this gives that 
$\,\Xi_{\tilde{\mathbf{b}}}^{c,\tilde{\mathbf{u}}}$
does not depend on the choice of the decompositions
$\tilde{\mathbf{b}}$ and $\tilde{\mathbf{u}}$ and
depends only on the sets $\mathbf{b},\mathbf{u}$ and $\,c$.

(ii) In particular, for  $\,\tilde{\mathbf{u}}=
\bigl((u_1),\ldots,(u_m)\bigr),\,
\tilde{b}=\bigl((b),(0),\ldots,(0)\bigr):$  
\begin{align}\label{ass-s-two}
&\Xi_{\tilde{b}}^{c,\tilde{\mathbf{u}}}=
\!\!\sum_{c_1,c_2,\ldots,c_{m-1}\in P_+}
C^{c_1,u_1}_{b,0} 
C^{c_2,u_2}_{c_1,0} 
C^{c_3,u_3}_{c_2,0}
\cdots 
C^{c,u_m}_{c_{m-1},0},
\end{align}
where $C^{b,u}_{a,0}$ are given by formula (\ref{lev-1-two}). 

More generally, let  
$\tilde{\mathbf{b}}=\bigl((b_0,b_1),(b_2),\ldots,(b_m)\bigr).$
Then  
\begin{align}\label{ass-s-three}
&\Xi_{\tilde{\mathbf{b}}}^{c,\tilde{\mathbf{u}}}=
\!\!\sum_{c_1,c_2,\ldots,c_{m-1}\in P_+}
C^{c_1,u_1}_{b_0,b_1} 
C^{c_2,u_2}_{c_1,b_2} 
C^{c_3,u_3}_{c_2,b_3}
\cdots
C^{c,u_m}_{c_{m-1},b_m}
\end{align}
depends only on the (unordered) sets $(b_i),(u_i)$ and $c$,
where $C^{c,u}_{a,b}$ are given by (\ref{lev-1-abc}).
For instance, the $3$-vertex associativity conditions hold:
\begin{align}\label{ass-s-4}
&\!\sum_{c_1\in P_+}
C^{c_1,u_1}_{b_0,b_1} 
C^{c,u_2}_{c_1,b_2}=
\!\sum_{c_1\in P_+}
C^{c_1,u_1}_{b_1,b_2} 
C^{c,u_2}_{c_1,b_0}=
\!\sum_{c_1\in P_+}
C^{c_1,u_1}_{b_0,b_2} 
C^{c,u_2}_{c_1,b_1}.
\end{align}
\end{theorem}
\vskip -1cm \sq

Note that the order of $(u_1,u_2)$ influences the
summations in (\ref{ass-s-4}), though the output does
not depend on it. Indeed,
$$
C^{c_1,u_1}_{b_0,b_1} 
C^{c,u_2}_{c_1,b_2}=C^{c_1,1'}_{b_0,b_1} 
C^{c,1'}_{c_1,b_2}/\bigl(u_1(b_0+b_1-c_1)u_2(c_1+b_2-c)\bigr),
$$
and the denominator here depends on $c_1$ unless $u_1=u_2$.
Permuting $\mathbf{u}^i$ in 
$\tilde{\mathbf{u}}$ from (\ref{ass-c-gen}) leads to quite 
non-trivial identities. This includes the $2$\~vertices from
(\ref{ass-s-two}); the corresponding identities 
considered under $t_\nu=0$ played 
an important role in \cite{ChF},  
with a direct link to the {\em level-rank duality\,} in the 
case of $A_n$. 

In a greater detail, formula (\ref{ass-s-two}) generalizes
Rogers-Ramanujan-type sums in \cite{ChF}, which were obtained
for $t_\nu=0$, $b=0$ and for minuscule weights $c$.
Then the corresponding sums will be {\em $q$\~modular functions\,}
under a certain normalization. The 
quantities $P_b(q^{c+\rho_k})$ disappear 
from (\ref{lev-1-two}) under $t_{\sht}=0=t_{\lng}$, which 
significantly simplifies the theory.
The modular properties of general
$\Xi_{\tilde{\mathbf{b}}}^{c,\tilde{\mathbf{u}}}$
are not known, though the case $t_\nu=q_\nu$ is
actually similar to $t_\nu=0$.

Relations (\ref{ass-s-4}) are basic for TQFT,
where they are (mainly) considered at roots of unity,
which makes the sums involved finite. They result,
for instance, in the WRT invariants of links
if the Quantum Group generalization of the $3j$-symbols 
is taken as the starting point. Our approach is compatible 
with the passage to roots of unity, but this will not be 
discussed in this work. 

Another important property of the {\em topological vertex\,}
needed in TQFT is its $\mathbf{S}_3$\~symmetry, though this is not 
always assumed (especially in the refined theories). 
See e.g., \cite{AKMV,GIKV,AFS}.
We obviously have that 
$C_{a,b}^{c,u}=C_{b,a}^{c,u}$. The invariance
with respect to $(a,b,c^\iota)\mapsto (a,c,b^\iota) $  
holds upon  multiplication
of our $C_{a,b}^{c,u}$ by 
$\lan P_c^\circ, P_c^\circ\ran$,
which readily follows from (\ref{lev-1-abc}). Thus we
essentially have the $\mathbf{S}_3$\~symmetry for the DAHA-vertex.

\medskip

\subsection{\bf The coinvariant approach}
Let us provide an alternative approach to level-one
$C_{a,b}^{c,u}$, which clarifies their $\mathbf{S}_3$\~invariance
and connects them with the DAHA polynomials associated 
with links (with the Hopf links, to be exact).

The key for such an approach is Lemma \ref{GAUS-EVAl}.
We will begin with the case of $E$\~polynomials with
two sets $\mathbf{b}=(b_i,\,1\le i\le k),
\mathbf{c}=(c_j,\,1\le j\le m)\in P$
and (a single) $u\in \Pi'$. The case of arbitrary sets 
$\mathbf{u}$ (i.e. arbitrary levels) requires general theory 
of links, which is beyond this section. Following the notation 
from (\ref{P-multi}), let $E^\circ_{\mathbf{b}}=
E^\circ_{b_1}\cdots E^\circ_{b_k}\,$
(and for $\mathbf{c})$. The same notation will be used
without $\circ$ and (later) 
for $P$\~polynomials. Recall from (\ref{zetauvrel1}):
$$
\tau_+^u=\tau_+\ze_{1,u}=\ze_{1,u}\tau_+,\,
\tau_-^u=\ze_{u,1}\tau_-=\ze_{u,1}\tau_-,\ 
\si^{u,v}=\ze_{uv^{-1},uv}\si.
$$
For any $g,h\in \C[X_a]$ (in particular, for 
$g=E_{\mathbf{b}}^\circ, h=E_{\mathbf{c}}^\circ\,$) and for the
{\em coinvariant\,} $\,\{.\}_{ev}$
from (\ref{evfunct}), we set:
\begin{align}\label{vert-coinv}
\e_{g}^{h,u}\equal&\bigr\{(\tau_+^{1/u})^{-1}\bigl(h(Y^{-1})\bigr)
(\tau_-^u)^{-1}\bigl(g(X)\bigr)\bigr\}_{ev}\!\!=
\bigr\{(\tau_-^u)^{-1}\bigl(h(X)g(X)\bigr)\bigr\}_{ev},\notag\\
&\e_{\mathbf{b}}^{\mathbf{c},u}\!=\!
\frac{\bigl\{\tau_+^{-1}\!
\bigl(E_{\mathbf{c}}^\circ(\!Y^{\!-1}\!)\bigr)
\tau_-^{-1}\!\bigl(E_{\mathbf{b}}^\circ(\!X\!)\bigr)\bigr\}_{ev}}
{u(\Si_ib_i+\Si_jc_j)}\!=\!
\frac{\bigl\{\tau_-^{-1}\!
\bigl(E_{\mathbf{c}}^\circ(\!X\!) 
E_{\mathbf{b}}^\circ(\!X\!) \bigr)\bigr\}_{ev}}
{u(\Si_ib_i+\Si_jc_j)}.
\end{align}
The action of $\tau_{\pm}$ here and below is that in $\HH$.
We use that 
$\si^{u,u}=\tau_+^u(\tau_-^u)^{-1}\tau_+^u=\ze_{1,u^{2\,}}\si$,
which gives that 
$$
(\tau_+^{1/u})^{-1}\bigl(h(Y^{-1})\bigr)\!=\!
(\tau_+^{1/u})^{-1}\si\bigl(h(X)\bigr)\!=\!
(\tau_+^{u})^{-1}\si^{u,u}\bigl(h\bigr)\!=\!
(\tau_-^{u})^{-1}\bigl(h\bigr).
$$
Then the action of $\ze_{u,1}$ on the $E$\~polynomials
from (\ref{zeonEP}) is applied. Note the $g\leftrightarrow h$\~
symmetry of (\ref{vert-coinv}) and that 
$\e_{\mathbf{b}}^{\mathbf{c},u}$ depends only on the
(unordered) set $\mathbf{b}\cup\mathbf{c}$.

It is worth mentioning that for any $g(X)$, 
\begin{align*}
&\{ \phi\Bigl(\tau_-^{-1}\bigl(g(X)\bigr)\Bigr)\}_{ev}=
\{ \tau_+^{-1}\phi\Bigl(\bigl(g(X)\bigr)\Bigr)\}_{ev}\\
=&\,\{ \tau_+^{-1}\si\Bigl(\bigl(g(X)\bigr)\Bigr)\}_{ev}=
\{ \tau_-^{-1}\tau_+\Bigl(\bigl(g(X)\bigr)\Bigr)\}_{ev}
=\{ \tau_-^{-1}\Bigl(\bigl(g(X)\bigr)\Bigr)\}_{ev}.
\end{align*}
This clarifies our using the inverses $\tau_{\pm}^{-1}$ for
$\e_{\mathbf{b}}^{\mathbf{c},u}$; if they are replaced by 
$\tau_{\pm}$, then  $Y^{-1}$ must be replaced there by 
$\si^{-1}(X)$. 
\smallskip

The same formulas
hold for the $P$\~polynomials; we set
$\p_{g}^{h,u}=\e_{g}^{h,u}$ for symmetric $g,h$.
Then for $b=(b_i)\subset P_+\supset (c_j)$,
\begin{align}\label{vert-coinv-p}
&\p_{\mathbf{b}}^{\mathbf{c},u}\!\equal\!
\frac{\bigl\{\tau_+^{-1}\!
\bigl(P_{\mathbf{c}}^\circ(\!Y^{\!-1}\!)\bigr)
\tau_-^{-1}\!\bigl(P_{\mathbf{b}}^\circ(\!X\!)\bigr)\bigr\}_{ev}}
{u(\Si_ib_i+\Si_jc_j)}\!=\!
\frac{\bigl\{\tau_-^{-1}\!
\bigl(P_{\mathbf{c}}^\circ(\!X\!) 
P_{\mathbf{b}}^\circ(\!X\!) \bigr)\bigr\}_{ev}}
{u(\Si_ib_i+\Si_jc_j)}.
\end{align}
\smallskip

There is a simplification in the case of a single 
$c\in P$ (or $b$) due to the following general identity based
on (\ref{tauminebx}). Recall that 
$b_\#^2/2-\rho_k^2/2=b_+^2/2+(b_+,\rho_k)$; see
(\ref{ebdual}), (\ref{tauminebx}) and Lemma \ref{GAUS-EVAl}.
For any $m\in \Z, H\in \HH$
and $g=g(X)$, one has:
\begin{align}\label{tau-m-c}
&\bigl\{H\cdot (\tau_-^u)^m\bigl(\!E_c^\circ(X)\bigr)
\bigr\}_{ev}\!=
u(mc)q^{-m(\frac{c_+^2}{2}+(c_+,\rho_k))}
\bigl\{H\cdot E_c^\circ(X)\bigr\}_{ev}, \\
&\bigl\{(\tau_-^u)^m\!\bigl(g\,E_c^\circ\bigr)\bigr\}_{ev}=
u(mc)\,q^{-m(\frac{c_+^2}{2}+(c_+,\rho_k))}
\bigl\{(\tau_-^u)^m\bigl(g\bigr)E_c^\circ(X)\bigr\}_{ev},\notag
\end{align}
where the second follows from the first
for $H=(\tau_-^u)^m(g(X))$.
\smallskip

\begin{proposition}\label{EV-VERT}
For $h=E_c^\circ (c\in P)$ and any $g(X)$, 
\begin{align}\label{e-single-c}
\e_{g}^{h,u}=
u(-c)q^{c_+^2/2+(c_+,\rho_k)}
\bigl\{E_c^\circ(Y^{-1})(\tau_-^u)^{-1}
\bigl(g(X)\bigr)\bigr\}_{ev}&,\\
\e_{\mathbf{b}}^{c,u}=
\frac{q^{c_+^2/2+(c_+,\rho_k)}}{u(c+\!b_1+\ldots+\!b_k)}
\bigl\{E_c^\circ(Y^{-1})\tau_-^{-1}
\bigl(E^\circ_{\mathbf{b}}(X)\bigr)\bigr\}_{ev}&\notag\\
=\frac{q^{c_+^2/2+(c_+,\rho_k)}
\Bigl(\dot{\tau}_-^{-1}\bigl(E^\circ_{\mathbf{b}}(X)\bigr)\Bigr)
(q^{c_\#})}{u(c+b_1+\ldots+b_k)}
=\frac{
\Bigl(\dot{\tau}_-^{-1}\bigl(E_c^\circ 
E^\circ_{\mathbf{b}}\bigr)\Bigr)
(q^{-\rho_k})}{u(c+b_1+\ldots+b_k)},&\notag\\
\e_{\mathbf{b}}^{c,u}\,=\,\frac{
\langle E^\circ_c\,E^\circ_{\mathbf{b}}\, 
\theta_u\mmu^\flat \rangle} 
{\langle \theta\mmu^\flat \rangle} \for 
\mathbf{b}=(b_1,\ldots,b_k)\subset P \ni c&,\label{coinvianorm}
\end{align}
where $\dot{\tau}^{-1}_-:\,g(X)\mapsto
\tau_-^{-1}\bigl(g(X)\bigr)(1)$ is the action of $\tau_-^{-1}$
in the polynomial representation $\v\ni g(X)$. In the symmetric
case, 
\begin{align}\label{e-single-c-p}
\p_{\mathbf{b}}^{c,u}=
\frac{q^{c^2/2+(c,\rho_k)}}{u(c+\!b_1+\ldots+\!b_k)}
\bigl\{P_c^\circ(Y^{-1})\tau_-^{-1}
\bigl(P_{\mathbf{b}}^\circ(X)\bigr)\bigr\}_{ev}&\\
=\frac{q^{c^2/2+(c,\rho_k)}
\Bigl(\dot{\tau}_-^{-1}\bigl(P^\circ_{\mathbf{b}}(X)\bigr)\Bigr)
(q^{c+\rho_k})}{u(c+b_1+\ldots+b_k)}
=\frac{
\Bigl(\dot{\tau}_-^{-1}\bigl(P_c^\circ 
P^\circ_{\mathbf{b}}\bigr)\Bigr)
(q^{-\rho_k})}{u(c+b_1+\ldots+b_k)}&\notag\\
=\frac{
\langle P^\circ_c\,P^\circ_{\mathbf{b}}\, 
\theta_u\mmu^\flat \rangle} 
{\langle \theta\mmu^\flat \rangle}=
C_{\mathbf{b}}^{\iota(c),u}
\frac{\langle P^\circ_c,\,P^\circ_c\rangle}
{\langle \theta\mmu^\flat \rangle}
\for \mathbf{b}=(b_i)\subset P_+ \ni c\,&,\notag
\end{align}
where $\iota(c)=c^\iota$, $C_{\mathbf{b}}^{c,u}$ is defined in
(\ref{P-multi}); see also formula (\ref{lev-1-abc}).
\end{proposition}

\comment{
\begin{align}\label{epepfy}
u(b+c)\langle E^\circ_b\, \tau_-(f_c)\, \theta_u\mmu^\flat  
\rangle\ & =\
q^{(b_\#,b_\#)/2-(\rho_k,\rho_k)/2}
f_c(q^{b_\#})\langle \theta\mmu^\flat  \rangle,\\
u(b +c)\langle P^\circ_{b}\, \dot{\tau}_-(g_c)\, 
\theta_u\mmu^\flat  
\rangle\ & =\ 
\end{align}}
{\it Proof.}
Applying $\phi$ inside the coinvariant (twice, back and
forth), we use 
(\ref{tau-m-c}) for $m=-1$ and 
that
$\phi \tau_-^u=\tau_+^{1/u}\phi$\,:
\begin{align*}
&\bigl\{E_c^\circ(Y^{-1})(\tau^u_-)^{-1}
\bigl(g(X)\bigr)\bigr\}_{ev}=
\bigl\{(\tau_+^{1/u})^{-1}\bigl(g(Y^{-1})\bigr)
E_c^\circ(X)\bigr\}_{ev}\\
&=(u(-c)q^{\frac{c_+^2}{2}+(c_+,\rho_k)})^{-1}
\bigl\{(\tau_+^{1/u})^{-1}\bigl(g(Y^{-1})\bigr)
(\tau_-^u)^{-1}\bigl(E_c^\circ(X)\bigr)\bigr\}_{ev}\\
&=u(c)\,q^{-(\frac{c_+^2}{2}+(c_+,\rho_k))}\,
\bigl\{(\tau_+^{1/u})^{-1}\bigl(E_c^\circ(Y^{-1})\bigr)\,
(\tau_-^{u})^{-1}\bigl(g(X)\bigr)
\bigr\}_{ev}.
\end{align*}
Then we use formula (\ref{epepf}) from 
Lemma \ref{GAUS-EVAl} for $b\!=\!0$ and
$f\!=\!\dot{\tau}_-^{-1}(E^\circ_c E^\circ_{\mathbf{b}})$,
where the $\dot{\tau}_-^{\pm 1}$ (here and in the next
formula)
is the action of $\tau_-^{\pm 1}$ in $\v$, which is $g(X)\mapsto
\tau_-^{\pm 1}\bigl(g(X)\bigr)(1)$.
Namely, one has\,:
\begin{align}
&u(c+\Si_i b_i)\,
\langle \dot{\tau}_-(f)\, \theta_u\mmu^\flat  
\rangle\ =
f(q^{-\rho_k})\langle \theta\mmu^\flat  \rangle\notag\\
=&\ u(c+\Si_i b_i)\,
\langle E^\circ_c\, E^\circ_{\mathbf{b}}\, \theta_u\mmu^\flat  
\rangle\ =\
\bigl(\dot{\tau}_-^{-1}
(E^\circ_{\mathbf{b}})\bigr)(q^{-\rho_k})\langle 
\theta\mmu^\flat  
\rangle.\notag
\end{align}
\vskip -0.7cm
\sq

This proposition readily results in the following important
general property of the evaluation pairing, which is 
necessary to connect the DAHA approach with 
the splice diagrams of the corresponding link.

\begin{theorem}\label{EVALTAUM}
For arbitrary $f=f(X),g=g(X)\in \v$ and the evaluation pairing
$\{f,g\}_{ev}=\{f(Y^{-1})\bigl(g(X)\bigr)\}_{ev}$,
\begin{align}\label{eval-tau}
\bigl\{\dot{\tau}_-^{-1}(f),\dot{\tau}_-^{-1}(g)\bigr\}_{ev}=
\bigl\{\dot{\tau}_-^{-1}(fg)\bigr\}_{ev}=
\bigl\{\tau_-^{-1}\bigl(f(X)g(X)\bigr)\bigr\}_{ev}.
\end{align}
Combining this relation with Lemma \ref{GAUS-EVAl} in the
notations there and under the
assumptions $\ze_{u,1}(f)=u(b) f$, 
$\ze_{u,1}(g)=u(c) g$ for certain $b,c\in P$ and all $u\in \Pi'$,
one has for $f,g$ as above:
\begin{align}\label{epepfgen}
u(b+c)\,\lan \dot{\tau}_-(f)\, \dot{\tau}_-(g)\, 
\theta_u\mmu^\flat  
\ran\  =\
\{f,g\}_{ev}\,\lan \theta\mmu^\flat  \ran.
\end{align}
Also, 
$\{\bigl(\dot{\tau}_-(f)\bigr)(Y)\,
\dot{\tau}_-(g)\}_{ev}
=\bigl\{\dot{\tau}_-(fg)\bigr\}_{ev}$\,
for $f\in \v^W$,
where $\v^W$ is the subsubspace of $W$\~invariants. \sq
\end{theorem}

{\it Proof.} 
We extend $\{f,g\}_{ev}$ to $A,B\in \HH$
to  $\{A,B\}_{ev}\equal\{\vph(A)B\}_{ev}$. Then
$\{A,B\}_{ev}\!=\!\{B,A\}_{ev}$\, and for $H(1)\in \v\subset\HH
\ni H:$
\begin{align}
\{A,B&\}_{ev}\!=\!
\{\vph(A)B(1)\}_{ev}\!=\!\{\vph\bigl(B(1)\bigr)A\}_{ev}\!=\!
\{A(1),B(1)\}_{ev}\,,\notag\\
\bigl\{\dot{\tau}_-^{-1}(f)&\,,\dot{\tau}_-^{-1}(g)\bigr\}_{ev}=
\bigl\{\dot{\tau}_-^{-1}(f),\tau_-^{-1}(g)\bigr\}_{ev}=
\bigl\{\tau_-^{-1}(g),\dot{\tau}_-^{-1}(f)\bigr\}_{ev}
\notag\\
=\!
\bigl\{\tau_-^{-1}(g)&\,,\tau_-^{-1}(f)\bigr\}_{ev}
\!\!=\!\bigl\{\tau_+^{-1}\bigl(\si(g)\bigr),
\tau_-^{-1}(f)\bigr\}_{ev}
\!\!=\!\bigl\{\tau_+^{-1}\bigl(\vph(g)\bigr),
\tau_-^{-1}(f)\bigr\}_{ev}
\notag\\
=&
\bigl\{\vph\bigl(\tau_-^{-1}(g)\bigr),\tau_-^{-1}(f)\bigr\}_{ev}=
\bigl\{\tau_-^{-1}\bigl(gf\bigr)\bigr\}_{ev}=
\bigl\{\dot{\tau}_-^{-1}\bigl(fg\bigr)\bigr\}_{ev}.
\label{eval-taux}
\end{align}

This gives (\ref{eval-tau}). 
Then use Lemma \ref{GAUS-EVAl} and 
$\,\si^2=\iota$  in $\v^W$. \sq

It is worth mentioning that the relation
\begin{align}\label{tilde-C}
\tilde{C}_{\mathbf{b}}^{\,c,u}\equal
\frac{C_{\mathbf{b}}^{\,c,u}}
{\langle \theta\mmu^\flat \rangle}=
\frac{\p_{\mathbf{b}}^{\,\iota(c),u}}
{\langle P^\circ_{c^{\iota}},\,P^\circ_{c^{\iota}}\rangle}
=\frac{\p_{\mathbf{b}}^{\,\iota(c),u}}
{\langle P^\circ_c,\,P^\circ_c\rangle}
\end{align}
is actually entirely conceptual, as well as its
nonsymmetric counterpart. However, we prefer to obtain
it via Lemma \ref{GAUS-EVAl}, which is of independent
interest. The right-hand side here is convenient to
analyze the action of $\star$ and $\eta$. Using 
(\ref{evsym}) and (\ref{eta-tau}):
\begin{align}\label{vert-coinv-px}
&(\p_{\mathbf{b}}^{\mathbf{c},u})^\star
=
\frac{\bigl\{\tau_-\!
\bigl(P_{\mathbf{c}^\iota}^\circ(\!X\!) 
P_{\mathbf{b}^\iota}^\circ(\!X\!) \bigr)\bigr\}_{ev}}
{u(\Si_ib_i^\iota+\Si_jc_j^\iota)},\
(\tilde{C}_{\mathbf{b}}^{\,c,u})^\star=
\tilde{C}_{\mathbf{b}^\iota}^{\,\iota(c),u},
\end{align}
where we apply $\iota=-w_0$ to $\mathbf{b}$ componentwise. 

We obviously used here the rationality of the right-hand side
of (\ref{tilde-C}); otherwise $q^\star=q^{-1}$ must be
addressed. Note that one can define 
$$
\tilde{C}_{\mathbf{b}}^{\,c,\mathbf{u}}
=C_{\mathbf{b}}^{\,c,\mathbf{u}}/
\langle \theta\mmu^\flat \rangle^l
$$
for any level $l\ge 0$ and $\mathbf{u}=(u_1,\ldots,u_l)$.
Then the relations from Theorem \ref{ASSOC-C} for
the $C$\~coefficients will hold
for the $\tilde{C}$\~coefficients; we mean those
based on the independence of
$\Xi_{\tilde{\mathbf{b}}}^{c,\tilde{\mathbf{u}}}$ on
the choice of the decompositions $\tilde{\mathbf{b}}$, 
$\tilde{\mathbf{u}}$.
However, such coefficients are generally infinite series 
in terms of (non-negative fractional powers of)
$q$ and applying $\star\,$, which includes 
$q\mapsto q^{-1}$, is more involved (though doable thanks 
to the $q,t$\~setting). 



\setcounter{equation}{0}
\section{\sc DAHA-Jones theory}
\subsection{\bf Iterated torus knots}\label{sec:ITER-KNOTS}
We will begin with torus knots and iterated torus knots; 
the iterated {\em links\,} will be considered next.

The torus knots $T(\rr,\ss)$ are defined for any
integers (including $0$ and negative ones)
assuming that \,gcd$(\rr,\ss)=1$.
One has the symmetries $\,T(\rr,\ss)=T(\ss,\rr)=T(-\rr,-\ss)$,
where we use ``$=$" for the ambient isotopy equivalence.
Also $\,T(\rr,\ss)=\unknot\,$ if $\,|\rr|\le 1$ or $|\ss|\le 1$
for the {\em unknot\,}, denoted by $\,\unknot\,$.
See e.g., \cite{RJ,EN,ChD} or Knot Atlas for the
details and the corresponding invariants. 
 
Following \cite{ChD}, the {\em \tax-presentation\,} 
of an
{\em iterated torus knots\,} will be 
$\t(\vec{\rr},\vec{\ss})$ for two sequences of
integers {\em of any signs\,}: 
\begin{align}\label{iterrss} 
\vec{\rr}=\{\rr_1,\ldots \rr_\ell\}, \ 
\vec{\ss}=\{\ss_1,\ldots \ss_\ell\} \hbox{\, such that\, 
gcd}(\rr_i,\ss_i)=1;
\end{align}
$\ell$ will be called the {\em length\,} of $\,\vec\rr,\vec\ss$.
The pairs $[\rr_i,\ss_i]$ are
{\em characteristic\,} or {\em Newton pairs\,}
for algebraic knots (such that $\rr_i,\ss_i>0$). We will 
call them Newton pairs too for arbitrary
(possibly negative) $\rr_i$ or $\ss_i$, which is allowed in our
considerations.

This presentation, referred to as the {\em $\tax$-presentation\,}, 
will be exactly the one needed 
in the DAHA approach. However it is not optimal for establishing
the symmetries of our polynomials and the justification 
that our construction depends only on the 
corresponding knot/link. We actually need 
the {\em cable presentation\,} for this, which is actually from
the definition of the corresponding iterated torus knots. 
It requires one more sequence of integers (possibly negative):
\begin{align}\label{Newtonpair}
\aa_1=\ss_1,\,\aa_{i}=\aa_{i-1}\rr_{i-1}\rr_{i}+\ss_{i}\,\ 
(i=2,\ldots,m). 
\end{align}
See  e.g., \cite{EN}. In terms of
the {\em cabling\,} discussed below, the corresponding knots
are as follows. First, $T(\rr,\ss)=C\!ab(\ss,\rr)(\unknot)$
(note that we transpose $\rr,\ss$ here), and then we set:
\begin{align}\label{Knotsiter}
\t(\vec{\rr},\!\vec{\ss})\rightsquigarrow
C\!ab(\vec{\aa},\!\vec{\rr})(\unknot)=
\bigl(C\!ab(\aa_\ell,\rr_\ell)\cdots C\!ab(\aa_2,\rr_2)\bigr)
\bigl(T(\rr_1,\!\ss_1)\bigr).
\end{align}

The first iteration here (application of 
$C\!ab$) is  $C\!ab(\aa_1,\rr_1)$, {\em not that for the last 
pair\,}, and then we proceed inverting the natural ordering.
\smallskip

{\sf Cabling}.
Knots and links will be always considered up to 
{\em ambient isotopy\,}; we use \,``$=$" for it.
The {\em cabling\,} $C\!ab(\aa,\bb)(K)$ of any oriented
knot $K$ in (oriented) $\S^3$ is defined as follows; 
see e.g., \cite{Mo,EN} and references therein. 
We consider a small $2$\~dimensional torus
around $K$ and put there the torus knot $T(\aa,\bb)$
in the direction of $K$, 
which is $C\!ab(\aa,\bb)(K)$ (up to ambient isotopy);
we set $C\!ab(\vec\aa,\vec\rr)\equal
C\!ab(\vec\aa,\vec\rr)(\unknot)$.   

This procedure depends on 
the order of $\aa,\bb$ and orientation of
$K$. We choose the latter in the standard way, matching  
the Mathematica package
Knot Atlas \cite{KA}; the parameter $\aa$ gives the number of
turns around $K$. This construction also depends
on the {\em framing\,} of the cable knots; we take the natural
one, associated with the parallel copy of the 
torus where a given cable knot sits (its parallel copy 
has zero linking number with this knot). It will be the 
same {\em standard framing\,} for iterated torus {\em links\,};
see below. Since the DAHA-invariants
are considered in this work up to powers of $q,t$, this 
will be sufficient. 

\smallskip

{\sf Topological symmetries.}
By construction,
$C\!ab(\aa,0)(K)=\unknot$ for any knot $K$
and  $C\!ab(\aa,1)(K)=K$ for any $\aa\neq 0$.
Accordingly, we have the following {\em reduction cases\,}:
\begin{align}\label{iter-triv}
&\hbox{If\, }
\rr_i=0 \hbox{\,\,for some\,} 1\le i\le \ell,\, \hbox{then}
\, \t(\vec{\rr},\vec{\ss})=\\
\t(\{\rr_{i+1},\cdots,\rr_{\ell}&\},
\{\ss_{i+1},\cdots,\ss_{\ell}\}),\  
\t(\vec{\rr},\vec{\ss})=\unknot
\for i=\ell.\notag\\
&\hbox{If\, } \rr_i=1,\, \ss_i\in \Z \hbox{\, for some\,\,} i, \,
\hbox{then}\, \t(\vec\rr,\vec\ss)=\label{iter-s1}\\
\t(\{\rr_1,\cdots,\rr_{i-1}&,
\rr_{i+1},\cdots,\rr_\ell\},\{\cdots,\ss_{i-1},
\,\ss'_{i+1}\,,\ss_{i+2},\cdots\}),\notag\\
\hbox{where\,\,\, } \ss_{i+1}'\!=\!\ss&_{i+1}\!+\!\ss_{i}\rr_{i+1} 
\hbox{\,\, if\,\, } i <\ell\hbox{\,\, (no\, $s'_{\ell+1}$\, 
for\, } i=\ell). \notag
\end{align}
\smallskip
Let us comment on the last relation; see
(\ref{Newtonpair}). Since $r_i=1$, one has:
$$
\aa_i\!=\!\aa_{i-1}\rr_{i-1}\!+\!\ss_i,\ \,
\aa_{i+1}\!=\!\aa_i \rr_{i+1}\!+\!\ss_{i+1}\!=\!
\aa_{i-1} \rr_{i-1}\rr_{i+1}\!+ (\ss_{i+1}\!+\ss_i \rr_{i+1}).$$
The pairs $\{\aa_i,\rr_i\}$ are sometimes called 
{\em topological or cable parameters\,}. Indeed,   
the isotopy equivalence of iterated knots generally can be
seen only at the level of these parameters (not at the level 
of Newton pairs).
\smallskip

Next, the symmetry $T(\rr,\ss)=T(\ss,\rr)$ results in
the following {\em transposition and reduction properties\,}.
For any $\vec{\rr},\vec{\ss},\vec{\aa}$ and for $m\ge 0$,
\begin{align}\label{iter-sym}
&\t(\vec{\rr},\vec{\ss})=
\t(\{\ss_1,\rr_{2},\ldots,\rr_{\ell}\},
\{\rr_{1},\ss_2\ldots,\ss_{\ell}\}),
\end{align}
\begin{align}\label{iter-symx}
&C\!ab(\aa_\ell,\rr_\ell)\!\cdots C\!ab(\aa_1,\rr_1)(T(m,1))\!=\!
C\!ab(\aa_\ell,\rr_\ell)\!\cdots C\!ab(\aa_1,\rr_1)(\unknot\,).
\end{align}

Then, switching from a knot $K$ to its {\em mirror image\,}, 
denoted by $K^\star$:
\begin{align}\label{iterorien}
C\!ab(\aa,\bb)(K^\star)&=
\bigl(C\!ab(-\aa,\bb)(K)\bigr)^\star \hbox{\,\, for any }
\aa,\bb \hbox{ with\,\,   gcd}(\aa,\bb)\!=\!1,\notag\\
\hbox{and\, }&C\!ab(-\vec\aa,\vec\rr)\!=\!
\bigl(C\!ab(\vec\aa,\vec\rr)\bigr)^\star,\ 
\t(\vec\rr, -\vec\ss)\!=\!
\bigl(\t(\vec\rr,\vec\ss)\bigr)^\star.
\end{align}
We note that the Jones and HOMFLY-PT polynomials for 
$K^\star$ are obtained from those for $K$ (can be a link)
by the formal conjugation of the parameters, which is 
$q\mapsto q^{-1}, a\mapsto a^{-1}$. This will hold
for the DAHA-Jones polynomials and
DAHA-superpolynomials, where the conjugation 
$t\mapsto t^{-1}$ must be added to that of $q,a$.

Furthermore,
changing the {\em orientation\,}, denoted by ``$-$", at the 
$i$\hbox{\tiny th}\, step, we obtain 
that for any $1\le i\le \ell$,
\begin{align}\label{iter2min}
-C\!ab(\vec\aa,\vec\rr)\!=\!
C\!ab&(\{\ldots\!,\aa_{i\!-\!1},\!-\aa_i,\aa_{i\!+\!1},\!\ldots\},
\{\ldots\!,\rr_{i\!-\!1},\!-\rr_i,\rr_{i\!+\!1},\!\ldots\}),
\\
(-)^{\ell-i+1}\,\t(\vec\rr,\vec\ss)\!=
\!\t&(\{\ldots,\rr_{i\!-\!1},
\!-\rr_i,\rr_{i\!+\!1},\ldots\},
\{\ldots\!,\ss_{i\!-\!1},\!-\ss_i,\ss_{i\!+\!1},\!\ldots\} ).
\notag
\end{align} 
The second transformation here results in the following cable: 
$$C\!ab(\{\ldots\!,\aa_{i\!-\!1},\!-\aa_i,-\aa_{i\!
+\!1},\!\ldots\},
\{\ldots\!,\rr_{i\!-\!1},\!-\rr_i,-\rr_{i\!+\!1},\!\ldots\}),$$
which explains the sign.

Changing the orientation of a knot or the simultaneous
change of the orientations of all components of
a link (equivalently, applying $\,\iota\,$ to the weights) does 
not influence its Jones and HOMFLY-PT polynomials;
so the DAHA-Jones polynomials must remain unchanged under
such a transformation (they are). We will mostly use
the symbol $\vee$ for the change of the orientation in this
work.
\medskip

\subsection{\bf From knots to links}\label{sec:knots-links}
Switching to links, we need to add 
{\em colors\,} to the cables above, which are  dominant
weights $b$. The \tax-presentation of a 
{\em torus iterated link\,} will be 
a union of $\kappa$ {\em colored knots\,} 
\begin{align}\label{tau-link}
\l_{(\vec\rr^{\,j},\vec\ss^{\,j})}^{\,\Up,\, (b^j)}=
\Bigl(\{\t(\vec\rr^{j},\vec\ss^{j}),\, b^j\in P_+\},\, 
j=1,\ldots,\kappa\Bigr) \hbox{\, together with}&\\
\hbox{\em the incidence} \hbox{\em\, matrix\,\,} \Up=(\up_{j,k}),
\hbox{\, where\, } 0\le \up_{j,k}=\up_{k,j}
\le \min\{\ell^j,\ell^{k}\}&,
\notag\\
\hbox{implie\,s that\,}\, [\rr_{i}^j,\ss_{i}^j]\!=\! 
[\rr_{i}^k,\ss_{i}^k]\hbox{\, for all\, }
1\le i \le \up_{j,k} \hbox{\, and any\, } 1\le j,k \le \kappa&. 
\notag
\end{align}
Here $\ell^j$ is the length of $\vec\rr^j=\{\rr_i^j\}$ and 
$\vec\ss^j=\{\ss_i^j\}$ for $1\le j\le \kappa$; we
naturally set $\up_{j,j}=\ell^j$. 
\smallskip

Subject to this above identification 
(for $i\le \up_{j,k}$),
the pairs $[\rr_i^j,\ss_i^j]$ will be treated as
{\em vertices\,} of a natural graph $\l$ determined by
$\Up$; these pairs we be called {\em \tax-labels\,}, of the 
corresponding vertices. The notation $[\,,\,]$
will be used exclusively for such labels $[\rr,\ss]$.

The {\em paths\,} are the sequences
of vertices with fixed $j$. The vertices 
for neighboring $i$ in the same path will be connected by 
the {\em edges}; the graph
has a natural orientation from $i$ to $i+1$ 
in each path.  Also, we will add the {\em arrowhead\,} 
at the end of each path, which is at $i=\ell^j$.
 
Such an {\em incidence graph\,}  (including the arrowheads)
is a union of trees, called
{\em subtrees\,}.
Every subtree has
at least one {\em base path\,}, the one that intersects all
other paths in this component. It also has  
a unique initial vertex (corresponding to $i=1$ in 
any base path). 
The colors $b^j$ will be assigned to
the arrowheads; thus the $j${\tiny th} path
corresponds to the knot $\t(\vec\rr^{j},\vec\ss^{j})$ 
colored by $b^j\in P_+$.  The graph can be empty (no vertices),
then it is a collection of paths that are pure arrowheads.
There can be several arrowheads from the same vertex, 
but one path always has one arrow. Topologically and in the DAHA 
construction, one can (technically) assume that the 
graph is a tree by adding an extra initial vertex with the 
label $[1,0]$ connected to its all subtrees. 


The $\aa$\~parameters can be calculated
along the paths exactly as we did for the knots
(starting from $i=1, \aa_1=\ss_1$); then
$\aa_i^j$ depends only on the corresponding vertex. 
The pairs 
$\{\aa_i^j,\rr_i^j\}$
will be called the {\em cab-labels} of the vertices. 
Actually only the \tax-labels will be needed in the DAHA 
constructions (we will mostly call them simply {\em labels} 
and use $[\,,\,]$ only for them). However
the cab-labels are necessary to explain the topological
symmetries (including the DAHA-Jones polynomials).

The torus knot colored by $b\in P_+$
will be denoted by $T_{\rr,\ss}^b$. Respectively,
$C\!ab_{\aa,\rr}^b(L)$, 
equivalently $C\!ab_{0,1}^b Cab_{\aa,\rr}(\mathbf L)$, will 
be the cable $C\!ab(\aa,\rr)(\mathbf L)$
of a link $\mathbf L$ colored by $b$.
The color can be attached only to the last $C\!ab$ in 
the sequence of cables. 
In the absence of vertices, the notation is $\unknot^{\,b}$ 
(the unknot colored by $b\in P_+$) or $C\!ab(0,1)^b$. {\em
We mostly use the same notation $\l$ for the graph and the
corresponding link $\mathbf{L}$.}

\smallskip 
The passage from the {\em \tax-presentation\,}
to the {\em cab-presentation\,}
is 
\begin{align}\label{Knotsiterx}
\l\bigl(\vec\rr^{j},\vec\ss^{j},\, 1\le j\le \kappa\bigr) 
\rightsquigarrow
\Bigl(\coprod_{j=1}^\kappa C\!ab(\vec{\aa}^j,\!\vec{\rr}^j)
\Bigr)(\unknot),
\end{align} 
where the composition and coproduct of cables is with respect to 
the tree structure and $C\!ab(\vec{\aa}^j,\!\vec{\rr}^j)=
\cdots C\!ab(\aa^j_2,\rr^j_2)T(\rr^j_1,\ss^j_1)$ is as in
(\ref{Knotsiter}).
In this work, the coproduct symbol $\,\hbox{\small$\coprod$}\,$
(sometimes omitted) 
is used when a union of cables is applied to the same link;
this union is disjoint but the result of its application will
generally have nonzero linking numbers. For a tree, 
the cab-presentation begins as follows:
\begin{align}\label{Knotsiterxx}
\Bigl(
\!\cdots\!\bigl(\!\coprod_{\ell^j=3}^{\hbox{\tiny tree}}
(C\!ab^{b^j}_{\aa_3^j,\rr_3^j} C\!ab_{\aa_2^j,\rr_2^j})\bigr)
\bigl(\!\coprod_{\ell^j=2}^{\hbox{\tiny tree}} 
C\!ab^{b^j}_{\aa_2^j,\rr_2^j}\bigr)
\bigl(\!\coprod_{\ell^j=1} C\!ab^{b^j}_{0,1}\bigr)
\Bigr) 
\bigl(T(\rr_1^1,\ss_1^1)\bigr).
\end{align} 
The coproduct for $\ell^j=1$ corresponds to pure arrowheads from
$i=1=j$, the next product is over single edges from the initial
point followed by arrowheads. The $2${\small nd} and
the $3${\small rd} coproducts are with respect to
the incidence tree, so must be understood as follows.
The third contains
$$
\bigl(C\!ab^{b^j}_{\aa_3^j,\rr_3^j}\coprod 
C\!ab^{b^k}_{\aa_3^k,\rr_3^k} 
\bigr)\bigl(C\!ab_{\aa_2^j,\rr_2^j}\bigr)
\hbox{\, if\, } \up_{j,k}=2,
$$
i.e. when the vertices labeled by
$[\rr_2^j,\ss_2^j]=[\rr_2^k,\ss_2^k]$
are identified in the incidence tree.
If $\up_{j,k}=3$, i.e. when the corresponding paths are 
different only by the arrowheads, this product contains
$$
\bigl(C\!ab^{b^j}_{0,1}\coprod C\!ab^{b^k}_{0,1}\bigr) 
\bigl(C\!ab_{\aa_3^j,\rr_3^j} C\!ab_{\aa_2^j,\rr_2^j}\bigr).
$$
This is similar for the product with
$\ell^j=2$ (when $\up_{j,k}=2$) and in general.

Here and generally $C\!ab(\aa_{i+1},\rr_{i+1})C\!ab(\aa_i,\rr_i)$
means the {\em composition\,} of cabling operations (for
the standard framing).  However, we will frequently
omit the symbol of coproduct between cables,
when it is clearly {\em not\,} the composition;
for instance the cables in (\ref{Knotsiterxx}) of the same
``level" $i=2,3$. 

Note that the property $C\!ab(\aa,1)(K)=K$ generally holds 
only for knots $K$. Applying $C\!ab(\aa,1)$ to a disjoint 
union of knots generally ties them up. For instance, 
$C\!ab(-1,1)$ of $m$ unknots $\,\coprod_{j=1}^m \unknot\ \,$
produces the {\,\em Hopf $\,m$\~link\,} for the standard
framing, which is $\mathbf{S}_m$\~symmetric
and with pairwise linking numbers all equal to $-1$. 
This corresponds to the tree with one vertex $[1,-1]$
and $m$ colored
arrows from it. Note that applying $C\!ab(0,1)$ here simply 
produces a union of $m$ arrowheads, without any vertices and 
edges. 
\smallskip

Using
the symmetry from (\ref{iter-s1}) requires recalculating
all $\ss$ after $[\rr^j_i=1,\, \ss^j_i\in \Z]$ in all 
paths through it. 
Relation (\ref{iter-triv}) holds for links, but now it must be
understood as deleting all vertices in the paths through
a vertex with $\rr_{i_o}^{j_o}=0$ from the first one
(i.e. for all $1\le i\le i_o$ in any paths containing this 
vertex).
The paths which share some vertices with those affected 
remain untouched, so the matrix $\Up$ and 
the incidence graph must be recalculated accordingly, which
can result in extra subtrees. 
Also, (\ref{iterorien}) and 
(\ref{iter2min}) hold if $\star$ and ``$-$" are applied to the 
whole link. 

\smallskip

\subsection{\bf Splice diagrams} \label{sec:Splice}
Let us extend the previous construction
to the case of a pair of incidence graphs
$\{\mathcal{L},\,'\!\mathcal{L}\}$, where the latter can be
with or without $\vee$. The {\em twisted union\,} of the 
corresponding links is defined
as follows. The cabling construction provides a canonical 
embedding of the iterated torus links into   
the solid torus. We put the links for $\l$ and 
$'\!\mathcal{L}$ into the horizontal solid torus
and the complementary vertical one. 
The presence of $\vee$
in $\,'\!\mathcal{L}^\vee$ means changing the orientation
of this component; for instance, $\{\square,\square\}$ and
$\{\square,\square^\vee\}$ represent uncolored Hopf
$2$\~links with the linking number $lk=-1$ and $+1$ 
correspondingly.


{\em Recall that the notation $\l$ or $\{\l,\,'\!\l\}$ 
is used for trees and also for
the corresponding link $\mathbf{L}$ (depending on the context).}
\smallskip

\thicklines  
Using the language of splice diagrams from \cite{EN},
the pair $\{\l,\,'\!\l^\vee\}$  corresponds to  
$\circ_1$ and $'\circ_1$ in  $\mathcal{L}$ and 
$'\!\mathcal{L}$ connected by an {\em arc\,}:
\begin{equation*}
\begin{picture}(60,20)
\put(8,11){\oval(10,15)[l]}
\put(10,0){$\,'\!\!\circ_1\cdots$}
\put(10,15){$\ \circ_1 \cdots$}
\end{picture}\!\!\!\!\!\! \hbox{\,\, or\,\,\,\,\,\,\,}
\begin{picture}(60,20)
\put(8,11){\oval(10,15)[l]}
\put(10,0){$\,'\!\!\circ_1\cdots^\vee$}
\put(10,15){$\ \circ_1\cdots$}
\end{picture}\!\!\!\!\!\!\!,
\end{equation*}
where adding $\vee$ changes  the 
orientation of the link with prime. For instance,
the Hopf links for $lk=\pm 1$ are represented as
$\Yboxdim7pt \yng(1)\longleftrightarrow\yng(1)$
(for $+$)
and  $\yng(1)\longleftrightarrow\yng(1)^\vee$.
Importantly, the transposition $\mathcal{L} \leftrightarrow 
\,'\!\mathcal{L}$ does not change the
output if $\vee$ is present or absent; the
total change of the orientation does not influence the link 
invariants we will consider.  

This is  topologically equivalent to combining the 
trees via an additional vertex $\circ\,$ labeled by 
$[1,\pm1]$, i.e. using the tree 
$\circ \rightrightarrows \genfrac{}{}{0pt}{1}{ \circ_1\,\cdots}
{ '\!\circ_1\cdots}$, where $\circ_1, \, '\!\circ_1$ from
$\l,\,'\!\l$ are connected
via such 
intermediate $\circ\,$. In the DAHA-construction, this 
is governed by Theorem \ref{EVALTAUM}.
The corresponding
$\aa$\~parameters in the {\em cab-labels\,}
must be then recalculated, 
since we now begin with a new vertex $\circ$. Thus the
twisted union (the arc-connection) actually results in
one tree; see (\ref{T2-3--1-2}) for a concrete example.

\smallskip

Let us translate more systematically our combinatorial data 
into the
language of splice diagrams. See \cite{EN} for details.
For one tree $\mathcal{L}$, the construction is as above
(including the directions and arrowheads), but 
we need to switch from the  labels $[\rr,\ss]$ in the 
$\tax$-presentation to the corresponding $\{\aa,\rr\}$ in 
the cab-presentation. Graphically, the passage is as follows:
\vskip -0.2cm

\begin{equation*}
\begin{picture}(40,20)

\put(0,10){$\rightarrow \circ\ \rightdotsarrow$}
\put(5,0){$[\rr,\ss]$}
\end{picture}
\begin{picture}(10,20)
\put(0,7){$\ \ \leftrightsquigarrow$}
\end{picture}\ \ \ \ \ 
\begin{picture}(65,40)
    \put(35,9){\line(2,1){20}}
    \put(35,9){\line(2,-1){20}}
    \put(25,9){\line(-1,0){20}}
    \put(30,11){\line(0,1){20}}
    \put(23.5,6){$\bigoplus$}
    \put(36,-2){$1$}
    \put(36,14){$1$}
    \put(56,4){$\vdots$}
    \put(20,0){$\aa$}
    \put(24,16){$\rr$}
    \put(27,30){$\circ$}
\end{picture}.
\end{equation*}

Note that we extend \cite{EN} by adding colors
to the arrowheads. Also, recall that $\aa_1=\ss_1$ for 
the initial vertices in any {\em paths\,}. 
\smallskip

In the case of two trees
$\{\mathcal{L},\,'\!\mathcal{L}^\vee\}$ (note $\vee$) 
the connection by an arc described above
corresponds to the
following splice diagram:
\begin{align}
\begin{picture}(60,85)
\put(0,20){
    \put(8,11){\oval(10,15)[l]}
    \put(10,0){$'\circ_1\rightarrow 
\,'\!\!\mathcal{M}\ \rightdotsarrow$}
    \put(10,15){$\,\circ_1\rightarrow 
\,\mathcal{M}\ \rightdotsarrow$}
    }
\put(0,45){$\![\rr_1,\ss_1]$}
\put(0,10){$\!['\rr_1,'\!\ss_1]$}
\end{picture}
\begin{picture}(20,35)
\put(5,27){$\ \ \ \ \,\leftrightsquigarrow$}
\end{picture}\ \ \ \ \ \ 
\begin{picture}(95,85)
\put(25,32){\oval(50,40)[l]}
\put(0,3){  
    \put(35,9){\line(1,0){20}}
    \put(75,11){\vector(2,1){20}}
    \put(75,7){\vector(2,-1){20}}
    \put(31.5,11){\line(0,1){20}}
    \put(25,6){$\bigoplus$}
    \put(40,-1){$1$}
    \put(73,-3){$1$}
    \put(73,16){$1$}
    \put(57,6){$\!'\!\widetilde{\mathcal{M}}$}
    \put(95,4){$\vdots$}
    \put(10,-1){$'\aa_1$}
    \put(21,16){$\!\!'\rr_1$}
    \put(29,30){$\circ$}
    }
\put(0,43){
    \put(35,9){\line(1,0){20}}
    \put(75,11){\vector(2,1){20}}
    \put(75,7){\vector(2,-1){20}}
    \put(31.5,11){\line(0,1){20}}
    \put(25,6){$\bigoplus$}
    \put(40,-1){$1$}
    \put(73,-3){$1$}
    \put(73,16){$1$}
    \put(57,6){$\widetilde{\mathcal{M}}$}
    \put(95,4){$\vdots$}
    \put(10,-1){$\,\aa_1$}
    \put(21,16){$\rr_1$}
    \put(29,30){$\circ$}
    }
\end{picture}\ ,
\label{splicefig}
\end{align}
where $\widetilde{\mathcal{M}}$ and 
$\,'\!\widetilde{\mathcal{M}}$ are 
splice diagrams made from
graphs 
$\mathcal{M}$ and $\,'\!\mathcal{M}$. For the
sake of simplicity of this figure, we assume 
that $\mathcal{L},\,'\!\mathcal{L}$ have only single
edges from their initial 
vertices $\circ_1,\,'\!\circ_1$ to $\mathcal{M},
\,'\!\mathcal{M}$, the remaining portions. The general 
correspondence is totally similar.

We note that adding $\vee$ to $\,'\!\mathcal{M}$
(the change of the orientation of this part) 
in the language of splice diagrams is as follows. One adds
a new {\em trivalent
node\,} to the arc
with weights $+1$ on the edges going to trees and  
$-1$ for the third edge, calling {\em a leaf\,}; see
\cite{EN}, Theorem 8.1, Statements 2 and 3 for the details.
The nodes with weights $+1$ (both) can be deleted 
from the diagram; 
in DAHA, this may result in some $q^m$, which is trivial due
to the normalization. 
%
%
Assigning $-1$ to the leaf as above
changes the orientation of one of the components; it can be
$\mathcal{M}$  or $\,'\!\mathcal{M}$,
up to the total change of the orientation. Algebraically,
this will correspond to the DAHA-relations (\ref{lsymlprime}),
(\ref{vphiotavee}).
See also (\ref{fig:splice}). 

\comment{
The relations with trees $\circ 
\rightrightarrows \genfrac{}{}{0pt}{1} { \circ_1\,\cdots}{ 
'\!\circ_1\cdots}$ becomes clear if one notes that in splice 
diagram language this corresponds to adding a new trivalent node 
with weights $1$ on the edges going to trees and an $\pm 1$ weight 
going to a leaf. If all weights are positive it can be removed 
without changing the link. If the leaf edge has weight $-1$ it can 
be removed with simultaneous change of all orientation on one of 
the trees (see \cite{EN}, Theorem 8.1, statements 2 and 3). These 
are exactly the two types of topological construction we presented 
above. 
}

\smallskip
We actually do not need the splice diagrams too much in this work.
However, they are an important and convenient tool
for the classification of iterated torus links, 
useful to analyze the topological  symmetries.  The latter
are one-to-one with the symmetries of our DAHA construction.
Also, the splice diagrams establish the connection with $3$-folds
and can be effectively used for calculating some invariants,
including the linking numbers, with applications to plain
curve singularities.
\smallskip

We note that our method can not give all 
splice diagrams of solvable type, though {\em all\,} 
algebraic links can be
reached by our construction. See Theorems 9.2, 9.4 in \cite{EN}). 
For instance, the DAHA-approach does not provide (so far)
the following diagram, corresponding to the {\em granny knot\,}:

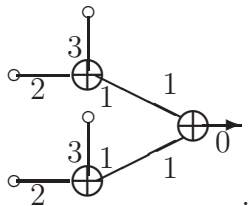
\begin{equation*}
\begin{picture}(90,80)
\put(0,40){
    \put(1,6){$\circ$}
    \put(35,9){\line(2,-1){32}}
    \put(25,9){\line(-1,0){20}}
    \put(32,11){\line(0,1){20}}
    \put(25,6){$\bigoplus$}
    \put(36,-3){$1$}
    \put(10,0){$2$}
    \put(24,16){$3$}
    \put(29,30){$\circ$}
    }
\put(60,20){
    \put(0,-10){$1$}
    \put(0,22){$1$}
    \put(20,0){$0$}
    \put(5,7){$\bigoplus$}
    \put(15,10){\vector(1,0){15}}
    }
\put(0,0){
    \put(1,6){$\circ$}
    \put(35,9){\line(2,1){33}}
    \put(25,9){\line(-1,0){20}}
    \put(32,11){\line(0,1){20}}
    \put(25,6){$\bigoplus$}
    \put(36,14){$1$}
    \put(10,0){$2$}
    \put(24,16){$3$}
    \put(29,30){$\circ$}
    }
\end{picture}.
\end{equation*}

It represents the connected sum of two trefoils with the
same orientation. 
The {\em connected sum\,} is the only one from the three
operations in \cite{EN} necessary to  obtain
arbitrary solvable links in $\S^3$ that is generally missing 
in the DAHA-approach.
 These three are  {\em disjoint sum}, 
{\em connected sum} and {\em cabling}.
We note that some $3$\~component connected sums
can be obtained, for instance the connected sum of two
Hopf $2$\~links (the chain of $3$ unknots, where the first
and the last are not linked); see (\ref{T10-10}) below.
\smallskip


{\sf Algebraic links.\,} 
We provide here only basic facts; see \cite{EN} for details and 
references, especially Theorem 9.4 there.
Generally, one begins with a polynomial equation $f(x,y)=0$
considered in a neighborhood of an isolated singularity
$0=(x=0,y=0)$. Its intersection 
with a small $3$-dimensional sphere in 
$\C^2$ around $0$ is called an {\em algebraic link\,}.
Assuming that $\rr^j_i,\ss^j_i>0\,$, 
any {\em tree\,} $\l^\Up_{(\vec\rr^j,\vec\ss^j)}$  
(in the \tax-presentation) corresponds to a germ of
{\em plane curve singularity\,}  at $0$. If these 
inequalities hold, the tree will be called {\em positive\,}.

Such germs are unions of  
{\em unibranch\,} components  for the paths of $\Up$
(numbered by $j$),  which are given as follows: 
\begin{align}\label{yxcurve}
y = c^j_1\,x^{\ss^j_1/\rr^j_1}
(1+c^j_2\,x^{\ss^j_2/(\rr^j_1\rr^j_2)}
\bigl(1+c^j_3\, 
x^{\ss^j_3/(\rr^j_1\rr^j_2\rr^j_3)} 
\Bigl(\ldots\Bigr)\bigr)) \hbox{\, at\, } 0. 
\end{align} 
The parameters $c_i^j\in \C$ are sufficiently general here.
The simplest example is the equation
$y^{\rr\kappa}= x^{\ss\kappa}$ under gcd$(\rr,\ss)=1$,
which corresponds to the torus link $T(\rr\kappa,\ss\kappa)$
with $\kappa$ knot components isotopic to $T(\rr,\ss)$; the
pairwise linking numbers here are all $\rr\ss$. See
also Section \ref{sec:Alex}.

The unibranch components and the corresponding 
pairwise linking numbers uniquely determine the corresponding
germ due to the Reeve theorem. The linking number between the 
branches corresponding to the paths with the indices $j\neq k$ is
\begin{align}\label{linkjk}
\hbox{lk}(j,k)=\aa^{j}_{i_o}\rr_{i_o}^j\,
\Bigl(\,\prod_{i=i_o+1}^{\ell^j} \rr^j_i\Bigr)
\Bigl(\,\prod_{i=i_o+1}^{\ell^k} \rr^k_i\Bigr),
\where i_\circ=\up(j,k).
\end{align}
This formula correctly gives lk$(j,k)=0$ if we
set $\aa^j_{i_\circ}=0$ when $\up(j,k)=0$ (i.e.
allow graphs $\l$ here, not only trees), but then 
the corresponding links will become non-algebraic.
All linking numbers must be strictly positive for 
algebraic links. 
\smallskip

Arbitrary algebraic links can be obtained using
this construction and
the {\em twisted union\,} (above) for the pairs
of positive algebraic trees with $\vee$ added to the second
subject to the inequality 
$\,'\ss_1 \ss_1>\,'\rr_1 \rr_1$ for the first vertices of
these two trees. Then $\{\l,\,'\!\l^\vee\}$ is
called a {\em positive pair\,}.

See \cite{EN} and (\ref{splicefig}) concerning the algebraic
links; the linking numbers between the branches  
of $\l$ and $\,'\!\l$ in their twisted union
are (full) products of $\rr^j_i$ and $\,'\rr^k_i$
over $i\,$ in the corresponding paths.

We note that the theory in \cite{EN} is without colors, as well as
that in \cite{ORS}. Attaching colors to the branches
requires more 
involved algebraic-geometric considerations; see \cite{Ma} for 
the case $t=q$. 
\medskip

\subsection{\bf DAHA-Jones polynomials}
They will be defined for an arbitrary (reduced,
irreducible) root system $R$ and its twisted
affinization $\tilde{R}$. The notations and
formulas are from the previous sections. The
combinatorial data will be the \tax-labeled graphs 
$\l_{(\vec\rr^{j},\vec\ss^{j})}^{\,\Up, (b^j)}$
from (\ref{tau-link}). Recall that 
$$
1\le j\le \kappa,\,\, \vec\rr^{\,j}=\{\rr_i^j\},\,\,
\vec\ss^{\,j}=\{\ss_i^j\},\,\, 1\le i\le \ell^j,
$$
and $\Up$ is the {\em incidence graph/matrix\,},  
and the {\em arrowheads\,} (at the ends of all {\em paths\,}) 
are colored by $b^j\in P_+$. The incidence graph is 
not supposed to be connected here and
the paths can contain no vertices; see 
(\ref{tau-link}). The construction below will be for
{\em two\,} arbitrary such graphs $\l$, $'\!\l$ (the second
can be empty).

The choice of the {\em integral form\,} of the $P$\~polynomials 
plays an important role in the following. Recall from 
Theorem \ref{REGTHM} that $\tilde{P}_b=
\tilde{\n}_b P_b^\circ$ for $b\in P_+$ is
$q,t$\~integral. This is not always the smallest 
$q,t$\~integral normalization of $P_b$. 
The following normalization is.

Let us consider the denominators of all coefficients of 
$P_b^\circ$ as polynomials in terms of $q$ with the coefficients 
in $\C(t_\nu)$ (actually in $\Z[t_\nu]$). Then let 
$\ddot{\n}_b$ be their least common multiple, {\em LCM\,}. 
We will assume that its constant term is $1$; then it is a  
polynomial in terms of $q$ with the coefficients
from $\Z[t_\nu]$. The polynomial 
$\ddot{P}_b\equal\ddot{\n}_bP_b^\circ$ is the 
minimal $q,t$\~integral
form of $P_b$. Accordingly, we set for $b^1,\ldots,b^m\in P_+$:
\begin{align}\label{lcmkrho}
&(b^1,\ldots,b^m)_{ev}^\sim=
LC\!M\bigl(\tilde{\n}_{b^1},\ldots,\tilde{\n}_{b^m}\bigr),\\
&(b^1,\ldots,b^m)_{ev}^{\cdot\cdot}=
LC\!M\bigl(\ddot{\n}_{b^1},\ldots,\ddot{\n}_{b^m}\bigr).\notag
\end{align}
Recall that the calculation of $(b^1,\ldots,b^m)_{ev}^\sim$
is entirely combinatorial; it is an explicit product of binomials
in the form of (\ref{binaffexp}).

In type $A$, one can also take here the $J$\~polynomials
$J_{\la^j}$ for $\la^j=\la(b^j)$; see (\ref{P-arms-legs}).
Then we set:
\begin{align}\label{P-arms-legs-LCM}
&(b^1,\ldots,b^m)_{ev}^J\!=\!(\la^1,\ldots,\la^m)_{ev}^J\!=
\!LC\!M^{^\odot}\bigl(J_{\la^1}(t^{\rho}),\ldots,
J_{\la^m}(t^{\rho})\bigr),
\end{align}
where ${}^\odot$ means here that we
normalize $LC\!M$ by the condition that it is a 
$q,t$\~polynomial with the constant term $1$. The latter
practically means that we 
ignore the factors in formula
(\ref{stabeval}) for $J_{\la}(t^{\rho})=J_{\la}(t^{-\rho})$
before the product there (upon $a=-t^{n+1}$). 

One has the following combinatorially transparent formula:
\begin{align}\label{P-arms-legs-union}
&(\la^1,\ldots,\la^m)_{ev}^J\ =\ 
(\la^1\!\!\vee\cdots\vee\!\la^m)_{ev}^J\,, \hbox{\,\,\, where}\\ 
&\la^1\!\!\vee\cdots\vee\!\la^m 
\hbox{\, is the union of diagrams\, }
\{\la^j\}.\notag
\end{align}
The latter union is by definition the smallest Young diagram 
containing all diagrams $\la^1,\ldots,\la^m$.

Going back to arbitrary root systems,
we can take
$P_b\, (b\in P_+)\,$ themselves when 
$k_{\lng}\!=\!1\!=\!k_{\sht}$, equivalently 
$t_{\nu}\!=\!q_{\nu}$. Then the corresponding
$P^{(\!1\!)}_b=P_b(X; t_\nu=q_\nu)$ do not depend on $q,t$ for any
(reduced, irreducible) root system $R$. We set in this case:
\begin{align}\label{P-arms-legs-k1}
&(b^1,\ldots,b^m)_{ev}^{(\!1\!)}= 
LC\!M^{{}^\odot}\bigl(P^{(\!1\!)}_{b^1}\!(q^{\rho_k}\!),\ldots,
P^{(\!1\!)}_{b^m}\!(q^{\rho_k}\!)\bigr).
\end{align}
\smallskip
Note that the $J$\~polynomials in the $A_n$\~case
are not minimal integral even when $t=q$.
However, using them vs. $P^{(\!1\!)}_\la$ for $t=q$
makes the corresponding HOMFLY-PT polynomials 
$a,q$\~integral and has other advantages. This is
important only for links; the choice of the integral
form does not influence the DAHA-construction for knots.
\smallskip

We represent torus knots $T(\rr,\ss)$    
by the matrices $\ga_{\rr,\ss}\in GL_{\,2}(\Z)$ with the 
first column $(\rr,\ss)^{tr}$ ($tr$ is the transposition)
for $\,\rr,\ss\in \Z$ provided \,gcd$(\rr,\ss)=1$. 
Let $\hat{\ga}_{\rr,\ss}\in GL_{\,2}^{\wedge}(\Z)$ be
{\em any\,} pullback of $\ga_{\rr,\ss}$.

Obviously $\,(\rr,\ss)\,$ can be lifted
to $\,\ga\,$ of determinant $1$ and, accordingly, 
to $\hat{\ga}\in PSL^\wedge_{\,2}(\Z)$
generated by $\{\tau_{\pm}\}$. I.e. the usage of $\eta$
can be avoided. However, $\eta$ results in important 
symmetries of the $J\!D$\~polynomials. Similarly, 
the {\em enhanced\,} $GL_{\,2}^{\wedge}(\Z)^\ze$ is actually
not needed in what will follow,
since all the quantities below (inside $\{\ldots\}_{ev}$)
will be $\ze_{u,v}$\~eigenvectors. Controlling 
the corresponding $\Pi'$\~characters is of importance, but 
not really significant in the construction below.


\comment{
For a polynomial $F$ in terms of 
fractional powers of $q$ and $t_\nu$, 
the {\em tilde-normalization\,}
$\tilde{F}$ will be the result of the division of $F$
by the lowest $q,t_\nu$\~monomial, assuming that it 
is well defined. We put $q^\bullet t^\bullet$ for
a monomial factor (possibly fractional)
in terms of $q,t_\nu$.}

\begin{definition}\label{JONITER}
{\sf Ingredients.}
Let $R$ be a reduced irreducible root system and $q$  not
a root of unity.
Recall that $H\mapsto H\!\!\Downarrow\, \equal H(1)$, where the
action of $H\in \HH$ in $\v$ is used.
The construction is in terms of two {\sf graphs/links\,}
from (\ref{tau-link})
\begin{align}\label{links-x-y}
\l=\l_{(\vec\rr^j,\vec\ss^j)}^{\Up,(b^j)},\
'\!\l=\, '\!\l_{('\vec\rr^j,'\vec\ss^j)}^{\,'\Up,('\!b^j)}
\where b^j,\,'\!b^j\in P_+,
\end{align} 
$$
1\le j\le \kappa,\,'\!\kappa\for \l,\, '\!\l,\   
\vec\rr^j= (\rr_i^j \mid 1\le i\le \ell^j),\,  
'\vec\rr^j=('\!\rr_i^j \mid 1\le i\le\, '\!\ell^j).
$$ 
The $j$\~ranges for \,$\vec\ss^j,'\!\vec\ss^j$
are those for \,$\vec\rr^j,'\!\vec\rr^j$.
Let us lift the columns $(\rr^j_i,\ss^j_i)^{tr}$,
 $('\!\rr^j_i,'\!\ss^j_i)^{tr}$\, to\,
$\ga_i^j,\,\!'\ga_i^j$ and then
to $\hat{\ga}_i^j,\, '\hat{\ga}_i^j\in GL_{\,2}^{\wedge}(\Z)$ 
as above.  

The construction also requires the choice of
the {\sf integral forms}
$\hat{P}_b$ for the polynomials $P_b$. Considering
the pairs $\{i,j\}$ as vertices
of the incidence graph for $\l$, we begin with letting
$$
\hat{\P}_{\ell^j+1}^j\equal\hat{P}_{b^j},\ 
\ga_{\ell^j+1}^j\equal\hbox{id} \for 1\le j\le \kappa;
$$ 
recall that $\ell^j=0$ when the $j$\hbox{\tiny th}
path contains only one arrowhead and $\up_{j,k}=0$
if the corresponding paths do not intersect.
\smallskip

{\sf Pre-polynomials.}
For a given {\sf path} with the index $j$, we 
define the polynomials $\hat{\P}_{i}^{j}$ 
by induction with respect
to $i$, starting with $i=\ell^{j}$ down to $i=0$:
\begin{align}\label{hatPij}
\hat{\P}_{i}^{j}=
\prod_{1\le k\le \kappa}^ {\up(k,j)=i} 
\bigl(\,\hat{\ga}_{i+1}^{k}(\hat{\P}_{i+1}^{k})
\!\!\Downarrow\,\bigr)
\end{align}
i.e. the last product is over all paths $k$ that enter
(intersect) the path for $j$ exactly at the index 
$i\ge 0$, including $k=j$\, when\, $i=\ell^j$.
Note that
$\hat{\P}_{\ell^{j}}^{j}=\prod_{1\le k\le 
\kappa}^ {\up(k,j)=\ell^j} \hat{P}_{b^k}$ for a {\sf base path}
$j$, where this product is over all 
{\sf arrowheads\,} from (originated at) the 
vertex $\{i=\ell^j,j\}$. 

Since $i=0$ is allowed in (\ref{hatPij}),
$\hat{\P}_{0}^{j}$ actually depends only on the corresponding
 subtree for 
{\sf any\,}
path $j$\, there. If $\l$ is the union of subtrees, then 
$\hat{\P}_{0}^{j}$
is the product of the corresponding polynomials
$\bigl(\,\hat{\ga}_{1}^{k}
(\hat{\P}_{1}^{k})\!\!\Downarrow\,\bigr)$
over all these subtrees (for any choices
of paths $k$ there). So we may set
$\,\hat{\P}_{0}^{tot}=\hat{\P}_{0}^{j}$. The 
polynomial $'\hat{\P}_{0}^{tot}$ for $'\!\l$
is defined in the same way.

{\sf Finale.}
Using the notations $\mathbf{b}=(b^j),\,\mathbf{'b}=('b^j)$,
the {\sf DAHA-Jones polynomial} for the integral forms
$\hat{P}_{b^j}$, $\hat{P}_{'\!b^j}$  and a certain 
fixed index $1\le j_o\le \kappa$ or $j_o=\emptyset$ is
(which determines the normalization):
\begin{align}\label{jones-hat}
\hat{J\!D}_{\,(\vec\rr^j,\,\vec\ss^j)\,,
\,('\!\vec\rr^j,\,'\!\vec\ss^j)}
^{R,\,j_o,\,\Up,\,'\!\Up}((b^j),('\!b^j)\, ; \,q,t)\, =&\, 
\hat{J\!D}^{j_o,\,\Up,\,'\!\Up}_{\,(\vec\rr^j,\,\vec\ss^j)
\, , (\,'\!\vec\rr^j,\,'\!\vec\ss^j)}
(\mathbf{b},\mathbf{'b}\, ; \,q,t)\\ 
=\,\hat{J\!D}^{j_o}_{\l,\,'\!\l}\equal 
\Bigl\{\, \vph(\,'\hat{\P}_0^{tot})&\,\hat{\P}_0^{tot}/
\hat{P}_{b^{j_o}}(q^{-\rho_k})\,\Bigr\}_{ev},
\notag
\end{align}
where there is no division by $\hat{P}_{b^{j_o}}(q^{-\rho_k})$
for $j_o=\emptyset$. Due to the definition of $\vph$ from 
(\ref{starphi}), $\vph(\,'\hat{\P}_0^{tot})$
is simply  $'\hat{\P}_0^{tot}(Y^{-1})=\,
'\hat{\P}_0^{tot}(X\!\mapsto\! Y^{-1})$. Taking $Y^{+1}$
here, we obtain:
\begin{align}\label{jones-hat-vee}
\hat{J\!D}^{j_o}_{\l,\,'\!\l^\vee}\equal 
\Bigl\{\,'\hat{\P}_0^{tot}(Y)&\,\hat{\P}_0^{tot}/
\hat{P}_{b^{j_o}}(q^{-\rho_k})\,\Bigr\}_{ev}.
\end{align}
\end{definition}
\vskip -1.2cm \sq
\smallskip

{\sf Using automorphisms.}
We note that the pre-polynomial
$\hat{\P}_0^{tot}$ is invariant 
with respect to the action of $\ze_{u,1}$ in
the polynomial representation; use
(\ref{zetauvrel}) and (\ref{zetheE}) to calculate
the corresponding character. 

{\em The role of $\vph$\,.} 
Applying $\vph$ inside the coinvariant gives that
\begin{align}\label{lsymlprime}
\hat{J\!D}^{\bullet}_{\l,\,'\!\l}=
\hat{J\!D}^{\bullet}_{'\!\l,\,\l}\for \bullet=j_o,\, min,
\, \emptyset,
\end{align}
where the normalizations $\bullet$ must be the same in both 
sides, i.e. the division must be
by the same $\hat{P}_{b^{j_o}}(q^{-\rho_k})$ for $j_o$.
Topologically, this relation means that adding the {\em meridian}
$\,'\!\l$ to (the link associated with) $\l$
is isotopic to adding  the meridian $\l$ to $\,'\!\l$.
\smallskip

The operator $\vph(\,'\hat{\P}_0^{tot})$
does not commute with $\hat{\P}_0^{tot}$, which is
an operator in terms of $X$. In type $A$,
the commutator relations
between such and similar operators 
are part of the theory of {\em elliptic Hall algebra\,}
(isomorphic to spherical DAHA); see \cite{SV}. This
connection is important for the 
{\em Hopf links\,} and the {\em DAHA-vertex\,}. See also
some remarks concerning the toric $q,t$\~skein in 
Section \ref{sec:skein}. 

The  {\em toric skein algebra}
was identify with the elliptic Hall algebra for $t=q$ in
\cite{MoS} (and with the spherical DAHA). 
This was used to check the coincidence of our superpolynomials
at $t=q$ with the HOMFLY-PT polynomials for torus iterated
{\em knots}, conjectured in \cite{ChD} (justified for any 
such knots for $A_1$ and torus knots in \cite{CJ}). 
Combining \cite{MoS} with 
our approach here proves Part $(i)$ of Conjecture \ref{CONCONJ}.

This implication is of
course important, but direct using 
{\em Rosso-Jones formula\,} and the DAHA-shift operators
is also quite relevant, as it was demonstrated in 
\cite{CJ,ChD,ChE}. We already have a sketch of the proof of
Part $(i)$ of Conjecture \ref{CONCONJ} below based on this
method, which can be generally extended to any root 
systems, i.e. applied to connect arbitrary WRT-invariants
with the corresponding DAHA-Jones polynomials at $t=q$. 
\smallskip

Continuing with the automorphisms, we note that 
$\vph$  maps $\tau_{\pm}\mapsto \tau_{\mp}$
in the construction of 
$\,'\hat{\P}_0^{tot}$, 
without changing
the order of the corresponding tau-matrices. This
can be used to establish the symmetry
$T(\rr,\ss)\leftrightarrow T(\ss,\rr)$,
but generally not inside the cables;
the application of $\vph$ is not compatible with the projection
$\Downarrow\,$ onto the polynomial presentation $\v$.
The following ones are. 

\smallskip
{\sf Automorphisms $\eta, \tau_-, \iota\,, \si^2$.} 
They are compatible 
with $\Downarrow\,$, which fact is the key in proving the
symmetries of the $J\!D$\~polynomials. For instance,
adding the vertex $[\rr_0=1,\ss_0=\mm]$ to a tree 
$\l$ results in the same $\hat{J\!D}$\~polynomial as
the change of $[\rr_1,\ss_1]$ in $\l$  by 
$[\rr_1, \mm\rr_1+\ss_1]$
(for every path, if there are several). 
This gives (the key) relations from (\ref{iter-s1});
see Part $(ii)$ of Theorem \ref{MAINTHM}. 

Concerning $\iota\,$, the  definition from
(\ref{jones-hat-vee}) is actually a particular
case of (\ref{jones-hat}), where all colors $'b^j$ are
replaced by $\iota(\,'b^j)$. We use that $\iota$ 
preserves the coinvariant and naturally
acts in $\HH$; see (\ref{evfunct},\ref{iotaXY}).
Also, $\,\iota\cdot \hat{\ga}=
\hat{\ga}\cdot \iota$ for $\ga\in GL_2(\Z)$. We 
use these properties of $\iota$ and similar
ones for $\eta$ in the theorems below.
\smallskip

{\em Adding $\vee$.} This changes 
$'\hat{\P}_0^{tot}(Y^{-1})=\si('\hat{\P}_0^{tot})$ in
the definition above by 
$'\hat{\P}_0^{tot}(Y)=\si^{-1}('\hat{\P}_0^{tot})$.
Equivalently,  we can send $'b\mapsto\, 'b^\iota$. 
Equivalently, this operation is
the change of the sign of the first label (or labels if
it is not a tree) in $\,'\!\l$:\,
$[\,'\rr_1^j,\,'\ss_1^j\,]\mapsto [-'\rr_1^j,-'\ss_1^j\,]$.

Since $\iota$ commutes with $\vph$, we have the relation:
\begin{align}\label{vphiotavee}
\hat{J\!D}^{\bullet}_{\l,\,'\!\l^\vee}=
\hat{J\!D}^{\bullet}_{\l^\vee,\,'\!\l}=
\hat{J\!D}^{\bullet}_{'\!\l,\,\l^\vee},
\end{align}
where the normalization $\bullet$ must be the same
in all three formulas (say, $\bullet=\emptyset$ or $\bullet=min$).
Indeed, the total change of the orientation of
the link represented by the pair $\{\l,\,'\!\l\}$ 
must not influence the $\hat{J\!D}$\~polynomial. 
Recall that topologically, adding $\vee$ to
$'\!\l$ is switching the orientation of the corresponding link.
\smallskip

{\sf Iterated knots.\,}
In the case of torus iterated {\em knots\,} (when there is 
only one path) and in the absence of $'\!\l$,
we arrive at formula (2.12) from \cite{ChD}:
\begin{align}\label{jones-ditx}
& J\!D_{\vec\rr,\vec\ss}
(b;q,t)\! =\!  
\Bigl\{\hat{\ga_{1}}\Bigr(
\cdots\Bigl(\hat{\ga}_{\ell-1}
\Bigl(\bigl(\hat{\ga}_\ell(P_b)/
P_b(q^{-\rho_k})\bigr)\!\Downarrow
\Bigr)\!\Downarrow\Bigr) \cdots\Bigr)\Bigr\}_{ev}.
\end{align}
It includes only one $b\in P_+$ and therefore does not depend
on the choice of the integral form $\hat{P}_b$ of $P_b$.

The simplest link then is obtained,
by adding any number of arrowheads 
colored by $b^1,\cdots, b^\kappa$, then $P_b/P_b(q^{-\rho_k})$
must be replaced here by the product 
\begin{align}\label{simp-link-P}
\hat{\P}_{\ell}^1=
\hat{P}_{\ell+1}^1\cdots \hat{P}_{\ell+1}^\kappa
/\hat{P}_{b^{j_o}}(q^{-\rho_k})=
\hat{P}_{b^1}\cdots \hat{P}_{b^\kappa}
/\hat{P}_{b^{j_o}}(q^{-\rho_k}).
\end{align}
This particular case is already quite interesting; 
expressing the products $P_{b^1}\cdots P_{b^\kappa}$
as linear combinations of Macdonald polynomials 
(generalizing the Pieri rules) is a challenge.
\smallskip

\subsection{\bf The polynomiality}
\begin{theorem}\label{THM-integr-Jones}
For any choice of the normalization
index $1\le j_o\le \kappa$, the DAHA-Jones polynomial \,
$\hat{J\!D}^{j_o}_{\l,\,'\!\l}\,$ defined above
is a polynomial in terms of $q,t_{\sht},t_{\lng}$
up to a factor $q^\bullet t^\bullet=
q^\bullet t_{\sht}^\bullet t_{\lng}^\bullet$,
where the powers $\bullet$ can be rational in the latter. 
It does not depend
on the particular choice of the lifts \,$\ga^j_i\in GL_2(\Z)$
and $\hat{\ga}^j_i\in GL_{\,2}^{\wedge}(\Z)$\, for $1\le  i
\le\ell^j$.
Also one can replace in formula (\ref{jones-hat}) the
ratio $\hat{P}_{b^{j_o}}/\hat{P}_{b^{j_o}}(q^{-\rho_k})$
by $E_{b^{j_o}}/E_{b^{j_o}}(q^{-\rho_k})$ without changing
the output. 
\end{theorem}

{\it Proof.} 
The justification of this and the next theorem almost exactly
follow those in Theorem 1.2 from \cite{CJJ} (for torus knots)
and in Theorem 2.1 \cite{ChD} (for iterated torus knots).

\smallskip
The quantity $\hat{J\!D}^{\,j_o}_{\,\l,\,'\!\l}$ 
can be non-integral with respect to $q,t_\nu$ only if 
$P_{b}^\circ\equal P_{b}/P_b(q^{-\rho_k})$ 
for $b=b^{j_o}$ has a pole in terms of $q$ when $\ep=0$,
where $\ep=(1-q^j t_{\sht}^r t_{\lng}^s)$ for 
certain $j>0,r,s\ge 0,r+s>0$.
We can assume that the binomial $\ep$ is maximal such, i.e. that 
$P_{b}^\circ$ has no singularity at
$(1-q^{jv} t_{\sht}^{rv}t_{\lng}^{sv})/\ep=0$ for any 
$\Z\ni v>1$.

Let us localize and complete the ring of coefficients of $\HH$
and the polynomial representation $\v$, which is 
$\Z_{q,t}=\Z[q^{\pm 1/(2\mm)},t_\nu^{\pm 1/2}]$, 
by such $\ep$,
i.e. with respect to the ideal $\bigl(1-q^v t_{\sht}^r
t_{\lng}^s\bigr)$; 
the notations 
will be $\,\Z_{q,t}^{(\ep)},\, \HH^{(\ep)},\,
\v^{(\ep)}$. Note that we added $q^{\pm 1/(2\mm)}$ to
$\Z_{q,t}$.
\smallskip

We will use
the {\em evaluation pairing}, defined as follows:
$$
\{E,F\}_{ev}\,=\, E(Y^{-1})(F(X))(q^{-\rho_k}),
\ E,F \in \v^{(\ep)}\,.
$$
See Theorem 11.8 from \cite{C103} and
Theorem 1.2 from \cite{CJJ}. We set
$Rad_{\ep,p}=\bigl\{F\in \v^{(\ep)}\,\mid\, 
\{F,\v\}_{ev}\in 
\ep^p\Z_{q,t}^{(\ep)}\bigr\}
\for p\in \Z_+\,. 
$

\smallskip
Switching from $\v$ to the whole $\HH$, we define  
\begin{align}\label{RADep}
R\!A\!D_{\ep,p}\equal\bigl\{ H\in \HH^{(\ep)}\, \mid\, 
\{\HH\, H\, \HH\}_{ev}\in \ep^p\Z_{q,t}^{(\ep)}\bigr\}
\for p\in \Z_+.
\end{align}
Equivalently, 
$
R\!A\!D_{\ep,p}=\bigr\{H\in \HH^{(\ep)} \mid 
H(\v^{(\ep)})\subset Rad_{\ep,p}\bigl\},
$ 
since $Rad_{\ep,p}\!=\!\bigl\{F\!\in\! \v^{(\ep)} 
\mid \{\HH(F)\}_{ev}\!\in 
\ep^\ell\Z_{q,t}^{(\ep)}\bigr\}$; see Lemma 11.3
from \cite{C103}.

Here $q$ is not a root of unity. Therefore any $Y$\~invariant
submodule of $\v$ is invariant with respect to the
natural action of $\dot{\tau}_-$ of $\tau_-$ in $\v$. 
We conclude that 
$\psi$ and $\tau_-$ preserve $R\!A\!D_{\ep,p}$ for any 
$p\in \Z_+$ (and for generic $q$). The same holds for $\eta$;
see (\ref{etatxpi}). Thus the whole  
$GL_{\,2}^{\wedge}(\Z)$ preserves each $R\!A\!D_{\ep,r}$. 

Let $\bar{P}_{b}\equal\ep^l 
P_{b}^\circ \in \v^{(\ep)}$ for $b=b^{j_o}$ and minimal
such $l\in \N$. Then 
$\bar{P}_{b}(q^{-\rho_k})\in 
\ep^l\Z_{q,t}^{(\ep)}$ due to the
normalization of $P_b^\circ$.  
Since $\bar{P}_{b}$  is an eigenfunction 
of $\{L_f\}$ from (\ref{macdopers}), one has that
$\bar{P}_{b} \in Rad_{\ep,l}$; see 
Lemmas 11.4-5 from \cite{C103}.

This implies that
$\hat{\ga}(\bar{P}_{b_{j_o}})\in R\!A\!D_{\ep,l}$ for any 
$\ga\in PSL_2(\Z)$, as well as for any
$\bar{P}=R\bar{P}_{b^{j_o}}$, where $R$ is {\em any\,} 
$q,t$\~integral Laurent polynomial.

The projection  $\bar{P}'=\hat{\ga}(\bar{P})\!\Downarrow=
\hat{\ga}(\bar{P})(1)$ then belongs to $Rad_{\ep,l}$.
Hence $\hat{\ga}'(\bar{P}')\in R\!A\!D_{\ep,l}$ for
any $\ga'\in  PSL_2(\Z)$ and we can continue this process.
One automatically obtains that
$\{\bar{P}''\}_{ev}$ is divisible by $\ep^l$ for 
$\bar{P}''\equal\hat{\ga}'(\bar{P}')\!\Downarrow$ 
and for any further such polynomials obtained by this procedure
continued for  $\ga'',\ga''',\ldots$ from  $PSL_2(\Z)$
or from $GL_2(\Z)$. 

Here we can multiply each
$\bar{P}',\bar{P}''$ by arbitrary $q,t$\~integral
$R',R''$ and so on. Furthermore, we can apply to the
{\em final\,} $\bar{P}$\~polynomial in this chain
any operator in the form $Q(Y^{-1})$ for
a $q,t$\~integral Laurent polynomial $Q(X)$ and the
output will have the same
divisibility by $\ep^l$ as above. 
Thus we conclude that 
$\hat{J\!D}^{j_o}_{\l,\,'\!\l}$ has actually
no singularity at $\ep=0$, which contradiction is sufficient to 
claim its $q,t$\~integrality. 
\smallskip

The independence of $\hat{J\!D}^{j_o}_{\l,\,'\!\l}$
of the lifts of the columns $(\rr_i^j,\ss_i^j)^{tr}$
to $\hat{\ga}_i^j$ and the possibility to replace
$P_{b^{j_o}}$ by $E_{b^{j_o}}$  
closely follow the corresponding claims
in Theorem 1.2 from \cite{CJJ}. \sq
\medskip

\subsection{\bf Major symmetries}
The following theorem is a link counterpart of Theorem 2.2 from
\cite{ChD}. We will comment on its proof (but omit the
details); it remains
essentially the same as for iterated {\em knots}.
All following claims 
hold when  $\,'\!\l$ is replaced by $\,'\!\l^\vee.$ 
Furthermore, formula
(\ref{eval-tau}) from Theorem \ref{EVALTAUM} can be used
to reduce the pairs of graphs to a single tree;
see, e.g. an example after (\ref{T2-3--1-2}). 

\begin{theorem} \label{MAINTHM} 
(i) {\sf Minimal normalization.}
The $q,t$\~integrality and other claims
from Theorem \ref{THM-integr-Jones} hold for the
following modifications of DAHA-Jones polynomials
(which does not require picking $j_o$):
\begin{align}\label{jones-bar}
\hat{J\!D}^{min}_{\,(\vec\rr^j,\,\vec\ss^j)
\, , (\,'\!\vec\rr^j,\,'\!\vec\ss^j)}
(\mathbf{b},\mathbf{'b}\, ; \,q,t)
=\,\hat{J\!D}^{min}_{\l,\,'\!\l}\equal 
\Bigl\{\frac{\, \vph(\,'\hat{\P}_0^{tot})\,
\hat{\P}_0^{tot}}{
(\mathbf{b},\mathbf{'b})_{ev}^\wedge}\Bigr\}_{ev},
\end{align}
where  the polynomials $\hat{P}_b (b\in P_+)$ 
are defined for one of the following
integral forms: 
\begin{align}\label{int-forms}
\hat{P}_b=
\tilde{P}_b,\, \ddot{P}_b,\ J_{\la(b)} (\hbox{for\,} A_n),
\hbox{\, or\, }
P^{(\!1\!)}_b (\hbox{when \, } k_{\sht}\!=\!1\!=\!k_{\lng});
\end{align}
$(b^1,\ldots,b^m)_{ev}^\wedge$ are defined correspondingly.
See (\ref{lcmkrho}),(\ref{P-arms-legs-union}),
(\ref{P-arms-legs-k1}). 
\smallskip

(ii) {\sf Topological symmetries.}
Following Section \ref{sec:knots-links},
the polynomial $\hat{J\!D}^{j_o}_{\l,\,'\!\l}$ from
Theorem \ref{THM-integr-Jones}
and $\hat{J\!D}^{min}_{\l,\,'\!\l}$ from $(i)$, 
considered up to a monomial factor 
$q^\bullet t_{\sht}^\bullet
t_{\lng}^\bullet$,  depends 
only on the topological link corresponding to the pair of graphs
$\{\l,\,'\!\l\}$.
For instance, (\ref{iter-triv}) holds and the reduction of the
vertices with $\rr=1$ from (\ref{iter-s1}) can be applied in
$\l$ or in $'\!\l$. Also, the transposition  
$[\rr_1^j,\ss_i^j]\mapsto $ $[\ss_1^j,\rr_i^j]$ from 
(\ref{iter-sym}) does not influence 
$\hat{J\!D}^{j_\circ}_{\l}$ or $\hat{J\!D}^{min}_{\l}$ 
if $\,'\!\l=\emptyset$.
 
Moreover, $\hat{J\!D}^{j_o}$ or  $\hat{J\!D}^{min}$
above become $\star-$conjugated 
if all $\ss^j_i,'\!\ss^j_i$ from $\l$ and $'\!\l$ change
their signs simultaneously; see (\ref{iterorien}), 
(\ref{iter2min}). Furthermore,
they do not change when $\,b^j\mapsto \iota^{\!\pi}(b^j)\,$ and
simultaneously 
\begin{align}\label{jones-sym-iter} 
&\{\rr^j_i\},\{\ss^j_i\}\mapsto
\{\ldots,\rr^j_{i\!-\!1},
-\rr^j_i,\rr^j_{i\!+\!1},\ldots\},
\{\ldots,\ss^j_{i\!-\!1},-\ss^j_i,\ss^j_{i\!+\!1},\ldots\}
\end{align}
in $\l$ and $\,'\!\l$ provided that total numbers of such 
transformations in every path of $\l,\,'\!\l\,$  have all the 
same parity $\pi\in \Z_2$.
\smallskip

(iii) {\sf Color exchange.} 
We assume that for given $\mathbf{b},'\mathbf{b}$,
generic $q$ and certain (special)
$t_{\sht},t_{\lng}$, there exist 
$\mathbf{c}=(c^j),\,'\mathbf{c}=('\!c^j)$ and
$w^j,\,'\!w^j\in W$ satisfying the relations   
\begin{align}\label{bcw-rel}
q^{\,(b+\rho_k-w(\rho_k)-w(c)\,,\,\al)\,}=1=
q_{\al}^{\,(b-w(c)\,,\,\al^{\!\vee})}\,
t_{\sht}{\!\!\!}^{(\rho^w_{\sht}\,,\,\al^{\!\vee})}\,
t_{\lng}{\!\!\!}^{\,(\rho^w_{\lng}\,,\,\al^{\!\vee})},
\end{align}
for all $\al\in R_+$,
where\, $w\!=\!w^j,b\!=\!b^j,c\!=\!c^j$\, or\, 
$w\!=\! '\!w, b\!=\! 'b^{j},c\!=\!'\!c^{j}$, and  
we set\, $\rho_{\nu}^w\equal w(\rho_{\nu})-\rho_{\nu}$. Then 
\begin{align}\label{jones-bc}
\hat{J\!D}^{min}_{\l,\,'\!\l}=\hat{J\!D}^{min}_{\m,\,'\!\!\m}
\hbox{\, for such \, } q,\{t_\nu\}, 
\end{align}
where $\m,\,'\!\m$ correspond to the same $(\rr_i^j,\ss_i^j)$,
$('\rr_i^j,'\!\ss_i^j)$, $\Up,\,'\Up$ but with the colors
$b^j,\,'\!b^j$ replaced by $c^j,\,'\!c^j$.
\smallskip

(iv) {\sf Specialization $q=1$.}
We now make $q=1$, assuming that $t_\nu$ are generic and
using the notation $(b^1,\ldots,b^m)_{ev}^\wedge$ from
(\ref{lcmkrho}). We switch here to the spherical
polynomials $P_b^\circ$, so the following does not depend 
on the choice of the integral forms $\hat{P}_b$: 
\begin{align}\label{q-1-prod}
&\frac{(b^1,\ldots,b^\kappa,\,'b^1,\ldots,
'\!b^{\,'\kappa})^\wedge_{ev}}
{(b^1)^\wedge_{ev}\cdots (b^\kappa)^\wedge_{ev}
\,('b^1)^\wedge_{ev}\cdots
('b^{\,'\kappa})^\wedge_{ev}}(q\!=\!1)\,\,
\hat{J\!D}^{min}_{\l,\,'\!\l}(q\!=\!1)\\
=\ \,&
\prod_{j=1}^{\kappa} J\!D_{\,\vec\rr^j,\,\vec\ss^j}
\,\bigl(\,b^j\,;\,q\!=\!1,\,t_\nu)\ 
\prod_{j=1}^{'\kappa} J\!D_{\,'\vec\rr^j,\,'\vec\ss^j}
\,\bigl(\,'\!b^j\,;\,q\!=\!1,\,t_\nu),\notag\\
&\hbox{where\,\,}
J\!D_{\,\vec\rr,\,\vec\ss}
\,\bigl(b;\,q\!=\!1,t_\nu
\bigr)\!=\!
\hbox{\small$\prod$}_{p=1}^n J\!D_{\,\vec\rr,\,\vec\ss}\,
(\om_p;\,q\!=\!1,t_\nu)^{b_p},\notag
\end{align}
for $b=\hbox{\small$\sum$}_{p=1}^n b_p \om_p\in P_+$,
where the knot $J\!D$\~polynomials from  
(\ref{jones-ditx}) are used. See formula (2.18) 
in \cite{ChD}. These specializations coincide for the
pairs $\{\l,\,'\!\l\}$ and $\{\l,\,'\!\l^\vee\}$.
\sq
 \end{theorem}
\smallskip

{\sf Comments.}
Concerning using ``min" in $\hat{J\!D}{}^{min}$ and 
for the superpolynomials later, they
can be reducible depending on the choice of the integral
form. We expect them to be irreducible for $\hat{P}_b=
\ddot{P}_b$, but this is not generally true for other
integral forms, including $\hat{P}_b=\tilde{P}_b$ and
$\hat{P}_b=J_b$ in the $A$\~case (which are not 
the smallest possible). 
\smallskip

The justification of the symmetries from Part $(ii)$ is
essentially parallel to Theorem 1.2 from
\cite{CJJ}. Let us comment on (\ref{iter-s1}).
Essentially, one needs to check here that the torus
knot  $T(m\rr+\ss,\rr)$ results in the same
DAHA-Jones polynomial as the ``2-cable" corresponding to
the Newton pairs $[m,1]\rightarrow [\rr,\ss]$. Topologically, 
$T(m\rr+\ss,\rr)$ is isotopic to 
$C\!ab(m\rr+\ss,\rr)T(m,1)$, since
$T(m,1)$ is unknot. As it was noted above, the corresponding
relation for the $J\!D$\~polynomials readily 
follows from the commutativity $\tau_-^m$
with $\Downarrow$, which simply means that $\tau_-$ acts 
in $\v$. Other symmetries are
based on applying $\eta, \iota, \si^2$ inside
the coinvariant; they are compatible with the 
projection $\Downarrow$ as well. Also, the case of the pair of
graphs is governed by Theorem \ref{EVALTAUM}.
\smallskip

The last part of the theorem states that the
DAHA-Jones polynomials for iterated torus
links become 
the products of those for the paths under $q=1$ 
(i.e. over the knots that constitute the link). 
This is compatible with our conjecture in \cite{ChD} 
concerning the Betti
numbers of the Jacobian factors for (the germs of)
unibranch planar singularities. The DAHA-superpolynomials
in type $A$ under $a=0,\,q=1$ occur here. Thus the passage
from knots to links does not add much to our conjecture
on the Betti numbers due to $(iii)$; 
we omit the discussion of the Jacobian factors and their
Betti numbers in this work. 
\smallskip

{\sf QG invariants.}
Extending the connection conjectures from \cite{CJ,ChD}, we
expect that for $k_{\sht}=1=k_{\lng}$ and the
integral form $\hat{P}_b=P_b^{(\!1\!)}$,
\begin{align}\label{jones-hatx}
\hat{J\!D}^{j_o}_{\l,\,'\!\l}(t_\nu\mapsto q_\nu)=
\hat{J\!D}_{\,(\vec\rr^j,\,\vec\ss^j)\,,
\,('\!\vec\rr^j,\,'\!\vec\ss^j)}
^{R,\Up,j_o}((b^j),('\!b^j)\, ; \,q,t_\nu\mapsto q_\nu)
\end{align}
coincide up to $q^\bullet$ 
with the {\em reduced Quantum Group (WRT) invariants\,} of the
corresponding  iterated torus colored links, where the
reduced normalization is for the weight $b^{j_o}$.
To obtain the {\em non-reduced\,} $QG$-invariants, 
one takes $j_o=\emptyset$, i.e.
omits the division by $P^{(\!1\!)}_{b^{j_o}}
(q^{-\rho})$
in $\hat{J\!D}^{j_o}_{\l,\,'\!\l}$ in this case
(for $k_\nu=1$).
 
We do not see at the moment any
topological meaning of the 
$\hat{J\!D}{}^{min}$\~polynomials from part $(i)$ in
Theorem \ref{MAINTHM} for such $k_\nu$ and
general $(b^j),('\!b^j)$. The division there 
is by $(\mathbf{b},\,'\mathbf{b})_{ev}^{(\!1\!)}$,
which provides the $q$\~integrality, but results
in the polynomials that seem too small for the
topological interpretation. This is unless all 
colors are the same.
\smallskip

We note that the Quantum Group (WRT) invariants
are associated with the root 
system $\tR$; see \cite{CJ}. The {\em shift operator\,} was used 
there to deduce this coincidence from \cite{Ste} in the case
of $A_n$ and torus knots; quite a few confirmations were provided 
for other root systems (including special ones).
The method was originated in  {\em CFT\,} Conformal Field Theory,
and {\em Verlinde 
algebras}. DAHA provides perfect tools to understand and
generalize the latter algebras.

As we discussed above, Part $(i)$ of 
Conjecture \ref{CONCONJ} can be checked using the Skein;
another approach, which is applicable to verifying
the connections with the $QG$\~invariants for any root systems,
is a combination of the DAHA
shift operators  with the {\em Rosso-Jones formula\,}. 

\section{\sc DAHA-superpolynomials}
\subsection{\bf Main theorem}
Following \cite{CJ,GoN,CJJ,ChD},
Theorem \ref{JONITER} and Theorem \ref{MAINTHM} can
be extended to the {\em DAHA- superpolynomials\,},
the result of the stabilization  of 
$\hat{J\!D}^{A_n,j_o}_{\l,\,'\!\l}$ 
(including $j_o=\emptyset$) 
or $\hat{J\!D}^{A_n,min}_{\l,\,'\!\l}$ 
with respect to $n\to \infty$. 

This stabilization was announced 
in \cite{CJ} for torus knots; its proof was published in 
\cite{GoN}. 
Both approaches use \cite{SV} and can be 
extended to arbitrary iterated knots and links;
the Duality Conjecture was proposed in \cite{CJ}
and proven in \cite{GoN} for torus knots; also see 
\cite{CJJ} for an alternative approach based on the
generalized level-rank duality. The justifications
of these claims for iterated torus knots \cite{ChD} and
iterated torus links is essentially parallel to the case
of torus knots, though there are some deviations. 
\smallskip

The main change here from knots to links is that the 
polynomiality of the superpolynomials for {\em links\,} is 
based on our
using the $J_\la$\~polynomials as the integral form of 
$\{P_\la\}$. 
Actually $\{J_\la\}$ were already employed in \cite{GoN}
for the stabilization and duality, but the construction
of (reduced) $J\!D$\~polynomials and superpolynomials for 
{\em knots\,} requires only spherical $\{P^\circ_\la\}$ and 
does not depend on picking their integral forms. For 
{\em links\,} vs. knots, the role of polynomials
$\{J_\la\}$ becomes the key; without using the $J$\~polynomials
the superpolynomials can become rational (not polynomial) in
terms of $t$.

\smallskip

The sequences $\vec\rr^j,\,\vec\ss^j$ of length $\ell^j$
for the graph $\l$ and $'\!\vec\rr^j,'\!\vec\ss^j$ of 
length $'\!\ell^j$ for the graph $'\!\l$ will be from the
previous sections. We consider only $A_n$ here,
setting $t=t_{\lng}=t_{\sht}=q^k$. We will always
use below the 
DAHA-Jones polynomials $\hat{J\!D}^{j_o},\, \hat{J\!D}^{min}$
in type $A_n$ (for $sl_{n+1}$)\, 
{\em defined in terms of $J_{\la}$\,}, unless
stated otherwise. The integral form $P_{\la}^{(\!1\!)}$ for $t=q$ 
(i.e. when $k=1$) will be needed only when discussing the 
connection
with the HOMFLY-PT polynomials.
Recall that $\la=\la(b)$ is the Young diagram
representing $b\in P_+$.
\smallskip

We will consider $P_+\ni b=$ $\sum_{i=1}^n b_i \om_i$ 
for $A_n $ as 
a (dominant) weight for {\em any\,} $A_m$ (for $sl_{m+1}$)
with $m\ge n-1$,
where we set $\om_{n}=0$ upon its restriction to $A_{n-1}$.
The integral form of $P_b$ in (\ref{jones-hat})
and (\ref{jones-bar}) will be $J_{\la}$ from 
(\ref{P-arms-legs})
for $\la=\la(b)$ in the next theorem. See 
\cite{GS,CJ,GoN,CJJ,ChD} concerning the version
of the following theorem for torus knots and
torus iterated {\em links}.

\begin{theorem}\label{STABILIZ}
(i) {\sf Stabilization.} Given the links $\l,\, '\!\l$ 
colored by $\mathbf{b}=(b^j),\, '\mathbf{b}=('\!b^j)$ and
the normalization index $1\le j_o\le \kappa$ (including
$j_o=\emptyset$), there exist
polynomials from $\Z[q,t^{\pm 1},a]$
\begin{align}\label{h-polynoms-hat}
&\hat{\h}^{j_o}_{\l,\,'\!\l}\ =\
\hat{\h}^{\,\Up,j_o}_{(\vec\rr^j,\,\vec\ss^j),
('\vec\rr^j,\,'\vec\ss^j)}(\mathbf{b},\,
'\mathbf{b};\,q,t,a),\\
&\hat{\h}^{min}_{\l,\ '\!\l}\, =\  
\hat{\h}^{\,\Up,min}_{(\vec\rr^j,\,\vec\ss^j),
('\vec\rr^j,\,'\vec\ss^j)\,}(\mathbf{b},\,
'\mathbf{b};\,q,t,a)\notag
\end{align}
such that for any  $m\!\ge\! n\!-\!1$ and proper powers
of $q,t$ (possibly rational)\,:
\begin{align}\label{jones-sup-hat}
&\hat{\h}^{j_o}_{\l,\,'\!\l}(q,t,a\!=\!-t^{m+1})\,=\,
\pm\, q^\bullet t^\bullet\, 
\hat{J\!D}^{A_m,j_o}_{\l,\,'\!\l}(q,t),
\\ 
&\hat{\h}^{min}_{\l,\,'\!\l}(q,t,a\!=\!-t^{m+1})\!=\!
\pm q^\bullet t^\bullet\, 
\hat{J\!D}^{A_m,min}_{\l,\,'\!\l}(q,\,t).\notag
\end{align} 
They are normalized as follows ({\sf the hat-normalization}).
 We multiply $\hat{\h}$
by $\pm q^\bullet t^\bullet$ to ensure that $\hat{\h}(a\!=\!0)$
is not divisible by 
$\,q\,$ and $\,t\,$ and that the coefficient of the 
minimal 
power of $t$ in $\hat{\h}(a\!=\!0)$ is positive. If 
$\hat{J\!D}^{A_m}$\~polynomials are considered under the
same hat-normalization, then 
relations (\ref{jones-sup-hat}) will automatically 
hold for sufficiently large $\,m$
without any correction factors $\pm q^\bullet t^\bullet$
(and one sufficiently large
$\,m$\, is actually sufficient to fix $\,\hat{\h}\,$ 
uniquely). The $a$\~stabilization and other claims here
hold for
$\hat{\h}_{\l,\,'\!\l^\vee}$; see (\ref{jones-hat-vee}).

(ii) {\sf Symmetries.}  
The polynomials $\hat{\h}$ depend only on the isotopy
class of the corresponding torus iterated links. For instance, 
the symmetries from Part $(ii)$ of Theorem \ref{MAINTHM}
hold for $\hat{\h}=\hat{\h}^{j_o}, \hat{\h}^{min}$ under the 
normalization from $(i)$ with the following reservation about 
(\ref{jones-sym-iter}); its $a$\~extension holds
upon $a^\star=a^{-1}$ 
(up to $\,a^\bullet q^\bullet t^\bullet).$ 
The same claims are for
$\,'\!\l^\vee$ instead of  $\,'\!\l$, and the
$a$\~extension of Theorem \ref{EVALTAUM} 
in the case of the pair of graphs can be reduced to  
a single graph.

For the {\sf color exchange\,}, we impose (\ref{bcw-rel}), 
and consider $\,w^j,\,'\!w^j\,$ 
as elements of $\mathbf{S}_{m+1}$ (the Weyl group for $A_m$)
for every
$m\ge n$. Then for $\,'\!\l^\diamond$
that is either $\,'\!\l$ 
or $\,'\!\l^\vee$ (for $\diamond=\emptyset,\vee$):
\begin{align}\label{jones-bca}
\hat{\h}_{\l,\,'\!\l^{\diamond}}
(\mathbf{b},\,'\mathbf{b}\,;\,q,t,a) = 
\hat{\h}_{\l,\,'\!\l^{\diamond}}
(\mathbf{c},\,'\mathbf{c}\,;\,q,t,a) = 
\hbox{\, for such \, } q,t.
\end{align}
Similarly, the {\sf specialization\, }
relations from (\ref{q-1-prod}) at $q=1$ in the case of $A_n$
hold for $\hat{\h}$ and such specializations coincide for
$\{\l,\,'\!\l\}$ and $\{\,'\!\l^\vee\}$. Recall that the product
formula there holds when spherical polynomials $P_b^\circ=
P_b/P_b(t^{\rho})$ are used in the formulas for $\hat{J\!D},
\hat{\h}$.

(iii) {\sf Super-duality and deg${}_a$}.
We switch from $\mathbf{b},\,'\mathbf{b}$ to the
corresponding sets of Young diagrams
$\lla$,\, $'\!\lla$. Let 
$\lla^{tr}$,\, $'\!\lla^{tr}$
be their transpositions. Up to powers of $q$ and $t$
denoted here an below by $q^\bullet t^\bullet$, one has:
\begin{align}\label{iter-duality}
\hat{\h}_{\l,\,'\!\l}(\lla,\,'\!\lla\,;\, q,t,a)=
q^{\bullet}t^{\bullet}
\hat{\h}_{\l,\,'\!\l}(\lla,\,'\!\lla\,;\,t^{-1},q^{-1},a);
\end{align}
the same super-duality holds with $\vee$,
i.e. for $\hat{\h}_{\l,\,'\!\l^\vee}$.

\comment{
We will assume that $|\rr^j_1|>|\ss^j_1|$,
 $|'\rr^j_1|>|'\ss^j_1|$ for any
$j$ such that there is only one arrowhead through 
the final vertex $[\rr^j_1,\ss^j_1]$ or
$\{'\rr^j_1,\,'\ss^j_1\}$ in this path. Then 
deg${}_a \hat{\h}^{min}_{\l,\,'\!\l^\vee}
(\,\lla,\,'\!\lla\,;\,q,t,a)\,$
(note $\vee$ here) is no greater than 
\begin{align}\label{deg-a-j}
\sum_{j=1}^\kappa
\max\{1,|\,\ss^j_1\rr^j_2\cdots\rr^j_{\ell^j}|\}\,|\la^j|\,+\,
\sum_{j=1}^{'\!\kappa}
\max\{1,|\,'\ss^j_1\,'\rr^j_2\cdots\,'\rr^j_{\,'\!\ell^j}|\}
\,|'\!\la^j|\, -\,\De,&\notag\\
\De\!=\!|\la^1\!\vee\!\ldots\!\vee\!\la^\kappa\!\vee\,
'\!\la^1\!\vee\!\ldots\!\vee\, '\!\la^{'\!\kappa}|
 \hbox{\, for \,} \hat{\h}^{min},\ 
\De\!=\!|\la^{j_o}| \hbox{\, for \,} \hat{\h}^{j_o},&
\end{align}
where $|\la|$ is the number of boxes in $\la$ 
and the products $(|\cdots|)$ in (\ref{deg-a-j})
equals $1$ if there are no vertices in the  $j${\tiny th}
path (only the arrowhead). Without $\vee$,
the estimate for deg${}_a \hat{\h}^{min}_{\l,\,'\!\l}$
is with $\rr^j_1,\,'\rr^j_1$ replacing
$\ss^j_1,\,'\ss^j_1$.
}

Let us assume that $\rr_i^j,\,'\!\rr_i^j\neq 0$ for $i\!>\! 1$.
Then \, 
deg${}_a \hat{\h}^{min}_{\l,\,'\!\l}
(\,\lla,\,'\!\lla\,;\,q,t,a)\,$  and \, 
deg${}_a \hat{\h}^{min}_{\l,\,'\!\l^\vee}
(\,\lla,\,'\!\lla\,;\,q,t,a)\,$ are no greater than
\vskip -0.7cm
\begin{align}\label{deg-a-j}
\sum_{j=1}^\kappa
\max\{1,|\,\rr^j_1|\}|\rr^j_2\cdots\rr^j_{\ell^j}|\,|\la^j|+
\sum_{j=1}^{'\!\kappa}
\max\{1,|\,'\rr^j_1|\}|\,'\rr^j_2\cdots\,'\rr^j_{\,'\!\ell^j}|
\,|'\!\la^j| -\De,&\notag\\
\De\!=\!|\la^1\!\vee\!\ldots\!\vee\!\la^\kappa\!\vee\,
'\!\la^1\!\vee\!\ldots\!\vee\, '\!\la^{'\!\kappa}|
 \hbox{\, for \,} \hat{\h}^{min},\ 
\De\!=\!|\la^{j_o}| \hbox{\, for \,} \hat{\h}^{j_o},&
\end{align}
where $|\la|$ is the number of boxes in $\la$.
This estimate holds for $\hat{\h}^\emptyset$ (when $\De=0$), 
i.e. without divisions by the evaluations at $t^\rho$.\sq
\end{theorem} 
\smallskip

Frequently for $\{\l,\,'\!\l^\vee\}$ with $\vee$
(conjecturally for {\em all\,} positive pairs) 
the following improvement of (\ref{deg-a-j})
holds (with the same $\De$):
\begin{align}\label{deg-a-jj}
&\hbox{deg}_a\hat{\h}^{min}_{\l,\,'\!\l^\vee}
=\hbox{\small $\sum$}_{j=1}^\kappa
\max\bigl\{1,\min\{|\rr^j_1|,|\ss^j_1|\}\, |\rr^j_2|\cdots
|\rr^j_{\ell^j}|\bigr\}\,|\la^j|\!\\
&+
\hbox{\small $\sum$}_{j=1}^{\,'\!\kappa}
\max\bigl\{1,\min\{|\,'\rr^j_1|,|\,'\ss^j_1|\}\, 
|\,'\rr^j_2|\cdots
|\,'\rr^j_{\,'\!\ell^j}|\bigr\} \,|\,'\!\la^j|+\De.\notag
\end{align}

\subsection{\bf Sketch of the proof}
Generally, the stabilization and duality are due to the 
switch from $P_{\la}$ to the {\em modified Macdonald 
polynomials\,} \cite{GoN}.  The 
projective action of $PSL_2(\Z)$ is compatible with such a 
switch (see there and \cite{SV}). The transition to the 
iterated links from torus knots is relatively straightforward;
only one $b$ will be discussed below.

It suffices to check the stabilization (up to $\pm q^\bullet
t^\bullet$) for $\hat{\h}^{min}$ and for the integral form
$J_{\la(b)}$. Then the stabilization of $\hat{\h}^{j_o}$ follows
from Proposition \ref{PROP:stab-values}. Note that
the stabilization holds for
$\hat{P}_b=P_b$, for instance for the
corresponding $\hat{\h}^\emptyset$ (i.e.
without the division by the evaluation at $t^\rho$). The
polynomials $P_b$ generally are not $q,t$\~integral, so
the corresponding $\hat{\h}^\emptyset$
is generally rational. We note that even if $t=q$,
the denominators in terms of $q$ do occur in $\hat{\h}^\emptyset$ 
for links.

The limit $\hat{\h}(a=0)$, used in the hat-normalization
condition,  
is the term-wise intersection (common part) of
$\hat{J\!D}{}^{A_m}$ for $m\gg 0$, under the normalization 
as in $(i)$. Indeed, the 
terms in $\hat{\h}(a=0)$ with deg$_{\,t}\le M$
can be extracted from a {\em single\,} $\,J\!D^{A_m}\,$
if $m\ge M$, since $a=-t^{m+1}$ is beyond
this range of $t$. This explains why the correction factors 
$\pm q^\bullet t^\bullet$ are not necessary for large $m$
for the hat-normalization.

Concerning the super-duality and deg${}_a\hat{\h}^{min}$, 
the simplest (and actually the key) particular case is when
there are no vertices at all in $\l$ and this graph contains
only $\kappa\,$ arrowheads. We set symbolically
$\Up=\,\{\rightdotsarrowtiny\!\}$ (any number of arrows).
Let us omit $'\!\l$; replace
$\lla=(\la^j)$ below by $\lla\cup\,'\!\lla\,$ if $\,'\!\l$
is present. Using 
(\ref{stabeval})
and taking into account the normalization from $(i)$,
\begin{align}\label{stabevaly}
&t^{(\rho,b)}J_\la(t^{-\rho})=\Pi_\la^\dag\equal
\prod_{p=1}^n \prod_{v=0}^{\la_{p\,}-1} 
\Bigl( 1\,+\,q^{v}\, a\, t^{-p+1}\Bigr),\\
&\hat{\h}^{\rightdotsarrowtiny\!,\,min}_{\l}(a,q,t)=
\Bigl(\Pi_{\la^1}^\dag\cdots\Pi_{\la^\kappa}^\dag\Bigr)/
\Pi_{\la^1\vee\ldots\vee\la^\kappa}^\dag.
\end{align}
The division here will be by $\Pi_{\la^{j_o}}^\dag$ for
$\hat{\h}^{j_o}$ (no division when $j_o=\emptyset$).

Then (\ref{stabevaly}) gives the
required formula for deg$_a\hat{\h}$ in this case.
The super-duality follows here from the relation
$$
\Pi_\la^\dag(a,q,t)=\Pi_{\la^{tr}}^\dag(a,t^{-1},q^{-1}).
$$

Let us briefly consider the
case of $\{\l,\,'\!\l\,\}$ for
$\,'\!\l=\rightdotsarrowtiny\,'\!\la=\mu$. Then
$\hat{\h}^{\emptyset}_{\l,\,'\!\l}$ is
essentially $J_{\mu}(q^{k\rho+\la})$. Indeed,
(\ref{ebdual}) and (\ref{macdopers}) give that
$f(Y^{-1})(P_b)\!=\!f(q^{-\rho_k-b})P_b$ \,for $b\in P_+$ and
symmetric $f$. Then the $a$\~stabilization is due to
Proposition \ref{PROP:stab-values}. The $a$\~degree will be
the same as that for $J_{\mu}(q^{k\rho})$, which is the same
with $\vee$, i.e. for $\,'\!\l^\vee$ instead of $\,'\!\l$; 
the corresponding coinvariant is $J_{\mu}(q^{-k\rho-\la})$ 
with $\vee$ (up to the proper renormalization).
These two cases, with $\vee$ and without, are
quite different; the $a$\~degree of
$\hat{\h}^{min}_{\l,\,'\!\l}$ is generally greater than
that with $\,\vee\,$. See 
(\ref{T4-2-y-3-2}), (\ref{6-3-T1-1-1}) and
(\ref{4-2-T21-1-1})  vs.
(\ref{T4-2-yv-3-2}), (\ref{6-3-T1-1-1v}) and (\ref{4-2-T21-1-1v}).
\smallskip

Generalizing, one needs to decompose the pre-polynomials 
$\hat{\P}_0^{tot},\,'\hat{\P}_0^{tot}$ from (\ref{jones-hat})
with respect to the Macdonald polynomial and determine the
leading one. We essentially follow \cite{GoN}.
Note that it was conjectured in \cite{ChD} that 
(\ref{deg-a-jj}) gives exact deg$_a$ for {\em algebraic knots\,};
there are examples of non-algebraic knots and links where this
formula for the $a$\~degrees does not work. Without $\,'\!\l$
and in the uncolored case, our formula for deg$_a$ is nothing but 
the {\em multiplicity\,} of the corresponding singularity (Rego). 
Also, the coincidence of $\hat{\h}(q=1)$ with $\vee$ and without 
is not difficult to see; the product formula from
Part $(iii)$ holds because  $q=1$ is the case of the trivial 
{\em center charge\,}.
\smallskip

Continuing the remarks after Theorem \ref{MAINTHM},
we note that the polynomial $\hat{\h}^{min}$ can be
reduced for certain iterated links, so they are not exactly
``minimal" $a,q,t^{\pm1}$\~integral. The importance of $J_\la$ 
(not sharp as integral forms) is mainly due to
the stabilization and super-duality. We omit details
here, but let us provide the following example.
If one makes $t=q$ and takes $\hat{P}_\la=P_\la^{(\!1\!)}$
(=Schur polynomials) as the integral form, then
the corresponding $\hat{\h}^{min}$ will be
{\em not\,} a polynomial in terms of $t\!=\!q$ even in the 
case of $2$ unknots, though the corresponding
$\hat{J\!D}^{min}_\la$ are of course always
$q$\~integral. This is important; the {\em polynomial\,}
$a$\~stabilization does require using ``non-minimal" 
integral forms.   
\smallskip

See \cite{CJJ} for a systematic consideration of the
{\em color exchange\,} condition (\ref{jones-bca})
from $(ii)$. In terms
of the Young diagrams, the procedure is as follows.
Recall that we associate with
$c^j=\sum_{i=1}^n c^j_i\om_i$ the Young diagram
$$
\la(c^j)=\{\la^j_1\!=\!c^j_1+\ldots+c^j_n,
\la^j_2\!=\!c^j_1+\ldots+c^j_{n-1},
\ldots, \la_n^j\!=\!c^j_n,0,0,\ldots\}.
$$
Then we switch to $\la(\bar{c}^j)=\{\bar{\la}^j_i=
\la^j_i-k(i-1)\}$ for $t=q^k$, apply $w^j\in W$ to
$\la(\bar{c}^j)$  and finally obtain 
\begin{align}\label{lamw}
\la(b^j)=\{\bar{\la}^j_{w^j(i)}+k(i-1)\}=
\{\la^j_{w^j(i)}+k(i-w^j(i))\}.
\end{align} 
Here $1\le j\le \kappa\,$ and 
$w^j$ transforms the rows of $\la^j$; we generally
set $w\{\la_1,\la_2,\ldots,\la_n\}=
\{\la_{w(1)},\la_{w(2)},\ldots,\la_{w(n)}\}$.
Given $k<0$ (it can be fractional), $\la(b^j)$ 
{\em must be a Young diagram\,}; this condition 
determines which $\,w^j\,$ can be used. Then we
repeat the same procedure for $'\!\l$ and obtain
$'\mathbf{b}$ from $'\mathbf{c}$ for proper
$'w^j\in W$, $1\le j\le\, '\!\kappa$.
\sq
\smallskip

\subsection{\bf Super-vertex}
Using the definition from (\ref{vert-coinv-p})
for trivial $u=1'$ (which will be omitted here), one has
for any root system $R$ and two sets
$
\ \mathbf{b}=(b^j,\,1\le j\le \kappa)\,\subset P_+ \,\supset\,
'\mathbf{b}=('b^j,\,1\le j\le\, '\!\kappa),
$
\begin{align}\label{vert-coinv-py}
&\p_{\mathbf{b}}^{'\mathbf{b}}\!\equal\!
\bigl\{\tau_+^{-1}\!
\bigl(P_{'\mathbf{b}}^\circ(\!Y^{\!-1}\!)\bigr)
\tau_-^{-1}\!\bigl(P_{\mathbf{b}}^\circ(\!X\!)\bigr)\bigr\}_{ev}
=\bigl\{
\tau_-^{-1}\!\bigl(P_{\mathbf{b}}^\circ
P_{'\mathbf{b}}^\circ\bigr)\bigr\}_{ev}\\
&=\frac{\hat{J\!D}{}^R_{\{1,-1\},\,\{1,-1\}}
(\mathbf{b},'\!\mathbf{b};q,t)}
{(b^1)_{ev}\cdots (b^\kappa)_{ev}
(\,'b^1)_{ev}\cdots (\,'b^{'\!\kappa})_{ev}}=
\frac{\hat{J\!D}{}^R_{\{1,-1\}}
(\mathbf{b}\cup\,'\mathbf{b})}
{(b^1)_{ev}\cdots (\,'b^{'\!\kappa})_{ev}}
\notag
\end{align}
for any integral form $\hat{P}_b$ in $\hat{J\!D}$ (and the 
corresponding $\{\,\}_{ev}$). I.e. $\p_{\mathbf{b}}^{'\mathbf{b}}$
can be interpreted as
$\hat{J\!D}_{\l,\,'\!\l}$ for $\l=\{[\rr_1=1,\ss_1=-1]\}=\,'\!\l$ 
with the arrowheads colored correspondingly by 
$\mathbf{b}$ and $\,'\mathbf{b}$. Note that we use here the
upper indices in $(b^j), ('\!b^j)$ in contrast to the lower ones 
in (\ref{vert-coinv-p}). Due to the division in 
(\ref{vert-coinv-py}), 
the choice of the integral form $\hat{P}_b$ does not really
matter, however the 
corresponding ratio is generally {\em not\,} a polynomial in terms 
of $q,t_\nu$. 
\smallskip

Thus the {\em level-one\,} DAHA-vertex appeared directly 
connected with our invariants for the Hopf links.
The case $'\mathbf{b}=c$ for (a single)
$c\in P_+$ is of particular
interest. Using (\ref{e-single-c-p}),
\begin{align}\label{e-single-c-py}
\p_{\mathbf{b}}^{c}=
C_{\mathbf{b}}^{\iota(c)}
\frac{\langle P^\circ_c,\,P^\circ_c\rangle}
{\langle \theta\mmu^\flat \rangle}
\for \mathbf{b}\subset P_+ \ni c\,,\,\iota(c)=c^\iota=-w_0(c),
\end{align}
where $C_{\mathbf{b}}^{c}\,$ is the 
{\em DAHA multi-vertex\,} for $u=1'$ from
(\ref{P-multi},\ref{lev-1-abc}). It is 
of key importance for $\mathbf{b}=(b^1,b^2)$
and for $\mathbf{b}=b$.
 

\begin{corollary}\label{DAHA-vertex-Hopf}
Let $R=A_n$. For $\hat{P}_\la=J_\la$, where $\la=\la(b)$, 
and for the evaluation
$\,(\la)^\dag_{ev}=t^{(\rho,b)}
J_\la(t^{-\rho})=\Pi_\la^\dag$, we set\,:
\begin{align}\label{h-vertex}
\hat{\h}_{\l,\,'\!\l}^{\emptyset,\dag}
=\frac{\hat{\h}^{\emptyset}_{\{1,-1\},\,\{1,-1\}}
(\lla,\,'\!\lla\,;\,q,t)}
{(\la^1)_{ev}^\dag\cdots (\la^\kappa)_{ev}^\dag
(\,'\!\la^1)_{ev}^\dag\cdots (\,'\!\la^{'\!\kappa})_{ev}^\dag}
\!=\!\frac{\hat{\h}^{\emptyset}_{\{1,-1\}}
(\lla\cup\,'\!\lla)}
{(\la^1)_{ev}^\dag\cdots (\,'\!\la^{'\!\kappa})_{ev}^\dag}.
\end{align}
It is a rational function in terms
of $a,q,t$, where the $a$\~degrees of the numerator and
denominator in (\ref{h-vertex}) coincide (they are 
$a,q,t^{\pm 1}$\~polynomials) and are
equal to $\sum_{j=1}^\kappa\,|\la^j|\,+\,
\sum_{j=1}^{'\!\kappa}\,|'\!\la^j|$. One has:
\begin{align}\label{h-vertexx}
\hat{\h}_{\l,\,'\!\l}^{\emptyset,\dag}= \pm \,q^\bullet t^\bullet 
\p_{\mathbf{b}}^{'\mathbf{b}} 
\for m\ge n-1,
\end{align}
where the extension of the weights $\,\mathbf{b},\,'\mathbf{b}$\, 
to any $\,A_m\, (m\ge n-1)$\, is as in (\ref{vert-coinv-py}) 
via the Young diagrams $\,\la^j=\la(b^j),\, '\!\la^j=\la('b^j)$;
recall that we put
 $\om_n=0$ upon its restriction to $A_{n-1}$. 

In particular, the existence of the $a$\~stabilization of
$\,C_{\mathbf{b}}^{\iota(c)}/\langle \theta\mmu^\flat \rangle=
\p_{\mathbf{b}}^{c}/\langle P^\circ_c,\,P^\circ_c\rangle$
up to $q^\bullet t^\bullet$ (depending on $m$) results
from Proposition \ref{MACD-NORMS-A}. Here
$1/\langle P^\circ_c,\,P^\circ_c\rangle$ 
up to $q^\bullet t^\bullet$ becomes an $a$\~polynomial of 
degree $2|\la(c)|$ with $q,t$\~rational coefficients
(up to $q^\bullet t^\bullet$).
Replacing here $\,'\!\l$ by $\,'\!\l^\vee$, we arrive at the
$a$\~stabilization of 
$\,C_{\mathbf{b}}^{c}/\langle \theta\mmu^\flat \rangle$
(without \,$\iota$), satisfying the associativity
from Part $(ii)$ of Theorem \ref{ASSOC-C}
with $u_i=1'$.
\sq
\end{corollary}

Let us also mention here the relations from
(\ref{h-vertex}): 
\begin{align}\label{rel-Hopf}
\hat{\h}^{\emptyset}_{\{1,-1\},\,\{1,-1\}}
(\lla,\,'\!\lla)=\hat{\h}^{\emptyset}_{\{1,-1\}}
(\lla\cup\,'\!\lla);
\end{align}
they are quite obvious {\em topologically\,}; see below.
We will provide quite a few examples of the
$a$\~stabilization of 
$\,C_{\la^1,\la^2}^{\la^3}/
\langle \theta\mmu^\flat \rangle$ (and that for
$\iota(\la^3)$), but will not discuss the 
``$a$\~associativity" from (\ref{ass-s-4}) in this work.
\medskip

\subsection{\bf HOMFLY-PT polynomials}
There are two approaches to
the unreduced {\em HOMFLY-PT polynomials\,} 
$\hbox{\small H\!O\!M}(\lla;q,\hbox{a})$, via QG (in type $A$) or 
using the corresponding {\em skein relations\,} and 
the corresponding Hecke algebras. Both are 
for any links and colors. See, e.g. \cite{QS} and references 
there. We provide here only a sketchy discussion.

Note that the symmetry 
from Theorem 4.8 in this work corresponds to our $q\mapsto q^{-1}, 
\hbox{a}\mapsto \hbox{a}, \la^j\mapsto (\la^j)^{tr}$. 
The passage to {\em reduced\,} HOMFLY-PT polynomials 
corresponds to our division by $P_{\la^{j_o}}(q^{\rho})$, so it 
requires picking one path-component $j_o$ from 
in a given graph-link. 
Let us impose the hat-normalization here;
the notation will be $\hat{\hbox{\small H\!O\!M}}(q,\hbox{a})$
or $\hat{\hbox{\small H\!O\!M}}_{\l,\,'\!\l}(\lla;q,\hbox{a})$
for the unreduced ones. {\em We will set $\hbox{\rm a}=-a$ below.}

For links, the 
$q$\~polynomiality of the unreduced
HOMFLY-PT polynomials does {\em not\,} 
hold. This is direct from the corresponding normalization of 
colored unknots. 
One has: 
\begin{align}\label{HOM-normal}
\hat{\hbox{\small H\!O\!M}}(q,\hbox{a})=
\Bigl(\frac{1-\hbox{\rm a}}{1-q}\Bigr)^{\kappa} 
\hbox{\,\, for
$\kappa$ uncolored unknots.}
\end{align}
With colors, $\hat{\hbox{\small H\!O\!M}}(q,\hbox{a})=
q^\bullet 
P^{(\!1\!)}_{\la^1}(q^{\rho})\cdots P^{(\!1\!)}_{\la^\kappa}
(q^{\rho})$ for $\kappa$ unknots 
upon the $a$\~stabilization with $\hbox{a}=q^{n+1}$ (for $A_n$).  
Note the absence of minus in the latter, in contrast 
to this substitution for $a$ in $\hat{\h}$. The power 
$q^\bullet$ here is adjusted 
to ensure the hat-normalization of the left-hand side.

From formulas (\ref{P-arms-legs}),(\ref{stabeval}) for
$t=q$ and with such $a$:
\begin{align}\label{stabevalz}
q^{\bullet}P^{(\!1\!)}_\la(q^{\rho})=
&\frac{\prod_{p=1}^n \prod_{j=0}^{\la_{p\,}-1} 
\bigl( 1\,-\,\hbox{a}\, q^{j-p+1}\bigr)}
{\prod_{\Box\in\la}
(1-q^{arm(\Box)+leg(\Box)+1})},\ \la=\{\la_p\},
\end{align}
where $n$ here is the number of (nonempty) rows in $\la$.
Thus the product on the right-hand sides in (\ref{stabevalz})
over the Young diagrams $\la^j$ in $\lla$ is the value of
$\hat{\hbox{\small H\!O\!M}}_{\l}(\lla;q,\hbox{a})$ for 
$\l=\coprod_{j=1}^\kappa \unknot$\,. 

We note that since we deal only with iterated torus links,  
the {\em Rosso-Jones cabling formula\,} is generally
sufficient  for calculating the corresponding
HOMFLY-PT polynomials; see e.g., \cite{RJ,Mo,ChE}. 
This can be actually done for arbitrary root systems
(and WRT-invariants of any colored iterated torus links).
See \cite{AM,AMM} and \cite{Ma} about using here
the HOMFLY-PT skein relations. Also, paper \cite{MoS} 
established the identification of the skein algebra of the
torus with the Elliptic Hall algebra for $t=q$ and therefore
with the corresponding spherical DAHA; see \cite{SV}. 
 
\smallskip

The algebraic-geometric interpretation of these relations
from \cite{Ma} presumably can be generalized to establish the
connection of our superpolynomials to papers \cite{ObS,ORS}.
These two papers are for arbitrary plane curve singularities,
but are restricted to the uncolored case; the main conjecture
from \cite{ObS} was extended in \cite{Ma} to the colored case
and proved; this was {\em unrefined\,}, i.e. without $t$. 

Our previous paper
\cite{ChD} was in the {\em unibranch\,} case; now
we can reach any multi-branch plane curve singularities 
(with arbitrary colors).
We will not discuss systematically the relations to 
\cite{ORS},\cite{Ma} and \cite{Pi} in the present work, though
provide quite a few examples.
\smallskip

The stable {\em Khovanov-Rozansky homology\,} is the
$sl_N$ homology from \cite{KhR1,KhR2} in the range of $N$ 
where the isomorphism in Theorem 1 from \cite{Ras} holds
(see also \cite{Kh}). Thus
they can be obtained from the 
{\em triply-graded HOMFLY-PT homology}, assuming that the
corresponding {\em differentials\,} are known (they are
generally involved).

Considering links and adding colors makes the KhR\~theory 
significantly more difficult. Even without such an
extension, the HOMFLY-PT  homology is known only 
for very limited number of examples and  
no formulas are known for torus iterated knots/links beyond
torus knots. The situation is much better for the {\em Khovanov 
homology\,} (for $sl_2$), though colors-links are 
a problem even in this case. The {\em categorification theory\,} 
can generally address arbitrary 
colors (dominant weights), but the HOMFLY-PT {\em homology\,} 
remains quite a challenge; see 
e.g., \cite{Kh,WW,Rou,Web} and references there. 

Thus we have to restrict Part $(iii)$ of the conjecture below
to the uncolored unreduced case (unless for $N=2$). 
The corresponding Poincar\`e series, {\em stable 
unreduced Khovanov-Rozansky series\,}, will be
denoted by  $K\!h\!R^{\hbox{\tiny stab}}_{\l,\,'\!\l}
(q_{st},t_{st},a_{st})$ in the (topologically)
standard parameters; see \cite{ORS} and below.
The passage to the Khovanov-Rozansky polynomials
for $sl_N$, denoted below by $K\!h\!R^N_{\l,\,'\!\l}$\,,  
is as follows:
\begin{align}\label{KhR-subs}
a_{st}\mapsto q_{st}^{N}, \hbox{\, equivalently, \, }
a\mapsto t^{N}\sqrt{q/t} \hbox{\ \,in our parameters}.
\end{align} 
This actually depends on the grading used in the theory
and there are variations here in different works.

These relations are applied as such only
for sufficiently large (stable) $N$,
otherwise the theory of differentials is necessary.
The differentials correspond to a different
substitution: $a=-t^N$ in DAHA parameters 
(see the conjecture below). 
Let us mention 
the relation to the 
{\em Heegaard-Floer homology\,}: $N=0$. 
Also, the 
{\em Alexander polynomial\,} 
of the corresponding singularity 
is 
$\hat{\h}^{min}_{\l}\,
(q,q,a=-1)/(1-q)^{\kappa-\de_{\kappa,1}}$                   
in the case of one {\em uncolored\,} tree $\l$
with $\kappa$ paths (the number of connected components in the
corresponding cable). This is {\em zeta-monodromy\,} from
\cite{DGPS} upon $t\mapsto q$ (unless for the unknot). 
\smallskip
 
We will always impose the {\em hat-normalization\,} 
from Part $(i)$ of Theorem \ref{STABILIZ} below.
The notation will be 
$\hat{K\!h\!R}{}^{\hbox{\tiny stab}}_{\l,\,'\!\l}$
and $\hat{K\!h\!R}{}^{N}_{\l,\,'\!\l}$. Namely, we divide the first 
polynomial
by the smallest power of $a_{st}$ and then divide (both of) them
by the greatest possible $q_{st}^\bullet t_{st}^\bullet\,$ 
making 
$\,\hat{K\!h\!R}{}^{\hbox{\tiny stab}}_{\l,\,'\!\l}(a_{st}\!=\!0)$
and $\hat{K\!h\!R}{}^{N}_{\l,\,'\!\l}$ from 
$\Z[q_{st},t_{st}]$ with the constant term $1$.

\smallskip
{\sf Other approaches.}
The Khovanov-Rozansky theory (unreduced or reduced)
is expected to be connected to the
physics {\em superpolynomials}  
based on the theory of {\em BPS states\,} \cite{DGR,AS,DMMSS,
FGS,GGS}.
This theory is not mathematically rigorous and the formulas
(for small knots/links) are mostly obtained 
via expected symmetries. We can {\em prove\,} 
them in  DAHA theory, so the coincidences of our formulas
with physical formulas is not surprising 
(unless they impose too many symmetries). 
For instance, the
approach in \cite{DMMSS} is very algebraic.
   
See e.g. \cite{Gor,GORS,GoN} for an important approach to 
superpolynomials
of uncolored torus knots based on {\em rational DAHA\,}.
It is expected that colors can be potentially added here 
(for torus knots); the  case of
symmetric powers of the fundamental representation is 
in progress, see \cite{GGS}. We will not touch this 
direction in this work. Using rational DAHA is connected
with the Hilbert schemes of plain curve singularities and
$C^2$ and with the {\em ORS-polynomials}  
from Part $(iv)$ below, so this approach is
related to our considerations.
\smallskip

{\sf Knot operators.}
Using the Macdonald polynomials
instead of Schur functions in the construction of the so-called
{\em knot opertators\,} was suggested in \cite{AS}, which triggered
paper \cite{CJ}. These operators for $t=q$ naturally appear in 
the approach to the invariants of torus knots via the Verlinde 
algebras. 
This method results in certain  
algorithms, but it is justified mathematically by now only 
for the root systems $A,D$. It requires using
the roots of unity $q$ ($t$ must be an integral
power of $q$) and the  
formula for the refined Verlinde $S$\~operator. In the refined
theory, with the Macdonald polynomials instead of
Schur functions, the formula for $S$ becomes
very involved (even for $A_1$). One must know {\em all}
Macdonald polynomials at roots of unity, which is 
almost impossible technically (though some formulas
for their coefficients are known).  
\smallskip

Furthermore,  
it must be justified that the final formulas 
are uniform in terms of $q$, i.e. can be lifted
to a generic $q,t$, which is quite a challenge. 
Formally, such a lift can be done only if the upper
bounds for the degrees of $q,t$ (and $a$) are known.
So the authors mainly compared their calculations with
known/expected formulas.
Only the simplest superpolynomials were 
discussed in \cite{AS} (mostly uncolored); they are 
reproduced in \cite{CJ} via DAHA. 

In more details, the refined $S$ operator is essentially the
matrix with the entries $P_b(q^{c+\rho_k})$ for all 
{\em admissible\,} (which depends on the root of unity) 
$b,c\in P_+$. Knowing {\em all\,}
$P_b(q^{c+\rho_k})$ is generally a {\em transcendental\,} problem
unless for $A_1$, and one need the whole projective
action of $PSL_2(\Z)$ due to \cite{Ki} and the works by I.Ch.
(see \cite{C101}). Furthermore, all $A_n$ are necessary for the 
stabilization.  It is the price of using 
the polynomial representation in \cite{Ki,AS} at
roots of unity vs. the direct usage  of $\HH$ suggested in 
\cite{CJ}. These problems were resolved there 
(actually bypassed) and the DAHA-Jones theory was extended to
any reduced root systems and dominant weights.

Having said this, DAHA at roots of unity is an important component
of the general theory. The {\em generalized Verlinde algebras\,}, 
also called {\em perfect DAHA representations\,} are one of the
main applications of DAHA. They are used in the theory of the 
$J\!D$\~polynomials and DAHA superpolynomials, 
especially toward $3$\~folds and applications in Number Theory.
The roots of unity were employed in \cite{CJJ} for the 
justification of the super-duality (a sketch).
\smallskip

\subsection{\bf Connection Conjecture}
From now on, let $\,\hat{\h}_{\l,\,'\!\l}\,(q,t,a)_{st}\,$
denote \,$\hat{\h}_{\l,\,'\!\l}\,(q,t,a)\,$ in
Theorem \ref{STABILIZ} expressed
in terms of the standard topological parameters 
(see \cite{CJ} and Section 1 in \cite{ORS}):
\begin{align}\label{qtareli}
&t=q^2_{st},\  q=(q_{st}t_{st})^2,\  a=a_{st}^2 t_{st},\notag\\
&q^2_{st}=t,\  t_{st}=\sqrt{q/t},\  a_{st}^2=a\sqrt{t/q}.
\end{align}
I.e. we use the substitutions from the first line here 
to obtain the polynomial $\,\hat{\h}_{\l,\,'\l}\,(q,t,a)_{st}$
from $\hat{\h}_{\l,\,'\l}\,(q,t,a).$

We will consider the integral forms $\hat{P}_\la=J_\la$ and
$\hat{P}_\la=P_\la^{(\!1\!)}$; the
latter is for $t=q$ when the Macdonald polynomials
coincide with the Schur functions. Also,  $\hat{\h}^{j_o}$
is the hat-normalization from Theorem \ref{STABILIZ}, where
$j_o=\emptyset$ means that there are no divisions 
by the evaluations at $q^\rho$. 

\begin{conjecture}\label{CONCONJ}
(i) For $t=q$ and $\hat{P}_\la=P_\la^{(\!1\!)}$, we conjecture 
that 
\begin{align}\label{conjhomfly}
&\hat{\h}^\emptyset_{\l,\,'\!\l}\,
(q,t\!\mapsto\! q,a\!\mapsto\! -a)_{st}=
\hat{\hbox{\small H\!O\!M}}_{\l,\,'\!\l}\,(q_{st},a_{st}),
\end{align}
where the latter is the hat-normalization 
of the unreduced HOMFLY-PT polynomial
for any pair of graphs $\{\l,\,'\!\l\}$ colored by an arbitrary 
sequences $\lla,\,'\!\lla$ of 
Young diagrams. Equivalently,
\begin{align}\label{conjhomflyy}
&\hat{\h}^{j_o}_{\l,\,'\!\l}\,
(q,t\!\mapsto\! q,a\!\mapsto\! - a)_{st}=
\hat{\hbox{\small H\!O\!M}}^{j_o}_{\l,\,'\!\l}\,
(q_{st}, a_{st}),
\end{align}
where 
$\hat{\hbox{\small H\!O\!M}}^{j_o}_{\l,\,'\!\l}\,
(q_{st}, a_{st})$
is the corresponding reduced HOMFLY-PT polynomial.
One can replace $\,'\!\l$ by $\,'\!\l^\vee$ in (\ref{conjhomfly})
and (\ref{conjhomflyy}).

(ii) Now let $\{\l,\,'\!\l^\vee\}$ be a pair of trees 
such that $\rr_i^j,\,\ss_i^j,\,'\rr_i^j,\,'\ss_i^j>0$, and 
$\,'\ss_1\ss_1>\,'\rr_1\rr_1$, where $1\le j\le \kappa$, $ 1\le 
i\le \ell^j$ and $1\le j\le \,'\!\kappa$, $ 1\le i\le 
\,'\!\ell^j$. Here and in $(iii)$ only 
$\,'\!\l^\vee$ with $\vee$ (can be empty) is considered and the 
integral form is $\hat{P}_\la=J_\la$. Then formula 
(\ref{deg-a-jj}) for deg$_a$ is conjectured to be exact. In the 
uncolored case, we also expect the positivity of the following
series: 
\begin{align}\label{posit-claim}
\hat{\h}^{min}_{\,\l,\,'\!\l^\vee}\,(q,t,a)/(1-t)^
{\kappa+'\!\kappa-1}\in
\Z_+[[q,t,a]]
\end{align}
upon the natural $t$\~expansion of this ratio. For any diagrams,
it is conjectured to hold for sufficiently
large powers of $(1-t)(1-q)$\, (provided $\rr_i^j,\,\ss_i^j,\,
'\rr_i^j,\,'\ss_i^j,\,
'\ss_1\ss_1\!-\!'\rr_1\rr_1>0$ and with the usage of $\vee$).

(iii) Furthermore, let $\la^j=\square=\la(\om_1)=\,'\!\la^j$ 
for all $j$ (the uncolored case).  Then $\{J_{\Box}\}_{ev}=
t^{1/2}(1+a)/(a^2)^{1/4}$ and 
we conjecture that for $\hat{P}_\la=J_\la$
and for the hat-normalization above:
\begin{align}\label{khrconj}
\Bigl(\frac{\hat{\h}^{\emptyset}_{\,\l,\,'\!\l^\vee}}
{(1-t)^{\kappa+'\!\kappa}}\Bigr)_{\!st}\!= 
\Bigl(\frac{(1\!+\!a)\,\hat{\h}^{min}_{\,\l,\,'\!\l^\vee}}
{(1-t)^{\kappa+'\!\kappa}}\Bigr)_{\!st}\!=
\hat{K\!h\!R}{}^{\hbox{\tiny stab}}_{\l,\,'\!\l^\vee}
(q_{st},t_{st},a_{st}).
\end{align}
The topological setting is unreduced here,
so $K\!h\!R{}^{\hbox{\tiny stab}}$ are polynomials in
terms of $a,q$ with the coefficients that are
generally infinite $t$\~series.

(iv) Conjecture 2 from \cite{ORS} states that \ 
$
K\!h\!R^{\hbox{\tiny stab}}_{\l,\,'\!\l}=
\overline{\mathscr{P}}_{\! \hbox{\tiny alg}},
$
where the latter series is defined there
for the corresponding germ of plane curve singularity
(see (\ref{yxcurve}))
in terms of the weight filtration in the cohomology of its nested
Hilbert scheme. Thus the series from
(\ref{khrconj}) can be also expected
to coincide under the  hat-normalization with
$\overline{\mathscr{P}}_{\! \hbox{\tiny alg}}$.
\sq
\end{conjecture}
\smallskip

Combining our paper with Section 7.1 from \cite{MoS}
proves Part (i) (we will post the details 
somewhere). Also we can follow here
Proposition 2.3 from \cite{CJ} (for torus knots),
where we used \cite{Ste}.
This approach is based on the DAHA shift operators and Verlinde 
algebras.  

Instead of using \cite{MoS} or
the knot operators from {\em CFT\,} 
and the Verlinde algebras, one can directly apply 
the {\em Rosso-Jones cabling formula} \cite{RJ,Mo,ChE} 
upon its relatively straightforward adjustment to iterated links. 

This formula used together with the theory of
DAHA {\em shift opertor\,} gives a relatively
straightforward way for the justification of Part $(i)$.
We have a sketch of a proof, which follows \cite{CJ} and
especially 
the case of $A_1$ (iterated torus knots) considered in detail in 
Proposition 4.2 from \cite{ChD}. This leads to an exact match,
not only to the coincidence up to proportionality
(under the hat-normalization). 

The advantage of this approach is that it
can be potentially extended to the 
{\em WRT invariants\,} for any root systems
and for any iterated torus links (following \cite{ChD}). 
For torus knots, the connection with the HOMFLY-PT
polynomials from \cite{CJ} at $t=q$ was extended now 
to the coincidence of the DAHA {\em hyperpolynomials\,}
of type $D$ at $t=q$ to the Kauffman polynomials; see there. 
\smallskip

In the uncolored case, the ratio
$(1\!+\!a)\,\hat{\h}^{min}_{\,\l,\,'\!\l}/
(1-t)^{\kappa+'\!\kappa}$
in Part $(iii)$ for $\hat{P}=J$ becomes exactly 
$\hat{\h}^\emptyset_{\l,\,'\!\l}$ for $\hat{P}=P^{(\!1\!)}$
from Part $(i)$ upon $t=q$. Accordingly, the reduced 
variant of the conjecture from Part $(i)$ becomes in the
uncolored case as follows:
\begin{align}\label{conjhomflyyred}
&\hat{\h}^{min}_{\l,\,'\!\l}\,
(q_{st}^2,q_{st}^2,-a_{st}^2)/(1-q_{st}^2)^{\kappa+'\!\kappa-1}=
\hat{\hbox{\small H\!O\!M}}^{reduced}_
{\l,\,'\!\l}\,(q_{st},a_{st}),
\end{align}
where the integral $J$\~form is used in $\hat{\h}^{min}$.  

\smallskip

Concerning the relation to the plane curve singularities,
Proposition 3 from \cite{ORS} matches the DAHA super-duality and
our estimate for \, deg$_a$\, from 
Part $(iii)$ of our Theorem \ref{STABILIZ} restricted to
the uncolored case. This is a confirmation
of Part $(iv)$ of the Connection Conjecture. However the relation 
of our construction
to $\overline{\mathscr{P}}_{\! \hbox{\tiny alg}}$  
the stable Khovanov-Rozansky 
polynomials can not be directly confirmed at the moment
beyond some cases of torus knots; see (\ref{khfromh}) below.
\smallskip

One of the reasons for such an uncertainty
is that the positivity of the polynomials 
$\hat{\h}^{min}_{\,\l,\,'\!\l^\vee}\,(q,t,a)$
does not hold for links (including uncolored ones) and for knots
if the corresponding Young diagrams are non-rectangle. We
address it in the positivity conjecture from Part $(ii)$,
which however corresponds to the {\em non-reduced theory}. 
Recall that (\ref{posit-claim}) is for any  
{\em positive pairs\,}  $\{\l,\,'\!\l^\vee\}$ and can be 
extended to any weights. 
 The experiments 
show that the positivity of the series there almost always  
fails when the inequalities from Part $(ii)$  are not satisfied 
(even in the uncolored case). This is in contrast to
the Khovanov-Rozansky theory, where knots are arbitrary.  

For quite a few 
colored algebraic links, the division by $(1-t)^\bullet$ is 
sufficient for the positivity in Part $(ii)$. However 
we have examples when the positivity holds with 
$(1-q)^\bullet$ and fails for the division by
any powers of $(1-t)$.
There are also examples when only proper powers of
$(1-q)(1-t)$ ensure the positivity. We did 
not reach any conjectures concerning the occurrence of
$t$ or $q$ here
and the minimal powers of these corresponding binomials.
\smallskip

One can also try to
replace the negative terms $-q^l t^m$ by $q^l t^{m\pm 1}$
following (\ref{khfromh}) below for small links (but this
is of experimental nature).  
An obvious problem is that there are practically no known 
formulas of the stable $K\!h\!R$\~polynomials and  
$\overline{\mathscr{P}}_{\! \hbox{\tiny alg}}$
beyond those for some (relatively simple) torus knots. 

Therefore Parts $(iii,iv)$ of Conjecture \ref{CONCONJ}
are not exactly verifiable conjectures at the moment, with a 
reservation concerning the positivity claim (\ref{posit-claim})
from $(ii)$, which is well confirmed in examples. We mention
that the relation with the stable $K\!h\!R$\~polynomials 
was conjectured in \cite{ChD} for {\em pseudo-algebraic\,} 
iterated 
torus knots (with positive DAHA-superpolynomials), not only for 
algebraic knots. We also suggested there some procedures
of experimental nature, which hopefully may work for any cables
(see below).

\smallskip
Let us address a bit using the {\em differentials}.
It is generally difficult to 
calculate $K\!h\!R{}^{n+1}$ unless for $n=1$ (the celebrated 
Khovanov polynomials). Not many formulas are known (and all
known ones are uncolored so far).
The polynomials $K\!h\!R{}^{\hbox{\tiny stab}}_{\l,\,'\!\l}$\, 
are of more algebraic nature, but are actually no simpler. 
If the polynomial $K\!h\!R{}^{\hbox{\tiny stab}}$\, 
is known, recovering all individual 
$K\!h\!R{}^{n+1}$ from it is generally provided by the theory of 
differentials $\partial_{n+1}$ from \cite{Kh,Ras}, but 
this is quite a challenge.

The corresponding homology 
$Ker(\partial_{n+1})/Im(\partial_{n+1})$ 
gives $K\!h\!R{}^{n+1}$ for any $n\ge 1$. These differentials  
are generally involved, but their certain algebraic 
simplifications,
suggested in \cite{Ras,DGR} and developed further in \cite{CJ},
work surprisingly well for sufficiently small links. 
The assumption 
in \cite{CJ} is that the actual $\partial_{n+1}$ are 
``as surjective as possible" beginning with $a=0$. See 
Conjecture 2.7 and 
Section 3.6 in \cite{CJ}; the ``smallest"  torus knot when 
(reduced) $K\!h\!R{}^{2}$ cannot be obtained this way is 
$T(12,7)$. Quite a few examples for torus iterated knots can 
be found in \cite{ChD}.

Applying this procedure is possible for links, but the 
{\em unreduced\,} setting create problems. Each particular
homology is finite but there are infinitely many non-trivial ones. 
We hope to address this in further works (at least in examples).
Let us consider now a special case.
\smallskip

{\sf Recovering $K\!h\!R^{n+1}$ from $J\!D^{n+1}$.} Generally, the
procedure of obtaining the $K\!h\!R$\~polynomials requires
knowing the whole superpolynomials and, moreover, all
differentials. However, the direct recovering 
{\em reduced\,} $K\!h\!R^{n+1}$ from 
$J\!D^{n+1}$ is not impossible for small links. This was noted
in \cite{CJ}. Due to the lack of other ways for justifying 
Part $(iii)$, this provides at least something,
and such an approach is  applicable to {\em any\,} torus iterated 
links (not only algebraic ones). 


Practically
it works as follows. For $\hat{P}=J$ and
modulo the operations of changing the sign explained below 
(we put $\doteq$), 
one can expect that   
\begin{align}\label{khfromh}
\Bigl(\frac{(1\!+\!a)\,\hat{\h}^{min}_{\,\l,\,'\!\l}}
{(1-t)^{\kappa}}\Bigr)_{q\,\mapsto (qt)^2,\, t\,\mapsto q^2,\,
a\,\mapsto -q^{2(n+1)}}\, \doteq\, 
\hat{K\!h\!R}{}^{n+1}_{\l,\,'\!\l}
(q,t)_{reduced}\,.
\end{align}

The latter polynomials can be calculated in few cases for $n>1$
and for many links  for $n=1$ using of the procedure 
$K\!hReduced[\cdot][q,t]$ from \cite{KA} 
(the hat-normalization must be applied). If there is no exact 
match, then the Connection Conjecture hints that the typical 
corrections are
the substitutions\, $-q^l t^m\mapsto +q^l t^{m+1}$,
$q^l t^m\mapsto +q^l t^{m+2}$ and so on. Here  
$-q^l t^m\mapsto +q^l t^{m-1}$ (etc.) can occur too
for {\em non-algebraic\,} links; see \cite{ChD}.
Generally the DAHA part in (\ref{khfromh}) is 
smaller than the right-hand side, but  in the examples we provide 
below they have the same number of terms.  Also, see the 
discussion after Conjecture 2.4 in \cite{ChD}. We will 
demonstrate
(\ref{khfromh}) for some basic knots.

\medskip

\comment{
The output of this reduction 
procedure will be denoted by $D\!A\!H\!A' K\!h\!R{}^{n+1}$. 

For example, if $\hat{\h}$ (in the DAHA-parameters) contains
$q^{i}(t^j a^{m}+ t^{j-n-1}a^{m+1}+t^{j-2n-2}a^{m+2})$, where
$q^{i}t^i a^m$ was not involved in the
reductions for smaller deg$_a$, then 
$q^i t^{j-2n-2}a^{m+2}$ will go to $D\!A\!H\!A' K\!h\!R_{n+1}$ 
from this triple (subject to possible further inductive 
reductions).
This is how the algebraic simplification of $\partial_{n+1}$ is
defined. 

However if the actual ({\em topological})  $\partial_{n+1}$ is not 
{\em onto\,} when acting from the space associated with  
$q^i t^{j-n-1}a^{m+1}$ to that for $q^i t^{j}a^m$, then 
$q^{i}t^j a^{m}$ is the right term to pick for 
$K\!h\!R^{n+1}$.  Thus there can be 
{\em topological corrections\,},
which makes the above recovery algorithm insufficient. 
Nevertheless, even if these corrections are present, one can
expect that in the reduced setting,
$\hat{K\!h\!R}{}^2-D\!A\!H\!A'K\!h\!R^{2}$ is a linear combination 
of  $q_{st}^i t_{st}^j(1+t_{st})$ (the case of $\partial_2=0$), 
$q_{st}^i t_{st}^j(1-t_{st}^2)$ and so on 
{\em with non-negative coefficients\,}. If this is the case,
then such a positivity can be considered a confirmation 
of Parts $(iii,iv)$.  Such positivity holds 
in all examples of 
{\em pseudo-algebraic\,} torus iterated knots (by definition,
with positive $\h$\~polynomials) we considered
in \cite{ChD}.

One can combinatorially extend this procedure to arbitrary 
iterated links (with non-positive $\h$)
changing all "\!$-1$" by $t^{n+1}/a$ in 
$\hat{\h}_{\l,\,'\!\l}(q,t,a)$ (in the uncolored case) 
followed by the reduction as above. Practically, one deletes
in $\hat{\h}$ all pairs 
$\pm t^i q^j a^m(t^{n+1}+a)$ 
and $t^i q^j a^m(t^{2n+2}-a^2)$ beginning with $m=0$,
then switches to the standard parameters and substitutes 
$-1\mapsto 1/t_{st}$. 

We will also use the variant
of this algorithm, denoted by 
$D\!A\!H\!A\!^+\!K\!h\!R$, where $-1$ is replaced 
by $a/t^{n+1}$ before
applying the reduction procedure.  
Presumably
$D\!A\!H\!A\!^+\!K\!h\!R$ must be used here when the $a$\~leading 
term of $\hat{\h}(q,t,a)$ is positive, but the 
evidence is limited so far. 
}


\setcounter{equation}{0}
\section{\sc Multiple torus knots}
\subsection{\bf Preliminary remarks}
We will consider examples (mainly of numerical kind) confirming
our theorems and the Connection Conjecture, including the
positivity claim (\ref{posit-claim}) for algebraic links. 
We selected only relatively simple and instructional ones; 
however, some formulas are long, which reflects the nature of
this theory. Explicit formulas are of obvious value;
they are expected to contain a 
lot of geometric information (beyond what we
discuss here and know now). We provide examples
only for links here;  see \cite{CJ,CJJ} for DAHA-superpolynomials 
for torus knots and \cite{ChD} for iterated torus knots. 

We present colored links $\l$ in the form of (\ref{tau-link})
and Theorem \ref{STABILIZ}:
\begin{align}\label{link-present}
\l=\l_{(\{\vec\rr^1,\vec\ss^1\},\ldots,
\{\vec\rr^\kappa,\vec\ss^\kappa\})}
^{\,\Up,\, (\la^1,\ldots,\la^\kappa)},\where \Up
\hbox{\, is presented graphically.}
\end{align}
Note that for a given $j$ (the index of the path),
we first collect $\rr^j_i$ in $\vec{\rr}^j$ and then
collect the corresponding $\ss^j_i$ in $\vec{\ss}^j$
(i.e. separate $\rr^j_i$ from $\ss^j_i$).
The labels
of vertices are the pairs $[\rr^j_i,\ss^j_i]$, identified with
respect to the incidence graph $\Up$.
Practically, we put $\bigl\{\vec\rr^j,\vec\ss^j\bigr\}$
as 
$\begin{pmatrix}
\rr_1^j,&\ldots,&\rr_\ell^j\\ 
\ss_1^j,&\ldots,&\ss_\ell^j
\end{pmatrix}$, but
such a 2-row presentation is graphically unreasonable in a paper,
especially in the indices of $\l$. {\em 
The square brackets $[\rr,\ss]$ is used only for labels 
in this work\,}; also, recall
that the first left vertex is $[\rr_1^j,\ss_1^j]$.

We will constantly use the natural diagrams for $\Up$. 
For instance,
$$
\{\Up,(\la^j)\}=
\{\circ_1\rightarrow\circ_2\,\, \raisebox{1pt}{\rightthreearrow},
\  
(\la^1,\la^2,\la^3)\}
$$
means that $\kappa=3$
and there are two vertices, shown as $\circ_{1,2}$ 
with the labels $[\rr^j_1,\ss^j_1],[\rr^j_2,\rr^j_2]$, 
which are the same for $j=1,2,3$. So the paths here are
different only by the arrowheads. Formally, we need
to repeat the same
$\,\vec\rr,\vec\ss\,$ three times (for each path) in
$\l$, but we will mostly omit coinciding\,
$\vec\rr^j,\vec\ss^j$\, if this is clear from the
graph. Thus the corresponding (labeled, colored) tree $\l$
will be presented as  
$$
\l_{(\bigl\{\{\rr_1,\rr_2\},\{\ss_1,\ss_2\}\bigr\},
\bigl\{\{\rr_1,\rr_2\},\{\ss_1,\ss_2\}\bigr\},
\bigl\{\{\rr_1,\rr_2\},\{\ss_1,\ss_2\}\bigr\})}
^{\,\circ\rightarrow\circ\rightthreearrow,\, 
(\la^1,\la^2,\la^3)} \hbox{\, or\, }
\l_{(\bigl\{\{\rr_1,\rr_2\},\{\ss_1,\ss_2\}\bigr\})}
^{\,\circ\rightarrow\circ\rightthreearrow,\, 
(\la^1,\la^2,\la^3)}.
$$

\smallskip
{\sf HOMFLY-PT polynomial.\,}
The relation
(\ref{conjhomflyy}) from the Connection
Conjecture, which is equivalent to
(\ref{conjhomfly}) there, is proven, so we give  
only $\hat{\h}{}^{min}_{\l}$ and (sometimes)  
$$
\hat{\h}{}^{min}_{\l}(q,t\mapsto q,a\mapsto -a), \where
\hat{P}_\la=J_\la.
$$
The corresponding denominators 
$\,den^{j_o}=den^{j_o}(\l,\,'\!\l)$ are provided :
\begin{align}\label{conjhomflyz}
&(\hat{\h}^{min}_{\l,\,'\!\l}\,
(q,q,-a)/den^{j_o})_{st}=
\hat{\hbox{\small H\!O\!M}}^{j_o}_{\l,\,'\!\l}\,(q_{st},a_{st}).
\end{align}

Here $()_{st}$ in the right-hand side means that
one needs to substitute $a\mapsto a_{st}^2,
q\mapsto q_{st}^2$ to compare our polynomials
with the usual presentation for the HOMFLY-PT polynomials. 
To calculate the latter, we mainly use 
the software by S.~Artamonov for 
colored HOMFLY-PT polynomials; see also \cite{AM,AMM}.

Generally, in notations from (\ref{stabevaly})
and (\ref{P-arms-legs}):
\begin{align}\label{den-hom}
\hat{\h}^{j_o}_{\l}=&
\hat{\h}^{min}_{\l}\,
\Pi_{\la^1\vee\ldots\vee\la^\kappa}^\dag/\Pi_{\la^{j_o}}^\dag,\,\
h_\la=J_{\la}/P_{\la}={\prod_{\Box\in\la}
(1-q^{arm(\Box)}t^{leg(\Box)+1})},\notag
\\
&den^{j_o}\!=\!
\Bigl(\frac{\Pi_{\la^{j_o}}^\dag}
{\Pi_{\la^1\vee\ldots\vee\la^\kappa}^\dag}\,
\prod_{j \neq j_o}\frac{J_{\la^j}}{P_{\la^j}}
\Bigr)_{t\mapsto q}\ ,\ \, 
\Pi_\la^\dag\!=\!
\prod_{p=1}^n \prod_{v=0}^{\la_{p\,}-1} 
\Bigl( 1\!-\!a\,q^{v}t^{1\!-\!p}\Bigr).
\end{align}
Note that $den^{j_o}=(1-q)^{\kappa-1}$ in the uncolored case
(for any $j_o$). We omitted $'\!\l$ here; use
(\ref{den-hom}) with $\lla$ replaced 
by $\lla\cup\,'\!\lla$ if it is present. Also, 
$den^{j_o}$ for $\l,\,'\!\l$ serve
$\l,\,'\!\l^\vee$ as well.
\smallskip

The connection of our construction
with the HOMFLY-PT polynomials is solid (it holds in all
examples we considered), as well as Part $(ii)$ of the
Connection Conjecture. Parts $(iii,iv)$ are much less 
convincing at the moment for torus
iterated {\em links\,} due to the lack (actually the absence) of
formulas for stable Khovanov-Rozansky polynomials and
their algebraic-geometric counterparts from \cite{ORS}
is such a generality. As for the iterated {\em knots\,},
especially for relatively small torus knots, there are
direct confirmations (not too many). The positivity property
(\ref{posit-claim}) and its generalizations for any weights
hold for algebraic links in all examples we considered
(and powers of $(1-t)$ are mostly sufficient); this is  
an indirect confirmation of the Connection Conjecture. 

The symmetries of DAHA
superpolynomials match well those (known and expected) in
geometry/topology/physics.  
The super-duality is important. It is proven in \cite{ORS} for 
uncolored algebraic links within their theory, which is
another (indirect) confirmation of the Connection Conjecture.
We note that this and other symmetries provide an 
{\em excellent\,} way 
to verify our numerical simulations (we always check all of
them).
\smallskip

\comment{
\begin{figure}[htbp]
\begin{center}
\vskip -0.0in
\hskip -0.5in
\includegraphics[scale=0.8]{5-braids.eps}
\vskip +0.in
\caption{
{Cables (11,2),(13,2),(15,2) of T(3,2),}
\newline 
{Cab(19,3)T(3,2) and Cab(25,2)T(4,3)}}
\label{bngt}
\end{center}
\end{figure}
\medskip
}

\subsection{\bf Uncolored 2-fold trefoil}
The corresponding links will be in this section uncolored
$\l=T(\kappa\, \rr,\kappa\,\ss)$. We will mostly provide only 
{\em minimal\,}
$\hat{\h}{}^{min}_{\l}$ (for $\hat{P}=J$); 
there will be no $\,'\!\l$ in this section. We begin with
the ``canonical" example of $T(6,4)$.

{\sf 2-fold trefoil.} The \tax-presentation  
and $\hat{\h}$ are as follows:
$$
1\le j\le \kappa=2,\ \vec\rr^j=3,\, \vec\ss^j=2,\, 
\Up=\{\circ\rightrightarrows\}\,,\ 
\la^1=\square=\la^2;
$$
\begin{align}\label{T6-4}
&T(6,4):\ 
\l=\l_{(\{3,2\},\{3,2\})}^{\,\circ\rightrightarrows,\, 
(\square,\square)},\ \ 
\hat{\h}{}^{min}_{\l}\,(q,t,a)=
\end{align}

\renewcommand{\baselinestretch}{0.5} 
\noindent
{\small
\(
1-t+q t+q^2 t+q^3 t-q t^2+2 q^4 t^2-q^2 t^3-q^4 t^3
+2 q^5 t^3-q^3 t^4-q^5 t^4+2 q^6 t^4-q^4 t^5+q^7 t^5-q^5 t^6
+q^7 t^6-q^6 t^7+q^7 t^7-q^7 t^8+q^8 t^8+
a^3 \bigl(q^6-q^6 t+q^7 t-q^7 t^2+q^8 t^2\bigr)+
a^2 \bigl(q^3+q^4+q^5-q^3 t+q^5 t+2 q^6 t-q^4 t^2-q^5 t^2
+2 q^7 t^2
-q^5 t^3-q^6 t^3+q^7 t^3+q^8 t^3-q^6 t^4+q^8 t^4-q^7 t^5
+q^8 t^5\bigr)+a \bigl(q+q^2+q^3-q t+q^3 t+3 q^4 t+q^5 t-q^2 t^2
-q^3 t^2-q^4 t^2+3 q^5 t^2+q^6 t^2-q^3 t^3-q^4 t^3-2 q^5 t^3
+3 q^6 t^3+q^7 t^3-q^4 t^4-q^5 t^4-q^6 t^4+3 q^7 t^4-q^5 t^5
-q^6 t^5
+q^7 t^5+q^8 t^5-q^6 t^6+q^8 t^6-q^7 t^7+q^8 t^7\bigr).
\)
}
\renewcommand{\baselinestretch}{1.0} 
\smallskip

The $a$\~degree of $\hat{\h}{}^{min}_{\l}\,(q,t,a)$
is $3$, which matches the formula deg$_a=\ss(2|\la|)-|\la|=3$
from (\ref{deg-a-jj}); see also  Part $(iii)$ of Theorem 
\ref{STABILIZ}. The self-duality and other claims
in this theorem hold. For instance, the transposition
$\rr\leftrightarrow \ss$ in
$T(\kappa \rr,\kappa \ss)$ here and below
does not influence the superpolynomial,
though the practical calculations can be quite different
depending on the order of $\rr$ and $\ss$.
For instance, the label $[3,2]$ corresponds to 
$\tau_+\tau_-^2$, but  $[2,3]$ is naturally represented 
by $\tau_-\tau_+\tau_-$.

Note that the series from the positivity claim from
(\ref{posit-claim}) in Conjecture \ref{CONCONJ} reads
$\hat{\h}{}^{min}_{\l}(q,t,a)/(1-t)$ for uncolored $2$\~links.
This positivity holds here and in all examples we considered.

Our superpolynomial matches that suggested in 
\cite{DMMSS} in Section 2.8. The main factor of
$-P_1^{T[4,6]}$ coincides with our one upon the
following substitution:
$ a\mapsto A^2, q\mapsto q^2, t\mapsto t^2$. I.e.
their {\em non-bold\,} $A,q,t$ are essentially the
DAHA parameters.

\smallskip

{\em HOMFLY-PT polynomial.\,}
Recall that Part $(i)$ of the Connection Conjecture
(which is a theorem) claims that                          
\begin{align*}
&(\hat{\h}^{min}_{\l}\,
(q,q,-a)/den)_{st}=
\hat{\hbox{\small H\!O\!M}}^{j_o}_{\l}\,(q_{st},a_{st}),
\end{align*}
with the denominators explained above. 

In this example,
$den=(1-q)$\, and\, \  $\hat{\h}{}^{min}_{\l}(q,q,-a)=$
\renewcommand{\baselinestretch}{0.5} 

\noindent
{\small
\(
1-q+q^2+q^4-q^5+2 q^6-2 q^7+2 q^8-2 q^9+2 q^{10}-q^{11}
+q^{12}+q^{14}-q^{15}+q^{16}-a^3 \bigl(q^6-q^7+q^8-q^9
+q^{10}\bigr)+a^2 \bigl(q^3+q^5+q^7-q^8+q^9+q^{11}+q^{13}\bigr)
-a \bigl(q+q^3+2 q^5-q^6+2 q^7
-2 q^8+2 q^9-q^{10}+2 q^{11}+q^{13}+q^{15}\bigr).
\)
}
\renewcommand{\baselinestretch}{1.0} 
\smallskip


{\em Khovanov polynomial.\,}
Let us apply the procedure from (\ref{khfromh})
to obtain the DAHA approximation to the reduced
Khovanov polynomial ($n=1$). 
We switch to the standard parameters:
$q\mapsto (qt)^2,\, t\mapsto q^2,\,a\mapsto -q^{2(n+1)}$
(we add $\{st,n\}$ but omit $st$ in $q_{st},t_{st},a_{st}$).
Then
\begin{align}\label{khfromh4-6}
\Bigl(\frac{(1\!+\!a)\,\hat{\h}^{min}_{\l}}
{(1-t)^{2}}\Bigr)_{st, n=1}=
\end{align}
{\small \(
1+q^2+q^4 t^2-q^8 t^2+q^6 t^4+q^8 t^4-q^{10} t^4
-q^{12} t^4+q^8 t^6+q^{10} t^6-q^{12} t^6-q^{14} t^6
+2 q^{12} t^8-2 q^{16} t^8+q^{16} t^{10}-q^{18} t^{10}
-q^{20} t^{10}+q^{22} t^{10}+q^{20} t^{12}+q^{22} t^{12},\)}
\begin{align}\label{khfromh4-6-kh}
(\hat{K\!h\!R}{}^{2}_{\l})_{reduced}=
\end{align}
{\small \(
1+q^2+q^4 t^2+q^8 t^3+q^6 t^4+q^8 t^4+q^{10} t^5+q^{12} t^5
+q^8 t^6+q^{10} t^6+q^{12} t^7+q^{14} t^7+2 q^{12} t^8
+2 q^{16} t^9+q^{16} t^{10}+q^{18} t^{11}
+q^{20} t^{11}+q^{20} t^{12}+ q^{22} t^{12}+ q^{22} t^{12}.
\)}

Their difference (the second minus the first)
is a sum of the terms $q^l(t^m+t^{m-1})$,
$q^l(t^m-t^{m-2})$, which exactly matches the
expectations concerning (\ref{khfromh}).

\smallskip

\subsection{\bf Similar links} 
Let us consider {\sf 2-fold T(2,1)}. One has:
$$
1\le j\le \kappa=2,\ \vec\rr^j=2,\, \vec\ss^j=1,\, 
\Up=\{\circ\rightrightarrows\}\,,\ 
\la^1=\square=\la^2;
$$
\begin{align}\label{T4-2}
&T(4,2):\ 
\l=\l_{(\{2,1\},\{2,1\})}^{\,\circ\rightrightarrows,\, 
(\square,\square)},\ \ 
\hat{\h}{}^{min}_{\l}\,(q,t,a)=
\end{align}
{\small
$$
1 - t + qt - q t^2 + q^2t^2 + a(q - q t + q^2 t).
$$
}

The $a$\~degree is $\ss(2|\la|)-|\la|=1$. 
The self-duality reads:
$$
q^2 t^2 \hat{\h}{}^{min}_{\l}\,(t^{-1},q^{-1},a)\,=\,
\hat{\h}{}^{min}_{\l}\,(q,t,a),
$$ 
and $\hat{\h}$ is obviously positive upon  division by $(1-t)$.

Upon the passage to the standard parameters, formula (\ref{T4-2})
becomes 
$$
1-q^2+q^4 t^2-q^6 t^2+q^8 t^4+a^2(q^2 t^3-q^4 t^3+q^6 t^5),
$$
which is the major factor of formula (136) from \cite{DMMSS}.

\medskip
{\em HOMFLY-PT polynomial.\,} One has:
\begin{align*}
&\hat{\h}{}^{min}_{\l}(q,q,-a)=
1 - q + q^2 - q^3 + q^4 - a(q - q^2 + q^3)),\\
&\hat{\h}^{min}_{\l}(q^2,q^2,-a^2)/(1-q^2)=
\hat{\hbox{\small H\!O\!M}}_{\l}\,(q,a).
\end{align*}
\medskip

{\em Khovanov polynomial.\,} Recall that this and other formulas 
for (reduced) Khovanov polynomials will be
in terms of the standard parameters; we add the suffix 
$\{st,n=1\}$ to remind this. One has:

\begin{align}\label{khfromh4-2}
\Bigl(\frac{(1\!+\!a)\,\hat{\h}^{min}_{\l}}
{(1-t)^{2}}\Bigr)_{st,n=1}=
\hbox{\small \(
1 + q^2 + q^4t^2 - q^8t^2 + q^8t^4 + q^{10}t^4,
\)}\\
\label{khfromh4-2-kh}
(\hat{K\!h\!R}{}^{2}_{\l})_{reduced}=
\hbox{\small \(
1+q^2+q^4 t^2+q^8 t^3+q^8 t^4+q^{10} t^4,
\)}
\end{align}
where the difference $q^8(t^3+t^2)$ matches our expectations.
\medskip

{\sf $2$-fold T(1,1)\,}.
This is the  simplest algebraic torus link. 
\begin{align}\label{T2-2}
&T(2,2):\ 
\l=\l_{(\{1,1\},\{1,1\})}^{\,\circ\rightrightarrows,\, 
(\square,\square)},\,
\hat{\h}{}^{min}_{\l}=1 - t + q t +a q,\\
&\hbox{where\, \,}
\hat{\hbox{\small H\!O\!M}}_{\l}\,(q,a)=(1 - q + q^2- aq)/(1-q).
\notag
\end{align}
This is for $\ga=\tau_-$. The hat-normalization 
of the reduced Khovanov polynomial is
$1 + q^2 + q^4t^2 + q^6t^2$, which exactly coincides with
that obtained via (\ref{khfromh}); no sign-corrections
are necessary. 

We note that the superpolynomial becomes
$1+a(q+1/t-q/t)$ for $\ga=\tau_-^{-1}$, i.e. for the tree
$\l_{(\{1,-1\},\{1,-1\})}$. 
\smallskip

{\sf $2$-fold T(m,1)\,}. Let us provide a general formula for
$\l=\l_{(\{1,m\},\{1,m\})}^{\,\circ\rightrightarrows,\, 
(\square,\square)}$ with $m>0$:
\begin{align}\label{T2+2m} 
&\hat{\h}{}^{min}_{\l}\,(q,t,a)=
(1-t)(1+q t+q^2t^2 +q^3 t^3+\ldots+ q^{m-1}t^{m-1})\\
&+q^m t^m+a 
\bigl((1-t)(q+q^2t+\ldots+ q^{m-1}t^{m-2})+q^ m t^{m-1}\bigr).
\notag
\end{align}

The corresponding link is isotopic to 
$\l_{(\{m,1\},\{m,1\})}^{\,\circ\rightrightarrows,\,
(\square,\square)}$,
so their superpolynomials must coincide. This is a simplest
example of the symmetries of DAHA superpolynomials
in Part $(ii)$ of Theorem 
\ref{STABILIZ}. Recall that generally
we can transpose $\rr$ and $\ss$
only for the {\em first\,} vertex $[\rr_1^j,\ss_1^j]$
in the presence of iterations. Topologically this is obvious.

Let us provide a counterpart of this formula for non-algebraic
$\l=\l_{(\{1,-m\},\{1,-m\})}$ with $m>0$:
\begin{align}\label{T2-2m} 
&\hat{\h}{}^{min}_{\l}\,(q,t,a)=
(1-q)(1+q t+\ldots+q^{m-2}t^{m-2})+q^{m-1} t^{m-1}\\
&+a \bigl((1-q)(\frac{1}{t}+q+q^2t+\ldots+ 
q^{m-1}t^{m-2})+q^ m t^{m-1}\bigr).
\notag
\end{align}
One can change here $\{1,-m\}$ by $\{-m,1\}$, the
superpolynomial will remain the same.
Both formulas, for $\pm\, m$, satisfy the self-duality (with proper
$q^\bullet t^\bullet$\~multipliers), which is simple 
to see directly.
\medskip

{\sf 2-fold T(5,2).} The \tax-parameters  
and $\hat{\h}$ are as follows:
$$
1\le j\le \kappa=2,\ \vec\rr^j=5,\, \vec\ss^j=2,\, 
\Up=\{\circ\rightrightarrows\}\,,\ 
\la^1=\square=\la^2;
$$
\begin{align}\label{T10-4}
&T(10,4):\ 
\l=\l_{(\{5,2\},\{5,2\})}^{\,\circ\rightrightarrows,\, 
(\square,\square)},\ \ 
\hat{\h}{}^{min}_{\l}\,(q,t,a)=
\end{align}

\renewcommand{\baselinestretch}{0.5} 
\noindent
{\small
\(
1-t+q t+q^2 t+q^3 t-q t^2+2 q^4 t^2+q^5 t^2+q^6 t^2
-q^2 t^3-q^4 t^3+q^5 t^3+q^6 t^3+2 q^7 t^3-q^3 t^4-q^5 t^4
+3 q^8 t^4-q^4 t^5-q^6 t^5-q^8 t^5+3 q^9 t^5-q^5 t^6-q^7 t^6
-q^9 t^6+3 q^{10} t^6-q^6 t^7-q^8 t^7+2 q^{11} t^7-q^7 t^8-q^9 t^8
+q^{11} t^8+q^{12} t^8-q^8 t^9-q^{10} t^9+q^{11} t^9+q^{12} t^9
-q^9 t^{10}-q^{11} t^{10}+2 q^{12} t^{10}-q^{10} t^{11}
+q^{13} t^{11}-q^{11} t^{12}+q^{13} t^{12}-q^{12} t^{13}
+q^{13} t^{13}-q^{13} t^{14}+q^{14} t^{14}+a^3 \bigl(q^6-q^6 t
+q^7 t+q^8 t+q^9 t-q^7 t^2+2 q^{10} t^2-q^8 t^3-q^{10} t^3
+2 q^{11} t^3-q^9 t^4-q^{11} t^4+2 q^{12} t^4-q^{10} t^5
+q^{13} t^5-q^{11} t^6+q^{13} t^6-q^{12} t^7+q^{13} t^7
-q^{13} t^8+q^{14} t^8\bigr)+a^2 \bigl(q^3+q^4+q^5-q^3 t
+q^5 t+3 q^6 t+2 q^7 t+q^8 t-q^4 t^2-q^5 t^2-q^6 t^2+2 q^7 t^2
+3 q^8 t^2+3 q^9 t^2-q^5 t^3-q^6 t^3-2 q^7 t^3+2 q^9 t^3
+4 q^{10} t^3-q^6 t^4-q^7 t^4-2 q^8 t^4-q^9 t^4+q^{10} t^4
+4 q^{11} t^4-q^7 t^5-q^8 t^5-2 q^9 t^5-q^{10} t^5+2 q^{11} t^5
+3 q^{12} t^5-q^8 t^6-q^9 t^6-2 q^{10} t^6+3 q^{12} t^6
+q^{13} t^6-q^9 t^7-q^{10} t^7-2 q^{11} t^7+2 q^{12} t^7
+2 q^{13} t^7-q^{10} t^8-q^{11} t^8-q^{12} t^8+3 q^{13} t^8
-q^{11} t^9-q^{12} t^9+q^{13} t^9+q^{14} t^9-q^{12} t^{10}
+q^{14} t^{10}-q^{13} t^{11}+q^{14} t^{11}\bigr)
+a \bigl(q+q^2+q^3-q t+q^3 t+3 q^4 t+2 q^5 t+q^6 t-q^2 t^2-q^3 t^2
-q^4 t^2+2 q^5 t^2+3 q^6 t^2+4 q^7 t^2+q^8 t^2-q^3 t^3-q^4 t^3
-2 q^5 t^3+q^7 t^3+5 q^8 t^3+2 q^9 t^3-q^4 t^4-q^5 t^4-2 q^6 t^4
-q^7 t^4-q^8 t^4+5 q^9 t^4+2 q^{10} t^4-q^5 t^5-q^6 t^5-2 q^7 t^5
-q^8 t^5-2 q^9 t^5+5 q^{10} t^5+2 q^{11} t^5-q^6 t^6-q^7 t^6
-2 q^8 t^6-q^9 t^6-q^{10} t^6+5 q^{11} t^6+q^{12} t^6-q^7 t^7
-q^8 t^7-2 q^9 t^7-q^{10} t^7+q^{11} t^7+4 q^{12} t^7-q^8 t^8
-q^9 t^8-2 q^{10} t^8+3 q^{12} t^8+q^{13} t^8-q^9 t^9-q^{10} t^9
-2 q^{11} t^9+2 q^{12} t^9+2 q^{13} t^9-q^{10} t^{10}-q^{11} t^{10}
-q^{12} t^{10}+3 q^{13} t^{10}-q^{11} t^{11}-q^{12} t^{11}
+q^{13} t^{11}+q^{14} t^{11}
-q^{12} t^{12}+q^{14} t^{12}-q^{13} t^{13}+q^{14} t^{13}\bigr).
\)
}
\renewcommand{\baselinestretch}{1.0} 
\smallskip

The $a$\~degree of $\hat{\h}{}^{min}_{\l}\,(q,t,a)$
is $3$. It is the same as for the trefoil and
remains unchanged for all uncolored  $2$-fold $T(2\mm+1,2)$
due to formula (\ref{deg-a-jj}).
The positivity of 
$\hat{\h}{}^{min}_{\l}(q,t,a)/(1-t)$  holds. Recall that
the division by $(1-t)^{\kappa-1}$ is presumably sufficient to
ensure the positivity for {\em uncolored\,} algebraic links.

\hbadness=10000
\medskip
{\em HOMFLY-PT polynomial.\,}
Here
$den=(1-q)$\, and\, \  $\hat{\h}{}^{min}_{\l}(q,q,-a)=$
\renewcommand{\baselinestretch}{0.5} 

\noindent
{\small
\(
1-q+q^2+q^4-q^5+2 q^6-q^7+2 q^8-q^9+2 q^{10}-2 q^{11}+3 q^{12}
-3 q^{13}+3 q^{14}-3 q^{15}+3 q^{16}-2 q^{17}+2 q^{18}-q^{19}
+2 q^{20}-q^{21}+2 q^{22}-q^{23}+q^{24}+q^{26}-q^{27}+q^{28}
-a^3 \bigl(q^6-q^7+q^8+q^{10}-q^{11}+2 q^{12}-2 q^{13}+2 q^{14}
-2 q^{15}+2 q^{16}-q^{17}+q^{18}+q^{20}-q^{21}+q^{22}\bigr)
-a^2 \bigl(-q^3-q^5-2 q^7-2 q^9-2 q^{11}+q^{12}-2 q^{13}+2 q^{14}
-2 q^{15}+q^{16}-2 q^{17}-2 q^{19}-2 q^{21}-q^{23}-q^{25}\bigr)
-a \bigl(q+q^3+2 q^5+2 q^7+3 q^9-q^{10}+3 q^{11}-2 q^{12}+3 q^{13}
-3 q^{14}+3 q^{15}-2 q^{16}+3 q^{17}
-q^{18}+3 q^{19}+2 q^{21}+2 q^{23}+q^{25}+q^{27}\bigr).
\)
}
\renewcommand{\baselinestretch}{1.0} 
\smallskip

Recall that one needs to divide this polynomial by
$den$ and then 
change our $a,q$ to the {\em standard\,}
$a_{st}^2,q_{st}^2$ to compare it with the usual presentation
for the reduced HOMFLY-PT polynomials. 
\smallskip

\subsection{\bf Uncolored 2-fold 
\texorpdfstring{{\mathversion{bold}$T(4,3)$}}{T(4,3)}}
The family $T(2\mm+1,2)$ is well known to be
quite special (topologically
and algebraically). Let us
provide two uncolored examples for 2-fold $T(3\mm\pm 1,3)$.

{\sf 2-fold T(4,3).} The natural notation is  $T(8,6)$.
The $a$\~degree of $\hat{\h}{}^{min}_{\l}\,(q,t,a)$
is deg$_a=\ss(2|\la|)-|\la|=5$.
The positivity of the series
$\hat{\h}{}^{min}_{\l}(q,t,a)/(1-t)$ from 
(\ref{posit-claim}) holds.
The \tax-parameters  
and $\hat{\h}$ are as follows:
$$
1\le j\le \kappa=2,\ \vec\rr^j=4,\, \vec\ss^j=3,\, 
\Up=\{\circ\rightrightarrows\}\,,\ 
\la^1=\square=\la^2;
$$
\begin{align}\label{T8-6}
&T(8,6):\ 
\l=\l_{(\{4,3\},\{4,3\})}^{\,\circ\rightrightarrows,\, 
(\square,\square)},\ \ 
\hat{\h}{}^{min}_{\l}\,(q,t,a)=
\end{align}

\renewcommand{\baselinestretch}{0.5} 
\noindent
{\small
\(
1-t+q t+q^2 t+q^3 t+q^4 t+q^5 t-q t^2+q^4 t^2+q^5 t^2
+3 q^6 t^2+q^7 t^2+q^8 t^2-q^2 t^3-q^4 t^3+3 q^7 t^3+2 q^8 t^3
+3 q^9 t^3-q^3 t^4-q^5 t^4-q^6 t^4-q^7 t^4+2 q^8 t^4+q^9 t^4
+5 q^{10} t^4+q^{11} t^4-q^4 t^5-q^6 t^5-q^7 t^5-2 q^8 t^5+q^9 t^5
+5 q^{11} t^5+q^{12} t^5-q^5 t^6-q^7 t^6-q^8 t^6-2 q^9 t^6
-q^{11} t^6+6 q^{12} t^6+q^{13} t^6-q^6 t^7-q^8 t^7-q^9 t^7
-2 q^{10} t^7-q^{12} t^7+5 q^{13} t^7+q^{14} t^7-q^7 t^8-q^9 t^8
-q^{10} t^8-2 q^{11} t^8+5 q^{14} t^8-q^8 t^9-q^{10} t^9-q^{11} t^9
-2 q^{12} t^9+q^{13} t^9+q^{14} t^9+3 q^{15} t^9-q^9 t^{10}
-q^{11} t^{10}-q^{12} t^{10}-2 q^{13} t^{10}+2 q^{14} t^{10}
+2 q^{15} t^{10}+q^{16} t^{10}-q^{10} t^{11}-q^{12} t^{11}
-q^{13} t^{11}-q^{14} t^{11}+3 q^{15} t^{11}+q^{16} t^{11}
-q^{11} t^{12}-q^{13} t^{12}-q^{14} t^{12}+3 q^{16} t^{12}
-q^{12} t^{13}-q^{14} t^{13}+q^{16} t^{13}+q^{17} t^{13}
-q^{13} t^{14}-q^{15} t^{14}+q^{16} t^{14}+q^{17} t^{14}
-q^{14} t^{15}+q^{17} t^{15}-q^{15} t^{16}+q^{17} t^{16}
-q^{16} t^{17}+q^{17} t^{17}-q^{17} t^{18}+q^{18} t^{18}
+a^5 \bigl(q^{15}-q^{15} t+q^{16} t-q^{16} t^2+q^{17} t^2
-q^{17} t^3+q^{18} t^3\bigr)
\)

\smallskip\noindent
\(
+a^4 \bigl(q^{10}+q^{11}+q^{12}
+q^{13}+q^{14}-q^{10} t+q^{12} t+q^{13} t+q^{14} t+2 q^{15} t
-q^{11} t^2-q^{12} t^2+q^{14} t^2+q^{15} t^2+2 q^{16} t^2
-q^{12} t^3-q^{13} t^3-q^{14} t^3+q^{16} t^3+2 q^{17} t^3
-q^{13} t^4-q^{14} t^4-q^{15} t^4+q^{16} t^4+q^{17} t^4+q^{18} t^4
-q^{14} t^5-q^{15} t^5+q^{17} t^5+q^{18} t^5-q^{15} t^6-q^{16} t^6
+q^{17} t^6+q^{18} t^6-q^{16} t^7+q^{18} t^7-q^{17} t^8
+q^{18} t^8\bigr)
\)

\smallskip\noindent
\(
+a^3 \bigl(q^6+q^7+2 q^8+2 q^9+2 q^{10}+q^{11}+q^{12}-q^6 t+2 q^9 t
+3 q^{10} t+5 q^{11} t+3 q^{12} t+3 q^{13} t+q^{14} t-q^7 t^2
-q^8 t^2-2 q^9 t^2+q^{11} t^2+6 q^{12} t^2+4 q^{13} t^2
+4 q^{14} t^2+q^{15} t^2-q^8 t^3-q^9 t^3-3 q^{10} t^3-2 q^{11} t^3
-2 q^{12} t^3+5 q^{13} t^3+4 q^{14} t^3+5 q^{15} t^3+q^{16} t^3
-q^9 t^4-q^{10} t^4-3 q^{11} t^4-3 q^{12} t^4-4 q^{13} t^4
+4 q^{14} t^4+4 q^{15} t^4+4 q^{16} t^4+q^{17} t^4-q^{10} t^5
-q^{11} t^5-3 q^{12} t^5-3 q^{13} t^5-4 q^{14} t^5+5 q^{15} t^5
+4 q^{16} t^5+3 q^{17} t^5-q^{11} t^6-q^{12} t^6-3 q^{13} t^6
-3 q^{14} t^6-2 q^{15} t^6+6 q^{16} t^6+3 q^{17} t^6+q^{18} t^6
-q^{12} t^7-q^{13} t^7-3 q^{14} t^7-2 q^{15} t^7+q^{16} t^7
+5 q^{17} t^7+q^{18} t^7-q^{13} t^8-q^{14} t^8-3 q^{15} t^8
+3 q^{17} t^8+2 q^{18} t^8-q^{14} t^9-q^{15} t^9-2 q^{16} t^9
+2 q^{17} t^9+2 q^{18} t^9-q^{15} t^{10}-q^{16} t^{10}
+2 q^{18} t^{10}-q^{16} t^{11}+q^{18} t^{11}-q^{17} t^{12}
+q^{18} t^{12}\bigr)
\)

\smallskip\noindent
\(
+a^2 \bigl(q^3+q^4+2 q^5+2 q^6+2 q^7+q^8
+q^9-q^3 t+2 q^6 t+4 q^7 t+6 q^8 t+5 q^9 t+5 q^{10} t+2 q^{11} t
+q^{12} t-q^4 t^2-q^5 t^2-2 q^6 t^2-q^7 t^2+q^8 t^2+6 q^9 t^2
+7 q^{10} t^2+9 q^{11} t^2+4 q^{12} t^2+2 q^{13} t^2-q^5 t^3
-q^6 t^3-3 q^7 t^3-3 q^8 t^3-3 q^9 t^3+2 q^{10} t^3+6 q^{11} t^3
+12 q^{12} t^3+6 q^{13} t^3+3 q^{14} t^3-q^6 t^4-q^7 t^4-3 q^8 t^4
-4 q^9 t^4-5 q^{10} t^4-2 q^{11} t^4+2 q^{12} t^4+12 q^{13} t^4
+6 q^{14} t^4+3 q^{15} t^4-q^7 t^5-q^8 t^5-3 q^9 t^5-4 q^{10} t^5
-6 q^{11} t^5-4 q^{12} t^5+q^{13} t^5+12 q^{14} t^5+6 q^{15} t^5
+2 q^{16} t^5-q^8 t^6-q^9 t^6-3 q^{10} t^6-4 q^{11} t^6
-6 q^{12} t^6-4 q^{13} t^6+2 q^{14} t^6+12 q^{15} t^6+4 q^{16} t^6
+q^{17} t^6-q^9 t^7-q^{10} t^7-3 q^{11} t^7-4 q^{12} t^7
-6 q^{13} t^7-2 q^{14} t^7+6 q^{15} t^7+9 q^{16} t^7+2 q^{17} t^7
-q^{10} t^8-q^{11} t^8-3 q^{12} t^8-4 q^{13} t^8-5 q^{14} t^8
+2 q^{15} t^8+7 q^{16} t^8+5 q^{17} t^8-q^{11} t^9-q^{12} t^9
-3 q^{13} t^9-4 q^{14} t^9-3 q^{15} t^9+6 q^{16} t^9+5 q^{17} t^9
+q^{18} t^9-q^{12} t^{10}-q^{13} t^{10}-3 q^{14} t^{10}
-3 q^{15} t^{10}+q^{16} t^{10}+6 q^{17} t^{10}+q^{18} t^{10}
-q^{13} t^{11}-q^{14} t^{11}-3 q^{15} t^{11}-q^{16} t^{11}
+4 q^{17} t^{11}+2 q^{18} t^{11}-q^{14} t^{12}-q^{15} t^{12}
-2 q^{16} t^{12}+2 q^{17} t^{12}+2 q^{18} t^{12}-q^{15} t^{13}
-q^{16} t^{13}+2 q^{18} t^{13}-q^{16} t^{14}+q^{18} t^{14}
-q^{17} t^{15}+q^{18} t^{15}\bigr)
\)

\smallskip\noindent
\(
+a \bigl(q+q^2+q^3+q^4+q^5-q t+q^3 t+2 q^4 t+3 q^5 t+5 q^6 t
+3 q^7 t+2 q^8 t+q^9 t-q^2 t^2-q^3 t^2-q^4 t^2+q^6 t^2+6 q^7 t^2
+6 q^8 t^2+6 q^9 t^2+3 q^{10} t^2+q^{11} t^2-q^3 t^3-q^4 t^3
-2 q^5 t^3-2 q^6 t^3-2 q^7 t^3+3 q^8 t^3+5 q^9 t^3+10 q^{10} t^3
+6 q^{11} t^3+2 q^{12} t^3-q^4 t^4-q^5 t^4-2 q^6 t^4-3 q^7 t^4
-4 q^8 t^4+q^{10} t^4+10 q^{11} t^4+8 q^{12} t^4+3 q^{13} t^4
-q^5 t^5-q^6 t^5-2 q^7 t^5-3 q^8 t^5-5 q^9 t^5-2 q^{10} t^5
-2 q^{11} t^5+9 q^{12} t^5+8 q^{13} t^5+3 q^{14} t^5-q^6 t^6
-q^7 t^6-2 q^8 t^6-3 q^9 t^6-5 q^{10} t^6-3 q^{11} t^6-3 q^{12} t^6
+9 q^{13} t^6+8 q^{14} t^6+2 q^{15} t^6-q^7 t^7-q^8 t^7-2 q^9 t^7
-3 q^{10} t^7-5 q^{11} t^7-3 q^{12} t^7-2 q^{13} t^7+10 q^{14} t^7
+6 q^{15} t^7+q^{16} t^7-q^8 t^8-q^9 t^8-2 q^{10} t^8-3 q^{11} t^8
-5 q^{12} t^8-2 q^{13} t^8+q^{14} t^8+10 q^{15} t^8+3 q^{16} t^8
-q^9 t^9-q^{10} t^9-2 q^{11} t^9-3 q^{12} t^9-5 q^{13} t^9
+5 q^{15} t^9+6 q^{16} t^9+q^{17} t^9-q^{10} t^{10}-q^{11} t^{10}
-2 q^{12} t^{10}-3 q^{13} t^{10}-4 q^{14} t^{10}+3 q^{15} t^{10}
+6 q^{16} t^{10}+2 q^{17} t^{10}-q^{11} t^{11}-q^{12} t^{11}
-2 q^{13} t^{11}-3 q^{14} t^{11}-2 q^{15} t^{11}+6 q^{16} t^{11}
+3 q^{17} t^{11}-q^{12} t^{12}-q^{13} t^{12}-2 q^{14} t^{12}
-2 q^{15} t^{12}+q^{16} t^{12}+5 q^{17} t^{12}-q^{13} t^{13}
-q^{14} t^{13}-2 q^{15} t^{13}+3 q^{17} t^{13}+q^{18} t^{13}
-q^{14} t^{14}-q^{15} t^{14}-q^{16} t^{14}+2 q^{17} t^{14}
+q^{18} t^{14}-q^{15} t^{15}-q^{16} t^{15}+q^{17} t^{15}
+q^{18} t^{15}
-q^{16} t^{16}+q^{18} t^{16}-q^{17} t^{17}+q^{18} t^{17}\bigr).
\)
}
\renewcommand{\baselinestretch}{1.0}

\medskip
{\em HOMFLY-PT polynomial.\,}
Here
$den=(1-q)$\, and\, \  $\hat{\h}{}^{min}_{\l}(q,q,-a)=$
\renewcommand{\baselinestretch}{0.5} 

\noindent
{\small
\(
1-q+q^2+q^4+2 q^6-q^7+3 q^8-q^9+3 q^{10}-q^{11}+4 q^{12}-3 q^{13}
+5 q^{14}-3 q^{15}+4 q^{16}-4 q^{17}+5 q^{18}-4 q^{19}+4 q^{20}
-3 q^{21}+5 q^{22}-3 q^{23}+4 q^{24}-q^{25}+3 q^{26}-q^{27}
+3 q^{28}-q^{29}+2 q^{30}+q^{32}+q^{34}-q^{35}+q^{36}
-a^5 \bigl(q^{15}-q^{16}+q^{17}-q^{18}+q^{19}-q^{20}+q^{21}\bigr)
-a^4 \bigl(-q^{10}-q^{12}-q^{13}-q^{14}-2 q^{16}+q^{17}-q^{18}
+q^{19}-2 q^{20}-q^{22}-q^{23}-q^{24}-q^{26}\bigr)
-a^3 \bigl(q^6+2 q^8+q^9+3 q^{10}+q^{11}+5 q^{12}+6 q^{14}-q^{15}
+5 q^{16}-3 q^{17}+5 q^{18}-3 q^{19}+5 q^{20}-q^{21}+6 q^{22}
+5 q^{24}+q^{25}+3 q^{26}+q^{27}+2 q^{28}+q^{30}\bigr)
-a^2 \bigl(-q^3-2 q^5-q^6-3 q^7-2 q^8-5 q^9-2 q^{10}-7 q^{11}
-2 q^{12}-7 q^{13}-q^{14}-7 q^{15}+2 q^{16}-6 q^{17}+3 q^{18}
-6 q^{19}+2 q^{20}-7 q^{21}-q^{22}-7 q^{23}-2 q^{24}-7 q^{25}
-2 q^{26}-5 q^{27}-2 q^{28}-3 q^{29}-q^{30}-2 q^{31}-q^{33}\bigr)
-a \bigl(q+q^3+q^4+2 q^5+q^6+4 q^7+q^8+5 q^9+2 q^{10}+5 q^{11}
+q^{12}+7 q^{13}-q^{14}+6 q^{15}-2 q^{16}+5 q^{17}-3 q^{18}
+5 q^{19}-2 q^{20}+6 q^{21}-q^{22}+7 q^{23}+q^{24}+5 q^{25}+2 q^{26}
+5 q^{27}+q^{28}
+4 q^{29}+q^{30}+2 q^{31}+q^{32}+q^{33}+q^{35}\bigr).
\)
}
\renewcommand{\baselinestretch}{1.0} 
\medskip

{\em Khovanov polynomial.\,}
Let us apply the procedure from (\ref{khfromh})
to obtain the DAHA approximation to the reduced
Khovanov polynomial ($n=1$). Recall that 
we switch to the standard parameters:
$q\mapsto (qt)^2,\, t\mapsto q^2,\,a\mapsto -q^{2(n+1)}$;
$st$ is omitted in $q_{st},t_{st},a_{st}$.
Then
\begin{align}\label{khfromh8-6}
\Bigl(\frac{(1\!+\!a)\,\hat{\h}^{min}_{\l}}
{(1-t)^{2}}\Bigr)_{st, n=1}=
\end{align}
{\small \(
1+q^2+q^4 t^2-q^8 t^2+q^6 t^4+q^8 t^4-q^{10} t^4-q^{12} t^4
+q^8 t^6+q^{10} t^6-q^{12} t^6-q^{14} t^6+q^{10} t^8+2 q^{12} t^8
-q^{14} t^8-2 q^{16} t^8+q^{12} t^{10}+2 q^{14} t^{10}-q^{16} t^{10}
-3 q^{18} t^{10}+q^{22} t^{10}+3 q^{16} t^{12}+q^{18} t^{12}
-4 q^{20} t^{12}-q^{22} t^{12}+q^{24} t^{12}+q^{18} t^{14}
+2 q^{20} t^{14}-3 q^{22} t^{14}-3 q^{24} t^{14}+2 q^{26} t^{14}
+q^{28} t^{14}+q^{20} t^{16}+2 q^{22} t^{16}-q^{24} t^{16}
-3 q^{26} t^{16}+q^{30} t^{16}+2 q^{24} t^{18}-3 q^{28} t^{18}
+q^{32} t^{18}+2 q^{28} t^{20}-q^{30} t^{20}-2 q^{32} t^{20}
+q^{34} t^{20}+q^{32} t^{22}-q^{34} t^{22}
-q^{36} t^{22}+q^{38} t^{22}+q^{36} t^{24}+q^{38} t^{24}
,\)}
\begin{align}\label{khfromh8-6-kh}
(\hat{K\!h\!R}{}^{2}_{\l})_{reduced}=
\end{align}
{\small \(
1+q^2+q^4 t^2+q^8 t^3+q^6 t^4+q^8 t^4+q^{10} t^5+q^{12} t^5
+q^8 t^6+q^{10} t^6+q^{12} t^7+q^{14} t^7+q^{10} t^8+2 q^{12} t^8
+q^{14} t^9+2 q^{16} t^9+q^{12} t^{10}+2 q^{14} t^{10}+q^{16} t^{11}
+3 q^{18} t^{11}+3 q^{16} t^{12}+q^{18} t^{12}+q^{22} t^{12}
+4 q^{20} t^{13}+q^{22} t^{13}+q^{18} t^{14}+2 q^{20} t^{14}
+q^{24} t^{14}+3 q^{22} t^{15}+3 q^{24} t^{15}+q^{20} t^{16}
+2 q^{22} t^{16}+2 q^{26} t^{16}+q^{28} t^{16}+q^{24} t^{17}
+3 q^{26} t^{17}+2 q^{24} t^{18}+q^{30} t^{18}+3 q^{28} t^{19}
+2 q^{28} t^{20}+q^{32} t^{20}+q^{30} t^{21}+2 q^{32} t^{21}
+q^{32} t^{22}+q^{34} t^{22}
+q^{34} t^{23}+q^{36} t^{23}+q^{36} t^{24}+2 q^{38} t^{24}
.\)}
\medskip

Their difference (the second minus the first)
is large, but still a sum of the terms $q^l(t^m+t^{m-1})$,
$q^l(t^m-t^{m-2})$ with positive coefficients, which confirms
the Connection Conjecture:
\medskip

\noindent
{\small
\(
\bigl(q^8 t^2+q^8 t^3\bigr)+\bigl(q^{10} t^4+q^{10} t^5\bigr)
+\bigl(q^{12} t^4+q^{12} t^5\bigr)+\bigl(q^{12} t^6
+q^{12} t^7\bigr)+\bigl(q^{14} t^6+q^{14} t^7\bigr)
+\bigl(q^{14} t^8+q^{14} t^9\bigr)+\bigl(2 q^{16} t^8+
2 q^{16} t^9\bigr)+\bigl(q^{16} t^{10}+q^{16} t^{11}\bigr)
+\bigl(3 q^{18} t^{10}+3 q^{18} t^{11}\bigr)
+\bigl(-q^{22} t^{10}+2 q^{22} t^{12}+q^{22} t^{13}\bigr)
+\bigl(4 q^{20} t^{12}+4 q^{20} t^{13}\bigr)
+\bigl(-q^{24} t^{12}+4 q^{24} t^{14}
+3 q^{24} t^{15}\bigr)
+\bigl(3 q^{22} t^{14}+3 q^{22} t^{15}\bigr)+
\bigl(q^{24} t^{16}+q^{24} t^{17}\bigr)
+\bigl(-2 q^{26} t^{14}+5 q^{26} t^{16}+3 q^{26} t^{17}\bigr)
+\bigl(-q^{28} t^{14}+q^{28} t^{16}\bigr)
+\bigl(-q^{30} t^{16}+q^{30} t^{18}\bigr)
+\bigl(3 q^{28} t^{18}+3 q^{28} t^{19}\bigr)+\bigl(-q^{32} t^{18}
+q^{32} t^{20}\bigr)+\bigl(q^{30} t^{20}+q^{30} t^{21}\bigr)
+\bigl(2q^{32} t^{20}+2 q^{32} t^{21}\bigr)+\bigl(-q^{34} t^{20}
+q^{34} t^{22}\bigr)+ \bigl(q^{34} t^{22}+q^{34} t^{23}\bigr)
+\bigl(q^{36} t^{22}+q^{36} t^{23}\bigr)
+\bigl(-q^{38} t^{22}+q^{38} t^{24}\bigr).
\)
}
\smallskip

\subsection{\bf Uncolored 2-fold 
\texorpdfstring{{\mathversion{bold}$T(5,3)$}}
{T(5,3)}} 
\begin{align}\label{T10-6}
&T(10,6):\ 
\l=\l_{(\{5,3\},\{5,3\})}^{\,\circ\rightrightarrows,\, 
(\square,\square)},\ \ 
\hat{\h}{}^{min}_{\l}\,(q,t,a)=
\end{align}

\renewcommand{\baselinestretch}{0.5} 
\noindent
{\small
\(
1-t+q t+q^2 t+q^3 t+q^4 t+q^5 t-q t^2+q^4 t^2+q^5 t^2
+3 q^6 t^2+2 q^7 t^2+2 q^8 t^2-q^2 t^3-q^4 t^3+2 q^7 t^3+2 q^8 t^3
+5 q^9 t^3+2 q^{10} t^3+q^{11} t^3-q^3 t^4-q^5 t^4-q^6 t^4-q^7 t^4
+q^8 t^4+5 q^{10} t^4+4 q^{11} t^4+4 q^{12} t^4-q^4 t^5-q^6 t^5
-q^7 t^5-2 q^8 t^5-q^{10} t^5+3 q^{11} t^5+3 q^{12} t^5+7 q^{13} t^5
+q^{14} t^5-q^5 t^6-q^7 t^6-q^8 t^6-2 q^9 t^6-q^{10} t^6
-2 q^{11} t^6+2 q^{12} t^6+q^{13} t^6+8 q^{14} t^6+2 q^{15} t^6
-q^6 t^7-q^8 t^7-q^9 t^7-2 q^{10} t^7-q^{11} t^7-3 q^{12} t^7
+q^{13} t^7+8 q^{15} t^7+2 q^{16} t^7-q^7 t^8-q^9 t^8-q^{10} t^8
-2 q^{11} t^8-q^{12} t^8-3 q^{13} t^8+8 q^{16} t^8+2 q^{17} t^8
-q^8 t^9-q^{10} t^9-q^{11} t^9-2 q^{12} t^9-q^{13} t^9-3 q^{14} t^9
+8 q^{17} t^9+q^{18} t^9-q^9 t^{10}-q^{11} t^{10}-q^{12} t^{10}
-2 q^{13} t^{10}-q^{14} t^{10}-3 q^{15} t^{10}+q^{16} t^{10}
+q^{17} t^{10}+7 q^{18} t^{10}-q^{10} t^{11}-q^{12} t^{11}
-q^{13} t^{11}-2 q^{14} t^{11}-q^{15} t^{11}-3 q^{16} t^{11}
+2 q^{17} t^{11}+3 q^{18} t^{11}+4 q^{19} t^{11}-q^{11} t^{12}
-q^{13} t^{12}-q^{14} t^{12}-2 q^{15} t^{12}-q^{16} t^{12}
-2 q^{17} t^{12}+3 q^{18} t^{12}+4 q^{19} t^{12}+q^{20} t^{12}
-q^{12} t^{13}-q^{14} t^{13}-q^{15} t^{13}-2 q^{16} t^{13}
-q^{17} t^{13}-q^{18} t^{13}+5 q^{19} t^{13}+2 q^{20} t^{13}
-q^{13} t^{14}-q^{15} t^{14}-q^{16} t^{14}-2 q^{17} t^{14}
+5 q^{20} t^{14}-q^{14} t^{15}-q^{16} t^{15}-q^{17} t^{15}
-2 q^{18} t^{15}+q^{19} t^{15}+2 q^{20} t^{15}+2 q^{21} t^{15}
-q^{15} t^{16}-q^{17} t^{16}-q^{18} t^{16}-q^{19} t^{16}
+2 q^{20} t^{16}+2 q^{21} t^{16}-q^{16} t^{17}-q^{18} t^{17}
-q^{19} t^{17}+3 q^{21} t^{17}-q^{17} t^{18}-q^{19} t^{18}
+q^{21} t^{18}+q^{22} t^{18}-q^{18} t^{19}-q^{20} t^{19}
+q^{21} t^{19}+q^{22} t^{19}-q^{19} t^{20}+q^{22} t^{20}
-q^{20} t^{21}+q^{22} t^{21}-q^{21} t^{22}+q^{22} t^{22}
-q^{22} t^{23}+q^{23} t^{23}
\)

\vfil
\smallskip\noindent
\(
+a^5 \bigl(q^{15}-q^{15} t+q^{16} t
+q^{17} t+q^{18} t-q^{16} t^2+2 q^{19} t^2-q^{17} t^3-q^{19} t^3
+2 q^{20} t^3-q^{18} t^4-q^{20} t^4+2 q^{21} t^4-q^{19} t^5
+q^{22} t^5-q^{20} t^6+q^{22} t^6-q^{21} t^7+q^{22} t^7-q^{22} t^8
+q^{23} t^8\bigr)
\)

\vfil
\smallskip\noindent
\(
+a^4 \bigl(q^{10}+q^{11}+q^{12}+q^{13}+q^{14}
-q^{10} t+q^{12} t+2 q^{13} t+3 q^{14} t+4 q^{15} t+2 q^{16} t
+q^{17} t-q^{11} t^2-q^{12} t^2-q^{13} t^2+2 q^{15} t^2+5 q^{16} t^2
+4 q^{17} t^2+3 q^{18} t^2-q^{12} t^3-q^{13} t^3-2 q^{14} t^3
-2 q^{15} t^3+4 q^{17} t^3+4 q^{18} t^3+4 q^{19} t^3-q^{13} t^4
-q^{14} t^4-2 q^{15} t^4-3 q^{16} t^4-2 q^{17} t^4+3 q^{18} t^4
+4 q^{19} t^4+4 q^{20} t^4-q^{14} t^5-q^{15} t^5-2 q^{16} t^5
-3 q^{17} t^5-3 q^{18} t^5+3 q^{19} t^5+4 q^{20} t^5+3 q^{21} t^5
-q^{15} t^6-q^{16} t^6-2 q^{17} t^6-3 q^{18} t^6-2 q^{19} t^6
+4 q^{20} t^6+4 q^{21} t^6+q^{22} t^6-q^{16} t^7-q^{17} t^7
-2 q^{18} t^7-3 q^{19} t^7+5 q^{21} t^7+2 q^{22} t^7-q^{17} t^8
-q^{18} t^8-2 q^{19} t^8-2 q^{20} t^8+2 q^{21} t^8+4 q^{22} t^8
-q^{18} t^9-q^{19} t^9-2 q^{20} t^9+3 q^{22} t^9+q^{23} t^9
-q^{19} t^{10}-q^{20} t^{10}-q^{21} t^{10}+2 q^{22} t^{10}
+q^{23} t^{10}-q^{20} t^{11}-q^{21} t^{11}+q^{22} t^{11}
+q^{23} t^{11}-q^{21} t^{12}+q^{23} t^{12}-q^{22} t^{13}
+q^{23} t^{13}\bigr)
\)

\vfil
\smallskip\noindent
\(
+a^3 \bigl(q^6+q^7+2 q^8+2 q^9+2 q^{10}
+q^{11}+q^{12}-q^6 t+2 q^9 t+4 q^{10} t+7 q^{11} t+6 q^{12} t
+6 q^{13} t+3 q^{14} t+q^{15} t-q^7 t^2-q^8 t^2-2 q^9 t^2
-q^{10} t^2+6 q^{12} t^2+8 q^{13} t^2+12 q^{14} t^2+8 q^{15} t^2
+5 q^{16} t^2+q^{17} t^2-q^8 t^3-q^9 t^3-3 q^{10} t^3-3 q^{11} t^3
-4 q^{12} t^3+q^{13} t^3+4 q^{14} t^3+14 q^{15} t^3+11 q^{16} t^3
+8 q^{17} t^3+2 q^{18} t^3-q^9 t^4-q^{10} t^4-3 q^{11} t^4
-4 q^{12} t^4-6 q^{13} t^4-3 q^{14} t^4-2 q^{15} t^4+12 q^{16} t^4
+12 q^{17} t^4+10 q^{18} t^4+2 q^{19} t^4-q^{10} t^5-q^{11} t^5
-3 q^{12} t^5-4 q^{13} t^5-7 q^{14} t^5-5 q^{15} t^5-6 q^{16} t^5
+9 q^{17} t^5+12 q^{18} t^5+10 q^{19} t^5+2 q^{20} t^5-q^{11} t^6
-q^{12} t^6-3 q^{13} t^6-4 q^{14} t^6-7 q^{15} t^6-6 q^{16} t^6
-7 q^{17} t^6+9 q^{18} t^6+12 q^{19} t^6+8 q^{20} t^6+q^{21} t^6
-q^{12} t^7-q^{13} t^7-3 q^{14} t^7-4 q^{15} t^7-7 q^{16} t^7
-6 q^{17} t^7-6 q^{18} t^7+12 q^{19} t^7+11 q^{20} t^7+5 q^{21} t^7
-q^{13} t^8-q^{14} t^8-3 q^{15} t^8-4 q^{16} t^8-7 q^{17} t^8
-5 q^{18} t^8-2 q^{19} t^8+14 q^{20} t^8+8 q^{21} t^8+q^{22} t^8
-q^{14} t^9-q^{15} t^9-3 q^{16} t^9-4 q^{17} t^9-7 q^{18} t^9
-3 q^{19} t^9+4 q^{20} t^9+12 q^{21} t^9+3 q^{22} t^9-q^{15} t^{10}
-q^{16} t^{10}-3 q^{17} t^{10}-4 q^{18} t^{10}-6 q^{19} t^{10}
+q^{20} t^{10}+8 q^{21} t^{10}+6 q^{22} t^{10}-q^{16} t^{11}
-q^{17} t^{11}-3 q^{18} t^{11}-4 q^{19} t^{11}-4 q^{20} t^{11}
+6 q^{21} t^{11}+6 q^{22} t^{11}+q^{23} t^{11}-q^{17} t^{12}
-q^{18} t^{12}-3 q^{19} t^{12}-3 q^{20} t^{12}+7 q^{22} t^{12}
+q^{23} t^{12}-q^{18} t^{13}-q^{19} t^{13}-3 q^{20} t^{13}
-q^{21} t^{13}+4 q^{22} t^{13}+2 q^{23} t^{13}-q^{19} t^{14}
-q^{20} t^{14}-2 q^{21} t^{14}+2 q^{22} t^{14}+2 q^{23} t^{14}
-q^{20} t^{15}-q^{21} t^{15}+2 q^{23} t^{15}-q^{21} t^{16}
+q^{23} t^{16}-q^{22} t^{17}+q^{23} t^{17}\bigr)
\)

\vfil
\smallskip\noindent
\(
+a^2 \bigl(q^3
+q^4+2 q^5+2 q^6+2 q^7+q^8+q^9-q^3 t+2 q^6 t+4 q^7 t+7 q^8 t
+7 q^9 t+7 q^{10} t+4 q^{11} t+2 q^{12} t-q^4 t^2-q^5 t^2-2 q^6 t^2
-q^7 t^2+5 q^9 t^2+9 q^{10} t^2+14 q^{11} t^2+12 q^{12} t^2
+10 q^{13} t^2+3 q^{14} t^2+q^{15} t^2-q^5 t^3-q^6 t^3-3 q^7 t^3
-3 q^8 t^3-4 q^9 t^3-q^{10} t^3+3 q^{11} t^3+13 q^{12} t^3
+16 q^{13} t^3+18 q^{14} t^3+8 q^{15} t^3+3 q^{16} t^3-q^6 t^4
-q^7 t^4-3 q^8 t^4-4 q^9 t^4-6 q^{10} t^4-5 q^{11} t^4-4 q^{12} t^4
+6 q^{13} t^4+14 q^{14} t^4+23 q^{15} t^4+12 q^{16} t^4
+5 q^{17} t^4-q^7 t^5-q^8 t^5-3 q^9 t^5-4 q^{10} t^5-7 q^{11} t^5
-7 q^{12} t^5-8 q^{13} t^5-q^{14} t^5+8 q^{15} t^5+24 q^{16} t^5
+14 q^{17} t^5+6 q^{18} t^5-q^8 t^6-q^9 t^6-3 q^{10} t^6
-4 q^{11} t^6-7 q^{12} t^6-8 q^{13} t^6-10 q^{14} t^6-5 q^{15} t^6
+4 q^{16} t^6+23 q^{17} t^6+14 q^{18} t^6+5 q^{19} t^6-q^9 t^7
-q^{10} t^7-3 q^{11} t^7-4 q^{12} t^7-7 q^{13} t^7-8 q^{14} t^7
-11 q^{15} t^7-6 q^{16} t^7+4 q^{17} t^7+24 q^{18} t^7+12 q^{19} t^7
+3 q^{20} t^7-q^{10} t^8-q^{11} t^8-3 q^{12} t^8-4 q^{13} t^8
-7 q^{14} t^8-8 q^{15} t^8-11 q^{16} t^8-5 q^{17} t^8+8 q^{18} t^8
+23 q^{19} t^8+8 q^{20} t^8+q^{21} t^8-q^{11} t^9-q^{12} t^9
-3 q^{13} t^9-4 q^{14} t^9-7 q^{15} t^9-8 q^{16} t^9-10 q^{17} t^9
-q^{18} t^9+14 q^{19} t^9+18 q^{20} t^9+3 q^{21} t^9-q^{12} t^{10}
-q^{13} t^{10}-3 q^{14} t^{10}-4 q^{15} t^{10}-7 q^{16} t^{10}
-8 q^{17} t^{10}-8 q^{18} t^{10}+6 q^{19} t^{10}+16 q^{20} t^{10}
+10 q^{21} t^{10}-q^{13} t^{11}-q^{14} t^{11}-3 q^{15} t^{11}
-4 q^{16} t^{11}-7 q^{17} t^{11}-7 q^{18} t^{11}-4 q^{19} t^{11}
+13 q^{20} t^{11}+12 q^{21} t^{11}+2 q^{22} t^{11}-q^{14} t^{12}
-q^{15} t^{12}-3 q^{16} t^{12}-4 q^{17} t^{12}-7 q^{18} t^{12}
-5 q^{19} t^{12}+3 q^{20} t^{12}+14 q^{21} t^{12}+4 q^{22} t^{12}
-q^{15} t^{13}-q^{16} t^{13}-3 q^{17} t^{13}-4 q^{18} t^{13}
-6 q^{19} t^{13}-q^{20} t^{13}+9 q^{21} t^{13}+7 q^{22} t^{13}
-q^{16} t^{14}-q^{17} t^{14}-3 q^{18} t^{14}-4 q^{19} t^{14}
-4 q^{20} t^{14}+5 q^{21} t^{14}+7 q^{22} t^{14}+q^{23} t^{14}
-q^{17} t^{15}-q^{18} t^{15}-3 q^{19} t^{15}-3 q^{20} t^{15}
+7 q^{22} t^{15}+q^{23} t^{15}-q^{18} t^{16}-q^{19} t^{16}
-3 q^{20} t^{16}-q^{21} t^{16}+4 q^{22} t^{16}+2 q^{23} t^{16}
-q^{19} t^{17}-q^{20} t^{17}-2 q^{21} t^{17}+2 q^{22} t^{17}
+2 q^{23} t^{17}-q^{20} t^{18}-q^{21} t^{18}+2 q^{23} t^{18}
-q^{21} t^{19}+q^{23} t^{19}-q^{22} t^{20}+q^{23} t^{20}\bigr)
\)

\vfil
\smallskip\noindent
\(+a \bigl(q+q^2+q^3+q^4+q^5-q t+q^3 t+2 q^4 t+3 q^5 t+5 q^6 t
+4 q^7 t+3 q^8 t+q^9 t-q^2 t^2-q^3 t^2-q^4 t^2+q^6 t^2+5 q^7 t^2
+7 q^8 t^2+10 q^9 t^2+7 q^{10} t^2+4 q^{11} t^2+q^{12} t^2-q^3 t^3
-q^4 t^3-2 q^5 t^3-2 q^6 t^3-2 q^7 t^3+q^8 t^3+3 q^9 t^3
+11 q^{10} t^3+12 q^{11} t^3+11 q^{12} t^3+5 q^{13} t^3+q^{14} t^3
-q^4 t^4-q^5 t^4-2 q^6 t^4-3 q^7 t^4-4 q^8 t^4-2 q^9 t^4
-2 q^{10} t^4+6 q^{11} t^4+11 q^{12} t^4+17 q^{13} t^4+10 q^{14} t^4
+3 q^{15} t^4-q^5 t^5-q^6 t^5-2 q^7 t^5-3 q^8 t^5-5 q^9 t^5
-4 q^{10} t^5-5 q^{11} t^5+q^{12} t^5+5 q^{13} t^5+18 q^{14} t^5
+14 q^{15} t^5+5 q^{16} t^5-q^6 t^6-q^7 t^6-2 q^8 t^6-3 q^9 t^6
-5 q^{10} t^6-5 q^{11} t^6-7 q^{12} t^6-2 q^{13} t^6+16 q^{15} t^6
+15 q^{16} t^6+6 q^{17} t^6-q^7 t^7-q^8 t^7-2 q^9 t^7-3 q^{10} t^7
-5 q^{11} t^7-5 q^{12} t^7-8 q^{13} t^7-4 q^{14} t^7-2 q^{15} t^7
+15 q^{16} t^7+15 q^{17} t^7+5 q^{18} t^7-q^8 t^8-q^9 t^8
-2 q^{10} t^8-3 q^{11} t^8-5 q^{12} t^8-5 q^{13} t^8-8 q^{14} t^8
-5 q^{15} t^8-2 q^{16} t^8+16 q^{17} t^8+14 q^{18} t^8+3 q^{19} t^8
-q^9 t^9-q^{10} t^9-2 q^{11} t^9-3 q^{12} t^9-5 q^{13} t^9
-5 q^{14} t^9-8 q^{15} t^9-4 q^{16} t^9+18 q^{18} t^9+10 q^{19} t^9
+q^{20} t^9-q^{10} t^{10}-q^{11} t^{10}-2 q^{12} t^{10}
-3 q^{13} t^{10}-5 q^{14} t^{10}-5 q^{15} t^{10}-8 q^{16} t^{10}
-2 q^{17} t^{10}+5 q^{18} t^{10}+17 q^{19} t^{10}+5 q^{20} t^{10}
-q^{11} t^{11}-q^{12} t^{11}-2 q^{13} t^{11}-3 q^{14} t^{11}
-5 q^{15} t^{11}-5 q^{16} t^{11}-7 q^{17} t^{11}+q^{18} t^{11}
+11 q^{19} t^{11}+11 q^{20} t^{11}+q^{21} t^{11}-q^{12} t^{12}
-q^{13} t^{12}-2 q^{14} t^{12}-3 q^{15} t^{12}-5 q^{16} t^{12}
-5 q^{17} t^{12}-5 q^{18} t^{12}+6 q^{19} t^{12}+12 q^{20} t^{12}
+4 q^{21} t^{12}-q^{13} t^{13}-q^{14} t^{13}-2 q^{15} t^{13}
-3 q^{16} t^{13}-5 q^{17} t^{13}-4 q^{18} t^{13}-2 q^{19} t^{13}
+11 q^{20} t^{13}+7 q^{21} t^{13}-q^{14} t^{14}-q^{15} t^{14}
-2 q^{16} t^{14}-3 q^{17} t^{14}-5 q^{18} t^{14}-2 q^{19} t^{14}
+3 q^{20} t^{14}+10 q^{21} t^{14}+q^{22} t^{14}-q^{15} t^{15}
-q^{16} t^{15}-2 q^{17} t^{15}-3 q^{18} t^{15}-4 q^{19} t^{15}
+q^{20} t^{15}+7 q^{21} t^{15}+3 q^{22} t^{15}-q^{16} t^{16}
-q^{17} t^{16}-2 q^{18} t^{16}-3 q^{19} t^{16}-2 q^{20} t^{16}
+5 q^{21} t^{16}+4 q^{22} t^{16}-q^{17} t^{17}-q^{18} t^{17}
-2 q^{19} t^{17}-2 q^{20} t^{17}+q^{21} t^{17}+5 q^{22} t^{17}
-q^{18} t^{18}-q^{19} t^{18}-2 q^{20} t^{18}+3 q^{22} t^{18}
+q^{23} t^{18}-q^{19} t^{19}-q^{20} t^{19}-q^{21} t^{19}
+2 q^{22} t^{19}+q^{23} t^{19}-q^{20} t^{20}-q^{21} t^{20}
+q^{22} t^{20}+q^{23} t^{20}
-q^{21} t^{21}+q^{23} t^{21}-q^{22} t^{22}+q^{23} t^{22}\bigr)
\)
}
\renewcommand{\baselinestretch}{1.0} 
\smallskip

The $a$\~degree of $\hat{\h}{}^{min}_{\l}\,(q,t,a)$
is $5$, the same as that for the $2$-fold $T(4,3).$
The positivity of the series
$\hat{\h}{}^{min}_{\l}(q,t,a)/(1-t)$  holds. 
\smallskip

\subsection{\bf Uncolored 3-links}
We will begin with the case of the uncolored 3-folded trefoil
defined as follows.
\smallskip
  
{\sf 3-fold trefoil.} The corresponding \tax-presentation  
and $\hat{\h}$ are: 
$$
1\le j\le \kappa=3,\ \vec\rr^j=3,\, \vec\ss^j=2,\, 
\Up=\{\circ\rightthreearrow\}\,,\ 
\la^1=\square=\la^2=\la^3;
$$

\begin{align}\label{T9-6}
&T(9,6):\ 
\l=\l_{(\{3,2\},\{3,2\},\{3,2\})}^{\,\circ\rightthreearrow,\, 
(\square,\square,\square)},\ \ 
\hat{\h}{}^{min}_{\l}\,(q,t,a)=
\end{align}

\renewcommand{\baselinestretch}{0.5} 
\noindent
{\small
\(
1-2 t+q t+q^2 t+q^3 t+q^4 t+q^5 t+t^2-2 q t^2-q^2 t^2
-q^3 t^2+3 q^6 t^2+2 q^7 t^2+q^8 t^2+q t^3-q^2 t^3-2 q^4 t^3-q^5 t^3
-3 q^6 t^3+2 q^8 t^3+4 q^9 t^3+q^{10} t^3+q^2 t^4-q^3 t^4+q^4 t^4
-q^5 t^4-q^6 t^4-3 q^7 t^4-2 q^8 t^4-3 q^9 t^4+4 q^{10} t^4
+4 q^{11} t^4+q^{12} t^4+q^3 t^5-q^4 t^5+q^5 t^5-3 q^8 t^5
-q^9 t^5-5 q^{10} t^5-q^{11} t^5+6 q^{12} t^5+3 q^{13} t^5
+q^4 t^6-q^5 t^6+q^6 t^6+q^8 t^6-2 q^9 t^6-q^{10} t^6-4 q^{11} t^6
-4 q^{12} t^6+5 q^{13} t^6+3 q^{14} t^6+q^{15} t^6+q^5 t^7-q^6 t^7
+q^7 t^7+q^9 t^7-q^{10} t^7-5 q^{12} t^7-5 q^{13} t^7+6 q^{14} t^7
+3 q^{15} t^7+q^6 t^8-q^7 t^8+q^8 t^8+q^{10} t^8-q^{11} t^8
+q^{12} t^8-5 q^{13} t^8-5 q^{14} t^8+5 q^{15} t^8+3 q^{16} t^8
+q^7 t^9-q^8 t^9+q^9 t^9+q^{11} t^9-q^{12} t^9+q^{13} t^9
-5 q^{14} t^9-4 q^{15} t^9+6 q^{16} t^9+q^{17} t^9+q^8 t^{10}
-q^9 t^{10}+q^{10} t^{10}+q^{12} t^{10}-q^{13} t^{10}
-4 q^{15} t^{10}-q^{16} t^{10}+4 q^{17} t^{10}+q^9 t^{11}
-q^{10} t^{11}+q^{11} t^{11}+q^{13} t^{11}-q^{14} t^{11}
-q^{15} t^{11}-5 q^{16} t^{11}+4 q^{17} t^{11}+q^{18} t^{11}
+q^{10} t^{12}-q^{11} t^{12}+q^{12} t^{12}+q^{14} t^{12}
-2 q^{15} t^{12}-q^{16} t^{12}-3 q^{17} t^{12}+4 q^{18} t^{12}
+q^{11} t^{13}-q^{12} t^{13}+q^{13} t^{13}+q^{15} t^{13}
-3 q^{16} t^{13}-2 q^{17} t^{13}+2 q^{18} t^{13}+q^{19} t^{13}
+q^{12} t^{14}-q^{13} t^{14}+q^{14} t^{14}-3 q^{17} t^{14}
+2 q^{19} t^{14}+q^{13} t^{15}-q^{14} t^{15}+q^{15} t^{15}
-q^{17} t^{15}-3 q^{18} t^{15}+3 q^{19} t^{15}+q^{14} t^{16}
-q^{15} t^{16}+q^{16} t^{16}-q^{17} t^{16}-q^{18} t^{16}
+q^{20} t^{16}+q^{15} t^{17}-q^{16} t^{17}+q^{17} t^{17}
-2 q^{18} t^{17}+q^{20} t^{17}+q^{16} t^{18}-q^{17} t^{18}
-q^{19} t^{18}+q^{20} t^{18}+q^{17} t^{19}-q^{18} t^{19}
-q^{19} t^{19}+q^{20} t^{19}+q^{18} t^{20}-2 q^{19} t^{20}
+q^{20} t^{20}+q^{19} t^{21}-2 q^{20} t^{21}+q^{21} t^{21}
+a^5 \bigl(q^{15}-2 q^{15} t+q^{16} t+q^{17} t+q^{15} t^2
-2 q^{16} t^2-q^{17} t^2+q^{18} t^2+q^{19} t^2+q^{16} t^3
-q^{17} t^3-q^{18} t^3+q^{19} t^3+q^{17} t^4-q^{18} t^4
-q^{19} t^4+q^{20} t^4+q^{18} t^5-2 q^{19} t^5+q^{20} t^5
+q^{19} t^6-2 q^{20} t^6+q^{21} t^6\bigr)
\)

\vfil
\smallskip\noindent
\(
+a^4 \bigl(q^{10}+q^{11}+q^{12}+q^{13}+q^{14}-2 q^{10} t-q^{11} t
+q^{13} t+q^{14} t+3 q^{15} t+q^{16} t+q^{10} t^2-q^{11} t^2
-2 q^{12} t^2-3 q^{13} t^2-q^{14} t^2-q^{15} t^2+3 q^{16} t^2
+3 q^{17} t^2+q^{18} t^2+q^{11} t^3-3 q^{14} t^3-3 q^{15} t^3
-q^{16} t^3+2 q^{17} t^3+3 q^{18} t^3+q^{19} t^3+q^{12} t^4
+q^{14} t^4-q^{15} t^4-4 q^{16} t^4-2 q^{17} t^4+2 q^{18} t^4
+3 q^{19} t^4+q^{13} t^5+q^{15} t^5-q^{16} t^5-4 q^{17} t^5
-q^{18} t^5+3 q^{19} t^5+q^{20} t^5+q^{14} t^6+q^{16} t^6
-q^{17} t^6-3 q^{18} t^6-q^{19} t^6+3 q^{20} t^6+q^{15} t^7
+q^{17} t^7-3 q^{18} t^7-q^{19} t^7+q^{20} t^7+q^{21} t^7
+q^{16} t^8-3 q^{19} t^8+q^{20} t^8+q^{21} t^8+q^{17} t^9
-2 q^{19} t^9+q^{21} t^9+q^{18} t^{10}-q^{19} t^{10}-q^{20} t^{10}
+q^{21} t^{10}+q^{19} t^{11}-2 q^{20} t^{11}+q^{21} t^{11}\bigr)
\)

\vfil
\smallskip\noindent
\(
+a^3 \bigl(q^6+q^7+2 q^8+2 q^9+2 q^{10}+q^{11}+q^{12}-2 q^6 t-q^7 t
-2 q^8 t+2 q^{10} t+5 q^{11} t+4 q^{12} t+4 q^{13} t+2 q^{14} t
+q^6 t^2-q^7 t^2-q^8 t^2-4 q^9 t^2-5 q^{10} t^2-5 q^{11} t^2
+q^{12} t^2+4 q^{13} t^2+6 q^{14} t^2+5 q^{15} t^2+2 q^{16} t^2
+q^7 t^3+q^9 t^3-2 q^{10} t^3-4 q^{11} t^3-9 q^{12} t^3-6 q^{13} t^3
+q^{14} t^3+7 q^{15} t^3+7 q^{16} t^3+3 q^{17} t^3+q^{18} t^3
+q^8 t^4+2 q^{10} t^4-q^{12} t^4-7 q^{13} t^4-10 q^{14} t^4
-3 q^{15} t^4+7 q^{16} t^4+8 q^{17} t^4+3 q^{18} t^4+q^9 t^5
+2 q^{11} t^5+q^{12} t^5+q^{13} t^5-5 q^{14} t^5-11 q^{15} t^5
-5 q^{16} t^5+7 q^{17} t^5+7 q^{18} t^5+2 q^{19} t^5+q^{10} t^6
+2 q^{12} t^6+q^{13} t^6+2 q^{14} t^6-4 q^{15} t^6-11 q^{16} t^6
-3 q^{17} t^6+7 q^{18} t^6+5 q^{19} t^6+q^{11} t^7+2 q^{13} t^7
+q^{14} t^7+2 q^{15} t^7-5 q^{16} t^7-10 q^{17} t^7+q^{18} t^7
+6 q^{19} t^7+2 q^{20} t^7+q^{12} t^8+2 q^{14} t^8+q^{15} t^8
+q^{16} t^8-7 q^{17} t^8-6 q^{18} t^8+4 q^{19} t^8+4 q^{20} t^8
+q^{13} t^9+2 q^{15} t^9+q^{16} t^9-q^{17} t^9-9 q^{18} t^9
+q^{19} t^9+4 q^{20} t^9+q^{21} t^9+q^{14} t^{10}+2 q^{16} t^{10}
-4 q^{18} t^{10}-5 q^{19} t^{10}+5 q^{20} t^{10}+q^{21} t^{10}
+q^{15} t^{11}+2 q^{17} t^{11}-2 q^{18} t^{11}-5 q^{19} t^{11}
+2 q^{20} t^{11}+2 q^{21} t^{11}+q^{16} t^{12}+q^{18} t^{12}
-4 q^{19} t^{12}+2 q^{21} t^{12}+q^{17} t^{13}-q^{19} t^{13}
-2 q^{20} t^{13}+2 q^{21} t^{13}+q^{18} t^{14}-q^{19} t^{14}
-q^{20} t^{14}+q^{21} t^{14}+q^{19} t^{15}-2 q^{20} t^{15}
+q^{21} t^{15}\bigr)
\)

\vfil
\smallskip\noindent
\(
+a^2 \bigl(q^3+q^4+2 q^5+2 q^6+2 q^7+q^8
+q^9-2 q^3 t-q^4 t-2 q^5 t+2 q^7 t+6 q^8 t+5 q^9 t+6 q^{10} t
+3 q^{11} t+q^{12} t+q^3 t^2-q^4 t^2-q^5 t^2-4 q^6 t^2-5 q^7 t^2
-7 q^8 t^2+3 q^{10} t^2+9 q^{11} t^2+8 q^{12} t^2+5 q^{13} t^2
+q^{14} t^2+q^4 t^3+q^6 t^3-2 q^7 t^3-3 q^8 t^3-10 q^9 t^3
-9 q^{10} t^3-7 q^{11} t^3+5 q^{12} t^3+11 q^{13} t^3+9 q^{14} t^3
+4 q^{15} t^3+q^{16} t^3+q^5 t^4+2 q^7 t^4-6 q^{10} t^4-9 q^{11} t^4
-16 q^{12} t^4-3 q^{13} t^4+12 q^{14} t^4+12 q^{15} t^4+6 q^{16} t^4
+q^{17} t^4+q^6 t^5+2 q^8 t^5+q^9 t^5+2 q^{10} t^5-3 q^{11} t^5
-4 q^{12} t^5-18 q^{13} t^5-11 q^{14} t^5+10 q^{15} t^5
+13 q^{16} t^5+6 q^{17} t^5+q^{18} t^5+q^7 t^6+2 q^9 t^6+q^{10} t^6
+3 q^{11} t^6-q^{12} t^6-2 q^{13} t^6-17 q^{14} t^6-13 q^{15} t^6
+10 q^{16} t^6+12 q^{17} t^6+4 q^{18} t^6+q^8 t^7+2 q^{10} t^7
+q^{11} t^7+3 q^{12} t^7-q^{14} t^7-17 q^{15} t^7-11 q^{16} t^7
+12 q^{17} t^7+9 q^{18} t^7+q^{19} t^7+q^9 t^8+2 q^{11} t^8
+q^{12} t^8+3 q^{13} t^8-2 q^{15} t^8-18 q^{16} t^8-3 q^{17} t^8
+11 q^{18} t^8+5 q^{19} t^8+q^{10} t^9+2 q^{12} t^9+q^{13} t^9
+3 q^{14} t^9-q^{15} t^9-4 q^{16} t^9-16 q^{17} t^9+5 q^{18} t^9
+8 q^{19} t^9+q^{20} t^9+q^{11} t^{10}+2 q^{13} t^{10}+q^{14} t^{10}
+3 q^{15} t^{10}-3 q^{16} t^{10}-9 q^{17} t^{10}-7 q^{18} t^{10}
+9 q^{19} t^{10}+3 q^{20} t^{10}+q^{12} t^{11}+2 q^{14} t^{11}
+q^{15} t^{11}+2 q^{16} t^{11}-6 q^{17} t^{11}-9 q^{18} t^{11}
+3 q^{19} t^{11}+6 q^{20} t^{11}+q^{13} t^{12}+2 q^{15} t^{12}
+q^{16} t^{12}-10 q^{18} t^{12}+5 q^{20} t^{12}+q^{21} t^{12}
+q^{14} t^{13}+2 q^{16} t^{13}-3 q^{18} t^{13}-7 q^{19} t^{13}
+6 q^{20} t^{13}+q^{21} t^{13}+q^{15} t^{14}+2 q^{17} t^{14}
-2 q^{18} t^{14}-5 q^{19} t^{14}+2 q^{20} t^{14}+2 q^{21} t^{14}
+q^{16} t^{15}+q^{18} t^{15}-4 q^{19} t^{15}+2 q^{21} t^{15}
+q^{17} t^{16}-q^{19} t^{16}-2 q^{20} t^{16}+2 q^{21} t^{16}
+q^{18} t^{17}-q^{19} t^{17}-q^{20} t^{17}+q^{21} t^{17}
+q^{19} t^{18}-2 q^{20} t^{18}+q^{21} t^{18}\bigr)
\)

\vfil
\smallskip\noindent
\(
+a \bigl(q+q^2+q^3+q^4+q^5-2 q t-q^2 t+q^4 t+2 q^5 t+5 q^6 t+4 q^7 t
+2 q^8 t+q^9 t+q t^2-q^2 t^2-2 q^3 t^2-3 q^4 t^2-3 q^5 t^2
-4 q^6 t^2+q^7 t^2+6 q^8 t^2+7 q^9 t^2+5 q^{10} t^2+2 q^{11} t^2
+q^2 t^3-2 q^5 t^3-3 q^6 t^3-7 q^7 t^3-7 q^8 t^3-3 q^9 t^3
+6 q^{10} t^3+9 q^{11} t^3+6 q^{12} t^3+2 q^{13} t^3+q^3 t^4
+q^5 t^4-q^7 t^4-5 q^8 t^4-7 q^9 t^4-11 q^{10} t^4-4 q^{11} t^4
+10 q^{12} t^4+11 q^{13} t^4+4 q^{14} t^4+q^{15} t^4+q^4 t^5+
q^6 t^5+q^7 t^5+q^8 t^5-3 q^9 t^5-4 q^{10} t^5-10 q^{11} t^5
-13 q^{12} t^5+5 q^{13} t^5+13 q^{14} t^5+7 q^{15} t^5+q^{16} t^5
+q^5 t^6+q^7 t^6+q^8 t^6+2 q^9 t^6-q^{10} t^6-2 q^{11} t^6
-8 q^{12} t^6-16 q^{13} t^6+2 q^{14} t^6+12 q^{15} t^6+7 q^{16} t^6
+q^{17} t^6+q^6 t^7+q^8 t^7+q^9 t^7+2 q^{10} t^7-8 q^{13} t^7
-16 q^{14} t^7+2 q^{15} t^7+13 q^{16} t^7+4 q^{17} t^7+q^7 t^8
+q^9 t^8+q^{10} t^8+2 q^{11} t^8+q^{13} t^8-8 q^{14} t^8
-16 q^{15} t^8+5 q^{16} t^8+11 q^{17} t^8+2 q^{18} t^8+q^8 t^9
+q^{10} t^9+q^{11} t^9+2 q^{12} t^9-8 q^{15} t^9-13 q^{16} t^9
+10 q^{17} t^9+6 q^{18} t^9+q^9 t^{10}+q^{11} t^{10}+q^{12} t^{10}
+2 q^{13} t^{10}-2 q^{15} t^{10}-10 q^{16} t^{10}-4 q^{17} t^{10}
+9 q^{18} t^{10}+2 q^{19} t^{10}+q^{10} t^{11}+q^{12} t^{11}
+q^{13} t^{11}+2 q^{14} t^{11}-q^{15} t^{11}-4 q^{16} t^{11}
-11 q^{17} t^{11}+6 q^{18} t^{11}+5 q^{19} t^{11}+q^{11} t^{12}
+q^{13} t^{12}+q^{14} t^{12}+2 q^{15} t^{12}-3 q^{16} t^{12}
-7 q^{17} t^{12}-3 q^{18} t^{12}+7 q^{19} t^{12}+q^{20} t^{12}
+q^{12} t^{13}+q^{14} t^{13}+q^{15} t^{13}+q^{16} t^{13}
-5 q^{17} t^{13}-7 q^{18} t^{13}+6 q^{19} t^{13}+2 q^{20} t^{13}
+q^{13} t^{14}+q^{15} t^{14}+q^{16} t^{14}-q^{17} t^{14}
-7 q^{18} t^{14}+q^{19} t^{14}+4 q^{20} t^{14}+q^{14} t^{15}
+q^{16} t^{15}-3 q^{18} t^{15}-4 q^{19} t^{15}+5 q^{20} t^{15}
+q^{15} t^{16}+q^{17} t^{16}-2 q^{18} t^{16}-3 q^{19} t^{16}
+2 q^{20} t^{16}+q^{21} t^{16}+q^{16} t^{17}-3 q^{19} t^{17}
+q^{20} t^{17}+q^{21} t^{17}+q^{17} t^{18}-2 q^{19} t^{18}
+q^{21} t^{18}+q^{18} t^{19}-q^{19} t^{19}-q^{20} t^{19}
+q^{21} t^{19}+q^{19} t^{20}-2 q^{20} t^{20}+q^{21} t^{20}\bigr).
\)
}
\renewcommand{\baselinestretch}{1.0} 
\smallskip

The $a$\~degree of $\hat{\h}{}^{min}_{\l}\,(q,t,a)$
is $5$, coinciding with that given by the formula
deg$_a=\ss(3|\la|)-|\la|=5$
from (\ref{deg-a-jj}). The self-duality and other claims
in this theorem hold. We note that the positivity claim from
Part $(ii)$ of the Connection Conjecture \ref{CONCONJ}
really requires here $\kappa-1=2$. Generally (\ref{posit-claim}) 
reads:
$$
\hat{\h}^{min}_{\,\l,\,'\!\l}\,(q,t,a)/(1-t)^{\kappa-1}\in
\Z_+[[q,t,a]].
$$
Taking here $(1-t)$ is insufficient for the
positivity; $(1-t)^2$ is necessary.
\smallskip

{\sf 3-fold T(2,1).} 
\begin{align}\label{T6-3}
&T(9,6):\ 
\l=\l_{(\{2,1\},\{2,1\},\{2,1\})}^{\,\circ\rightthreearrow,\, 
(\square,\square,\square)},\ \ 
\hat{\h}{}^{min}_{\l}\,(q,t,a)=
\end{align}

\renewcommand{\baselinestretch}{0.5} 
\noindent
{\small
\(
1-2 t+q t+q^2 t+t^2-2 q t^2-q^2 t^2+q^3 t^2+q^4 t^2+q t^3-q^2 t^3
-q^3 t^3+q^4 t^3+q^2 t^4-q^3 t^4-q^4 t^4+q^5 t^4+q^3 t^5-2 q^4 t^5
+q^5 t^5+q^4 t^6-2 q^5 t^6+q^6 t^6+a^2 \bigl(q^3-2 q^3 t+q^4 t
+q^5 t+q^3 t^2-2 q^4 t^2+q^5 t^2+q^4 t^3-2 q^5 t^3+q^6 t^3\bigr)
+a \bigl(q+q^2-2 q t-q^2 t+2 q^3 t+q^4 t+q t^2-q^2 t^2-3 q^3 t^2
+2 q^4 t^2+q^5 t^2+q^2 t^3-3 q^4 t^3+2 q^5 t^3+q^3 t^4-q^4 t^4
-q^5 t^4+q^6 t^4+q^4 t^5-2 q^5 t^5+q^6 t^5\bigr).
\)
}
\renewcommand{\baselinestretch}{1.0} 
\smallskip

Here again the division by $(1-t)$ is insufficient for the
positivity.
\medskip

\setcounter{equation}{0}
\section{\sc Hopf links}\label{sec:torus-hopf}
\subsection{\bf Basic constructions}
We will provide basic superpolynomials for
the Hopf links. They correspond to taking $\tau_-^{-1}$
inside the coinvariant, i.e. this is the
case of multiple $T(1,-1)$. 
These examples are directly related to the DAHA-vertex 
$\,C_{\mathbf{b}}^{c}/\langle \theta\mmu^\flat \rangle=
\p_{\mathbf{b}}^{\iota(c)}/\langle P^\circ_c,\,P^\circ_c\rangle$
from Corollary \ref{DAHA-vertex-Hopf}. Namely, 
\begin{align}\label{h-vertexz}
\hat{\h}^{\emptyset,\dag}
\!=\!\frac{\hat{\h}^{\emptyset}_{\{1,-1\}}
(\lla)}
{(\la^1)_{ev}^\dag\cdots (\la^{\kappa})_{ev}^\dag}=
\frac{(\la^1\!\vee\!\ldots\!\vee\!\la^{\kappa})_{ev}^\dag\,
\hat{\h}^{min}_{\{1,-1\}}(\lla)}
{(\la^1)_{ev}^\dag\cdots (\la^{\kappa})_{ev}^\dag}
\end{align}
is the stabilization of $\p_{\hbox{\tiny $\lla$}}=
\p_{(\la^1\!,\ldots,\la^{\kappa-1})}^{\la^\kappa}$ (up
to \,$q^\bullet t^\bullet$)\, for $\lla=(\la^j)$.

The following identity from (\ref{e-single-c-p})
is the key for this connection: 
\begin{align*}
&\hat{\h}^{\circ\rightthreearrow,min}_{\{1,-1\}}(\la^1,\la^2,\la^3)=
\hat{\h}^{min}_{\l,\,'\!\l}(\la^1,\la^2;\la^3)
\for \l=\circ\rightrightarrows,\,'\!\l=\circ\rightarrow,
\end{align*}
where the trees $\l,\,'\!\l$ are colored by
$\{\la^1,\la^2\}$ and $\{\la^3\}$ correspondingly.
 
Upon adding $\vee$ to $\,'\!\l$ the last formula (which is
the switch from $Y^{-1}$ 
to $Y$ in $\,'\!\l$), the super-polynomials
\begin{align}\label{vert-3-c}
&\hat{\h}^{min}_{\{1,-1\},\{1,0\}}(\la^1,\la^2;\la^3)^\vee\equal
\hat{\h}^{\rightrightarrows,\rightarrow, min}_{\l,\,'\!\l^\vee}
\end{align}
provide the $a$\~stabilizations
of $C_{\la^1,\la^2}^{\la^3}\,\lan P_{\la^3},P_{\la^3}\ran/
\langle \theta\mmu^\flat \rangle$, which are of obvious interest
since they contain no $\iota$ and therefore 
satisfy straight associativity (with certain multiplies)
from (\ref{ass-s-4}).
We will always take
$\la^3=\square$ in the examples below.

Since the recalculation of the $\hat{\h}$\~polynomials to
$\p_{\hbox{\tiny $\lla$}}$ and the 
$C$\~coefficients is straightforward, we will provide
only $\hat{\h}$. Recall that 
Proposition \ref{MACD-NORMS-A} and (\ref{stabevaly})
are used for this, where 
\begin{align}\label{stabevalyz}
&\Pi_\la^\dag=(\la)_{ev}^\dag=
\prod_{p=1}^n \prod_{v=0}^{\la_{p\,}-1} 
\Bigl( 1\,+\,q^{v}\, a\, t^{-p+1}\Bigr).
\end{align}
Also, $den^{j_o}$ will be provided, which are needed in
(\ref{conjhomflyz}):
\begin{align*}
&\hat{\h}^{min}_{\l}\,
(q,q,-a)/den^{j_o}=
\hat{\hbox{\small H\!O\!M}}^{j_o}_{\l}\,(q,a).
\end{align*}
\smallskip

\subsection{\bf Colored Hopf 3-links}
The first example will be uncolored; it is a direct continuation
of the previous section, namely the case of $T(3,-3)$ in the
notations there. It is self-dual with respect
to the super-duality.
\Yboxdim5pt 
 
\begin{align}\label{T1-1-1}
&\Yboxdim7pt \yng(1)\,,\yng(1)\,,\yng(1)\, :\ \,\Yboxdim5pt 
\l=\l_{(\{1,-1\},\{1,-1\},\{1,-1\})}^{\,\circ\rightthreearrow,\, 
(\yng(1)\,,\yng(1)\,,\yng(1)\,)},\ \ 
\hat{\h}{}^{min}_{\l}\,(q,t,a)=
\end{align}

\renewcommand{\baselinestretch}{0.5} 
\noindent
{\small
\(
1+a \bigl(q+q^2+\frac{1}{t^2}-\frac{2 q}{t^2}
+\frac{q^2}{t^2}+\frac{1}{t}+\frac{q}{t}-\frac{2 q^2}{t}\bigr)
+a^2 \bigl(q^3+\frac{1}{t^3}-\frac{2 q}{t^3}+\frac{q^2}{t^3}+
\frac{q}{t^2}-\frac{2 q^2}{t^2}+\frac{q^3}{t^2}+
\frac{q}{t}+\frac{q^2}{t}-\frac{2 q^3}{t}\bigr).
\)
}
\renewcommand{\baselinestretch}{1.0} 

Let us provide its $\vee$\~variant:
\begin{align}\label{T1-1-1v}
&\Yboxdim7pt \yng(1)\,,\yng(1)\,;\yng(1)\,{}^\vee\, :\ \Yboxdim5pt 
\l=\l_{(\{1,-1\},\{1,-1\})}^{\,\circ\rightrightarrows,
(\yng(1)\,,\yng(1)\,)}, 
\,'\!\l=\l_{\{1,0\}}^{\,\circ\rightarrow,\,\yng(1)\,},\ \ 
\hat{\h}{}^{min}_{\l,\,'\!\l^\vee}=
\end{align}

\renewcommand{\baselinestretch}{0.5} 
\noindent
\centerline{\small
\(
2-q+a^2 \bigl(q^2+\frac{q}{t}-\frac{q^2}{t}\bigr)-2 t
+2 q t+a \bigl(3 q-q^2+\frac{1}{t}-\frac{q}{t}-t+q^2 t\bigr).
\)
}
\renewcommand{\baselinestretch}{1.0} 
\smallskip

We will provide the formulas with $\vee$
for all further examples:
\smallskip

\begin{align}\label{T11-1-1}
&\Yboxdim7pt \yng(2)\,,\yng(1)\,,\yng(1)\, :\ \,\ \Yboxdim5pt
\l=\l_{(\{1,-1\},\{1,-1\},\{1,-1\})}^{\,\circ\rightthreearrow,\, 
(\yng(2)\,,\yng(1)\,,\yng(1)\,)},\ \ 
\hat{\h}{}^{min}_{\l}\,(q,t,a)=
\end{align}

\renewcommand{\baselinestretch}{0.5} 
\noindent
{\small
\(
1+a \bigl(q^2+q^3+\frac{1}{t^2}-\frac{q}{t^2}-\frac{q^2}{t^2}
+\frac{q^3}{t^2}+\frac{1}{t}+\frac{q}{t}-\frac{2 q^3}{t}\bigr)
+a^2 \bigl(q^5+\frac{1}{t^3}-\frac{q}{t^3}-\frac{q^2}{t^3}
+\frac{q^3}{t^3}+\frac{q}{t^2}-\frac{2 q^3}{t^2}
+\frac{q^5}{t^2}+\frac{q^2}{t}+\frac{q^3}{t}-\frac{2 q^5}{t}\bigr).
\)
}
\renewcommand{\baselinestretch}{1.0} 

\begin{align}\label{T11-1-1v}
&\Yboxdim7pt \yng(2)\,,\yng(1)\,;\yng(1)\,{}^\vee\, :\ \Yboxdim5pt 
\l=\l_{(\{1,-1\},\{1,-1\})}^{\,\circ\rightrightarrows,
(\yng(2)\,,\yng(1)\,)}, 
\,'\!\l=\l_{\{1,0\}}^{\,\circ\rightarrow,\,\yng(1)\,},\ \ 
\hat{\h}{}^{min}_{\l,\,'\!\l^\vee}=
\end{align}

\renewcommand{\baselinestretch}{0.5} 
\noindent
{\small
\(
1+q-q^2+a^2 \bigl(q^4+\frac{q^2}{t}-\frac{q^4}{t}\bigr)
-t-q t+2 q^2 t
+a \bigl(-1+q+3 q^2-q^4+\frac{1}{t}-\frac{q^2}{t}-q t+q^4 t\bigr).
\)
}
\renewcommand{\baselinestretch}{1.0} 
\smallskip

\begin{align}\label{T20-1-1}
&\Yboxdim7pt \yng(1,1)\,,\yng(1)\,,\yng(1)\, :\ \,\ \Yboxdim5pt
\l=\l_{(\{1,-1\},\{1,-1\},\{1,-1\})}^{\,\circ\rightthreearrow,\, 
(\yng(1,1)\,,\yng(1)\,,\yng(1)\,)},\ \ 
\hat{\h}{}^{min}_{\l}\,(q,t,a)=
\end{align}

\renewcommand{\baselinestretch}{0.5} 
\noindent
{\small
\(
1+a \bigl(q+q^2+\frac{1}{t^3}-\frac{2 q}{t^3}+\frac{q^2}{t^3}
+\frac{1}{t^2}-\frac{q^2}{t^2}+\frac{q}{t}-\frac{q^2}{t}\bigr)
+a^2 \bigl(q^3+\frac{1}{t^5}-\frac{2 q}{t^5}+\frac{q^2}{t^5}
+\frac{q}{t^3}-\frac{2 q^2}{t^3}+\frac{q^3}{t^3}
+\frac{q}{t^2}-\frac{q^3}{t^2}+\frac{q^2}{t}-\frac{q^3}{t}\bigr).
\)
}
\renewcommand{\baselinestretch}{1.0} 
\smallskip

The last example is dual to (\ref{T11-1-1}) under
$t\leftrightarrow q^{-1}$. Usually the super-duality is up to 
certain $q^\bullet t^\bullet$, but for the Hopf link
{\em without $\vee$\,}, such multipliers are 
not needed since the $a$\~constant term is always $1$.
The constant term is not $1$ with $\vee$.
We omit the corresponding $\vee$\~variant for
(\ref{T20-1-1}), since it is
super-dual to (\ref{T11-1-1v}). The next case will
be self-dual, since the diagram 
$\Yboxdim5pt \yng(2,1)$ is transposition-symmetric.

\begin{align}\label{T21-1-1}
&\Yboxdim7pt \yng(2,1)\,,\yng(1)\,,\yng(1)\, :\ \,\ \Yboxdim5pt
\l=\l_{(\{1,-1\},\{1,-1\},\{1,-1\})}^{\,\circ\rightthreearrow,\, 
(\yng(2,1)\,,\yng(1)\,,\yng(1)\,)},\ \ 
\hat{\h}{}^{min}_{\l}\,(q,t,a)=
\end{align}

\renewcommand{\baselinestretch}{0.5} 
\noindent
{\small
\(
1+a \bigl(q^2+q^3+\frac{1}{t^3}-\frac{2 q}{t^3}+\frac{q^2}{t^3}
+\frac{1}{t^2}+\frac{q}{t^2}-\frac{3 q^2}{t^2}+\frac{q^3}{t^2}
+\frac{q}{t}+\frac{q^2}{t}-\frac{2 q^3}{t}\bigr)
+a^2 \bigl(q^5+\frac{1}{t^5}-\frac{2 q}{t^5}+\frac{q^2}{t^5}
+\frac{q}{t^4}-\frac{3 q^2}{t^4}+\frac{3 q^3}{t^4}-\frac{q^4}{t^4}
+\frac{q}{t^3}+\frac{q^2}{t^3}-\frac{5 q^3}{t^3}+\frac{3 q^4}{t^3}
+\frac{q^2}{t^2}+\frac{q^3}{t^2}-\frac{3 q^4}{t^2}
+\frac{q^5}{t^2}+\frac{q^3}{t}+\frac{q^4}{t}-\frac{2 q^5}{t}\bigr).
\)
}
\renewcommand{\baselinestretch}{1.0} 

\begin{align}\label{T21-1-1v}
&\Yboxdim7pt \yng(2,1)\,,\yng(1)\,;\yng(1)\,{}^\vee\, :\ \Yboxdim5pt 
\l=\l_{(\{1,-1\},\{1,-1\})}^{\,\circ\rightrightarrows,
(\yng(2,1)\,,\yng(1)\,)}, 
\,'\!\l=\l_{\{1,0\}}^{\,\circ\rightarrow,\,\yng(1)\,},\ \ 
\hat{\h}{}^{min}_{\l,\,'\!\l^\vee}=
\end{align}

\renewcommand{\baselinestretch}{0.5} 
\noindent
{\small
\(
1+q-q^2+a^2 \bigl(q^4+\frac{q^2}{t^2}-\frac{q^3}{t^2}+\frac{q^3}{t}
-\frac{q^4}{t}\bigr)-t+q^2 t-q t^2+q^2 t^2+a \bigl(-q+5 q^2-2 q^3
+\frac{1}{t^2}-\frac{q}{t^2}-\frac{1}{t}+\frac{3 q}{t}
-\frac{2 q^2}{t}-q t-q^2 t+3 q^3 t-q^4 t-q^3 t^2+q^4 t^2\bigr).
\)
}
\renewcommand{\baselinestretch}{1.0} 
\smallskip

Note that 
$\hat{\h}^{\emptyset,\dag}=
\hat{\h}^{\min}/(1+a)^2$ and  $den^1=(1-q)^2$
in all these examples, since 
$\la^1\!\vee\la^2\!\vee\la^{3}$ is $\la^1$
(it contains the remaining two). In the next two examples,
$den^1=(1-q)^2(1-q^2)$. 
\smallskip

\begin{align}\label{T11-11-1}
&\Yboxdim7pt \yng(2)\,,\yng(2)\,,\yng(1)\, :\ \,\ \Yboxdim5pt
\l=\l_{(\{1,-1\},\{1,-1\},\{1,-1\})}^{\,\circ\rightthreearrow,\, 
(\yng(2)\,,\yng(2)\,,\yng(1)\,)},\ \ 
\hat{\h}{}^{min}_{\l}\,(q,t,a)=
\end{align}

\renewcommand{\baselinestretch}{0.5} 
\noindent
{\small
\(
1+a \bigl(q^2+q^3+q^4+\frac{1}{t^2}-\frac{2 q^2}{t^2}
+\frac{q^4}{t^2}+\frac{1}{t}+\frac{q}{t}+\frac{q^2}{t}-\frac{q^3}{t}
-\frac{2 q^4}{t}\bigr)+a^2 \bigl(q^5+q^6+q^7+\frac{1}{t^3}
+\frac{q}{t^3}-\frac{3 q^2}{t^3}-\frac{2 q^3}{t^3}+\frac{3 q^4}{t^3}
+\frac{q^5}{t^3}-\frac{q^6}{t^3}+\frac{q}{t^2}+\frac{2 q^2}{t^2}
-\frac{5 q^4}{t^2}-\frac{2 q^5}{t^2}+\frac{3 q^6}{t^2}
+\frac{q^7}{t^2}+\frac{q^2}{t}+\frac{2 q^3}{t}+\frac{2 q^4}{t}
-\frac{3 q^6}{t}-\frac{2 q^7}{t}\bigr)+a^3 \bigl(q^9+\frac{q}{t^4}
-\frac{q^2}{t^4}-\frac{2 q^3}{t^4}+\frac{2 q^4}{t^4}+\frac{q^5}{t^4}
-\frac{q^6}{t^4}+\frac{q^2}{t^3}+\frac{q^3}{t^3}-\frac{3 q^4}{t^3}
-\frac{2 q^5}{t^3}+\frac{3 q^6}{t^3}+\frac{q^7}{t^3}-\frac{q^8}{t^3}
+\frac{q^3}{t^2}+\frac{q^4}{t^2}-\frac{3 q^6}{t^2}-\frac{2 q^7}{t^2}
+\frac{2 q^8}{t^2}+\frac{q^9}{t^2}+\frac{q^5}{t}
+\frac{q^6}{t}+\frac{q^7}{t}-\frac{q^8}{t}-\frac{2 q^9}{t}\bigr).
\)
}
\renewcommand{\baselinestretch}{1.0} 

\begin{align}\label{T11-11-1v}
&\Yboxdim7pt \yng(2)\,,\yng(2)\,;\yng(1)\,{}^\vee\, :\ \Yboxdim5pt 
\l=\l_{(\{1,-1\},\{1,-1\})}^{\,\circ\rightrightarrows,
(\yng(2)\,,\yng(2)\,)}, 
\,'\!\l=\l_{\{1,0\}}^{\,\circ\rightarrow,\,\yng(1)\,},\ \ 
\hat{\h}{}^{min}_{\l,\,'\!\l^\vee}=
\end{align}

\renewcommand{\baselinestretch}{0.5} 
\noindent
{\small
\(
2-q^2+a^3 \bigl(q^7+\frac{q^3}{t^2}-\frac{q^4}{t^2}-\frac{q^5}{t^2}
+\frac{q^6}{t^2}+\frac{q^4}{t}+\frac{q^5}{t}-\frac{q^6}{t}
-\frac{q^7}{t}\bigr)-2 t+2 q^2 t+a \bigl(-2+7 q^2+2 q^3-3 q^4-q^5
+\frac{2}{t}+\frac{q}{t}-\frac{3 q^2}{t}-\frac{q^3}{t}+\frac{q^4}{t}
-q t-2 q^2 t+2 q^4 t+q^5 t\bigr)+a^2 \bigl(-q-q^2+2 q^3+3 q^4+2 q^5
-q^6-q^7+\frac{q}{t^2}-\frac{q^2}{t^2}-\frac{q^3}{t^2}
+\frac{q^4}{t^2}+\frac{3 q^2}{t}+\frac{2 q^3}{t}
-\frac{4 q^4}{t}-\frac{2 q^5}{t}+\frac{q^6}{t}-q^3 t+q^7 t\bigr).
\)
}
\renewcommand{\baselinestretch}{1.0} 
\smallskip

The next 2 cases are dual to the previous ones. We provide them, 
but then will omit the examples dual to those already given.
Note that 
\begin{align*}
\hat{\h}_{\l}^{\emptyset,\dag}
(\Yboxdim5pt \yng(2)\,,\yng(2)\,,\yng(1))\!=\!
\frac{\hat{\h}_{\l}^{\min}(\Yboxdim5pt \yng(2)\,,\yng(2)\,,\yng(1))}
{(1+a)^2(1+qa)},\  
\hat{\h}_{\l}^{\emptyset,\dag}
(\Yboxdim5pt \yng(1,1)\,,\yng(1,1)\,,\yng(1))\!=\!
\frac{\hat{\h}_{\l}^{\min}
(\Yboxdim5pt \yng(1,1)\,,\yng(1,1)\,,\yng(1))}
{(1+a)^2(1+a/t)}.
\end{align*}
\smallskip

\begin{align}\label{T20-20-1}
&\Yboxdim7pt \yng(1,1)\,,\yng(1,1)\,,\yng(1)\, :\ \,\ \Yboxdim5pt
\l=\l_{(\{1,-1\},\{1,-1\},\{1,-1\})}^{\,\circ\rightthreearrow,\, 
(\yng(1,1)\,,\yng(1,1)\,,\yng(1)\,)},\ \ 
\hat{\h}{}^{min}_{\l}\,(q,t,a)=
\end{align}

\renewcommand{\baselinestretch}{0.5} 
\noindent
{\small
\(
1+a \bigl(q+q^2+\frac{1}{t^4}-\frac{2 q}{t^4}
+\frac{q^2}{t^4}+\frac{1}{t^3}-\frac{q}{t^3}+\frac{1}{t^2}
+\frac{q}{t^2}-\frac{2 q^2}{t^2}+\frac{q}{t}\bigr)
+a^2 \bigl(q^3+\frac{1}{t^7}-\frac{2 q}{t^7}+\frac{q^2}{t^7}
+\frac{1}{t^6}-\frac{3 q}{t^6}+\frac{3 q^2}{t^6}-\frac{q^3}{t^6}
+\frac{1}{t^5}-\frac{2 q^2}{t^5}+\frac{q^3}{t^5}+\frac{2 q}{t^4}
-\frac{5 q^2}{t^4}+\frac{3 q^3}{t^4}+\frac{2 q}{t^3}
-\frac{2 q^3}{t^3}+\frac{q}{t^2}+\frac{2 q^2}{t^2}-\frac{3 q^3}{t^2}
+\frac{q^2}{t}+\frac{q^3}{t}\bigr)+a^3 \bigl(\frac{1}{t^9}
-\frac{2 q}{t^9}+\frac{q^2}{t^9}-\frac{q}{t^8}+\frac{2 q^2}{t^8}
-\frac{q^3}{t^8}+\frac{q}{t^7}-\frac{2 q^2}{t^7}+\frac{q^3}{t^7}
+\frac{q}{t^6}-\frac{3 q^2}{t^6}+\frac{3 q^3}{t^6}-\frac{q^4}{t^6}
+\frac{q}{t^5}-\frac{2 q^3}{t^5}+\frac{q^4}{t^5}+\frac{q^2}{t^4}
-\frac{3 q^3}{t^4}+\frac{2 q^4}{t^4}+\frac{q^2}{t^3}+\frac{q^3}{t^3}
-\frac{2 q^4}{t^3}
+\frac{q^3}{t^2}-\frac{q^4}{t^2}+\frac{q^4}{t}\bigr).
\)
}
\renewcommand{\baselinestretch}{1.0} 

\begin{align}\label{T20-20-1v}
&\Yboxdim7pt \yng(1,1)\,,\yng(1,1)\,;\yng(1)\,{}^\vee\, :\ 
\Yboxdim5pt 
\l=\l_{(\{1,-1\},\{1,-1\})}^{\,\circ\rightrightarrows,
(\yng(1,1)\,,\yng(1,1)\,)}, 
\,'\!\l=\l_{\{1,0\}}^{\,\circ\rightarrow,\,\yng(1)\,},\ \ 
\hat{\h}{}^{min}_{\l,\,'\!\l^\vee}=
\end{align}

\renewcommand{\baselinestretch}{0.5} 
\noindent
{\small
\(
2-q+a^3 \bigl(\frac{q}{t^5}-\frac{q^2}{t^5}-\frac{q^2}{t^4}
+\frac{q^3}{t^4}+\frac{q^2}{t^3}-\frac{q^3}{t^3}+\frac{q^2}{t^2}
-\frac{q^3}{t^2}+\frac{q^3}{t}\bigr)-2 t^2+2 q t^2
+a^2 \bigl(-q+3 q^2-q^3+\frac{1}{t^5}-\frac{q}{t^5}-\frac{q}{t^4}
+\frac{q^2}{t^4}+\frac{2 q}{t^3}-\frac{2 q^2}{t^3}+\frac{3 q}{t^2}
-\frac{4 q^2}{t^2}+\frac{q^3}{t^2}-\frac{1}{t}+\frac{2 q}{t}
+\frac{2 q^2}{t}-\frac{q^3}{t}-q t+q^3 t\bigr)+a \bigl(-2+7 q-3 q^2
+\frac{1}{t^3}-\frac{q}{t^3}+\frac{2}{t^2}-\frac{3 q}{t^2}
+\frac{q^2}{t^2}
+\frac{2 q}{t}-\frac{q^2}{t}-t+q^2 t-2 q t^2+2 q^2 t^2\bigr).
\)
}
\renewcommand{\baselinestretch}{1.0} 
\smallskip

\begin{align}\label{T20-11-1}
&\Yboxdim7pt \yng(2)\,,\yng(1,1)\,,\yng(1)\, :\ \,\ \Yboxdim5pt
\l=\l_{(\{1,-1\},\{1,-1\},\{1,-1\})}^{\,\circ\rightthreearrow,\, 
(\yng(2)\,,\yng(1,1)\,,\yng(1)\,)},\ \ 
\hat{\h}{}^{min}_{\l}\,(q,t,a)=
\end{align}

\renewcommand{\baselinestretch}{0.5} 
\noindent
{\small
\(
1+a \bigl(q^2+q^3+\frac{1}{t^3}-\frac{q}{t^3}
-\frac{q^2}{t^3}+\frac{q^3}{t^3}+\frac{1}{t^2}-\frac{q^3}{t^2}
+\frac{q}{t}-\frac{q^3}{t}\bigr)+a^2 \bigl(q^5+\frac{1}{t^5}
-\frac{q}{t^5}-\frac{q^2}{t^5}+\frac{q^3}{t^5}+\frac{q}{t^3}
-\frac{2 q^3}{t^3}+\frac{q^5}{t^3}
+\frac{q^2}{t^2}-\frac{q^5}{t^2}+\frac{q^3}{t}-\frac{q^5}{t}\bigr).
\)
}
\renewcommand{\baselinestretch}{1.0} 

\begin{align}\label{T20-11-1v}
&\Yboxdim7pt \yng(2)\,,\yng(1,1)\,;\yng(1)\,{}^\vee\, :\ 
\Yboxdim5pt 
\l=\l_{(\{1,-1\},\{1,-1\})}^{\,\circ\rightrightarrows,
(\yng(2)\,,\yng(1,1)\,)}, 
\,'\!\l=\l_{\{1,0\}}^{\,\circ\rightarrow,\,\yng(1)\,},\ \ 
\hat{\h}{}^{min}_{\l,\,'\!\l^\vee}=
\end{align}

\renewcommand{\baselinestretch}{0.5} 
\noindent
{\small
\(
1+q-q^2+a^2 \bigl(q^4+\frac{q^2}{t^2}-\frac{q^4}{t^2}\bigr)
-t+q^2 t-q t^2+q^2 t^2+a \bigl(2 q^2+q^3-q^4+\frac{1}{t^2}
-\frac{q^2}{t^2}-\frac{1}{t}+\frac{q}{t}
+\frac{q^2}{t}-\frac{q^3}{t}-q t+q^3 t-q^3 t^2+q^4 t^2\bigr).
\)
}
\renewcommand{\baselinestretch}{1.0} 
\medskip

These 2 examples are self-dual. The degree here is
$|\ss|(|\la^1|+|\la^2|+|\la^3|)-|\la^1\!\vee\la^2\!\vee\la^{3}|=
2+2+1-3=2,$
since $\la^1\!\vee\la^2\!\vee\la^{3}\,$ is the $3$-hook. Also:
$
\,\hat{\h}_{\l}^{\emptyset,\dag}
(\Yboxdim5pt \yng(2)\,,\yng(1,1)\,,\yng(1))=
\hat{\h}_{\l}^{\min}(\Yboxdim5pt \yng(2)\,,\yng(1,1)\,,\yng(1))/
(1+a)^2,\ 
$
$den^1=\frac{(1-q)^2(1-q^2)}{(1-a/q)}.$

\smallskip
\begin{align}\label{T21-20-1}
&\Yboxdim7pt \yng(2,1),\,\yng(1,1)\,,\,\yng(1)\, :\ \,\ \Yboxdim5pt
\l=\l_{(\{1,-1\},\{1,-1\},\{1,-1\})}^{\,\circ\rightthreearrow,\, 
(\yng(2,1),\,\yng(1,1)\,,\,\yng(1)\,)},\ \ 
\hat{\h}{}^{min}_{\l}\,(q,t,a)=
\end{align}

\renewcommand{\baselinestretch}{0.5} 
\noindent
{\small
\(
1+a \bigl(q^2+q^3+\frac{1}{t^4}-\frac{2 q}{t^4}
+\frac{q^2}{t^4}+\frac{1}{t^3}-\frac{2 q^2}{t^3}
+\frac{q^3}{t^3}+\frac{1}{t^2}+\frac{q}{t^2}-\frac{q^2}{t^2}
-\frac{q^3}{t^2}+\frac{q}{t}+\frac{q^2}{t}-\frac{q^3}{t}\bigr)
+a^2 \bigl(q^5+\frac{1}{t^7}-\frac{2 q}{t^7}+\frac{q^2}{t^7}
+\frac{1}{t^6}-\frac{2 q}{t^6}+\frac{2 q^3}{t^6}-\frac{q^4}{t^6}
+\frac{1}{t^5}+\frac{q}{t^5}-\frac{4 q^2}{t^5}+\frac{2 q^3}{t^5}
+\frac{2 q}{t^4}-\frac{q^2}{t^4}-\frac{5 q^3}{t^4}
+\frac{4 q^4}{t^4}
+\frac{q}{t^3}+\frac{3 q^2}{t^3}-\frac{3 q^3}{t^3}
-\frac{2 q^4}{t^3}
+\frac{q^5}{t^3}+\frac{q^2}{t^2}+\frac{3 q^3}{t^2}
-\frac{3 q^4}{t^2}
-\frac{q^5}{t^2}+\frac{q^3}{t}+\frac{2 q^4}{t}-\frac{q^5}{t}\bigr)
+a^3 \bigl(\frac{1}{t^9}-\frac{2 q}{t^9}+\frac{q^2}{t^9}
-\frac{q^2}{t^8}+\frac{2 q^3}{t^8}-\frac{q^4}{t^8}+\frac{q}{t^7}
-\frac{2 q^2}{t^7}+\frac{q^3}{t^7}+\frac{q}{t^6}-\frac{4 q^3}{t^6}
+\frac{4 q^4}{t^6}-\frac{q^5}{t^6}+\frac{2 q^2}{t^5}
-\frac{2 q^3}{t^5}-\frac{2 q^4}{t^5}+\frac{3 q^5}{t^5}
-\frac{q^6}{t^5}+\frac{2 q^3}{t^4}-\frac{3 q^4}{t^4}
-\frac{q^5}{t^4}
+\frac{2 q^6}{t^4}+\frac{q^3}{t^3}+\frac{2 q^4}{t^3}
-\frac{3 q^5}{t^3}
+\frac{2 q^5}{t^2}-\frac{2 q^6}{t^2}+\frac{q^6}{t}\bigr).
\)
}
\renewcommand{\baselinestretch}{1.0} 

\begin{align}\label{T21-20-1v}
&\Yboxdim7pt \yng(2,1)\,,\yng(1,1)\,;\yng(1)\,{}^\vee\, :\ 
\Yboxdim5pt 
\l=\l_{(\{1,-1\},\{1,-1\})}^{\,\circ\rightrightarrows,
(\yng(2,1)\,,\yng(1,1)\,)}, 
\,'\!\l=\l_{\{1,0\}}^{\,\circ\rightarrow,\,\yng(1)\,},\ \ 
\hat{\h}{}^{min}_{\l,\,'\!\l^\vee}=
\end{align}

\renewcommand{\baselinestretch}{0.5} 
\noindent
{\small
\(
1+q-q^2+a^3 \bigl(\frac{q^2}{t^5}-\frac{q^3}{t^5}-\frac{q^4}{t^4}
+\frac{q^5}{t^4}+\frac{q^3}{t^3}-\frac{q^5}{t^3}+\frac{q^4}{t^2}
-\frac{q^5}{t^2}+\frac{q^5}{t}\bigr)-t+q t-2 q t^2+2 q^2 t^2
+a^2 \bigl(-q^2+3 q^4-q^5+\frac{1}{t^5}-\frac{q}{t^5}-\frac{1}{t^4}
+\frac{2 q}{t^4}-\frac{2 q^2}{t^4}+\frac{q^3}{t^4}+\frac{q}{t^3}
+\frac{2 q^2}{t^3}-\frac{4 q^3}{t^3}+\frac{q^4}{t^3}-\frac{q}{t^2}
+\frac{4 q^2}{t^2}-\frac{q^3}{t^2}-\frac{3 q^4}{t^2}
+\frac{q^5}{t^2}
-\frac{q}{t}+\frac{5 q^3}{t}-\frac{q^4}{t}-\frac{q^5}{t}-q^3 t
+q^5 t\bigr)+a \bigl(-3 q+6 q^2+q^3-2 q^4+\frac{1}{t^3}
-\frac{q}{t^3}+\frac{3 q}{t^2}-\frac{4 q^2}{t^2}+\frac{q^3}{t^2}
-\frac{1}{t}+\frac{2 q}{t}+\frac{3 q^2}{t}-\frac{4 q^3}{t}
+\frac{q^4}{t}-q t-2 q^2 t+4 q^3 t-q^4 t-2 q^3 t^2+2 q^4 t^2\bigr).
\)
}
\renewcommand{\baselinestretch}{1.0} 
\smallskip

The $a$\~degree is
$|\ss|(|\la^1|+|\la^2|+|\la^3|)-|\la^1\!\vee\la^2\!\vee\la^{3}|=
3+2+1-3=3.$ Correspondingly,
$
\hat{\h}_{\l}^{\emptyset,\dag}
(\Yboxdim5pt \yng(2,1),\,\yng(1,1)\,,\,\yng(1))=
\hat{\h}_{\l}^{\min}(\Yboxdim5pt \yng(2,1),\,
\yng(1,1)\,,\,\yng(1))/
\bigl((1+a)^2(1+aq)\bigr),\,
$
$den^1=(1-q)^2(1-q^2).$ 
\smallskip

\subsection{\bf Hopf 2-links}
The previous examples were for the Hopf $3$\~links. Let us
discuss the Hopf $2$\~links.

\begin{align}\label{T11-1}
&\Yboxdim7pt \yng(2)\,,\yng(1)\,:\ \Yboxdim5pt
\l\!=\!\l_{(\{1,-1\},\{1,-1\})}^{\,\circ\rightrightarrows,\, 
(\yng(2)\,,\yng(1)\,)},\,
\hat{\h}{}^{min}_{\l}\,(q,t,a)\!=\!\hbox{\small
\(1+a \bigl(q^2+\frac{1}{t}-\frac{q^2}{t}\bigr).
\)}
\end{align}

\smallskip

\begin{align}\label{T11-20}
&\Yboxdim7pt \yng(2)\,,\yng(1,1)\, :\ 
\hat{\h}{}^{min}_{\,\circ\rightrightarrows,\,}
(\Yboxdim5pt \yng(2)\,,\yng(1,1)\,;\,q,t,a)=\hbox{\small
\(1+a \bigl(q^2+\frac{1}{t^2}-\frac{q^2}{t^2}\bigr)\)},\\
\label{T21-1}
&\Yboxdim7pt \yng(2,1)\,,\yng(1)\, :\ 
\hat{\h}_{\,\circ\rightrightarrows\,}^{min}\,
(\Yboxdim5pt \yng(2,1)\,,\yng(1)\,;\,q,t,a)\!=\!\hbox{\small
\(
1\!+\!a \bigl(q^2\!+\!\frac{1}{t^2}\!-\!\frac{q}{t^2}\!+\!
\frac{q}{t}\!-\!\frac{q^2}{t}\bigr)
\)}.
\end{align}

The last two examples are self-dual, as well as
the simplest superpolynomial
$\hat{\h}_{\,\circ\rightrightarrows,\,}^{\min}
(\Yboxdim5pt \yng(1)\,,\yng(1)\,)=
1+a(q+1/t-q/t)$, which is actually uncolored $T(2,-2)$, already 
considered in the previous section. 
The next $2$ ones are not self-dual:

\begin{align}\label{T11-11}
&\Yboxdim7pt \yng(2)\,,\yng(2)\, :\ \,\ \Yboxdim5pt
\l=\l_{(\{1,-1\},\{1,-1\})}^{\,\circ\rightrightarrows,\, 
(\yng(2)\,,\yng(2)\,)},\ \ 
\hat{\h}{}^{min}_{\l}\,(q,t,a)=
\end{align}

\renewcommand{\baselinestretch}{0.5} 
\noindent
{\small
\(
1+a \bigl(q^2+q^3+\frac{1}{t}+\frac{q}{t}-\frac{q^2}{t}
-\frac{q^3}{t}\bigr)+a^2 \bigl(q^5+\frac{q}{t^2}-\frac{q^2}{t^2}
-\frac{q^3}{t^2}+\frac{q^4}{t^2}
+\frac{q^2}{t}+\frac{q^3}{t}-\frac{q^4}{t}-\frac{q^5}{t}\bigr).
\)
}
\renewcommand{\baselinestretch}{1.0} 
\smallskip

\begin{align}\label{T21-20}
&\Yboxdim7pt \yng(2,1)\,,\yng(2)\, :\ \,\ \Yboxdim5pt
\l=\l_{(\{1,-1\},\{1,-1\})}^{\,\circ\rightrightarrows,\, 
(\yng(2,1)\,,\yng(2)\,)},\ \ 
\hat{\h}{}^{min}_{\l}\,(q,t,a)=
\end{align}

\renewcommand{\baselinestretch}{0.5} 
\noindent
{\small
\(
1+a \bigl(q^2+q^3+\frac{1}{t^2}-\frac{q^2}{t^2}
+\frac{q}{t}-\frac{q^3}{t}\bigr)+a^2 \bigl(q^5+\frac{q}{t^3}
-\frac{q^2}{t^3}-\frac{q^3}{t^3}+\frac{q^4}{t^3}+\frac{q^2}{t^2}
-\frac{q^4}{t^2}
+\frac{q^3}{t}-\frac{q^5}{t}\bigr).
\)
}
\renewcommand{\baselinestretch}{1.0} 
\smallskip

Here deg$_a=2$ and 
$
\hat{\h}_{\l}^{\emptyset,\dag}
(\Yboxdim5pt \yng(2,1)\,,\yng(2)\,)=
\hat{\h}_{\l}^{\min}(\Yboxdim5pt \yng(2,1)\,,\yng(2)\,)/
\bigl((1+a)(1+aq)\bigr).
$
The same denominator serves 
$\hat{\h}_{\l}^{\emptyset,\dag}
(\Yboxdim5pt \yng(2)\,,\yng(2)\,)$. Otherwise, deg$_a=1$  
in the examples above 
and this denominator is $(1+a)$. Let us provide $den^1$
from (\ref{den-hom}), needed for obtaining 
the reduced hat-normalized HOMFLY-PT polynomials. One has
$den^1(\yng(2)\,,\yng(1,1))=
\frac{(1-q)(1-q^2)}{(1-a/q)}$; otherwise:
\begin{align}\label{den1-two}
den^1(\la^1,\raisebox{1pt}{\yng(1)})\!=\!(1\!-\!q),\,
den^1(\la^1,\raisebox{1pt}{\yng(2)})\!=\!(1\!-\!q)(1\!-\!q^2)\!=\!
den^1(\la^1,\raisebox{1pt}{\yng(1,1)}),
\end{align}
since $\la^1$ contains $\la^2$ in (\ref{den1-two}).

\medskip
{\sf Connection to the topological vertex.}
Let $(\hat{\h}_{\,\circ\rightrightarrows\,}^{min})^{\infty}$
be the leading $a$\~coefficient of  
$\hat{\h}_{\,\circ\rightrightarrows\,}^{min}$ and 
$
C_{D\!A\!H\!\!A}^{num}(\la^1,\la^2)=
\bigl(\hat{\h}_{\,\circ\rightrightarrows\,}^{min}(\la^1,\la^2) 
\bigr)^{\infty}(t\mapsto q),$ 
which is proportional to $C_{\la^1,\la^2}^0$ and
$C_{\,\la^1}^{\,\iota(\la^2)} 
\lan P_{\la^2}^\circ,P_{\la^2}^\circ\ran\,$
in this limit. Let us provide $C_{D\!A\!H\!\!A}^{num}$
in the following 2 examples.

\begin{align}\label{T30-11}
&\Yboxdim7pt \yng(1,1,1)\,,\yng(2)\, :\ 
\hat{\h}{}^{min}=
\hbox{\small
\(
1+a \bigl(q^2+\frac{1}{t^3}-\frac{q^2}{t^3}\bigr);\,
den^1=\frac{(1-q)(1-q^2)}{1-aq}\).}
\end{align}
{\small
\begin{align*}
\ \ \ \ \ \ \ \hat{\hbox{\small H\!O\!M}}\!=\!
\frac{
\bigl(1-a (\frac{1}{q^3}-\frac{1}{q}+q^2)\bigr)(1-aq)}
{(1-q)(1-q^2)}, \ C_{D\!A\!H\!\!A}^{num}=1-q^2+q^5,
\end{align*}
}
\vskip -0.1cm

\begin{align}\label{T30-20}
&\Yboxdim7pt \yng(1,1,1)\,,\yng(1,1)\, :\ 
\hat{\h}_{\,\circ\rightrightarrows\,}^{min}\,
(\Yboxdim5pt \yng(1,1,1)\,,\yng(1,1)\,;\,q,t,a)\!=\!
\end{align}
{\small
\(
1+a \bigl(q\!+\frac{1}{t^4}\!-\frac{q}{t^4}\!+\frac{1}{t^3}
\!-\frac{q}{t^3}\!+\frac{q}{t}\bigr)
+a^2 \bigl(\frac{1}{t^7}\!-\frac{q}{t^7}\!-\frac{q}{t^6}
\!+\frac{q^2}{t^6}\!+\frac{q}{t^4}\!-\frac{q^2}{t^4}
\!+\frac{q}{t^3}\!-\frac{q^2}{t^3}\!+\frac{q^2}{t}\bigr),\,
\)
\(\hat{\h}
(\Yboxdim5pt \yng(1,1,1)\,,\yng(1,1)\,;\,q,q,-a)\!=\!
1\!-a \bigl(1\!+\frac{1}{q^4}\!-\frac{1}{q^2}\!+q\bigr)
\!+a^2 \bigl(\frac{1}{q^7}\!-\frac{1}{q^6}
\!-\frac{1}{q^5}\!+\frac{1}{q^4}\!+\frac{1}{q^3}
\!-\frac{1}{q}\!+q\bigr),
\)

\(
den^1=(1-q)(1-q^2),\ \ C_{D\!A\!H\!\!A}^{num}=
(\Yboxdim5pt \yng(1,1,1)\,,\yng(1,1)\,;\,q,a)=
1-q-q^2+q^3+q^4-q^6+q^8.
\)
}

\medskip

{\em Vertex amplitudes.}
We claim that  
$C_{D\!A\!H\!\!A}^{num}(\la^1,\la^2)$ 
coincide up to $q^\bullet$ with the 
{\em numerators\,} 
of $C_{V\!\!A}(\la^1,(\la^2)^{tr})$, from Section ``The 
vertex amplitudes" of \cite{AKMV}. We denote them by $C_{V\!\!A}$ 
to 
avoid confusion with our own $C$. Note the transposition here; by
the way, if $tr$ is applied to $\la^1$ instead of $\la^2$, 
then $q$ must be 
changed to $q^{-1}$. This coincidence includes  $\la^2=\emptyset$, 
when the corresponding weight is $b=0$ (denoted by dot in 
\cite{AKMV}); they numerators are pure powers
$q^\bullet$ in this case. 
\smallskip

This claim can be deduced from 
(\ref{lev-1-two}) and we checked it 
numerically in all examples we considered.
The denominators of $C_{V\!\!A}(\la^1,(\la^2)^{\tr})$ are the 
products of binomials in the form $(1-q^\bullet)$ and are actually
the matter of normalization of the $2$\~vertex. 
Such vertex amplitudes are connected with the stable Macdonald
theory, which generally corresponds to taking $a=0$ in our 
approach. Here we take the top power of $a$ instead of $a=0$, 
since the conjugation is actually
needed to connect our $C$ with $C_{V\!\!A}.$
See e.g., \cite{AFS}.

However, {\em the connection of our $C$ and $C_{V\!\!A}$
fails for the $3$\~vertex,} with $\vee$ or without
(i.e. with $\iota$ in $C$ or without).
Our formulas for (the numerators of) 
$\hat{\h}_{\,\circ\rightrightarrows\,}^{min}(\la^1,\la^2)^{top}
(t\mapsto q)$ seem different from the corresponding ones in
\cite{AKMV} for $3$\~links. It is possible that we have two 
different theories (based on different kinds of
Mehta-Macdonald identities). Recall that 
our approach is for arbitrary root systems 
$\iota=-w_0$ naturally emerges. It is in contrast to 
the theory from \cite{AKMV,AFS}, 
including the $K$\~theoretical
Nekrasov partition functions \cite{Nek}; it heavily involves
the transposition $\,\{\,\}^{tr}\,$ of diagrams, which
certainly requires {\em stable\,} $A_n$.  Thus at this moment, 
such a coincidence of two theories seems a special feature
of the $2$\~vertex.
\smallskip

Actually the Macdonald polynomials are
needed for the refined topological (physics) vertex
only a little; the {\em skew Schur functions\,}
are the main ingredient. See e.g., \cite{AFS}; the appearance
of stable Macdonald polynomials in the Nekrasov instanton 
sums from
formulas (5.1-3) there is via the evaluation of these polynomials
at $t^{\rho}$, which is a very explicit (and simple) product. 
Replacing the skew Schur function
by the corresponding Macdonald polynomials does not actually
influence the partition sums. 

The Hopf links naturally emerge and play an important
role in \cite{AKMV}, \cite{GIKV} and other papers
on topological vertex. However we do not quite 
understand  the connection of the superpolynomials 
for Hopf links with the (refined) topological vertex
considered in \cite{GIKV};  these directions look like two 
parallel but {\em different\,} theories to us (though both
are from the same physics source). In DAHA theory, these
two  are really closely related, the DAHA-vertex is part of
the DAHA-superpolynomials for {\em links}. However
our approach seems (so far) different from
the physics theory of the topological vertex. 
Mathematically, the DAHA-vertex seems almost inevitable, 
to combine the topological vertex and superpolynomials
in one theory (what we do).  
\medskip
 
{\sf Khovanov polynomials for Hopf links.}
Let us comment a bit on Section 3.1 from \cite{GIKV}
on the Hopf $2$\~link for the
fundamental representation, denoted by $L2a1$ in \cite{KA}.
We will use (\ref{khfromh}). Recall that the standard parameters 
there were modified vs. those in (\ref{qtareli}) as follows:
\begin{align}\label{newaqt}
t\mapsto q^2,\, q\mapsto (q t)^2,\, a\mapsto -a^2\,\
(a\mapsto -q^{2(n+1)} \for A_n). 
\end{align}
Only the substitution for $a$ is different here; it is
exactly the one from the HOMFLY-PT reduction vs. 
$a\mapsto a^2 t$ in the Connection Conjecture. Such
a modification of $a$ is very natural, since it automatically
provides the divisibility by  $(1-t)^{\kappa}$ (see below). 
It allows to deal directly with $\hat{\h}^{min}$, which is a 
polynomial in terms of $a,q,t^{\pm 1}$.

Let us switch in what will follow  to {\em new\,}
$q,t,a$ from (\ref{newaqt}).
Then the uncolored superpolynomial suggested in \cite{GIKV} 
is essentially obtained from our one by the 
conjugation $q,t,a\mapsto q^{-1},t^{-1},a^{-1}$
up to some simple factor. Upon 
the substitution $a=q^2$, it results in 
the Khovanov polynomial of this link ($L2a1$), 
which is $(1+q^2)(1+q^4 t^2)=1+q^2+q^4 t^2+q^6 t^2$
up to $q^\bullet t^\bullet$\~proportionality.

More systematically, in the {\em new parameters\,} above,
our superpolynomial  reads: 
$\hat{\h}^{min}_{L2a1}=1 - a^2 (1/q^2 - t^2 + q^2 t^2)$.
We obtain, which matches (41) in \cite{GIKV},
that: 
\begin{align}\label{hopf-2-gukov}
\Bigl(\frac{1\!-\!a^2}{(1\!-\!q^2)^2} 
\hat{\h}^{min}_{L2a1}\Bigr)_{a\,\mapsto q^N}\!=\! 
q^\bullet t^\bullet \hat{K\!h\!R}{}^N_{L2a1}(q,t) \for N=2,3,4,5.
\end{align}
Let us provide the left-hand side for $N=5$:
{\small
\(1+2 q^2+3 q^4+4 q^6+4 q^8+3 q^{10}+2 q^{12}+q^{14}+q^{10} t^2
+q^{12} t^2+q^{14} t^2+q^{16} t^2+q^{18} t^2.
\)
}

\medskip
Let us apply the same procedure to the Hopf $3$\~link,
denoted  by $L6n1\,=\, Link[\,6,\,NonAlternating,\,1\,]$
in \cite{KA}; we use 
the procedure $K\!hReduced[\cdot]$.
We generally expect that 
\begin{align}\label{hopf-3-gukov}
\Bigl(\frac{1-a^2}
{(1-q^2)^3} \hat{\h}^{min}_{L6n1}\Bigr)_{a\,\mapsto q^N}
\doteq 
q^\bullet t^\bullet \hat{K\!h\!R}{}^N_{L6n1}(q,t) \for N\gg 0.
\end{align}
For $N=2$, we obtain the following\,:
\begin{align*}
\frac{1-a^2}{(1-q^2)^3} \hat{\h}^{min}_{L6n1}(a=q^2)
\!=&(2+3 q^2+q^4-\underline{q^2 t^2}+q^6 t^2+q^8 t^4
+q^{10} t^4)q^2t^2,\\ 
\hat{K\!h\!R}{}^2_{L6n1}(q,t)\!=&
(2+3 q^2+q^4+\,\underline{q^2 t}\,+q^6 t^2+q^8 t^4+q^{10} t^4)/q.
\end{align*}
The change $-q^2t^2\mapsto q^2t$ (they are outlined)
is a typical correction
due to the nontrivial differential $\partial_2$.
Note that the $t$\~degree in this correction 
{\em diminishes\,} (by $1$); 
it always {\em increases\,} (here by $1$) in such corrections in 
the examples of algebraic links we considered.
\smallskip

We see that it is not impossible to obtain the 
reduced Khovanov polynomials $\hat{K\!h\!R}{}^2$ and other 
Khovanov- Rozansky polynomials directly
from $\hat{\h}^{min}$, i.e. without the division by 
$(1-t)^{\kappa-1}$.  
However we do not know how far this and similar procedures can 
go for links. In the examples we provide, managing the negative
terms in $\hat{\h}^{min}$ is always by 
$-q^lt^m\mapsto q^lt^{m\pm1}$,
$q^lt^m\mapsto q^lt^{m\pm2}$ and so on.
This is always (in examples) in the same direction
for a given link, which matches
the discussion of the usage of differentials
after the Connection Conjecture on \cite{ChD}. Employing the 
differentials provides sharper match here, but, as we have
seen, 
even a ``recovery" of reduced $\hat{K\!h\!R}{}^{n+1}$ from a 
single $\hat{\h}^{min}$ for the corresponding $A_n$ is not
impossible (for small links).
\smallskip

\subsection{\bf Algebraic Hopf links}
The matrix $\tau_-$, which played the key role above,
will be replaced by $\tau_+$ in this section. Recall that 
according to Part $(ii)$ of
the Connection Conjecture \ref{CONCONJ}, we expect the properties 
of the superpolynomials for {\em positive pairs\,} of trees
$\{\l,\,'\!\l^\vee\}$  
(with $\vee$ and satisfying the inequalities for $\rr,\ss$
in Part $(ii)$) to be the same as for a single
positive tree $\l$. Recall that the positivity of  
$\{\l,\,'\!\l^\vee\}$ is sufficient for the corresponding
link to be algebraic. It is necessary if the
{\em trees\,} $\l,\,'\!\l$ are {\em reduced\,}(minimal
modulo the equivalence); 
see \cite{EN}.
\smallskip

Before considering the algebraic links,
let us provide the simplest non-algebraic  {\em twisted union\,}, 
which is that for two disjoint (unlinked) 
unknots. I.e. $\Up$ is $\,\rightrightarrows\,$ or, 
equivalently, $\circ\rightrightarrows$ with the vertex 
labeled by $\{1,0\}$. Such a union is the result of 
adding the meridian having the linking numbers $+1$ with both 
unknots. This is the {\em 3-chain}, which is obviously
non-algebraic. Its superpolynomial is as follows:

\begin{align}\label{T10-10}
&
\Yboxdim5pt
\hbox{3-chain}\, :\ \,
\l=\l_{(\{1,0\},\{1,0\})}^{\,\circ\rightrightarrows,\,
(\yng(1)\,,\yng(1)\,)},
\,'\!\l= \l_{\{1,0\}}^{\,\circ\rightarrow,\,\yng(1)\,},\ \ 
\hat{\h}{}^{min}_{\l,\,'\!\l^\vee}\,(q,t,a)=
\end{align}

\renewcommand{\baselinestretch}{0.5} 
\noindent
{\small
\(
1+a^2 q^2-2 t+2 q t+t^2-2 q t^2+q^2 t^2+
a \bigl(2 q-2 q t+2 q^2 t\bigr)
\!=\!(1+a q-t+q t)^2.
\)
}
\renewcommand{\baselinestretch}{1.0} 
\smallskip

The square here is not surprising, since the
3-chain is the {\em connected sum\,} of two (positive) 
Hopf $2$\~links, with the DAHA-superpolynomials $(1+aq-t+qt)$.
Such a product formula holds if the meridian, which is
a colored unknot, is added to an arbitrary 
{\em disjoint union\,} of links. Let us justify this.

Recall the general formula $\{P,Q\}_{ev}=
\{\tau_-\bigl(\dot{\tau}_-^{-1}(P)
\dot{\tau}_-^{-1}(Q)\bigr)\}_{ev}$ from Theorem \ref{EVALTAUM}
for $P,Q\in \v^W$. 
Let $Q=\hat{\P}_{0}^{tot}$ and
$P=\,'\hat{\P}_{0}^{tot}$ be from (\ref{jones-hat})
for any {\em graphs\,} $\l,\,'\!\l$. Due to
the symmetry (\ref{iter-s1}), the components of the link 
for $\{\l,\,'\!\l^\vee\}$ is a certain union of  those for 
$\l$ and $'\!\l$. 
\smallskip

\begin{proposition}\label{CONECSUM}
Let $Q=\hat{\P}_{0}^{\l,tot}$,  $R=\hat{\P}_{0}^{\mathcal{M},tot}$
for arbitrary graphs  
$\l,\mathcal{M}$ and $\mathcal{N}$ be the disjoint union
of $\l,\mathcal{M}$. We take $\, \rightarrow\,$ as 
$\,'\!\l$; equivalently, $\,'\!\l=\,
'\!\circ\!\!\rightarrow\ $ where $\circ$
has the label $[1,0]$; it will be colored by 
any Young diagram $\la$. The corresponding 
$\hat{\h}_{\{\mathcal{N},\,'\!\l^\vee\}}$ will be normalized  
for $\la$, i.e. for $j_o$ being the index of 
$\,'\!\l$ (we put below $\la$ instead). So this is not
the standard $min$\~normalization. Then
\begin{align}\label{con-sum}
\hat{\h}^{\la}_{\{\mathcal{N},\,'\!\l^\vee\}}=
\hat{\h}^{\la}_{\{\mathcal{L},\,'\!\l^\vee\}}
\hat{\h}^{\la}_{\{\mathcal{M},\,'\!\l^\vee\}}.
\end{align}
\end{proposition}
{\it Proof.} 
$\hat{\h}^{\la}_{\{\mathcal{N},\,'\!\l^\vee\}}\!=\!
\{QR,\iota(P_{\la}^\circ)\}_{ev}\!=\!
QR(q^{\la+k\rho})\!=\!
\hat{\h}^{\la}_{\{\mathcal{L},\,'\!\l^\vee\}}
\hat{\h}^{\la}_{\{\mathcal{M},\,'\!\l^\vee\}}.
$
\vskip -0.6cm \sq

\smallskip

At the level of HOMFLY-PT polynomials (upon $a\!\mapsto\! -a, 
t\!\mapsto\! q$), relation
(\ref{con-sum}) becomes a well-known property of the {\em connected
sum\,} of $\tilde{\l},\tilde{\mathcal{M}}$, where by tilde,
we mean the twist with unknot. This sum
is $\tilde{\mathcal{N}}$ and (\ref{con-sum}) holds.
We cannot comment on the validity of such a product formula
in Khovanov-Rozansky theory.
\smallskip

{\sf Unusual positivity factors.} 
Let us demonstrate that the division by 
powers $(1-t)^\bullet$ can be insufficient to ensure
the positivity in Part $(ii)$  for colored algebraic links.
We begin with 

\begin{align}\label{T2-0+2x1-1}
&\Yboxdim7pt \yng(1,1)\,,2\!\times\! \yng(2)\, :\ \,\ \Yboxdim5pt
\l=\l_{\{1,1\}}^{\,\circ\rightthreearrow,\, 
(\yng(1,1)\,,\yng(2)\,,\yng(2)\,)},\ \ 
\hat{\h}{}^{min}_{\l}\,(q,t,a)=
\end{align}

\renewcommand{\baselinestretch}{0.5} 
\noindent
{\small
\(
1+a^3 q^9-t-q t+q^2 t+q^3 t-t^2+q t^2-q^2 t^2-q^3 t^2
+2 q^4 t^2+t^3+q t^3-q^2 t^3-2 q^3 t^3-q^4 t^3+q^5 t^3
+q^6 t^3-q t^4+q^2 t^4+2 q^3 t^4-2 q^4 t^4-q^5 t^4+q^6 t^4
+q^3 t^5-2 q^5 t^5+q^7 t^5-q^3 t^6+q^4 t^6+2 q^5 t^6-2 q^6 t^6
-q^7 t^6+q^8 t^6+a^2 \bigl(q^5+q^6+q^7-q^5 t-q^6 t+q^7 t+q^8 t
-q^5 t^2+q^9 t^2+q^5 t^3-2 q^7 t^3+q^9 t^3\bigr)+a \bigl(q^2
+q^3+q^4-q^2 t-2 q^3 t+2 q^5 t+q^6 t-q^2 t^2-q^4 t^2-q^5 t^2
+2 q^6 t^2+q^7 t^2+q^2 t^3+2 q^3 t^3-q^4 t^3-4 q^5 t^3-q^6 t^3
+2 q^7 t^3+q^8 t^3-q^3 t^4+q^4 t^4+2 q^5 t^4
-2 q^6 t^4-q^7 t^4+q^8 t^4+q^5 t^5-2 q^7 t^5+q^9 t^5\bigr).
\)
}
\renewcommand{\baselinestretch}{1.0} 
\smallskip

\Yboxdim5pt
Here deg$_a=3=|\yng(1,1)\,|+2|\yng(2)|-|\yng(2,1)|$ as
it is conjectured. This degree remains the same \Yboxdim5pt
for the triple $(\yng(1)\,,2\!\times\! \yng(2)\,)$; we will
omit the formula for $\hat{\h}^{min}$ for the latter.
The corresponding  series $\hat{\h}^{min}/(1-t)^p$ are
non-positive for any $p\ge 0$ (for both). However, the
series $\hat{\h}^{min}/(1-q)^3$ are positive and $3$
is the smallest here (for both). Thus $(1-t)^\bullet$
is sufficient to ensure the positivity
for $(\yng(1)\,,2\!\times\! \yng(1,1)\,)$, but any power
of $(1\!-\!t)$ is insufficient for $2\times \yng(2)\,$.
\smallskip

The next case demonstrates that  
$(1\!-\!t)(1\!-\!q)$ may be really needed.
\begin{align}\label{T2-1+2-1}
&\Yboxdim7pt \yng(2,1)\,,\yng(2,1)\, :\ \,\ \Yboxdim5pt
\l=\l_{\{1,1\}}^{\,\circ\rightrightarrows,\, 
(\yng(2,1)\,,\yng(2,1)\,)},\ \ 
\hat{\h}{}^{min}_{\l}\,(q,t,a)=
\end{align}

\renewcommand{\baselinestretch}{0.5} 
\noindent
{\small
\(
1+\frac{a^3 q^6}{t}-2 t+2 q t+t^2-4 q t^2+2 q^2 t^2+q^3 t^2
+3 q t^3-5 q^2 t^3+2 q^3 t^3-q t^4+3 q^2 t^4-4 q^3 t^4+2 q^4 t^4
+q^3 t^5-2 q^4 t^5+q^5 t^5+a^2 \bigl(-q^3+q^4+q^5+\frac{q^3}{t}
+\frac{q^4}{t}-2 q^4 t+q^5 t+q^6 t-q^5 t^2+q^6 t^2\bigr)
+a \bigl(-q+2 q^2+q^3+\frac{q}{t}-4 q^2 t+3 q^3 t+q^4 t+q^2 t^2
-5 q^3 t^2+3 q^4 t^2+q^5 t^2
+q^2 t^3+q^3 t^3-4 q^4 t^3+2 q^5 t^3-q^5 t^4+q^6 t^4\bigr).
\)
}
\renewcommand{\baselinestretch}{1.0} 
\smallskip

\Yboxdim5pt
The $a$\~degree remains $3=|\yng(2,1)\,|+|\yng(2,1)|-|\yng(2,1)|$.
The corresponding  series $\hat{\h}^{min}/(1-t)^p$ is
non-positive for any $p\ge 0$. Since $\hat{\h}^{min}$
is self-dual, the same non-positivity holds for $(1-q)^\bullet$. 
However, the series $\hat{\h}^{min}/((1\!-\!t)(1\!-\!q))^2$ 
is positive (as it was conjectured).
\smallskip

\subsection{\bf Twisted 
\texorpdfstring{{\mathversion{bold}$T(2\kappa,\kappa)$}}
{T(2k,k)}}\label{sec:T(2k,k)}
Let us discuss now in greater
detail the case of two trees when $'\!\l$ is  
{\em uncolored\,} $\circ\rightarrow$ with trivial
$\,'\rr=1,\, '\ss=0$. I.e. this will be the {\em twisted union\,}
with the uncolored unknot. In our notations, it is for
\Yboxdim5pt
$'\!\l=\l_{1,0}^{\circ\rightarrow,\,\yng(1)\,}$.
We will add ``prime" to the corresponding notation/name
of the link
without showing $\,'\!\l$ (with or without $\vee$).
Only the multiples of $T(2,1)$ will be considered in this section.
\smallskip
  
{\sf T(4,2)-prime for 3-hook and 1-box.}
\begin{align}\label{4-2-T21-1-1}
&T(4,2)_{\yng(2,1)\,,\,\yng(1)\,}^{\prime}\, :\ \,\Yboxdim5pt
\l=\l_{2,1}^{\,\circ\rightrightarrows,\, 
(\yng(2,1)\,,\,\yng(1)\,)},\ 
'\!\l=\l_{1,0}^{\circ\rightarrow,\,\yng(1)\,}\,,\ 
\hat{\h}{}^{min}_{\l,\,'\!\l}\,(q,t,a)=
\end{align}

\renewcommand{\baselinestretch}{0.5} 
\noindent
{\small
\(
1-t+q t-2 q t^2+2 q^2 t^2+q t^3-2 q^2 t^3+q^3 t^3-q^3 t^4+q^4 t^4
+a^3 \bigl(q^4-3 q^5+4 q^6-3 q^7+q^8-\frac{q^3}{t^4}+\frac{q^4}{t^4}
+\frac{2 q^3}{t^3}-\frac{3 q^4}{t^3}+\frac{q^6}{t^3}-\frac{q^3}{t^2}
+\frac{4 q^4}{t^2}-\frac{5 q^5}{t^2}+\frac{q^6}{t^2}+\frac{q^7}{t^2}
-\frac{3 q^4}{t}+\frac{8 q^5}{t}-\frac{5 q^6}{t}-q^6 t+2 q^7 t
-q^8 t\bigr)+a^2 \bigl(q^2-7 q^3+11 q^4-3 q^5-2 q^6+q^7
-\frac{q}{t^4}+\frac{q^2}{t^4}+\frac{2 q}{t^3}-\frac{3 q^2}{t^3}
+\frac{q^4}{t^3}-\frac{q}{t^2}+\frac{5 q^2}{t^2}-\frac{4 q^3}{t^2}
-\frac{2 q^4}{t^2}+\frac{q^5}{t^2}+\frac{q^6}{t^2}-\frac{4 q^2}{t}
+\frac{10 q^3}{t}-\frac{3 q^4}{t}-\frac{4 q^5}{t}+\frac{q^6}{t}
+q^3 t-7 q^4 t+10 q^5 t-4 q^6 t+q^4 t^2-4 q^5 t^2+5 q^6 t^2
-3 q^7 t^2+q^8 t^2-q^6 t^3+2 q^7 t^3-q^8 t^3\bigr)+a \bigl(-3 q
+6 q^2-q^4+\frac{1}{t^2}-\frac{q}{t^2}-\frac{1}{t}+\frac{3 q}{t}
-\frac{q^2}{t}-\frac{q^3}{t}+q t-6 q^2 t+6 q^3 t-q^5 t+q^2 t^2
-6 q^3 t^2+6 q^4 t^2
-q^5 t^2+q^3 t^3-3 q^4 t^3+3 q^5 t^3-q^6 t^3-q^5 t^4+q^6 t^4\bigr).
\)
}
\renewcommand{\baselinestretch}{1.0} 

\begin{align}\label{4-2-T21-1-1v}
&T(4,2)_{\yng(2,1)\,,\,\yng(1)\,}^{\prime,\vee}\, :\ \,\Yboxdim5pt
\l=\l_{2,1}^{\,\circ\rightrightarrows,\, 
(\yng(2,1)\,,\,\yng(1)\,)},\ 
'\!\l=\l_{1,0}^{\circ\rightarrow,\,\yng(1)\,}\,,\ 
\hat{\h}{}^{min}_{\l,\,'\!\l^\vee}\,(q,t,a)=
\end{align}

\renewcommand{\baselinestretch}{0.5} 
\noindent
{\small
\(
1-2 t+q t+q^2 t+t^2-3 q t^2-q^2 t^2+2 q^3 t^2+q^4 t^2+3 q t^3
-2 q^2 t^3-4 q^3 t^3+2 q^4 t^3+q^5 t^3-q t^4+3 q^2 t^4+q^3 t^4
-7 q^4 t^4+2 q^5 t^4+2 q^6 t^4-q^2 t^5+2 q^3 t^5+3 q^4 t^5
-7 q^5 t^5
+2 q^6 t^5+q^7 t^5-q^3 t^6+2 q^4 t^6+3 q^5 t^6-7 q^6 t^6+2 q^7 t^6
+q^8 t^6-q^4 t^7+2 q^5 t^7+q^6 t^7-4 q^7 t^7+2 q^8 t^7-q^5 t^8
+3 q^6 t^8-2 q^7 t^8-q^8 t^8+q^9 t^8-q^6 t^9+3 q^7 t^9-3 q^8 t^9
+q^9 t^9+q^8 t^{10}-2 q^9 t^{10}+q^{10} t^{10}+
a^2 \bigl(q^5-2 q^5 t
+q^6 t+q^7 t+q^5 t^2-3 q^6 t^2+q^7 t^2+q^8 t^2+3 q^6 t^3-5 q^7 t^3
+q^8 t^3+q^9 t^3-q^6 t^4+3 q^7 t^4-3 q^8 t^4+q^9 t^4+q^8 t^5
-2 q^9 t^5+q^{10} t^5\bigr)+a \bigl(q^2+q^3-2 q^2 t-q^3 t+2 q^4 t
+q^5 t+q^2 t^2-2 q^3 t^2-4 q^4 t^2+3 q^5 t^2+2 q^6 t^2+3 q^3 t^3
+q^4 t^3-9 q^5 t^3+3 q^6 t^3+2 q^7 t^3-q^3 t^4+2 q^4 t^4+4 q^5 t^4
-10 q^6 t^4+3 q^7 t^4+2 q^8 t^4-q^4 t^5+2 q^5 t^5+4 q^6 t^5
-9 q^7 t^5+3 q^8 t^5+q^9 t^5-q^5 t^6+2 q^6 t^6+q^7 t^6-4 q^8 t^6
+2 q^9 t^6-q^6 t^7+3 q^7 t^7
-2 q^8 t^7-q^9 t^7+q^{10} t^7+q^8 t^8-2 q^9 t^8+q^{10} t^8\bigr).
\)
}
\renewcommand{\baselinestretch}{1.0} 
\smallskip

The $a$\~degree in the first formula is $3$, i.e. 
greater than that in the second (with $\vee$):
deg$_a=\ss(|\yng(2,1)\,|+|\yng(1)\,|)+
\max\{\,'\ss,1\}|\yng(1)\,|-|\yng(2,1)\,|=2.$ This is
an important confirmation of the role of $\vee$
in our theorems/conjectures.
With $\vee$ here, the positivity
of $\hat{\h}{}^{min}_{\l,\,'\!\l^\vee}(q,t,a)/(1-t)^3$ holds,
but it fails for $(1-t)^2$. The term $(1-t)^3$ agrees
with the informal observation
above that the proper power of $(1-t)$ here can be 
$(1-t)^{\kappa+\,'\kappa -1}$ times $(1-t)$, which is
the extra correction factor for $\yng(2,1)\,$.
\medskip

{\sf T(6,3)-prime for 3 boxes.}  
\begin{align}\label{6-3-T1-1-1}
&T(6,3)_{\yng(1)\,,\,\yng(1)\,,\,\yng(1)}^{\prime}\, 
: \,\Yboxdim5pt
\l\!=\!\l_{2,1}^{\,\circ\rightthreearrow,\, 
(\yng(1)\,,\,\yng(1)\,,\,\yng(1)\,)},\, 
'\!\l=\l_{1,0}^{\circ\rightarrow,\,\yng(1)\,}\,,\, 
\hat{\h}{}^{min}_{\l,\,'\!\l}\,(q,t,a)\!=
\end{align}

\renewcommand{\baselinestretch}{0.5} 
\noindent
{\small
\(
1+a^4 \bigl(-3 q^6+6 q^7-3 q^8-\frac{3 q^6}{t^2}+\frac{3 q^7}{t^2}
+\frac{6 q^6}{t}-\frac{9 q^7}{t}+\frac{3 q^8}{t}\bigr)-2 t+q t
+q^2 t+t^2-2 q t^2-q^2 t^2+q^3 t^2+q^4 t^2+q t^3-q^2 t^3-q^3 t^3
+q^4 t^3+q^2 t^4-q^3 t^4-q^4 t^4+q^5 t^4+q^3 t^5-2 q^4 t^5+q^5 t^5
+q^4 t^6-2 q^5 t^6+q^6 t^6+a^3 \bigl(-2 q^3-2 q^4+3 q^5+12 q^6
-13 q^7+3 q^8-\frac{q^4}{t^2}-\frac{3 q^5}{t^2}+\frac{3 q^6}{t^2}
+\frac{q^7}{t^2}+\frac{q^3}{t}+\frac{3 q^4}{t}+\frac{5 q^5}{t}
-\frac{13 q^6}{t}+\frac{3 q^7}{t}+\frac{q^8}{t}+q^3 t-q^4 t-5 q^5 t
+3 q^6 t+5 q^7 t-3 q^8 t+q^4 t^2-q^5 t^2-2 q^6 t^2+3 q^7 t^2
-q^8 t^2
+q^5 t^3-2 q^6 t^3+q^7 t^3\bigr)+a^2 \bigl(-2 q-q^2+9 q^4+4 q^5
-8 q^6+q^8-\frac{q^3}{t^2}+\frac{q^6}{t^2}+\frac{q}{t}
+\frac{q^2}{t}
+\frac{3 q^3}{t}-\frac{2 q^4}{t}-\frac{4 q^5}{t}+\frac{q^7}{t}
+q t-q^2 t-4 q^3 t-6 q^4 t+10 q^5 t+4 q^6 t-4 q^7 t+q^2 t^2
+q^3 t^2-3 q^4 t^2-6 q^5 t^2+9 q^6 t^2-2 q^7 t^2+q^3 t^3+q^4 t^3
-4 q^5 t^3+3 q^7 t^3-q^8 t^3+q^4 t^4-q^5 t^4-q^6 t^4+q^7 t^4
+q^5 t^5-2 q^6 t^5+q^7 t^5\bigr)+a \bigl(-2+2 q+2 q^2+3 q^3-q^4
-q^5+\frac{1}{t}-\frac{q^3}{t}+t-4 q t-2 q^2 t+7 q^4 t-q^5 t
-q^6 t+2 q t^2-2 q^2 t^2-3 q^3 t^2-3 q^4 t^2+7 q^5 t^2-q^6 t^2
+2 q^2 t^3-q^3 t^3-3 q^4 t^3
+3 q^6 t^3-q^7 t^3+2 q^3 t^4-2 q^4 t^4-2 q^5 t^4+2 q^6 t^4
+2 q^4 t^5-4 q^5 t^5+2 q^6 t^5+q^5 t^6-2 q^6 t^6+q^7 t^6\bigr).
\)
}
\renewcommand{\baselinestretch}{1.0} 

\begin{align}\label{6-3-T1-1-1v}
&T(6,3)_{\yng(1)\,,\,\yng(1)\,,\,\yng(1)}^{\prime,\vee}\, 
:\,\Yboxdim5pt
\l\!=\!\l_{2,1}^{\,\circ\rightthreearrow,\, 
(\yng(1)\,,\,\yng(1)\,,\,\yng(1)\,)},\, 
'\!\l=\l_{1,0}^{\circ\rightarrow,\,\yng(1)\,}\,,\,
\hat{\h}{}^{min}_{\l,\,'\!\l^\vee}\,(q,t,a)\!=
\end{align}

\renewcommand{\baselinestretch}{0.5} 
\noindent
{\small
\(
1-3 t+q t+q^2 t+q^3 t+3 t^2-3 q t^2-2 q^2 t^2-2 q^3 t^2+2 q^4 t^2
+q^5 t^2+q^6 t^2-t^3+3 q t^3+q^3 t^3-5 q^4 t^3-q^5 t^3+2 q^6 t^3
+q^7 t^3-q t^4+2 q^2 t^4-q^3 t^4+4 q^4 t^4-2 q^5 t^4-6 q^6 t^4
+2 q^7 t^4+2 q^8 t^4-q^2 t^5+2 q^3 t^5-2 q^4 t^5+3 q^5 t^5+q^6 t^5
-6 q^7 t^5+2 q^8 t^5+q^9 t^5-q^3 t^6+2 q^4 t^6-2 q^5 t^6+3 q^6 t^6
+q^7 t^6-6 q^8 t^6+2 q^9 t^6+q^{10} t^6-q^4 t^7+2 q^5 t^7-2 q^6 t^7
+3 q^7 t^7-2 q^8 t^7-q^9 t^7+q^{10} t^7-q^5 t^8+2 q^6 t^8-2 q^7 t^8
+4 q^8 t^8-5 q^9 t^8+2 q^{10} t^8-q^6 t^9+2 q^7 t^9-q^8 t^9+q^9 t^9
-2 q^{10} t^9+q^{11} t^9-q^7 t^{10}+2 q^8 t^{10}-2 q^{10} t^{10}
+q^{11} t^{10}-q^8 t^{11}+3 q^9 t^{11}-3 q^{10} t^{11}
+q^{11} t^{11}-q^9 t^{12}+3 q^{10} t^{12}-3 q^{11} t^{12}
+q^{12} t^{12}+a^3 \bigl(q^6-3 q^6 t+q^7 t+q^8 t+q^9 t+3 q^6 t^2
-3 q^7 t^2-2 q^8 t^2+q^9 t^2+q^{10} t^2-q^6 t^3+3 q^7 t^3
-4 q^9 t^3+q^{10} t^3
+q^{11} t^3-q^7 t^4+2 q^8 t^4-2 q^{10} t^4+q^{11} t^4-q^8 t^5
+3 q^9 t^5-3 q^{10} t^5+q^{11} t^5-q^9 t^6+3 q^{10} t^6
-3 q^{11} t^6
+q^{12} t^6\bigr)+a^2 \bigl(q^3+q^4+q^5-3 q^3 t-2 q^4 t-q^5 t
+3 q^6 t+2 q^7 t+q^8 t+3 q^3 t^2-2 q^5 t^2-7 q^6 t^2+q^7 t^2
+3 q^8 t^2+2 q^9 t^2-q^3 t^3+2 q^4 t^3+2 q^5 t^3+4 q^6 t^3
-9 q^7 t^3
-4 q^8 t^3+4 q^9 t^3+2 q^{10} t^3-q^4 t^4+q^5 t^4+6 q^7 t^4
-6 q^8 t^4-4 q^9 t^4+3 q^{10} t^4+q^{11} t^4-q^5 t^5+q^6 t^5
+6 q^8 t^5-9 q^9 t^5+q^{10} t^5+2 q^{11} t^5-q^6 t^6+q^7 t^6
+4 q^9 t^6-7 q^{10} t^6+3 q^{11} t^6-q^7 t^7+q^8 t^7+2 q^9 t^7
-2 q^{10} t^7-q^{11} t^7+q^{12} t^7-q^8 t^8+2 q^9 t^8-2 q^{11} t^8
+q^{12} t^8-q^9 t^9+3 q^{10} t^9-3 q^{11} t^9+q^{12} t^9\bigr)
+a \bigl(q+q^2+q^3-3 q t-2 q^2 t-q^3 t+3 q^4 t+2 q^5 t+q^6 t
+3 q t^2
-2 q^3 t^2-7 q^4 t^2-2 q^5 t^2+4 q^6 t^2+3 q^7 t^2+q^8 t^2-q t^3
+2 q^2 t^3+2 q^3 t^3+4 q^4 t^3-4 q^5 t^3-10 q^6 t^3+2 q^7 t^3
+4 q^8 t^3+q^9 t^3-q^2 t^4+q^3 t^4+5 q^5 t^4+2 q^6 t^4-12 q^7 t^4
+4 q^9 t^4+q^{10} t^4-q^3 t^5+q^4 t^5-q^5 t^5+4 q^6 t^5+4 q^7 t^5
-12 q^8 t^5+2 q^9 t^5+3 q^{10} t^5-q^4 t^6+q^5 t^6-q^6 t^6
+4 q^7 t^6
+2 q^8 t^6-10 q^9 t^6+4 q^{10} t^6+q^{11} t^6-q^5 t^7+q^6 t^7
-q^7 t^7+5 q^8 t^7-4 q^9 t^7-2 q^{10} t^7+2 q^{11} t^7-q^6 t^8
+q^7 t^8+4 q^9 t^8-7 q^{10} t^8+3 q^{11} t^8-q^7 t^9+q^8 t^9
+2 q^9 t^9-2 q^{10} t^9-q^{11} t^9+q^{12} t^9-q^8 t^{10}
+2 q^9 t^{10}-2 q^{11} t^{10}+q^{12} t^{10}
-q^9 t^{11}+3 q^{10} t^{11}-3 q^{11} t^{11}+q^{12} t^{11}\bigr).
\)
}
\renewcommand{\baselinestretch}{1.0} 
\smallskip

The $a$\~degree in the first formula is $4$, which is
(again) greater than that in the second :
deg$_a=\ss(3|\yng(1)\,|)+
\max\{\,'\ss,1\}|\yng(1)\,|-|\yng(1)\,|=3.$ 
The positivity of
$\hat{\h}{}^{min}_{\l,\,'\!\l}(q,t,a)/(1-t)^3$ holds,
but fails for $(1-t)^2$, this agrees with 
$\kappa+\,'\kappa-1=3$
from the Connection Conjecture.
Both polynomials are self-dual and become $(1+a)^3$ as $q=1$.
Using the spherical polynomials would result at $1$ for $q=1$
(since all components are unknots); the factor $(1+a)^{4-1}$
is due to using $J$\~polynomials and the $min$\~normalization.
\smallskip

The corresponding plane curve singularity is:

$$
\c=\{y^8+x^4=0\}, \,\, \hbox{Link}(\c)=
\bigl\{\,4\times \unknot\,,\,\, \hbox{lk}_{i\neq j}=2\,\bigr\}.
$$

Its Alexander polynomial is $1+q^4-2 q^8-2 q^{12}+q^{16}+q^{20}$.
See Section \ref{sec:Alex} below for details.
\smallskip

This link is isotopic to uncolored 
$(4\times C\!ab(2,1))T(1,1)$. The
corresponding $\hat{\h}{}^{min}_{\m}\,(q,t,a)$ for
$\m=\l^{\circ\rightarrow\circ\,\rightfourarrow,\, 
(4\times\yng(1)\,)}_{\{1,1\},\{1,1\}}$ coincides with that
from (\ref{6-3-T1-1-1v}), which also coincides with  
$\hat{\h}{}^{min}_{\n}\,(q,t,a)$ for 
$\n=\l^{\circ\,\rightfourarrow,\, 
(4\times\yng(1)\,)}_{\{2,1\}}$. Its link is
$(4\times C\!ab(0,1))T(2,1)$, isotopic to the previous two.

\smallskip
Let us check directly
the coincidence of the superpolynomials 
for the original $\{\l,\,'\!\l\}$ and $\n$ above. For $E=E_\la,$
$P=P_\la$, where $\la\in P,P_+$ respectively, one has:
\begin{align}\label{iden2-1}
&\{E(Y^{-1})\, \tau_+^{-1}\tau_-^{-1}(Q)\}_{ev}\simeq
\{\tau_+^{-1}\tau_-^{-1}(EQ)\}_{ev},\ Q\in \v,\\ 
&\{\,P(Y)\, \tau_+\tau_-(Q)\,\}_{ev}\simeq
\{\,\tau_+\tau_-(P\,Q)\,\}_{ev} \for Q\in \v^W\!\!, \notag
\end{align}
though one can assume that $Q\in \v$ in both relations.
By $\simeq$\,, we mean the equality up to a factor 
$q^\bullet t^\bullet$.
The second relation follows from the first upon applying 
$\,\eta\,$ inside the coinvariant.
See (\ref{eta-tau}), (\ref{evsym}) and  Section \ref{sect:Aut}
for this and other properties we will need.

Let us check the first equality. We use the identity
$\tau_+\tau_-^{-1}\tau_+=\si=\tau_-^{-1}\tau_+\tau_-^{-1}$, 
its corollary 
$\tau_+^{-1}\tau_-^{-1}\tau_+=\tau_-^{-1}\tau_+\tau_-$,
the $\vph$\~invariance of $\{\cdot\}_{ev}$, and also
the relation  $\vph(Q)=\si(Q)$. One has:
\begin{align*}
&\{E(Y^{-1})\,(\tau_+^{-1}\tau_-^{-1})(Q)\}_{ev}=
\{\vph\Bigl(E(Y^{-1})\,
(\tau_+^{-1}\tau_-^{-1}\tau_+)(Q)\Bigr)\}_{ev}\\
=
&\{(\tau_+^{-1}\tau_-\tau_+)(\tau_-^{-1}\tau_+\tau_-^{-1})
(Q)\,E\}_{ev}\!=\!
\{(\tau_+^{-1}\tau_-(\tau_+\tau_-^{-1}\tau_+)\tau_-^{-1})
(Q)\,E\}_{ev}\\
=&\{\tau_-^{-2}(Q)\,E\}_{ev}\,\simeq\,\{\tau_-^{-2}(QE)\}_{ev}\,=\,
\{(\tau_+^{-1}\tau_-^{-1})(QE)\}_{ev}.
\end{align*}
\smallskip

\setcounter{equation}{0}
\section{\sc Further examples}

\subsection{\bf On cable notations}
Let us somewhat simplify the notations. If the tree
is one path extended by several arrowheads at its end, then
the notation without
the duplication of the sets $(\vec\rr^j,\vec\ss^j)$
can be used, since the corresponding paths are different
only by the arrowheads (i.e. by colors). 
More generally, in the case of one {\em base path\,} with
several arrowheads going from its vertices (not only from
its very last vertex), we will provide 
the \tax-parameters only for this path. 
The other paths will be its portions that
and the arrowheads before the final vertex.
\smallskip

For example,
for the tree $\Up=\{\twoone\}$ with 2 vertices and
2 arrowheads, the link
$\l^{\twoonetiny}_{\{3,2\},\{-2,1\}}$ means 
combinatorially and topologically:
$$
\l^{\twoonetiny,\, (\la^1,\la^2)}_{(\bigl\{\{3,2\},\{-2,1\}\bigr\},
\bigl\{\{3,2\}\bigr\})}\, \sim\,
\bigl(C\!ab(-11,2)^{\la^1}\coprod\, C\!ab(0,1)^{\la^2}\bigr)T(3,2),
$$
where we count the paths from the bottom to the top. Recall that 
the labels $[\rr,\ss]$ corresponding here
to $\vec\rr=\{3,2\},\vec\ss=\{-2,1\}$ are $[3,-2]$ and $[2,1]$,
so they can be naturally considered as the columns of 
{\tiny $\begin{pmatrix}
3&\!\!2\\ 
-2&\!\!1
\end{pmatrix}$}; these columns become 
the first columns of the corresponding $\ga$. 
To avoid any  misunderstanding in some examples,
we will continue
duplicating $\{\vec{\rr},\vec{\ss}\}$ for the
paths different only by the arrowheads. 
Also, the symbol $\,\hbox{\small$\coprod$}\,$
will be mainly omitted below if the confusion with the 
composition of cables is impossible.
\smallskip

\subsection{\bf Colored 2-links}
Actually we have already considered quite a few colored
examples in Section \ref{sec:torus-hopf}. Let us add
some simple colors to the $2$\~fold trefoil.
\smallskip
 
{\sf 2-column and 2-row for trefoil.} The \tax-presentation  
and $\hat{\h}^{min}$ are as follows:
$$
1\le j\le \kappa=2,\ \vec\rr^j=3,\, \vec\ss^j=2,\, 
\Up=\{\circ\rightrightarrows\}\,,\ 
\Yboxdim7pt \la^1=\yng(2)\,,\la^2=\yng(1,1)\,;
$$
\Yboxdim5pt
\begin{align}\label{T6-4-11-20}
&T(6,4):\ 
\l=\l_{(\{3,2\},\{3,2\})}^{\,\circ\rightrightarrows,\, 
(\yng(2)\,,\yng(1,1)\,)},\ \ 
\hat{\h}{}^{min}_{\l}\,(q,t,a)=
\end{align}

\renewcommand{\baselinestretch}{0.5} 
\noindent
{\small
\(
1+q t+q^2 t+q^3 t-t^2+q^3 t^2+3 q^4 t^2+q^5 t^2-q t^3-q^2 t^3
-q^3 t^3+3 q^5 t^3+2 q^6 t^3+q^7 t^3-q^3 t^4-3 q^4 t^4-q^5 t^4
+3 q^6 t^4+3 q^7 t^4+q^8 t^4-3 q^5 t^5-2 q^6 t^5+q^7 t^5+3 q^8 t^5
+3 q^9 t^5-3 q^6 t^6-3 q^7 t^6+2 q^8 t^6+3 q^9 t^6+2 q^{10} t^6
-2 q^7 t^7-3 q^8 t^7-q^9 t^7+3 q^{10} t^7+3 q^{11} t^7-3 q^8 t^8
-3 q^9 t^8+q^{10} t^8+3 q^{11} t^8+2 q^{12} t^8-2 q^9 t^9
-3 q^{10} t^9-q^{11} t^9+3 q^{12} t^9+3 q^{13} t^9-3 q^{10} t^{10}
-3 q^{11} t^{10}+2 q^{12} t^{10}+3 q^{13} t^{10}+q^{14} t^{10}
-2 q^{11} t^{11}-3 q^{12} t^{11}+q^{13} t^{11}+3 q^{14} t^{11}
+q^{15} t^{11}-3 q^{12} t^{12}-2 q^{13} t^{12}+3 q^{14} t^{12}
+2 q^{15} t^{12}-3 q^{13} t^{13}-q^{14} t^{13}+3 q^{15} t^{13}
+q^{16} t^{13}-3 q^{14} t^{14}+3 q^{16} t^{14}-q^{14} t^{15}
-q^{15} t^{15}+q^{16} t^{15}+q^{17} t^{15}-q^{15} t^{16}
+q^{17} t^{16}-q^{15} t^{17}+q^{17} t^{17}-q^{16} t^{18}
+q^{18} t^{18}
\)
\vfil

\noindent
\(
+a^5 \bigl(\frac{q^{15}}{t}-q^{15} t+q^{17} t
-q^{17} t^3+q^{19} t^3\bigr)+a^4 \bigl(q^{12}+q^{13}+q^{14}
+\frac{q^{10}}{t}+\frac{q^{11}}{t}+\frac{q^{12}}{t}
+\frac{q^{13}}{t}-q^{10} t-q^{11} t-q^{12} t+2 q^{14} t+2 q^{15} t
-q^{12} t^2-q^{13} t^2+2 q^{15} t^2+q^{16} t^2-q^{13} t^3
-2 q^{14} t^3-q^{15} t^3+2 q^{16} t^3+2 q^{17} t^3-q^{14} t^4
-2 q^{15} t^4+2 q^{17} t^4+q^{18} t^4-q^{15} t^5-q^{16} t^5
+q^{18} t^5+q^{19} t^5-q^{16} t^6-q^{17} t^6+q^{18} t^6+q^{19} t^6
-q^{17} t^7+q^{19} t^7-q^{17} t^8+q^{19} t^8\bigr)
\)
\vfil

\noindent
\(
+a^3 \bigl(q^8
+3 q^9+3 q^{10}+2 q^{11}+q^{12}+\frac{q^6}{t}+\frac{q^7}{t}
+\frac{2 q^8}{t}+\frac{q^9}{t}+\frac{q^{10}}{t}-q^6 t-q^7 t
-2 q^8 t-q^9 t+3 q^{10} t+4 q^{11} t+4 q^{12} t+q^{13} t-q^8 t^2
-3 q^9 t^2-3 q^{10} t^2+q^{11} t^2+5 q^{12} t^2+5 q^{13} t^2
+2 q^{14} t^2-4 q^{10} t^3-4 q^{11} t^3-q^{12} t^3+4 q^{13} t^3
+6 q^{14} t^3+q^{15} t^3-3 q^{11} t^4-6 q^{12} t^4-2 q^{13} t^4
+4 q^{14} t^4+6 q^{15} t^4+2 q^{16} t^4-3 q^{12} t^5-5 q^{13} t^5
-2 q^{14} t^5+4 q^{15} t^5+5 q^{16} t^5+q^{17} t^5-3 q^{13} t^6
-6 q^{14} t^6-q^{15} t^6+5 q^{16} t^6+4 q^{17} t^6+q^{18} t^6
-3 q^{14} t^7-4 q^{15} t^7+q^{16} t^7+4 q^{17} t^7+2 q^{18} t^7
-4 q^{15} t^8-3 q^{16} t^8+3 q^{17} t^8+3 q^{18} t^8+q^{19} t^8
-3 q^{16} t^9-q^{17} t^9+3 q^{18} t^9+q^{19} t^9-q^{16} t^{10}
-2 q^{17} t^{10}+q^{18} t^{10}+2 q^{19} t^{10}-q^{17} t^{11}
+q^{19} t^{11}-q^{17} t^{12}+q^{19} t^{12}\bigr)
\)
\vfil

\noindent
\(
+a^2 \bigl(2 q^5
+3 q^6+4 q^7+2 q^8+q^9+\frac{q^3}{t}+\frac{q^4}{t}+\frac{q^5}{t}
+\frac{q^6}{t}-q^3 t-q^4 t-q^5 t+4 q^7 t+7 q^8 t+4 q^9 t+2 q^{10} t
-2 q^5 t^2-3 q^6 t^2-4 q^7 t^2+q^8 t^2+8 q^9 t^2+8 q^{10} t^2
+4 q^{11} t^2+q^{12} t^2-q^6 t^3-4 q^7 t^3-7 q^8 t^3-2 q^9 t^3
+7 q^{10} t^3+9 q^{11} t^3+5 q^{12} t^3-3 q^8 t^4-9 q^9 t^4
-6 q^{10} t^4+5 q^{11} t^4+10 q^{12} t^4+6 q^{13} t^4+q^{14} t^4
-2 q^9 t^5-9 q^{10} t^5-7 q^{11} t^5+4 q^{12} t^5+10 q^{13} t^5
+6 q^{14} t^5-2 q^{10} t^6-9 q^{11} t^6-9 q^{12} t^6+4 q^{13} t^6
+10 q^{14} t^6+5 q^{15} t^6+q^{16} t^6-2 q^{11} t^7-9 q^{12} t^7
-7 q^{13} t^7+5 q^{14} t^7+9 q^{15} t^7+4 q^{16} t^7-2 q^{12} t^8
-9 q^{13} t^8-6 q^{14} t^8+7 q^{15} t^8+8 q^{16} t^8+2 q^{17} t^8
-2 q^{13} t^9-9 q^{14} t^9-2 q^{15} t^9+8 q^{16} t^9+4 q^{17} t^9
+q^{18} t^9-3 q^{14} t^{10}-7 q^{15} t^{10}+q^{16} t^{10}
+7 q^{17} t^{10}+2 q^{18} t^{10}-4 q^{15} t^{11}-4 q^{16} t^{11}
+4 q^{17} t^{11}+4 q^{18} t^{11}-q^{15} t^{12}-3 q^{16} t^{12}
+3 q^{18} t^{12}+q^{19} t^{12}-2 q^{16} t^{13}-q^{17} t^{13}
+2 q^{18} t^{13}+q^{19} t^{13}-q^{17} t^{14}+q^{19} t^{14}
-q^{17} t^{15}+q^{19} t^{15}\bigr)
\)
\vfil

\noindent
\(
+a \bigl(q^2+2 q^3+2 q^4+q^5
+\frac{q}{t}-q t+q^3 t+2 q^4 t+5 q^5 t+3 q^6 t+q^7 t-q^2 t^2
-2 q^3 t^2-2 q^4 t^2+6 q^6 t^2+7 q^7 t^2+3 q^8 t^2+q^9 t^2-q^3 t^3
-2 q^4 t^3-5 q^5 t^3-3 q^6 t^3+4 q^7 t^3+8 q^8 t^3+6 q^9 t^3
+q^{10} t^3-q^5 t^4-6 q^6 t^4-7 q^7 t^4+2 q^8 t^4+8 q^9 t^4
+7 q^{10} t^4+2 q^{11} t^4-5 q^7 t^5-8 q^8 t^5-2 q^9 t^5
+8 q^{10} t^5+9 q^{11} t^5+2 q^{12} t^5-5 q^8 t^6-9 q^9 t^6
-2 q^{10} t^6+7 q^{11} t^6+8 q^{12} t^6+2 q^{13} t^6-4 q^9 t^7
-9 q^{10} t^7-5 q^{11} t^7+7 q^{12} t^7+9 q^{13} t^7+2 q^{14} t^7
-5 q^{10} t^8-9 q^{11} t^8-2 q^{12} t^8+8 q^{13} t^8+7 q^{14} t^8
+q^{15} t^8-4 q^{11} t^9-9 q^{12} t^9-2 q^{13} t^9+8 q^{14} t^9
+6 q^{15} t^9+q^{16} t^9-5 q^{12} t^{10}-8 q^{13} t^{10}
+2 q^{14} t^{10}+8 q^{15} t^{10}+3 q^{16} t^{10}-5 q^{13} t^{11}
-7 q^{14} t^{11}+4 q^{15} t^{11}+7 q^{16} t^{11}+q^{17} t^{11}
-6 q^{14} t^{12}-3 q^{15} t^{12}+6 q^{16} t^{12}+3 q^{17} t^{12}
-q^{14} t^{13}-5 q^{15} t^{13}+5 q^{17} t^{13}+q^{18} t^{13}
-2 q^{15} t^{14}-2 q^{16} t^{14}+2 q^{17} t^{14}+2 q^{18} t^{14}
-q^{15} t^{15}-2 q^{16} t^{15}+q^{17} t^{15}+2 q^{18} t^{15}
-q^{16} t^{16}+q^{18} t^{16}-q^{17} t^{17}+q^{19} t^{17}\bigr).
\)
}
\renewcommand{\baselinestretch}{1.0} 

The $a$\~degree of $\hat{\h}{}^{min}_{\l}\,(q,t,a)$
is $5$, which matches the formula 
deg$_a=\ss(|\yng(2)|+|\yng(1,1)|)-|\yng(2,1)|=5$
from (\ref{deg-a-jj}). The self-duality and other claims
in this theorem hold. For instance, the transposition
$6\leftrightarrow 4$ in
$T(6,4)$ does not influence the superpolynomial.

The positivity of 
$\hat{\h}{}^{min}_{\l}(q,t,a)/(1-t)$ from 
(\ref{posit-claim}) holds;
though it was conjectured with such a power of $(1-t)$
only in the uncolored case.

Concerning the HOMFLY-PT polynomials, the following holds:
\begin{align*}
&\hat{\h}^{min}_{\l}\,
(q,q,-a)/den^1=
\hat{\hbox{\small H\!O\!M}}^{j_o=1}_{\l}\,(q,a).
\end{align*}
Here 
$den^1(\yng(2)\,,\yng(1,1))=
\frac{(1-q)(1-q^2)}{(1-a/q)}$. Otherwise in this section:
\Yboxdim5pt
$$
den^1(\la^1,\raisebox{1pt}{\yng(1)})=(1-q),\,
den^1(\la^1,\raisebox{1pt}{\yng(2)})=(1-q)(1-q^2)=
den^1(\la^1,\raisebox{1pt}{\yng(1,1)}).
$$
Note that $den^1$ depends only on the colors and 
therefore is the same 
as for the colored Hopf links in (\ref{den1-two}).
\medskip

{\sf 2-fold T(2,1)\,}. One has:
$$
1\le j\le \kappa=2,\ \vec\rr^j=2,\, \vec\ss^j=1,\, 
\Up=\{\circ\rightrightarrows\}\,,\ 
\Yboxdim7pt \la^1=\yng(2)\,,\la^2=\yng(1,1)\,;
$$
\Yboxdim5pt
\begin{align}\label{T4-2-11-20}
&T(4,2):\ 
\l=\l_{(\{2,1\},\{2,1\})}^{\,\circ\rightrightarrows,\, 
(\yng(2)\,,\yng(1,1)\,)},\ \ 
\hat{\h}{}^{min}_{\l}\,(q,t,a)=
\end{align}
{\small
$$
1-t^2+q^2 t^2-q^2 t^4+q^4 t^4+a(q^2-q^2 t^2+q^4 t^2).
$$
}

The $a$\~degree is $\ss(|\la^1|+|\la^2|)-|\la^1\!\vee\!\la^2|=1$. 
The self-duality and other claims hold, including the 
positivity upon the division by $(1-t)$. Interestingly,
(\ref{T4-2-11-20}) is the result of the substitution 
$q\mapsto q^2, t\mapsto t^2$ in the superpolynomial at
(\ref{T4-2}) for \Yboxdim7pt $\la^1=\yng(1)=\la^2$:
\Yboxdim5pt
$$
\hat{\h}^{min}_
{\yng(1)\,,\yng(1)\,}=
\hbox{\small
$1 - t + qt - q t^2 + q^2t^2 + a(q - q t + q^2 t)$}.
$$

For the same tree and \Yboxdim7pt
$\la^1=\yng(1,1)\,, \la^2=\yng(1)\,$:\Yboxdim5pt
\begin{align}\label{T4-2-11-10}
\hat{\h}{}^{min}_{\yng(1,1)\,,\yng(1)}\,(q,t,a)=
\hbox{\small
$
1 - t^2 + qt^2 - qt^4 + q^2t^4 + a(q - qt^2 + q^2t^2).
$
}
\end{align}
This is the result of the substitution $t\mapsto t^2$
in the previous formula. 

The superpolynomial from (\ref{T4-2-11-10}) is the
major factor of $-P_{[1,1],1}^{T[4,2]}$ in (237) from
\cite{DMMSS} upon the substitution
$q\mapsto q^2,t\mapsto t^2,a\mapsto A^2$. Colored superpolynomials
for $T(6,4)$ are not given in \cite{DMMSS}.

\medskip

{\sf 2-column and 1 box for trefoil.} Let us somewhat
reduce the previous example; we will consider the 
diagrams 
\Yboxdim7pt $\la^1\!=\yng(1,1)\,, \la^2\!=\yng(1)\,$\,\Yboxdim5pt
\,\, for the trefoil:

\begin{align}\label{T6-4-11-10}
&T(6,4):\ 
\l=\l_{(\{3,2\},\{3,2\})}^{\,\circ\rightrightarrows,\, 
(\yng(1,1)\,,\yng(1)\,)},\ \ 
\hat{\h}{}^{min}_{\l}\,(q,t,a)=
\end{align}

\renewcommand{\baselinestretch}{0.5} 
\noindent
{\small
\(
1+2 q t-t^2+2 q^2 t^2+q^3 t^2-2 q t^3+q^2 t^3+2 q^3 t^3-2 q^2 t^4
+3 q^4 t^4-q^2 t^5-q^3 t^5+2 q^4 t^5+q^5 t^5-q^3 t^6-2 q^4 t^6
+3 q^5 t^6-q^3 t^7-q^4 t^7+q^5 t^7+q^6 t^7-q^4 t^8-2 q^5 t^8
+3 q^6 t^8-q^4 t^9-q^5 t^9+q^6 t^9+q^7 t^9-q^5 t^{10}-q^6 t^{10}
+2 q^7 t^{10}-q^5 t^{11}-q^6 t^{11}+2 q^7 t^{11}-q^6 t^{12}
+q^8 t^{12}-q^6 t^{13}+q^8 t^{13}-q^7 t^{14}+q^8 t^{14}-q^7 t^{15}
+q^8 t^{15}-q^8 t^{17}+q^9 t^{17}+a^4 \bigl(\frac{q^7}{t}
-q^7 t+q^8 t-q^8 t^3+q^9 t^3\bigr)+a^3 \bigl(q^5+q^6+\frac{q^4}{t}
+\frac{q^5}{t}+\frac{q^6}{t}-q^4 t-q^5 t+q^6 t+2 q^7 t-q^5 t^2
+2 q^7 t^2-2 q^6 t^3+2 q^8 t^3-q^6 t^4-q^7 t^4+2 q^8 t^4-2 q^7 t^5
+q^8 t^5+q^9 t^5-q^7 t^6+q^9 t^6-q^8 t^7+q^9 t^7-q^8 t^8
+q^9 t^8\bigr)+a^2 \bigl(2 q^3+2 q^4+q^5+\frac{q^2}{t}
+\frac{q^3}{t}+\frac{q^4}{t}-q^2 t-q^3 t+2 q^4 t+3 q^5 t
+q^6 t-2 q^3 t^2-q^4 t^2+4 q^5 t^2+3 q^6 t^2-3 q^4 t^3-q^5 t^3
+4 q^6 t^3+q^7 t^3-q^4 t^4-4 q^5 t^4+2 q^6 t^4+4 q^7 t^4-2 q^5 t^5
-3 q^6 t^5+4 q^7 t^5+q^8 t^5-q^5 t^6-4 q^6 t^6+2 q^7 t^6+3 q^8 t^6
-2 q^6 t^7-2 q^7 t^7+4 q^8 t^7-q^6 t^8-3 q^7 t^8+3 q^8 t^8+q^9 t^8
-2 q^7 t^9+q^8 t^9+q^9 t^9-q^7 t^{10}-q^8 t^{10}+2 q^9 t^{10}
-q^8 t^{11}+q^9 t^{11}-q^8 t^{12}+q^9 t^{12}\bigr)
+a \bigl(q+2 q^2+q^3+\frac{q}{t}-q t+2 q^2 t+3 q^3 t+q^4 t-q t^2
-2 q^2 t^2+2 q^3 t^2+5 q^4 t^2+q^5 t^2-2 q^2 t^3-2 q^3 t^3
+3 q^4 t^3+3 q^5 t^3-3 q^3 t^4-3 q^4 t^4+6 q^5 t^4+2 q^6 t^4
-q^3 t^5-3 q^4 t^5+q^5 t^5+4 q^6 t^5-2 q^4 t^6-5 q^5 t^6+5 q^6 t^6
+2 q^7 t^6-q^4 t^7-3 q^5 t^7+4 q^7 t^7-2 q^5 t^8-4 q^6 t^8
+5 q^7 t^8+q^8 t^8-q^5 t^9-3 q^6 t^9+2 q^7 t^9+2 q^8 t^9
-2 q^6 t^{10}-2 q^7 t^{10}+4 q^8 t^{10}-q^6 t^{11}-2 q^7 t^{11}
+3 q^8 t^{11}-2 q^7 t^{12}+q^8 t^{12}+q^9 t^{12}-q^7 t^{13}
+q^9 t^{13}-q^8 t^{14}+q^9 t^{14}-q^8 t^{15}+q^9 t^{15}\bigr).
\)
}
\renewcommand{\baselinestretch}{1.0} 
\smallskip

The $a$\~degree of $\hat{\h}{}^{min}_{\l}\,(q,t,a)$ is
deg$_a=\ss(|\yng(2)\,|+|\yng(1)\,|)-|\yng(2)\,|=4.$
This superpolynomial is
dual to that for
\Yboxdim7pt $\la^1=\yng(2)\,, \la^2=\yng(1)\,,$
which we will omit. \Yboxdim5pt
The positivity of the series
$\hat{\h}{}^{min}_{\l}(q,t,a)/(1-t)$ holds.
\medskip

\subsection{\bf 3-hook for trefoil} Let us conclude this
set of examples with 
\Yboxdim7pt $\la^1\!=\yng(2,1)\,, \la^2\!=\yng(1)\,$\,\Yboxdim5pt
\, for the trefoil:

\begin{align}\label{T6-4-21-10}
&T(6,4):\ 
\l=\l_{(\{3,2\},\{3,2\})}^{\,\circ\rightrightarrows,\, 
(\yng(2,1)\,,\yng(1)\,)},\ \ 
\hat{\h}{}^{min}_{\l}\,(q,t,a)=
\end{align}

\renewcommand{\baselinestretch}{0.5} 
\noindent
{\small
\(
1-t+2 q t+q^2 t+q^3 t-3 q t^2+q^2 t^2+q^3 t^2+4 q^4 t^2+q^5 t^2
+q t^3-4 q^2 t^3-q^3 t^3-3 q^4 t^3+6 q^5 t^3+3 q^6 t^3+q^7 t^3
+2 q^2 t^4-3 q^3 t^4-q^4 t^4-9 q^5 t^4+6 q^6 t^4+5 q^7 t^4
+2 q^8 t^4
+2 q^3 t^5-2 q^4 t^5+2 q^5 t^5-13 q^6 t^5+2 q^7 t^5+6 q^8 t^5
+5 q^9 t^5+2 q^4 t^6-2 q^5 t^6+4 q^6 t^6-14 q^7 t^6+5 q^9 t^6
+6 q^{10} t^6+2 q^5 t^7-2 q^6 t^7+6 q^7 t^7-14 q^8 t^7-3 q^9 t^7
+4 q^{10} t^7+7 q^{11} t^7+2 q^6 t^8-2 q^7 t^8+6 q^8 t^8-14 q^9 t^8
-2 q^{10} t^8+4 q^{11} t^8+6 q^{12} t^8+2 q^7 t^9-2 q^8 t^9
+7 q^9 t^9-14 q^{10} t^9-3 q^{11} t^9+5 q^{12} t^9+5 q^{13} t^9
+2 q^8 t^{10}-2 q^9 t^{10}+6 q^{10} t^{10}-14 q^{11} t^{10}
+6 q^{13} t^{10}+2 q^{14} t^{10}+2 q^9 t^{11}-2 q^{10} t^{11}
+6 q^{11} t^{11}-14 q^{12} t^{11}+2 q^{13} t^{11}+5 q^{14} t^{11}
+q^{15} t^{11}+2 q^{10} t^{12}-2 q^{11} t^{12}+4 q^{12} t^{12}
-13 q^{13} t^{12}+6 q^{14} t^{12}+3 q^{15} t^{12}+2 q^{11} t^{13}
-2 q^{12} t^{13}+2 q^{13} t^{13}-9 q^{14} t^{13}+6 q^{15} t^{13}
+q^{16} t^{13}+2 q^{12} t^{14}-2 q^{13} t^{14}-q^{14} t^{14}
-3 q^{15} t^{14}+4 q^{16} t^{14}+2 q^{13} t^{15}-3 q^{14} t^{15}
-q^{15} t^{15}+q^{16} t^{15}+q^{17} t^{15}+2 q^{14} t^{16}
-4 q^{15} t^{16}+q^{16} t^{16}+q^{17} t^{16}+q^{15} t^{17}
-3 q^{16} t^{17}+2 q^{17} t^{17}-q^{17} t^{18}+q^{18} t^{18}
\)

\smallskip
\noindent
\(
+a^5 \bigl(-q^{15}+q^{16}+\frac{q^{15}}{t}-2 q^{16} t+2 q^{17} t
+q^{16} t^2-2 q^{17} t^2+q^{18} t^2-q^{18} t^3+q^{19} t^3\bigr)
+a^4 \bigl(-q^{10}+q^{12}+q^{13}+2 q^{14}+\frac{q^{10}}{t}
+\frac{q^{11}}{t}+\frac{q^{12}}{t}+\frac{q^{13}}{t}-2 q^{11} t
-2 q^{12} t+q^{14} t+4 q^{15} t+q^{11} t^2-q^{12} t^2-3 q^{13} t^2
+4 q^{16} t^2+q^{12} t^3-4 q^{14} t^3-q^{15} t^3+4 q^{17} t^3
+q^{13} t^4-4 q^{15} t^4+q^{17} t^4+2 q^{18} t^4+q^{14} t^5
-3 q^{16} t^5+q^{18} t^5+q^{19} t^5+q^{15} t^6-q^{16} t^6
-2 q^{17} t^6+q^{18} t^6+q^{19} t^6+q^{16} t^7-2 q^{17} t^7
+q^{19} t^7-q^{18} t^8+q^{19} t^8\bigr)
\)

\smallskip
\noindent
\(
+a^3 \bigl(-q^6+4 q^9
+3 q^{10}+3 q^{11}+q^{12}+\frac{q^6}{t}+\frac{q^7}{t}
+\frac{2 q^8}{t}+\frac{q^9}{t}+\frac{q^{10}}{t}-2 q^7 t-2 q^8 t
-5 q^9 t+4 q^{10} t+4 q^{11} t+6 q^{12} t+2 q^{13} t+q^7 t^2
-q^8 t^2-q^9 t^2-11 q^{10} t^2+q^{11} t^2+5 q^{12} t^2+8 q^{13} t^2
+4 q^{14} t^2+q^8 t^3+2 q^{10} t^3-14 q^{11} t^3-4 q^{12} t^3
+4 q^{13} t^3+9 q^{14} t^3+4 q^{15} t^3+q^9 t^4+5 q^{11} t^4
-15 q^{12} t^4-7 q^{13} t^4+4 q^{14} t^4+9 q^{15} t^4+4 q^{16} t^4
+q^{10} t^5+6 q^{12} t^5-14 q^{13} t^5-7 q^{14} t^5+4 q^{15} t^5
+8 q^{16} t^5+2 q^{17} t^5+q^{11} t^6+6 q^{13} t^6-15 q^{14} t^6
-4 q^{15} t^6+5 q^{16} t^6+6 q^{17} t^6+q^{18} t^6+q^{12} t^7
+5 q^{14} t^7-14 q^{15} t^7+q^{16} t^7+4 q^{17} t^7+3 q^{18} t^7
+q^{13} t^8+2 q^{15} t^8-11 q^{16} t^8+4 q^{17} t^8+3 q^{18} t^8
+q^{19} t^8+q^{14} t^9-q^{16} t^9-5 q^{17} t^9+4 q^{18} t^9
+q^{19} t^9+q^{15} t^{10}-q^{16} t^{10}-2 q^{17} t^{10}
+2 q^{19} t^{10}+q^{16} t^{11}-2 q^{17} t^{11}+q^{19} t^{11}
-q^{18} t^{12}+q^{19} t^{12}\bigr)
\)

\smallskip
\noindent
\(
+a^2 \bigl(-q^3+2 q^5+3 q^6
+5 q^7+2 q^8+q^9+\frac{q^3}{t}+\frac{q^4}{t}+\frac{q^5}{t}
+\frac{q^6}{t}-2 q^4 t-3 q^5 t-q^6 t+2 q^7 t+10 q^8 t+5 q^9 t
+3 q^{10} t+q^4 t^2-q^5 t^2-5 q^6 t^2-7 q^7 t^2-6 q^8 t^2
+13 q^9 t^2+10 q^{10} t^2+7 q^{11} t^2+q^{12} t^2+q^5 t^3+q^6 t^3
-3 q^7 t^3-10 q^8 t^3-17 q^9 t^3+11 q^{10} t^3+13 q^{11} t^3
+10 q^{12} t^3+q^{13} t^3+q^6 t^4+2 q^7 t^4+q^8 t^4-8 q^9 t^4
-27 q^{10} t^4+6 q^{11} t^4+15 q^{12} t^4+12 q^{13} t^4
+2 q^{14} t^4+q^7 t^5+2 q^8 t^5+3 q^9 t^5-4 q^{10} t^5
-31 q^{11} t^5+2 q^{12} t^5
+16 q^{13} t^5+12 q^{14} t^5+q^{15} t^5+q^8 t^6+2 q^9 t^6
+4 q^{10} t^6-2 q^{11} t^6-33 q^{12} t^6+2 q^{13} t^6+15 q^{14} t^6
+10 q^{15} t^6+q^{16} t^6+q^9 t^7+2 q^{10} t^7+4 q^{11} t^7
-2 q^{12} t^7-31 q^{13} t^7+6 q^{14} t^7+13 q^{15} t^7+7 q^{16} t^7
+q^{10} t^8+2 q^{11} t^8+4 q^{12} t^8-4 q^{13} t^8-27 q^{14} t^8
+11 q^{15} t^8+10 q^{16} t^8+3 q^{17} t^8+q^{11} t^9+2 q^{12} t^9
+3 q^{13} t^9-8 q^{14} t^9-17 q^{15} t^9+13 q^{16} t^9+5 q^{17} t^9
+q^{18} t^9+q^{12} t^{10}+2 q^{13} t^{10}+q^{14} t^{10}
-10 q^{15} t^{10}-6 q^{16} t^{10}+10 q^{17} t^{10}+2 q^{18} t^{10}
+q^{13} t^{11}+2 q^{14} t^{11}-3 q^{15} t^{11}-7 q^{16} t^{11}
+2 q^{17} t^{11}+5 q^{18} t^{11}+q^{14} t^{12}+q^{15} t^{12}
-5 q^{16} t^{12}-q^{17} t^{12}+3 q^{18} t^{12}+q^{19} t^{12}
+q^{15} t^{13}-q^{16} t^{13}-3 q^{17} t^{13}+2 q^{18} t^{13}
+q^{19} t^{13}+q^{16} t^{14}-2 q^{17} t^{14}+q^{19} t^{14}
-q^{18} t^{15}+q^{19} t^{15}\bigr)
\)

\smallskip
\noindent
\(
+a \bigl(-q+2 q^2+2 q^3
+2 q^4+q^5+\frac{q}{t}-3 q^2 t+q^3 t+2 q^4 t+6 q^5 t+4 q^6 t
+q^7 t+q^2 t^2-5 q^3 t^2-3 q^4 t^2-4 q^5 t^2+8 q^6 t^2+10 q^7 t^2
+4 q^8 t^2+q^9 t^2+2 q^3 t^3-4 q^4 t^3-4 q^5 t^3-14 q^6 t^3
+3 q^7 t^3+14 q^8 t^3+9 q^9 t^3+2 q^{10} t^3+3 q^4 t^4-2 q^5 t^4
-21 q^7 t^4-6 q^8 t^4+14 q^9 t^4+13 q^{10} t^4+4 q^{11} t^4
+3 q^5 t^5-q^6 t^5+5 q^7 t^5-22 q^8 t^5-15 q^9 t^5+12 q^{10} t^5
+17 q^{11} t^5+5 q^{12} t^5+3 q^6 t^6-q^7 t^6+8 q^8 t^6-21 q^9 t^6
-19 q^{10} t^6+9 q^{11} t^6+17 q^{12} t^6+5 q^{13} t^6+3 q^7 t^7
-q^8 t^7+10 q^9 t^7-20 q^{10} t^7-22 q^{11} t^7+9 q^{12} t^7
+17 q^{13} t^7+4 q^{14} t^7+3 q^8 t^8-q^9 t^8+10 q^{10} t^8
-20 q^{11} t^8-19 q^{12} t^8+12 q^{13} t^8+13 q^{14} t^8
+2 q^{15} t^8+3 q^9 t^9-q^{10} t^9+10 q^{11} t^9-21 q^{12} t^9
-15 q^{13} t^9+14 q^{14} t^9+9 q^{15} t^9+q^{16} t^9
+3 q^{10} t^{10}-q^{11} t^{10}+8 q^{12} t^{10}-22 q^{13} t^{10}
-6 q^{14} t^{10}+14 q^{15} t^{10}+4 q^{16} t^{10}+3 q^{11} t^{11}
-q^{12} t^{11}+5 q^{13} t^{11}-21 q^{14} t^{11}+3 q^{15} t^{11}
+10 q^{16} t^{11}+q^{17} t^{11}+3 q^{12} t^{12}-q^{13} t^{12}
-14 q^{15} t^{12}+8 q^{16} t^{12}+4 q^{17} t^{12}+3 q^{13} t^{13}
-2 q^{14} t^{13}-4 q^{15} t^{13}-4 q^{16} t^{13}+6 q^{17} t^{13}
+q^{18} t^{13}+3 q^{14} t^{14}-4 q^{15} t^{14}-3 q^{16} t^{14}
+2 q^{17} t^{14}+2 q^{18} t^{14}+2 q^{15} t^{15}-5 q^{16} t^{15}
+q^{17} t^{15}+2 q^{18} t^{15}+q^{16} t^{16}
-3 q^{17} t^{16}+2 q^{18} t^{16}-q^{18} t^{17}+q^{19} t^{17}\bigr).
\)
}
\renewcommand{\baselinestretch}{1.0} 

\smallskip

The $a$\~degree of $\hat{\h}{}^{min}_{\l}\,(q,t,a)$ is\, 
$\ss(|\yng(2,1)\,|+|\yng(1)\,|)-|\yng(2,1)\,|=5.$
This superpolynomial is self-dual.
The positivity of the series
$\hat{\h}{}^{min}_{\l}(q,t,a)/(1-t)^2$ holds.
Here $(1-t)^2$ is minimal necessary, since the 3-hook
alone requires the division by $(1-t)$ for the positivity
and we have another $(1-t)$ because we do a $2$\~link.
This is of course an informal argument, not a justification.

Generally, Part $(ii)$ of the Connection Conjecture claims 
that proper powers of $(1-t)(1-q)$ are necessary for 
arbitrary colors for {\em positive pairs\,} of trees
$\{\l,\,'\!\l^\vee\}$, i.e. those satisfying the
inequalities $\rr,\,'\!\rr, \ss,\,'\!\ss>0$ and
$\,'\ss_1\ss_1>\,'\rr_1\rr_1$. 
This is supposed to be similar to Corollary 4.10 from 
\cite{EGL} concerning the positivity
of {\em unreduced\,} HOMFLY-PT polynomials. See also
Example 4.22 there concerning the $3$\~hook and
Section 5.4 devoted to the $3$\~hook in \cite{CJ}.
However the powers $(1-t)^\bullet$ are insufficient for 
algebraic links, which cannot be seen in the
theory of HOMFLY-PT polynomials (where $t=q$). 
We do not
have any conjectures at the moment about the minimal powers of
$(1-t)$ and $(1-q)$ necessary for the positivity (assuming
that Part $(ii)$ is true).
\smallskip

\subsection{\bf Alexander polynomials}\label{sec:Alex}
The uncolored reduced hat-normalized
HOMFLY-PT polynomial is  
$\hat{\h}{}^{min}_{\l}\,(q,q,-a)/
(1-q)^{\kappa-1}$, recalculated to the standard $a,q$. 
Therefore the corresponding
{\em Alex\-an\-der polynomial\,} must be 
$\hat{\h}{}^{min}_{\l}\,(q,q,-1)/(1-q)^{\kappa-\de_{\kappa,1}}$,
where $\de$ is the Kronecker symbol.
Without $\de_{\kappa,1}$, this
will be the multivariable Alexander polynomial for coinciding 
parameters. 
For algebraic links, this polynomial is directly connected 
with the zeta function of monodromy of the Milnor fiber
of the corresponding singularity. Here
$\kappa$ is the number of connected component of the link, 
which is $2$ in the example above.
Generally, the substitution
$a\mapsto -1$ is the passage to the 
{\em Heegaard-Floer homology}. We will provide in this section
the equations of the corresponding singularities  
(at $x=0,y=0$) and their Alexander polynomials, the
{\em zeta-monodromy\,} from \cite{DGPS} upon $t\mapsto q$
(unless for the unknot). 
\smallskip

Note that the value
of the Alexander polynomial at $q=1$ for any $2$\~link 
(a link with $2$ components) is $\pm\,lk$ 
for the linking number $\,lk$. Thus this polynomial 
uniquely determines the (topological) 
type of singularity in this 
case, assuming that the corresponding knot components are known. 
Therefore $\hat{\h}$ must uniquely determine such singularities. 

\smallskip

Let us consider the uncolored $C\!ab(7,1)T(6,4)$ for 
$\Up=\{\circ\rightarrow\circ\rightrightarrows\}$
with $2$ arrowheads.
We set  $\vec\rr=\{3,1\},\vec\ss=\{2,1\}$.
The \tax-labels of the vertices are
$[3,2]$ and $[1,1]$ (the first and the second).  
As above, we put $T(6,4)$ instead of $(C\!ab(2,3) C\!ab(2,3))
(\unknot\,)$.

\begin{align}\label{T7-1-6-4}
&C\!ab(7,1)T(6,4):\ 
\l=\l_{\{3,1\},\{2,1\}}^
{\,\circ\rightarrow\circ\rightrightarrows,\, 
(\yng(1)\,,\yng(1)\,)},\ \ 
\hat{\h}{}^{min}_{\l}\,(q,t,a)=
\end{align}

\renewcommand{\baselinestretch}{0.5} 
\noindent
{\small
\(
1-t+q t+q^2 t+q^3 t-q t^2+2 q^4 t^2-q^2 t^3-q^4 t^3+2 q^5 t^3
-q^3 t^4-q^5 t^4+2 q^6 t^4-q^4 t^5-q^6 t^5+2 q^7 t^5-q^5 t^6
+q^8 t^6-q^6 t^7+q^8 t^7-q^7 t^8+q^8 t^8-q^8 t^9+q^9 t^9+a^3 
\bigl(q^6-q^6 t+q^7 t-q^7 t^2+q^8 t^2-q^8 t^3+q^9 t^3\bigr)
+a^2 \bigl(q^3+q^4+q^5-q^3 t+q^5 t+2 q^6 t-q^4 t^2-q^5 t^2
+2 q^7 t^2-q^5 t^3-q^6 t^3+2 q^8 t^3-q^6 t^4-q^7 t^4+q^8 t^4
+q^9 t^4-q^7 t^5+q^9 t^5-q^8 t^6+q^9 t^6\bigr)+a \bigl(q+q^2
+q^3-q t+q^3 t+3 q^4 t+q^5 t-q^2 t^2-q^3 t^2-q^4 t^2+3 q^5 t^2
+q^6 t^2-q^3 t^3-q^4 t^3-2 q^5 t^3+3 q^6 t^3+q^7 t^3-q^4 t^4
-q^5 t^4-2 q^6 t^4+3 q^7 t^4+q^8 t^4-q^5 t^5-q^6 t^5-q^7 t^5
+3 q^8 t^5-q^6 t^6
-q^7 t^6+q^8 t^6+q^9 t^6-q^7 t^7+q^9 t^7-q^8 t^8+q^9 t^8\bigr).
\)
}
\renewcommand{\baselinestretch}{1.0} 

\smallskip

The $a$\~degree of $\hat{\h}{}^{min}_{\l}\,(q,t,a)$ is
deg$_a=\ss_1\rr_2|\yng(1)\,|+
\ss_1|\yng(1)\,|-|\yng(1)\,|= 2\times 1+2-1=3$.
The corresponding singularity is:
$$ \c: \{(x^3\!+\!y^2)^2\!+\!x^7\!=\!0\},\, 
\hbox{Link}(\c)=\bigl\{ T(3,2),T(3,2),
\hbox{lk}\!=\!7\bigr\},
$$
The Alexander polynomial is 
$1+q^4+q^6+q^8+q^{10}+q^{12}+q^{16}$; use 
\cite{DGPS}. This example is of interest
because the link of $ (x^3+y^2)^2+x^7+x^5y=0$ is
$C\!ab(13,2)T(3,2)$, which is the simplest 
non-torus algebraic {\em knot\,}; its Alexander polynomial 
is
$1-q+q^4-q^5+q^6-q^7+q^8-q^9+q^{10}-q^{11}+q^{12}-q^{15}+q^{16}$. 

We note that the term  $x^5y$ here is  minimal in the following 
sense. If it is replaced by  $x^4y$, then the 
corresponding
$\hbox{Link}((x^3\!+\!y^2)^2\!+\!x^7\!+\!x^4y\!=\!0)$ becomes 
$\bigl\{ T(4,3),T(1,0),
\hbox{lk}\!=\!4\bigr\}$, i.e. not a knot.
\smallskip
 
Now let us discuss the simplest algebraic link obtained
by linking $T(3,2)$ and $T(5,3)$.
It will be $C\!ab(1,1)(T(3,2)T(5,3))$ for
$\Up=\{\twotwo\}$ with 
$\vec\rr^1=\{1,2\},\vec\ss^1=\{1,1\}$ and 
 $\vec\rr^2=\{1,3\},\vec\ss^2=\{1,2\}$.
The \tax-labels of the vertices are
$[1,1]$ and $\{[2,1],[3,2]\} $ (for the first $\circ$ and 
for the remaining two vertices).  

\begin{align}\label{T3-2-5-3-1-1}
&\bigl(C\!ab(3,2)C\!ab(5,3)\bigr)T(1,1):\ 
\l_{\bigl\{\{1,2\},\{1,1\}\bigr\},
\bigl\{\{1,3\},\{1,2\}\bigr\}}
^{\,\twotwotiny,\, 
(\yng(1)\,,\yng(1)\,)},\,
\hat{\h}{}^{min}_{\l}=
\end{align}

\renewcommand{\baselinestretch}{0.5} 
\noindent
{\small
\(
1-t+q t+q^2 t+q^3 t+q^4 t-q t^2+2 q^4 t^2+2 q^5 t^2+q^6 t^2-q^2 t^3
-2 q^4 t^3+q^5 t^3+3 q^6 t^3+2 q^7 t^3-q^3 t^4-2 q^5 t^4+3 q^7 t^4
+2 q^8 t^4-q^4 t^5-2 q^6 t^5+3 q^8 t^5+q^9 t^5-q^5 t^6-2 q^7 t^6
+q^8 t^6+2 q^9 t^6-q^6 t^7-2 q^8 t^7+2 q^9 t^7+q^{10} t^7-q^7 t^8
+q^{10} t^8-q^8 t^9+q^{10} t^9-q^9 t^{10}+q^{10} t^{10}
-q^{10} t^{11}+q^{11} t^{11}+a^4 \bigl(q^{10}+q^{11} t\bigr)
+a^3 \bigl(q^6+q^7+q^8+q^9-q^6 t+q^7 t+2 q^8 t+2 q^9 t+2 q^{10} t
-2 q^7 t^2+2 q^9 t^2+2 q^{10} t^2+q^{11} t^2-2 q^8 t^3+2 q^{10} t^3
+q^{11} t^3-2 q^9 t^4+q^{10} t^4+q^{11} t^4-q^{10} t^5
+q^{11} t^5\bigr)+a^2 \bigl(q^3+q^4+2 q^5+q^6+q^7-q^3 t+q^5 t
+4 q^6 t+4 q^7 t+3 q^8 t+q^9 t-q^4 t^2-2 q^5 t^2-2 q^6 t^2
+4 q^7 t^2+6 q^8 t^2+4 q^9 t^2+q^{10} t^2-q^5 t^3-2 q^6 t^3
-4 q^7 t^3+3 q^8 t^3+6 q^9 t^3+3 q^{10} t^3-q^6 t^4-2 q^7 t^4
-4 q^8 t^4
+4 q^9 t^4+4 q^{10} t^4+q^{11} t^4-q^7 t^5-2 q^8 t^5-2 q^9 t^5
+4 q^{10} t^5+q^{11} t^5-q^8 t^6-2 q^9 t^6+q^{10} t^6+2 q^{11} t^6
-q^9 t^7+q^{11} t^7-q^{10} t^8+q^{11} t^8\bigr)+a \bigl(q+q^2+q^3
+q^4-q t+q^3 t+3 q^4 t+4 q^5 t+2 q^6 t+q^7 t-q^2 t^2-q^3 t^2
-2 q^4 t^2+2 q^5 t^2+6 q^6 t^2+5 q^7 t^2+2 q^8 t^2-q^3 t^3-q^4 t^3
-4 q^5 t^3-q^6 t^3+6 q^7 t^3+6 q^8 t^3+2 q^9 t^3-q^4 t^4-q^5 t^4
-4 q^6 t^4-2 q^7 t^4+6 q^8 t^4+5 q^9 t^4+q^{10} t^4-q^5 t^5-q^6 t^5
-4 q^7 t^5-q^8 t^5+6 q^9 t^5+2 q^{10} t^5-q^6 t^6-q^7 t^6-4 q^8 t^6
+2 q^9 t^6+4 q^{10} t^6-q^7 t^7-q^8 t^7-2 q^9 t^7+3 q^{10} t^7
+q^{11} t^7-q^8 t^8-q^9 t^8+q^{10} t^8+q^{11} t^8-q^9 t^9
+q^{11} t^9
-q^{10} t^{10}+q^{11} t^{10}\bigr).
\)
}
\renewcommand{\baselinestretch}{1.0} 
\smallskip

The $a$\~degree is
$\ss^1_1\rr^1_2+
\ss^2_1\rr^2_2-1= 3+2-1=4$.
The positivity of 
$\hat{\h}{}^{min}_{\l}(q,t,a)/(1-t)$ holds.
The corresponding
singularity is given by
the equation $(y^5+x^3)(y^2+x^3)=0$ with
the linking number $lk=6$; its
Alexander polynomial is
$1+q^6+q^7+q^{13}+q^{14}+q^{20}.$
\smallskip
 
Our last example will be for  the uncolored tree
$\Up=\{\twoone\}$
with the labeled vertices
$[\rr^1_1=2,\ss^1_1=1],[\rr^1_2=3,\ss^1_2=2]$
for the long path and $[\rr^2_1=2,\ss^2_1=1]$ for
the short one (equal to $[\rr_1^1,\ss^1_1]$ in the tree).
\begin{align}\label{T8-3-2-1}
&\bigl(C\!ab(8,3)C\!ab(0,1)\bigr)T(2,1):\ 
\l=\l_{\{2,3\},\{1,2\}}^
{\,\twoonetiny,\, 
(\yng(1)\,,\yng(1)\,)},\ \ 
\hat{\h}{}^{min}_{\l}\,(q,t,a)=
\end{align}

\renewcommand{\baselinestretch}{0.5} 
\noindent
{\small
\(
1-t+q t+q^2 t+q^3 t-q t^2+2 q^4 t^2+q^5 t^2+q^6 t^2-q^2 t^3
-q^4 t^3+q^5 t^3+2 q^6 t^3+2 q^7 t^3-q^3 t^4-q^5 t^4+q^7 t^4
+3 q^8 t^4-q^4 t^5-q^6 t^5+q^8 t^5+3 q^9 t^5-q^5 t^6-q^7 t^6
+q^9 t^6+2 q^{10} t^6-q^6 t^7-q^8 t^7+2 q^{10} t^7+q^{11} t^7
-q^7 t^8-q^9 t^8+q^{10} t^8+q^{11} t^8-q^8 t^9-q^{10} t^9
+2 q^{11} t^9-q^9 t^{10}+q^{12} t^{10}-q^{10} t^{11}+q^{12} t^{11}
-q^{11} t^{12}+q^{12} t^{12}-q^{12} t^{13}+q^{13} t^{13}
+a^3 \bigl(q^6-q^6 t+q^7 t+q^8 t+q^9 t-q^7 t^2+q^9 t^2+2 q^{10} t^2
-q^8 t^3+2 q^{11} t^3-q^9 t^4+q^{11} t^4+q^{12} t^4-q^{10} t^5
+q^{12} t^5-q^{11} t^6+q^{12} t^6-q^{12} t^7+q^{13} t^7\bigr)
+a^2 \bigl(q^3+q^4+q^5-q^3 t+q^5 t+3 q^6 t+2 q^7 t+q^8 t-q^4 t^2
-q^5 t^2-q^6 t^2+3 q^7 t^2+4 q^8 t^2+3 q^9 t^2-q^5 t^3-q^6 t^3
-2 q^7 t^3+q^8 t^3+4 q^9 t^3+4 q^{10} t^3-q^6 t^4-q^7 t^4-2 q^8 t^4
+q^9 t^4+4 q^{10} t^4+3 q^{11} t^4-q^7 t^5-q^8 t^5-2 q^9 t^5
+q^{10} t^5+4 q^{11} t^5+q^{12} t^5-q^8 t^6-q^9 t^6-2 q^{10} t^6
+3 q^{11} t^6+2 q^{12} t^6-q^9 t^7-q^{10} t^7-q^{11} t^7
+3 q^{12} t^7-q^{10} t^8-q^{11} t^8+q^{12} t^8+q^{13} t^8
-q^{11} t^9+q^{13} t^9-q^{12} t^{10}+q^{13} t^{10}\bigr)
+a \bigl(q+q^2+q^3-q t+q^3 t+3 q^4 t+2 q^5 t+q^6 t-q^2 t^2
-q^3 t^2-q^4 t^2+2 q^5 t^2+4 q^6 t^2+4 q^7 t^2+q^8 t^2-q^3 t^3
-q^4 t^3-2 q^5 t^3+3 q^7 t^3+6 q^8 t^3+2 q^9 t^3-q^4 t^4-q^5 t^4
-2 q^6 t^4-q^7 t^4+2 q^8 t^4+6 q^9 t^4+2 q^{10} t^4-q^5 t^5
-q^6 t^5-2 q^7 t^5-q^8 t^5+2 q^9 t^5+6 q^{10} t^5+q^{11} t^5
-q^6 t^6-q^7 t^6-2 q^8 t^6-q^9 t^6+3 q^{10} t^6+4 q^{11} t^6
-q^7 t^7-q^8 t^7-2 q^9 t^7+4 q^{11} t^7+q^{12} t^7-q^8 t^8
-q^9 t^8-2 q^{10} t^8+2 q^{11} t^8+2 q^{12} t^8-q^9 t^9
-q^{10} t^9-q^{11} t^9+3 q^{12} t^9-q^{10} t^{10}-q^{11} t^{10}
+q^{12} t^{10}+q^{13} t^{10}
-q^{11} t^{11}+q^{13} t^{11}-q^{12} t^{12}+q^{13} t^{12}\bigr).
\)
}
\renewcommand{\baselinestretch}{1.0} 
\medskip

The $a$\~degree of $\hat{\h}{}^{min}_{\l}\,(q,t,a)$ is
deg$_a=\ss_1\rr_2|\yng(1)\,|+
\ss_1|\yng(1)\,|-|\yng(1)\,|= 3+1-1=3$.
The positivity of 
$\hat{\h}{}^{min}_{\l}(q,t,a)/(1-t)$ holds. The corresponding
singularity is $(y^8+x^3)(y^2+x)=0$ with $lk=6$ and the
Alexander polynomial $1+q^4+q^{10}+q^{14}+q^{20}+q^{24}$.
Recall that the latter is {\em zeta-monodromy\,} from
\cite{DGPS} as $t\mapsto q$. 
\smallskip

\subsection{\bf Two different paths}
Let us extend the previous example. 
The tree will be now uncolored
$\Up=\twotwo$ with the vertices labeled by 
$$
[\rr^1_1=2,\ss^1_1=1], [\rr^1_2=3,\ss^1_2=2], \hbox{\, and\, }
[\rr^2_1=2,\ss^2_1=1],[\rr^2_2=2,\ss^2_2=1].
$$
The first vertices $[\rr_1^i,\ss_1^i]$ from the
corresponding paths are identified in this tree.
Since we have two different paths, the simplified
notations for $\l$ from the previous example cannot
be used now; both sets of labels (for both paths)
must be shown.

\begin{align}\label{T5-2-8-3-2-1}
&\bigl(C\!ab(8,3)C\!ab(5,2)\bigr)T(2,1):\ 
\l_{\bigl\{\{2,3\},\{1,2\}\bigr\},
\bigl\{\{2,2\},\{1,1\}\bigr\}}
^{\,\twotwotiny,\, 
(\yng(1)\,,\yng(1)\,)},\,
\hat{\h}{}^{min}_{\l}=
\end{align}

{\vbadness=10000 \hbadness=10000
\renewcommand{\baselinestretch}{0.5} 
\noindent
{\small
\(
1-t+q t+q^2 t+q^3 t+q^4 t-q t^2+q^4 t^2+2 q^5 t^2+2 q^6 t^2+q^7 t^2
+q^8 t^2-q^2 t^3-q^4 t^3+q^6 t^3+2 q^7 t^3+3 q^8 t^3+3 q^9 t^3
+q^{10} t^3-q^3 t^4-q^5 t^4-q^6 t^4+2 q^9 t^4+5 q^{10} t^4
+3 q^{11} t^4+q^{12} t^4-q^4 t^5-q^6 t^5-q^7 t^5-q^8 t^5-q^9 t^5
+5 q^{11} t^5+5 q^{12} t^5+2 q^{13} t^5-q^5 t^6-q^7 t^6-q^8 t^6
-q^9 t^6-2 q^{10} t^6-2 q^{11} t^6+4 q^{12} t^6+6 q^{13} t^6
+3 q^{14} t^6-q^6 t^7-q^8 t^7-q^9 t^7-q^{10} t^7-2 q^{11} t^7
-3 q^{12} t^7+3 q^{13} t^7+7 q^{14} t^7+3 q^{15} t^7-q^7 t^8
-q^9 t^8
-q^{10} t^8-q^{11} t^8-2 q^{12} t^8-3 q^{13} t^8+3 q^{14} t^8
+6 q^{15} t^8+2 q^{16} t^8-q^8 t^9-q^{10} t^9-q^{11} t^9
-q^{12} t^9
-2 q^{13} t^9-3 q^{14} t^9+4 q^{15} t^9+5 q^{16} t^9+q^{17} t^9
-q^9 t^{10}-q^{11} t^{10}-q^{12} t^{10}-q^{13} t^{10}
-2 q^{14} t^{10}-2 q^{15} t^{10}+5 q^{16} t^{10}+3 q^{17} t^{10}
-q^{10} t^{11}-q^{12} t^{11}-q^{13} t^{11}-q^{14} t^{11}
-2 q^{15} t^{11}+5 q^{17} t^{11}+q^{18} t^{11}-q^{11} t^{12}
-q^{13} t^{12}-q^{14} t^{12}-q^{15} t^{12}-q^{16} t^{12}
+2 q^{17} t^{12}+3 q^{18} t^{12}-q^{12} t^{13}-q^{14} t^{13}
-q^{15} t^{13}-q^{16} t^{13}+3 q^{18} t^{13}+q^{19} t^{13}
-q^{13} t^{14}-q^{15} t^{14}-q^{16} t^{14}+2 q^{18} t^{14}
+q^{19} t^{14}-q^{14} t^{15}-q^{16} t^{15}-q^{17} t^{15}
+q^{18} t^{15}+2 q^{19} t^{15}-q^{15} t^{16}-q^{17} t^{16}
+2 q^{19} t^{16}-q^{16} t^{17}-q^{18} t^{17}+q^{19} t^{17}
+q^{20} t^{17}-q^{17} t^{18}+q^{20} t^{18}-q^{18} t^{19}
+q^{20} t^{19}-q^{19} t^{20}+q^{20} t^{20}-q^{20} t^{21}
+q^{21} t^{21}
\)

\smallskip
\noindent
\(+a^4 \bigl(q^{10}-q^{10} t+q^{11} t+q^{12} t
+q^{13} t+q^{14} t-q^{11} t^2+2 q^{14} t^2+2 q^{15} t^2+q^{16} t^2
-q^{12} t^3-2 q^{14} t^3+q^{15} t^3+3 q^{16} t^3+2 q^{17} t^3
-q^{13} t^4-2 q^{15} t^4+3 q^{17} t^4+2 q^{18} t^4-q^{14} t^5
-2 q^{16} t^5+3 q^{18} t^5+q^{19} t^5-q^{15} t^6-2 q^{17} t^6
+q^{18} t^6+2 q^{19} t^6-q^{16} t^7-2 q^{18} t^7+2 q^{19} t^7
+q^{20} t^7-q^{17} t^8+q^{20} t^8-q^{18} t^9+q^{20} t^9
-q^{19} t^{10}+q^{20} t^{10}-q^{20} t^{11}+q^{21} t^{11}\bigr)
\)

\smallskip
\noindent
\(
+a^3 \bigl(q^6+q^7+q^8+q^9-q^6 t+q^8 t+2 q^9 t+4 q^{10} t
+3 q^{11} t+2 q^{12} t+q^{13} t-q^7 t^2-q^8 t^2-q^9 t^2
+4 q^{11} t^2+6 q^{12} t^2+6 q^{13} t^2+4 q^{14} t^2+q^{15} t^2
-q^8 t^3-q^9 t^3-2 q^{10} t^3-3 q^{11} t^3+4 q^{13} t^3
+9 q^{14} t^3+8 q^{15} t^3+3 q^{16} t^3-q^9 t^4-q^{10} t^4
-2 q^{11} t^4-4 q^{12} t^4
-3 q^{13} t^4+q^{14} t^4+9 q^{15} t^4+10 q^{16} t^4+4 q^{17} t^4
-q^{10} t^5-q^{11} t^5-2 q^{12} t^5-4 q^{13} t^5-5 q^{14} t^5
-q^{15} t^5+9 q^{16} t^5+10 q^{17} t^5+3 q^{18} t^5-q^{11} t^6
-q^{12} t^6-2 q^{13} t^6-4 q^{14} t^6-5 q^{15} t^6-q^{16} t^6
+9 q^{17} t^6+8 q^{18} t^6+q^{19} t^6-q^{12} t^7-q^{13} t^7
-2 q^{14} t^7-4 q^{15} t^7-5 q^{16} t^7+q^{17} t^7+9 q^{18} t^7
+4 q^{19} t^7-q^{13} t^8-q^{14} t^8-2 q^{15} t^8-4 q^{16} t^8
-3 q^{17} t^8+4 q^{18} t^8+6 q^{19} t^8+q^{20} t^8-q^{14} t^9
-q^{15} t^9-2 q^{16} t^9-4 q^{17} t^9+6 q^{19} t^9+2 q^{20} t^9
-q^{15} t^{10}-q^{16} t^{10}-2 q^{17} t^{10}-3 q^{18} t^{10}
+4 q^{19} t^{10}+3 q^{20} t^{10}-q^{16} t^{11}-q^{17} t^{11}
-2 q^{18} t^{11}+4 q^{20} t^{11}-q^{17} t^{12}-q^{18} t^{12}
-q^{19} t^{12}+2 q^{20} t^{12}+q^{21} t^{12}-q^{18} t^{13}
-q^{19} t^{13}+q^{20} t^{13}+q^{21} t^{13}-q^{19} t^{14}
+q^{21} t^{14}-q^{20} t^{15}+q^{21} t^{15}\bigr)
\)

\smallskip
\noindent
\(
+a^2 \bigl(q^3+q^4+2 q^5+q^6+q^7-q^3 t+3 q^6 t+4 q^7 t+5 q^8 t
+4 q^9 t+2 q^{10} t+q^{11} t-q^4 t^2-q^5 t^2-2 q^6 t^2+2 q^8 t^2
+7 q^9 t^2+9 q^{10} t^2+9 q^{11} t^2+5 q^{12} t^2+2 q^{13} t^2
-q^5 t^3-q^6 t^3-3 q^7 t^3-2 q^8 t^3-3 q^9 t^3+2 q^{10} t^3
+8 q^{11} t^3+14 q^{12} t^3+12 q^{13} t^3+6 q^{14} t^3+q^{15} t^3
-q^6 t^4-q^7 t^4-3 q^8 t^4-3 q^9 t^4-5 q^{10} t^4-4 q^{11} t^4
+2 q^{12} t^4+14 q^{13} t^4+17 q^{14} t^4+10 q^{15} t^4
+2 q^{16} t^4-q^7 t^5-q^8 t^5-3 q^9 t^5-3 q^{10} t^5-6 q^{11} t^5
-7 q^{12} t^5-4 q^{13} t^5+12 q^{14} t^5+20 q^{15} t^5
+11 q^{16} t^5+2 q^{17} t^5-q^8 t^6-q^9 t^6-3 q^{10} t^6
-3 q^{11} t^6-6 q^{12} t^6-8 q^{13} t^6-7 q^{14} t^6+11 q^{15} t^6
+20 q^{16} t^6+10 q^{17} t^6+q^{18} t^6-q^9 t^7-q^{10} t^7
-3 q^{11} t^7-3 q^{12} t^7-6 q^{13} t^7-8 q^{14} t^7-7 q^{15} t^7
+12 q^{16} t^7+17 q^{17} t^7+6 q^{18} t^7
-q^{10} t^8-q^{11} t^8-3 q^{12} t^8-3 q^{13} t^8-6 q^{14} t^8
-8 q^{15} t^8-4 q^{16} t^8+14 q^{17} t^8+12 q^{18} t^8+2 q^{19} t^8
-q^{11} t^9-q^{12} t^9-3 q^{13} t^9-3 q^{14} t^9-6 q^{15} t^9
-7 q^{16} t^9+2 q^{17} t^9+14 q^{18} t^9+5 q^{19} t^9-q^{12} t^{10}
-q^{13} t^{10}-3 q^{14} t^{10}-3 q^{15} t^{10}-6 q^{16} t^{10}
-4 q^{17} t^{10}+8 q^{18} t^{10}+9 q^{19} t^{10}+q^{20} t^{10}
-q^{13} t^{11}-q^{14} t^{11}-3 q^{15} t^{11}-3 q^{16} t^{11}
-5 q^{17} t^{11}+2 q^{18} t^{11}+9 q^{19} t^{11}+2 q^{20} t^{11}
-q^{14} t^{12}-q^{15} t^{12}-3 q^{16} t^{12}-3 q^{17} t^{12}
-3 q^{18} t^{12}+7 q^{19} t^{12}+4 q^{20} t^{12}-q^{15} t^{13}
-q^{16} t^{13}-3 q^{17} t^{13}-2 q^{18} t^{13}+2 q^{19} t^{13}
+5 q^{20} t^{13}-q^{16} t^{14}-q^{17} t^{14}-3 q^{18} t^{14}
+4 q^{20} t^{14}+q^{21} t^{14}-q^{17} t^{15}-q^{18} t^{15}
-2 q^{19} t^{15}+3 q^{20} t^{15}+q^{21} t^{15}-q^{18} t^{16}
-q^{19} t^{16}+2 q^{21} t^{16}-q^{19} t^{17}+q^{21} t^{17}
-q^{20} t^{18}+q^{21} t^{18}\bigr)
\)

\smallskip
\noindent
\(
+a \bigl(q+q^2+q^3+q^4
-q t+q^3 t+2 q^4 t+4 q^5 t+3 q^6 t+2 q^7 t+q^8 t-q^2 t^2-q^3 t^2
-q^4 t^2+3 q^6 t^2+5 q^7 t^2+7 q^8 t^2+6 q^9 t^2+3 q^{10} t^2
+q^{11} t^2-q^3 t^3-q^4 t^3-2 q^5 t^3-2 q^6 t^3+q^8 t^3+7 q^9 t^3
+11 q^{10} t^3+9 q^{11} t^3+4 q^{12} t^3+q^{13} t^3-q^4 t^4-q^5 t^4
-2 q^6 t^4-3 q^7 t^4-2 q^8 t^4-3 q^9 t^4+q^{10} t^4+10 q^{11} t^4
+14 q^{12} t^4+9 q^{13} t^4+3 q^{14} t^4-q^5 t^5-q^6 t^5-2 q^7 t^5
-3 q^8 t^5-3 q^9 t^5-5 q^{10} t^5-4 q^{11} t^5+6 q^{12} t^5
+16 q^{13} t^5+13 q^{14} t^5+4 q^{15} t^5-q^6 t^6-q^7 t^6-2 q^8 t^6
-3 q^9 t^6-3 q^{10} t^6-6 q^{11} t^6-7 q^{12} t^6+2 q^{13} t^6
+16 q^{14} t^6+15 q^{15} t^6+4 q^{16} t^6-q^7 t^7-q^8 t^7-2 q^9 t^7
-3 q^{10} t^7-3 q^{11} t^7-6 q^{12} t^7-8 q^{13} t^7+q^{14} t^7
+16 q^{15} t^7+13 q^{16} t^7+3 q^{17} t^7-q^8 t^8-q^9 t^8
-2 q^{10} t^8-3 q^{11} t^8-3 q^{12} t^8-6 q^{13} t^8-8 q^{14} t^8
+2 q^{15} t^8+16 q^{16} t^8+9 q^{17} t^8+q^{18} t^8-q^9 t^9
-q^{10} t^9-2 q^{11} t^9-3 q^{12} t^9-3 q^{13} t^9-6 q^{14} t^9
-7 q^{15} t^9+6 q^{16} t^9+14 q^{17} t^9+4 q^{18} t^9-q^{10} t^{10}
-q^{11} t^{10}-2 q^{12} t^{10}-3 q^{13} t^{10}-3 q^{14} t^{10}
-6 q^{15} t^{10}-4 q^{16} t^{10}+10 q^{17} t^{10}+9 q^{18} t^{10}
+q^{19} t^{10}-q^{11} t^{11}-q^{12} t^{11}-2 q^{13} t^{11}
-3 q^{14} t^{11}-3 q^{15} t^{11}-5 q^{16} t^{11}+q^{17} t^{11}
+11 q^{18} t^{11}+3 q^{19} t^{11}-q^{12} t^{12}-q^{13} t^{12}
-2 q^{14} t^{12}-3 q^{15} t^{12}-3 q^{16} t^{12}-3 q^{17} t^{12}
+7 q^{18} t^{12}+6 q^{19} t^{12}-q^{13} t^{13}-q^{14} t^{13}
-2 q^{15} t^{13}-3 q^{16} t^{13}-2 q^{17} t^{13}+q^{18} t^{13}
+7 q^{19} t^{13}+q^{20} t^{13}-q^{14} t^{14}-q^{15} t^{14}
-2 q^{16} t^{14}-3 q^{17} t^{14}+5 q^{19} t^{14}+2 q^{20} t^{14}
-q^{15} t^{15}-q^{16} t^{15}-2 q^{17} t^{15}-2 q^{18} t^{15}
+3 q^{19} t^{15}+3 q^{20} t^{15}-q^{16} t^{16}-q^{17} t^{16}
-2 q^{18} t^{16}+4 q^{20} t^{16}-q^{17} t^{17}-q^{18} t^{17}
-q^{19} t^{17}+2 q^{20} t^{17}+q^{21} t^{17}-q^{18} t^{18}
-q^{19} t^{18}+q^{20} t^{18}+q^{21} t^{18}
-q^{19} t^{19}+q^{21} t^{19}-q^{20} t^{20}+q^{21} t^{20}\bigr).
\)
}
\renewcommand{\baselinestretch}{1.0} 
\smallskip

The $a$\~degree of $\hat{\h}{}^{min}_{\l}\,(q,t,a)$ is
deg$_a=\ss^1_1\rr^1_2|\yng(1)\,|+
\ss^2_1\rr^2_2|\yng(1)\,|-|\yng(1)\,|= 3+2-1=4$.
The positivity of 
$\hat{\h}{}^{min}_{\l}(q,t,a)/(1-t)$ holds.
\medskip
}
 
The last one will be uncolored $C\!ab(5,2)T(4,2)$ for
$\Up=\{\circ\rightarrow\circ\rightrightarrows\}$.
The \tax-presentation is as follows: 
$$
1\le j\le \kappa=2,\ \vec\rr^j=\{2,2\},\, \vec\ss^j=\{1,1\},\, 
\Up=\{\circ\rightarrow\circ\rightrightarrows\}\,,\ 
\Yboxdim7pt \la^1\!=\!\yng(1)=\!\la^2.
$$
Now we can the avoid duplication of
notation; $\{2,2\},\{1,1\}$ in $\l$ will mean
 $(\bigl\{\{2,2\},\{1,1\}\bigr\},
\bigl\{\{2,2\},\{1,1\}\bigr\})$ for $2$ paths. 
So the labels are $[2,1]$ and $[2,1]$. The cable is
$C\!ab(5,2)T(4,2)$ for each path.

\begin{align}\label{T5-2-4-2}
&C\!ab(5,2)T(4,2):\ 
\l=\l_{\{2,2\},\{1,1\}}^
{\,\circ\rightarrow\circ\rightrightarrows,\, 
(\yng(1)\,,\yng(1)\,)},\ \ 
\hat{\h}{}^{min}_{\l}\,(q,t,a)=
\end{align}
\smallskip

\renewcommand{\baselinestretch}{0.5} 
\noindent
{\small
\(
1-t+q t+q^2 t+q^3 t-q t^2+2 q^4 t^2+q^5 t^2+q^6 t^2-q^2 t^3
-q^4 t^3+q^5 t^3+q^6 t^3+2 q^7 t^3-q^3 t^4-q^5 t^4+3 q^8 t^4
-q^4 t^5-q^6 t^5-q^8 t^5+3 q^9 t^5-q^5 t^6-q^7 t^6-q^9 t^6
+3 q^{10} t^6-q^6 t^7-q^8 t^7+2 q^{11} t^7-q^7 t^8-q^9 t^8
+q^{11} t^8+q^{12} t^8-q^8 t^9-q^{10} t^9+q^{11} t^9+q^{12} t^9
-q^9 t^{10}-q^{11} t^{10}+2 q^{12} t^{10}-q^{10} t^{11}
+q^{13} t^{11}-q^{11} t^{12}+q^{13} t^{12}-q^{12} t^{13}
+q^{13} t^{13}-q^{13} t^{14}+q^{14} t^{14}+a^3 \bigl(q^6
-q^6 t+q^7 t+q^8 t+q^9 t-q^7 t^2+2 q^{10} t^2-q^8 t^3
-q^{10} t^3+2 q^{11} t^3-q^9 t^4-q^{11} t^4+2 q^{12} t^4
-q^{10} t^5+q^{13} t^5-q^{11} t^6+q^{13} t^6-q^{12} t^7
+q^{13} t^7-q^{13} t^8+q^{14} t^8\bigr)+a^2 \bigl(q^3+q^4+q^5
-q^3 t+q^5 t+3 q^6 t+2 q^7 t+q^8 t-q^4 t^2-q^5 t^2-q^6 t^2
+2 q^7 t^2+3 q^8 t^2+3 q^9 t^2-q^5 t^3-q^6 t^3-2 q^7 t^3+2 q^9 t^3
+4 q^{10} t^3-q^6 t^4-q^7 t^4-2 q^8 t^4-q^9 t^4+q^{10} t^4
+4 q^{11} t^4-q^7 t^5-q^8 t^5-2 q^9 t^5-q^{10} t^5+2 q^{11} t^5
+3 q^{12} t^5-q^8 t^6-q^9 t^6-2 q^{10} t^6+3 q^{12} t^6+q^{13} t^6
-q^9 t^7-q^{10} t^7-2 q^{11} t^7+2 q^{12} t^7+2 q^{13} t^7
-q^{10} t^8-q^{11} t^8-q^{12} t^8+3 q^{13} t^8-q^{11} t^9
-q^{12} t^9+q^{13} t^9+q^{14} t^9-q^{12} t^{10}+q^{14} t^{10}
-q^{13} t^{11}+q^{14} t^{11}\bigr)+a \bigl(q+q^2+q^3-q t+q^3 t
+3 q^4 t+2 q^5 t+q^6 t-q^2 t^2-q^3 t^2-q^4 t^2+2 q^5 t^2+3 q^6 t^2
+4 q^7 t^2+q^8 t^2-q^3 t^3-q^4 t^3-2 q^5 t^3+q^7 t^3+5 q^8 t^3
+2 q^9 t^3-q^4 t^4-q^5 t^4-2 q^6 t^4-q^7 t^4-q^8 t^4+5 q^9 t^4
+2 q^{10} t^4-q^5 t^5-q^6 t^5-2 q^7 t^5-q^8 t^5-2 q^9 t^5
+5 q^{10} t^5+2 q^{11} t^5-q^6 t^6-q^7 t^6-2 q^8 t^6-q^9 t^6
-q^{10} t^6+5 q^{11} t^6+q^{12} t^6-q^7 t^7-q^8 t^7-2 q^9 t^7
-q^{10} t^7+q^{11} t^7+4 q^{12} t^7-q^8 t^8-q^9 t^8-2 q^{10} t^8
+3 q^{12} t^8+q^{13} t^8-q^9 t^9-q^{10} t^9-2 q^{11} t^9
+2 q^{12} t^9+2 q^{13} t^9-q^{10} t^{10}-q^{11} t^{10}
-q^{12} t^{10}+3 q^{13} t^{10}-q^{11} t^{11}-q^{12} t^{11}
+q^{13} t^{11}+q^{14} t^{11}
-q^{12} t^{12}+q^{14} t^{12}-q^{13} t^{13}+q^{14} t^{13}\bigr).
\)
}
\renewcommand{\baselinestretch}{1.0} 
\smallskip

The $a$\~degree of $\hat{\h}{}^{min}_{\l}\,(q,t,a)$ is\, 
$\ss_1\rr_2(|\yng(1)\,|+|\yng(1)\,|)-|\yng(1)\,|=3.$
This polynomial is
self-dual; all uncolored ones are self-dual.
The positivity of the series
$\hat{\h}{}^{min}_{\l}(q,t,a)/(1-t)^{\kappa-1}$ holds here 
and in further (uncolored) examples; recall that $\kappa$ is 
the number of paths (components of the link), which is $2$ in the 
considered case. 
\smallskip

The plane curve singularity for (\ref{T5-2-4-2}) is
$$
\c=\{(y^5+x^2)(y^5-x^2)=0\},\, \
\hbox{Link}(\c)=\bigl\{T(5,2),T(5,2),\, \hbox{lk}=10\bigr\}.
$$
The knots here  (and the corresponding branches of this
singularity)  
and the linking number $lk$ between them (equivalently, the 
intersection multiplicity) uniquely determine the germ of 
this curve. The  Alexander polynomial is
$1+q^4+q^8+q^{10}+q^{12}+q^{14}+q^{16}+q^{18}+q^{22}+q^{26}$,
which coincides with 
$\hat{\h}{}^{min}_{\l}\,(q,q,-1)/(1-q)^2$. 

\setcounter{equation}{0}
\section{\sc Generalized twisted union}
Continuing Section \ref{sec:T(2k,k)}, 
we will provide further examples of the pairs
of trees $\{\l,\,'\!\l\}$ and $\{\l,\,'\!\l^\vee\}$. 
They will be mostly for the {\em uncolored\,}
unknot taken as $\,'\!\l$, i.e. for 
\Yboxdim5pt
$$'\!\l=\l_{1,0}^{\circ\rightarrow,\,\yng(1)\,}.$$

We will simply add ``prime" to the corresponding knot/link
$\l$ (followed by $\vee$ if it is present), to indicate 
using this very $'\!\l$.
Recall that 
$\hat{\h}{}^{min}_{\l,\,'\!\l^\vee}$
is obtained by the substitution $X\mapsto Y$ instead of 
$X\mapsto Y^{-1}$  in $\hat{\h}{}^{min}_{\l,\,'\!\l}$.
Topologically, this is adding the meridian to $\l$. 
Note that all examples of Hopf links from Section 
\ref{sec:torus-hopf} can be restated in terms of such 
pairs due to (\ref{vert-coinv-p}). 
\smallskip

\subsection{\bf Uncolored trefoil-prime} 
The \tax-presentation  
for such $\l$ is as follows:
$
1\le j\le \kappa=2,\ \vec\rr^j=3,\, \vec\ss^j=2,\, 
\Up=\{\circ\rightarrow\}\,,\ 
\Yboxdim7pt \la^1=\yng(1)\,;
$

\Yboxdim5pt
\begin{align}\label{T6-4-y}
&T(3,2)'_{\yng(1)\,}:\ 
\l=\l_{\{3,2\}}^{\,\circ\rightarrow,\, 
\yng(1)\,}\,,\ \ 
\hat{\h}{}^{min}_{\l,\,'\!\l}\,(q,t,a)=
\end{align}

\renewcommand{\baselinestretch}{0.5} 
\noindent
\centerline{\small
\(
1+a^2 \bigl(q^2-\frac{q^2}{t^2}+\frac{q^3}{t^2}+\frac{q}{t}
-\frac{q^3}{t}\bigr)
+q t+a \bigl(2 q+\frac{1}{t}-\frac{q^2}{t}+q^2 t\bigr).
\)
}
\renewcommand{\baselinestretch}{1.0} 
\smallskip

{\sf Adding $\vee$.}
The positivity of the series
$\hat{\h}{}^{min}_{\l,\,'\!\l}(q,t,a)/(1-t)$ from 
(\ref{posit-claim}) in Conjecture \ref{CONCONJ} 
fails (and for any power of $(1-t)$), but holds 
for the $\vee$\~variant of this superpolynomial (as conjectured):
\smallskip

\Yboxdim5pt
\begin{align}\label{T6-4-yv}
&T(3,2)^{\prime,\vee}_{\yng(1)\,}:\ 
\l=\l_{\{3,2\}}^{\,\circ\rightarrow,\, 
\yng(1)\,}\,,\ \ 
\hat{\h}{}^{min}_{\l,\,'\!\l^\vee}\,(q,t,a)=
\end{align}

\renewcommand{\baselinestretch}{0.5} 
\noindent
{\small
\(
1-t+q t+q^2 t-q t^2+q^3 t^2-q^2 t^3+q^3 t^3-q^3 t^4
+q^4 t^4+a^2 \bigl(q^3-q^3 t+q^4 t\bigr)
+a \bigl(q+q^2-q t+2 q^3 t-q^2 t^2+q^4 t^2-q^3 t^3+q^4 t^3\bigr).
\)
}
\renewcommand{\baselinestretch}{1.0} 
\smallskip

\Yboxdim7pt
The $a$\~degree of\, 
$\hat{\h}{}^{min}_{\l,\,'\!\l^\vee}\,(q,t,a)$\, 
is given by the formula deg$_a=\max\{\ss,1\}
|\yng(1)\,|+\max\{\,'\ss,1\}|\yng(1)\,|-|\yng(1)\,|=2+1-1=2.$
It is the same as in the previous example. Both polynomials 
result at $q=1$ in $(1+a)(1+a+t)$; the second factor is the 
specialization at $q=1$ for the uncolored trefoil.

The series
$\hat{\h}{}^{min}_{\l,\,'\!\l^\vee}(q,t,a)/(1-t)$ is positive;
the corresponding germ of the singularity is
$$
\c=\{(y^3+x^2)(y^3+x)=0\},\, \
\hbox{Link}(\c)=\bigl\{T(3,2),T(3,1),\, \hbox{lk}=3\bigr\},
$$
where $T(3,1)=\unknot$\,,
and the Alexander polynomial is $\,1+q^3+q^6$.
\Yboxdim5pt
\smallskip

We note that the superpolynomial for
(\ref{T6-4-yv}) does not coincide with 
\Yboxdim5pt
\begin{align}\label{T3-2-1-1}
&\bigl(C\!ab(3,2)C\!ab(1,1)\bigr)T(1,1):\, 
\l\!=\!\l_{\{1,1\},\{1,2\}}^{\,\twoonetiny,\, 
(\yng(1)\,,\yng(1)\,)},\ 
\hat{\h}{}^{min}_{\l}\,(q,t,a)=
\end{align}
{\small
$$1+a^2 q^3-t+q t+q^2 t-q t^2+q^2 t^2-q^2 t^3+q^3 t^3+a 
(q+q^2-q t+q^2 t+q^3 t-q^2 t^2+q^3 t^2).
$$
}
The corresponding singularity is
$\c=\{(y^3+x^2)(y+x)=0\}$ with the same link components
$\{T(3,2),\, \unknot\,\}$ 
as above but with $\,lk=2$
and the Alexander polynomial $1+q^4$. 
The self-duality and other claims 
in Theorem \ref{STABILIZ} hold in these cases.
\medskip

{\sf Splice interpretation.}
Let $k\in \Z$.  
Recall that formula (\ref{eval-tau}) provides the
coincidence of $\hat{\h}^{min}$ from (\ref{T6-4-yv}) 
with that for
\begin{align}\label{T2-3--1-2}
&\bigl(C\!ab(2,3)C\!ab(k,\!-1)\bigr)T(1,\!-1): \,
\l_{\bigl\{\{1,3\},\{-1,5\}\bigr\},
\bigl\{\{1,0\},\{-1,k-\!1\}\bigr\}}  
^{\,\twotwotiny,\, 
(\yng(1)\,,\yng(1)\,)}.
\end{align}
Explicitly, the following DAHA-identity establishes this
coincidence:
\begin{align*}
&\bigl\{\tau_-^{-1}\bigr(\tau_-^{1-k}\si^2(P_{\yng(1)})
\!\Downarrow 
\tau_-\tau_+\tau_-^{2}(P_{\yng(1)})\!\Downarrow\bigr)\bigr\}_{ev}
\!=\!
q^\bullet t^\bullet\bigl\{P^\iota_{\yng(1)}\!\Downarrow\,,\,
\tau_+\tau_-^{2}(P_{\yng(1)})\!\Downarrow\bigr\}_{ev},
\end{align*}
where we use the commutativity of $\tau_-$ with the projection
$H\!\!\Downarrow=H(1)$ of $H\in \HH\,$ onto $\v$, i.e. the
fundamental fact that $\tau_-$ acts in $\v$; this action was
denoted by  $\dot{\tau}_-$ in  (\ref{eval-tau}). Recall that
$P^\iota=P(X^{-1}),\,'\!\l$ is $\circ\rightarrow$, and the 
colors are trivial (shown by  
\Yboxdim7pt
$\yng(1)\,$ in the diagrams below).
\Yboxdim5pt

Let us translate the latter identity into the language of 
splice diagrams. The corresponding diagrams for $k=1$ are: 
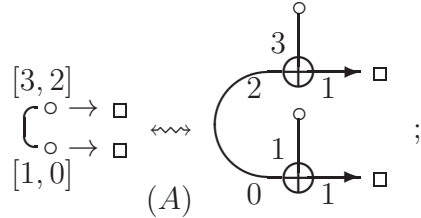
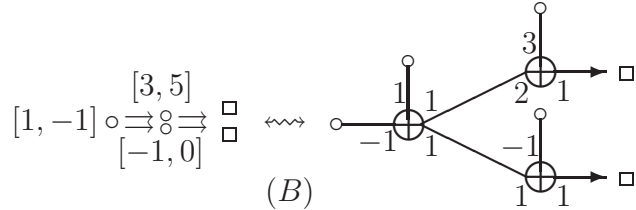
\begin{align*}
\begin{picture}(45,85)
\put(0,20){
    \put(8,11){\oval(10,15)[l]}
    \put(10,0){$\circ\rightarrow 
\,\yng(1)$}
    \put(10,15){$\circ\rightarrow 
\,\yng(1)$}
    }
\put(0,45){$\![3,2]$}
\put(0,10){$\![1,0]$}
\end{picture}
&\begin{picture}(20,35)
\put(5,27){$\leftrightsquigarrow$}$\  (A)$
\end{picture} 
\begin{picture}(80,85)
\put(35,32){\oval(50,40)[l]}
\put(10,3){  
    \put(35,9){\vector(1,0){20}}
    \put(31.5,11){\line(0,1){20}}
    \put(25,6){$\bigoplus$}
    \put(40,-1){$1$}
    \put(60,6){$\yng(1)$}
    \put(12,-1){$0$}
    \put(21,16){$1$}
    \put(29,30){$\circ$}
    }
\put(10,43){
    \put(35,9){\vector(1,0){20}}
    \put(31.5,11){\line(0,1){20}}
    \put(25,6){$\bigoplus$}
    \put(40,-1){$1$}
    \put(60,6){$\yng(1)$}
    \put(12,-1){$2$}
    \put(21,16){$3$}
    \put(29,30){$\circ$}
    }
\end{picture}
\begin{picture}(20,35)
\put(5,27){$;$}
\end{picture}\\
\begin{picture}(90,85)
\put(0,15){
    \put(35,12){$\twotwo$}
    \put(0,12){$[1,-1]$}
    \put(40,0){$[-1,0]$}   
    \put(45,25){$[3,5]$}
    \put(80,18){$\yng(1)$}
    \put(80,8){$\yng(1)$}
    }
\end{picture}
&\begin{picture}(20,35)
\put(5,27){$\leftrightsquigarrow$}$\ \,(B)$
\end{picture} 
\begin{picture}(130,40)
\put(30,20)
{
    \put(15,9){\line(2,1){40}}
    \put(15,9){\line(2,-1){40}}
    \put(5,9){\line(-1,0){20}}
    \put(10,11){\line(0,1){20}}
    \put(3.5,6){$\bigoplus$}
    \put(4,17){$1$}
    \put(16,-3){$1$}
    \put(16,14){$1$}
    \put(-9,-1){$-1$}
    \put(-20,6){$\circ$}
    \put(7,30){$\circ$}
    }
\put(80,40){
    \put(15,9){\vector(1,0){20}}
    \put(10,11){\line(0,1){20}}
    \put(3.5,6){$\bigoplus$}
    \put(16,-2){$1$}
    \put(40,6){$\yng(1)$}
    \put(0,-3){$2$}
    \put(3,16){$3$}
    \put(7,30){$\circ$}
    }
\put(80,0){
    \put(15,9){\vector(1,0){20}}
    \put(10,11){\line(0,1){20}}
    \put(3.5,6){$\bigoplus$}
    \put(16,-2){$1$}
    \put(40,6){$\yng(1)$}
    \put(0,-2){$1$}
    \put(-5,17){$-1$}
    \put(7,30){$\circ$}
    }
\end{picture}
\begin{picture}(10,45)
\put(5,22){$.$}
\end{picture}
\end{align*}

Recall that the colors are assigned to the arrowheads; also
we will use $\vee$ for the change of orientation of
the corresponding component.
Generally, any {\em leaves\,} of weight $1$ can be
deleted without changing 
the link (see here and below \cite{EN}, Theorem 8.1, 
Statements 2,3). The situation with leaves of weight $-1$ is 
somewhat more involved. 
Namely, one can change weights of two edges from the same 
node to their negatives together with the simultaneous change of 
orientations of all components starting with these 
edges. Starting with the diagram $(B)$, this results in the 
following (isotopic) equivalence
$\thickapprox\hskip -2.5pt\thickapprox$\,:

\begin{equation*}
\begin{picture}(130,80)
\put(30,30)
{
    \put(15,9){\line(2,1){40}}
    \put(15,9){\line(2,-1){40}}
    \put(5,9){\line(-1,0){20}}
    \put(10,11){\line(0,1){20}}
    \put(3.5,6){$\bigoplus$}
    \put(4,17){$1$}
    \put(16,-3){$1$}
    \put(16,14){$1$}
    \put(-9,-1){$-1$}
    \put(-20,6){$\circ$}
    \put(7,30){$\circ$}
    }
\put(80,50){
    \put(15,9){\vector(1,0){20}}
    \put(10,11){\line(0,1){20}}
    \put(3.5,6){$\bigoplus$}
    \put(16,-2){$1$}
    \put(40,6){$\yng(1)$}
    \put(0,-3){$2$}
    \put(3,16){$3$}
    \put(7,30){$\circ$}
    }
\put(80,10){
    \put(15,9){\vector(1,0){20}}
    \put(10,11){\line(0,1){20}}
    \put(3.5,6){$\bigoplus$}
    \put(16,-2){$1$}
    \put(40,6){$\yng(1)$}
    \put(0,-2){$1$}
    \put(-5,17){$-1$}
    \put(7,30){$\circ$}
    }
\end{picture}
\begin{picture}(20,45)
\put(5,37){$\thickapprox\hskip -2.5pt\thickapprox$}
\end{picture} 
\begin{picture}(130,80)
\put(30,30)
{
    \put(15,9){\line(2,1){40}}
    \put(15,9){\line(2,-1){40}}
    \put(5,9){\line(-1,0){20}}
    \put(10,11){\line(0,1){20}}
    \put(3.5,6){$\bigoplus$}
    \put(4,17){$1$}
    \put(11,-7){$-1$}
    \put(16,14){$1$}
    \put(0,-1){$1$}
    \put(-20,6){$\circ$}
    \put(7,30){$\circ$}
    }
\put(80,50){
    \put(15,9){\vector(1,0){20}}
    \put(10,11){\line(0,1){20}}
    \put(3.5,6){$\bigoplus$}
    \put(16,-2){$1$}
    \put(40,6){$\yng(1)$}
    \put(0,-3){$2$}
    \put(3,16){$3$}
    \put(7,30){$\circ$}
    }
\put(80,10){
    \put(15,9){\vector(1,0){20}}
    \put(10,11){\line(0,1){20}}
    \put(3.5,6){$\bigoplus$}
    \put(16,-2){$1$}
    \put(40,6){$\yng(1)^\vee$}
    \put(0,-2){$1$}
    \put(-5,17){$-1$}
    \put(7,30){$\circ$}
    }
\end{picture}
\end{equation*}
\begin{align}
\begin{picture}(20,45)
\put(5,37){$\thickapprox\hskip -2.5pt\thickapprox$}
\end{picture}
\begin{picture}(130,80)
\put(30,30)
{
    \put(15,9){\line(2,1){40}}
    \put(15,9){\line(2,-1){40}}
    \put(5,9){\line(-1,0){20}}
    \put(10,11){\line(0,1){20}}
    \put(3.5,6){$\bigoplus$}
    \put(4,17){$1$}
    \put(11,-7){$-1$}
    \put(16,14){$1$}
    \put(0,-1){$1$}
    \put(-20,6){$\circ$}
    \put(7,30){$\circ$}
    }
\put(80,50){
    \put(15,9){\vector(1,0){20}}
    \put(10,11){\line(0,1){20}}
    \put(3.5,6){$\bigoplus$}
    \put(16,-2){$1$}
    \put(40,6){$\yng(1)$}
    \put(0,-3){$2$}
    \put(3,16){$3$}
    \put(7,30){$\circ$}
    }
\put(80,10){
    \put(15,9){\vector(1,0){20}}
    \put(10,11){\line(0,1){20}}
    \put(3.5,6){$\bigoplus$}
    \put(9,-2){$-1$}
    \put(40,6){$\yng(1)$}
    \put(0,-2){$1$}
    \put(4,17){$1$}
    \put(7,30){$\circ$}
    }
\end{picture}
\begin{picture}(20,45)
\put(5,37){$\thickapprox\hskip -2.5pt\thickapprox$}
\end{picture} \ \ 
\begin{picture}(75,80)
\put(30,30){
    \put(5,9){\vector(-1,0){20}}
    \put(15,9){\vector(1,0){20}}
    \put(10,11){\line(0,1){20}}
    \put(3.5,6){$\bigoplus$}
    \put(16,-2){$1$}
    \put(40,6){$\yng(1)$}
    \put(0,-3){$2$}
    \put(3,16){$3$}
    \put(7,30){$\circ$}
    \put(-25,6){$\yng(1)$}
    }
\end{picture}\label{fig:splice}
\begin{picture}(10,45)
\put(5,32){$.$}
\end{picture}
\end{align}

At the last step here, we used Statement 3 of Theorem 8.1 from 
\cite{EN} three times. Applying it one more time, we conclude that  
the diagram/link in $(A)$  is equivalent to that in $(B)$. 
\smallskip

\subsection{\bf Colored/iterated examples}
Let us begin with the  trefoil-prime for the 2-row:
\Yboxdim5pt
\begin{align}\label{T6-4-y-11}
&T(3,2)'_{\yng(2)\,}\, :\ 
\l=\l_{\{3,2\}}^{\,\circ\rightarrow,\, 
\yng(2)\,}\,,\ \ 
\hat{\h}{}^{min}_{\l,\,'\!\l}\,(q,t,a)=
\end{align}

\renewcommand{\baselinestretch}{0.5} 
\noindent
{\small
\(
1+a^3 \bigl(q^7-\frac{q^6}{t^2}+\frac{q^8}{t^2}
+\frac{q^5}{t}+\frac{q^6}{t}-\frac{q^7}{t}
-\frac{q^8}{t}\bigr)+q^2 t+q^3 t+q^4 t^2
+a^2 \bigl(q^4+3 q^5+q^6-q^7-q^8-\frac{q^4}{t^2}
+\frac{q^6}{t^2}+\frac{q^2}{t}+\frac{q^3}{t}+\frac{q^4}{t}
-\frac{q^5}{t}-\frac{3 q^6}{t}+\frac{q^8}{t}+q^6 t+q^7 t\bigr)
+a \bigl(2 q^2+2 q^3+q^4-q^5-q^6
+\frac{1}{t}-\frac{q^4}{t}+2 q^4 t+2 q^5 t+q^6 t^2\bigr).
\)
}
\renewcommand{\baselinestretch}{1.0} 
\smallskip

\Yboxdim7pt
It is super-dual to the superpolynomial for $\yng(1,1)\,$.
\Yboxdim5pt
The positivity of the series
$\hat{\h}{}^{min}_{\l,\,'\!\l}(q,t,a)/(1-t)^p$ fails
for $p=1$ and any $p>1$
(due to the absence of $\vee$).
At $q=1$, it equals $(1+a)(1+a+t)^2$, which also coincides
with the specialization for the corresponding $\vee$\~polynomial: 
\smallskip

\Yboxdim5pt
\begin{align}\label{T6-4-yv-11}
&T(3,2)^{\prime,\vee}_{\yng(2)\,}\, :\ 
\l=\l_{\{3,2\}}^{\,\circ\rightarrow,\, 
\yng(2)\,}\,,\ \ 
\hat{\h}{}^{min}_{\l,\,'\!\l^\vee}\,(q,t,a)=
\end{align}

\renewcommand{\baselinestretch}{0.5} 
\noindent
{\small
\(
1-t+q^2 t+q^3 t+q^4 t-q^2 t^2-q^3 t^2+q^5 t^2+2 q^6 t^2
-q^4 t^3-q^5 t^3+q^7 t^3+q^8 t^3-q^6 t^4-q^7 t^4+q^8 t^4+q^9 t^4
-q^8 t^5+q^{10} t^5+a^3 \bigl(q^9-q^9 t+q^{11} t\bigr)
+a^2 \bigl(q^5+q^6+q^7-q^5 t-q^6 t+2 q^8 t+2 q^9 t-q^7 t^2-q^8 t^2
+q^{10} t^2+q^{11} t^2-q^9 t^3+q^{11} t^3\bigr)+a \bigl(q^2+q^3
+q^4-q^2 t-q^3 t+2 q^5 t+3 q^6 t+q^7 t-q^4 t^2-2 q^5 t^2-q^6 t^2
+q^7 t^2+3 q^8 t^2+q^9 t^2-q^6 t^3-2 q^7 t^3
+2 q^9 t^3+q^{10} t^3-q^8 t^4-q^9 t^4+q^{10} t^4+q^{11} t^4\bigr).
\)
}
\renewcommand{\baselinestretch}{1.0} 
\smallskip

\Yboxdim7pt
The $a$\~degree of\, 
$\hat{\h}{}^{min}_{\l,\,'\!\l^\vee}\,(q,t,a)$\, 
is given by the formula deg$_a=\max\{\ss,1\}
|\yng(2)\,|+\max\{\,'\ss,1\}|\yng(1)\,|-|\yng(2)\,|=4+1-2=3.$
See ((\ref{deg-a-j}) and (\ref{deg-a-jj}). 
This polynomial is dual to the superpolynomial for $\,\yng(1,1)\,$.
\Yboxdim5pt
The positivity of the series
$\hat{\h}{}^{min}_{\l,\,'\!\l^\vee}(q,t,a)/(1-t)$ holds.
\medskip

{\sf Iterated T(2,1)-prime.}
Let us replace the additional arrowhead in
the example from (\ref{T8-3-2-1}) by the prime-construction
(twisted union with one box)
for the starting $\vec{\rr}=\{2,3\},\vec{\ss}=\{1,2\}$.
Recall that all prime examples are for \Yboxdim5pt
$'\!\l=\l_{1,0}^{\circ\rightarrow,\,\yng(1)\,}\,$.
 
\Yboxdim5pt
\begin{align}\label{T4-2-y-3-2}
&C\!ab(8,3)(T(2,1))'_{\yng(1)\,}\, :\ 
\l=\l_{\{2,3\},\{1,2\}}^{\,\circ\rightarrow\circ\rightarrow,\, 
\yng(1)\,}\,,\ \ 
\hat{\h}{}^{min}_{\l,\,'\!\l}\,(q,t,a)=
\end{align}

\renewcommand{\baselinestretch}{0.5} 
\noindent
{\small
\(
1+q t+q^2 t+q^2 t^2+q^3 t^2+q^4 t^2+q^3 t^3+q^4 t^3+q^5 t^3+q^4 t^4
+q^5 t^4+q^5 t^5+q^6 t^5+q^6 t^6+q^7 t^7+a^4 \bigl(q^7-q^8
-\frac{q^6}{t^2}+\frac{q^8}{t^2}+\frac{q^6}{t}-\frac{q^7}{t}
-\frac{q^8}{t}+\frac{q^9}{t}+q^8 t-q^9 t\bigr)+a^3 \bigl(q^4+2 q^5
+2 q^6-4 q^7-2 q^8+2 q^9-\frac{q^4}{t^2}-\frac{q^5}{t^2}
+\frac{2 q^7}{t^2}+\frac{q^3}{t}+\frac{q^4}{t}-\frac{3 q^6}{t}
-\frac{2 q^7}{t}+\frac{3 q^8}{t}+q^5 t+2 q^6 t+2 q^7 t-3 q^8 t
+q^6 t^2+2 q^7 t^2-q^9 t^2+q^7 t^3+q^8 t^3-q^9 t^3+q^8 t^4\bigr)
+a^2 \bigl(q^2+3 q^3+3 q^4-4 q^6-2 q^7+2 q^8-\frac{q^3}{t^2}
+\frac{q^6}{t^2}+\frac{q}{t}+\frac{q^2}{t}+\frac{q^3}{t}
-\frac{2 q^4}{t}-\frac{2 q^5}{t}-\frac{q^6}{t}+\frac{2 q^7}{t}
+q^3 t+3 q^4 t+5 q^5 t+q^6 t-4 q^7 t-q^8 t+q^9 t+q^4 t^2+3 q^5 t^2
+5 q^6 t^2-2 q^8 t^2+q^5 t^3+3 q^6 t^3+3 q^7 t^3-2 q^8 t^3+q^6 t^4
+3 q^7 t^4+q^8 t^4-q^9 t^4+q^7 t^5+q^8 t^5+q^8 t^6\bigr)
+a \bigl(2 q+2 q^2+q^3-q^4-q^5+\frac{1}{t}-\frac{q^3}{t}+2 q^2 t
+3 q^3 t+3 q^4 t-2 q^6 t+2 q^3 t^2+3 q^4 t^2+4 q^5 t^2-q^7 t^2
+2 q^4 t^3+3 q^5 t^3+3 q^6 t^3-q^7 t^3+2 q^5 t^4
+3 q^6 t^4+q^7 t^4-q^8 t^4
+2 q^6 t^5+2 q^7 t^5+2 q^7 t^6+q^8 t^7\bigr).
\)
}
\renewcommand{\baselinestretch}{1.0} 
\smallskip

\Yboxdim7pt
It is self-dual; the positivity of the series
$\hat{\h}{}^{min}_{\l,\,'\!\l}(q,t,a)/(1-t)^p$ fails
for any $p\ge 0$, 
which failure is generally expected without $\vee$.
As in (\ref{4-2-T21-1-1}), 
the $a$\~degree of\, 
$\hat{\h}{}^{min}_{\l,\,'\!\l}\,(q,t,a)$, which is $4$,
is greater than 
deg$_a=\ss_1\rr_2
|\yng(1)\,|+\max\{\,'\ss,1\}|\yng(1)\,|-|\yng(1)\,|=3+1-1=3.$
\smallskip

The $\vee$\~version is generally  
with $X$ replaced by $Y$ in the {\em pre-polynomial}
$\,'\P_0^{tot}$ for $\,'\!\l$. This pre-polynomial
is \Yboxdim5pt $J_{\yng(1)\,}$ (up to a renormalization). 
We will omit the corresponding 
\begin{align}\label{T4-2-yv-3-2}
&\hat{\h}{}^{min}_{\l,\,'\!\l^\vee}\,(q,t,a)
\for C\!ab(8,3)(T(2,1))^{\prime,\vee}_{\yng(1)\,}\, :\ 
\l=\l_{\{2,3\},\{1,2\}}^{\,\circ\rightarrow\circ\rightarrow,\, 
\yng(1)\,}\,.\ \ 
\end{align}
It coincides with 
$\hat{\h}{}^{min}_{\mathcal{M}}\,(q,t,a)$ for
the polynomial from (\ref{T8-3-2-1}). The latter was
defined (and calculated) for
\begin{align}\label{T8-3-2-1X}
&\bigl(C\!ab(8,3)C\!ab(0,1)\bigr)T(2,1):\ 
\mathcal{M}=\l_{\{2,3\},\{1,2\}}^
{\,\twoonetiny,\, 
(\yng(1)\,,\yng(1)\,)}.
\end{align}

Let us demonstrate this coincidence theoretically. 
Setting $P=P_{\yng(1)\,}/\{P\}_{ev},\, P^\iota=\iota(P),\,
Q=\tau_+\tau_-^2(P_{\yng(1)\,})\!\Downarrow=
\tau_+\tau_-^2(P_{\yng(1)\,})(1)\in\v$, 
\begin{align*}
&\tau_+\tau_-(PQ)=\tau_+\tau_-\bigl(\vph\si^{-1}(P^\iota)Q\bigr)
=\vph\tau_-^2(P^\iota) \tau_+\tau_-(Q),\\
\{&\tau_+\tau_-(PQ)\}_{ev}
=\{\vph\bigl(\dot{\tau}^2_-(P^\iota)\bigr) 
\bigl(\tau_+\tau_-(Q)\!\Downarrow\bigr)\}_{ev}
\simeq \{P^\iota, \tau_+\tau_-(Q)\!\Downarrow\}_{ev},
\end{align*}
where $\dot{\tau}_-$ is the action of $\tau_-$ in $\v$
and $\simeq$ is the equality up to $q^\bullet t^\bullet$.
Modulo this equivalence, the last quantity coincides with
$\hat{\h}^{min}_{\l,\,'\!\l^\vee}$. Recall that 
$P^\iota$ can be replaced by
$E_{\om_n}/\{E_{\om_n}\,\}_{ev}$ here and the
$E$\~polynomials are $\dot{\tau}_-$\~invariant 
up to $\simeq$\,; see (\ref{taumineb}) and also
(\ref{iden2-1}).

\comment{
\renewcommand{\baselinestretch}{0.5} 
\noindent
{\small
\(
1-t+q t+q^2 t+q^3 t-q t^2+2 q^4 t^2+q^5 t^2+q^6 t^2-q^2 t^3-q^4 t^3
+q^5 t^3+2 q^6 t^3+2 q^7 t^3-q^3 t^4-q^5 t^4+q^7 t^4+3 q^8 t^4
-q^4 t^5-q^6 t^5+q^8 t^5+3 q^9 t^5-q^5 t^6-q^7 t^6+q^9 t^6
+2 q^{10} t^6-q^6 t^7-q^8 t^7+2 q^{10} t^7+q^{11} t^7-q^7 t^8
-q^9 t^8+q^{10} t^8+q^{11} t^8-q^8 t^9-q^{10} t^9+2 q^{11} t^9
-q^9 t^{10}+q^{12} t^{10}-q^{10} t^{11}+q^{12} t^{11}-q^{11} t^{12}
+q^{12} t^{12}-q^{12} t^{13}+q^{13} t^{13}+a^3 \bigl(q^6-q^6 t
+q^7 t
+q^8 t+q^9 t-q^7 t^2+q^9 t^2+2 q^{10} t^2-q^8 t^3+2 q^{11} t^3
-q^9 t^4+q^{11} t^4+q^{12} t^4-q^{10} t^5+q^{12} t^5-q^{11} t^6
+q^{12} t^6-q^{12} t^7+q^{13} t^7\bigr)+a^2 \bigl(q^3+q^4+q^5
-q^3 t
+q^5 t+3 q^6 t+2 q^7 t+q^8 t-q^4 t^2-q^5 t^2-q^6 t^2+3 q^7 t^2
+4 q^8 t^2+3 q^9 t^2-q^5 t^3-q^6 t^3-2 q^7 t^3+q^8 t^3+4 q^9 t^3
+4 q^{10} t^3-q^6 t^4-q^7 t^4-2 q^8 t^4+q^9 t^4+4 q^{10} t^4
+3 q^{11} t^4-q^7 t^5-q^8 t^5-2 q^9 t^5+q^{10} t^5+4 q^{11} t^5
+q^{12} t^5-q^8 t^6-q^9 t^6-2 q^{10} t^6+3 q^{11} t^6+2 q^{12} t^6
-q^9 t^7-q^{10} t^7-q^{11} t^7+3 q^{12} t^7-q^{10} t^8-q^{11} t^8
+q^{12} t^8+q^{13} t^8-q^{11} t^9+q^{13} t^9-q^{12} t^{10}
+q^{13} t^{10}\bigr)+a \bigl(q+q^2+q^3-q t+q^3 t+3 q^4 t+2 q^5 t
+q^6 t-q^2 t^2-q^3 t^2-q^4 t^2+2 q^5 t^2+4 q^6 t^2+4 q^7 t^2
+q^8 t^2-q^3 t^3-q^4 t^3-2 q^5 t^3+3 q^7 t^3+6 q^8 t^3+2 q^9 t^3
-q^4 t^4-q^5 t^4-2 q^6 t^4-q^7 t^4+2 q^8 t^4+6 q^9 t^4+2 q^{10} t^4
-q^5 t^5-q^6 t^5-2 q^7 t^5-q^8 t^5+2 q^9 t^5+6 q^{10} t^5
+q^{11} t^5
-q^6 t^6-q^7 t^6-2 q^8 t^6-q^9 t^6+3 q^{10} t^6+4 q^{11} t^6
-q^7 t^7
-q^8 t^7-2 q^9 t^7+4 q^{11} t^7+q^{12} t^7-q^8 t^8-q^9 t^8
-2 q^{10} t^8+2 q^{11} t^8+2 q^{12} t^8-q^9 t^9-q^{10} t^9
-q^{11} t^9+3 q^{12} t^9-q^{10} t^{10}-q^{11} t^{10}+q^{12} t^{10}
+q^{13} t^{10}
-q^{11} t^{11}+q^{13} t^{11}-q^{12} t^{12}+q^{13} t^{12}\bigr).
\)
}
\renewcommand{\baselinestretch}{1.0} 
\smallskip
}

\Yboxdim7pt

\smallskip
Continuing with
(\ref{T4-2-y-3-2}),(\ref{T4-2-yv-3-2}), let us provide
their values at $q=1$: 
$$\ \ \ \ \hat{\h}{}^{min}_{\l,\,'\!\l}\,(q=1,t,a)=
\hat{\h}{}^{min}_{\l,\,'\!\l^\vee}\,(q=1,t,a)=$$
{\small
\begin{align*}
&\bigl(1+a\bigr)\,\bigr(1+2 t+3 t^2+3 t^3+2 t^4
+2 t^5+t^6+t^7\\+a^2(1+2 t
+&\,2 t^2+t^3+t^4)+a (2+4 t+5 t^2+4 t^3+3 t^4+2 t^5+t^6)\bigr)
\end{align*}}
\vskip -0.9cm
$$=\bigl(1+a\bigr)\hat{\h}{}^{min}_{\l}\,(q=1,t,a).
$$
\Yboxdim5pt
Recall that factor $(1+a)$ here is due to our using
$J_{\yng(1)}$ instead of spherical $P_{\yng(1)}^\circ$
(which would give $1$ for $'\!\l=\unknot$); the division
by the $LCM$ of the evaluations of $J$\~polynomials 
in $\h^{min}$ is only by one $J_{\yng(1)}(t^\rho)$ here.
\Yboxdim7pt

The second factor is relatively long.
It is the evaluation at $q=1$ of 
$\hat{\h}^{min}_{8,3}(\yng(1)\,)$ for the torus knot $T(8,3)$. 
Indeed,
$C\!ab(8,3)(T(2,1))$ is isotopic to $T(8,3)$ or $T(3,8)$, where
the latter presentation  is somewhat  more convenient practically; 
the corresponding $\ga_{3,8}$ corresponds to 
$\tau_-^2\tau_+\tau_-^2$. For
the sake of completeness,  
\ $\hat{\h}^{min}_{8,3}(\yng(1)\,)=$

\renewcommand{\baselinestretch}{0.5} 
\noindent
{\small
\(
1+q t+q^2 t+q^2 t^2+q^3 t^2+q^4 t^2+q^3 t^3+q^4 t^3+q^5 t^3+q^4 t^4
+q^5 t^4+q^5 t^5+q^6 t^5+q^6 t^6+q^7 t^7+a^2 \bigl(q^3+q^4 t+q^5 t
+q^5 t^2+q^6 t^2+q^6 t^3+q^7 t^4\bigr)+a \bigl(q+q^2+q^2 t+2 q^3 t
+q^4 t+q^3 t^2+2 q^4 t^2+2 q^5 t^2+q^4 t^3
+2 q^5 t^3+q^6 t^3+q^5 t^4+2 q^6 t^4+q^6 t^5+q^7 t^5+q^7 t^6\bigr).
\)
}
\renewcommand{\baselinestretch}{1.0} 
\smallskip

\Yboxdim5pt
{\sf Double-prime T(2,1).}
Let us provide at least one example of twisted union
with $\,'\!\l$ colored by $\yng(1,1)\,$. We will use
double-prime then:

\begin{align}\label{4-2-T2011-1-1v}
&T(4,2)_{\yng(2)\,,\,\yng(1)\,}^{\prime\prime,\vee}\, 
:\ \,\Yboxdim5pt
\l=\l_{2,1}^{\,\circ\rightrightarrows,\, 
(\yng(2)\,,\,\yng(1)\,)},\ 
'\!\l=\l_{1,0}^{\circ\rightarrow,\,\yng(1,1)\,}\,,\ 
\hat{\h}{}^{min}_{\l,\,'\!\l^\vee}\,(q,t,a)=
\end{align}

\renewcommand{\baselinestretch}{0.5} 
\noindent
{\small
\(
1-t+q^2 t-t^2-q^2 t^2+q^3 t^2+q^4 t^2+t^3-q^2 t^3-q^3 t^3+q^5 t^3
+q^2 t^4-q^3 t^4-2 q^4 t^4+2 q^6 t^4+q^3 t^5-2 q^5 t^5+q^7 t^5
+q^4 t^6-2 q^6 t^6+q^8 t^6+q^5 t^7-q^6 t^7-q^7 t^7+q^8 t^7+q^6 t^8
-q^7 t^8-q^8 t^8+q^9 t^8+q^7 t^{10}-q^8 t^{10}-q^9 t^{10}
+q^{10} t^{10}+a^2 \bigl(q^5-q^5 t+q^7 t-q^5 t^2+q^8 t^2+q^5 t^3
-2 q^7 t^3+q^9 t^3+q^7 t^5-q^8 t^5-q^9 t^5+q^{10} t^5\bigr)
+a \bigl(q^2+q^3-q^2 t-q^3 t+q^4 t+q^5 t-q^2 t^2-q^3 t^2-q^4 t^2
+q^5 t^2+2 q^6 t^2+q^2 t^3+q^3 t^3-q^4 t^3-3 q^5 t^3+2 q^7 t^3
+q^4 t^4-q^5 t^4-2 q^6 t^4+2 q^8 t^4+2 q^5 t^5-q^6 t^5-3 q^7 t^5
+q^8 t^5+q^9 t^5+q^6 t^6-q^7 t^6-q^8 t^6+q^9 t^6+q^7 t^7
-q^8 t^7-q^9 t^7+q^{10} t^7+q^7 t^8-q^8 t^8-q^9 t^8
+q^{10} t^8\bigr).
\)
}
\renewcommand{\baselinestretch}{1.0} 
\smallskip

The $a$\~degree is $2=|\yng(1,1)\,|+|\yng(2)\,|+|\yng(1)\,|
-|\yng(2,1)\,|$. The positivity
of $\hat{\h}{}^{min}_{\l,\,'\!\l}(q,t,a)/(1-t)^2$ holds
(as for all uncolored diagrams).
The super-duality ($q\leftrightarrow t^{-1}$) is true.
\medskip

\subsection{\bf Generalized twisting}
\label{sec:twist}
The following generalization of  
twisted unions  $\{\l,\,'\!\l\}$ and  $\{\l,\,'\!\l^\vee\}$
is very natural algebraically.

Recall that the $\hat{\h}$\~polynomials for the standard
twisted union 
are essentially $\{P(Y^{-1})Q(X)\}_{ev}$
for the {\em pre-polynomials\,} $P=\,'\hat{\P}_{0}^{tot}$ and
$Q=\hat{\P}_{0}^{tot}$ from (\ref{jones-hat}) defined
via the graphs $\,'\!\l$ and $\l$.

The division by the evaluations is used here; it can be
for $j_o$ as in (\ref{jones-hat})
or the one used in the construction of $\hat{\h}^{min}$ 
from Theorem \ref{MAINTHM}, where the division is
by the $LC\!M$ of the evaluations 
of all $J$\~polynomials involved in both, $\,'\!\l$ and $\l$.

When $Y^{-1}$ is changed here to $Y$, we add $\vee$ to
$\,'\!\l$. The {\em twisting\,} here is the
application  $\si^{\pm 1}: X\mapsto Y^{\mp 1}$ to
$P$ (corresponding to $\,'\!\l$).
Let us replace $\si^{\pm 1}$ 
by $\hat{\xi}$ for an arbitrary \,$\xi\in PSL_2(\Z)$. 
Let $(\al,\be)^{tr}$ be the first column of $\xi$.
Here and below we will omit hat in $\hat{\xi}$ (which
means the standard lift of $\xi$ to an automorphism 
of $\HH$).

Given $\xi\in PSL_2(\Z)$ and graphs $\l,\,'\!\l$
(labeled, with arrowheads colored by Young diagrams), 
we define  
$\ \lxi\hat{\h}_{\l,\,'\!\l}^{min}(q,t,a)$ as 
$\{\xi\bigl(P(X)\bigr)Q(X)\}_{ev}$ in the
normalization of Theorem  \ref{MAINTHM}.
It is a polynomial
in terms of $q,t^{\pm1},a\,$
depending only on $\{\al,\be\}$.
We say that $\l$ is {\em $\xi$\~twisted\,} by 
$\,'\!\l$.
Note that 
$
\lxi\hat{\h}_{\emptyset,\,\l}^{min}(q,t,a)=
\lxi\hat{\h}_{\l}^{min}(q,t,a)= 
\lxi\hat{\h}_{\l,\,\emptyset}^{min}(q,t,a).
$

This construction can be reduced to 
the $\si^{-1}$\~twisted union. Namely,
\begin{align}\label{xitwist}
\ \lxi\hat{\h}_{\l,\,'\!\l}^{min}(q,t,a)=
\hat{\h}_{\l,\,'\!\mathcal{M}^\vee}^{min}(q,t,a)
\for '\!\mathcal{M}=[\be,\al]\,\rightdotsarrow'\!\l,
\end{align} 
i.e. for $'\!\mathcal{M}$ obtained from 
$\,'\!\l$ by adding $[\be,\al]$ as the first vertex 
(connected with all paths). If $\al,\be>0$ and 
$\al\ss_1>\be\rr_1$ or 
$\be=1,\al=0$  and $\, '\ss_1\ss_1>\,'\rr_1\rr_1$,  
then we call 
$\{\xi,\l,\,'\!\tilde{\l}^\vee\}$ {\em positive\,} if
the trees $\l,\,'\!\l$  are positive.
This is sufficient
for the corresponding link to be algebraic
(necessary if $\l,\,'\!\l$ are reduced \cite{EN}). Note 
that we obtain $\{\l,\,'\!\l^\vee\}$ for $\be=1,\al=0$.
Also, if $\vee$ is omitted, then $[\be,\al]$ in (\ref{xitwist})
must be changed to $[-\be,-\al]$ and the
resulting link becomes non-algebraic. 

To justify (\ref{xitwist}), let 
$\tilde{\xi}=\vph(\xi)=\vph\,\xi\vph$, which is the
conjugation of $\xi$ by the matrix
{\tiny $\begin{pmatrix} 0 & 1 \\ 1 & 0\\ \end{pmatrix}$}.
Then for any $P,Q\in \v^W$ and $P^\iota=\iota(P)$,
\begin{align}\label{xiPQ}
\xi(P)(Q)\!=\!\vph\bigl(\tilde{\xi}\si^{-1}(P^\iota)\bigr)
(Q),\ 
\{\xi(P)\,Q\}_{ev}\!=\!
\{\vph\bigl(\tilde{\xi}\si^{-1}(P^\iota)\bigr)Q\}_{ev},\\
=\{\tilde{\xi}\si^{-1}(P^\iota),Q\}_{ev}
=\{Q,\tilde{\xi}\si^{-1}(P^\iota)\!\Downarrow\}_{ev}=
\{\tilde{\xi}\si^{-1}(P^\iota)\!\Downarrow,Q\}_{ev},\notag
\end{align}
where the first column of $\tilde{\xi}\si^{-1}$ is 
$(\be,\al)^{tr}$. This gives the required.
\smallskip

Thus we can use Theorem \ref{STABILIZ}.
For instance, the super-duality from (\ref{iter-duality}) holds 
for such $\ \lxi\hat{\h}$. Upon the switch to the 
spherical $P_\la^\circ$, the product formula
(\ref{q-1-prod}) for the 
specializations at $q=1$ is also true (in the case of $A_n$),
where we add $[\be,\al]$ as the first vertex to (all paths in) 
$\,'\!\l$. Upon such a  modification of $\,'\!\l$, the
estimate (\ref{deg-a-j}) for deg$_a$ holds too.
Also, the positivity claim from Part $(ii)$
of Conjecture \ref{CONCONJ} can be extended to positive
 $\{\xi,\,\l,\,'\!\l^\vee\}$.
\smallskip

Let us combine (\ref{xiPQ}) with Theorem \ref{EVALTAUM}. 
Recall that $\dot{\tau}_-(P)=
\tau_-(P)\!\Downarrow\, =\bigl(\tau_-(P)\bigr)(1)\in \v^W$.
For any $P,Q\in \v^W$
and for $\ze\equal\tau_-^{-1}\tilde{\xi}\si^{-1}=$
{\tiny $\begin{pmatrix} \be & \ast \\ \al\!-\!\be & \ast\\
\end{pmatrix}$}, 
\begin{align}\label{pre-tau-min}
&\dot{\tau}_-(P^\iota Q)=
\tau_-\bigl(\vph\,\si^{-1}(P)\,Q\bigr)\!\Downarrow\,\,=\,
\bigl(\vph(\tau_-(P))\bigr)\bigl(\dot{\tau}_-(Q)\bigr),\\
\label{xitaum}
&\dot{\tau}_-\Bigl(\bigl(\ze(P^\iota)\bigr)\!\!\Downarrow
\,\dot{\tau}_-^{-1}(Q)\Bigr)\ =\ 
\vph\bigl(\vph\,\xi(P)\!\Downarrow\bigr)(Q),\
\notag\\
&\Bigl\{\dot{\tau}_-\Bigl(\bigl(\ze(P)\!\Downarrow\bigr)
\,\dot{\tau}_-^{-1}(Q^\iota)\Bigr)\Bigr\}_{ev}
=\{\xi(P)(Q)\}_{ev}\,.\
\end{align}

Therefore
 $\xi(P^\emptyset)(\iota(Q^\emptyset))$ serves as 
a {\em pre-polynomail\,} for the tree
\begin{align}\label{treeunion}
[1,1]\!\rightrightarrows 
\bigl(\, [\be,\al\!-\!\be]\,\rightdotsarrow '\!\l,\ 
[1,-1]\,\rightdotsarrow\l\, 
\bigr)\,, 
\ P^\emptyset\sim\, '\!\l,\, Q^\emptyset\sim\, \l,
\end{align}
where $\sim$ means that $P^\emptyset,Q^\emptyset$ are 
pre-polynomials for $\,'\!\l,\l$. However, it does not
coincide with the {\em standard pre-polynomial\,} associated 
with this tree in (\ref{jones-hat}) (and through this work).


\Yboxdim5pt
\subsection{\bf Some examples} 
The first will be for $\xi$
with the first column $(3,2)^{tr}$. We will write
$\xi=\ga_{3,2}$ then; $\tau_+\tau_-^2$ will do. We take
two unknots $\l,\,'\!\l$ 
colored by $\yng(2)\,,\,\yng(1,1,1)\,$ ; the second is
where $\xi$ is applied. 

\Yboxdim5pt
\begin{align}\label{T3-2-111-20K}
&\xi=\ga_{3,2},\ \ 
\l=\, \rightarrow\!\raisebox{1.5pt}{$\yng(2)\,$},\ \  
'\!\l=\, \rightarrow\!\yng(1,1,1)\, ,\  \
\lxi\hat{\h}{}^{min}_{\l,\,'\!\l}\,(q,t,a)=
\end{align}

\renewcommand{\baselinestretch}{0.5} 
\noindent
{\small
\(
1+\frac{a^4 q^7}{t^3}+q t+q t^2-t^3+q^2 t^3+q^3 t^3-q t^4+q^2 t^4
+q^3 t^4-q t^5+q^3 t^5+q^4 t^5-q^2 t^6+2 q^4 t^6-q^2 t^7-q^3 t^7
+q^4 t^7+q^5 t^7-q^3 t^8+q^5 t^8-q^3 t^9+2 q^5 t^9-q^4 t^{10}
+q^6 t^{10}-q^4 t^{11}+q^6 t^{11}-q^4 t^{12}+q^6 t^{12}-q^5 t^{15}
+q^7 t^{15}+a^3 \bigl(-q^4+q^6+q^7+\frac{q^4}{t^3}+\frac{q^5}{t^3}
+\frac{q^6}{t^2}+\frac{q^6}{t}+q^7 t+q^7 t^2-q^5 t^3+q^7 t^3\bigr)
+a^2 \bigl(-q^2+2 q^4+2 q^5+\frac{q^2}{t^3}+\frac{q^3}{t^2}
+\frac{q^4}{t^2}+\frac{q^3}{t}+\frac{q^4}{t}-q^3 t+q^4 t+2 q^5 t
+q^6 t-q^3 t^2+2 q^5 t^2+2 q^6 t^2-2 q^4 t^3+3 q^6 t^3-2 q^4 t^4
+2 q^6 t^4+q^7 t^4-q^4 t^5+q^6 t^5+q^7 t^5-q^5 t^6+2 q^7 t^6
-q^5 t^7
+q^7 t^7-q^5 t^8+q^7 t^8\bigr)+a \bigl(2 q^2+q^3+\frac{q}{t^2}
+\frac{q}{t}-q t+q^2 t+2 q^3 t+q^4 t-q t^2+2 q^3 t^2+2 q^4 t^2
-2 q^2 t^3+3 q^4 t^3+q^5 t^3-q^2 t^4-2 q^3 t^4+2 q^4 t^4+3 q^5 t^4
-2 q^3 t^5+3 q^5 t^5+q^6 t^5-q^3 t^6-q^4 t^6+2 q^5 t^6+2 q^6 t^6
-3 q^4 t^7+3 q^6 t^7-2 q^4 t^8+2 q^6 t^8-q^4 t^9+q^6 t^9+q^7 t^9
-q^5 t^{10}+q^7 t^{10}-q^5 t^{11}+q^7 t^{11}-q^5 t^{12}
+q^7 t^{12}\bigr).
\)
}
\renewcommand{\baselinestretch}{1.0} 
\smallskip

\Yboxdim5pt
The positivity of the series
$\lxi\hat{\h}{}^{min}_{\l,\,'\!\l}(q,t,a)/(1-t)$ holds.
This polynomial is super-dual to the polynomial 
$\lxi\hat{\h}^{min}$
for the pair of diagrams 
$\yng(1,1)\,,\,\yng(3)\,$. Using the evaluation formula at
$q=1$ for the latter 
and the $q\leftrightarrow 1/t$-duality,
we obtain that $\lxi\hat{\h}{}^{min}_{\l,\,'\!\l}(q,t=1,a)=
(1+a)(1+a+aq)^3$, which is true indeed for the polynomial
in (\ref{T3-2-111-20K}).
The value $\lxi\hat{\h}{}^{min}_{\l,\,'\!\l}(q=1,t,a)$
is $(1+a)$ times the evaluation at $q=1$ for the 
trefoil and $\om_3$; the latter is irreducible. 
\smallskip
 
The next example will be for
the uncolored trefoil taken as $\l$.
\Yboxdim5pt
\begin{align}\label{T3-2-10-11K}
&\xi=\ga_{3,2},\ \ 
\l=\, \l^{\circ\!\rightarrow,\, \yng(1)}_{\{3,2\}},\ \ 
'\!\l=\l^{\circ\!\rightarrow,\, \yng(2)}_{\{1,0\}},\ \  
\lxi\hat{\h}^{min}_{\l,\,'\!\l}\,(q,t,a)=
\end{align}

\renewcommand{\baselinestretch}{0.5} 
\noindent
{\small
\(
1-t+q^2 t+q^3 t+q^4 t+q^5 t-q^2 t^2-q^3 t^2+2 q^6 t^2+2 q^7 t^2
+q^8 t^2-q^4 t^3-q^5 t^3-q^6 t^3-q^7 t^3+q^8 t^3+3 q^9 t^3
+q^{10} t^3-q^6 t^4-q^7 t^4-q^8 t^4-2 q^9 t^4+q^{10} t^4
+3 q^{11} t^4+q^{12} t^4-q^8 t^5-q^9 t^5-q^{10} t^5-2 q^{11} t^5
+2 q^{12} t^5+3 q^{13} t^5-q^{10} t^6-q^{11} t^6-q^{12} t^6
+2 q^{14} t^6+q^{15} t^6-q^{12} t^7-2 q^{13} t^7+q^{14} t^7
+2 q^{15} t^7-2 q^{14} t^8+2 q^{16} t^8-q^{15} t^9+q^{17} t^9
+a^4 \bigl(q^{14}-q^{14} t+q^{16} t-q^{16} t^2+q^{18} t^2\bigr)
+a^3 \bigl(q^9+q^{10}+q^{11}+q^{12}-q^9 t-q^{10} t+q^{12} t
+2 q^{13} t+2 q^{14} t-q^{11} t^2-2 q^{12} t^2-q^{13} t^2
+2 q^{15} t^2+2 q^{16} t^2-q^{13} t^3-2 q^{14} t^3+q^{16} t^3
+q^{17} t^3+q^{18} t^3-q^{15} t^4-q^{16} t^4+q^{17} t^4+q^{18} t^4
-q^{16} t^5+q^{18} t^5\bigr)+a^2 \bigl(q^5+q^6+2 q^7+q^8+q^9-q^5 t
-q^6 t-q^7 t+q^8 t+3 q^9 t+4 q^{10} t+3 q^{11} t+q^{12} t-q^7 t^2
-2 q^8 t^2-3 q^9 t^2-2 q^{10} t^2+2 q^{11} t^2+4 q^{12} t^2
+4 q^{13} t^2+q^{14} t^2-q^9 t^3-2 q^{10} t^3-4 q^{11} t^3
-3 q^{12} t^3+2 q^{13} t^3+4 q^{14} t^3+3 q^{15} t^3+q^{16} t^3
-q^{11} t^4-2 q^{12} t^4-4 q^{13} t^4-q^{14} t^4+4 q^{15} t^4
+3 q^{16} t^4+q^{17} t^4-q^{13} t^5-3 q^{14} t^5-q^{15} t^5
+2 q^{16} t^5+2 q^{17} t^5+q^{18} t^5-2 q^{15} t^6-q^{16} t^6
+2 q^{17} t^6+q^{18} t^6-q^{16} t^7+q^{18} t^7\bigr)+a \bigl(q^2
+q^3+q^4+q^5-q^2 t-q^3 t+q^5 t+3 q^6 t+4 q^7 t+2 q^8 t+q^9 t
-q^4 t^2-2 q^5 t^2-2 q^6 t^2-2 q^7 t^2+2 q^8 t^2+5 q^9 t^2
+4 q^{10} t^2+2 q^{11} t^2-q^6 t^3-2 q^7 t^3-3 q^8 t^3-4 q^9 t^3
+5 q^{11} t^3+4 q^{12} t^3+2 q^{13} t^3-q^8 t^4-2 q^9 t^4
-3 q^{10} t^4-5 q^{11} t^4+q^{12} t^4+6 q^{13} t^4+3 q^{14} t^4
+q^{15} t^4-q^{10} t^5-2 q^{11} t^5-3 q^{12} t^5-3 q^{13} t^5
+3 q^{14} t^5+5 q^{15} t^5+q^{16} t^5-q^{12} t^6-3 q^{13} t^6
-2 q^{14} t^6+2 q^{15} t^6+3 q^{16} t^6+q^{17} t^6-2 q^{14} t^7
-2 q^{15} t^7+2 q^{16} t^7
+2 q^{17} t^7-q^{15} t^8-q^{16} t^8+q^{17} t^8+q^{18} t^8\bigr).
\)
}
\renewcommand{\baselinestretch}{1.0} 
\smallskip

\Yboxdim5pt
The super-duality is with  
the $\lxi\hat{\h}^{min}$\~polynomial
for the diagrams 
$\yng(1)\,,\yng(1,1)\,$. The positivity is for
$\lxi\hat{\h}{}^{min}_{\l,\,'\!\l}(q,t,a)/(1-t)$.
One has: $\lxi\hat{\h}{}^{min}_{\l,\,'\!\l}(1,t,a)=
(1+a)(1+a+t)^3$. 
Concerning the latter product,
the switch to polynomials $P_\la^\circ$ 
in Part $(iv)$ of Theorem \ref{MAINTHM} 
means the division of the latter product
by $(1+a)$, which corresponds
to the intersection of $\yng(1)\,,\yng(2)\,$. 
Since the tree $\l$ coincides with $\,\!\l$ extended by $[3,2]$,
and $\yng(2)=2\om_1$, we indeed arrive at $(1+a+t)^3$.
\smallskip

{\sf Examples of iterated type.} 
It will be $\lxi\h=
\{\hat{\xi}\bigl(\hat{\ga}(J_{\yng(1)\,})\!\!\Downarrow\bigr)\,
J^\circ_{\yng(1)\,}\}_{ev}$ 
for $\ga=\ga_{2,1}$ upon the hat-normalization. Recall that
$J_\la^\circ=P_\la^\circ=P_\la/P_\la(t^\rho)$.
\Yboxdim5pt
\begin{align}\label{T1-0-3-2-2-1K}
&\xi=\ga_{3,2},\ \ 
\l=\, \l^{\circ\!\rightarrow,\, \yng(1)}_{\{1,0\}},\ \ 
'\!\l=\l^{\circ\!\rightarrow,\, \yng(1)}_{\{2,1\}},\ \  
\lxi\hat{\h}^{min}_{\l,\,'\!\l}\,(q,t,a)=
\end{align}
\renewcommand{\baselinestretch}{0.5} 
\noindent
{\small
\(
1-t+q t+q^2 t+q^3 t+q^4 t-q t^2+2 q^4 t^2+2 q^5 t^2+q^6 t^2
-q^2 t^3-q^4 t^3+q^5 t^3+3 q^6 t^3+2 q^7 t^3-q^3 t^4-q^5 t^4
+3 q^7 t^4+3 q^8 t^4-q^4 t^5-q^6 t^5+3 q^8 t^5+2 q^9 t^5-q^5 t^6
-q^7 t^6+3 q^9 t^6+q^{10} t^6-q^6 t^7-q^8 t^7+q^9 t^7
+2 q^{10} t^7-q^7 t^8-q^9 t^8+2 q^{10} t^8+q^{11} t^8-q^8 t^9
+q^{11} t^9-q^9 t^{10}+q^{11} t^{10}-q^{10} t^{11}+q^{11} t^{11}
-q^{11} t^{12}+q^{12} t^{12}+a^4 \bigl(q^{10}+q^{11} t
+q^{12} t^2\bigr)+a^3 \bigl(q^6+q^7+q^8+q^9-q^6 t+q^7 t+2 q^8 t
+2 q^9 t+2 q^{10} t-q^7 t^2+2 q^9 t^2+3 q^{10} t^2+2 q^{11} t^2
-q^8 t^3+2 q^{10} t^3+2 q^{11} t^3+q^{12} t^3-q^9 t^4+2 q^{11} t^4
+q^{12} t^4-q^{10} t^5+q^{11} t^5+q^{12} t^5-q^{11} t^6
+q^{12} t^6\bigr)+a^2 \bigl(q^3+q^4+2 q^5+q^6+q^7-q^3 t+q^5 t
+4 q^6 t+4 q^7 t+3 q^8 t+q^9 t-q^4 t^2-q^5 t^2-q^6 t^2+4 q^7 t^2
+6 q^8 t^2+5 q^9 t^2+q^{10} t^2-q^5 t^3-q^6 t^3-2 q^7 t^3
+3 q^8 t^3+7 q^9 t^3+5 q^{10} t^3+q^{11} t^3-q^6 t^4-q^7 t^4
-2 q^8 t^4+3 q^9 t^4+6 q^{10} t^4+3 q^{11} t^4-q^7 t^5-q^8 t^5
-2 q^9 t^5+4 q^{10} t^5+4 q^{11} t^5+q^{12} t^5-q^8 t^6-q^9 t^6
-q^{10} t^6+4 q^{11} t^6+q^{12} t^6-q^9 t^7-q^{10} t^7+q^{11} t^7
+2 q^{12} t^7-q^{10} t^8+q^{12} t^8-q^{11} t^9+q^{12} t^9\bigr)
+a \bigl(q+q^2+q^3+q^4-q t+q^3 t+3 q^4 t+4 q^5 t+2 q^6 t+q^7 t
-q^2 t^2-q^3 t^2-q^4 t^2+2 q^5 t^2+6 q^6 t^2+5 q^7 t^2+2 q^8 t^2
-q^3 t^3-q^4 t^3-2 q^5 t^3+6 q^7 t^3+7 q^8 t^3+3 q^9 t^3-q^4 t^4
-q^5 t^4-2 q^6 t^4-q^7 t^4+6 q^8 t^4+7 q^9 t^4+2 q^{10} t^4
-q^5 t^5-q^6 t^5-2 q^7 t^5-q^8 t^5+6 q^9 t^5+5 q^{10} t^5
+q^{11} t^5-q^6 t^6-q^7 t^6-2 q^8 t^6+6 q^{10} t^6+2 q^{11} t^6
-q^7 t^7-q^8 t^7-2 q^9 t^7+2 q^{10} t^7+4 q^{11} t^7-q^8 t^8
-q^9 t^8-q^{10} t^8+3 q^{11} t^8+q^{12} t^8-q^9 t^9-q^{10} t^9
+q^{11} t^9+q^{12} t^9
-q^{10} t^{10}+q^{12} t^{10}-q^{11} t^{11}+q^{12} t^{11}\bigr).
\)
}
\renewcommand{\baselinestretch}{1.0} 
\smallskip

\Yboxdim7pt
This polynomial is self-dual. The positivity of
$\lxi\hat{\h}{}^{min}_{\l,\,'\!\l}(q,t,a)/(1-t)$ holds.
One has: $\lxi\hat{\h}{}^{min}_{\l,\,'\!\l}(q=1,t,a)=$
\renewcommand{\baselinestretch}{0.5} 
\noindent
{\small
\(
(1+a)\bigl(1+3 a+3 a^2+a^3+3 t+7 a t+5 a^2 t+a^3 t
+4 t^2+8 a t^2+5 a^2 t^2+a^3 t^2+4 t^3+8 a t^3+4 a^2 t^3+4 t^4
+6 a t^4+2 a^2 t^4+3 t^5+4 a t^5+a^2 t^5+2 t^6+2 a t^6+t^7
+a t^7+t^8\bigr),
\)
}
\renewcommand{\baselinestretch}{1.0} 
where the second factor is that for $C\!ab(13,2)T(3,2)$,
which appears when making $\l$ empty.

\Yboxdim5pt
Let us provide the identity from (\ref{xitwist}) in this
case. The polynomial
$\lxi\hat{\h}^{min}_{\l,\,'\!\l}\,(q,t,a)$ from
(\ref{T1-0-3-2-2-1K}) coincides with 

\begin{align}\label{T1-0-3-2-2-1KI}
\hat{\h}^{min}_{\l,\,'\!\tilde{\l}^\vee}\,(q,t,a) \for
\l=\, \l^{\circ\!\rightarrow\, \yng(1)}_{\{1,0\}},\ \ 
'\!\tilde{\l}=\l^{\circ\!\rightarrow\circ\!\rightarrow\,
\yng(1)}_{\{\{2,2\},\{3,1\}\}}.
\end{align}
\Yboxdim7pt
The latter corresponds to $C\!ab(13,2)T(2,3)$
for $[2,3]\!\rightarrow\![2,1]\!\rightarrow\!\yng(1)\,$.
Note $[2,3]$ here. The label $[\,3,2\,]$ would
result in a different 
$\hat{\h}^{min}$\~polynomial due
to applying $P_{\square\,}(Y)$ to the corresponding 
pre-polynomial.
\smallskip

\Yboxdim5pt
Making now $\xi=\ga_{1,1}$, say taking
$\xi=\tau_-$, we obtain:
\Yboxdim5pt
\begin{align}\label{T1-0-2-1K}
&\xi=\ga_{1,1},\ \ 
\l=\, \l^{\circ\!\rightarrow,\, \yng(1)}_{\{1,0\}},\ \ 
'\!\l=\l^{\circ\!\rightarrow,\, \yng(1)}_{\{2,1\}},\ \  
\lxi\hat{\h}^{min}_{\l,\,'\!\l}\,(q,t,a)=
\hbox{\small $1+a^2 q^3$}
\end{align}
\renewcommand{\baselinestretch}{0.5} 
\noindent
{\small
\(
-t+q t+q^2 t-q t^2+q^2 t^2-q^2 t^3+q^3 t^3
+a \bigl(q+q^2-q t+q^2 t+q^3 t-q^2 t^2+q^3 t^2\bigr),
\)
}
\renewcommand{\baselinestretch}{1.0} 

\Yboxdim5pt
\noindent
which is self-dual. The positivity of
$\lxi\hat{\h}{}^{min}_{\l,\,'\!\l}(q,t,a)/(1-t)$ holds.
One has: $\lxi\hat{\h}{}^{min}_{\l,\,'\!\l}(1,t,a)=
(1+a)(1+a+t)$. Taking empty $\l$ here will result in
$C\!ab(3,2)T(1,1)$, i.e. in the trefoil; this explains the 
factor $(1+a+t)$. 
\smallskip

This polynomial can be obtained by our usual construction.
Namely, we have the following DAHA identities:
\begin{align*}
{}^{\tau_-}\hat{\h}^{min}_{[1,0]\!\rightarrow 
\yng(1)\,,\ [2,1]\!\rightarrow \yng(1)\, }=
\hat{\h}^{min}_{[2,3]\!\rightarrow 
\yng(1)\,,\ [1,0]\!\rightarrow \yng(1)\,^\vee }=
\hat{\h}^{min}_{[1,1]\!\rightarrow[2,1]\!\rightarrow 
\yng(1)\,}\,.
\end{align*}
Recall that the first polynomial
is $\{\tau_-\bigl(\tau_+\tau_-(J_{\yng(1)\,})\!\!
\Downarrow\bigr)J^\circ_{\yng(1)\,}\}_{ev}$ upon 
the hat-normalization.
Up to $q^\bullet t^\bullet$, the remaining two correspond to:
$$
\{J_{\yng(1)\,}(Y)\bigl(\tau_-\tau_+\tau_-
(J^\circ_{\yng(1)})\!\Downarrow\bigr)\}_{ev}\, 
=\,\{J_{\yng(1)\,}(Y)\Bigl(\tau_-\bigl(
\tau_+\tau_-(J^\circ_{\yng(1)})\!\!\Downarrow\bigr)
\!\!\Downarrow\Bigr)\}_{ev}\,,
$$
which are identical due the commutativity of 
$\tau_-$ with $\Downarrow$. Their coincidence with
the first one is essentially the verification of the
fact that our twisted construction depends only
on the topological type of the corresponding link.  
\smallskip

\subsection{\bf Toward the Skein}\label{sec:skein}
The generalized twisting can be used for
the following topological characterization of the 
pre-polynomials. 
Let $P^\emptyset$ be a pre-polynomial, where $\emptyset$
means that it is
obtained without the division by the evaluations of the
$J$\~polynomials involved. We mostly need $P=P^{min}$
in this work, which is $P^\emptyset$ divided by the total
$LC\!M$ of these evaluations. The corresponding superpolynomial 
is $\h^{min}\!\!=\{P^{min}\}_{ev}$; we will mostly drop the
hat-normalization in this section. Due to this normalization,
$P^\emptyset$ and $P^{min}$ are actually needed only up 
to $q,t$\~monomial factors in the rest of this work.

\begin{theorem}\label{GENUNION}
(i) Let $P^{\emptyset}$ be the standard pre-polynomial
for a tree $\l$ from (\ref{jones-hat})
and $\xi\in PSL_2(\Z)$ be as above. Using the
non-degeneracy of $\{\,,\,\}_{ev}$,
the pre-polynomial\, $\xi(P^\emptyset)\!\Downarrow$
for the tree $[\al,\be]\rightarrow \l$  can be uniquely 
determined via the coinvariants
$\{\xi(P^\emptyset)(J_\la^\iota)\}_{ev}$ 
for all diagrams $\la$. The latter coinvariants are
the superpolynomials  
$\h^\emptyset_{\mathcal{M}}$ for the trees
$\ \mathcal{M}\,=\,[1,1]\!\rightrightarrows
\bigl(\, [\be,\al\!-\!\be]\,\rightarrow \l,\ 
[1,-1]\,\rightarrow\la\, \bigr)\,.$ 

(ii) Therefore all pre-polynomials (for
any trees and colors) can be uniquely recovered 
if all $\h$\~invariants are known. Topologically, this 
gives that the standard ore-polynomials (\ref{jones-hat})
are invariants of the 
corresponding links 
considered (naturally) in the solid torus. This means that the
symmetries from Sections \ref{sec:ITER-KNOTS}, 
\ref{sec:knots-links} hold unless they involve the first
vertex $[\rr_1,\ss_1]$. For instance,
(\ref{iter-symx}) holds, but (\ref{iter-sym}) does not. \sq
\end{theorem}

Note that the class of pre-polynomials is closed
with respect to the (usual) multiplication, where the 
renormalization (the division by the $LC\!M$ of all
evaluations)
is necessary if $\h^{min}$ are considered. Indeed,
the pure product
corresponds to the union of the trees, which can be considered
as one tree by connecting them to the additional 
(initial) vertex $[1,0]$. Generally, 
the reduction to a single
tree is not always reasonable. For instance,
this may result in non-positive labels for algebraic 
pairs of trees $\{\l,\,'\!\l^\vee\}$, as it is clear from 
(\ref{treeunion}).

\smallskip

The final step here would be a similar
topological understanding of the {\em knot operators\,}
themselves in terms of the topological
interpretation of their matrix elements. These operators are 
$K=\xi(P)$ (and then by induction) in the notation from 
Theorem \ref{GENUNION}; the matrix elements for such $K$ are 
$\{J_\mu,K(J_\la^\iota)\}_{ev}$ for any $\la,\mu$. Similar to
the pre-polynomials, the knot operators can be expected 
invariants of the links in the solid torus, but now we have 
to remove the middle circle too. Equivalently, they correspond
to links  in the torus times an interval, $T\times I$. 

Indeed, considering  $\{J_\mu, K (J_\la^\iota)\}_{ev}$ with
``free ends" $\la,\mu$
means topologically that the corresponding trees 
are subtrees of undermined larger ones, 
extending them in {\em both\,}
directions (below $i=0$ and beyond $i=\ell^j$). Accordingly,
the twisted unions 
can be applied in both directions of the corresponding links.  
This requires fixing the
meridian as for pre-polynomials and also the middle circle.
This topological restriction
gives that the symmetries from (\ref{iter-sym}) and 
(\ref{iter-symx}) must not generally
hold (now both) for the knot operators.
\smallskip

This approach can be expected to provide 
an identification
of the {\em spherical DAHA\,} with 
the {\em toric $q,t$\~skein algebra\,}, the Skein,
which is for $T\times I$. This will hopefully 
help to obtain $\h$\~polynomials for 
arbitrary links (not only cables of $\unknot\,$),
though this problem can be approached independently.

Our way to understand the Skein via twisted unions 
(following Theorem \ref{GENUNION}) is {\em global\,} (not via
{\em local\,} skein relations), but there is a relatively
direct connection with the relations of the
Elliptic Hall Algebra from \cite{SV}. 
Such an approach to the Skein is not at all restricted to $A_n$ 
and can be readily extended to any root systems. The  
$a$\~stabilization for exceptional root systems will/may be lost 
(though we found some in \cite{ChE}), but anything else is 
expected to hold. The {\em refined\,}  WRT-invariants (for
any Quantum Groups) will play the role of DAHA-superpolynomials. 
\smallskip

This program is partially connected with paper \cite{ChE}
devoted to the {\em $q,t$\~composite theory\,} in  
\,annulus$\,\times I$. The Macdonald polynomials for 
{\em composite partitions\,}  $[\la,\mu]\equal
\iota(\la)+\mu$  
and their $a$\~stabilization were used, instead of 
the products $P_{\la}P_{\mu}$ here. This construction 
is actually for links too
(for double torus knots), but the corresponding superpolynomials 
are (naturally) smaller than those in the present work. For 
instance, the uncolored composite case is for the weight 
$\om_1+\om_n$ upon the stabilization at $n\to \infty$. We do
$P_{\Box}^2$ or  $P_{\Box}P_{\Box}^\iota$ here.  
\smallskip

{\sf Seifert 3-folds.}
It seems that the expected interpretation of the 
spherical DAHA as Skein($T\times I$) will require 
{\em Seifert manifolds}. 

The twisted unions (with $\vee$) can be 
naturally considered as links in the lens spaces.
We place two links for $\l$ and $\,'\l$ in the standard
horizontal solid torus and its vertical complement in $\S^3$ 
and then perform the corresponding twist at their common
boundary, the 2-dimensional torus. The lens spaces and
Seifert manifolds seem really necessary if 
arbitrary sequences
$\vec{\xi}=(\xi_i,1\le i\le m)$ of
matrices from $PGL_2(\Z)$ are considered for $m>1$. 
Given a sequence of pre-polynomials
$\vec{P^\emptyset}=(P_i^\emptyset)$, the corresponding
twisted knot operator and pre-polynomial are 
\begin{align}\label{knotoper}
\k^{\,\emptyset}=
\xi_1(P^\emptyset_1)\cdots 
\xi_m(P^\emptyset_m),\ \,
\p^{\,\emptyset}= \bigl(\xi_1(P^\emptyset_1)\cdots 
\xi_m(P^\emptyset_m)\bigr)\!\Downarrow.
\end{align}
Accordingly, we set
$\h^\emptyset={}^{\vec{\xi}}\h^{{}^\emptyset}_{\vec{P}}=
\{\p^{\emptyset}\}_{ev}$
and define $\k^{min},\p^{min}, \h^{min}$ upon the division by
the $LC\!M$ of all $J$\~polynomials involved.

The theory of Khovanov-Rozansky polynomials of the 
links in the lens spaces and Seifert $3$\~folds
is not developed. However see \cite{Ste1} and
references there for the mathematics and physics approaches
to the HOMFLY-PT theory there. The generalized 
{\em ORS polynomials} can be expected then
for the germs of curve singularities in the toric surfaces 
associated with these manifolds (and their plumbing). 

Using Seifert manifolds is highly desirable
anyway, especially if $\S^3$ is really
insufficient for a topological interpretation 
of the knot operators $\k$,   
pre-polynomials $\p$ for $m>1$ and 
$\h$\~invariants for $m>2$.

Recall that the latter can be reduced to the invariants 
in $\S^3$ for $m=2$ due to (\ref{xitwist}). Namely 
$\{\hat{\xi}(Q),\hat{\xi}'(P)\}_{ev}$ for 
pre-polynomials $P,Q$ and $\xi,\xi'\in PSL(2,\Z)$
are associated with proper links in $\S^3$. We
note that this bilinear form (in $P,Q$) is symmetric when 
$\xi'=\vph(\xi)$, which intersects
the {\em higher-level
DAHA coinvariants\,} from \cite{ChM}.
\smallskip

The DAHA theory of $\k,\p,\h$ from
(\ref{knotoper}) is similar to that for $m=2$.
They depend only on the
first columns $(\al_i,\be_i)^{tr}$ of $\,\xi_i$ 
and the hat-normalized
$\hat{\h}^{min}\,$ are $a,q,t^{\pm1}$\~polynomials 
satisfying the 
super-duality. Also, if $P_i$ are for positive trees
such that $\ss_1>\rr_1>0$,
then the conditions $\al_i>\be_i>0$ seem sufficient 
(not always necessary) for the positivity 
of $\hat{\h}^{min}/((1\!-\!q)(1\!-\!t))^M$
for $M$ large enough.

Conceptually, such $\hat{\h}^{min}\,$ is the
{\em dynamical\,} $m$\~point correlation function 
in the DAHA theory; dynamical, because the parameters are weights. 
The theory of the corresponding KZ-type difference equations, 
closely related to the so-called {\em $A$\~polynomials\,}, is in 
progress. An interpretation of generalized
$\hat{\h}^{min}\,$ (for any $m$) within the $\S^3$\~theory is
not impossible, as well as that for the toric $q,t$\~skein 
algebra, but it seems 
that Seifert manifolds will be needed here. Anyway, the
interpretation of the 
knot operators from (\ref{knotoper}) and their iterations
topologically and physically is a challenge for the new theory 
we present in this work. 
\medskip

\comment{ 
Using (\ref{treeunion}), we can continue now by induction
For instance,
the matrix element  $\{J_\la,\,\xi(J_\mu)J_\nu\}_{ev}$
$=\{\xi(J_\mu)\bigl(J_\la J_\nu\bigr)\}_{ev}$ of $\k=\xi(J_\mu)$
for the (pure) Macdonald $J$\~polynomials and $\xi=
\hbox{\tiny $\begin{pmatrix} \al & \ast \\ \be & \ast\\ 
\end{pmatrix}$}\in PSL_2(\Z)$
is the evaluation of the pre-polynomial for the following tree: 
\begin{align}\label{tree3union}
[1,1]\rightthreearrow\bigl([1,-1]\rightarrow \la\,,\ 
[\be,\al\!-\!\be]\,\rightarrow \mu\,,\ 
[1,-1]\,\rightarrow \nu\,\bigr).
\end{align}

Generally, we obtain that the class of
{\em \,pre-polynomials\,} for arbitrary graphs associated with
torus iterated links is closed with respect to
the following operations.

\begin{definition}
(i) Let $\vec{\xi}=(\xi_i,1\le i\le m)$ be a sequence of
matrices from $PGL_2(\Z)$ and 
$\vec{P^\emptyset}=(P_i^\emptyset)$ be any sequence of 
pre-polynomials, where by $\emptyset$, we mean that
no division by the evaluations of the $J$\~polynomials is
performed. Then we set:
\begin{align}\label{prepol}
\p^{\,\emptyset}\equal \bigl(\xi_1(P^\emptyset_1)
\xi_2(P^\emptyset_2)\cdots 
\xi_m(P^\emptyset_m)\bigr)\!\Downarrow.
\end{align}
Accordingly,
${}^{\vec{\xi}}\hat{\h}^{{}^\emptyset}_{\vec{P}}$ 
will be the hat-normalization
of \, $\{\p^{\,\emptyset}\}_{ev}$. Then 
$\p^{min}$ is  $\p^{\,\emptyset}$
upon the division by the $LC\!M$ of all $J$\~polynomials 
involved and we define 
${}^{\vec{\xi}}\hat{\h}^{{}^{min}}_{\p}$ as the
hat-normalization of $\,\{\p^{min}\}_{ev}$. We also
call $\p_\la=\{J_\la,\p\}_{ev},\,
\p^\vee_\la=\{\iota(J_\la),\p\}_{ev}$ the components of $\p$
(with respect to the pairing $\{\,,\,\}_{ev}$).

(ii) Furthermore, the operator
\begin{align}\label{knotoper}
{}^{\xi\!}\vec{K}^{{}^\bullet}_{\vec{P}}=
\xi_1(P^\bullet_1)
\xi_2(P^\bullet_2)\cdots 
\xi_m(P^\bullet_m)
\end{align}
will be called
the knot operator corresponding to $\p$, where 
{\small $\bullet$} is {\small $\emptyset$} or $min$.
It belongs to the spherical DAHA, which is
 $\, \s\!\HH\equal
\{H\in \HH\,\mid\, T_i H=HT_i \for i=1,2,\ldots, n\}$.
Finally, $\k_{\la,\mu}\equal
\{J_\la, \vec{K}(
J_\mu)\}_{ev}^{\bullet}$ and  $\k^\vee_{\la,\mu}\equal
\{\iota(J_\la), \vec{K}(
J_\mu)\}_{ev}^{\bullet}$

for Young diagrams $\la,\mu$ will be called 
matrix elements of $\k=
{}^{\vec{\xi}\!}\vec{K}^{{}^\bullet}_{\vec{P}},$ 
where {\small $\bullet$} is
{\small $\emptyset$}  or $min$; the latter means 
the division by the total $LC\!M$ of evaluations of all 
$J$\~polynomials involved, including $J_\la,J_\mu$.
\sq
\end{definition}

\begin{theorem}\label{GENUNION}
(i) The polynomial $\p^{\,min}$ defined above depends 
only on the
first columns $(\al_i,\be_i)^{tr}$ of $\,\xi_i$ 
(and the initial pre-polynomials). It can be obtained
as a pre-polynomial for one single tree $\mathscr{M}$ constructed
in terms of the corresponding graphs $\l_i$ for $P_i$ and
$(\al_i,\be_i)$ and depends only on the isotopic type of the
corresponding link (naturally considered) inside the solid torus.
Accordingly, ${}^{\vec{\xi}}\hat{\h}^{min}(\vec{P}\,)$ 
are $a,q,t^{\pm1}$\~polynomials satisfying the super-duality.
Also, if $P_i$ are for positive trees
such that $\ss_1>\rr_1>0$,
then the conditions $\al_i>\be_i>0$ are sufficient 
(not always necessary) for the resulting
link to be algebraic. 

(ii) The operator $\k^{min}=
{}^{\vec{\xi}\!}\vec{K}^{{}^{min}}_{\vec{P}}$ above is
uniquely determined by its matrix elements  $\k_{\la,\mu}$
(assuming that they are for all $\la,\mu$), which is immediate
since the evaluation pairing $\{\cdot,\cdot\}_{em}$ is 
non-degenerate (for generic $q,t$). These matrix elements
are $\hat{\h}$\~polynomials for the corresponding trees,
naturally generalizing (\ref{tree3union}), so they are 
isotopic invariants of the corresponding links. Therefore
$\vec{K}$ depends only on the tree $\l$ associated with 
the pre-polynomial $\p=\vec{K}\!\Downarrow$ up to the symmetries
from (\ref{iter-triv}) and (\ref{iter-s1}). Topologically,
it depends only on the isotopy class of this link considered 
inside the solid torus without the middle circle (parallel), which 
corresponds to adding the meridian and the middle circle
(two unknots) colored by $\la,\mu$.
\sq
\end{theorem}
} 

\appendix
\setcounter{equation}{0}
\section{\sc Links and splice diagrams}
\label{TopologySection}

In this appendix, we remind the main topological constructions 
from \cite{EN}, especially the definition and basic properties
of the splice diagrams. We reduce the generality and always 
assume that the links are in $\mathbb{S}^3$, not in an arbitrary 
homology 3-sphere.

The key operation there is {\em splicing\,}; it provides a large 
family of links. Other operations can be mainly considered 
as its special cases. This includes the cabling, unions 
and twisted unions, which play the major role in the DAHA
approach to torus iterated links.  

\subsection{\bf Links, cables and splices}\label{LinkSection}
A \textit{link} in $\mathbb{S}^3$ is a disjoint union  
$\{ S_i| i\in I \}$  of oriented simple closed (smooth) 
curves $S_i \colon \mathbb{S}^1\hookrightarrow \mathbb{S}^3$ 
called {\em components}. It will be denoted by $\mathbf{L} = 
\left(\mathbb{S}^3, \bigcup_{i\in I} S_i\right)$. Since we deal 
only with the
links in $\mathbb{S}^3$, we mostly omit $\mathbb{S}^3$ in this 
notation and write $\mathbf{L} = \bigcup_{i\in I} S_i$. A 
\textit{knot} is a link with a single component. We consider 
links  up to isotopies (smooth homotopies), ignoring the order of 
components. However the orientation of the components will
matter. \textit{The unknot} will be denoted by  $\unknot\,$,
the border of a (standard) disk $\mathbb{D}^2\subset \mathbb{S}^3$.

\medskip

For each component $S_i$ of a link $\mathbf{L} = 
\bigcup_{i\in I} S_i$, there is a system of
solid tori $\mathcal{N}(S_i)\subset 
\mathbb{S}^3$ 
such that $S_i\subset\mathcal{N}(S_i)$, 
$\mathcal{N}(S_i)\cap S_j = \emptyset\:(j\neq i)$ and 
there is 
a homeomorphism $\mathcal{N}(S_i) \to 
\mathbb{D}^2\times \mathbb{S}^1$ for each $i$
mapping $S_i$ to ``the middle circle'' 
$\{a\}\times \mathbb{S}^1 \subset \mathbb{D}^2\times 
\mathbb{S}^1$, $a\in \mathbb{D}^2$. Such 
$\mathcal{N}(S_i)$ is called a \textit{tubular neighborhood} 
of the component $S_i$.

Since $\partial \mathcal{N}(S_i)$ is homeomorphic to the 
torus $\mathbb{S}^1\times \mathbb{S}^1$, the first
homology is $H_1(\partial \mathcal{N}(S_i), 
\mathbb{Z}) = \mathbb{Z}^2$ and  there is a topologically 
distinguished basis there. 
Let $\mathit{l}(\cdot,\cdot)$ be the linking 
number. There is a single pair 
$(\mathsf{M},\mathsf{L})$, $\mathsf{M},\mathsf{L}
\in H_1(\partial \mathcal{N}(S_i), \mathbb{Z})$ such that\, $i$) 
$\mathsf{M}$ is contractible in $\mathcal{N}(S_i)$,\, $ii$) 
$\mathit{l}(S_i,\mathsf{M})=1$,\, $iii$) $\mathsf{L}$ is 
homotopically 
equivalent to $S_i$ in $\mathcal{N}(S_i)$, and $iv$) 
$\mathit{l}(S_i,\mathsf{L})=0$. We call $\mathsf{M}$ the 
{\em meridian\,} and $\mathsf{L}$ the {\em longitude\,}.

\smallskip

A \textit{framed link} is a link with a choice 
 of an element $\mathsf{F}_i\in H_1(\partial 
\mathcal{N}(S_i), \mathbb{Z})$ for each component 
$S_i$ such that $\mathsf{F}_i$ is 
homotopically equivalent $S_i$ in $\mathcal{N}(S_i)$. This 
element is called \textit{framing}. One can always equip a 
component with the longitude as its framing and we call such
a choice \textit{topological framing}. However in the
DAHA-approach, a different framing is generally needed 
(see below).

Given two framings $\mathsf{F}_i$ and 
$\tilde{\mathsf{F}}_i$ of the same component $S_i$,
one can find $k\in \mathbb{Z}$ such that 
$\mathsf{F}_i-\tilde{\mathsf{F}}_i = k \mathsf{M}_i$, where 
$\mathsf{M}_i$ is the meridian of $S_i$. Let us call $k$ 
the \textit{difference} of framings $\mathsf{F}_i$ and 
$\tilde{\mathsf{F}}_i$.

\medskip

We now briefly present the operations on links.

{\sf Erasing components.}
Given a link $\mathbf{L} = \bigcup_{i\in I} S_i$ and a 
subset of its components $J\subset I$, one can consider a 
new link $\bigcup_{i\in I\setminus J} S_i$.

{\sf Orientation reversion.}
For a link $\mathbf{L} = \bigcup_{i\in I} S_i$, one 
can reverse the orientation for a subset $J\subset I$ 
of components. This link will be denoted by
$\bigcup_{i\in I\setminus J} 
S_i \cup \bigcup_{j\in J} (S_j^\vee)$.
\smallskip

We will need the following general operation on 
manifolds. Let us define 
the \textit{connected sum} $\Sigma\, \sharp\, \Sigma'$
of two (connected) manifolds $\Sigma$ and $\Sigma'$
of the same dimension. We select 
two open balls $U\subset\Sigma$ and $U'\subset\Sigma'$ and then
paste 
$\Sigma\setminus U$ and $\Sigma'\setminus U'$ along 
$\partial U$ and $\partial U'$. This  
does not depend (up to a homeomorphism) on the choice of 
$U$ and $U'$. Note that 
$\mathbb{S}^3 \sharp\, \mathbb{S}^3 = 
\mathbb{S}^3$.
\smallskip

{\sf Disjoint sum.}
For two links $\mathbf{L} = \left(\mathbb{S}^3, 
\bigcup_{i\in I} 
S_i\right)$ and $\mathbf{L}' = \left(\mathbb{S}^3, 
\bigcup_{j\in J} 
S'_j\right)$, their \textit{disjoint sum} is the link
\begin{equation*}
\mathbf{L} + \mathbf{L}' = \left( \mathbb{S}^3 \sharp\, 
\mathbb{S}^3,\; 
\bigcup_{i\in I} S_i \cup \bigcup_{j\in J} S'_j \right).
\end{equation*}

The components $S_i$ and $S'_j$ here are the images of the 
corresponding components of the original links under 
the natural 
maps $\mathbb{S}^3\setminus U \rightarrow \mathbb{S}^3 \sharp\, 
\mathbb{S}^3$ and 
$\mathbb{S}^3\setminus U' \rightarrow \mathbb{S}^3 \sharp \,
\mathbb{S}^3$ from the 
definition of $\mathbb{S}^3 \sharp\, \mathbb{S}^3$. The open
balls $U$ and $U'$ must 
not intersect the components of $\mathbf{L}$ and $\mathbf{L}'$;
$\mathbf{L} + \mathbf{L}'$ does not depend on choice 
of these domains.

{\sf Connected sum.}
Now let us pick two components  $S_{i_0}$, $S'_{j_0}$ in
 $\mathbf{L}$, $\mathbf{L}'$ ($i_0\in I$, $j_0\in J$) and 
define the \textit{connected sum} of links 
$\mathbf{L}$ and $\mathbf{L}'$ \textit{along} $S_{i_0}$ and 
$S'_{j_0}:$ 
\begin{equation}\label{Connesum}
\mathbf{L} \sharp\, \mathbf{L}'(S_{i_0},\,S'_{j_0}) = 
\Bigl( \mathbb{S}^3 \sharp\, \mathbb{S}^3,\; \left( S_{i_0}
\sharp\, S'_{j_0} \right) 
\cup \bigcup_{i\in I\setminus\{i_0\}} S_i \cup 
\bigcup_{j\in J\setminus\{j_0\}} S'_j \Bigr).
\end{equation}

We identify $S_i$ and $S'_j$ for $i\ne i_0$, $j\ne j_0$ 
with the images of the corresponding components of the original 
links under the natural maps $\mathbb{S}^3\setminus U 
\hookrightarrow \mathbb{S}^3 
\sharp\, \mathbb{S}^3$ and $\mathbb{S}^3\setminus U' 
\hookrightarrow \mathbb{S}^3 
\sharp\, \mathbb{S}^3$. The definition of $\mathbb{S}^3 
\sharp\, \mathbb{S}^3$ is above.
Let us pick the open balls $U$ and $U'$ intersecting
only with\, $S_i,S'_j$ and such that 
the intervals $U\cap S_{i_0} = 
\mathcal{I}$ and $U'\cap S'_{j_0} = \mathcal{I}'$ are non-empty. 
Then we glue 
$S_{i_0}\setminus U$ 
to $S'_{i_0}\setminus U'$ with respect to the orientation 
of $S_{i_0}$ and $S'_{j_0}$ in $\mathbb{S}^3 \sharp \,
\mathbb{S}^3$. This gives $S_{i_0}\sharp\, S'_{j_0}$
needed in  (\ref{Connesum}). 
One can check that $\mathbf{L} \sharp\, \mathbf{L}'(S_{i_0},
\,S'_{j_0})$ does not depend on the choice of $U$ and $U'$; 
however it of course depends on the choice of components  
$S_{i_0}$ and $S'_{j_0}$. 

Note that both operations, $\mathbf{L}+\mathbf{L}'$ and
$\mathbf{L} \sharp\, \mathbf{L}'$, result in links in 
$\mathbb{S}^3$ since 
$\mathbb{S}^3 \sharp\, \mathbb{S}^3 = \mathbb{S}^3$.
The next operation is the key in DAHA theory.
\smallskip

{\sf Cabling.}
For coprime $k,l \in \mathbb{Z}$, there exists a unique 
up to homotopy oriented simple curve 
$C\!ab_{k,l}(S_i)\subset \partial \mathcal{N}(S_i)$ such that 
it is  
homotopy equivalent to $k\mathsf{M} + l\mathsf{L}$. 
Here $k,l \in \mathbb{Z}$ can be arbitrary,
then  $C\!ab_{k,l}(S_i)\subset 
\partial \mathcal{N}(S_i)$ will be a union (unique 
up to homotopy) of $d= \hbox{gcd}(k,l)$ non-intersecting oriented 
connected closed components; each of them is
homology equivalent to 
$(k/d)\mathsf{M} + (l/d)\mathsf{L}.$

The \textit{cable} of type $(k,l)$, $k,l\in\mathbb{Z}$ 
of a link $\mathbf{L} = \bigcup_{i\in I} S_i$ along the
component $S_{i_0}$ is the link 
\begin{equation*}
C\!ab^{S_{i_0}}_{k,l}\mathbf{L} = \;C\!ab_{k,l}(S_i) 
\cup \bigcup_{i\in I\setminus\{i_0\} } S_i.
\end{equation*}

The knot $C\!ab_{k,l}(\,\unknot\,)$ is the
\textit{torus knot} $(k,l)$ denoted by $T(k,l)$
(by $T(l,k)$ in the paper).

The $(k,1)$-cabling preserves a link (up to isotopy)
and $(k,-1)$-cabling
reverses the orientation of the cabled component:
\begin{equation*}
C\!ab^{S_{i_0}}_{k,1}\mathbf{L} = \mathbf{L}, \  
C\!ab^{S_{i_0}}_{k,-1} \Bigl(\bigcup_{i\in I} S_i\Bigr) = 
S_{i_0}^\vee \cup \bigcup_{i\in I\setminus \{ i_0\} } S_i.
\end{equation*}

The cabling operations along different components commute 
with each other because different components have 
non-intersecting tubular neighborhoods.
\smallskip

{\sf Splice.}
This is actually the most general
operation on links; all previous ones
can be obtained as its proper specializations.

We take two links $\mathbf{L} = \bigcup_{i\in I} S_i$ and 
$\mathbf{L}' = \bigcup_{j\in J} S'_j$ and fix one component
in each: $S_{i_0}$ and $S'_{j_0}$, where $i_0\in I$ and 
$j_0\in J$. \textit{The splice} of the links $\mathbf{L}$ 
and $\mathbf{L}'$ \textit{along} $S_{i_0}$ and $S'_{j_0}$,
denoted by
\begin{equation*}
\mathbf{L} \stackrel{\underline{S_{i_0} \; S'_{j_0}} }
{\phantom{\longrightarrow}} \mathbf{L}' = 
\left( \Sigma,\; \bigcup_{i\in I\setminus \{i_0\} } S_i 
\cup \bigcup_{j\in J\setminus \{j_0\} } S'_j \right),
\end{equation*}
is as follows. We
glue together  $\mathbb{S}^3 \setminus 
\mathcal{N}(S_{i_0}) $ 
and $\mathbb{S}^3 \setminus \mathcal{N}(S'_{j_0}) $ along 
$\partial \mathcal{N}(S_{i_0})$ and $\partial 
\mathcal{N}(S'_{j_0})$, where the \textit{meridian} 
of $\partial \mathcal{N}(S_{i_0})$ is identified
with the \textit{longitude} 
of $\partial \mathcal{N}(S'_{j_0})$ and \textit{vice versa}. 
As above,
$S_i$ and $S'_j$ ($i\ne i_0$, $j\ne j_0$) are considered as
the images of the corresponding components under the inclusions 
$\mathbb{S}^3 \setminus \mathcal{N}(S_{i_0}) \hookrightarrow 
\Sigma$ 
and $\mathbb{S}^3 \setminus \mathcal{N}(S'_{j_0}) 
\hookrightarrow \Sigma$. We will need only cases when 
$\Sigma$ is diffeomorphic to $\mathbb{S}^3$, so the splice
(the union of $S_i,S_j'$ above) will be a link in $\S^3$. 
\medskip

This operation is 
$\mathbf{L}\leftrightarrow\mathbf{L}'$\~symmetric. 
We will need its non-symmetric generalization. 
It will depend now on the framing $\mathsf{F}_{i_0}$
of $S_{i_0}$. The   \textit{splice of $\mathbf{L}$ 
and $\mathbf{L}'$ along $S_{i_0}$ and $S'_{j_0}$ deformed by 
framing $\mathsf{F}_{i_0}$} is as follows. Now only the 
longitude of $S'_{j_0}$ will be identified with the meridian of 
$S_{i_0}$. The meridian of $S'_{j_0}$ will not be identified 
with the corresponding longitude; it will be now identified with 
the framing $\mathsf{F}_{i_0}$. We set:
\begin{equation}\label{deformfram}
\mathbf{L} \stackrel{\lefteqn{
\underrightarrow{\phantom{S_{i_0\;}}}}
\underline{S_{i_0} \; S'_{j_0}} }{} \mathbf{L}'.
\end{equation}

Here $\mathsf{F}_{i_0}$ is determined from the context.
If $\mathsf{F}_{i_0}$ is the topological framing, then the 
deformed splice coincides with the standard one. The
difference can be only because of non-trivial choices of 
the framing. Note that the splice operations (deformed or 
non-deformed) commute with each other.
\smallskip

\subsection{\bf Splice diagrams}
Splice diagrams provide a convenient way to deal with links 
made by splices of ``canonical'' links in Seifert manifolds
(here in $\mathbb{S}^3$).

{\sf Seifert fibrations.}
A \textit{standard fibered solid torus} corresponding to a pair 
of coprime integers $(k,l)$ ($k>0$) is the topological space 
$\mathbb{D}^2\times[0,1]$ with its (two)
borders identified by relations 
$(x,1)\sim(\rho_{l/k}x,1)$, where $\rho_{\alpha}$ is the 
rotation by angle $2\pi \alpha$ in $\mathbb{D}^2$ around the 
center 
$O\in \mathbb{D}^2$. The fibration is inherited from 
$\mathbb{D}^2\times[0,1]$ 
in the following way. For each $x\in \mathbb{D}^2\setminus O$, 
the segments $\{\rho_{i/k}x\} \times [0,1]$ form a circle. 
These circles will be called \textit{regular fibers}. The segment 
$O\times [0,1]$ is a circle too, called a
\textit{singular fiber} of multiplicity $k$ if $k\ne 1$ 
(if $k=1$ it is regular). Thus we have a
fibration with the $\mathbb{S}^1$\~fibers.
\smallskip

A \textit{Seifert fibration} of a 3-manifold $\Sigma$ is a 
continuous map $\pi\colon \Sigma \to B$ onto 2-manifold
$B$ satisfying the following.
For each point $b\in B$, there exists 
its neighborhood $U_b$  such that the fibration 
$\pi^{-1}(U_b)$ is isomorphic to the interior of the standard 
{\em fibered\,} solid torus. 
Singular fibers of multiplicity $k\geq 2$
correspond then the singular fibers above. 
This definition does not depend on the choice of the isomorphism.

The following fibrations of $\mathbb{S}^3$ with 
$\mathbb{S}^1$-fibers are important for us. 
For each unordered pair of coprime integers 
$\alpha_1,\alpha_2$ ($\alpha_i\geq 2$), there exists a 
projection $\pi_{\alpha_1,\alpha_2} \colon \mathbb{S}^3 \to 
\mathbb{S}^2$ with exactly two  
singular fibers with multiplicities $\alpha_1,\alpha_2$. 
This property uniquely determines such a fibration.
Following \cite{EN}, let us provide some details.

{\sf Topological description.}
Let $\mathbb{B} = \mathbb{S}^2\setminus \bigcup_{i=1}^2 U_i$ be 
the 2-sphere 
without two disjoint domains $U_i$, 
$E=\mathbb{S}^1\times \mathbb{B} \to \mathbb{B}$ be 
the trivial $\mathbb{S}^1$-bundle over it. The border $E$ 
consists of 
two tori, one for each 
$U_i$: $\partial E = \bigcup_{i=1}^2 T^i$. 
The required space is obtained by gluing $E$ with
two solid tori $\mathbb{D}^2\times \mathbb{S}^1$ 
(one for each $T_i$) along border.
Fixing a section $\sigma\colon \mathbb{B} \to E$,
the fiber $H_i$ in $T_i\subset E$ and 
$\mathsf{Q}_i = T_i\cap \sigma (\mathbb{B})$ form a basis in 
$H_1(T_i,\mathbb{Z})$ for $i=1,2$.

The fibers are as follows. Since $\alpha_i$ are pairwise 
coprime, there exist 
$\beta_i$ ($1\leq i \leq 2$) such that
\begin{equation*}
\beta_1 \alpha_2 + \beta_2 \alpha_1 = 1.
\end{equation*}

For each $i$, we then glue $T_i$ with the border of 
$\mathbb{D}^2\times 
\mathbb{S}^1$ in 
such a way that the homology class of $\alpha_i \mathsf{Q}_i 
+ \beta_i \mathsf{H}_i $ becomes zero. The resulting manifold is 
$\mathbb{S}^3$ with the  
inherited $\mathbb{S}^1$-fibration.
Note that there is a flexibility with the choices of 
$\sigma$ and $\beta_i$ here, but the output does not 
depend on this up to a diffeomorphism.
\smallskip

{\sf Analytic description.}
Let $\mathbb{S}^3 \subset \mathbb{C}^2$ be the standard unit 
sphere $|Z_1|^2+|Z_2|^2=1$ for $Z_i\in \mathbb{C}$. The fibers 
will be then the orbits of the following $\mathbb{S}^1$-action:
\begin{equation*}
t(Z_1, Z_2) = (t^{\alpha_2}Z_1, t^{\alpha_1}Z_2),
\end{equation*}
where $t\in \mathbb{C}$, $|t|=1$.
The singular fiber of multiplicity $\alpha_i$ are exactly 
the intersections of the sphere with the hyperplanes $Z_i = 0$.
This gives the required. 
\smallskip

{\sf The cases with $\alpha_i\leq 2$.}
In the construction of the fibration
$\pi_{\alpha_1,\alpha_2}$, we
imposed the conditions  $\alpha_i\geq 2$ (mainly to ensure 
that the fibrations have two singular fibers).
This restriction can be omitted. We have the following 3 cases.

$(i)$ 
If $\alpha_i = 1$ for one of $i=1,2$ or both,
then the corresponding fiber is regular
and such an extension is obvious. 

$(ii)$ 
If $\alpha_1 < 0$ or $\alpha_2<0$, then the resulting 
fibration of $\mathbb{S}^3$ will be the one for  
$\{|\alpha_1|,|\alpha_2|\}$  upon the
inversion of the orientation of $\mathbb{S}^3$ 
(as a $3$-manifold) when $\al_1\al_2<0$.

$(iii)$ 
If $\alpha_i = 0$, then $\alpha_j=\pm 1$ for
the remaining index $j$
due to gcd$(\al_i,\al_j)=1$, so the fiber for $\al_j$ is 
regular. The multiplicity of the
singular fiber is then $0$ and it is the unknot; 
the regular fibres are the meridians of the 
singular one. The $\mathbb{S}^1$-action with respect to 
this fibration is not free, but there is another
$\mathbb{S}^1$-action on $\mathbb{S}^3$ whose fixed point are 
precisely the points of the singular fiber.
\smallskip

{\sf Seifert links.}
Let us provide a series of ``canonical'' links, which will 
be then used as starting points for the operations above.
The \textit{Seifert link} is a link with the exterior that 
admits the Seifert fibration. It is known \cite{EN}, 
that every 
Seifert link in $\mathbb{S}^3$ is a collection of fibers in 
the Seifert fibration 
of $\mathbb{S}^3$. Let $(\alpha_1, \dots, \alpha_k)$
be a $k$\~tuple of non-negative integers
$(k\geq 2)$ with $\al_j=1$ for $j\neq i,i'$
for a certain pair $1\le i\neq i'\le k$ such that 
gcd$(\al_i,\al_{i'})=1$. Then it defines the link
\begin{equation*}
\mathbf{L}_{\alpha_1, \dots, \alpha_k} = \bigcup^{k}_{l=1} S_l,
\end{equation*}
where $S_l$ are the fibers of the Seifert fibration 
$\pi_{\alpha_i, \alpha_{i'}} \colon \mathbb{S}^3 \to 
\mathbb{S}^2$. Namely, if $l = i,i'$, they are the singular 
fibers of multiplicity $\alpha_i,\alpha_{i'}$ and 
regular fibers otherwise. 
Note that if there are several choices of such a pair 
$\alpha_i, \alpha_{i'}$ (then one of them must be  $1$), the
link $\mathbf{L}_{\alpha_1, \dots, \alpha_k}$ does not depend 
on such a choice.

We will represent 
$\mathbf{L}_{\alpha_1, \dots, \alpha_k}$ as follows:
\vskip -0.2cm

\begin{equation}
\begin{picture}(60,30)
\put(27,8){\vector(-1,1){20}}
\put(31,10){\vector(0,1){19}}
\put(34,9){\vector(1,1){20}}
\put(2,18){$\alpha_1$}
\put(19,22){$\alpha_i$}
\put(33,22){$\alpha_k$}
\put(11,28){$\dots$}
\put(36,28){$\dots$}
\put(25,2){$\pmb{\bigcirc}$}
\put(28,2){$\epsilon$}
\end{picture}
\text{with $\epsilon = +1$ or simply}
\begin{picture}(60,30)
\put(27,8){\vector(-1,1){20}}
\put(31,10){\vector(0,1){19}}
\put(34,9){\vector(1,1){20}}
\put(2,18){$\alpha_1$}
\put(19,22){$\alpha_i$}
\put(33,22){$\alpha_k$}
\put(11,28){$\dots$}
\put(36,28){$\dots$}
\put(25,2){$\bigoplus$}
\end{picture}. \label{SeifertLink}
\end{equation}
The reflection 
$-\mathbf{L}_{\alpha_1, \dots, \alpha_k}$ of
$\mathbf{L}_{\alpha_1, \dots, \alpha_k}$ 
will then correspond to $\epsilon=-1$.
 
\comment{
$\mathbf{L}_{\alpha_1, \dots, \alpha_k}$ (i.e. the same 
link in $\mathbb{S}^3$ with different orientation) as
\vskip -0.2cm
\begin{equation*}
\begin{picture}(60,30)
\put(27,8){\vector(-1,1){20}}
\put(31,10){\vector(0,1){19}}
\put(34,9){\vector(1,1){20}}
\put(2,18){$\alpha_1$}
\put(19,22){$\alpha_i$}
\put(33,22){$\alpha_k$}
\put(11,28){$\dots$}
\put(36,28){$\dots$}
\put(25,2){$\pmb{\bigcirc}$\ \.}
\put(28,2){$\epsilon$}
\end{picture}
\text{with } \epsilon = -1
\end{equation*}
}

If we do not assume that all $\alpha_i$ are non-negative, 
the corresponding link can be defined by the relations
\vskip -0.2cm

\begin{equation}\label{minus-orient}
\begin{picture}(60,30)
\put(27,8){\vector(-1,1){20}}
\put(31,10){\vector(0,1){19}}
\put(34,9){\vector(1,1){20}}
\put(2,18){$\alpha_1$}
\put(19,22){$\alpha_i$}
\put(33,22){$\alpha_k$}
\put(11,28){$\dots$}
\put(36,28){$\dots$}
\put(25,2){$\pmb{\bigcirc}$}
\put(28,2){$\epsilon$}
\end{picture}
\begin{picture}(10,30)
\put(0,12){$\sim$}
\end{picture}
\begin{picture}(60,30)
\put(27,8){\vector(-1,1){20}}
\put(31,10){\vector(0,1){19}}
\put(34,9){\vector(1,1){20}}
\put(2,18){$\alpha_1$}
\put(10,22){$\,-\!\alpha_i\!$}
\put(33,22){$\alpha_k$}
\put(11,29){$\dots$}
\put(36,29){$\dots$}
\put(25,2){$\pmb{\bigcirc}$}
\put(27,31){$\vee$}
\put(28,2){$\delta$}
\end{picture},
\end{equation}
where $\epsilon\delta = -1$; by $\vee$, we mean the
reversion of the orientation of the corresponding
component. Here $i$ can be 
arbitrary, so (\ref{minus-orient}) allows to
define links for any
$k$-tuples $(\alpha_1, \dots, \alpha_k)$ such that
\begin{equation}\label{Condition}
\alpha_j\!=\!\pm 1 \hbox{\, for\, }
 j\neq i,i', \hbox{\, where }
 1\!\le i\neq i'\!\le k,\,\, \hbox{gcd}
(\alpha_i,\alpha_{i'})=1.
\end{equation}
\smallskip

Let us discuss the framing. In  DAHA theory, only
{\em regular components\,} are needed, which
are regular fibers (those of multiplicity 
$\pm 1$). Any fibers have tubular neighborhoods with unions of 
Seifert fibres as their borders. The regular fibers are
homotopic to the fibers belonging to the border of the 
corresponding tubular neighborhood. This supplies them with
a canonical framing, called \textit{Seifert framing}. For 
irregular components, one may choose any framing (say, 
topological), but such components are not actually
needed in our work. 

The topological framing is the most common choice. 
However the Seifert framing appeared the one  
serving the DAHA invariants. Let us comment on 
the difference.

The regular components $S_i$ 
and $S_j$ of a given link $\pm \mathbf{L}_{\alpha_1, 
\dots, \alpha_k}$ have the linking number 
$lk=\pm \alpha_1 
\dots \hat{\alpha}_i \dots \hat{\alpha}_j \dots \alpha_k$ 
(it results from the formula for 
linking numbers in \cite{EN}, Chapter III, Section 10). 
Here $\hat{\alpha}_i$ means that $\alpha_i$ is omitted.
Let us add a (regular) component $S_l'$ to a regular component
$S_l$ in  $\pm\mathbf{L}_{\alpha_1, \dots, \alpha_k}$
using the Seifert framing. The resulting link will be  
$\pm \mathbf{L}_{\alpha_1, \dots, \alpha_k} \cup S'_l$,  
which is $\pm \mathbf{L}_{\alpha_1, \dots, \alpha_k, 1}$ for
$\alpha_l=1$ and 
$\mp \mathbf{L}_{\alpha_1, \dots, \alpha_k, -1}$ for
$\alpha_l=-1$
by definition. The linking number between
$S_l$ and $S_l'$ will be
$\pm\alpha_l\cdot\alpha_1\dots\hat{\alpha}_l\dots\alpha_k =
\pm\alpha_1\dots\alpha_k$, which is
the difference between the Seifert framing
and the topological framing
(when the linking number between $S$ and $S'$ would be zero).
\smallskip

{\sf Splice diagrams.}
They are trees or disjoint unions of trees, where the 
vertices can be {\em nodes\,} labeled by $\pm1$, 
unlabeled nodes, called {\em leaves}, and 
{\em arrowheads\,};
the pairs $\{$a labeled node, an edge from it$\}$ are
decorated by integers, called {\em weights}. 
In greater detail, the vertices are as follows.

\textit{An arrowhead}, at the end of an edge:
\begin{equation*}
\begin{picture}(15,10)
\put(0,5){\vector(1,0){13}}
\put(5,7){$S$}
\end{picture}
\;\mbox{ or }\;
\begin{picture}(19,10)
\put(0,5){\vector(1,0){13}}
\put(5,7){$S$}
\put(14,1){$\vee$\ .}
\end{picture}
\end{equation*}
It is a vertex of valence $1$.
The arrowheads are topologically interpreted as
the components of the resulting links. We mostly
put the names of the components over the
arrowheads, $S$ or $S_\vee$ in these 
examples. Here $\vee$ stands for the 
change of orientation of the component $S$; see below.
The $S$\~labels are necessary in the operations on links.

\textit{A leaf}, an unlabeled vertex of valence $1$:
\begin{equation*}
\begin{picture}(15,10)
\put(0,5){\line(1,0){13}}
\put(12,2){$\circ$\ .}
\end{picture}
\end{equation*}

\textit{A general vertex}, a node labeled by $\pm1$:
\begin{equation*}
\begin{picture}(60,30)
\put(27,8){\line(-1,1){20}}
\put(31,10){\line(0,1){19}}
\put(34,9){\line(1,1){20}}
\put(2,18){$\alpha_1$}
\put(19,22){$\alpha_i$}
\put(33,22){$\alpha_k$}
\put(11,28){$\dots$}
\put(36,28){$\dots$}
\put(25,2){$\pmb{\bigcirc}\ \ \ \ \ \ .$}
\put(28,2){$\epsilon$}
\end{picture}
\end{equation*}
This must be a vertex of valence $\geq 3$ with the label 
$\epsilon = \pm 1$ in it. Every adjacent edge has an integer 
weight $\alpha_i$ such that $(\alpha_1, \dots, \alpha_k)$ 
satisfies \ref{Condition}.
As in (\ref{SeifertLink}), we replace the node
with $\epsilon = +1$ in it by $\bigoplus$.
\smallskip

Given a splice diagram $\Gamma$, the corresponding  
link $\mathbf{L}(\Gamma)$ is constructed as follows.

A general vertex with the label $\epsilon$ and the
weights $\alpha_1, \dots, \alpha_k$ at the adjacent edges 
is interpreted as the 
link $\epsilon\mathbf{L}_{\alpha_1, \dots, \alpha_k}$.
It has the standard orientation for $\epsilon = 1$, which must
be reversed for $\epsilon = -1$. The adjacent edges give 
the components of this link; the weights show their 
multiplicities (as singular fibers).

An edge that does not have the arrowhead is interpreted 
as the {\em splice\,} of two links associated with the
vertices at these ends (along the corresponding components).
The arrowheads provide the components of 
$\mathbf{L}(\Gamma)$; 
if there is $\vee$ there, the orientation 
of the corresponding component must be reversed. 

A  leaf alone is interpreted as the unknot. Replacing
an arrowhead by a leaf topologically means {\em deleting\,} the
corresponding component. The {\em disjoint sum\,} of graphs 
corresponds to the disjoint sum of links. See also below.
\smallskip

The links that are constructed by splice 
diagrams form the class of \textit{graph links}. Note that 
we restrict our definitions to a subclass sufficient to work 
with links in $\mathbb{S}^3$. See \cite{EN}
for arbitrary graph links.

\subsection{\bf Operations on links}
The operations from section \ref{LinkSection} 
for the links associated with such diagrams
naturally result in operations on splice diagrams. 
The graphs of links $\mathbf{L}$ will be presented as boxes 
with $\mathbf{L}$ inside in the figures below. Sometimes
we show the arrows from this box corresponding to
the components of the link labeled by $S,S'$ and so on. 
We will use the operations $S^{\vee}$, $S\sharp\, S'$ etc., 
defined above and all previous notations for the link operations 
will be used. We will show the name of the initial component 
(before the transformation) in the figures below; the
corresponding arrowhead is replaced as follows.

\textsf{The erasure of a component}. We replace the arrow 
labeled by $S$ in the graph of $\mathbf{L}$ by the leaf. This
corresponds to erasing component $S$: 
\vskip -0.5cm
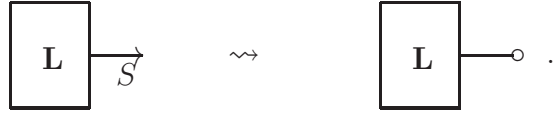
\begin{equation*}
\begin{picture}(140,45)
\put(5,0){\line(1,0){30}}
\put(5,0){\line(0,1){40}}
\put(5,40){\line(1,0){30}}
\put(35,0){\line(0,1){40}}
\put(17,15){$\mathbf{L}$}
\put(45,9){$S$}
\put(35,20){\line(1,0){20}}
\put(44,17){$\rightarrow$ $\ \ \ \ \ \ \ \rightsquigarrow$}
\end{picture}
\begin{picture}(80,45)
\put(5,0){\line(1,0){30}}
\put(5,0){\line(0,1){40}}
\put(5,40){\line(1,0){30}}
\put(35,0){\line(0,1){40}}
\put(17,15){$\mathbf{L}$}
\put(37,9){$$}
\put(35,20){\line(1,0){20}}
\put(54,17){$\circ$\ \ .}
\end{picture}
\end{equation*}

\textsf{The orientation reversion} of a component $S$ 
of a link $\mathbf{L}$. The new (changed) component is 
denoted by  $S^\vee$:
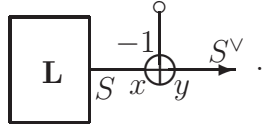
\begin{equation*}
\begin{picture}(100,45)
\put(5,0){\line(1,0){30}}
\put(5,0){\line(0,1){40}}
\put(5,40){\line(1,0){30}}
\put(35,0){\line(0,1){40}}
\put(17,15){$\mathbf{L}$}
\put(37,9){$S$}
\put(80,23){$S^\vee$}
\put(30,15){
    \put(35,5){\vector(1,0){25}\ \ .}
    \put(26,5){\line(-1,0){21}}
    \put(31.5,7){\line(0,1){20}}
    \put(25,2){$\bigoplus$}
    \put(37,-4){$y$}
    \put(20,-4){$x$}
    \put(15,12){$-1$}
    \put(28.5,26){$\circ$}
    }
\end{picture}
\end{equation*}
The weights $x$ and $y$ here are arbitrary; the link does not 
depend on them, only the Seifert framing does. So one can 
choose them to obtain the desired Seifert framing. 

\textsf{The connected sum} of two links $\mathbf{L}_1$ and 
$\mathbf{L}_2$ along $S$ and $S'$. The new (changed)
component is denoted by $S\sharp\, S'$.
\vskip -0.1cm
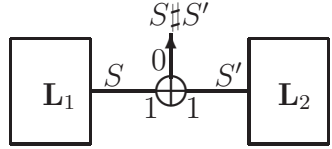
\begin{equation*}
\begin{picture}(140,60)
\put(0,0){
    \put(5,5){\line(1,0){30}}
    \put(5,5){\line(0,1){40}}
    \put(5,45){\line(1,0){30}}
    \put(35,5){\line(0,1){40}}
    \put(17,20){$\mathbf{L}_1$}
    \put(40,27){$S$}
    }
\put(80,0){
    \put(15,5){\line(1,0){30}}
    \put(15,5){\line(0,1){40}}
    \put(15,45){\line(1,0){30}}
    \put(45,5){\line(0,1){40}}
    \put(26,20){$\mathbf{L}_2$}
    \put(2,27){$S'$}
    }
\put(35,25){\line(1,0){25}}
\put(71,25){\line(1,0){24}}
\put(59,22){$\bigoplus$}
\put(55,15){$1$}
\put(71,15){$1$}
\put(65.5,30){\vector(0,1){17}}
\put(58,33){$0$}
\put(57,50){$S\sharp S'$}
\end{picture}.
\end{equation*}

\textsf{The cable} of type $(a,r)$ along any component 
$S$ of a link $\mathbf{L}$. The new component is called 
$C\!ab_{a,r}(S)$:
\vskip -0.3cm
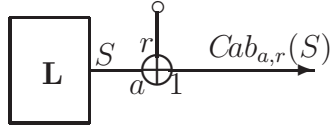
\begin{equation*}
\begin{picture}(130,60)
\put(5,5){\line(1,0){30}}
\put(5,5){\line(0,1){40}}
\put(5,45){\line(1,0){30}}
\put(35,5){\line(0,1){40}}
\put(17,20){$\mathbf{L}$}
\put(37,27){$S$}
\put(80,30){$C\!ab_{a,r}(S)$}
\put(30,20){
    \put(35,5){\vector(1,0){55}}
    \put(25,5){\line(-1,0){20}}
    \put(30.5,10){\line(0,1){17}}
    \put(24,2){$\bigoplus$}
    \put(35,-4){$1$}
    \put(20,-4){$a$}
    \put(24,12){$r$}
    \put(28,26){$\circ$}
    }
\end{picture}.
\end{equation*}

From the definition of the splice diagrams, it is obvious that 
the spice of two graph links is a graph link. The same is true 
for the disjoint sum. Moreover, the following theorem holds.

{\theorem The set of all graph links in $\mathbb{S}^3$ is the 
set of all 
solvable links, i.e. all links which can be constructed 
from the unknot $\unknot\,$ by any number of cablings, disjoint 
and connected sums. }

The proof is in \cite{EN} (Theorem 9.2).\sq

\subsection{\bf Equivalent diagrams}
Now we provide all relations between splice diagrams. We 
denote the link made from $\mathbf{L}$ by reversing of all 
orientations of components as $\mathbf{L}^\vee$, which is
adding $\vee$ to the components that have none or deleting it 
for the components with $\vee$.

{\theorem (\cite{EN}, Theorem 8.1) The following 
relations hold:

\noindent
$(i)$
\vskip -1.0cm
\begin{equation*}
\mathbf{L} = \mathbf{L}^\vee
\end{equation*}
(i.e. when the orientations of all components are reversed);
\smallskip

\noindent
$(ii)$
\vskip -1.0cm
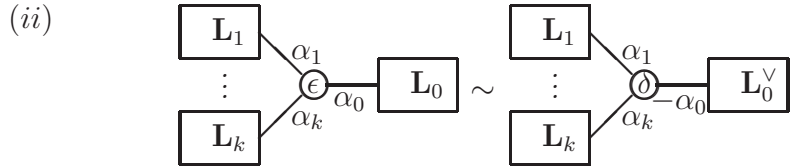
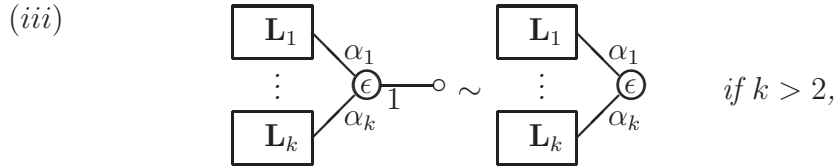
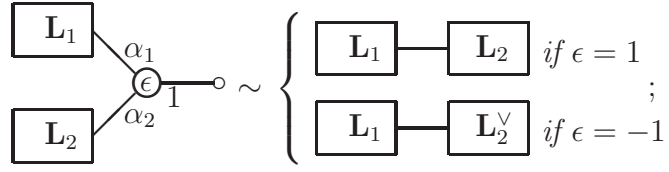
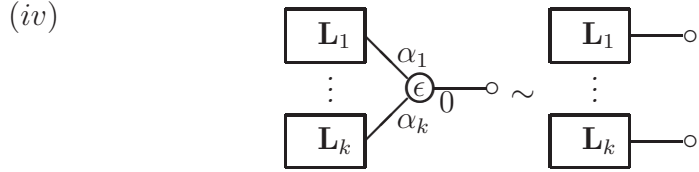
\begin{equation*}
\begin{picture}(115,65)
\put(5,0){\line(1,0){30}}
\put(5,0){\line(0,1){20}}
\put(5,20){\line(1,0){30}}
\put(35,0){\line(0,1){20}}
\put(5,40){\line(1,0){30}}
\put(5,40){\line(0,1){20}}
\put(5,60){\line(1,0){30}}
\put(35,40){\line(0,1){20}}
\put(17,47){$\mathbf{L}_1$}
\put(17,7){$\mathbf{L}_k$}
\put(47,40){$\alpha_1$}
\put(47,15){$\alpha_k$}
\put(35,50){\line(1,-1){16}}
\put(35,10){\line(1,1){16}}
\put(20,25){$\vdots$}
\put(50,27){$\pmb{\bigcirc}$}
\put(53,27){$\epsilon$}
\put(63,22){$\alpha_0$}
\put(60,30){\line(1,0){20}}
\put(80,20){\line(1,0){30}}
\put(80,20){\line(0,1){20}}
\put(80,40){\line(1,0){30}}
\put(110,20){\line(0,1){20}}
\put(92,27){$\mathbf{L}_0$}
\end{picture}
\begin{picture}(10,65)
\put(0,25){$\sim$}
\end{picture}
\begin{picture}(115,65)
\put(5,0){\line(1,0){30}}
\put(5,0){\line(0,1){20}}
\put(5,20){\line(1,0){30}}
\put(35,0){\line(0,1){20}}
\put(5,40){\line(1,0){30}}
\put(5,40){\line(0,1){20}}
\put(5,60){\line(1,0){30}}
\put(35,40){\line(0,1){20}}
\put(17,47){$\mathbf{L}_1$}
\put(17,7){$\mathbf{L}_k$}
\put(47,40){$\alpha_1$}
\put(47,15){$\alpha_k$}
\put(35,50){\line(1,-1){16}}
\put(35,10){\line(1,1){16}}
\put(20,25){$\vdots$}
\put(50,27){$\pmb{\bigcirc}$}
\put(53,27){$\delta$}
\put(58,21){$-\alpha_0$}
\put(60,30){\line(1,0){20}}
\put(80,20){\line(1,0){30}}
\put(80,20){\line(0,1){20}}
\put(80,40){\line(1,0){30}}
\put(110,20){\line(0,1){20}}
\put(92,27){$\mathbf{L}_0^\vee$}
\end{picture}
\end{equation*}
for $\delta = - \epsilon$;
\smallskip

\noindent
$(iii)$
\vskip -1.0cm
\begin{equation*}
\begin{picture}(90,65)
\put(5,0){\line(1,0){30}}
\put(5,0){\line(0,1){20}}
\put(5,20){\line(1,0){30}}
\put(35,0){\line(0,1){20}}
\put(5,40){\line(1,0){30}}
\put(5,40){\line(0,1){20}}
\put(5,60){\line(1,0){30}}
\put(35,40){\line(0,1){20}}
\put(17,47){$\mathbf{L}_1$}
\put(17,7){$\mathbf{L}_k$}
\put(47,40){$\alpha_1$}
\put(47,15){$\alpha_k$}
\put(35,50){\line(1,-1){16}}
\put(35,10){\line(1,1){16}}
\put(20,25){$\vdots$}
\put(50,27){$\pmb{\bigcirc}$}
\put(53,27){$\epsilon$}
\put(63,21){$1$}
\put(61,30){\line(1,0){20}}
\put(80,27){$\circ$}
\end{picture}
\begin{picture}(10,65)
\put(0,25){$\sim$}
\end{picture}
\begin{picture}(90,65)
\put(5,0){\line(1,0){30}}
\put(5,0){\line(0,1){20}}
\put(5,20){\line(1,0){30}}
\put(35,0){\line(0,1){20}}
\put(5,40){\line(1,0){30}}
\put(5,40){\line(0,1){20}}
\put(5,60){\line(1,0){30}}
\put(35,40){\line(0,1){20}}
\put(17,47){$\mathbf{L}_1$}
\put(17,7){$\mathbf{L}_k$}
\put(47,40){$\alpha_1$}
\put(47,15){$\alpha_k$}
\put(35,50){\line(1,-1){16}}
\put(35,10){\line(1,1){16}}
\put(20,25){$\vdots$}
\put(50,27){$\pmb{\bigcirc}$}
\put(53,27){$\epsilon$}
\end{picture}
\begin{picture}(10,65)
\put(0,25){if $k>2$,}
\end{picture}
\end{equation*}
\vskip -0.5cm
\begin{equation*}
\begin{picture}(90,65)
\put(5,0){\line(1,0){30}}
\put(5,0){\line(0,1){20}}
\put(5,20){\line(1,0){30}}
\put(35,0){\line(0,1){20}}
\put(5,40){\line(1,0){30}}
\put(5,40){\line(0,1){20}}
\put(5,60){\line(1,0){30}}
\put(35,40){\line(0,1){20}}
\put(17,47){$\mathbf{L}_1$}
\put(17,7){$\mathbf{L}_2$}
\put(47,40){$\alpha_1$}
\put(47,15){$\alpha_2$}
\put(35,50){\line(1,-1){16}}
\put(35,10){\line(1,1){16}}
\put(50,27){$\pmb{\bigcirc}$}
\put(53,27){$\epsilon$}
\put(63,21){$1$}
\put(61,30){\line(1,0){20}}
\put(80,27){$\circ$}
\end{picture}
\begin{picture}(10,65)
\put(0,25){$\sim$}
\end{picture}
\begin{picture}(90,25)
\put(0,25){
    $\begin{cases}
    \begin{picture}(115,25)
        \put(5,0){\line(1,0){30}}
        \put(5,0){\line(0,1){20}}
        \put(5,20){\line(1,0){30}}
        \put(35,0){\line(0,1){20}}
        \put(55,0){\line(1,0){30}}
        \put(55,0){\line(0,1){20}}
        \put(55,20){\line(1,0){30}}
        \put(85,0){\line(0,1){20}}
        \put(17,7){$\mathbf{L}_1$}
        \put(65,7){$\mathbf{L}_2$}
        \put(35,10){\line(1,0){20}}
        \put(90,5){if $\epsilon=1$}
    \end{picture} \\
    \begin{picture}(115,25)
        \put(5,0){\line(1,0){30}}
        \put(5,0){\line(0,1){20}}
        \put(5,20){\line(1,0){30}}
        \put(35,0){\line(0,1){20}}
        \put(55,0){\line(1,0){30}}
        \put(55,0){\line(0,1){20}}
        \put(55,20){\line(1,0){30}}
        \put(85,0){\line(0,1){20}}
        \put(17,7){$\mathbf{L}_1$}
        \put(65,7){$\mathbf{L}_2^\vee$}
        \put(35,10){\line(1,0){20}}
        \put(90,5){if $\epsilon=-1$}
    \end{picture}
    \end{cases};$
    }
\end{picture}
\end{equation*}
\smallskip

\noindent
$(iv)$
\vskip -1.0cm
\begin{equation*}
\begin{picture}(90,65)
\put(5,0){\line(1,0){30}}
\put(5,0){\line(0,1){20}}
\put(5,20){\line(1,0){30}}
\put(35,0){\line(0,1){20}}
\put(5,40){\line(1,0){30}}
\put(5,40){\line(0,1){20}}
\put(5,60){\line(1,0){30}}
\put(35,40){\line(0,1){20}}
\put(17,47){$\mathbf{L}_1$}
\put(17,7){$\mathbf{L}_k$}
\put(47,40){$\alpha_1$}
\put(47,15){$\alpha_k$}
\put(35,50){\line(1,-1){16}}
\put(35,10){\line(1,1){16}}
\put(20,25){$\vdots$}
\put(50,27){$\pmb{\bigcirc}$}
\put(53,27){$\epsilon$}
\put(63,21){$0$}
\put(61,30){\line(1,0){20}}
\put(80,27){$\circ$}
\end{picture}
\begin{picture}(10,65)
\put(0,25){$\sim$}
\end{picture}
\begin{picture}(60,65)
\put(5,0){\line(1,0){30}}
\put(5,0){\line(0,1){20}}
\put(5,20){\line(1,0){30}}
\put(35,0){\line(0,1){20}}
\put(5,40){\line(1,0){30}}
\put(5,40){\line(0,1){20}}
\put(5,60){\line(1,0){30}}
\put(35,40){\line(0,1){20}}
\put(17,47){$\mathbf{L}_1$}
\put(17,7){$\mathbf{L}_k$}
\put(35,50){\line(1,0){20}}
\put(35,10){\line(1,0){20}}
\put(55,47){$\circ$}
\put(55,7){$\circ$}
\put(20,25){$\vdots$}
\end{picture}
\end{equation*}
(the disjoint sum in the right-hand side);
\smallskip

\noindent
$(v)$
\vskip -1.0cm
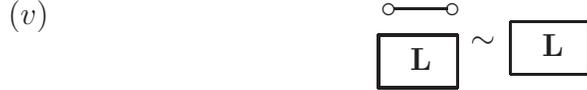
\begin{equation*}
\begin{picture}(40,35)
\put(5,0){\line(1,0){30}}
\put(5,0){\line(0,1){20}}
\put(5,20){\line(1,0){30}}
\put(35,0){\line(0,1){20}}
\put(11,30){\line(1,0){20}}
\put(6,27){$\circ$}
\put(30,27){$\circ$}
\put(17,7){$\mathbf{L}$}
\end{picture}
\begin{picture}(10,35)
\put(0,15){$\sim$}
\end{picture}
\begin{picture}(40,35)
\put(5,5){\line(1,0){30}}
\put(5,5){\line(0,1){20}}
\put(5,25){\line(1,0){30}}
\put(35,5){\line(0,1){20}}
\put(17,12){$\mathbf{L}$}
\end{picture}
\end{equation*}
(the disjoint sum in the left-hand side);
\smallskip

\noindent
$(vi)$
If $\alpha_0 \alpha_0' = \gamma\delta\alpha_1\dots
\alpha_k\alpha'_1\dots\alpha'_l$ define $\epsilon = \pm 1$ 
such that
\begin{equation*}
\alpha_ 0 = \epsilon \gamma \alpha'_1\dots\alpha'_l,
\end{equation*}
\vskip -0.5cm
\begin{equation*}
\alpha'_ 0 = \epsilon \delta \alpha_1\dots\alpha_k.
\end{equation*}
$(vi'\,)$ If $\gamma\delta\epsilon = 1$, then
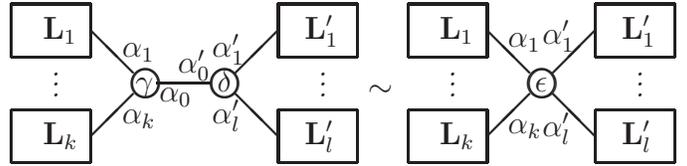
\begin{equation*}
\begin{picture}(140,65)
\put(5,0){\line(1,0){30}}
\put(5,0){\line(0,1){20}}
\put(5,20){\line(1,0){30}}
\put(35,0){\line(0,1){20}}
\put(5,40){\line(1,0){30}}
\put(5,40){\line(0,1){20}}
\put(5,60){\line(1,0){30}}
\put(35,40){\line(0,1){20}}
\put(17,47){$\mathbf{L}_1$}
\put(17,7){$\mathbf{L}_k$}
\put(47,40){$\alpha_1$}
\put(47,15){$\alpha_k$}
\put(35,50){\line(1,-1){16}}
\put(35,10){\line(1,1){16}}
\put(20,25){$\vdots$}
\put(50,27){$\pmb{\bigcirc}$}
\put(52,27){$\gamma$}
\put(61,23){$\alpha_0$}
\put(68,34){$\alpha'_0$}
\put(60,30){\line(1,0){20}}
\put(80,27){$\pmb{\bigcirc}$}
\put(83,27){$\delta$}
\put(81,40){$\alpha'_1$}
\put(81,15){$\alpha'_l$}
\put(1,0){
    \put(105,50){\line(-1,-1){16}}
    \put(105,10){\line(-1,1){16}}
    \put(105,0){\line(1,0){30}}
    \put(105,0){\line(0,1){20}}
    \put(105,20){\line(1,0){30}}
    \put(135,0){\line(0,1){20}}
    \put(105,40){\line(1,0){30}}
    \put(105,40){\line(0,1){20}}
    \put(105,60){\line(1,0){30}}
    \put(135,40){\line(0,1){20}}
    \put(115,47){$\mathbf{L}'_1$}
    \put(115,7){$\mathbf{L}'_l$}
    \put(120,25){$\vdots$}
    }
\end{picture}
\begin{picture}(10,65)
\put(0,25){$\sim$}
\end{picture}
\begin{picture}(110,65)
\put(5,0){\line(1,0){30}}
\put(5,0){\line(0,1){20}}
\put(5,20){\line(1,0){30}}
\put(35,0){\line(0,1){20}}
\put(5,40){\line(1,0){30}}
\put(5,40){\line(0,1){20}}
\put(5,60){\line(1,0){30}}
\put(35,40){\line(0,1){20}}
\put(17,47){$\mathbf{L}_1$}
\put(17,7){$\mathbf{L}_k$}
\put(43,44){$\alpha_1$}
\put(43,11){$\alpha_k$}
\put(35,50){\line(1,-1){16}}
\put(35,10){\line(1,1){16}}
\put(20,25){$\vdots$}
\put(50,27){$\pmb{\bigcirc}$}
\put(53,27){$\epsilon$}
\put(56,45){$\alpha'_1$}
\put(56,10){$\alpha'_l$}
\put(1,0){
    \put(75,50){\line(-1,-1){16}}
    \put(75,10){\line(-1,1){16}}
    \put(75,0){\line(1,0){30}}
    \put(75,0){\line(0,1){20}}
    \put(75,20){\line(1,0){30}}
    \put(105,0){\line(0,1){20}}
    \put(75,40){\line(1,0){30}}
    \put(75,40){\line(0,1){20}}
    \put(75,60){\line(1,0){30}}
    \put(105,40){\line(0,1){20}}
    \put(85,47){$\mathbf{L}'_1$}
    \put(85,7){$\mathbf{L}'_l$}
    \put(90,25){$\vdots$}
    }
\end{picture}.
\end{equation*}
$(vi'')$ If $\gamma\delta\epsilon = -1$, then
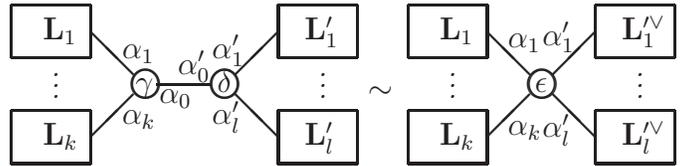
\begin{equation*}
\begin{picture}(140,65)
\put(5,0){\line(1,0){30}}
\put(5,0){\line(0,1){20}}
\put(5,20){\line(1,0){30}}
\put(35,0){\line(0,1){20}}
\put(5,40){\line(1,0){30}}
\put(5,40){\line(0,1){20}}
\put(5,60){\line(1,0){30}}
\put(35,40){\line(0,1){20}}
\put(17,47){$\mathbf{L}_1$}
\put(17,7){$\mathbf{L}_k$}
\put(47,40){$\alpha_1$}
\put(47,15){$\alpha_k$}
\put(35,50){\line(1,-1){16}}
\put(35,10){\line(1,1){16}}
\put(20,25){$\vdots$}
\put(50,27){$\pmb{\bigcirc}$}
\put(52,27){$\gamma$}
\put(61,23){$\alpha_0$}
\put(68,34){$\alpha'_0$}
\put(60,30){\line(1,0){20}}
\put(80,27){$\pmb{\bigcirc}$}
\put(83,27){$\delta$}
\put(81,40){$\alpha'_1$}
\put(81,15){$\alpha'_l$}
\put(1,0){
    \put(105,50){\line(-1,-1){16}}
    \put(105,10){\line(-1,1){16}}
    \put(105,0){\line(1,0){30}}
    \put(105,0){\line(0,1){20}}
    \put(105,20){\line(1,0){30}}
    \put(135,0){\line(0,1){20}}
    \put(105,40){\line(1,0){30}}
    \put(105,40){\line(0,1){20}}
    \put(105,60){\line(1,0){30}}
    \put(135,40){\line(0,1){20}}
    \put(115,47){$\mathbf{L}'_1$}
    \put(115,7){$\mathbf{L}'_l$}
    \put(120,25){$\vdots$}
    }
\end{picture}
\begin{picture}(10,65)
\put(0,25){$\sim$}
\end{picture}
\begin{picture}(110,65)
\put(5,0){\line(1,0){30}}
\put(5,0){\line(0,1){20}}
\put(5,20){\line(1,0){30}}
\put(35,0){\line(0,1){20}}
\put(5,40){\line(1,0){30}}
\put(5,40){\line(0,1){20}}
\put(5,60){\line(1,0){30}}
\put(35,40){\line(0,1){20}}
\put(17,47){$\mathbf{L}_1$}
\put(17,7){$\mathbf{L}_k$}
\put(43,44){$\alpha_1$}
\put(43,11){$\alpha_k$}
\put(35,50){\line(1,-1){16}}
\put(35,10){\line(1,1){16}}
\put(20,25){$\vdots$}
\put(50,27){$\pmb{\bigcirc}$}
\put(53,27){$\epsilon$}
\put(56,45){$\alpha'_1$}
\put(56,10){$\alpha'_l$}
\put(1,0){
    \put(75,50){\line(-1,-1){16}}
    \put(75,10){\line(-1,1){16}}
    \put(75,0){\line(1,0){30}}
    \put(75,0){\line(0,1){20}}
    \put(75,20){\line(1,0){30}}
    \put(105,0){\line(0,1){20}}
    \put(75,40){\line(1,0){30}}
    \put(75,40){\line(0,1){20}}
    \put(75,60){\line(1,0){30}}
    \put(105,40){\line(0,1){20}}
    \put(85,47){$\mathbf{L}_1^{\prime\vee}$}
    \put(85,7){$\mathbf{L}_l^{\prime\vee}$}
    \put(90,25){$\vdots$}
    }
\end{picture}.
\end{equation*}
(such $\epsilon$ exist, since $\alpha_0$ is prime to 
$\alpha_1\dots\alpha_k$ and $\alpha'_0$ is prime to 
$\alpha'_1\dots\alpha'_k$.)
\sq}

Splice diagrams are called equivalent if they can be obtained from
each other by these relations \cite{EN}. 
We have the following corollary, one of the key in
the DAHA approach; see below.

{\corollary For any general vertex and arbitrary  
$\epsilon_i = \pm 1$, 
\vskip -1.0cm
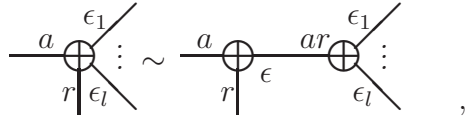
\begin{equation*}
\begin{picture}(50,45)
\put(11,23){$a$}
\put(0,20){\line(1,0){21}}
\put(20,17){$\bigoplus$}
\put(48,40){\line(-1,-1){17}}
\put(48,0){\line(-1,1){17}}
\put(40,15){$\vdots$}
\put(28,32){$\epsilon_1$}
\put(30,3){$\epsilon_l$}
\put(26.5,15){\line(0,-1){17}}
\put(20,3){$r$}
\end{picture}
\begin{picture}(10,65)
\put(0,15){$\sim$}
\end{picture}
\begin{picture}(110,45)
\put(11,23){$a$}
\put(5,20){\line(1,0){16}}
\put(20,17){$\bigoplus$}
\put(32,20){\line(1,0){29}}
\put(50,23){$ar$}
\put(60,17){$\bigoplus$}
\put(88,40){\line(-1,-1){17}}
\put(88,0){\line(-1,1){17}}
\put(80,15){$\vdots$}
\put(68,32){$\epsilon_1$}
\put(70,3){$\epsilon_l$}
\put(26.5,15){\line(0,-1){17}}
\put(20,3){$r$}
\put(35,10){$\epsilon$}
\end{picture},
\end{equation*}
where $\epsilon = \epsilon_1\dots \epsilon_l$. \label{CorDec}
\sq}

\smallskip
We call a splice diagram \textit{minimal} if no equivalent 
diagram has fewer edges. We call a splice diagram 
$\Gamma$ \textit{normal} if

$(i)$ $\Gamma$ is minimal;\ \ 
$(ii)$ all edge weights are non-negative;

$(iii)$ if an edge weight is zero, the adjacent vertex has 
the label $+1$.

{\theorem (Corollary 8.3 in \cite{EN}). Up to the transformation
$\Gamma\mapsto \Gamma^\vee$ (see above), there is a unique normal 
form for any splice diagram.\sq}
\medskip

\subsection{\bf Connection with DAHA}
We begin with the general description
of algebraic links from Theorem 9.4 in \cite{EN}.
All algebraic links are graph links; see Appendix to 
Chapter I in \cite{EN}. 
The following theorem describes all of them.

{\theorem The solvable link 
$\mathbf{L}(\Gamma)$ is an algebraic link if 

(a) all edge weights and vertex labels in $\Gamma$
are positive; 

(b) $\alpha_0 \alpha'_0 > \alpha_1 \dots \alpha_k\alpha'_1 
\dots \alpha'_l$ for every edge as follows: 
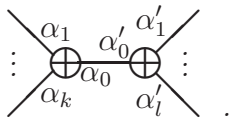
\begin{equation*}
\begin{picture}(80,45)
\put(17,30){$\alpha_1$}
\put(17,5){$\alpha_k$}
\put(5,40){\line(1,-1){17}}
\put(5,0){\line(1,1){17}}
\put(5,15){$\vdots$}
\put(20,17){$\bigoplus$}
\put(32,13){$\alpha_0$}
\put(39,25){$\alpha'_0$}
\put(31,20){\line(1,0){20}}
\put(50,17){$\bigoplus$}
\put(77,40){\line(-1,-1){16}}
\put(77,0){\line(-1,1){16}\,\ \ \ .}
\put(70,15){$\vdots$}
\put(53,33){$\alpha'_1$}
\put(53,3){$\alpha'_l$}
\end{picture}
\end{equation*}

Vice versa, if $\mathbf{L}$ is an algebraic 
link then its normal form graph satisfies the above conditions.
\sq}

\smallskip
These conditions automatically holds for the
{\em positive trees\,} from our paper.
The {\em twisted union\,} for $\{\l,\,'\!\l^\vee\}$
there corresponds to 
 
\vskip -0.5cm
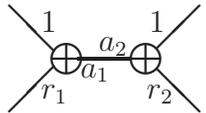
\begin{equation*}
\begin{picture}(80,45)
\put(17,30){$1$}
\put(17,5){$r_1$}
\put(5,40){\line(1,-1){17}}
\put(5,0){\line(1,1){17}}
\put(20,17){$\bigoplus$}
\put(32,13){$a_1$}
\put(39,23){$a_2$}
\put(31,20){\line(1,0){20}}
\put(50,17){$\bigoplus$}
\put(77,40){\line(-1,-1){16}}
\put(77,0){\line(-1,1){16}}
\put(58,30){$1$}
\put(57,5){$r_2$}
\end{picture}.
\end{equation*}
\vskip -0.2cm

Condition $(b)$ from the theorem then becomes
$a_1 a_2 > r_1 r_2$ (if the minimality holds).
This is the positivity condition
for the pairs $\{\l,\,'\!\l^\vee\}$ in our work;
the positive pairs are sufficient to
obtain {\em arbitrary\,} algebraic links. 
\medskip

{\sf Seifert framing.}
We will briefly discuss the framing needed in the DAHA
approach.  
 We need not only splices of links, but also 
{\em deformed\,} splices, defined in (\ref{deformfram})
at the end of 
Section \ref{LinkSection}; they are denoted by an edge marked 
by an arrow in the middle:
 $\begin{picture}(40,10)
 \put(5,5){\line(1,0){30}}
 \put(15,5){\vector(1,0){10}}
 \end{picture}$, not related to the arrowheads we use for
components.
This is a certain extension of the definition
of the splice diagram.  
They arrowed edges correspond to deformed splices in the
same way as the ordinary edges correspond to usual splices.

Let  $\mathbf{L}$ be a link with a chosen component $S$ in it. 
If the framing on $S$ differs by $D$ from the topological one, 
then the following links are equivalent:

\vskip -0.8cm
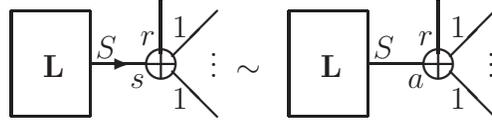
\begin{equation*}
\begin{picture}(90,60)
\put(5,5){\line(1,0){30}}
\put(5,5){\line(0,1){40}}
\put(5,45){\line(1,0){30}}
\put(35,5){\line(0,1){40}}
\put(17,20){$\mathbf{L}$}
\put(37,27){$S$}
\put(35,25){\vector(1,0){15}}
\put(30,20){
    \put(26,5){\line(-1,0){20}}
    \put(31.5,9){\line(0,1){20}}
    \put(25,2){$\bigoplus$}
    \put(53,25){\line(-1,-1){17}}
    \put(53,-15){\line(-1,1){17}}
    \put(50,0){$\vdots$}
    \put(36,15){$1$}
    \put(36,-12){$1$}
    \put(20,-4){$s$}
    \put(24,12){$r$}
    }
\end{picture}
\begin{picture}(15,60)
\put(0,20){$\sim$}
\end{picture}
\begin{picture}(90,60)
\put(5,5){\line(1,0){30}}
\put(5,5){\line(0,1){40}}
\put(5,45){\line(1,0){30}}
\put(35,5){\line(0,1){40}}
\put(17,20){$\mathbf{L}$}
\put(37,27){$S$}
\put(30,20){
    \put(26,5){\line(-1,0){21}}
    \put(31.5,9){\line(0,1){20}}
    \put(25,2){$\bigoplus$}
    \put(53,25){\line(-1,-1){17}}
    \put(53,-15){\line(-1,1){17}}
    \put(50,0){$\vdots$}
    \put(36,15){$1$}
    \put(36,-12){$1$}
    \put(20,-4){$a$}
    \put(24,12){$r$}
    }
\end{picture},
\end{equation*}
\vskip -0.3cm
\noindent
where $a = s+rD$. This is clear from the definitions.

If $S$ is a regular component of a Seifert link, then
it has the Seifert framing. One can use it to simplify the 
relation for the decomposition in Corollary \ref{CorDec},
which can be presented as follows:

\vskip -1.0cm
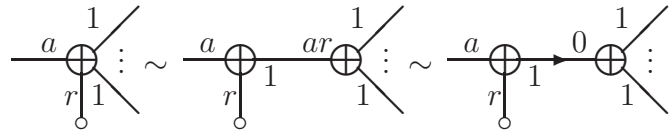
\begin{equation*}
\begin{picture}(50,45)
\put(11,23){$a$}
\put(0,20){\line(1,0){21}}
\put(20,17){$\bigoplus$}
\put(48,40){\line(-1,-1){17}}
\put(48,0){\line(-1,1){17}}
\put(40,15){$\vdots$}
\put(28,32){$1$}
\put(30,3){$1$}
\put(26.5,15){\line(0,-1){17}}
\put(20,3){$r$}
\put(23.5,-7){$\circ$}
\end{picture}
\begin{picture}(10,65)
\put(0,15){$\sim$}
\end{picture}
\begin{picture}(90,45)
\put(11,23){$a$}
\put(5,20){\line(1,0){16}}
\put(20,17){$\bigoplus$}
\put(31,20){\line(1,0){30}}
\put(50,23){$ar$}
\put(60,17){$\bigoplus$}
\put(88,40){\line(-1,-1){17}}
\put(88,0){\line(-1,1){17}}
\put(80,15){$\vdots$}
\put(68,32){$1$}
\put(70,3){$1$}
\put(26.5,15){\line(0,-1){17}}
\put(20,3){$r$}
\put(35,10){$1$}
\put(23.5,-7){$\circ$}
\end{picture}
\begin{picture}(10,65)
\put(0,15){$\sim$}
\end{picture}
\begin{picture}(90,45)
\put(11,23){$a$}
\put(5,20){\line(1,0){16}}
\put(20,17){$\bigoplus$}
\put(31,20){\vector(1,0){20}}
\put(31,20){\line(1,0){30}}
\put(52,23){$0$}
\put(60,17){$\bigoplus$}
\put(88,40){\line(-1,-1){17}}
\put(88,0){\line(-1,1){17}}
\put(80,15){$\vdots$}
\put(68,32){$1$}
\put(70,3){$1$}
\put(26.5,15){\line(0,-1){17}}
\put(20,3){$r$}
\put(35,10){$1$}
\put(23.5,-7){$\circ$}
\end{picture}
\end{equation*}

This is important for the DAHA-approach. Recall
the $(a,r,1)$-vertex corresponds to the projective
$PSL_2(\mathbb{Z})$\~action and 
$(0,1,\dots 1)$\~vertex corresponds
to the multiplication of pre-polynomials.
The usage of the Seifert framing completely clarifies the 
passage from the $(r,s)$-pairs to $(a,r)$-pairs in the 
DAHA-approach. Namely,
\vskip -0.5cm
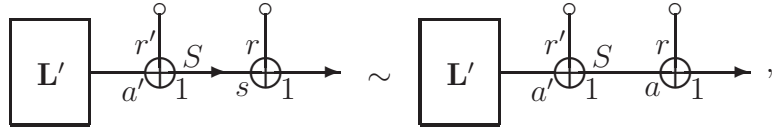
\begin{equation*}
\begin{picture}(140,60)
\put(5,5){\line(1,0){30}}
\put(5,5){\line(0,1){40}}
\put(5,45){\line(1,0){30}}
\put(35,5){\line(0,1){40}}
\put(15,20){$\mathbf{L}'$}
\put(70,27){$S$}
\put(30,20){
    \put(35,5){\line(1,0){35}}
    \put(26,5){\line(-1,0){21}}
    \put(31.5,7){\line(0,1){20}}
    \put(25,2){$\bigoplus$}
    \put(37,-5){$1$}
    \put(17,-6){$a'$}
    \put(22,12){$r'$}
    \put(28.5,26){$\circ$}
    }
\put(70,20){
    \put(35,5){\vector(1,0){25}}
    \put(5,5){\vector(1,0){12}}
    \put(31.5,7){\line(0,1){20}}
    \put(25,2){$\bigoplus$}
    \put(37,-5){$1$}
    \put(20,-4){$s$}
    \put(24,12){$r$}
    \put(28.5,26){$\circ$}
    }
\end{picture}
\begin{picture}(15,60)
\put(0,20){$\sim$}
\end{picture}
\begin{picture}(140,60)
\put(5,5){\line(1,0){30}}
\put(5,5){\line(0,1){40}}
\put(5,45){\line(1,0){30}}
\put(35,5){\line(0,1){40}}
\put(15,20){$\mathbf{L}'$}
\put(70,27){$S$}
\put(30,20){
    \put(35,5){\line(1,0){35}}
    \put(26,5){\line(-1,0){21}}
    \put(31.5,7){\line(0,1){20}}
    \put(25,2){$\bigoplus$}
    \put(37,-5){$1$}
    \put(17,-6){$a'$}
    \put(22,12){$r'$}
    \put(28.5,26){$\circ$}
    }
\put(70,20){
    \put(35,5){\vector(1,0){25}\,\ ,}
    \put(31.5,7){\line(0,1){20}}
    \put(25,2){$\bigoplus$}
    \put(37,-5){$1$}
    \put(20,-4){$a$}
    \put(24,12){$r$}
    \put(28.5,26){$\circ$}
    }
\end{picture}
\end{equation*}

\noindent
where $a = s+a'r'r$ for the Seifert framing on $S$. This 
is precisely the recursive formula for the switch from the
$(r,s)$-pairs to $(a,r)$-pairs in the DAHA-approach. The
arrows at the right ends in this figure are {\em arrowheads\,}
showing the components.

Thus applying  $\gamma_{r,s}$ to a pre-polynomial in DAHA theory
corresponds to splicing $\mathbf{L}\mapsto \mathbf{L}_{s,r,1}$ 
in topology.  See the main body of our work for this and other 
connections. The Seifert framing is exactly the one which
corresponds to the DAHA-approach. This results in the exact
topological interpretation of the DAHA-superpolynomials
of torus iterated links, i.e. without using the 
hat-normalization (ignoring factors 
$q^\bullet t^\bullet$). This will be discussed somewhere.
\medskip

\hbadness=10000
\vbadness=10000

\setcounter{equation}{0}
\section{\sc Double 
\texorpdfstring{$C\!ab(13,2)T(3,2)$}
{Cab(13,2)T(3,2)}}


Let us provide the DAHA-superpolynomial for the duplication of
the simplest non-torus algebraic knot $C\!ab(13,2)T(3,2)$.
It is long but, we think, it is of importance for
the theory of algebraic links.

\begin{align}\label{T12-2-6-4}
&C\!ab(13,2)T(6,4):\ 
\l=\l_{\{3,2\},\{2,1\}}^
{\,\circ\rightarrow\circ\rightrightarrows,\, 
(\yng(1)\,,\yng(1)\,)},\ \ 
\hat{\h}{}^{min}_{\l}\,(q,t,a)=
\end{align}

\renewcommand{\baselinestretch}{0.5} 
\noindent
{\small
\(
1-t+q t+q^2 t+q^3 t+q^4 t+q^5 t+q^6 t+q^7 t-q t^2+q^4 t^2+q^5 t^2
+2 q^6 t^2+2 q^7 t^2+4 q^8 t^2+3 q^9 t^2+3 q^{10} t^2+q^{11} t^2
+q^{12} t^2-q^2 t^3-q^4 t^3+q^7 t^3+q^8 t^3+4 q^9 t^3+4 q^{10} t^3
+7 q^{11} t^3+6 q^{12} t^3+6 q^{13} t^3+3 q^{14} t^3+2 q^{15} t^3
-q^3 t^4-q^5 t^4-q^6 t^4-q^7 t^4-q^9 t^4+2 q^{10} t^4+2 q^{11} t^4
+6 q^{12} t^4+7 q^{13} t^4+11 q^{14} t^4+9 q^{15} t^4+10 q^{16} t^4
+4 q^{17} t^4+2 q^{18} t^4-q^4 t^5-q^6 t^5-q^7 t^5-2 q^8 t^5-q^9t^5
-2 q^{10} t^5-q^{12} t^5+3 q^{13} t^5+4 q^{14} t^5+10 q^{15} t^5
+10 q^{16} t^5+16 q^{17} t^5+13 q^{18} t^5+11 q^{19} t^5
+3 q^{20} t^5+q^{21} t^5-q^5 t^6-q^7 t^6-q^8 t^6-2 q^9 t^6
-2 q^{10} t^6-3 q^{11} t^6-q^{12} t^6-3 q^{13} t^6+6 q^{16} t^6
+7 q^{17} t^6+15 q^{18} t^6+15 q^{19} t^6+21 q^{20} t^6
+12 q^{21} t^6+7 q^{22} t^6+q^{23} t^6-q^6 t^7-q^8 t^7-q^9 t^7
-2 q^{10} t^7-2 q^{11} t^7-4 q^{12} t^7-2 q^{13} t^7
-4 q^{14} t^7-2 q^{15} t^7-3 q^{16} t^7+2 q^{17} t^7+2 q^{18} t^7
+11 q^{19} t^7+12 q^{20} t^7+23 q^{21} t^7+19 q^{22} t^7
+17 q^{23} t^7+5 q^{24} t^7+q^{25} t^7-q^7 t^8-q^9 t^8-q^{10} t^8
-2 q^{11} t^8-2 q^{12} t^8-4 q^{13} t^8-3 q^{14} t^8-5 q^{15} t^8
-3 q^{16} t^8-5 q^{17} t^8-q^{18} t^8-2 q^{19} t^8+6 q^{20} t^8
+7 q^{21} t^8+20 q^{22} t^8+20 q^{23} t^8+27 q^{24} t^8
+11 q^{25} t^8+3 q^{26} t^8-q^8 t^9-q^{10} t^9-q^{11} t^9
-2 q^{12} t^9-2 q^{13} t^9-4 q^{14} t^9-3 q^{15} t^9-6 q^{16} t^9
-4 q^{17} t^9-6 q^{18} t^9-3 q^{19} t^9-5 q^{20} t^9+2 q^{21} t^9
+2 q^{22} t^9+15 q^{23} t^9+16 q^{24} t^9+32 q^{25} t^9
+18 q^{26} t^9+7 q^{27} t^9-q^9 t^{10}-q^{11} t^{10}-q^{12} t^{10}
-2 q^{13} t^{10}-2 q^{14} t^{10}-4 q^{15} t^{10}-3 q^{16} t^{10}
-6 q^{17} t^{10}-5 q^{18} t^{10}-7 q^{19} t^{10}-4 q^{20} t^{10}
-7 q^{21} t^{10}-q^{22} t^{10}-2 q^{23} t^{10}+10 q^{24} t^{10}
+11 q^{25} t^{10}+33 q^{26} t^{10}+22 q^{27} t^{10}+11 q^{28}t^{10}
-q^{10} t^{11}-q^{12} t^{11}-q^{13} t^{11}-2 q^{14} t^{11}
-2 q^{15} t^{11}-4 q^{16} t^{11}-3 q^{17} t^{11}-6 q^{18} t^{11}
-5 q^{19} t^{11}-8 q^{20} t^{11}-5 q^{21} t^{11}-8 q^{22} t^{11}
-3 q^{23} t^{11}-5 q^{24} t^{11}+6 q^{25} t^{11}+7 q^{26} t^{11}
+31 q^{27} t^{11}+24 q^{28} t^{11}+14 q^{29} t^{11}-q^{11} t^{12}
-q^{13} t^{12}-q^{14} t^{12}-2 q^{15} t^{12}-2 q^{16} t^{12}
-4 q^{17} t^{12}-3 q^{18} t^{12}-6 q^{19} t^{12}-5 q^{20} t^{12}
-8 q^{21} t^{12}-6 q^{22} t^{12}-9 q^{23} t^{12}-4 q^{24} t^{12}
-7 q^{25} t^{12}+3 q^{26} t^{12}+4 q^{27} t^{12}+30 q^{28} t^{12}
+25 q^{29} t^{12}+15 q^{30} t^{12}-q^{12} t^{13}-q^{14} t^{13}
-q^{15} t^{13}-2 q^{16} t^{13}-2 q^{17} t^{13}-4 q^{18} t^{13}
-3 q^{19} t^{13}-6 q^{20} t^{13}-5 q^{21} t^{13}-8 q^{22} t^{13}
-6 q^{23} t^{13}-10 q^{24} t^{13}-5 q^{25} t^{13}-8 q^{26} t^{13}
+2 q^{27} t^{13}+2 q^{28} t^{13}+29 q^{29} t^{13}+25 q^{30} t^{13}
+14 q^{31} t^{13}-q^{13} t^{14}-q^{15} t^{14}-q^{16} t^{14}
-2 q^{17} t^{14}-2 q^{18} t^{14}-4 q^{19} t^{14}-3 q^{20} t^{14}
-6 q^{21} t^{14}-5 q^{22} t^{14}-8 q^{23} t^{14}-6 q^{24} t^{14}
-10 q^{25} t^{14}-6 q^{26} t^{14}-9 q^{27} t^{14}+2 q^{28} t^{14}
+2 q^{29} t^{14}+30 q^{30} t^{14}+24 q^{31} t^{14}+11 q^{32} t^{14}
-q^{14} t^{15}-q^{16} t^{15}-q^{17} t^{15}-2 q^{18} t^{15}
-2 q^{19} t^{15}-4 q^{20} t^{15}-3 q^{21} t^{15}-6 q^{22} t^{15}
-5 q^{23} t^{15}-8 q^{24} t^{15}-6 q^{25} t^{15}-10 q^{26} t^{15}
-6 q^{27} t^{15}-9 q^{28} t^{15}+2 q^{29} t^{15}+4 q^{30} t^{15}
+31 q^{31} t^{15}+22 q^{32} t^{15}+7 q^{33} t^{15}-q^{15} t^{16}
-q^{17} t^{16}-q^{18} t^{16}-2 q^{19} t^{16}-2 q^{20} t^{16}
-4 q^{21} t^{16}-3 q^{22} t^{16}-6 q^{23} t^{16}-5 q^{24} t^{16}
-8 q^{25} t^{16}-6 q^{26} t^{16}-10 q^{27} t^{16}-6 q^{28} t^{16}
-8 q^{29} t^{16}+3 q^{30} t^{16}+7 q^{31} t^{16}+33 q^{32} t^{16}
+18 q^{33} t^{16}+3 q^{34} t^{16}-q^{16} t^{17}-q^{18} t^{17}
-q^{19} t^{17}-2 q^{20} t^{17}-2 q^{21} t^{17}-4 q^{22} t^{17}
-3 q^{23} t^{17}-6 q^{24} t^{17}-5 q^{25} t^{17}-8 q^{26} t^{17}
-6 q^{27} t^{17}-10 q^{28} t^{17}-5 q^{29} t^{17}-7 q^{30} t^{17}
+6 q^{31} t^{17}+11 q^{32} t^{17}+32 q^{33} t^{17}+11 q^{34} t^{17}
+q^{35} t^{17}-q^{17} t^{18}-q^{19} t^{18}-q^{20} t^{18}
-2 q^{21} t^{18}-2 q^{22} t^{18}-4 q^{23} t^{18}-3 q^{24} t^{18}
-6 q^{25} t^{18}-5 q^{26} t^{18}-8 q^{27} t^{18}-6 q^{28} t^{18}
-10 q^{29} t^{18}-4 q^{30} t^{18}-5 q^{31} t^{18}+10 q^{32} t^{18}
+16 q^{33} t^{18}+27 q^{34} t^{18}+5 q^{35} t^{18}-q^{18} t^{19}
-q^{20} t^{19}-q^{21} t^{19}-2 q^{22} t^{19}-2 q^{23} t^{19}
-4 q^{24} t^{19}-3 q^{25} t^{19}-6 q^{26} t^{19}-5 q^{27} t^{19}
-8 q^{28} t^{19}-6 q^{29} t^{19}-9 q^{30} t^{19}-3 q^{31} t^{19}
-2 q^{32} t^{19}+15 q^{33} t^{19}+20 q^{34} t^{19}+17 q^{35} t^{19}
+q^{36} t^{19}-q^{19} t^{20}-q^{21} t^{20}-q^{22} t^{20}
-2 q^{23} t^{20}-2 q^{24} t^{20}-4 q^{25} t^{20}-3 q^{26} t^{20}
-6 q^{27} t^{20}-5 q^{28} t^{20}-8 q^{29} t^{20}-6 q^{30} t^{20}
-8 q^{31} t^{20}-q^{32} t^{20}+2 q^{33} t^{20}+20 q^{34} t^{20}
+19 q^{35} t^{20}+7 q^{36} t^{20}-q^{20} t^{21}-q^{22} t^{21}
-q^{23} t^{21}-2 q^{24} t^{21}-2 q^{25} t^{21}-4 q^{26} t^{21}
-3 q^{27} t^{21}-6 q^{28} t^{21}-5 q^{29} t^{21}-8 q^{30} t^{21}
-5 q^{31} t^{21}-7 q^{32} t^{21}+2 q^{33} t^{21}+7 q^{34} t^{21}
+23 q^{35} t^{21}+12 q^{36} t^{21}+q^{37} t^{21}-q^{21} t^{22}
-q^{23} t^{22}-q^{24} t^{22}-2 q^{25} t^{22}-2 q^{26} t^{22}
-4 q^{27} t^{22}-3 q^{28} t^{22}-6 q^{29} t^{22}-5 q^{30} t^{22}
-8 q^{31} t^{22}-4 q^{32} t^{22}-5 q^{33} t^{22}+6 q^{34} t^{22}
+12 q^{35} t^{22}+21 q^{36} t^{22}+3 q^{37} t^{22}-q^{22} t^{23}
-q^{24} t^{23}-q^{25} t^{23}-2 q^{26} t^{23}-2 q^{27} t^{23}
-4 q^{28} t^{23}-3 q^{29} t^{23}-6 q^{30} t^{23}-5 q^{31} t^{23}
-7 q^{32} t^{23}-3 q^{33} t^{23}-2 q^{34} t^{23}+11 q^{35} t^{23}
+15 q^{36} t^{23}+11 q^{37} t^{23}-q^{23} t^{24}-q^{25} t^{24}
-q^{26} t^{24}-2 q^{27} t^{24}-2 q^{28} t^{24}-4 q^{29} t^{24}
-3 q^{30} t^{24}-6 q^{31} t^{24}-5 q^{32} t^{24}-6 q^{33} t^{24}
-q^{34} t^{24}+2 q^{35} t^{24}+15 q^{36} t^{24}+13 q^{37} t^{24}
+2 q^{38} t^{24}-q^{24} t^{25}-q^{26} t^{25}-q^{27} t^{25}
-2 q^{28} t^{25}-2 q^{29} t^{25}-4 q^{30} t^{25}-3 q^{31} t^{25}
-6 q^{32} t^{25}-4 q^{33} t^{25}-5 q^{34} t^{25}+2 q^{35} t^{25}
+7 q^{36} t^{25}+16 q^{37} t^{25}+4 q^{38} t^{25}-q^{25} t^{26}
-q^{27} t^{26}-q^{28} t^{26}-2 q^{29} t^{26}-2 q^{30} t^{26}
-4 q^{31} t^{26}-3 q^{32} t^{26}-6 q^{33} t^{26}-3 q^{34} t^{26}
-3 q^{35} t^{26}+6 q^{36} t^{26}+10 q^{37} t^{26}+10 q^{38} t^{26}
-q^{26} t^{27}-q^{28} t^{27}-q^{29} t^{27}-2 q^{30} t^{27}
-2 q^{31} t^{27}-4 q^{32} t^{27}-3 q^{33} t^{27}-5 q^{34} t^{27}
-2 q^{35} t^{27}+10 q^{37} t^{27}+9 q^{38} t^{27}+2 q^{39} t^{27}
-q^{27} t^{28}-q^{29} t^{28}-q^{30} t^{28}-2 q^{31} t^{28}
-2 q^{32} t^{28}-4 q^{33} t^{28}-3 q^{34} t^{28}-4 q^{35} t^{28}
+4 q^{37} t^{28}+11 q^{38} t^{28}+3 q^{39} t^{28}-q^{28} t^{29}
-q^{30} t^{29}-q^{31} t^{29}-2 q^{32} t^{29}-2 q^{33} t^{29}
-4 q^{34} t^{29}-2 q^{35} t^{29}-3 q^{36} t^{29}+3 q^{37} t^{29}
+7 q^{38} t^{29}+6 q^{39} t^{29}-q^{29} t^{30}-q^{31} t^{30}
-q^{32} t^{30}-2 q^{33} t^{30}-2 q^{34} t^{30}-4 q^{35} t^{30}
-q^{36} t^{30}-q^{37} t^{30}+6 q^{38} t^{30}+6 q^{39} t^{30}
+q^{40} t^{30}-q^{30} t^{31}-q^{32} t^{31}-q^{33} t^{31}
-2 q^{34} t^{31}-2 q^{35} t^{31}-3 q^{36} t^{31}+2 q^{38} t^{31}
+7 q^{39} t^{31}+q^{40} t^{31}-q^{31} t^{32}-q^{33} t^{32}
-q^{34} t^{32}-2 q^{35} t^{32}-2 q^{36} t^{32}-2 q^{37} t^{32}
+2 q^{38} t^{32}+4 q^{39} t^{32}+3 q^{40} t^{32}-q^{32} t^{33}
-q^{34} t^{33}-q^{35} t^{33}-2 q^{36} t^{33}-q^{37} t^{33}
-q^{38} t^{33}+4 q^{39} t^{33}+3 q^{40} t^{33}-q^{33} t^{34}
-q^{35} t^{34}-q^{36} t^{34}-2 q^{37} t^{34}+q^{39} t^{34}
+4 q^{40} t^{34}-q^{34} t^{35}-q^{36} t^{35}-q^{37} t^{35}
-q^{38} t^{35}+q^{39} t^{35}+2 q^{40} t^{35}+q^{41} t^{35}
-q^{35} t^{36}-q^{37} t^{36}-q^{38} t^{36}+2 q^{40} t^{36}
+q^{41} t^{36}-q^{36} t^{37}-q^{38} t^{37}+q^{40} t^{37}
+q^{41} t^{37}-q^{37} t^{38}-q^{39} t^{38}+q^{40} t^{38}
+q^{41} t^{38}-q^{38} t^{39}+q^{41} t^{39}-q^{39} t^{40}
+q^{41} t^{40}-q^{40} t^{41}+q^{41} t^{41}-q^{41} t^{42}
+q^{42} t^{42}+a^7 \bigl(q^{28}-q^{28} t+q^{29} t+q^{30} t
+q^{31} t-q^{29} t^2+2 q^{32} t^2+q^{33} t^2+q^{34} t^2-q^{30} t^3
-q^{32} t^3+q^{33} t^3+q^{34} t^3+2 q^{35} t^3-q^{31} t^4-q^{33}t^4
+3 q^{36} t^4-q^{32} t^5-q^{34} t^5-q^{36} t^5+3 q^{37} t^5
-q^{33} t^6-q^{35} t^6-q^{37} t^6+3 q^{38} t^6-q^{34} t^7
-q^{36} t^7+2 q^{39} t^7-q^{35} t^8-q^{37} t^8+q^{39} t^8
+q^{40} t^8-q^{36} t^9-q^{38} t^9+q^{39} t^9+q^{40} t^9
-q^{37} t^{10}-q^{39} t^{10}+2 q^{40} t^{10}-q^{38} t^{11}
+q^{41} t^{11}-q^{39} t^{12}+q^{41} t^{12}-q^{40} t^{13}
+q^{41} t^{13}-q^{41} t^{14}+q^{42} t^{14}\bigr)
+a^6 \bigl(q^{21}+q^{22}+q^{23}+q^{24}+q^{25}+q^{26}
+q^{27}-q^{21} t+q^{23} t+2 q^{24} t+3 q^{25} t+3 q^{26} t
+3 q^{27} t+4 q^{28} t+2 q^{29} t+q^{30} t-q^{22} t^2-q^{23} t^2
-q^{24} t^2+2 q^{26} t^2+4 q^{27} t^2+5 q^{28} t^2+7 q^{29} t^2
+6 q^{30} t^2+5 q^{31} t^2+2 q^{32} t^2+q^{33} t^2-q^{23} t^3
-q^{24} t^3-2 q^{25} t^3-2 q^{26} t^3-q^{27} t^3+q^{28} t^3
+4 q^{29} t^3+8 q^{30} t^3+8 q^{31} t^3+8 q^{32} t^3+5 q^{33} t^3
+3 q^{34} t^3-q^{24} t^4-q^{25} t^4-2 q^{26} t^4-3 q^{27} t^4
-3 q^{28} t^4-2 q^{29} t^4+6 q^{31} t^4+8 q^{32} t^4+9 q^{33} t^4
+7 q^{34} t^4+5 q^{35} t^4-q^{25} t^5-q^{26} t^5-2 q^{27} t^5
-3 q^{28} t^5-4 q^{29} t^5-4 q^{30} t^5-3 q^{31} t^5+3 q^{32} t^5
+7 q^{33} t^5+9 q^{34} t^5+7 q^{35} t^5+6 q^{36} t^5-q^{26} t^6
-q^{27} t^6-2 q^{28} t^6-3 q^{29} t^6-4 q^{30} t^6-5 q^{31} t^6
-5 q^{32} t^6+5 q^{34} t^6+9 q^{35} t^6+7 q^{36} t^6+6 q^{37} t^6
-q^{27} t^7-q^{28} t^7-2 q^{29} t^7-3 q^{30} t^7-4 q^{31} t^7
-5 q^{32} t^7-6 q^{33} t^7-2 q^{34} t^7+5 q^{35} t^7+9 q^{36} t^7
+7 q^{37} t^7+5 q^{38} t^7-q^{28} t^8-q^{29} t^8-2 q^{30} t^8
-3 q^{31} t^8-4 q^{32} t^8-5 q^{33} t^8-6 q^{34} t^8
-2 q^{35} t^8+5 q^{36} t^8+9 q^{37} t^8+7 q^{38} t^8+3 q^{39} t^8
-q^{29} t^9-q^{30} t^9-2 q^{31} t^9-3 q^{32} t^9-4 q^{33} t^9
-5 q^{34} t^9-6 q^{35} t^9+7 q^{37} t^9+9 q^{38} t^9+5 q^{39} t^9
+q^{40} t^9-q^{30} t^{10}-q^{31} t^{10}-2 q^{32} t^{10}
-3 q^{33} t^{10}-4 q^{34} t^{10}-5 q^{35} t^{10}-5 q^{36} t^{10}
+3 q^{37} t^{10}+8 q^{38} t^{10}+8 q^{39} t^{10}+2 q^{40} t^{10}
-q^{31} t^{11}-q^{32} t^{11}-2 q^{33} t^{11}-3 q^{34} t^{11}
-4 q^{35} t^{11}-5 q^{36} t^{11}-3 q^{37} t^{11}+6 q^{38} t^{11}
+8 q^{39} t^{11}+5 q^{40} t^{11}-q^{32} t^{12}-q^{33} t^{12}
-2 q^{34} t^{12}-3 q^{35} t^{12}-4 q^{36} t^{12}-4 q^{37} t^{12}
+8 q^{39} t^{12}+6 q^{40} t^{12}+q^{41} t^{12}-q^{33} t^{13}
-q^{34} t^{13}-2 q^{35} t^{13}-3 q^{36} t^{13}-4 q^{37} t^{13}
-2 q^{38} t^{13}+4 q^{39} t^{13}+7 q^{40} t^{13}+2 q^{41} t^{13}
-q^{34} t^{14}-q^{35} t^{14}-2 q^{36} t^{14}-3 q^{37} t^{14}
-3 q^{38} t^{14}+q^{39} t^{14}+5 q^{40} t^{14}+4 q^{41} t^{14}
-q^{35} t^{15}-q^{36} t^{15}-2 q^{37} t^{15}-3 q^{38} t^{15}
-q^{39} t^{15}+4 q^{40} t^{15}+3 q^{41} t^{15}+q^{42} t^{15}
-q^{36} t^{16}-q^{37} t^{16}-2 q^{38} t^{16}-2 q^{39} t^{16}
+2 q^{40} t^{16}+3 q^{41} t^{16}+q^{42} t^{16}-q^{37} t^{17}
-q^{38} t^{17}-2 q^{39} t^{17}+3 q^{41} t^{17}+q^{42} t^{17}
-q^{38} t^{18}-q^{39} t^{18}-q^{40} t^{18}+2 q^{41} t^{18}
+q^{42} t^{18}-q^{39} t^{19}-q^{40} t^{19}+q^{41} t^{19}
+q^{42} t^{19}-q^{40} t^{20}+q^{42} t^{20}-q^{41} t^{21}
+q^{42} t^{21}\bigr)+
\)

\smallskip
\noindent
\(
a^5 \bigl(q^{15}+q^{16}+2 q^{17}+2 q^{18}+3 q^{19}+3 q^{20}+3q^{21}
+2 q^{22}+2 q^{23}+q^{24}+q^{25}-q^{15} t+2 q^{18} t+3 q^{19} t
+6 q^{20} t+8 q^{21} t+11 q^{22} t+10 q^{23} t+10 q^{24} t
+7 q^{25} t+6 q^{26} t+3 q^{27} t+q^{28} t-q^{16} t^2-q^{17} t^2
-2 q^{18} t^2-q^{19} t^2-q^{20} t^2+3 q^{21} t^2+6 q^{22} t^2
+14 q^{23} t^2+17 q^{24} t^2+22 q^{25} t^2+19 q^{26} t^2
+18 q^{27} t^2+12 q^{28} t^2+8 q^{29} t^2+3 q^{30} t^2+q^{31} t^2
-q^{17} t^3-q^{18} t^3-3 q^{19} t^3-3 q^{20} t^3-5 q^{21} t^3
-3 q^{22} t^3-2 q^{23} t^3+7 q^{24} t^3+13 q^{25} t^3
+26 q^{26} t^3+28 q^{27} t^3+32 q^{28} t^3+25 q^{29} t^3
+19 q^{30} t^3+10 q^{31} t^3+5 q^{32} t^3+q^{33} t^3-q^{18} t^4
-q^{19} t^4-3 q^{20} t^4-4 q^{21} t^4-7 q^{22} t^4-7 q^{23} t^4
-9 q^{24} t^4-3 q^{25} t^4+q^{26} t^4+18 q^{27} t^4+26 q^{28} t^4
+39 q^{29} t^4+35 q^{30} t^4+30 q^{31} t^4+17 q^{32} t^4
+10 q^{33} t^4+2 q^{34} t^4-q^{19} t^5-q^{20} t^5-3 q^{21} t^5
-4 q^{22} t^5-8 q^{23} t^5-9 q^{24} t^5-13 q^{25} t^5-10 q^{26} t^5
-10 q^{27} t^5+5 q^{28} t^5+15 q^{29} t^5+38 q^{30} t^5
+39 q^{31} t^5+38 q^{32} t^5+23 q^{33} t^5+14 q^{34} t^5
+3 q^{35} t^5-q^{20} t^6-q^{21} t^6-3 q^{22} t^6-4 q^{23} t^6
-8 q^{24} t^6-10 q^{25} t^6-15 q^{26} t^6-14 q^{27} t^6
-17 q^{28} t^6-6 q^{29} t^6+q^{30} t^6+31 q^{31} t^6+37 q^{32} t^6
+41 q^{33} t^6+26 q^{34} t^6+16 q^{35} t^6+3 q^{36} t^6-q^{21} t^7
-q^{22} t^7-3 q^{23} t^7-4 q^{24} t^7-8 q^{25} t^7-10 q^{26} t^7
-16 q^{27} t^7-16 q^{28} t^7-21 q^{29} t^7-13 q^{30} t^7
-10 q^{31} t^7+23 q^{32} t^7+34 q^{33} t^7+42 q^{34} t^7
+27 q^{35} t^7+16 q^{36} t^7+3 q^{37} t^7-q^{22} t^8-q^{23} t^8
-3 q^{24} t^8-4 q^{25} t^8-8 q^{26} t^8-10 q^{27} t^8-16 q^{28} t^8
-17 q^{29} t^8-23 q^{30} t^8-17 q^{31} t^8-16 q^{32} t^8
+18 q^{33} t^8+32 q^{34} t^8+42 q^{35} t^8+26 q^{36} t^8
+14 q^{37} t^8+2 q^{38} t^8-q^{23} t^9-q^{24} t^9-3 q^{25} t^9
-4 q^{26} t^9-8 q^{27} t^9-10 q^{28} t^9-16 q^{29} t^9-17 q^{30}t^9
-24 q^{31} t^9-19 q^{32} t^9-18 q^{33} t^9+18 q^{34} t^9
+34 q^{35} t^9+41 q^{36} t^9+23 q^{37} t^9+10 q^{38} t^9+q^{39} t^9
-q^{24} t^{10}-q^{25} t^{10}-3 q^{26} t^{10}-4 q^{27} t^{10}
-8 q^{28} t^{10}-10 q^{29} t^{10}-16 q^{30} t^{10}-17 q^{31} t^{10}
-24 q^{32} t^{10}-19 q^{33} t^{10}-16 q^{34} t^{10}+23 q^{35}t^{10}
+37 q^{36} t^{10}+38 q^{37} t^{10}+17 q^{38} t^{10}+5 q^{39} t^{10}
-q^{25} t^{11}-q^{26} t^{11}-3 q^{27} t^{11}-4 q^{28} t^{11}
-8 q^{29} t^{11}-10 q^{30} t^{11}-16 q^{31} t^{11}-17 q^{32} t^{11}
-24 q^{33} t^{11}-17 q^{34} t^{11}-10 q^{35} t^{11}+31 q^{36}t^{11}
+39 q^{37} t^{11}+30 q^{38} t^{11}+10 q^{39} t^{11}+q^{40} t^{11}
-q^{26} t^{12}-q^{27} t^{12}-3 q^{28} t^{12}-4 q^{29} t^{12}
-8 q^{30} t^{12}-10 q^{31} t^{12}-16 q^{32} t^{12}-17 q^{33} t^{12}
-23 q^{34} t^{12}-13 q^{35} t^{12}+q^{36} t^{12}+38 q^{37} t^{12}
+35 q^{38} t^{12}+19 q^{39} t^{12}+3 q^{40} t^{12}-q^{27} t^{13}
-q^{28} t^{13}-3 q^{29} t^{13}-4 q^{30} t^{13}-8 q^{31} t^{13}
-10 q^{32} t^{13}-16 q^{33} t^{13}-17 q^{34} t^{13}-21 q^{35}t^{13}
-6 q^{36} t^{13}+15 q^{37} t^{13}+39 q^{38} t^{13}+25 q^{39} t^{13}
+8 q^{40} t^{13}-q^{28} t^{14}-q^{29} t^{14}-3 q^{30} t^{14}
-4 q^{31} t^{14}-8 q^{32} t^{14}-10 q^{33} t^{14}-16 q^{34} t^{14}
-16 q^{35} t^{14}-17 q^{36} t^{14}+5 q^{37} t^{14}+26 q^{38} t^{14}
+32 q^{39} t^{14}+12 q^{40} t^{14}+q^{41} t^{14}-q^{29} t^{15}
-q^{30} t^{15}-3 q^{31} t^{15}-4 q^{32} t^{15}-8 q^{33} t^{15}
-10 q^{34} t^{15}-16 q^{35} t^{15}-14 q^{36} t^{15}-10 q^{37}t^{15}
+18 q^{38} t^{15}+28 q^{39} t^{15}+18 q^{40} t^{15}+3 q^{41} t^{15}
-q^{30} t^{16}-q^{31} t^{16}-3 q^{32} t^{16}-4 q^{33} t^{16}
-8 q^{34} t^{16}-10 q^{35} t^{16}-15 q^{36} t^{16}-10 q^{37} t^{16}
+q^{38} t^{16}+26 q^{39} t^{16}+19 q^{40} t^{16}+6 q^{41} t^{16}
-q^{31} t^{17}-q^{32} t^{17}-3 q^{33} t^{17}-4 q^{34} t^{17}
-8 q^{35} t^{17}-10 q^{36} t^{17}-13 q^{37} t^{17}-3 q^{38} t^{17}
+13 q^{39} t^{17}+22 q^{40} t^{17}+7 q^{41} t^{17}+q^{42} t^{17}
-q^{32} t^{18}-q^{33} t^{18}-3 q^{34} t^{18}-4 q^{35} t^{18}
-8 q^{36} t^{18}-9 q^{37} t^{18}-9 q^{38} t^{18}+7 q^{39} t^{18}
+17 q^{40} t^{18}+10 q^{41} t^{18}+q^{42} t^{18}-q^{33} t^{19}
-q^{34} t^{19}-3 q^{35} t^{19}-4 q^{36} t^{19}-8 q^{37} t^{19}
-7 q^{38} t^{19}-2 q^{39} t^{19}+14 q^{40} t^{19}+10 q^{41} t^{19}
+2 q^{42} t^{19}-q^{34} t^{20}-q^{35} t^{20}-3 q^{36} t^{20}
-4 q^{37} t^{20}-7 q^{38} t^{20}-3 q^{39} t^{20}+6 q^{40} t^{20}
+11 q^{41} t^{20}+2 q^{42} t^{20}-q^{35} t^{21}-q^{36} t^{21}
-3 q^{37} t^{21}-4 q^{38} t^{21}-5 q^{39} t^{21}+3 q^{40} t^{21}
+8 q^{41} t^{21}+3 q^{42} t^{21}-q^{36} t^{22}-q^{37} t^{22}
-3 q^{38} t^{22}-3 q^{39} t^{22}-q^{40} t^{22}+6 q^{41} t^{22}
+3 q^{42} t^{22}-q^{37} t^{23}-q^{38} t^{23}-3 q^{39} t^{23}
-q^{40} t^{23}+3 q^{41} t^{23}+3 q^{42} t^{23}-q^{38} t^{24}
-q^{39} t^{24}-2 q^{40} t^{24}+2 q^{41} t^{24}+2 q^{42} t^{24}
-q^{39} t^{25}-q^{40} t^{25}+2 q^{42} t^{25}
-q^{40} t^{26}+q^{42} t^{26}-q^{41} t^{27}+q^{42} t^{27}\bigr)+
\)

\smallskip
\noindent
\(
a^4 \bigl(q^{10}+q^{11}+2 q^{12}+3 q^{13}+4 q^{14}+4 q^{15}
+5 q^{16}+4 q^{17}+4 q^{18}+3 q^{19}+2 q^{20}+q^{21}+q^{22}
-q^{10} t+q^{13} t+3 q^{14} t+7 q^{15} t+10 q^{16} t+15 q^{17} t
+17 q^{18} t+19 q^{19} t+18 q^{20} t+16 q^{21} t+11 q^{22} t
+8 q^{23} t+4 q^{24} t+2 q^{25} t-q^{11} t^2-q^{12} t^2
-2 q^{13} t^2-2 q^{14} t^2-2 q^{15} t^2+q^{16} t^2+5 q^{17} t^2
+14 q^{18} t^2+22 q^{19} t^2+33 q^{20} t^2+38 q^{21} t^2
+42 q^{22} t^2+37 q^{23} t^2+31 q^{24} t^2+20 q^{25} t^2
+13 q^{26} t^2+5 q^{27} t^2+2 q^{28} t^2-q^{12} t^3-q^{13} t^3
-3 q^{14} t^3-4 q^{15} t^3-6 q^{16} t^3-6 q^{17} t^3-6 q^{18} t^3
+q^{19} t^3+9 q^{20} t^3+26 q^{21} t^3+42 q^{22} t^3+59 q^{23} t^3
+64 q^{24} t^3+65 q^{25} t^3+51 q^{26} t^3+38 q^{27} t^3
+21 q^{28} t^3+11 q^{29} t^3+3 q^{30} t^3+q^{31} t^3-q^{13} t^4
-q^{14} t^4-3 q^{15} t^4-5 q^{16} t^4-8 q^{17} t^4-10 q^{18} t^4
-14 q^{19} t^4-12 q^{20} t^4-9 q^{21} t^4+4 q^{22} t^4
+22 q^{23} t^4+51 q^{24} t^4+71 q^{25} t^4+89 q^{26} t^4
+83 q^{27} t^4+70 q^{28} t^4+44 q^{29} t^4+26 q^{30} t^4
+9 q^{31} t^4+3 q^{32} t^4-q^{14} t^5-q^{15} t^5-3 q^{16} t^5
-5 q^{17} t^5-9 q^{18} t^5-12 q^{19} t^5-18 q^{20} t^5
-20 q^{21} t^5-23 q^{22} t^5-16 q^{23} t^5-5 q^{24} t^5
+23 q^{25} t^5+53 q^{26} t^5+89 q^{27} t^5+101 q^{28} t^5
+99 q^{29} t^5+68 q^{30} t^5+43 q^{31} t^5+17 q^{32} t^5
+6 q^{33} t^5-q^{15} t^6-q^{16} t^6-3 q^{17} t^6-5 q^{18} t^6
-9 q^{19} t^6-13 q^{20} t^6-20 q^{21} t^6-24 q^{22} t^6
-31 q^{23} t^6-30 q^{24} t^6-26 q^{25} t^6-6 q^{26} t^6
+22 q^{27} t^6+69 q^{28} t^6+99 q^{29} t^6+115 q^{30} t^6
+86 q^{31} t^6+57 q^{32} t^6+23 q^{33} t^6+8 q^{34} t^6-q^{16} t^7
-q^{17} t^7-3 q^{18} t^7-5 q^{19} t^7-9 q^{20} t^7-13 q^{21} t^7
-21 q^{22} t^7-26 q^{23} t^7-35 q^{24} t^7-38 q^{25} t^7
-40 q^{26} t^7-27 q^{27} t^7-7 q^{28} t^7+41 q^{29} t^7
+85 q^{30} t^7+118 q^{31} t^7+96 q^{32} t^7+65 q^{33} t^7
+26 q^{34} t^7+9 q^{35} t^7-q^{17} t^8-q^{18} t^8-3 q^{19} t^8
-5 q^{20} t^8-9 q^{21} t^8-13 q^{22} t^8-21 q^{23} t^8-27 q^{24}t^8
-37 q^{25} t^8-42 q^{26} t^8-48 q^{27} t^8-41 q^{28} t^8
-28 q^{29} t^8+17 q^{30} t^8+69 q^{31} t^8+116 q^{32} t^8
+100 q^{33} t^8+68 q^{34} t^8+26 q^{35} t^8+8 q^{36} t^8-q^{18} t^9
-q^{19} t^9-3 q^{20} t^9-5 q^{21} t^9-9 q^{22} t^9-13 q^{23} t^9
-21 q^{24} t^9-27 q^{25} t^9-38 q^{26} t^9-44 q^{27} t^9
-52 q^{28} t^9-49 q^{29} t^9-41 q^{30} t^9+q^{31} t^9+59 q^{32} t^9
+114 q^{33} t^9+100 q^{34} t^9+65 q^{35} t^9+23 q^{36} t^9
+6 q^{37} t^9-q^{19} t^{10}-q^{20} t^{10}-3 q^{21} t^{10}
-5 q^{22} t^{10}-9 q^{23} t^{10}-13 q^{24} t^{10}-21 q^{25} t^{10}
-27 q^{26} t^{10}-38 q^{27} t^{10}-45 q^{28} t^{10}-54 q^{29}t^{10}
-53 q^{30} t^{10}-47 q^{31} t^{10}-4 q^{32} t^{10}+59 q^{33} t^{10}
+116 q^{34} t^{10}+96 q^{35} t^{10}+57 q^{36} t^{10}
+17 q^{37} t^{10}+3 q^{38} t^{10}-q^{20} t^{11}-q^{21} t^{11}
-3 q^{22} t^{11}-5 q^{23} t^{11}-9 q^{24} t^{11}-13 q^{25} t^{11}
-21 q^{26} t^{11}-27 q^{27} t^{11}-38 q^{28} t^{11}-45 q^{29}t^{11}
-55 q^{30} t^{11}-54 q^{31} t^{11}-47 q^{32} t^{11}+q^{33} t^{11}
+69 q^{34} t^{11}+118 q^{35} t^{11}+86 q^{36} t^{11}
+43 q^{37} t^{11}+9 q^{38} t^{11}+q^{39} t^{11}-q^{21} t^{12}
-q^{22} t^{12}-3 q^{23} t^{12}-5 q^{24} t^{12}-9 q^{25} t^{12}
-13 q^{26} t^{12}-21 q^{27} t^{12}-27 q^{28} t^{12}-38 q^{29}t^{12}
-45 q^{30} t^{12}-55 q^{31} t^{12}-53 q^{32} t^{12}-41 q^{33}t^{12}
+17 q^{34} t^{12}+85 q^{35} t^{12}+115 q^{36} t^{12}
+68 q^{37} t^{12}+26 q^{38} t^{12}+3 q^{39} t^{12}-q^{22} t^{13}
-q^{23} t^{13}-3 q^{24} t^{13}-5 q^{25} t^{13}-9 q^{26} t^{13}
-13 q^{27} t^{13}-21 q^{28} t^{13}-27 q^{29} t^{13}-38 q^{30}t^{13}
-45 q^{31} t^{13}-54 q^{32} t^{13}-49 q^{33} t^{13}-28 q^{34}t^{13}
+41 q^{35} t^{13}+99 q^{36} t^{13}+99 q^{37} t^{13}+44 q^{38}t^{13}
+11 q^{39} t^{13}-q^{23} t^{14}-q^{24} t^{14}-3 q^{25} t^{14}
-5 q^{26} t^{14}-9 q^{27} t^{14}-13 q^{28} t^{14}-21 q^{29} t^{14}
-27 q^{30} t^{14}-38 q^{31} t^{14}-45 q^{32} t^{14}-52 q^{33}t^{14}
-41 q^{34} t^{14}-7 q^{35} t^{14}+69 q^{36} t^{14}+101 q^{37}t^{14}
+70 q^{38} t^{14}+21 q^{39} t^{14}+2 q^{40} t^{14}-q^{24} t^{15}
-q^{25} t^{15}-3 q^{26} t^{15}-5 q^{27} t^{15}-9 q^{28} t^{15}
-13 q^{29} t^{15}-21 q^{30} t^{15}-27 q^{31} t^{15}-38 q^{32}t^{15}
-44 q^{33} t^{15}-48 q^{34} t^{15}-27 q^{35} t^{15}+22 q^{36}t^{15}
+89 q^{37} t^{15}+83 q^{38} t^{15}+38 q^{39} t^{15}+5 q^{40} t^{15}
-q^{25} t^{16}-q^{26} t^{16}-3 q^{27} t^{16}-5 q^{28} t^{16}
-9 q^{29} t^{16}-13 q^{30} t^{16}-21 q^{31} t^{16}-27 q^{32} t^{16}
-38 q^{33} t^{16}-42 q^{34} t^{16}-40 q^{35} t^{16}-6 q^{36} t^{16}
+53 q^{37} t^{16}+89 q^{38} t^{16}+51 q^{39} t^{16}+13 q^{40}t^{16}
-q^{26} t^{17}-q^{27} t^{17}-3 q^{28} t^{17}-5 q^{29} t^{17}
-9 q^{30} t^{17}-13 q^{31} t^{17}-21 q^{32} t^{17}-27 q^{33} t^{17}
-37 q^{34} t^{17}-38 q^{35} t^{17}-26 q^{36} t^{17}+23 q^{37}t^{17}
+71 q^{38} t^{17}+65 q^{39} t^{17}+20 q^{40} t^{17}+2 q^{41} t^{17}
-q^{27} t^{18}-q^{28} t^{18}-3 q^{29} t^{18}-5 q^{30} t^{18}
-9 q^{31} t^{18}-13 q^{32} t^{18}-21 q^{33} t^{18}-27 q^{34} t^{18}
-35 q^{35} t^{18}-30 q^{36} t^{18}-5 q^{37} t^{18}+51 q^{38} t^{18}
+64 q^{39} t^{18}+31 q^{40} t^{18}+4 q^{41} t^{18}-q^{28} t^{19}
-q^{29} t^{19}-3 q^{30} t^{19}-5 q^{31} t^{19}-9 q^{32} t^{19}
-13 q^{33} t^{19}-21 q^{34} t^{19}-26 q^{35} t^{19}-31 q^{36}t^{19}
-16 q^{37} t^{19}+22 q^{38} t^{19}+59 q^{39} t^{19}+37 q^{40}t^{19}
+8 q^{41} t^{19}-q^{29} t^{20}-q^{30} t^{20}-3 q^{31} t^{20}
-5 q^{32} t^{20}-9 q^{33} t^{20}-13 q^{34} t^{20}-21 q^{35} t^{20}
-24 q^{36} t^{20}-23 q^{37} t^{20}+4 q^{38} t^{20}+42 q^{39} t^{20}
+42 q^{40} t^{20}+11 q^{41} t^{20}+q^{42} t^{20}-q^{30} t^{21}
-q^{31} t^{21}-3 q^{32} t^{21}-5 q^{33} t^{21}-9 q^{34} t^{21}
-13 q^{35} t^{21}-20 q^{36} t^{21}-20 q^{37} t^{21}-9 q^{38} t^{21}
+26 q^{39} t^{21}+38 q^{40} t^{21}+16 q^{41} t^{21}+q^{42} t^{21}
-q^{31} t^{22}-q^{32} t^{22}-3 q^{33} t^{22}-5 q^{34} t^{22}
-9 q^{35} t^{22}-13 q^{36} t^{22}-18 q^{37} t^{22}-12 q^{38} t^{22}
+9 q^{39} t^{22}+33 q^{40} t^{22}+18 q^{41} t^{22}+2 q^{42} t^{22}
-q^{32} t^{23}-q^{33} t^{23}-3 q^{34} t^{23}-5 q^{35} t^{23}
-9 q^{36} t^{23}-12 q^{37} t^{23}-14 q^{38} t^{23}+q^{39} t^{23}
+22 q^{40} t^{23}+19 q^{41} t^{23}+3 q^{42} t^{23}-q^{33} t^{24}
-q^{34} t^{24}-3 q^{35} t^{24}-5 q^{36} t^{24}-9 q^{37} t^{24}
-10 q^{38} t^{24}-6 q^{39} t^{24}+14 q^{40} t^{24}+17 q^{41} t^{24}
+4 q^{42} t^{24}-q^{34} t^{25}-q^{35} t^{25}-3 q^{36} t^{25}
-5 q^{37} t^{25}-8 q^{38} t^{25}-6 q^{39} t^{25}+5 q^{40} t^{25}
+15 q^{41} t^{25}+4 q^{42} t^{25}-q^{35} t^{26}-q^{36} t^{26}
-3 q^{37} t^{26}-5 q^{38} t^{26}-6 q^{39} t^{26}+q^{40} t^{26}
+10 q^{41} t^{26}+5 q^{42} t^{26}-q^{36} t^{27}-q^{37} t^{27}
-3 q^{38} t^{27}-4 q^{39} t^{27}-2 q^{40} t^{27}+7 q^{41} t^{27}
+4 q^{42} t^{27}-q^{37} t^{28}-q^{38} t^{28}-3 q^{39} t^{28}
-2 q^{40} t^{28}+3 q^{41} t^{28}+4 q^{42} t^{28}-q^{38} t^{29}
-q^{39} t^{29}-2 q^{40} t^{29}+q^{41} t^{29}+3 q^{42} t^{29}
-q^{39} t^{30}-q^{40} t^{30}+2 q^{42} t^{30}
-q^{40} t^{31}+q^{42} t^{31}-q^{41} t^{32}+q^{42} t^{32}\bigr)+
\)

\smallskip
\noindent
\(
a^3 \bigl(q^6+q^7+2 q^8+3 q^9+4 q^{10}+4 q^{11}+5 q^{12}+4 q^{13}
+4 q^{14}+3 q^{15}+2 q^{16}+q^{17}+q^{18}-q^6 t+q^9 t+3 q^{10} t
+7 q^{11} t+10 q^{12} t+16 q^{13} t+18 q^{14} t+21 q^{15} t
+21 q^{16} t+19 q^{17} t+14 q^{18} t+11 q^{19} t+6 q^{20} t
+3 q^{21} t+q^{22} t-q^7 t^2-q^8 t^2-2 q^9 t^2-2 q^{10} t^2
-2 q^{11} t^2+q^{12} t^2+4 q^{13} t^2+14 q^{14} t^2+22 q^{15} t^2
+35 q^{16} t^2+43 q^{17} t^2+51 q^{18} t^2+48 q^{19} t^2
+44 q^{20} t^2+32 q^{21} t^2+22 q^{22} t^2+11 q^{23} t^2
+5 q^{24} t^2+q^{25} t^2-q^8 t^3-q^9 t^3-3 q^{10} t^3-4 q^{11} t^3
-6 q^{12} t^3-6 q^{13} t^3-7 q^{14} t^3+6 q^{16} t^3+24 q^{17} t^3
+41 q^{18} t^3+65 q^{19} t^3+78 q^{20} t^3+89 q^{21} t^3
+79 q^{22} t^3+66 q^{23} t^3+44 q^{24} t^3+26 q^{25} t^3
+11 q^{26} t^3+4 q^{27} t^3+q^{28} t^3-q^9 t^4-q^{10} t^4
-3 q^{11} t^4-5 q^{12} t^4-8 q^{13} t^4-10 q^{14} t^4-15 q^{15} t^4
-13 q^{16} t^4-13 q^{17} t^4-q^{18} t^4+14 q^{19} t^4
+45 q^{20} t^4+72 q^{21} t^4+110 q^{22} t^4+120 q^{23} t^4
+120 q^{24} t^4+96 q^{25} t^4+68 q^{26} t^4+35 q^{27} t^4
+16 q^{28} t^4+5 q^{29} t^4+q^{30} t^4-q^{10} t^5-q^{11} t^5
-3 q^{12} t^5-5 q^{13} t^5-9 q^{14} t^5-12 q^{15} t^5-19 q^{16} t^5
-21 q^{17} t^5-27 q^{18} t^5-22 q^{19} t^5-16 q^{20} t^5
+9 q^{21} t^5+35 q^{22} t^5+88 q^{23} t^5+124 q^{24} t^5
+155 q^{25} t^5+145 q^{26} t^5+120 q^{27} t^5+72 q^{28} t^5
+36 q^{29} t^5+13 q^{30} t^5+3 q^{31} t^5-q^{11} t^6-q^{12} t^6
-3 q^{13} t^6-5 q^{14} t^6-9 q^{15} t^6-13 q^{16} t^6-21 q^{17} t^6
-25 q^{18} t^6-35 q^{19} t^6-36 q^{20} t^6-38 q^{21} t^6
-23 q^{22} t^6-6 q^{23} t^6+44 q^{24} t^6+89 q^{25} t^6
+153 q^{26} t^6+171 q^{27} t^6+164 q^{28} t^6+110 q^{29} t^6
+60 q^{30} t^6+23 q^{31} t^6+6 q^{32} t^6-q^{12} t^7-q^{13} t^7
-3 q^{14} t^7-5 q^{15} t^7-9 q^{16} t^7-13 q^{17} t^7-22 q^{18} t^7
-27 q^{19} t^7-39 q^{20} t^7-44 q^{21} t^7-52 q^{22} t^7
-45 q^{23} t^7-39 q^{24} t^7+2 q^{25} t^7+42 q^{26} t^7
+122 q^{27} t^7+168 q^{28} t^7+190 q^{29} t^7+138 q^{30} t^7
+79 q^{31} t^7+31 q^{32} t^7+8 q^{33} t^7-q^{13} t^8-q^{14} t^8
-3 q^{15} t^8-5 q^{16} t^8-9 q^{17} t^8-13 q^{18} t^8-22 q^{19} t^8
-28 q^{20} t^8-41 q^{21} t^8-48 q^{22} t^8-60 q^{23} t^8
-59 q^{24} t^8-61 q^{25} t^8-31 q^{26} t^8-q^{27} t^8+83 q^{28} t^8
+148 q^{29} t^8+200 q^{30} t^8+156 q^{31} t^8+92 q^{32} t^8
+35 q^{33} t^8+9 q^{34} t^8-q^{14} t^9-q^{15} t^9-3 q^{16} t^9
-5 q^{17} t^9-9 q^{18} t^9-13 q^{19} t^9-22 q^{20} t^9-28 q^{21}t^9
-42 q^{22} t^9-50 q^{23} t^9-64 q^{24} t^9-67 q^{25} t^9
-75 q^{26} t^9-53 q^{27} t^9-33 q^{28} t^9+49 q^{29} t^9
+125 q^{30} t^9+199 q^{31} t^9+164 q^{32} t^9+96 q^{33} t^9
+35 q^{34} t^9+8 q^{35} t^9-q^{15} t^{10}-q^{16} t^{10}
-3 q^{17} t^{10}-5 q^{18} t^{10}-9 q^{19} t^{10}-13 q^{20} t^{10}
-22 q^{21} t^{10}-28 q^{22} t^{10}-42 q^{23} t^{10}-51 q^{24}t^{10}
-66 q^{25} t^{10}-71 q^{26} t^{10}-83 q^{27} t^{10}-67 q^{28}t^{10}
-53 q^{29} t^{10}+27 q^{30} t^{10}+111 q^{31} t^{10}
+198 q^{32} t^{10}+164 q^{33} t^{10}+92 q^{34} t^{10}
+31 q^{35} t^{10}+6 q^{36} t^{10}-q^{16} t^{11}-q^{17} t^{11}
-3 q^{18} t^{11}-5 q^{19} t^{11}-9 q^{20} t^{11}-13 q^{21} t^{11}
-22 q^{22} t^{11}-28 q^{23} t^{11}-42 q^{24} t^{11}-51 q^{25}t^{11}
-67 q^{26} t^{11}-73 q^{27} t^{11}-87 q^{28} t^{11}-74 q^{29}t^{11}
-63 q^{30} t^{11}+20 q^{31} t^{11}+111 q^{32} t^{11}
+199 q^{33} t^{11}+156 q^{34} t^{11}+79 q^{35} t^{11}
+23 q^{36} t^{11}+3 q^{37} t^{11}-q^{17} t^{12}-q^{18} t^{12}
-3 q^{19} t^{12}-5 q^{20} t^{12}-9 q^{21} t^{12}-13 q^{22} t^{12}
-22 q^{23} t^{12}-28 q^{24} t^{12}-42 q^{25} t^{12}-51 q^{26}t^{12}
-67 q^{27} t^{12}-74 q^{28} t^{12}-89 q^{29} t^{12}-76 q^{30}t^{12}
-63 q^{31} t^{12}+27 q^{32} t^{12}+125 q^{33} t^{12}
+200 q^{34} t^{12}+138 q^{35} t^{12}+60 q^{36} t^{12}
+13 q^{37} t^{12}+q^{38} t^{12}-q^{18} t^{13}-q^{19} t^{13}
-3 q^{20} t^{13}-5 q^{21} t^{13}-9 q^{22} t^{13}-13 q^{23} t^{13}
-22 q^{24} t^{13}-28 q^{25} t^{13}-42 q^{26} t^{13}-51 q^{27}t^{13}
-67 q^{28} t^{13}-74 q^{29} t^{13}-89 q^{30} t^{13}-74 q^{31}t^{13}
-53 q^{32} t^{13}+49 q^{33} t^{13}+148 q^{34} t^{13}
+190 q^{35} t^{13}+110 q^{36} t^{13}+36 q^{37} t^{13}
+5 q^{38} t^{13}-q^{19} t^{14}-q^{20} t^{14}-3 q^{21} t^{14}
-5 q^{22} t^{14}-9 q^{23} t^{14}-13 q^{24} t^{14}-22 q^{25} t^{14}
-28 q^{26} t^{14}-42 q^{27} t^{14}-51 q^{28} t^{14}-67 q^{29}t^{14}
-74 q^{30} t^{14}-87 q^{31} t^{14}-67 q^{32} t^{14}-33 q^{33}t^{14}
+83 q^{34} t^{14}+168 q^{35} t^{14}+164 q^{36} t^{14}
+72 q^{37} t^{14}+16 q^{38} t^{14}+q^{39} t^{14}-q^{20} t^{15}
-q^{21} t^{15}-3 q^{22} t^{15}-5 q^{23} t^{15}-9 q^{24} t^{15}
-13 q^{25} t^{15}-22 q^{26} t^{15}-28 q^{27} t^{15}-42 q^{28}t^{15}
-51 q^{29} t^{15}-67 q^{30} t^{15}-73 q^{31} t^{15}-83 q^{32}t^{15}
-53 q^{33} t^{15}-q^{34} t^{15}+122 q^{35} t^{15}+171 q^{36} t^{15}
+120 q^{37} t^{15}+35 q^{38} t^{15}+4 q^{39} t^{15}-q^{21} t^{16}
-q^{22} t^{16}-3 q^{23} t^{16}-5 q^{24} t^{16}-9 q^{25} t^{16}
-13 q^{26} t^{16}-22 q^{27} t^{16}-28 q^{28} t^{16}-42 q^{29}t^{16}
-51 q^{30} t^{16}-67 q^{31} t^{16}-71 q^{32} t^{16}-75 q^{33}t^{16}
-31 q^{34} t^{16}+42 q^{35} t^{16}+153 q^{36} t^{16}
+145 q^{37} t^{16}+68 q^{38} t^{16}+11 q^{39} t^{16}-q^{22} t^{17}
-q^{23} t^{17}-3 q^{24} t^{17}-5 q^{25} t^{17}-9 q^{26} t^{17}
-13 q^{27} t^{17}-22 q^{28} t^{17}-28 q^{29} t^{17}-42 q^{30}t^{17}
-51 q^{31} t^{17}-66 q^{32} t^{17}-67 q^{33} t^{17}-61 q^{34}t^{17}
+2 q^{35} t^{17}+89 q^{36} t^{17}+155 q^{37} t^{17}+96 q^{38}t^{17}
+26 q^{39} t^{17}+q^{40} t^{17}-q^{23} t^{18}-q^{24} t^{18}
-3 q^{25} t^{18}-5 q^{26} t^{18}-9 q^{27} t^{18}-13 q^{28} t^{18}
-22 q^{29} t^{18}-28 q^{30} t^{18}-42 q^{31} t^{18}-51 q^{32}t^{18}
-64 q^{33} t^{18}-59 q^{34} t^{18}-39 q^{35} t^{18}+44 q^{36}t^{18}
+124 q^{37} t^{18}+120 q^{38} t^{18}+44 q^{39} t^{18}
+5 q^{40} t^{18}-q^{24} t^{19}-q^{25} t^{19}-3 q^{26} t^{19}
-5 q^{27} t^{19}-9 q^{28} t^{19}-13 q^{29} t^{19}-22 q^{30} t^{19}
-28 q^{31} t^{19}-42 q^{32} t^{19}-50 q^{33} t^{19}-60 q^{34}t^{19}
-45 q^{35} t^{19}-6 q^{36} t^{19}+88 q^{37} t^{19}+120 q^{38}t^{19}
+66 q^{39} t^{19}+11 q^{40} t^{19}-q^{25} t^{20}-q^{26} t^{20}
-3 q^{27} t^{20}-5 q^{28} t^{20}-9 q^{29} t^{20}-13 q^{30} t^{20}
-22 q^{31} t^{20}-28 q^{32} t^{20}-42 q^{33} t^{20}-48 q^{34}t^{20}
-52 q^{35} t^{20}-23 q^{36} t^{20}+35 q^{37} t^{20}+110q^{38}t^{20}
+79 q^{39} t^{20}+22 q^{40} t^{20}+q^{41} t^{20}-q^{26} t^{21}
-q^{27} t^{21}-3 q^{28} t^{21}-5 q^{29} t^{21}-9 q^{30} t^{21}
-13 q^{31} t^{21}-22 q^{32} t^{21}-28 q^{33} t^{21}-41 q^{34}t^{21}
-44 q^{35} t^{21}-38 q^{36} t^{21}+9 q^{37} t^{21}+72 q^{38} t^{21}
+89 q^{39} t^{21}+32 q^{40} t^{21}+3 q^{41} t^{21}-q^{27} t^{22}
-q^{28} t^{22}-3 q^{29} t^{22}-5 q^{30} t^{22}-9 q^{31} t^{22}
-13 q^{32} t^{22}-22 q^{33} t^{22}-28 q^{34} t^{22}-39 q^{35}t^{22}
-36 q^{36} t^{22}-16 q^{37} t^{22}+45 q^{38} t^{22}+78 q^{39}t^{22}
+44 q^{40} t^{22}+6 q^{41} t^{22}-q^{28} t^{23}-q^{29} t^{23}
-3 q^{30} t^{23}-5 q^{31} t^{23}-9 q^{32} t^{23}-13 q^{33} t^{23}
-22 q^{34} t^{23}-27 q^{35} t^{23}-35 q^{36} t^{23}-22 q^{37}t^{23}
+14 q^{38} t^{23}+65 q^{39} t^{23}+48 q^{40} t^{23}+11 q^{41}t^{23}
-q^{29} t^{24}-q^{30} t^{24}-3 q^{31} t^{24}-5 q^{32} t^{24}
-9 q^{33} t^{24}-13 q^{34} t^{24}-22 q^{35} t^{24}-25 q^{36} t^{24}
-27 q^{37} t^{24}-q^{38} t^{24}+41 q^{39} t^{24}+51 q^{40} t^{24}
+14 q^{41} t^{24}+q^{42} t^{24}-q^{30} t^{25}-q^{31} t^{25}
-3 q^{32} t^{25}-5 q^{33} t^{25}-9 q^{34} t^{25}-13 q^{35} t^{25}
-21 q^{36} t^{25}-21 q^{37} t^{25}-13 q^{38} t^{25}+24 q^{39}t^{25}
+43 q^{40} t^{25}+19 q^{41} t^{25}+q^{42} t^{25}-q^{31} t^{26}
-q^{32} t^{26}-3 q^{33} t^{26}-5 q^{34} t^{26}-9 q^{35} t^{26}
-13 q^{36} t^{26}-19 q^{37} t^{26}-13 q^{38} t^{26}+6 q^{39} t^{26}
+35 q^{40} t^{26}+21 q^{41} t^{26}+2 q^{42} t^{26}-q^{32} t^{27}
-q^{33} t^{27}-3 q^{34} t^{27}-5 q^{35} t^{27}-9 q^{36} t^{27}
-12 q^{37} t^{27}-15 q^{38} t^{27}+22 q^{40} t^{27}+21 q^{41}t^{27}
+3 q^{42} t^{27}-q^{33} t^{28}-q^{34} t^{28}-3 q^{35} t^{28}
-5 q^{36} t^{28}-9 q^{37} t^{28}-10 q^{38} t^{28}-7 q^{39} t^{28}
+14 q^{40} t^{28}+18 q^{41} t^{28}+4 q^{42} t^{28}-q^{34} t^{29}
-q^{35} t^{29}-3 q^{36} t^{29}-5 q^{37} t^{29}-8 q^{38} t^{29}
-6 q^{39} t^{29}+4 q^{40} t^{29}+16 q^{41} t^{29}+4 q^{42} t^{29}
-q^{35} t^{30}-q^{36} t^{30}-3 q^{37} t^{30}-5 q^{38} t^{30}
-6 q^{39} t^{30}+q^{40} t^{30}+10 q^{41} t^{30}+5 q^{42} t^{30}
-q^{36} t^{31}-q^{37} t^{31}-3 q^{38} t^{31}-4 q^{39} t^{31}
-2 q^{40} t^{31}+7 q^{41} t^{31}+4 q^{42} t^{31}-q^{37} t^{32}
-q^{38} t^{32}-3 q^{39} t^{32}-2 q^{40} t^{32}+3 q^{41} t^{32}
+4 q^{42} t^{32}-q^{38} t^{33}-q^{39} t^{33}-2 q^{40} t^{33}
+q^{41} t^{33}+3 q^{42} t^{33}-q^{39} t^{34}-q^{40} t^{34}
+2 q^{42} t^{34}
-q^{40} t^{35}+q^{42} t^{35}-q^{41} t^{36}+q^{42} t^{36}\bigr)+
\)

\smallskip
\noindent
\(
a^2 \bigl(q^3+q^4+2 q^5+2 q^6+3 q^7+3 q^8+3 q^9+2 q^{10}+2 q^{11}
+q^{12}+q^{13}-q^3 t+2 q^6 t+3 q^7 t+6 q^8 t+9 q^9 t+13 q^{10} t
+14 q^{11} t+15 q^{12} t+13 q^{13} t+12 q^{14} t+8 q^{15} t
+5 q^{16} t+2 q^{17} t+q^{18} t-q^4 t^2-q^5 t^2-2 q^6 t^2-q^7 t^2
-q^8 t^2+2 q^9 t^2+5 q^{10} t^2+13 q^{11} t^2+20 q^{12} t^2
+30 q^{13} t^2+34 q^{14} t^2+39 q^{15} t^2+35 q^{16} t^2
+31 q^{17} t^2+20 q^{18} t^2+13 q^{19} t^2+5 q^{20} t^2
+2 q^{21} t^2-q^5 t^3-q^6 t^3-3 q^7 t^3-3 q^8 t^3-5 q^9 t^3
-4 q^{10} t^3-4 q^{11} t^3+2 q^{12} t^3+8 q^{13} t^3+23 q^{14} t^3
+36 q^{15} t^3+55 q^{16} t^3+64 q^{17} t^3+71 q^{18} t^3
+62 q^{19} t^3+51 q^{20} t^3+30 q^{21} t^3+17 q^{22} t^3
+6 q^{23} t^3+2 q^{24} t^3-q^6 t^4-q^7 t^4-3 q^8 t^4-4 q^9 t^4
-7 q^{10} t^4-8 q^{11} t^4-11 q^{12} t^4-9 q^{13} t^4-8 q^{14} t^4
+3 q^{15} t^4+15 q^{16} t^4+40 q^{17} t^4+62 q^{18} t^4
+91 q^{19} t^4+101 q^{20} t^4+105 q^{21} t^4+80 q^{22} t^4
+57 q^{23} t^4+29 q^{24} t^4+13 q^{25} t^4+3 q^{26} t^4+q^{27} t^4
-q^7 t^5-q^8 t^5-3 q^9 t^5-4 q^{10} t^5-8 q^{11} t^5-10 q^{12} t^5
-15 q^{13} t^5-16 q^{14} t^5-20 q^{15} t^5-15 q^{16} t^5
-10 q^{17} t^5+11 q^{18} t^5+32 q^{19} t^5+71 q^{20} t^5
+103 q^{21} t^5+137 q^{22} t^5+132 q^{23} t^5+115 q^{24} t^5
+72 q^{25} t^5+39 q^{26} t^5+13 q^{27} t^5+4 q^{28} t^5-q^8 t^6
-q^9 t^6-3 q^{10} t^6-4 q^{11} t^6-8 q^{12} t^6-11 q^{13} t^6
-17 q^{14} t^6-20 q^{15} t^6-27 q^{16} t^6-27 q^{17} t^6
-29 q^{18} t^6-16 q^{19} t^6-2 q^{20} t^6+33 q^{21} t^6
+70 q^{22} t^6+126 q^{23} t^6+151 q^{24} t^6+165 q^{25} t^6
+123 q^{26} t^6+77 q^{27} t^6+31 q^{28} t^6+11 q^{29} t^6
+q^{30} t^6-q^9 t^7-q^{10} t^7-3 q^{11} t^7-4 q^{12} t^7
-8 q^{13} t^7-11 q^{14} t^7-18 q^{15} t^7-22 q^{16} t^7
-31 q^{17} t^7-34 q^{18} t^7-41 q^{19} t^7-35 q^{20} t^7
-30 q^{21} t^7-3 q^{22} t^7+28 q^{23} t^7+89 q^{24} t^7
+135 q^{25} t^7+185 q^{26} t^7+163 q^{27} t^7+116 q^{28} t^7
+51 q^{29} t^7+19 q^{30} t^7+2 q^{31} t^7-q^{10} t^8-q^{11} t^8
-3 q^{12} t^8-4 q^{13} t^8-8 q^{14} t^8-11 q^{15} t^8-18 q^{16}t^8
-23 q^{17} t^8-33 q^{18} t^8-38 q^{19} t^8-48 q^{20} t^8
-47 q^{21} t^8-49 q^{22} t^8-31 q^{23} t^8-9 q^{24} t^8
+46 q^{25} t^8+100 q^{26} t^8+180 q^{27} t^8+184 q^{28} t^8
+148 q^{29} t^8+69 q^{30} t^8+26 q^{31} t^8+3 q^{32} t^8
-q^{11} t^9-q^{12} t^9-3 q^{13} t^9-4 q^{14} t^9-8 q^{15} t^9
-11 q^{16} t^9-18 q^{17} t^9-23 q^{18} t^9-34 q^{19} t^9
-40 q^{20} t^9-52 q^{21} t^9-54 q^{22} t^9-61 q^{23} t^9
-50 q^{24} t^9-37 q^{25} t^9+9 q^{26} t^9+62 q^{27} t^9
+160 q^{28} t^9+190 q^{29} t^9+167 q^{30} t^9+82 q^{31} t^9
+30 q^{32} t^9+3 q^{33} t^9-q^{12} t^{10}-q^{13} t^{10}
-3 q^{14} t^{10}-4 q^{15} t^{10}-8 q^{16} t^{10}-11 q^{17} t^{10}
-18 q^{18} t^{10}-23 q^{19} t^{10}-34 q^{20} t^{10}
-41 q^{21} t^{10}-54 q^{22} t^{10}-58 q^{23} t^{10}-68 q^{24}t^{10}
-62 q^{25} t^{10}-56 q^{26} t^{10}-18 q^{27} t^{10}+31 q^{28}t^{10}
+139 q^{29} t^{10}+188 q^{30} t^{10}+175 q^{31} t^{10}
+86 q^{32} t^{10}+30 q^{33} t^{10}+3 q^{34} t^{10}-q^{13} t^{11}
-q^{14} t^{11}-3 q^{15} t^{11}-4 q^{16} t^{11}-8 q^{17} t^{11}
-11 q^{18} t^{11}-18 q^{19} t^{11}-23 q^{20} t^{11}-34 q^{21}t^{11}
-41 q^{22} t^{11}-55 q^{23} t^{11}-60 q^{24} t^{11}-72 q^{25}t^{11}
-69 q^{26} t^{11}-68 q^{27} t^{11}-35 q^{28} t^{11}+12 q^{29}t^{11}
+127 q^{30} t^{11}+187 q^{31} t^{11}+175 q^{32} t^{11}
+82 q^{33} t^{11}+26 q^{34} t^{11}+2 q^{35} t^{11}-q^{14} t^{12}
-q^{15} t^{12}-3 q^{16} t^{12}-4 q^{17} t^{12}-8 q^{18} t^{12}
-11 q^{19} t^{12}-18 q^{20} t^{12}-23 q^{21} t^{12}-34 q^{22}t^{12}
-41 q^{23} t^{12}-55 q^{24} t^{12}-61 q^{25} t^{12}-74 q^{26}t^{12}
-73 q^{27} t^{12}-74 q^{28} t^{12}-43 q^{29} t^{12}+5 q^{30} t^{12}
+127 q^{31} t^{12}+188 q^{32} t^{12}+167 q^{33} t^{12}
+69 q^{34} t^{12}+19 q^{35} t^{12}+q^{36} t^{12}-q^{15} t^{13}
-q^{16} t^{13}-3 q^{17} t^{13}-4 q^{18} t^{13}-8 q^{19} t^{13}
-11 q^{20} t^{13}-18 q^{21} t^{13}-23 q^{22} t^{13}-34 q^{23}t^{13}
-41 q^{24} t^{13}-55 q^{25} t^{13}-61 q^{26} t^{13}-75 q^{27}t^{13}
-75 q^{28} t^{13}-76 q^{29} t^{13}-43 q^{30} t^{13}+12 q^{31}t^{13}
+139 q^{32} t^{13}+190 q^{33} t^{13}+148 q^{34} t^{13}
+51 q^{35} t^{13}+11 q^{36} t^{13}-q^{16} t^{14}-q^{17} t^{14}
-3 q^{18} t^{14}-4 q^{19} t^{14}-8 q^{20} t^{14}-11 q^{21} t^{14}
-18 q^{22} t^{14}-23 q^{23} t^{14}-34 q^{24} t^{14}-41 q^{25}t^{14}
-55 q^{26} t^{14}-61 q^{27} t^{14}-75 q^{28} t^{14}-75 q^{29}t^{14}
-74 q^{30} t^{14}-35 q^{31} t^{14}+31 q^{32} t^{14}
+160 q^{33} t^{14}+184 q^{34} t^{14}+116 q^{35} t^{14}
+31 q^{36} t^{14}+4 q^{37} t^{14}-q^{17} t^{15}-q^{18} t^{15}
-3 q^{19} t^{15}-4 q^{20} t^{15}-8 q^{21} t^{15}-11 q^{22} t^{15}
-18 q^{23} t^{15}-23 q^{24} t^{15}-34 q^{25} t^{15}-41 q^{26}t^{15}
-55 q^{27} t^{15}-61 q^{28} t^{15}-75 q^{29} t^{15}-73 q^{30}t^{15}
-68 q^{31} t^{15}-18 q^{32} t^{15}+62 q^{33} t^{15}
+180 q^{34} t^{15}+163 q^{35} t^{15}+77 q^{36} t^{15}
+13 q^{37} t^{15}+q^{38} t^{15}-q^{18} t^{16}-q^{19} t^{16}
-3 q^{20} t^{16}-4 q^{21} t^{16}-8 q^{22} t^{16}-11 q^{23} t^{16}
-18 q^{24} t^{16}-23 q^{25} t^{16}-34 q^{26} t^{16}-41 q^{27}t^{16}
-55 q^{28} t^{16}-61 q^{29} t^{16}-74 q^{30} t^{16}-69 q^{31}t^{16}
-56 q^{32} t^{16}+9 q^{33} t^{16}+100 q^{34} t^{16}
+185 q^{35} t^{16}+123 q^{36} t^{16}+39 q^{37} t^{16}
+3 q^{38} t^{16}-q^{19} t^{17}-q^{20} t^{17}-3 q^{21} t^{17}
-4 q^{22} t^{17}-8 q^{23} t^{17}-11 q^{24} t^{17}-18 q^{25} t^{17}
-23 q^{26} t^{17}-34 q^{27} t^{17}-41 q^{28} t^{17}-55 q^{29}t^{17}
-61 q^{30} t^{17}-72 q^{31} t^{17}-62 q^{32} t^{17}-37 q^{33}t^{17}
+46 q^{34} t^{17}+135 q^{35} t^{17}+165 q^{36} t^{17}
+72 q^{37} t^{17}+13 q^{38} t^{17}-q^{20} t^{18}-q^{21} t^{18}
-3 q^{22} t^{18}-4 q^{23} t^{18}-8 q^{24} t^{18}-11 q^{25} t^{18}
-18 q^{26} t^{18}-23 q^{27} t^{18}-34 q^{28} t^{18}-41 q^{29}t^{18}
-55 q^{30} t^{18}-60 q^{31} t^{18}-68 q^{32} t^{18}-50 q^{33}t^{18}
-9 q^{34} t^{18}+89 q^{35} t^{18}+151 q^{36} t^{18}
+115 q^{37} t^{18}+29 q^{38} t^{18}+2 q^{39} t^{18}-q^{21} t^{19}
-q^{22} t^{19}-3 q^{23} t^{19}-4 q^{24} t^{19}-8 q^{25} t^{19}
-11 q^{26} t^{19}-18 q^{27} t^{19}-23 q^{28} t^{19}-34 q^{29}t^{19}
-41 q^{30} t^{19}-55 q^{31} t^{19}-58 q^{32} t^{19}-61 q^{33}t^{19}
-31 q^{34} t^{19}+28 q^{35} t^{19}+126 q^{36} t^{19}
+132 q^{37} t^{19}+57 q^{38} t^{19}+6 q^{39} t^{19}-q^{22} t^{20}
-q^{23} t^{20}-3 q^{24} t^{20}-4 q^{25} t^{20}-8 q^{26} t^{20}
-11 q^{27} t^{20}-18 q^{28} t^{20}-23 q^{29} t^{20}-34 q^{30}t^{20}
-41 q^{31} t^{20}-54 q^{32} t^{20}-54 q^{33} t^{20}-49 q^{34}t^{20}
-3 q^{35} t^{20}+70 q^{36} t^{20}+137 q^{37} t^{20}+80 q^{38}t^{20}
+17 q^{39} t^{20}-q^{23} t^{21}-q^{24} t^{21}-3 q^{25} t^{21}
-4 q^{26} t^{21}-8 q^{27} t^{21}-11 q^{28} t^{21}-18 q^{29} t^{21}
-23 q^{30} t^{21}-34 q^{31} t^{21}-41 q^{32} t^{21}-52 q^{33}t^{21}
-47 q^{34} t^{21}-30 q^{35} t^{21}+33 q^{36} t^{21}+103q^{37}t^{21}
+105 q^{38} t^{21}+30 q^{39} t^{21}+2 q^{40} t^{21}-q^{24} t^{22}
-q^{25} t^{22}-3 q^{26} t^{22}-4 q^{27} t^{22}-8 q^{28} t^{22}
-11 q^{29} t^{22}-18 q^{30} t^{22}-23 q^{31} t^{22}-34 q^{32}t^{22}
-40 q^{33} t^{22}-48 q^{34} t^{22}-35 q^{35} t^{22}-2 q^{36} t^{22}
+71 q^{37} t^{22}+101 q^{38} t^{22}+51 q^{39} t^{22}+5 q^{40}t^{22}
-q^{25} t^{23}-q^{26} t^{23}-3 q^{27} t^{23}-4 q^{28} t^{23}
-8 q^{29} t^{23}-11 q^{30} t^{23}-18 q^{31} t^{23}-23 q^{32} t^{23}
-34 q^{33} t^{23}-38 q^{34} t^{23}-41 q^{35} t^{23}-16 q^{36}t^{23}
+32 q^{37} t^{23}+91 q^{38} t^{23}+62 q^{39} t^{23}+13 q^{40}t^{23}
-q^{26} t^{24}-q^{27} t^{24}-3 q^{28} t^{24}-4 q^{29} t^{24}
-8 q^{30} t^{24}-11 q^{31} t^{24}-18 q^{32} t^{24}-23 q^{33} t^{24}
-33 q^{34} t^{24}-34 q^{35} t^{24}-29 q^{36} t^{24}+11 q^{37}t^{24}
+62 q^{38} t^{24}+71 q^{39} t^{24}+20 q^{40} t^{24}+q^{41} t^{24}
-q^{27} t^{25}-q^{28} t^{25}-3 q^{29} t^{25}-4 q^{30} t^{25}
-8 q^{31} t^{25}-11 q^{32} t^{25}-18 q^{33} t^{25}-23 q^{34} t^{25}
-31 q^{35} t^{25}-27 q^{36} t^{25}-10 q^{37} t^{25}+40 q^{38}t^{25}
+64 q^{39} t^{25}+31 q^{40} t^{25}+2 q^{41} t^{25}-q^{28} t^{26}
-q^{29} t^{26}-3 q^{30} t^{26}-4 q^{31} t^{26}-8 q^{32} t^{26}
-11 q^{33} t^{26}-18 q^{34} t^{26}-22 q^{35} t^{26}-27 q^{36}t^{26}
-15 q^{37} t^{26}+15 q^{38} t^{26}+55 q^{39} t^{26}+35 q^{40}t^{26}
+5 q^{41} t^{26}-q^{29} t^{27}-q^{30} t^{27}-3 q^{31} t^{27}
-4 q^{32} t^{27}-8 q^{33} t^{27}-11 q^{34} t^{27}-18 q^{35} t^{27}
-20 q^{36} t^{27}-20 q^{37} t^{27}+3 q^{38} t^{27}+36 q^{39} t^{27}
+39 q^{40} t^{27}+8 q^{41} t^{27}-q^{30} t^{28}-q^{31} t^{28}
-3 q^{32} t^{28}-4 q^{33} t^{28}-8 q^{34} t^{28}-11 q^{35} t^{28}
-17 q^{36} t^{28}-16 q^{37} t^{28}-8 q^{38} t^{28}+23 q^{39} t^{28}
+34 q^{40} t^{28}+12 q^{41} t^{28}-q^{31} t^{29}-q^{32} t^{29}
-3 q^{33} t^{29}-4 q^{34} t^{29}-8 q^{35} t^{29}-11 q^{36} t^{29}
-15 q^{37} t^{29}-9 q^{38} t^{29}+8 q^{39} t^{29}+30 q^{40} t^{29}
+13 q^{41} t^{29}+q^{42} t^{29}-q^{32} t^{30}-q^{33} t^{30}
-3 q^{34} t^{30}-4 q^{35} t^{30}-8 q^{36} t^{30}-10 q^{37} t^{30}
-11 q^{38} t^{30}+2 q^{39} t^{30}+20 q^{40} t^{30}+15 q^{41} t^{30}
+q^{42} t^{30}-q^{33} t^{31}-q^{34} t^{31}-3 q^{35} t^{31}
-4 q^{36} t^{31}-8 q^{37} t^{31}-8 q^{38} t^{31}-4 q^{39} t^{31}
+13 q^{40} t^{31}+14 q^{41} t^{31}+2 q^{42} t^{31}-q^{34} t^{32}
-q^{35} t^{32}-3 q^{36} t^{32}-4 q^{37} t^{32}-7 q^{38} t^{32}
-4 q^{39} t^{32}+5 q^{40} t^{32}+13 q^{41} t^{32}+2 q^{42} t^{32}
-q^{35} t^{33}-q^{36} t^{33}-3 q^{37} t^{33}-4 q^{38} t^{33}
-5 q^{39} t^{33}+2 q^{40} t^{33}+9 q^{41} t^{33}+3 q^{42} t^{33}
-q^{36} t^{34}-q^{37} t^{34}-3 q^{38} t^{34}-3 q^{39} t^{34}
-q^{40} t^{34}+6 q^{41} t^{34}+3 q^{42} t^{34}-q^{37} t^{35}
-q^{38} t^{35}-3 q^{39} t^{35}-q^{40} t^{35}+3 q^{41} t^{35}
+3 q^{42} t^{35}-q^{38} t^{36}-q^{39} t^{36}-2 q^{40} t^{36}
+2 q^{41} t^{36}+2 q^{42} t^{36}-q^{39} t^{37}-q^{40} t^{37}
+2 q^{42} t^{37}
-q^{40} t^{38}+q^{42} t^{38}-q^{41} t^{39}+q^{42} t^{39}\bigr)+
\)

\smallskip
\noindent
\(
a \bigl(q+q^2+q^3+q^4+q^5+q^6+q^7-q t+q^3 t+2 q^4 t+3 q^5 t+4 q^6 t
+5 q^7 t+7 q^8 t+6 q^9 t+5 q^{10} t+3 q^{11} t+2 q^{12} t+q^{13} t
-q^2 t^2-q^3 t^2-q^4 t^2+q^6 t^2+3 q^7 t^2+5 q^8 t^2+10 q^9 t^2
+13 q^{10} t^2+17 q^{11} t^2+16 q^{12} t^2+15 q^{13} t^2
+11 q^{14} t^2+7 q^{15} t^2+3 q^{16} t^2+q^{17} t^2-q^3 t^3-q^4 t^3
-2 q^5 t^3-2 q^6 t^3-2 q^7 t^3-q^8 t^3+5 q^{10} t^3+9 q^{11} t^3
+18 q^{12} t^3+24 q^{13} t^3+31 q^{14} t^3+31 q^{15} t^3
+29 q^{16} t^3+20 q^{17} t^3+12 q^{18} t^3+5 q^{19} t^3+q^{20} t^3
-q^4 t^4-q^5 t^4-2 q^6 t^4-3 q^7 t^4-4 q^8 t^4-4 q^9 t^4
-5 q^{10} t^4-2 q^{11} t^4+8 q^{13} t^4+16 q^{14} t^4+30 q^{15} t^4
+39 q^{16} t^4+50 q^{17} t^4+48 q^{18} t^4+41 q^{19} t^4
+25 q^{20} t^4+12 q^{21} t^4+4 q^{22} t^4+q^{23} t^4-q^5 t^5
-q^6 t^5-2 q^7 t^5-3 q^8 t^5-5 q^9 t^5-6 q^{10} t^5-8 q^{11} t^5
-7 q^{12} t^5-8 q^{13} t^5-3 q^{14} t^5+2 q^{15} t^5+16 q^{16} t^5
+28 q^{17} t^5+48 q^{18} t^5+60 q^{19} t^5+71 q^{20} t^5
+59 q^{21} t^5+40 q^{22} t^5+20 q^{23} t^5+7 q^{24} t^5+q^{25} t^5
-q^6 t^6-q^7 t^6-2 q^8 t^6-3 q^9 t^6-5 q^{10} t^6-7 q^{11} t^6
-10 q^{12} t^6-10 q^{13} t^6-13 q^{14} t^6-11 q^{15} t^6
-10 q^{16} t^6+10 q^{18} t^6+31 q^{19} t^6+48 q^{20} t^6
+76 q^{21} t^6+82 q^{22} t^6+75 q^{23} t^6+49 q^{24} t^6
+23 q^{25} t^6+6 q^{26} t^6+q^{27} t^6-q^7 t^7-q^8 t^7-2 q^9 t^7
-3 q^{10} t^7-5 q^{11} t^7-7 q^{12} t^7-11 q^{13} t^7-12 q^{14} t^7
-16 q^{15} t^7-16 q^{16} t^7-18 q^{17} t^7-12 q^{18} t^7
-7 q^{19} t^7+11 q^{20} t^7+27 q^{21} t^7+61 q^{22} t^7
+81 q^{23} t^7+97 q^{24} t^7+79 q^{25} t^7+46 q^{26} t^7
+16 q^{27} t^7+3 q^{28} t^7-q^8 t^8-q^9 t^8-2 q^{10} t^8
-3 q^{11} t^8-5 q^{12} t^8-7 q^{13} t^8-11 q^{14} t^8-13 q^{15} t^8
-18 q^{16} t^8-19 q^{17} t^8-23 q^{18} t^8-20 q^{19} t^8
-19 q^{20} t^8-6 q^{21} t^8+6 q^{22} t^8+39 q^{23} t^8
+64 q^{24} t^8+102 q^{25} t^8+100 q^{26} t^8+69 q^{27} t^8
+28 q^{28} t^8+6 q^{29} t^8-q^9 t^9-q^{10} t^9-2 q^{11} t^9
-3 q^{12} t^9-5 q^{13} t^9-7 q^{14} t^9-11 q^{15} t^9-13 q^{16} t^9
-19 q^{17} t^9-21 q^{18} t^9-26 q^{19} t^9-25 q^{20} t^9
-27 q^{21} t^9-18 q^{22} t^9-11 q^{23} t^9+18 q^{24} t^9
+41 q^{25} t^9+93 q^{26} t^9+109 q^{27} t^9+88 q^{28} t^9
+40 q^{29} t^9+9 q^{30} t^9-q^{10} t^{10}-q^{11} t^{10}
-2 q^{12} t^{10}-3 q^{13} t^{10}-5 q^{14} t^{10}-7 q^{15} t^{10}
-11 q^{16} t^{10}-13 q^{17} t^{10}-19 q^{18} t^{10}-22 q^{19}t^{10}
-28 q^{20} t^{10}-28 q^{21} t^{10}-32 q^{22} t^{10}-26 q^{23}t^{10}
-23 q^{24} t^{10}+q^{25} t^{10}+21 q^{26} t^{10}+79 q^{27} t^{10}
+108 q^{28} t^{10}+99 q^{29} t^{10}+47 q^{30} t^{10}
+11 q^{31} t^{10}-q^{11} t^{11}-q^{12} t^{11}-2 q^{13} t^{11}
-3 q^{14} t^{11}-5 q^{15} t^{11}-7 q^{16} t^{11}-11 q^{17} t^{11}
-13 q^{18} t^{11}-19 q^{19} t^{11}-22 q^{20} t^{11}-29 q^{21}t^{11}
-30 q^{22} t^{11}-35 q^{23} t^{11}-31 q^{24} t^{11}-31 q^{25}t^{11}
-11 q^{26} t^{11}+6 q^{27} t^{11}+67 q^{28} t^{11}+105 q^{29}t^{11}
+104 q^{30} t^{11}+50 q^{31} t^{11}+11 q^{32} t^{11}-q^{12} t^{12}
-q^{13} t^{12}-2 q^{14} t^{12}-3 q^{15} t^{12}-5 q^{16} t^{12}
-7 q^{17} t^{12}-11 q^{18} t^{12}-13 q^{19} t^{12}-19 q^{20} t^{12}
-22 q^{21} t^{12}-29 q^{22} t^{12}-31 q^{23} t^{12}-37 q^{24}t^{12}
-34 q^{25} t^{12}-36 q^{26} t^{12}-18 q^{27} t^{12}-3 q^{28} t^{12}
+60 q^{29} t^{12}+104 q^{30} t^{12}+104 q^{31} t^{12}
+47 q^{32} t^{12}+9 q^{33} t^{12}-q^{13} t^{13}-q^{14} t^{13}
-2 q^{15} t^{13}-3 q^{16} t^{13}-5 q^{17} t^{13}-7 q^{18} t^{13}
-11 q^{19} t^{13}-13 q^{20} t^{13}-19 q^{21} t^{13}-22 q^{22}t^{13}
-29 q^{23} t^{13}-31 q^{24} t^{13}-38 q^{25} t^{13}-36 q^{26}t^{13}
-39 q^{27} t^{13}-21 q^{28} t^{13}-6 q^{29} t^{13}+60 q^{30} t^{13}
+105 q^{31} t^{13}+99 q^{32} t^{13}+40 q^{33} t^{13}+6 q^{34}t^{13}
-q^{14} t^{14}-q^{15} t^{14}-2 q^{16} t^{14}-3 q^{17} t^{14}
-5 q^{18} t^{14}-7 q^{19} t^{14}-11 q^{20} t^{14}-13 q^{21} t^{14}
-19 q^{22} t^{14}-22 q^{23} t^{14}-29 q^{24} t^{14}-31 q^{25}t^{14}
-38 q^{26} t^{14}-37 q^{27} t^{14}-40 q^{28} t^{14}-21 q^{29}t^{14}
-3 q^{30} t^{14}+67 q^{31} t^{14}+108 q^{32} t^{14}+88 q^{33}t^{14}
+28 q^{34} t^{14}+3 q^{35} t^{14}-q^{15} t^{15}-q^{16} t^{15}
-2 q^{17} t^{15}-3 q^{18} t^{15}-5 q^{19} t^{15}-7 q^{20} t^{15}
-11 q^{21} t^{15}-13 q^{22} t^{15}-19 q^{23} t^{15}-22 q^{24}t^{15}
-29 q^{25} t^{15}-31 q^{26} t^{15}-38 q^{27} t^{15}-37 q^{28}t^{15}
-39 q^{29} t^{15}-18 q^{30} t^{15}+6 q^{31} t^{15}+79 q^{32} t^{15}
+109 q^{33} t^{15}+69 q^{34} t^{15}+16 q^{35} t^{15}+q^{36} t^{15}
-q^{16} t^{16}-q^{17} t^{16}-2 q^{18} t^{16}-3 q^{19} t^{16}
-5 q^{20} t^{16}-7 q^{21} t^{16}-11 q^{22} t^{16}-13 q^{23} t^{16}
-19 q^{24} t^{16}-22 q^{25} t^{16}-29 q^{26} t^{16}-31 q^{27}t^{16}
-38 q^{28} t^{16}-36 q^{29} t^{16}-36 q^{30} t^{16}-11 q^{31}t^{16}
+21 q^{32} t^{16}+93 q^{33} t^{16}+100 q^{34}t^{16}+46 q^{35}t^{16}
+6 q^{36} t^{16}-q^{17} t^{17}-q^{18} t^{17}-2 q^{19} t^{17}
-3 q^{20} t^{17}-5 q^{21} t^{17}-7 q^{22} t^{17}-11 q^{23} t^{17}
-13 q^{24} t^{17}-19 q^{25} t^{17}-22 q^{26} t^{17}-29 q^{27}t^{17}
-31 q^{28} t^{17}-38 q^{29} t^{17}-34 q^{30} t^{17}-31 q^{31}t^{17}
+q^{32} t^{17}+41 q^{33} t^{17}+102 q^{34} t^{17}+79 q^{35} t^{17}
+23 q^{36} t^{17}+q^{37} t^{17}-q^{18} t^{18}-q^{19} t^{18}
-2 q^{20} t^{18}-3 q^{21} t^{18}-5 q^{22} t^{18}-7 q^{23} t^{18}
-11 q^{24} t^{18}-13 q^{25} t^{18}-19 q^{26} t^{18}-22 q^{27}t^{18}
-29 q^{28} t^{18}-31 q^{29} t^{18}-37 q^{30} t^{18}-31 q^{31}t^{18}
-23 q^{32} t^{18}+18 q^{33} t^{18}+64 q^{34} t^{18}+97 q^{35}t^{18}
+49 q^{36} t^{18}+7 q^{37} t^{18}-q^{19} t^{19}-q^{20} t^{19}
-2 q^{21} t^{19}-3 q^{22} t^{19}-5 q^{23} t^{19}-7 q^{24} t^{19}
-11 q^{25} t^{19}-13 q^{26} t^{19}-19 q^{27} t^{19}-22 q^{28}t^{19}
-29 q^{29} t^{19}-31 q^{30} t^{19}-35 q^{31} t^{19}-26 q^{32}t^{19}
-11 q^{33} t^{19}+39 q^{34} t^{19}+81 q^{35} t^{19}+75 q^{36}t^{19}
+20 q^{37} t^{19}+q^{38} t^{19}-q^{20} t^{20}-q^{21} t^{20}
-2 q^{22} t^{20}-3 q^{23} t^{20}-5 q^{24} t^{20}-7 q^{25} t^{20}
-11 q^{26} t^{20}-13 q^{27} t^{20}-19 q^{28} t^{20}-22 q^{29}t^{20}
-29 q^{30} t^{20}-30 q^{31} t^{20}-32 q^{32} t^{20}-18 q^{33}t^{20}
+6 q^{34} t^{20}+61 q^{35} t^{20}+82 q^{36} t^{20}+40 q^{37} t^{20}
+4 q^{38} t^{20}-q^{21} t^{21}-q^{22} t^{21}-2 q^{23} t^{21}
-3 q^{24} t^{21}-5 q^{25} t^{21}-7 q^{26} t^{21}-11 q^{27} t^{21}
-13 q^{28} t^{21}-19 q^{29} t^{21}-22 q^{30} t^{21}-29 q^{31}t^{21}
-28 q^{32} t^{21}-27 q^{33} t^{21}-6 q^{34} t^{21}+27 q^{35} t^{21}
+76 q^{36} t^{21}+59 q^{37} t^{21}+12 q^{38} t^{21}-q^{22} t^{22}
-q^{23} t^{22}-2 q^{24} t^{22}-3 q^{25} t^{22}-5 q^{26} t^{22}
-7 q^{27} t^{22}-11 q^{28} t^{22}-13 q^{29} t^{22}-19 q^{30} t^{22}
-22 q^{31} t^{22}-28 q^{32} t^{22}-25 q^{33} t^{22}-19 q^{34}t^{22}
+11 q^{35} t^{22}+48 q^{36} t^{22}+71 q^{37} t^{22}+25 q^{38}t^{22}
+q^{39} t^{22}-q^{23} t^{23}-q^{24} t^{23}-2 q^{25} t^{23}
-3 q^{26} t^{23}-5 q^{27} t^{23}-7 q^{28} t^{23}-11 q^{29} t^{23}
-13 q^{30} t^{23}-19 q^{31} t^{23}-22 q^{32} t^{23}-26 q^{33}t^{23}
-20 q^{34} t^{23}-7 q^{35} t^{23}+31 q^{36} t^{23}+60 q^{37} t^{23}
+41 q^{38} t^{23}+5 q^{39} t^{23}-q^{24} t^{24}-q^{25} t^{24}
-2 q^{26} t^{24}-3 q^{27} t^{24}-5 q^{28} t^{24}-7 q^{29} t^{24}
-11 q^{30} t^{24}-13 q^{31} t^{24}-19 q^{32} t^{24}-21 q^{33}t^{24}
-23 q^{34} t^{24}-12 q^{35} t^{24}+10 q^{36} t^{24}+48 q^{37}t^{24}
+48 q^{38} t^{24}+12 q^{39} t^{24}-q^{25} t^{25}-q^{26} t^{25}
-2 q^{27} t^{25}-3 q^{28} t^{25}-5 q^{29} t^{25}-7 q^{30} t^{25}
-11 q^{31} t^{25}-13 q^{32} t^{25}-19 q^{33} t^{25}-19 q^{34}t^{25}
-18 q^{35} t^{25}+28 q^{37} t^{25}+50 q^{38} t^{25}+20 q^{39}t^{25}
+q^{40} t^{25}-q^{26} t^{26}-q^{27} t^{26}-2 q^{28} t^{26}
-3 q^{29} t^{26}-5 q^{30} t^{26}-7 q^{31} t^{26}-11 q^{32} t^{26}
-13 q^{33} t^{26}-18 q^{34} t^{26}-16 q^{35} t^{26}-10 q^{36}t^{26}
+16 q^{37} t^{26}+39 q^{38} t^{26}+29 q^{39} t^{26}+3 q^{40} t^{26}
-q^{27} t^{27}-q^{28} t^{27}-2 q^{29} t^{27}-3 q^{30} t^{27}
-5 q^{31} t^{27}-7 q^{32} t^{27}-11 q^{33} t^{27}-13 q^{34} t^{27}
-16 q^{35} t^{27}-11 q^{36} t^{27}+2 q^{37} t^{27}+30 q^{38} t^{27}
+31 q^{39} t^{27}+7 q^{40} t^{27}-q^{28} t^{28}-q^{29} t^{28}
-2 q^{30} t^{28}-3 q^{31} t^{28}-5 q^{32} t^{28}-7 q^{33} t^{28}
-11 q^{34} t^{28}-12 q^{35} t^{28}-13 q^{36} t^{28}-3 q^{37} t^{28}
+16 q^{38} t^{28}+31 q^{39} t^{28}+11 q^{40} t^{28}-q^{29} t^{29}
-q^{30} t^{29}-2 q^{31} t^{29}-3 q^{32} t^{29}-5 q^{33} t^{29}
-7 q^{34} t^{29}-11 q^{35} t^{29}-10 q^{36} t^{29}-8 q^{37} t^{29}
+8 q^{38} t^{29}+24 q^{39} t^{29}+15 q^{40} t^{29}+q^{41} t^{29}
-q^{30} t^{30}-q^{31} t^{30}-2 q^{32} t^{30}-3 q^{33} t^{30}
-5 q^{34} t^{30}-7 q^{35} t^{30}-10 q^{36} t^{30}-7 q^{37} t^{30}
+18 q^{39} t^{30}+16 q^{40} t^{30}+2 q^{41} t^{30}-q^{31} t^{31}
-q^{32} t^{31}-2 q^{33} t^{31}-3 q^{34} t^{31}-5 q^{35} t^{31}
-7 q^{36} t^{31}-8 q^{37} t^{31}-2 q^{38} t^{31}+9 q^{39} t^{31}
+17 q^{40} t^{31}+3 q^{41} t^{31}-q^{32} t^{32}-q^{33} t^{32}
-2 q^{34} t^{32}-3 q^{35} t^{32}-5 q^{36} t^{32}-6 q^{37} t^{32}
-5 q^{38} t^{32}+5 q^{39} t^{32}+13 q^{40} t^{32}+5 q^{41} t^{32}
-q^{33} t^{33}-q^{34} t^{33}-2 q^{35} t^{33}-3 q^{36} t^{33}
-5 q^{37} t^{33}-4 q^{38} t^{33}+10 q^{40} t^{33}+6 q^{41} t^{33}
-q^{34} t^{34}-q^{35} t^{34}-2 q^{36} t^{34}-3 q^{37} t^{34}
-4 q^{38} t^{34}-q^{39} t^{34}+5 q^{40} t^{34}+7 q^{41} t^{34}
-q^{35} t^{35}-q^{36} t^{35}-2 q^{37} t^{35}-3 q^{38} t^{35}
-2 q^{39} t^{35}+3 q^{40} t^{35}+5 q^{41} t^{35}+q^{42} t^{35}
-q^{36} t^{36}-q^{37} t^{36}-2 q^{38} t^{36}-2 q^{39} t^{36}
+q^{40} t^{36}+4 q^{41} t^{36}+q^{42} t^{36}-q^{37} t^{37}
-q^{38} t^{37}-2 q^{39} t^{37}+3 q^{41} t^{37}+q^{42} t^{37}
-q^{38} t^{38}-q^{39} t^{38}-q^{40} t^{38}+2 q^{41} t^{38}
+q^{42} t^{38}-q^{39} t^{39}-q^{40} t^{39}+q^{41} t^{39}
+q^{42} t^{39}-q^{40} t^{40}+q^{42} t^{40}-q^{41} t^{41}
+q^{42} t^{41}\bigr).
\)
}
\renewcommand{\baselinestretch}{1.0} 
\smallskip

The $a$\~degree of $\hat{\h}{}^{min}_{\l}\,(q,t,a)$ is
deg$_a=\ss_1\rr_2(|\yng(1)\,|+|\yng(1)\,|)-|\yng(1)\,|=7.$
This polynomial is
self-dual; all uncolored ones are self-dual.
The positivity of the series
$\hat{\h}{}^{min}_{\l}(q,t,a)/(1-t)$ holds.

Let us provide an equation of the corresponding 
plane curve singularity at $x=0,y=0$:
$$
(y^4+2x^3y^2+x^6+x^5y)(y^4+2x^3y^2+x^6-x^5y)=0.
$$

The linking number between its two component is $26$.
Its Alexander polynomial 
coincides with  $\hat{\h}{}^{min}_{\l}\,(q,q,-a)/(1-q)^2=$ 
{\small
\begin{align*}
1&+q^8+q^{12}+q^{16}+q^{20}+q^{24}+q^{26}+q^{28}+q^{32}+q^{34}
+q^{36}+q^{38}+q^{40}+q^{42}\\
&+q^{44}+q^{46}+q^{48}+q^{50}+q^{54}
+q^{56}+q^{58}+q^{62}+q^{66}+q^{70}+q^{74}+q^{82}. 
\end{align*}
}

\setcounter{equation}{0}
\section{\sc Intermediate twisting}

Let us extend the twisting construction
from Section \ref{sec:twist} in the following
natural direction. In 
$\bigl\{\hat{\xi}\bigl(
\hat{\ga}(\,'\hat{\P}_{0}^{tot})\!\Downarrow\bigr)\,
\hat{\P}_{0}^{tot}\bigr\}_{ev}$,
instead of the twist at the end (as we did),
one can do this {\em before\,} taking $\ga$. Namely, we
can consider $\bigl\{\hat{\ga}\bigl(
\hat{\xi}(\,'\hat{\P}_{0}^{tot})\,
\hat{\P}_{0}^{tot}\!\!\Downarrow\bigr)
\bigr\}_{ev}$
followed by the hat-normalization. Generally, this can be
done at any place of an arbitrary 
iterated link (after taking the corresponding 
$\Downarrow$).  This demonstrates what can be expected
if we extend the class of pre-polynomials by adding
those in the form $P(Y)(Q(X))\,R(X)$ 
for pre-polynomials $P,Q,R$ with further applying $\ga$
(and so on); this is related to the generalized coinvariants from
(\ref{knotoper}). 

Let us do
this for the example in (\ref{T1-0-3-2-2-1K}),
where  $\xi=\ga_{3,2}, \ga=\ga_{2,1}$. Instead of
$\lxi\h=
\{\hat{\xi}(J_{\yng(1)})\,(\hat{\ga}
(J^\circ_{\yng(1)\,})\!\Downarrow)\}_{ev}$ 
considered there, we will switch to more involved
${}_\xi\!\h=
\bigl\{\hat{\ga}\Bigl(\bigl(\hat{\xi}(J_{\yng(1)\,})
J^\circ_{\yng(1)\,}\bigr)\!\Downarrow\Bigr)\bigr\}_{ev}\,,$ 
and even more involved
${}_\xi^{\!\Box}\!\h=
\bigl\{\hat{\ga}\Bigl(\bigl(\hat{\xi}(J_{\yng(1)\,})
J_{\yng(1)\,}\bigr)\!\Downarrow\Bigr)
J^\circ_{\yng(1)}\bigr\}_{ev}\,.$ Recall that
$J_\la^\circ=P_\la^\circ=P_\la/P_\la(t^\rho).$ 

We obtain:
\Yboxdim5pt
\begin{align}\label{T3-2-2-1KK}
&\xi=\ga_{3,2},\ \ 
\l=\, \l^{\circ\!\rightarrow,\, \yng(1)}_{\{1,0\}},\ \  
'\!\l=\l^{\circ\!\rightarrow,\, \yng(1)}_{\{2,1\}},\ \  
{}_\xi\!\hat{\h}{}^{min}_{\l,\,'\!\l}\,(q,t,a)=
\end{align}

\renewcommand{\baselinestretch}{0.5} 
\noindent
{\small
\(
1-t+q t+q^2 t+q^3 t+q^4 t+q^5 t-q t^2+q^4 t^2+q^5 t^2+3 q^6 t^2
+2 q^7 t^2+q^8 t^2-q^2 t^3-q^4 t^3+2 q^7 t^3+3 q^8 t^3+4 q^9 t^3
+q^{10} t^3-q^3 t^4-q^5 t^4-q^6 t^4-q^7 t^4+q^8 t^4+q^9 t^4
+5 q^{10} t^4+4 q^{11} t^4+q^{12} t^4-q^4 t^5-q^6 t^5-q^7 t^5
-2 q^8 t^5+3 q^{11} t^5+5 q^{12} t^5+4 q^{13} t^5-q^5 t^6-q^7 t^6
-q^8 t^6-2 q^9 t^6-q^{10} t^6-q^{11} t^6+2 q^{12} t^6+4 q^{13} t^6
+5 q^{14} t^6+q^{15} t^6-q^6 t^7-q^8 t^7-q^9 t^7-2 q^{10} t^7
-q^{11} t^7-2 q^{12} t^7+q^{13} t^7+5 q^{14} t^7+4 q^{15} t^7
+q^{16} t^7-q^7 t^8-q^9 t^8-q^{10} t^8-2 q^{11} t^8-q^{12} t^8
-2 q^{13} t^8+5 q^{15} t^8+5 q^{16} t^8-q^8 t^9-q^{10} t^9
-q^{11} t^9-2 q^{12} t^9-q^{13} t^9-2 q^{14} t^9+q^{15} t^9
+4 q^{16} t^9+4 q^{17} t^9-q^9 t^{10}-q^{11} t^{10}-q^{12} t^{10}
-2 q^{13} t^{10}-q^{14} t^{10}-2 q^{15} t^{10}+2 q^{16} t^{10}
+5 q^{17} t^{10}+q^{18} t^{10}-q^{10} t^{11}-q^{12} t^{11}
-q^{13} t^{11}-2 q^{14} t^{11}-q^{15} t^{11}-q^{16} t^{11}
+3 q^{17} t^{11}+4 q^{18} t^{11}-q^{11} t^{12}-q^{13} t^{12}
-q^{14} t^{12}-2 q^{15} t^{12}-q^{16} t^{12}+5 q^{18} t^{12}
+q^{19} t^{12}-q^{12} t^{13}-q^{14} t^{13}-q^{15} t^{13}
-2 q^{16} t^{13}+q^{18} t^{13}+4 q^{19} t^{13}-q^{13} t^{14}
-q^{15} t^{14}-q^{16} t^{14}-2 q^{17} t^{14}+q^{18} t^{14}
+3 q^{19} t^{14}+q^{20} t^{14}-q^{14} t^{15}-q^{16} t^{15}
-q^{17} t^{15}-q^{18} t^{15}+2 q^{19} t^{15}+2 q^{20} t^{15}
-q^{15} t^{16}-q^{17} t^{16}-q^{18} t^{16}+3 q^{20} t^{16}
-q^{16} t^{17}-q^{18} t^{17}+q^{20} t^{17}+q^{21} t^{17}
-q^{17} t^{18}-q^{19} t^{18}+q^{20} t^{18}+q^{21} t^{18}
-q^{18} t^{19}+q^{21} t^{19}-q^{19} t^{20}+q^{21} t^{20}
-q^{20} t^{21}+q^{21} t^{21}-q^{21} t^{22}+q^{22} t^{22}
+a^5 \bigl(q^{15}-q^{15} t+q^{16} t+q^{17} t-q^{16} t^2
+q^{18} t^2+q^{19} t^2-q^{17} t^3+q^{20} t^3-q^{18} t^4
+q^{20} t^4-q^{19} t^5+q^{21} t^5-q^{20} t^6+q^{21} t^6
-q^{21} t^7+q^{22} t^7\bigr)+a^4 \bigl(q^{10}+q^{11}+q^{12}
+q^{13}+q^{14}-q^{10} t+q^{12} t+2 q^{13} t+2 q^{14} t+3 q^{15} t
+q^{16} t-q^{11} t^2-q^{12} t^2-q^{13} t^2+q^{14} t^2+2 q^{15} t^2
+4 q^{16} t^2+3 q^{17} t^2+q^{18} t^2-q^{12} t^3-q^{13} t^3
-2 q^{14} t^3-q^{15} t^3+q^{16} t^3+4 q^{17} t^3+3 q^{18} t^3
+2 q^{19} t^3-q^{13} t^4-q^{14} t^4-2 q^{15} t^4-2 q^{16} t^4
+4 q^{18} t^4+3 q^{19} t^4+q^{20} t^4-q^{14} t^5-q^{15} t^5
-2 q^{16} t^5-2 q^{17} t^5+4 q^{19} t^5+3 q^{20} t^5-q^{15} t^6
-q^{16} t^6-2 q^{17} t^6-2 q^{18} t^6+q^{19} t^6+4 q^{20} t^6
+q^{21} t^6-q^{16} t^7-q^{17} t^7-2 q^{18} t^7-q^{19} t^7
+2 q^{20} t^7+3 q^{21} t^7-q^{17} t^8-q^{18} t^8-2 q^{19} t^8
+q^{20} t^8+2 q^{21} t^8+q^{22} t^8-q^{18} t^9-q^{19} t^9
-q^{20} t^9+2 q^{21} t^9+q^{22} t^9-q^{19} t^{10}-q^{20} t^{10}
+q^{21} t^{10}+q^{22} t^{10}-q^{20} t^{11}+q^{22} t^{11}
-q^{21} t^{12}+q^{22} t^{12}\bigr)+a^3 \bigl(q^6+q^7+2 q^8
+2 q^9+2 q^{10}+q^{11}+q^{12}-q^6 t+2 q^9 t+4 q^{10} t+6 q^{11} t
+5 q^{12} t+4 q^{13} t+2 q^{14} t-q^7 t^2-q^8 t^2-2 q^9 t^2
-q^{10} t^2+q^{11} t^2+6 q^{12} t^2+8 q^{13} t^2+8 q^{14} t^2
+5 q^{15} t^2+2 q^{16} t^2-q^8 t^3-q^9 t^3-3 q^{10} t^3
-3 q^{11} t^3-3 q^{12} t^3+2 q^{13} t^3+7 q^{14} t^3
+11 q^{15} t^3+8 q^{16} t^3+4 q^{17} t^3+q^{18} t^3-q^9 t^4
-q^{10} t^4-3 q^{11} t^4-4 q^{12} t^4-5 q^{13} t^4-2 q^{14} t^4
+4 q^{15} t^4+11 q^{16} t^4+10 q^{17} t^4+4 q^{18} t^4+q^{19} t^4
-q^{10} t^5-q^{11} t^5-3 q^{12} t^5-4 q^{13} t^5-6 q^{14} t^5
-4 q^{15} t^5+2 q^{16} t^5+11 q^{17} t^5+10 q^{18} t^5
+4 q^{19} t^5-q^{11} t^6-q^{12} t^6-3 q^{13} t^6-4 q^{14} t^6
-6 q^{15} t^6-5 q^{16} t^6+2 q^{17} t^6+11 q^{18} t^6
+8 q^{19} t^6+2 q^{20} t^6-q^{12} t^7-q^{13} t^7-3 q^{14} t^7
-4 q^{15} t^7-6 q^{16} t^7-4 q^{17} t^7+4 q^{18} t^7
+11 q^{19} t^7+5 q^{20} t^7-q^{13} t^8-q^{14} t^8-3 q^{15} t^8
-4 q^{16} t^8-6 q^{17} t^8-2 q^{18} t^8+7 q^{19} t^8+8 q^{20} t^8
+2 q^{21} t^8-q^{14} t^9-q^{15} t^9-3 q^{16} t^9-4 q^{17} t^9
-5 q^{18} t^9+2 q^{19} t^9+8 q^{20} t^9+4 q^{21} t^9
-q^{15} t^{10}-q^{16} t^{10}-3 q^{17} t^{10}-4 q^{18} t^{10}
-3 q^{19} t^{10}+6 q^{20} t^{10}+5 q^{21} t^{10}+q^{22} t^{10}
-q^{16} t^{11}-q^{17} t^{11}-3 q^{18} t^{11}-3 q^{19} t^{11}
+q^{20} t^{11}+6 q^{21} t^{11}+q^{22} t^{11}-q^{17} t^{12}
-q^{18} t^{12}-3 q^{19} t^{12}-q^{20} t^{12}+4 q^{21} t^{12}
+2 q^{22} t^{12}-q^{18} t^{13}-q^{19} t^{13}-2 q^{20} t^{13}
+2 q^{21} t^{13}+2 q^{22} t^{13}-q^{19} t^{14}-q^{20} t^{14}
+2 q^{22} t^{14}-q^{20} t^{15}+q^{22} t^{15}-q^{21} t^{16}
+q^{22} t^{16}\bigr)+a^2 \bigl(q^3+q^4+2 q^5+2 q^6+2 q^7+q^8
+q^9-q^3 t+2 q^6 t+4 q^7 t+7 q^8 t+6 q^9 t+6 q^{10} t+3 q^{11} t
+q^{12} t-q^4 t^2-q^5 t^2-2 q^6 t^2-q^7 t^2+6 q^9 t^2
+9 q^{10} t^2+12 q^{11} t^2+9 q^{12} t^2+5 q^{13} t^2+q^{14} t^2
-q^5 t^3-q^6 t^3-3 q^7 t^3-3 q^8 t^3-4 q^9 t^3+5 q^{11} t^3
+14 q^{12} t^3+14 q^{13} t^3+11 q^{14} t^3+4 q^{15} t^3
+q^{16} t^3-q^6 t^4-q^7 t^4-3 q^8 t^4-4 q^9 t^4-6 q^{10} t^4
-4 q^{11} t^4-2 q^{12} t^4+10 q^{13} t^4+17 q^{14} t^4
+15 q^{15} t^4+7 q^{16} t^4+2 q^{17} t^4-q^7 t^5-q^8 t^5
-3 q^9 t^5-4 q^{10} t^5-7 q^{11} t^5-6 q^{12} t^5-6 q^{13} t^5
+4 q^{14} t^5+16 q^{15} t^5+18 q^{16} t^5+8 q^{17} t^5
+2 q^{18} t^5-q^8 t^6-q^9 t^6-3 q^{10} t^6-4 q^{11} t^6
-7 q^{12} t^6-7 q^{13} t^6-8 q^{14} t^6+2 q^{15} t^6
+14 q^{16} t^6+18 q^{17} t^6+7 q^{18} t^6+q^{19} t^6-q^9 t^7
-q^{10} t^7-3 q^{11} t^7-4 q^{12} t^7-7 q^{13} t^7-7 q^{14} t^7
-9 q^{15} t^7+2 q^{16} t^7+16 q^{17} t^7+15 q^{18} t^7
+4 q^{19} t^7-q^{10} t^8-q^{11} t^8-3 q^{12} t^8-4 q^{13} t^8
-7 q^{14} t^8-7 q^{15} t^8-8 q^{16} t^8+4 q^{17} t^8
+17 q^{18} t^8+11 q^{19} t^8+q^{20} t^8-q^{11} t^9-q^{12} t^9
-3 q^{13} t^9-4 q^{14} t^9-7 q^{15} t^9-7 q^{16} t^9-6 q^{17} t^9
+10 q^{18} t^9+14 q^{19} t^9+5 q^{20} t^9-q^{12} t^{10}
-q^{13} t^{10}-3 q^{14} t^{10}-4 q^{15} t^{10}-7 q^{16} t^{10}
-6 q^{17} t^{10}-2 q^{18} t^{10}+14 q^{19} t^{10}+9 q^{20} t^{10}
+q^{21} t^{10}-q^{13} t^{11}-q^{14} t^{11}-3 q^{15} t^{11}
-4 q^{16} t^{11}-7 q^{17} t^{11}-4 q^{18} t^{11}+5 q^{19} t^{11}
+12 q^{20} t^{11}+3 q^{21} t^{11}-q^{14} t^{12}-q^{15} t^{12}
-3 q^{16} t^{12}-4 q^{17} t^{12}-6 q^{18} t^{12}+9 q^{20} t^{12}
+6 q^{21} t^{12}-q^{15} t^{13}-q^{16} t^{13}-3 q^{17} t^{13}
-4 q^{18} t^{13}-4 q^{19} t^{13}+6 q^{20} t^{13}+6 q^{21} t^{13}
+q^{22} t^{13}-q^{16} t^{14}-q^{17} t^{14}-3 q^{18} t^{14}
-3 q^{19} t^{14}+7 q^{21} t^{14}+q^{22} t^{14}-q^{17} t^{15}
-q^{18} t^{15}-3 q^{19} t^{15}-q^{20} t^{15}+4 q^{21} t^{15}
+2 q^{22} t^{15}-q^{18} t^{16}-q^{19} t^{16}-2 q^{20} t^{16}
+2 q^{21} t^{16}+2 q^{22} t^{16}-q^{19} t^{17}-q^{20} t^{17}
+2 q^{22} t^{17}-q^{20} t^{18}+q^{22} t^{18}-q^{21} t^{19}
+q^{22} t^{19}\bigr)+a \bigl(q+q^2+q^3+q^4+q^5-q t+q^3 t
+2 q^4 t+3 q^5 t+5 q^6 t+4 q^7 t+2 q^8 t+q^9 t-q^2 t^2-q^3 t^2
-q^4 t^2+q^6 t^2+5 q^7 t^2+8 q^8 t^2+8 q^9 t^2+5 q^{10} t^2
+2 q^{11} t^2-q^3 t^3-q^4 t^3-2 q^5 t^3-2 q^6 t^3-2 q^7 t^3
+q^8 t^3+5 q^9 t^3+11 q^{10} t^3+11 q^{11} t^3+6 q^{12} t^3
+2 q^{13} t^3-q^4 t^4-q^5 t^4-2 q^6 t^4-3 q^7 t^4-4 q^8 t^4
-2 q^9 t^4+7 q^{11} t^4+14 q^{12} t^4+12 q^{13} t^4+5 q^{14} t^4
+q^{15} t^4-q^5 t^5-q^6 t^5-2 q^7 t^5-3 q^8 t^5-5 q^9 t^5
-4 q^{10} t^5-3 q^{11} t^5+2 q^{12} t^5+11 q^{13} t^5
+16 q^{14} t^5+8 q^{15} t^5+2 q^{16} t^5-q^6 t^6-q^7 t^6
-2 q^8 t^6-3 q^9 t^6-5 q^{10} t^6-5 q^{11} t^6-5 q^{12} t^6
-q^{13} t^6+9 q^{14} t^6+15 q^{15} t^6+10 q^{16} t^6+2 q^{17} t^6
-q^7 t^7-q^8 t^7-2 q^9 t^7-3 q^{10} t^7-5 q^{11} t^7-5 q^{12} t^7
-6 q^{13} t^7-3 q^{14} t^7+9 q^{15} t^7+15 q^{16} t^7
+8 q^{17} t^7+q^{18} t^7-q^8 t^8-q^9 t^8-2 q^{10} t^8
-3 q^{11} t^8-5 q^{12} t^8-5 q^{13} t^8-6 q^{14} t^8-3 q^{15} t^8
+9 q^{16} t^8+16 q^{17} t^8+5 q^{18} t^8-q^9 t^9-q^{10} t^9
-2 q^{11} t^9-3 q^{12} t^9-5 q^{13} t^9-5 q^{14} t^9-6 q^{15} t^9
-q^{16} t^9+11 q^{17} t^9+12 q^{18} t^9+2 q^{19} t^9
-q^{10} t^{10}-q^{11} t^{10}-2 q^{12} t^{10}-3 q^{13} t^{10}
-5 q^{14} t^{10}-5 q^{15} t^{10}-5 q^{16} t^{10}+2 q^{17} t^{10}
+14 q^{18} t^{10}+6 q^{19} t^{10}-q^{11} t^{11}-q^{12} t^{11}
-2 q^{13} t^{11}-3 q^{14} t^{11}-5 q^{15} t^{11}-5 q^{16} t^{11}
-3 q^{17} t^{11}+7 q^{18} t^{11}+11 q^{19} t^{11}+2 q^{20} t^{11}
-q^{12} t^{12}-q^{13} t^{12}-2 q^{14} t^{12}-3 q^{15} t^{12}
-5 q^{16} t^{12}-4 q^{17} t^{12}+11 q^{19} t^{12}+5 q^{20} t^{12}
-q^{13} t^{13}-q^{14} t^{13}-2 q^{15} t^{13}-3 q^{16} t^{13}
-5 q^{17} t^{13}-2 q^{18} t^{13}+5 q^{19} t^{13}+8 q^{20} t^{13}
+q^{21} t^{13}-q^{14} t^{14}-q^{15} t^{14}-2 q^{16} t^{14}
-3 q^{17} t^{14}-4 q^{18} t^{14}+q^{19} t^{14}+8 q^{20} t^{14}
+2 q^{21} t^{14}-q^{15} t^{15}-q^{16} t^{15}-2 q^{17} t^{15}
-3 q^{18} t^{15}-2 q^{19} t^{15}+5 q^{20} t^{15}+4 q^{21} t^{15}
-q^{16} t^{16}-q^{17} t^{16}-2 q^{18} t^{16}-2 q^{19} t^{16}
+q^{20} t^{16}+5 q^{21} t^{16}-q^{17} t^{17}-q^{18} t^{17}
-2 q^{19} t^{17}+3 q^{21} t^{17}+q^{22} t^{17}-q^{18} t^{18}
-q^{19} t^{18}-q^{20} t^{18}+2 q^{21} t^{18}+q^{22} t^{18}
-q^{19} t^{19}-q^{20} t^{19}+q^{21} t^{19}+q^{22} t^{19}
-q^{20} t^{20}+q^{22} t^{20}-q^{21} t^{21}+q^{22} t^{21}\bigr).
\)
}
\renewcommand{\baselinestretch}{1.0} 
\smallskip

\Yboxdim7pt
This polynomial is self-dual.  
The positivity of
${}_\xi\!\hat{\h}{}^{min}_{\l,\,'\!\l}(q,t,a)/(1-t)$ holds.
One has: ${}_\xi\!\hat{\h}{}^{min}_{\l,\,'\!\l}(1,t,a)=$
\renewcommand{\baselinestretch}{0.5} 
\noindent
{\small
\(
(1+a) (1+a+t) \bigl(1+3 a+3 a^2+a^3+3 t+7 a t+5 a^2 t+a^3 t
+4 t^2+8 a t^2+5 a^2 t^2+a^3 t^2+4 t^3+8 a t^3+4 a^2 t^3+4 t^4
+6 a t^4+2 a^2 t^4+3 t^5+4 a t^5+a^2 t^5+2 t^6+2 a t^6+t^7
+a t^7+t^8\bigr),
\)
}
\renewcommand{\baselinestretch}{1.0} 
where the last factor is that for $C\!ab(13,2)T(3,2)$.
This cable appears if we make $\l$ empty in the twisting
construction.

\smallskip
\Yboxdim5pt
Switching to
${}_\xi^{\!\Box}\!\h=
\bigl\{\hat{\ga}\Bigl(\bigl(\hat{\xi}(J_{\yng(1)\,})
J_{\yng(1)\,}\bigr)\!\Downarrow\Bigr)
J^\circ_{\yng(1)}\bigr\}_{ev}\,,$ 

\begin{align}\label{T3-2-2-1KK2}
&\xi=\ga_{3,2},\ \ 
\l=\, \l^{\circ\!\rightarrow,\, \yng(1)}_{\{1,0\}},\ \  
'\!\l=\l^{\circ\!\rightarrow,\, \yng(1)}_{\{2,1\}},\ \  
{}_\xi^{\!\Box}\!\hat{\h}{}^{min}_{\l,\,'\!\l}\,(q,t,a)=
\end{align}

\renewcommand{\baselinestretch}{0.5} 
\noindent
{\small
\(
1-2 t+q t+q^2 t+q^3 t+q^4 t+q^5 t+q^6 t+t^2-2 q t^2-q^2 t^2
-q^3 t^2+2 q^6 t^2+3 q^7 t^2+3 q^8 t^2+q^9 t^2+q^{10} t^2+q t^3
-q^2 t^3-2 q^4 t^3-q^5 t^3-3 q^6 t^3-q^7 t^3+5 q^9 t^3
+4 q^{10} t^3+4 q^{11} t^3+2 q^{12} t^3+q^2 t^4-q^3 t^4+q^4 t^4
-q^5 t^4-q^6 t^4-3 q^7 t^4-2 q^8 t^4-5 q^9 t^4+3 q^{11} t^4
+7 q^{12} t^4+6 q^{13} t^4+3 q^{14} t^4+q^3 t^5-q^4 t^5+q^5 t^5
-3 q^8 t^5-q^9 t^5-5 q^{10} t^5-4 q^{11} t^5-4 q^{12} t^5
+5 q^{13} t^5+8 q^{14} t^5+8 q^{15} t^5+2 q^{16} t^5+q^4 t^6
-q^5 t^6+q^6 t^6+q^8 t^6-2 q^9 t^6-q^{10} t^6-4 q^{11} t^6
-3 q^{12} t^6-8 q^{13} t^6-2 q^{14} t^6+6 q^{15} t^6
+11 q^{16} t^6+5 q^{17} t^6+q^{18} t^6+q^5 t^7-q^6 t^7+q^7 t^7
+q^9 t^7-q^{10} t^7-4 q^{12} t^7-2 q^{13} t^7-7 q^{14} t^7
-7 q^{15} t^7+q^{16} t^7+12 q^{17} t^7+7 q^{18} t^7+2 q^{19} t^7
+q^6 t^8-q^7 t^8+q^8 t^8+q^{10} t^8-q^{11} t^8+q^{12} t^8
-3 q^{13} t^8-2 q^{14} t^8-6 q^{15} t^8-8 q^{16} t^8-3 q^{17} t^8
+11 q^{18} t^8+9 q^{19} t^8+2 q^{20} t^8+q^7 t^9-q^8 t^9+q^9 t^9
+q^{11} t^9-q^{12} t^9+q^{13} t^9-2 q^{14} t^9-q^{15} t^9
-6 q^{16} t^9-9 q^{17} t^9-4 q^{18} t^9+10 q^{19} t^9
+9 q^{20} t^9+2 q^{21} t^9+q^8 t^{10}-q^9 t^{10}+q^{10} t^{10}
+q^{12} t^{10}-q^{13} t^{10}+q^{14} t^{10}-2 q^{15} t^{10}
-5 q^{17} t^{10}-10 q^{18} t^{10}-4 q^{19} t^{10}
+11 q^{20} t^{10}+7 q^{21} t^{10}+q^{22} t^{10}+q^9 t^{11}
-q^{10} t^{11}+q^{11} t^{11}+q^{13} t^{11}-q^{14} t^{11}
+q^{15} t^{11}-2 q^{16} t^{11}-5 q^{18} t^{11}-9 q^{19} t^{11}
-3 q^{20} t^{11}+12 q^{21} t^{11}+5 q^{22} t^{11}+q^{10} t^{12}
-q^{11} t^{12}+q^{12} t^{12}+q^{14} t^{12}-q^{15} t^{12}
+q^{16} t^{12}-2 q^{17} t^{12}-6 q^{19} t^{12}-8 q^{20} t^{12}
+q^{21} t^{12}+11 q^{22} t^{12}+2 q^{23} t^{12}+q^{11} t^{13}
-q^{12} t^{13}+q^{13} t^{13}+q^{15} t^{13}-q^{16} t^{13}
+q^{17} t^{13}-2 q^{18} t^{13}-q^{19} t^{13}-6 q^{20} t^{13}
-7 q^{21} t^{13}+6 q^{22} t^{13}+8 q^{23} t^{13}+q^{12} t^{14}
-q^{13} t^{14}+q^{14} t^{14}+q^{16} t^{14}-q^{17} t^{14}
+q^{18} t^{14}-2 q^{19} t^{14}-2 q^{20} t^{14}-7 q^{21} t^{14}
-2 q^{22} t^{14}+8 q^{23} t^{14}+3 q^{24} t^{14}+q^{13} t^{15}
-q^{14} t^{15}+q^{15} t^{15}+q^{17} t^{15}-q^{18} t^{15}
+q^{19} t^{15}-3 q^{20} t^{15}-2 q^{21} t^{15}-8 q^{22} t^{15}
+5 q^{23} t^{15}+6 q^{24} t^{15}+q^{14} t^{16}-q^{15} t^{16}
+q^{16} t^{16}+q^{18} t^{16}-q^{19} t^{16}+q^{20} t^{16}
-4 q^{21} t^{16}-3 q^{22} t^{16}-4 q^{23} t^{16}+7 q^{24} t^{16}
+2 q^{25} t^{16}+q^{15} t^{17}-q^{16} t^{17}+q^{17} t^{17}
+q^{19} t^{17}-q^{20} t^{17}-4 q^{22} t^{17}-4 q^{23} t^{17}
+3 q^{24} t^{17}+4 q^{25} t^{17}+q^{16} t^{18}-q^{17} t^{18}
+q^{18} t^{18}+q^{20} t^{18}-q^{21} t^{18}-q^{22} t^{18}
-5 q^{23} t^{18}+4 q^{25} t^{18}+q^{26} t^{18}+q^{17} t^{19}
-q^{18} t^{19}+q^{19} t^{19}+q^{21} t^{19}-2 q^{22} t^{19}
-q^{23} t^{19}-5 q^{24} t^{19}+5 q^{25} t^{19}+q^{26} t^{19}
+q^{18} t^{20}-q^{19} t^{20}+q^{20} t^{20}+q^{22} t^{20}
-3 q^{23} t^{20}-2 q^{24} t^{20}+3 q^{26} t^{20}+q^{19} t^{21}
-q^{20} t^{21}+q^{21} t^{21}-3 q^{24} t^{21}-q^{25} t^{21}
+3 q^{26} t^{21}+q^{20} t^{22}-q^{21} t^{22}+q^{22} t^{22}
-q^{24} t^{22}-3 q^{25} t^{22}+2 q^{26} t^{22}+q^{27} t^{22}
+q^{21} t^{23}-q^{22} t^{23}+q^{23} t^{23}-q^{24} t^{23}
-q^{25} t^{23}+q^{27} t^{23}+q^{22} t^{24}-q^{23} t^{24}
+q^{24} t^{24}-2 q^{25} t^{24}+q^{27} t^{24}+q^{23} t^{25}
-q^{24} t^{25}-q^{26} t^{25}+q^{27} t^{25}+q^{24} t^{26}
-q^{25} t^{26}-q^{26} t^{26}+q^{27} t^{26}+q^{25} t^{27}
-2 q^{26} t^{27}+q^{27} t^{27}+q^{26} t^{28}-2 q^{27} t^{28}
+q^{28} t^{28}
\)

\smallskip
\noindent
\(
+a^6 \bigl(q^{21}-q^{21} t+q^{22} t+q^{23} t
-q^{22} t^2+q^{24} t^2+q^{25} t^2-q^{23} t^3+q^{26} t^3
-q^{24} t^4+q^{26} t^4-q^{25} t^5+q^{27} t^5-q^{26} t^6
+q^{27} t^6-q^{27} t^7+q^{28} t^7\bigr)+a^5 \bigl(q^{15}
+q^{16}+q^{17}+q^{18}+q^{19}+q^{20}-2 q^{15} t+q^{17} t
+2 q^{18} t+2 q^{19} t+2 q^{20} t+3 q^{21} t+q^{22} t+q^{15} t^2
-2 q^{16} t^2-2 q^{17} t^2-2 q^{18} t^2+q^{19} t^2+2 q^{20} t^2
+3 q^{21} t^2+4 q^{22} t^2+3 q^{23} t^2+q^{24} t^2+q^{16} t^3
-q^{17} t^3-q^{18} t^3-4 q^{19} t^3-2 q^{20} t^3+3 q^{22} t^3
+4 q^{23} t^3+3 q^{24} t^3+2 q^{25} t^3+q^{17} t^4-q^{18} t^4
-3 q^{20} t^4-4 q^{21} t^4-2 q^{22} t^4+3 q^{23} t^4+4 q^{24} t^4
+3 q^{25} t^4+q^{26} t^4+q^{18} t^5-q^{19} t^5-2 q^{21} t^5
-4 q^{22} t^5-3 q^{23} t^5+3 q^{24} t^5+4 q^{25} t^5+3 q^{26} t^5
+q^{19} t^6-q^{20} t^6-2 q^{22} t^6-4 q^{23} t^6-2 q^{24} t^6
+3 q^{25} t^6+4 q^{26} t^6+q^{27} t^6+q^{20} t^7-q^{21} t^7
-2 q^{23} t^7-4 q^{24} t^7+3 q^{26} t^7+3 q^{27} t^7+q^{21} t^8
-q^{22} t^8-3 q^{24} t^8-2 q^{25} t^8+2 q^{26} t^8+2 q^{27} t^8
+q^{28} t^8+q^{22} t^9-q^{23} t^9-4 q^{25} t^9+q^{26} t^9
+2 q^{27} t^9+q^{28} t^9+q^{23} t^{10}-q^{24} t^{10}-q^{25} t^{10}
-2 q^{26} t^{10}+2 q^{27} t^{10}+q^{28} t^{10}+q^{24} t^{11}
-q^{25} t^{11}-2 q^{26} t^{11}+q^{27} t^{11}+q^{28} t^{11}
+q^{25} t^{12}-2 q^{26} t^{12}+q^{28} t^{12}+q^{26} t^{13}
-2 q^{27} t^{13}+q^{28} t^{13}\bigr)
\)

\smallskip
\noindent
\(
+a^4 \bigl(q^{10}+q^{11}
+2 q^{12}+2 q^{13}+3 q^{14}+2 q^{15}+2 q^{16}+q^{17}+q^{18}
-2 q^{10} t-q^{11} t-q^{12} t+q^{13} t+2 q^{14} t+6 q^{15} t
+7 q^{16} t+7 q^{17} t+5 q^{18} t+4 q^{19} t+2 q^{20} t
+q^{10} t^2-q^{11} t^2-2 q^{12} t^2-4 q^{13} t^2-4 q^{14} t^2
-4 q^{15} t^2+2 q^{16} t^2+7 q^{17} t^2+12 q^{18} t^2
+10 q^{19} t^2+9 q^{20} t^2+5 q^{21} t^2+2 q^{22} t^2+q^{11} t^3
-3 q^{14} t^3-5 q^{15} t^3-10 q^{16} t^3-6 q^{17} t^3-q^{18} t^3
+11 q^{19} t^3+13 q^{20} t^3+14 q^{21} t^3+8 q^{22} t^3
+4 q^{23} t^3+q^{24} t^3+q^{12} t^4+q^{14} t^4-q^{15} t^4
-3 q^{16} t^4-10 q^{17} t^4-11 q^{18} t^4-9 q^{19} t^4
+4 q^{20} t^4+13 q^{21} t^4+16 q^{22} t^4+10 q^{23} t^4
+4 q^{24} t^4+q^{25} t^4+q^{13} t^5+q^{15} t^5-q^{17} t^5
-8 q^{18} t^5-12 q^{19} t^5-14 q^{20} t^5-q^{21} t^5
+12 q^{22} t^5+17 q^{23} t^5+10 q^{24} t^5+4 q^{25} t^5
+q^{14} t^6+q^{16} t^6-6 q^{19} t^6-11 q^{20} t^6-17 q^{21} t^6
-3 q^{22} t^6+12 q^{23} t^6+16 q^{24} t^6+8 q^{25} t^6
+2 q^{26} t^6+q^{15} t^7+q^{17} t^7-5 q^{20} t^7-10 q^{21} t^7
-17 q^{22} t^7-q^{23} t^7+13 q^{24} t^7+14 q^{25} t^7
+5 q^{26} t^7+q^{16} t^8+q^{18} t^8-5 q^{21} t^8-11 q^{22} t^8
-14 q^{23} t^8+4 q^{24} t^8+13 q^{25} t^8+9 q^{26} t^8
+2 q^{27} t^8+q^{17} t^9+q^{19} t^9-6 q^{22} t^9-12 q^{23} t^9
-9 q^{24} t^9+11 q^{25} t^9+10 q^{26} t^9+4 q^{27} t^9
+q^{18} t^{10}+q^{20} t^{10}-8 q^{23} t^{10}-11 q^{24} t^{10}
-q^{25} t^{10}+12 q^{26} t^{10}+5 q^{27} t^{10}+q^{28} t^{10}
+q^{19} t^{11}+q^{21} t^{11}-q^{23} t^{11}-10 q^{24} t^{11}
-6 q^{25} t^{11}+7 q^{26} t^{11}+7 q^{27} t^{11}+q^{28} t^{11}
+q^{20} t^{12}+q^{22} t^{12}-3 q^{24} t^{12}-10 q^{25} t^{12}
+2 q^{26} t^{12}+7 q^{27} t^{12}+2 q^{28} t^{12}+q^{21} t^{13}
+q^{23} t^{13}-q^{24} t^{13}-5 q^{25} t^{13}-4 q^{26} t^{13}
+6 q^{27} t^{13}+2 q^{28} t^{13}+q^{22} t^{14}+q^{24} t^{14}
-3 q^{25} t^{14}-4 q^{26} t^{14}+2 q^{27} t^{14}+3 q^{28} t^{14}
+q^{23} t^{15}-4 q^{26} t^{15}+q^{27} t^{15}+2 q^{28} t^{15}
+q^{24} t^{16}-2 q^{26} t^{16}-q^{27} t^{16}+2 q^{28} t^{16}
+q^{25} t^{17}-q^{26} t^{17}-q^{27} t^{17}+q^{28} t^{17}
+q^{26} t^{18}-2 q^{27} t^{18}+q^{28} t^{18}\bigr)
\)

\smallskip
\noindent
\(
+a^3 \bigl(q^6+q^7+2 q^8+3 q^9+3 q^{10}+3 q^{11}+3 q^{12}
+2 q^{13}+q^{14}+q^{15}-2 q^6 t-q^7 t-2 q^8 t-q^9 t+2 q^{10} t
+6 q^{11} t+8 q^{12} t+11 q^{13} t+11 q^{14} t+8 q^{15} t
+6 q^{16} t+3 q^{17} t+q^{18} t+q^6 t^2-q^7 t^2-q^8 t^2-4 q^9 t^2
-6 q^{10} t^2-8 q^{11} t^2-3 q^{12} t^2+2 q^{13} t^2
+12 q^{14} t^2+18 q^{15} t^2+20 q^{16} t^2+16 q^{17} t^2
+10 q^{18} t^2+5 q^{19} t^2+q^{20} t^2+q^7 t^3+q^9 t^3
-2 q^{10} t^3-4 q^{11} t^3-11 q^{12} t^3-13 q^{13} t^3
-14 q^{14} t^3-q^{15} t^3+12 q^{16} t^3+25 q^{17} t^3
+28 q^{18} t^3+20 q^{19} t^3+12 q^{20} t^3+4 q^{21} t^3
+q^{22} t^3+q^8 t^4+2 q^{10} t^4-q^{12} t^4-7 q^{13} t^4
-12 q^{14} t^4-23 q^{15} t^4-18 q^{16} t^4-4 q^{17} t^4
+18 q^{18} t^4+32 q^{19} t^4+30 q^{20} t^4+17 q^{21} t^4
+7 q^{22} t^4+2 q^{23} t^4+q^9 t^5+2 q^{11} t^5+q^{12} t^5
+q^{13} t^5-4 q^{14} t^5-7 q^{15} t^5-21 q^{16} t^5-27 q^{17} t^5
-21 q^{18} t^5+6 q^{19} t^5+30 q^{20} t^5+36 q^{21} t^5
+21 q^{22} t^5+8 q^{23} t^5+2 q^{24} t^5+q^{10} t^6+2 q^{12} t^6
+q^{13} t^6+2 q^{14} t^6-2 q^{15} t^6-4 q^{16} t^6-16 q^{17} t^6
-27 q^{18} t^6-31 q^{19} t^6-6 q^{20} t^6+28 q^{21} t^6
+37 q^{22} t^6+21 q^{23} t^6+7 q^{24} t^6+q^{25} t^6+q^{11} t^7
+2 q^{13} t^7+q^{14} t^7+2 q^{15} t^7-q^{16} t^7-2 q^{17} t^7
-13 q^{18} t^7-25 q^{19} t^7-33 q^{20} t^7-11 q^{21} t^7
+28 q^{22} t^7+36 q^{23} t^7+17 q^{24} t^7+4 q^{25} t^7
+q^{12} t^8+2 q^{14} t^8+q^{15} t^8+2 q^{16} t^8-q^{17} t^8
-q^{18} t^8-11 q^{19} t^8-25 q^{20} t^8-33 q^{21} t^8
-6 q^{22} t^8+30 q^{23} t^8+30 q^{24} t^8+12 q^{25} t^8
+q^{26} t^8+q^{13} t^9+2 q^{15} t^9+q^{16} t^9+2 q^{17} t^9
-q^{18} t^9-q^{19} t^9-11 q^{20} t^9-25 q^{21} t^9-31 q^{22} t^9
+6 q^{23} t^9+32 q^{24} t^9+20 q^{25} t^9+5 q^{26} t^9
+q^{14} t^{10}+2 q^{16} t^{10}+q^{17} t^{10}+2 q^{18} t^{10}
-q^{19} t^{10}-q^{20} t^{10}-13 q^{21} t^{10}-27 q^{22} t^{10}
-21 q^{23} t^{10}+18 q^{24} t^{10}+28 q^{25} t^{10}
+10 q^{26} t^{10}+q^{27} t^{10}+q^{15} t^{11}+2 q^{17} t^{11}
+q^{18} t^{11}+2 q^{19} t^{11}-q^{20} t^{11}-2 q^{21} t^{11}
-16 q^{22} t^{11}-27 q^{23} t^{11}-4 q^{24} t^{11}
+25 q^{25} t^{11}+16 q^{26} t^{11}+3 q^{27} t^{11}
+q^{16} t^{12}+2 q^{18} t^{12}+q^{19} t^{12}+2 q^{20} t^{12}
-q^{21} t^{12}-4 q^{22} t^{12}-21 q^{23} t^{12}-18 q^{24} t^{12}
+12 q^{25} t^{12}+20 q^{26} t^{12}+6 q^{27} t^{12}+q^{17} t^{13}
+2 q^{19} t^{13}+q^{20} t^{13}+2 q^{21} t^{13}-2 q^{22} t^{13}
-7 q^{23} t^{13}-23 q^{24} t^{13}-q^{25} t^{13}+18 q^{26} t^{13}
+8 q^{27} t^{13}+q^{28} t^{13}+q^{18} t^{14}+2 q^{20} t^{14}
+q^{21} t^{14}+2 q^{22} t^{14}-4 q^{23} t^{14}-12 q^{24} t^{14}
-14 q^{25} t^{14}+12 q^{26} t^{14}+11 q^{27} t^{14}+q^{28} t^{14}
+q^{19} t^{15}+2 q^{21} t^{15}+q^{22} t^{15}+q^{23} t^{15}
-7 q^{24} t^{15}-13 q^{25} t^{15}+2 q^{26} t^{15}+11 q^{27} t^{15}
+2 q^{28} t^{15}+q^{20} t^{16}+2 q^{22} t^{16}+q^{23} t^{16}
-q^{24} t^{16}-11 q^{25} t^{16}-3 q^{26} t^{16}+8 q^{27} t^{16}
+3 q^{28} t^{16}+q^{21} t^{17}+2 q^{23} t^{17}-4 q^{25} t^{17}
-8 q^{26} t^{17}+6 q^{27} t^{17}+3 q^{28} t^{17}+q^{22} t^{18}
+2 q^{24} t^{18}-2 q^{25} t^{18}-6 q^{26} t^{18}+2 q^{27} t^{18}
+3 q^{28} t^{18}+q^{23} t^{19}+q^{25} t^{19}-4 q^{26} t^{19}
-q^{27} t^{19}+3 q^{28} t^{19}+q^{24} t^{20}-q^{26} t^{20}
-2 q^{27} t^{20}+2 q^{28} t^{20}+q^{25} t^{21}-q^{26} t^{21}
-q^{27} t^{21}+q^{28} t^{21}+q^{26} t^{22}-2 q^{27} t^{22}
+q^{28} t^{22}\bigr)
\)

\smallskip
\noindent
\(
+a^2 \bigl(q^3+q^4+2 q^5+2 q^6+3 q^7+2 q^8
+2 q^9+q^{10}+q^{11}-2 q^3 t-q^4 t-2 q^5 t+q^7 t+6 q^8 t+8 q^9 t
+10 q^{10} t+9 q^{11} t+8 q^{12} t+5 q^{13} t+2 q^{14} t+q^{15} t
+q^3 t^2-q^4 t^2-q^5 t^2-4 q^6 t^2-5 q^7 t^2-8 q^8 t^2-4 q^9 t^2
+q^{10} t^2+11 q^{11} t^2+16 q^{12} t^2+21 q^{13} t^2
+17 q^{14} t^2+12 q^{15} t^2+5 q^{16} t^2+2 q^{17} t^2+q^4 t^3
+q^6 t^3-2 q^7 t^3-3 q^8 t^3-10 q^9 t^3-12 q^{10} t^3
-16 q^{11} t^3-5 q^{12} t^3+5 q^{13} t^3+24 q^{14} t^3
+29 q^{15} t^3+28 q^{16} t^3+16 q^{17} t^3+7 q^{18} t^3
+2 q^{19} t^3+q^5 t^4+2 q^7 t^4-6 q^{10} t^4-9 q^{11} t^4
-21 q^{12} t^4-21 q^{13} t^4-17 q^{14} t^4+7 q^{15} t^4
+28 q^{16} t^4+40 q^{17} t^4+30 q^{18} t^4+16 q^{19} t^4
+5 q^{20} t^4+q^{21} t^4+q^6 t^5+2 q^8 t^5+q^9 t^5+2 q^{10} t^5
-3 q^{11} t^5-4 q^{12} t^5-16 q^{13} t^5-23 q^{14} t^5
-33 q^{15} t^5-16 q^{16} t^5+12 q^{17} t^5+42 q^{18} t^5
+40 q^{19} t^5+25 q^{20} t^5+8 q^{21} t^5+2 q^{22} t^5+q^7 t^6
+2 q^9 t^6+q^{10} t^6+3 q^{11} t^6-q^{12} t^6-q^{13} t^6
-11 q^{14} t^6-17 q^{15} t^6-36 q^{16} t^6-33 q^{17} t^6
-7 q^{18} t^6+36 q^{19} t^6+46 q^{20} t^6+29 q^{21} t^6
+10 q^{22} t^6+2 q^{23} t^6+q^8 t^7+2 q^{10} t^7+q^{11} t^7
+3 q^{12} t^7+q^{14} t^7-8 q^{15} t^7-12 q^{16} t^7-32 q^{17} t^7
-40 q^{18} t^7-20 q^{19} t^7+31 q^{20} t^7+47 q^{21} t^7
+29 q^{22} t^7+8 q^{23} t^7+q^{24} t^7+q^9 t^8+2 q^{11} t^8
+q^{12} t^8+3 q^{13} t^8+2 q^{15} t^8-6 q^{16} t^8-9 q^{17} t^8
-29 q^{18} t^8-42 q^{19} t^8-24 q^{20} t^8+31 q^{21} t^8
+46 q^{22} t^8+25 q^{23} t^8+5 q^{24} t^8+q^{10} t^9+2 q^{12} t^9
+q^{13} t^9+3 q^{14} t^9+2 q^{16} t^9-5 q^{17} t^9-7 q^{18} t^9
-28 q^{19} t^9-42 q^{20} t^9-20 q^{21} t^9+36 q^{22} t^9
+40 q^{23} t^9+16 q^{24} t^9+2 q^{25} t^9+q^{11} t^{10}
+2 q^{13} t^{10}+q^{14} t^{10}+3 q^{15} t^{10}+2 q^{17} t^{10}
-5 q^{18} t^{10}-7 q^{19} t^{10}-29 q^{20} t^{10}
-40 q^{21} t^{10}-7 q^{22} t^{10}+42 q^{23} t^{10}
+30 q^{24} t^{10}+7 q^{25} t^{10}+q^{12} t^{11}+2 q^{14} t^{11}
+q^{15} t^{11}+3 q^{16} t^{11}+2 q^{18} t^{11}-5 q^{19} t^{11}
-9 q^{20} t^{11}-32 q^{21} t^{11}-33 q^{22} t^{11}
+12 q^{23} t^{11}+40 q^{24} t^{11}+16 q^{25} t^{11}
+2 q^{26} t^{11}+q^{13} t^{12}+2 q^{15} t^{12}+q^{16} t^{12}
+3 q^{17} t^{12}+2 q^{19} t^{12}-6 q^{20} t^{12}-12 q^{21} t^{12}
-36 q^{22} t^{12}-16 q^{23} t^{12}+28 q^{24} t^{12}
+28 q^{25} t^{12}+5 q^{26} t^{12}+q^{14} t^{13}+2 q^{16} t^{13}
+q^{17} t^{13}+3 q^{18} t^{13}+2 q^{20} t^{13}-8 q^{21} t^{13}
-17 q^{22} t^{13}-33 q^{23} t^{13}+7 q^{24} t^{13}
+29 q^{25} t^{13}+12 q^{26} t^{13}+q^{27} t^{13}+q^{15} t^{14}
+2 q^{17} t^{14}+q^{18} t^{14}+3 q^{19} t^{14}+q^{21} t^{14}
-11 q^{22} t^{14}-23 q^{23} t^{14}-17 q^{24} t^{14}
+24 q^{25} t^{14}+17 q^{26} t^{14}+2 q^{27} t^{14}+q^{16} t^{15}
+2 q^{18} t^{15}+q^{19} t^{15}+3 q^{20} t^{15}-q^{22} t^{15}
-16 q^{23} t^{15}-21 q^{24} t^{15}+5 q^{25} t^{15}
+21 q^{26} t^{15}+5 q^{27} t^{15}+q^{17} t^{16}+2 q^{19} t^{16}
+q^{20} t^{16}+3 q^{21} t^{16}-q^{22} t^{16}-4 q^{23} t^{16}
-21 q^{24} t^{16}-5 q^{25} t^{16}+16 q^{26} t^{16}+8 q^{27} t^{16}
+q^{18} t^{17}+2 q^{20} t^{17}+q^{21} t^{17}+3 q^{22} t^{17}
-3 q^{23} t^{17}-9 q^{24} t^{17}-16 q^{25} t^{17}+11 q^{26} t^{17}
+9 q^{27} t^{17}+q^{28} t^{17}+q^{19} t^{18}+2 q^{21} t^{18}
+q^{22} t^{18}+2 q^{23} t^{18}-6 q^{24} t^{18}-12 q^{25} t^{18}
+q^{26} t^{18}+10 q^{27} t^{18}+q^{28} t^{18}+q^{20} t^{19}
+2 q^{22} t^{19}+q^{23} t^{19}-10 q^{25} t^{19}-4 q^{26} t^{19}
+8 q^{27} t^{19}+2 q^{28} t^{19}+q^{21} t^{20}+2 q^{23} t^{20}
-3 q^{25} t^{20}-8 q^{26} t^{20}+6 q^{27} t^{20}+2 q^{28} t^{20}
+q^{22} t^{21}+2 q^{24} t^{21}-2 q^{25} t^{21}-5 q^{26} t^{21}
+q^{27} t^{21}+3 q^{28} t^{21}+q^{23} t^{22}+q^{25} t^{22}
-4 q^{26} t^{22}+2 q^{28} t^{22}+q^{24} t^{23}-q^{26} t^{23}
-2 q^{27} t^{23}+2 q^{28} t^{23}+q^{25} t^{24}-q^{26} t^{24}
-q^{27} t^{24}+q^{28} t^{24}+q^{26} t^{25}-2 q^{27} t^{25}
+q^{28} t^{25}\bigr)
\)

\smallskip
\noindent
\(
+a \bigl(q+q^2+q^3+q^4+q^5+q^6-2 q t-q^2 t
+q^4 t+2 q^5 t+4 q^6 t+6 q^7 t+5 q^8 t+3 q^9 t+2 q^{10} t
+q^{11} t+q t^2-q^2 t^2-2 q^3 t^2-3 q^4 t^2-3 q^5 t^2-4 q^6 t^2
-q^7 t^2+4 q^8 t^2+10 q^9 t^2+11 q^{10} t^2+10 q^{11} t^2
+7 q^{12} t^2+3 q^{13} t^2+q^{14} t^2+q^2 t^3-2 q^5 t^3-3 q^6 t^3
-7 q^7 t^3-8 q^8 t^3-8 q^9 t^3+9 q^{11} t^3+17 q^{12} t^3
+18 q^{13} t^3+12 q^{14} t^3+6 q^{15} t^3+q^{16} t^3+q^3 t^4
+q^5 t^4-q^7 t^4-5 q^8 t^4-7 q^9 t^4-13 q^{10} t^4-13 q^{11} t^4
-7 q^{12} t^4+9 q^{13} t^4+22 q^{14} t^4+24 q^{15} t^4
+16 q^{16} t^4+6 q^{17} t^4+q^{18} t^4+q^4 t^5+q^6 t^5+q^7 t^5
+q^8 t^5-3 q^9 t^5-4 q^{10} t^5-10 q^{11} t^5-15 q^{12} t^5
-20 q^{13} t^5-9 q^{14} t^5+12 q^{15} t^5+29 q^{16} t^5
+26 q^{17} t^5+13 q^{18} t^5+4 q^{19} t^5+q^5 t^6+q^7 t^6+q^8 t^6
+2 q^9 t^6-q^{10} t^6-2 q^{11} t^6-7 q^{12} t^6-11 q^{13} t^6
-21 q^{14} t^6-23 q^{15} t^6-5 q^{16} t^6+24 q^{17} t^6
+32 q^{18} t^6+20 q^{19} t^6+6 q^{20} t^6+q^{21} t^6+q^6 t^7
+q^8 t^7+q^9 t^7+2 q^{10} t^7-5 q^{13} t^7-8 q^{14} t^7
-17 q^{15} t^7-27 q^{16} t^7-18 q^{17} t^7+15 q^{18} t^7
+34 q^{19} t^7+24 q^{20} t^7+6 q^{21} t^7+q^{22} t^7+q^7 t^8
+q^9 t^8+q^{10} t^8+2 q^{11} t^8+q^{13} t^8-3 q^{14} t^8
-6 q^{15} t^8-14 q^{16} t^8-27 q^{17} t^8-24 q^{18} t^8
+9 q^{19} t^8+35 q^{20} t^8+24 q^{21} t^8+6 q^{22} t^8+q^8 t^9
+q^{10} t^9+q^{11} t^9+2 q^{12} t^9+q^{14} t^9-2 q^{15} t^9
-4 q^{16} t^9-12 q^{17} t^9-27 q^{18} t^9-26 q^{19} t^9
+9 q^{20} t^9+34 q^{21} t^9+20 q^{22} t^9+4 q^{23} t^9+q^9 t^{10}
+q^{11} t^{10}+q^{12} t^{10}+2 q^{13} t^{10}+q^{15} t^{10}
-2 q^{16} t^{10}-3 q^{17} t^{10}-11 q^{18} t^{10}-27 q^{19} t^{10}
-24 q^{20} t^{10}+15 q^{21} t^{10}+32 q^{22} t^{10}
+13 q^{23} t^{10}+q^{24} t^{10}+q^{10} t^{11}+q^{12} t^{11}
+q^{13} t^{11}+2 q^{14} t^{11}+q^{16} t^{11}-2 q^{17} t^{11}
-3 q^{18} t^{11}-12 q^{19} t^{11}-27 q^{20} t^{11}
-18 q^{21} t^{11}+24 q^{22} t^{11}+26 q^{23} t^{11}
+6 q^{24} t^{11}+q^{11} t^{12}+q^{13} t^{12}+q^{14} t^{12}
+2 q^{15} t^{12}+q^{17} t^{12}-2 q^{18} t^{12}-4 q^{19} t^{12}
-14 q^{20} t^{12}-27 q^{21} t^{12}-5 q^{22} t^{12}
+29 q^{23} t^{12}+16 q^{24} t^{12}+q^{25} t^{12}+q^{12} t^{13}
+q^{14} t^{13}+q^{15} t^{13}+2 q^{16} t^{13}+q^{18} t^{13}
-2 q^{19} t^{13}-6 q^{20} t^{13}-17 q^{21} t^{13}-23 q^{22} t^{13}
+12 q^{23} t^{13}+24 q^{24} t^{13}+6 q^{25} t^{13}+q^{13} t^{14}
+q^{15} t^{14}+q^{16} t^{14}+2 q^{17} t^{14}+q^{19} t^{14}
-3 q^{20} t^{14}-8 q^{21} t^{14}-21 q^{22} t^{14}-9 q^{23} t^{14}
+22 q^{24} t^{14}+12 q^{25} t^{14}+q^{26} t^{14}+q^{14} t^{15}
+q^{16} t^{15}+q^{17} t^{15}+2 q^{18} t^{15}+q^{20} t^{15}
-5 q^{21} t^{15}-11 q^{22} t^{15}-20 q^{23} t^{15}+9 q^{24} t^{15}
+18 q^{25} t^{15}+3 q^{26} t^{15}+q^{15} t^{16}+q^{17} t^{16}
+q^{18} t^{16}+2 q^{19} t^{16}-7 q^{22} t^{16}-15 q^{23} t^{16}
-7 q^{24} t^{16}+17 q^{25} t^{16}+7 q^{26} t^{16}+q^{16} t^{17}
+q^{18} t^{17}+q^{19} t^{17}+2 q^{20} t^{17}-2 q^{22} t^{17}
-10 q^{23} t^{17}-13 q^{24} t^{17}+9 q^{25} t^{17}
+10 q^{26} t^{17}+q^{27} t^{17}+q^{17} t^{18}+q^{19} t^{18}
+q^{20} t^{18}+2 q^{21} t^{18}-q^{22} t^{18}-4 q^{23} t^{18}
-13 q^{24} t^{18}+11 q^{26} t^{18}+2 q^{27} t^{18}+q^{18} t^{19}
+q^{20} t^{19}+q^{21} t^{19}+2 q^{22} t^{19}-3 q^{23} t^{19}
-7 q^{24} t^{19}-8 q^{25} t^{19}+10 q^{26} t^{19}+3 q^{27} t^{19}
+q^{19} t^{20}+q^{21} t^{20}+q^{22} t^{20}+q^{23} t^{20}
-5 q^{24} t^{20}-8 q^{25} t^{20}+4 q^{26} t^{20}+5 q^{27} t^{20}
+q^{20} t^{21}+q^{22} t^{21}+q^{23} t^{21}-q^{24} t^{21}
-7 q^{25} t^{21}-q^{26} t^{21}+6 q^{27} t^{21}+q^{21} t^{22}
+q^{23} t^{22}-3 q^{25} t^{22}-4 q^{26} t^{22}+4 q^{27} t^{22}
+q^{28} t^{22}+q^{22} t^{23}+q^{24} t^{23}-2 q^{25} t^{23}
-3 q^{26} t^{23}+2 q^{27} t^{23}+q^{28} t^{23}+q^{23} t^{24}
-3 q^{26} t^{24}+q^{27} t^{24}+q^{28} t^{24}+q^{24} t^{25}
-2 q^{26} t^{25}+q^{28} t^{25}+q^{25} t^{26}-q^{26} t^{26}
-q^{27} t^{26}
+q^{28} t^{26}+q^{26} t^{27}-2 q^{27} t^{27}+q^{28} t^{27}\bigr).
\)
}
\renewcommand{\baselinestretch}{1.0} 
\smallskip

\Yboxdim7pt
This polynomial is self-dual.  
The positivity of
${}_\xi^{\!\Box}\!\hat{\h}{}^{min}_{\l,\,'\!\l}(q,t,a)/(1-t)^2$ 
holds.
One has: ${}_\xi^{\!\Box}\!\hat{\h}{}^{min}_{\l,\,'\!\l}(1,t,a)=$
\renewcommand{\baselinestretch}{0.5} 
\noindent
{\small
\(
(1+a)^2 (1+a+t) \bigl(1+3 a+3 a^2+a^3+3 t+7 a t+5 a^2 t+a^3 t
+4 t^2+8 a t^2+5 a^2 t^2+a^3 t^2+4 t^3+8 a t^3+4 a^2 t^3+4 t^4
+6 a t^4+2 a^2 t^4+3 t^5+4 a t^5+a^2 t^5+2 t^6+2 a t^6+t^7
+a t^7+t^8\bigr).
\)
}
\renewcommand{\baselinestretch}{1.0} 
\medskip

\Yboxdim5pt
{\sf Two examples with $\ga_{1,1}$.}
Let us simplify here $\xi=\ga_{3,2}$ 
as much as possible. We will take now $\xi=\ga_{1,1}$, say 
$\xi=\tau_-$:

\Yboxdim5pt
\begin{align}\label{T1-0-2-1KK}
&\xi=\ga_{1,1},\ \ 
\l=\, \l^{\circ\!\rightarrow,\, \yng(1)}_{\{1,0\}},\ \ 
'\!\l=\l^{\circ\!\rightarrow,\, \yng(1)}_{\{2,1\}},\ \  
{}_\xi\!\hat{\h}^{min}_{\l,\,'\!\l}\,(q,t,a)=
\end{align}

\renewcommand{\baselinestretch}{0.5} 
\noindent
{\small
\(
1-t+q t+q^2 t-q t^2+q^3 t^2-q^2 t^3+q^3 t^3-q^3 t^4+q^4 t^4
+a^2 \bigl(q^3-q^3 t+q^4 t\bigr)
+a \bigl(q+q^2-q t+2 q^3 t-q^2 t^2+q^4 t^2-q^3 t^3+q^4 t^3\bigr).
\)
}
\renewcommand{\baselinestretch}{1.0} 
\smallskip

\Yboxdim5pt
This polynomial coincides with $\hat{\h}^{min}$ for
$T(3,2)^{\,\prime,\vee}_{\yng(1)\,}$ from
(\ref{T6-4-yv}). 
Indeed, the latter is the 
hat-normalization of
$$
\{\tau_+\tau_-^2(J_{\yng(1)}),\iota(J^\circ_{\yng(1)})\}_{ev}=
\{\vph(\tau_+\tau_-^2(\vph\si^{-1}J_{\yng(1)})\,
J^\circ_{\yng(1)})\}_{ev}=
\{\tau_-\tau_+\tau_-(J_{\yng(1)})\,J^\circ_{\yng(1)}\}_{ev},
$$
which coincides with 
$
\{\tau_-\bigl(\tau_+\tau_-(J_{\yng(1)})\,
J^\circ_{\yng(1)}\bigr)\}_{ev}\!=\!
\{\tau_-\tau_+\tau_-(J_{\yng(1)})\,J^\circ_{\yng(1)}\}_{ev}
$
for (\ref{T1-0-2-1KK}).
\smallskip

Generally, we do not see the reasons for such a connection with
our core link-invariants,
but this can be expected if some of the  $PSL_2(\Z)$\~matrices 
involved are powers of $\tau_-$.
Let us add another $\yng(1)$\, to such $\h$, 
reducing (\ref{T3-2-2-1KK2}) to $\ga_{1,1}$:

\Yboxdim5pt
\begin{align}\label{T1-0-2-1KK2}
&\xi=\ga_{1,1},\ \ 
\l=\, \l^{\circ\!\rightarrow,\, \yng(1)}_{\{1,0\}},\ \ 
'\!\l=\l^{\circ\!\rightarrow,\, \yng(1)}_{\{2,1\}},\ \  
{}_\xi^{\!\Box}\!\hat{\h}^{min}_{\l,\,'\!\l}\,(q,t,a)=
\end{align}

\renewcommand{\baselinestretch}{0.5} 
\noindent
{\small
\(
1-2 t+q t+q^2 t+q^3 t+t^2-2 q t^2-q^2 t^2+2 q^4 t^2
+q t^3-q^2 t^3-q^3 t^3-q^4 t^3+2 q^5 t^3+q^2 t^4-q^3 t^4
-q^4 t^4+q^6 t^4+q^3 t^5-q^4 t^5-q^5 t^5+q^6 t^5+q^4 t^6
-2 q^5 t^6+q^6 t^6+q^5 t^7-2 q^6 t^7+q^7 t^7
+a^3 \bigl(q^6-q^6 t+q^7 t\bigr)+a^2 \bigl(q^3+q^4+q^5-2 q^3 t
+q^5 t+2 q^6 t+q^3 t^2-2 q^4 t^2-q^5 t^2+q^6 t^2+q^7 t^2+q^4 t^3
-2 q^5 t^3+q^7 t^3+q^5 t^4-2 q^6 t^4+q^7 t^4\bigr)+a \bigl(q+q^2
+q^3-2 q t-q^2 t+q^3 t+3 q^4 t+q^5 t+q t^2-q^2 t^2-3 q^3 t^2
-q^4 t^2+3 q^5 t^2+q^6 t^2+q^2 t^3-3 q^4 t^3-q^5 t^3+3 q^6 t^3
+q^3 t^4-3 q^5 t^4+q^6 t^4+q^7 t^4+q^4 t^5
-q^5 t^5-q^6 t^5+q^7 t^5+q^5 t^6-2 q^6 t^6+q^7 t^6\bigr).
\)
}
\renewcommand{\baselinestretch}{1.0} 
\smallskip

\Yboxdim7pt
It is self-dual. The positivity of
${}_\xi^{\!\Box}\!\hat{\h}{}^{min}_{\l,\,'\!\l}(q,t,a)/(1-t)^2$ 
holds.
One has: ${}_\xi^{\!\Box}\!\hat{\h}{}^{min}_{\l,\,'\!\l}(1,t,a)=
(1+a)^2(1+a+t)$. 
\smallskip

{\bf Acknowledgements.}
We are grateful to  Mikhail Khovanov, Andrei Negut, 
David Rose for useful discussions, Semen Artamonov 
for our using his software for calculating colored 
HOMFLY-PT polynomials, and the referee for important
suggestions. Special thanks to Aaron Lauda
for his attention to this work and helpful discussions
on paper \cite{MoS}.
The first author thanks Andras Szenes and University of
Geneva for the invitation and hospitality. 
I.D. acknowledges partial support from the RFBR grants 13-02-00478, 
14-02-31446-mol-a, NSh-1500.2014.2 and the common grant 
14-01-92691-Ind-a.

\vskip -2cm
\bibliographystyle{unsrt}

\medskip
\end{document}